
\input amstex
\documentstyle{amsppt}
\loadbold

\magnification=\magstep1
\hsize=6.5truein
\vsize=9truein

\document

\baselineskip=14pt

\font \smallrm=cmr10 at 9pt
\font\smallit=cmti10 at 9pt
\font\smallsl=cmsl10 at 9pt
\font\sc=cmcsc10

\def \loongrightarrow {\relbar\joinrel\relbar\joinrel\rightarrow}
\def \llongrightarrow
{\relbar\joinrel\relbar\joinrel\relbar\joinrel\rightarrow}
\def \longtwoheadrightarrow
{\relbar\joinrel\relbar\joinrel\twoheadrightarrow}
\def \llongtwoheadrightarrow
{\relbar\joinrel\relbar\joinrel\relbar\joinrel\twoheadrightarrow}
\def \gerg {\frak g}
\def \gere {\frak e}
\def \gerh {\frak h}
\def \gerk {\frak k}
\def \germ {\frak m}
\def \gersl {{\frak{sl}}}
\def \und1 {\underline{1}}
\def \undepsilon {\underline{\epsilon}}
\def \u {{\text{\bf u}}}
\def \hyp {{\text{\sl Hyp}\,}}
\def \h {\hbar}
\def \a {{\bold a}}
\def \boldalpha {{\boldsymbol\alpha}}
\def \b {{\bold b}}
\def \x {{\bold x}}
\def \id {\text{\rm id}}
\def \dif {{\text{dif}}}
\def \G {\Cal{G}}
\def \L {\Cal{L}}
\def \U {\Cal{U}}
\def \calM {\Cal{M}}
\def \calA {\Cal{A}}
\def \calB {\Cal{B}}
\def \calC {\Cal{C}}
\def \calH {\Cal{H}}
\def \calK {\Cal{K}}
\def \calX {\Cal{X}}
\def \HA {\Cal{HA}}
\def \A {\Bbb{A}}
\def \H {\Bbb{H}}
\def \N {\Bbb{N}}
\def \Z {\Bbb{Z}}
\def \C {\Bbb{C}}
\def \F {\Bbb{F}}
\def \Rhat {\widehat{R}}

\def \Char {\hbox{\it Char}\,}
\def \QrUEA {\Cal{Q}{\hskip1pt}r\Cal{UE{\hskip-1,1pt}A}}
\def \QFA {\Cal{QF{\hskip-1,7pt}A}}
\def \PrUEA {\Cal{P}{\hskip0,8pt}r\Cal{UE{\hskip-1,1pt}A}}
\def \PFA {\Cal{PF{\hskip-1,7pt}A}}


\topmatter

{\ } 

\vskip-47pt  

\hfill   {\smallrm {\smallsl preprint \  math.QA/0303019\/}  (2003)}  
\hskip19pt   {\ }  

\vskip61pt  

\title
  The global quantum duality principle:  \\
  theory, examples, and applications
\endtitle

\author
       Fabio Gavarini
\endauthor

\rightheadtext{ The global quantum duality principle: theory,
examples, applications {\ \ } }

\affil
  Universit\`a degli Studi di Roma ``Tor Vergata'' ---
Dipartimento di Matematica  \\
  Via della Ricerca Scientifica 1, I-00133 Roma --- ITALY  \\
\endaffil

\address\hskip-\parindent
  Fabio Gavarini   \newline
   \indent   Universit\`a degli Studi di Roma
``Tor Vergata''   \newline  
   \indent   Dipartimento di Matematica   \newline   
   \indent   Via della Ricerca Scientifica 1,
I-00133 Roma, ITALY   \newline   
   \indent   e-mail: gavarini\@{}mat.uniroma2.it
\endaddress

\abstract
   Let  $ R $  be an integral domain, let  $ \h \! \in \! R \! \setminus
\! \{0\} $  be such that  $ R \big/ \h R $  is a field, and  $ \HA $  the
category of torsionless (or flat) Hopf algebras over  $ R $.  We call 
$ H \! \in \HA $  a  {\it ``quantized func\-tion
algebra'' (=QFA)}, resp.~{\it ``quantized restricted universal enveloping
algebras'' (=QrUEA)}, at  $ \h \, $  if  $ \, H \big/ \h \, H \, $ 
is the function algebra of a connected Poisson group, resp.~the
(restricted, if  $ \, R \big/ \h \, R \, $  has positive characteristic)
universal enveloping algebra of a (restricted) Lie bialgebra.
%
%
                                               \par
   We establish an ``inner'' Galois correspondence on  $ \HA $,
via the definition of two endofunctors,  $ (\ )^\vee $  and
$ (\ )' $,  of  $ \HA $  such that:  {\it (a)} \, the image of
$ (\ )^\vee $,  resp.~of  $ (\ )' $,  is the full subcategory
of all QrUEAs, resp.~QFAs,  at  $ \h \, $;  {\it (b)} \, if
$ \, p := \text{\it Char}\big(R \big/ \h \, R \big) = 0 \, $,
the restrictions  $ (\ )^\vee{\big|}_{\text{QFAs}} $  and
$ (\ )'{\big|}_{\text{QrUEAs}} $  yield equivalences inverse
to each  other;  {\it (c)} \, if  $ \, p = 0 \, $,  starting from
a QFA over a Poisson group  $ G $,  resp.~from a QrUEA over a Lie
bialgebra  $ \gerg $,  the functor  $ (\ )^\vee $,  resp.~$ (\ )' $,
gives a QrUEA, resp.~a QFA, over the dual Lie bialgebra, resp.~a dual
Poisson group.  In particular,  {\it (a)\/}  yields a machine to
produce quantum groups of both types (either QFAs or QrUEAs),  {\it
(b)\/}  gives a characterization of them among objects of  $ \HA $,
and  {\it (c)\/}  gives a ``global'' version of the so-called
``quantum duality principle'' (after Drinfeld's, cf.~[Dr]).
                                               \par
   We then apply our result to Hopf algebras of the form  $ \, \Bbbk[\h]
\otimes_\Bbbk H \, $  where  $ H $  is a Hopf algebra over the field
$ \Bbbk $:  this yields quantum groups, hence ``classical'' geometrical
symmetries of Poisson type (Poisson groups or Lie bialgebras, via
specialization) associated to the ``generalized'' symmetry encoded
by  $ H $.  Both our main result and the above mentioned application
are illustrated by means of several examples, which are studied in
some detail.
\endabstract

\endtopmatter

\footnote""{Keywords: \ {\sl Hopf algebras, Quantum Groups}.}

\footnote""{ 2000 {\it Mathematics Subject Classification:} \
Primary 16W30, 17B37, 20G42; Secondary 81R50. }


\hfill  \hbox{\vbox{ \baselineskip=10pt
                     \hbox{\smallit  \   ``Dualitas dualitatum }
                     \hbox{\smallit \ \ \;\, et omnia dualitas'' }
                     \vskip4pt
                     \hbox{\smallsl    N.~Barbecue, ``Scholia'' } }
\hskip1truecm }

\vskip3pt


\centerline{ \bf INDEX }

\vskip3pt

  \indent  \phantom{\S~1} \  Introduction   \dotfill
\ pag.\ \  2
\vskip1pt

  \indent \S~1 \  Notation and terminology   \dotfill
\ pag.\ \  5
\vskip1pt

  \indent \S~2 \  The global quantum duality principle
\dotfill  \ pag.\ \  7
\vskip1pt

  \indent  \S~3 \  General properties of Drinfeld's functors
\dotfill  \ pag.\ \  10
\vskip1pt

  \indent  \S~4 \  Drinfeld's functors on quantum groups
\dotfill  \ pag.\  19
\vskip1pt

  \indent  \S~5 \  Application to trivial deformations:
the Crystal Duality Principle
\dotfill  \ pag.\  30   
\vskip1pt

  \indent  \S~6 \  First example: the Kostant-Kirillov structure
\dotfill  \ pag.\  55  
\vskip1pt

  \indent  \S~7 \  Second example: quantum  $ {SL}_2 $,  $ {SL}_n $,
finite and affine Kac-Moody groups   \dotfill  \ pag.\  62  
\vskip1pt

  \indent  \S~8 \  Third example: quantum three-dimensional Euclidean
group   \dotfill  \ pag.\  75  
\vskip1pt

  \indent  \S~9 \  Fourth example: quantum Heisenberg group
\dotfill  \ pag.\  82  
\vskip1pt

  \indent  \S~10 \, Fifth example:
   \hbox{non-commutative Hopf algebra of formal diffeomorphisms
\dotfill  \ pag.\  89}  

\vskip2,1truecm

\centerline {\bf Introduction }

\vskip10pt

   The most general notion of ``symmetry'' in mathematics is encoded
in the notion of Hopf algebra.  Among Hopf algebras  $ H $  over a
field, the commutative and the cocommutative ones encode ``geometrical''
symmetries, in that they correspond, under some technical conditions, to
algebraic groups and to (restricted, if the ground field has positive
characteristic) Lie algebras respectively: in the first case  $ H $
is the algebra  $ F[G] $  of regular functions over an algebraic group
$ G $,  whereas in the second case it is the (restricted) universal
enveloping algebra  $ U(\gerg) $  ($ \u(\gerg) $  in the restricted case)
of a (restricted) Lie algebra  $ \gerg \, $.  A popular generalization
of these two types of ``geometrical symmetry'' is given by quantum groups:
roughly, these are Hopf algebras  $ H $  depending on a parameter  $ \h $
such that setting  $ \h = 0 \, $  the Hopf algebra one gets is either of
the type  $ F[G] $   --- hence  $ H $  is a  {\sl quantized function
algebra}, in short QFA ---   or of the type  $ U(\gerg) $  or
$ \u(\gerg) $  (according to the characteristic of the ground field)
--- hence  $ H $  is a  {\sl quantized (restricted) universal enveloping
algebra}, in short QrUEA.  When a QFA exists whose specialization
(i.e.~its ``value'' at  $ \, \h = 0 \, $)  is  $ F[G] $,  the group
$ G $  inherits from this ``quantization'' a Poisson bracket, which
makes it a Poisson (algebraic) group; similarly, if a QrUEA exists
whose specialization is  $ U(\gerg) $  or  $ \u(\gerg) $,  the
(restricted) Lie algebra  $ \gerg $  inherits a Lie cobracket
which makes it a Lie bialgebra.  Then by Poisson group theory
one has Poisson groups  $ G^* $  dual to  $ G $  and a Lie bialgebra
$ \gerg^* $  dual to  $ \gerg \, $,  so other geometrical symmetries
are related to the initial ones.
                                            \par
   The dependence of a Hopf algebra on  $ \h $  can be described
as saying that it is defined over a ring  $ R $  and  $ \, \h \in
R \, $:  \, so one is lead to dwell upon the category  $ \HA $  of
Hopf  $ R $--algebras  (maybe with some further conditions), and
then raises three basic questions:

\vskip3pt

   {\bf --- (1)}  {\sl How can we produce quantum groups?}
                                            \par
   {\bf --- (2)}  {\sl How can we characterize quantum groups
(of either kind) within  $ \HA $?}
                                            \par
   {\bf --- (3)}  {\sl What kind of relationship, if any, does
exist between quantum groups over mutually dual Poisson groups,
or mutually dual Lie bialgebras?}

\vskip3pt

   A first answer to question  {\bf (1)}  and  {\bf (3)}  together
is given, in characteristic zero, by the so-called ``quantum duality
principle'', known in literature with at least two formulations.  One
claims that quantum function algebras associated to dual Poisson groups
can be taken to be dual   --- in the Hopf sense ---   to each other;
and similarly for quantum enveloping algebras (cf.~[FRT1] and [Se]).
The second one, formulated by Drinfeld in local terms (i.e., using
formal groups, rather than algebraic groups, and complete topological
Hopf algebras; cf.~[Dr], \S 7, and see [Ga4] for a proof), gives a
recipe to get, out of a QFA over  $ G $,  a QrUEA over  $ \gerg^* $,
and, conversely, to get a QFA over  $ G^* $  out of a QrUEA over
$ \gerg \, $.  More precisely, Drinfeld defines two functors, inverse
to each other, from the category of quantized universal enveloping
algebras (in his sense) to the category of quantum formal series
Hopf algebras (his analogue of QFAs) and viceversa, such that  $ \,
U_\h(\gerg) \mapsto F_\h[[G^*]] \, $  and  $ \, F_\h[[G]] \mapsto
U_\h(\gerg^*) \, $  (in his notation, where the subscript  $ \h $
stands as a reminder for ``quantized'' and the double square
brackets stand for ``formal series Hopf algebra'').
                                            \par
   In this paper we establish a  {\sl global\/}  version of the quantum
duality principle which gives a complete answer to questions  {\bf (1)}
through  {\bf (3)}.  The idea is to push as far as possible Drinfeld's
original method so to apply it to the category  $ \HA $  of all Hopf
algebras which are torsion-free   --- or flat, if one prefers this
narrower setup ---   modules over some (integral)
domain, say  $ R $,  and to do it for each non-zero element  $ \h $ 
in  $ R $  such that  $ R \big/ \h \, R $  be a field.
                                            \par
   To be precise, we extend Drinfeld's recipe so to define functors
from  $ \HA $  to itself.  We show that the image of these
``generalized'' Drinfeld's functors is contained in a category
of quantum groups   --- one gives QFAs, the other QrUEAs ---   so
we answer question  {\bf (1)}.  Then, in the characteristic zero case,
we prove that when restricted to quantum groups these functors yield
equivalences inverse to each other.  Moreover, we show that these
equivalences exchange the types of quantum group (switching QFA with
QrUEA) and the underlying Poisson symmetries (interchanging  $ G $
or  $ \gerg $  with  $ G^* $  or  $ \gerg^* $),  thus solving  {\bf
(3)}.  Other details enter the picture to show that these functors
endow  $ \HA $  with sort of a (inner) ``Galois correspondence'',
in which QFAs on one side and QrUEAs on the other side are the
subcategories (in  $ \HA $)  of ``fixed points'' for the composition
of both Drinfeld's functors (in the suitable order): in particular,
this answers question  {\bf (2)}.  It is worth stressing that, since
our ``Drinfeld's functors'' are defined for each non-trivial point
$ (\h) $  of  $ \text{\it Spec}_{\,\text{\it max}}\,(R) $,  \, for
any such  $ (\h) $  and for any  $ H $  in  $ \HA $  they yield two
quantum groups, namely a QFA and a QrUEA,  w.r.t.~$ \h $  itself.
   \hbox{Thus we have a method to get, out of  any single
$ H \hskip-0pt \in \hskip-0pt \HA \, $, several quantum groups.}
                                             \par
   Therefore the ``global'' in the title is meant in several respects:
geometrically, we consider global objects (Poisson groups rather
than Poisson  {\sl formal\/}  groups, as in Drinfeld's approach);
algebraically we consider quantum groups over any domain  $ R $, 
so there may be several different ``semiclassical limits''
(=specializations) to consider, one for each non-trivial point
of type  $ (\h) $  in the maximal spectrum of  $ R $  (while Drinfeld
has  $ \, R = \Bbb\Bbbk[[\h]] \, $  so one can specialise only at  $ \,
\h = 0 \, $);  more generally, our recipe
applies to  {\sl any\/}  Hopf algebra, i.e.~not only to quantum groups;
finally, most of our results are characteristic-free, i.e.~they hold
not only in characteristic zero (as in Drinfeld's case) but also in
positive characteristic.  Furthermore, this ``global version'' of the
quantum duality principle opens the way to formulate a ``quantum
duality principle for subgroups and homogeneous spaces'', see [CG].
                                             \par
   A key, long-ranging application of our  {\sl global quantum duality
principle\/}  (GQDP) is the following.  Take as  $ R $  the polynomial
ring  $ \, R = \Bbbk[\h\,] \, $,  where  $ \Bbbk $  is a field: then
for any Hopf algebra over  $ \Bbbk $  we have that  $ \, H[\h\,] := R
\otimes_\Bbbk H \, $  is a torsion-free Hopf  $ R $--algebra,  hence
we can apply Drinfeld's functors to it.  The outcome of this procedure
is the  {\sl crystal duality principle\/}  (CDP), whose statement strictly
resembles that of the GQDP: now Hopf  $ \Bbbk $--algebras  are looked at
instead of torsionless Hopf  $ R $--algebras,  and quantum groups are
replaced by Hopf algebras with canonical filtrations such that the
associated graded Hopf algebra is either commutative or cocommutative.
Correspondingly, we have a method to associate to  $ H $  a Poisson
group  $ G $  and a Lie bialgebra  $ \gerk $  such that  $ G $  is an
affine space (as an algebraic variety) and  $ \gerk $  is graded (as
a Lie algebra); in both cases, the ``geometrical'' Hopf algebra can be
attained   --- roughly ---   through a continuous 1-parameter deformation
process.  This result can also be formulated in purely classical   ---
i.e.~``non-quantum'' ---   terms and proved by purely classical means.
However, the approach via the GQDP also yields further possibilities to
deform  $ H $  into other Hopf algebras of geometrical type, which is
out of reach of any classical approach.
                                             \par
   The paper is organized as follows.  In \S 1 we fix notation and
terminology, while \S 2 is devoted to define Drinfeld's functors and
state our main result, the GQDP (Theorem 2.2).  In \S 3 we extend
Drinfeld's functors to a broader framework, that of  {\sl (co)augmented
(co)algebras},  and study their properties in general.  \S 4 instead is
devoted to the analysis of the effect of such functors on quantum groups,
and prove Theorem 2.2, i.e.~the GQDP.  In \S 5 we explain the CDP, which
is deduced as an application of the CDP to trivial deformations of Hopf
$ \Bbbk $--algebras:  in particular, we study in detail the case of group
algebras.  In the last part of the paper we illustrate our results by
studying in full detail several relevant examples.  First we dwell upon
some well-known quantum groups: the standard quantization of the
Kostant-Kirillov structure on a Lie algebra (\S 6), the standard
Drinfeld-Jimbo quantization of semisimple groups (\S 7), the quantization
of the Euclidean group (\S 8) and the quantization of the Heisenberg group
(\S 9).  Then we study a key example of non-commutative, non-cocommutative
Hopf algebra   --- a non-commutative version of the Hopf algebra of formal
diffeomorphisms ---   as a nice application of the CDP (\S 10).

\vskip3pt

   {\sl $ \underline{\hbox{\it Warning}} \, $:  \, this paper is
$ \underline{\hbox{\it not\/}}  $  meant for publication!  The
results presented here will be published in separate articles;
therefore, any reader willing to quote anything from the present
preprint is kindly invited to ask the author for the precise
reference(s).}

\vskip11pt

\centerline{ \sc acknowledgements }

  The author thanks P.~Baumann, G.~Carnovale, N.~Ciccoli, A.~D'Andrea,
I.~Damiani, B.~Di Blasio, L.~Foissy, A.~Frabetti, C.~Gasbarri and
E.~Taft for many helpful discussions.

\vskip1,9truecm

\centerline {\bf \S \; 1 \ Notation and terminology }

\vskip10pt

  {\bf 1.1 The classical setting.} \, Let  $ \Bbbk $  be a fixed
field of any characteristic.  We call ``algebraic group'' the
maximal spectrum  $ G $  associated to any commutative Hopf
$ \Bbbk $--algebra  $ H $  (in particular, we deal with  {\sl
proaffine\/}  as well as  {\sl affine\/}  algebraic groups); then
$ H $  is called the algebra of regular function on  $ G $,  denoted
with  $ F[G] $.  We say that  $ G $  is connected if  $ F[G] $  has
no non-trivial idempotents; this is equivalent to the classical
topological notion when  $ F[G] $  is of finite type, i.e.~$ \dim(G) $
is finite.

  If  $ G $  is an algebraic group, we denote by
$ \germ_e $  the defining ideal of the unit element  $ \, e
\in G \, $  (it is the augmentation ideal of  $ F[G] \, $);
the cotangent space of  $ G $  at  $ e $  is  $ \, \gerg^\times
:= \germ_e \Big/ {\germ_e}^{\!2} \, $,  \, 
%
%
which is naturally a Lie coalgebra.  The tangent space of  $ G $ 
at  $ e $  is the dual space  $ \, \gerg := {\big( \gerg^\times
\big)}^* \, $  to  $ \gerg^\times \, $:  \, this is a Lie algebra,
which coincides with the set of all left-(or right-)invariant 
derivations of  $ F[G] \, $.  By  $ U(\gerg) $  we mean the
universal enveloping algebra of  $ \gerg \, $:  it is a connected
cocommutative Hopf algebra, and there is a natural Hopf pairing
(see  \S 1.2{\it (a)\/})  between  $ F[G] $  and  $ U(\gerg) $.
If  $ \, \hbox{\it Char} \,(\Bbbk) = p > 0 \, $,  \, then
$ \gerg $  is a restricted Lie algebra, and  $ \, \u(\gerg)
:= U(\gerg) \Big/ \big( \big\{\, x^p - x^{[p\hskip0,7pt]}
\,\big|\, x \in \gerg \,\big\} \big) \, $  is the restricted
universal enveloping algebra of  $ \gerg \, $.  In the sequel,
to unify notation and terminology, when  $ \, \Char(\Bbbk) = 0 \, $ 
we shall call any Lie algebra  $ \gerg $  ``restricted'', and its
``restricted universal enveloping algebra'' will be  $ U(\gerg) $, 
and we shall write  $ \, \U(\gerg) := U(\gerg) \, $  if  $ \,
\Char(\Bbbk) = 0 \, $  and  $ \, \U(\gerg) := \u(\gerg) \, $ 
if  $ \, \Char(\Bbbk) > 0 \, $.
                                            \par
   We shall also consider  $ \, \hyp(G) := {\big( {F[G]}^\circ
\big)}_\varepsilon = \big\{\, f \in {F[G]}^\circ \,\big|\,
f({\germ_e}^{\!n}) = 0 \; \forall \, n \geqq 0 \,\big\} \, $,  \,
i.e.~the connected component of the Hopf algebra  $ {F[G]}^\circ $
dual to  $ F[G] $,  which is called the  {\sl hyperalgebra\/}  of
$ G $.  By construction  $ \hyp(G) $  is a connected Hopf algebra,
containing  $ \, \gerg = \text{\sl Lie\/}(G) \, $;  if  $ \, \text{\it
Char}\,(\Bbbk) = 0 \, $  one has  $ \, \hyp(G) = U(\gerg) \, $,  \,
whereas if  $ \, \text{\it Char}\,(\Bbbk) > 0 \, $  one has a sequence
of Hopf algebra morphisms  $ \; U(\gerg) \longtwoheadrightarrow \u(\gerg)
\; {\lhook\joinrel\relbar\joinrel\relbar\joinrel\rightarrow} \, \hyp(G)
\; $.  In any case, there exists a natural perfect (= non-degenerate)
Hopf pairing between  $ F[G] $  and  $ \hyp(G) $.
                                            \par
   Now assume  $ G $  is a Poisson group (for this and other notions
hereafter see, e.g., [CP], but within an  {\sl algebraic geometry\/}
setting): then  $ F[G] $  is a Poisson Hopf algebra, and its Poisson
bracket induces on  $ \gerg^\times $  a Lie bracket which makes it
into a Lie bialgebra, and so  $ U(\gerg^\times) $  and 
$ \U(\gerg^\times) $  are co-Poisson Hopf algebras too.  On the
other hand,  $ \gerg $  turns into a Lie bialgebra   --- maybe in
topological sense, if  $ G $  is infinite dimensional ---   and 
$ U(\gerg) $  and  $ \U(\gerg) $  are (maybe topological) co-Poisson
Hopf algebras.  The Hopf pairing above between  $ F[G] $  and 
$ \U(\gerg) $  then is compatible with these additional co-Poisson
and Poisson structures.  Similarly,  $ \hyp(G) $  is a co-Poisson Hopf
algebra as well and the Hopf pairing between  $ F[G] $  and  $ \hyp(G) $ 
is compatible with the additional structures.  Moreover, the perfect
(=non-degenerate) pairing  $ \, \gerg \times \gerg^\times \!
\longrightarrow \Bbbk \, $  given by evaluation is compatible
with the Lie bialgebra structure on either side (see  \S 1.2{\it
(b)\/}):  so  $ \gerg $  and  $ \gerg^\times $  are Lie bialgebras 
{\sl dual to each other}.  In the sequel, we denote by  $ G^\star $ 
any connected algebraic Poisson group with  $ \gerg $  as cotangent
Lie bialgebra, and say it is  {\sl dual\/}  to  $ G \, $.
                                            \par
   Let  $ H $  be a Hopf algebra over an integral domain  $ D \, $. 
We call  $ H $  a  {\sl ``function algebra''\/}  (FA in short) if it
is commutative, with no non-trivial idempotents, and such that, if 
$ \, p := \Char(\Bbbk) > 0 \, $,  \, then  $ \, \eta^p = 0 \, $  for
all  $ \eta $  in the kernel of the counit of  $ H \, $.  If  $ D $ 
is a field, an FA is the algebra of regular functions of an algebraic
group-scheme over  $ D $  which is connected and, if  $ \, \Char(\Bbbk)
> 0 \, $,  \, is zero-dimensional of height 1; conversely, if  $ G $  is
such a group-scheme then  $ F[G] $  has these properties.  Instead, we
call  $ H $  a ``restricted universal enveloping algebra'' (=rUEA) if it
is cocommutative, connected, and generated by its primitive part.  If 
$ D $  is a field, an rUEA is the restricted universal enveloping algebra
of some (restricted) Lie algebra over  $ D \, $;  \, conversely, if 
$ \gerg $  is such a Lie algebra, then  $ \U(\gerg) $  has these
properties   
%
%
(see, e.g., [Mo], Theorem 5.6.5, and references therein).   
                                            \par
%
%
   For the Hopf operations in any Hopf algebra we shall use standard
notation, as in [Ab].

\vskip7pt

\proclaim {Definition 1.2}
                                   \hfill\break
   \indent   (a) \, Let  $ H $,  $ K $  be Hopf algebras (in any
category).  A pairing  $ \; \langle \,\ , \,\ \rangle \, \colon \, H
\times K \loongrightarrow R \; $  (where  $ R $  is the ground ring)
is a  {\sl Hopf (algebra) pairing\/}  if  $ \;\; \big\langle x, y_1
\cdot y_2 \big\rangle = \big\langle \Delta(x), y_1 \otimes y_2
\big\rangle := \sum_{(x)} \big\langle x_{(1)}, y_1 \big\rangle \cdot
\big\langle x_{(2)}, y_2 \big\rangle \, $,  $ \; \big\langle x_1
\cdot x_2, y \big\rangle = \big\langle x_1 \otimes x_2, \Delta(y)
\big\rangle := \sum_{(y)} \big\langle x_1, y_{(1)} \big\rangle
\cdot \big\langle x_2, y_{(2)} \big\rangle  \, $,  $ \, \langle
x, 1 \rangle = \epsilon(x) \, $,  $ \; \langle 1, y \rangle =
\epsilon(y) \, $,  $ \; \big\langle S(x), y \big\rangle =
\big\langle x, S(y) \big\rangle \, $,  for all  $ \, x, x_1,
x_2 \in H $,  $ \, y, y_1, y_2 \in K $.
                                   \hfill\break
   \indent   (b) \, Let  $ \gerg $,  $ \gerh $  be Lie bialgebras (in
any category).  A pairing  $ \; \langle \,\ , \,\ \rangle \, \colon \,
\gerg \times \gerh \loongrightarrow \Bbbk \; $  (where  $ \Bbbk $  is the
ground ring) is  called a  {\sl Lie bialgebra pairing\/}  if  $ \;\;
\big\langle x, [y_1,y_2] \big\rangle = \big\langle \delta(x), y_1
\otimes y_2 \big\rangle := \sum_{[x]} \big\langle x_{[1]}, y_1
\big\rangle \cdot \big\langle x_{[2]}, y_2 \big\rangle  \, $,
$ \; \big\langle [x_1,x_2], y \big\rangle = \big\langle x_1 \otimes
x_2, \delta(y) \big\rangle := \sum_{[y]} \big\langle x_1, y_{[1]}
\big\rangle \cdot \big\langle x_2, y_{[2]} \big\rangle \, $,  \,
for all  $ \, x, x_1, x_2 \in \gerg \, $  and  $ \, y, y_1, y_2
\in \gerh $,  \, with  $ \, \delta(x) = \sum_{[x]} x_{[1]} \otimes
x_{[2]} \, $  and  $ \, \delta(y) = \sum_{[x]} y_{[1]} \otimes
y_{[2]} \, $.
\endproclaim

\vskip7pt

  {\bf 1.3 The quantum setting.} \, Let  $ R $  be a (integral) domain,
and let  $ \, F = F(R) \, $  be its quotient field.  Denote by  $ \calM $ 
the category of torsion-free  $ R $--modules,  and by  $ \HA $  the
category of all Hopf algebras in  $ \calM \, $;  \, note that  {\sl
flat\/}  modules form a full subcategory of  $ \calM \, $.  Let 
$ \calM_F $  be the category of  $ F $--vector  spaces, and  $ \HA_F $ 
be the category of all Hopf algebras in  $ \calM_F \, $.  For any  $ \,
M \in \calM \, $,  \, set  $ \, M_F := F(R) \otimes_R M \, $.  Scalar
extension gives a functor  $ \; \calM \longrightarrow \calM_F \, $, 
$ \, M \mapsto M_F \, $,  \, which restricts to a functor  $ \; \HA
\longrightarrow \HA_F \, $  as well.
                                            \par
   Let  $ \, \h \in R \, $  be a non-zero element (which will be fixed
throughout), and let  $ \, \Bbbk := R \big/ (\h) = R \big/ \h \, R \, $ 
be the quotient ring.  For any  $ R $--module  $ M $,  we set  $ \,
M_\h{\Big|}_{\h=0} \!\! := M \big/ \! \h \, M = \Bbbk \otimes_R M \, $: 
this is a  $ \Bbbk $--module  (via scalar restriction  $ \, R \rightarrow
R \big/ \h \, R =: \Bbbk \, $),  which we call the  {\sl specialization\/} 
of  $ M $  at  $ \, \h = 0 \, $;  we use also notation  $ \, M
\,{\buildrel \, \h \rightarrow 0 \, \over \llongrightarrow}\,
\overline{M} \, $  to mean shortly that  $ \, M_\h{\Big|}_{\h=0}
\hskip-3pt \cong \overline{M} \, $.  Moreover, set  $ \, M_\infty
:= \bigcap_{n=0}^{+\infty} \h^n M \, $
    \hbox{(this is the closure of  $ \{0\} $  in the  $ \h $--adic}
topology of  $ M $).
   For any  $ H \! \in \! \HA \, $,  let  $ \, I_{\!
\scriptscriptstyle H} \! := \! \text{\sl Ker} \Big( H
\,{\buildrel \epsilon \over \twoheadrightarrow}\,
R \,{\buildrel {\h \mapsto 0} \over
{\relbar\joinrel\relbar\joinrel\twoheadrightarrow}}\,
\Bbbk \Big) \, $
       \hbox{and set  $ \, {I_{\!
\scriptscriptstyle H}}^{\!\!\infty} \!\! :=
\bigcap_{n=0}^{+\infty} {I_{\! \scriptscriptstyle H}}^{\!\!n} $.}
                                            \par
   Finally, given  $ \Bbb{H} $  in  $ \HA_F $,  a subset
$ \overline{H} $  of  $ \Bbb{H} $  is called  {\it an
$ R $--integer  form}  (or simply  {\it an  $ R $--form})  {\it of}
$ \, \Bbb{H} $  iff  $ \; \overline{H} \, $  is a Hopf  $ R $--subalgebra
of  $ \, \Bbb{H} \, $  (so  $ \, \overline{H} \, $  is torsion-free as an
$ R $--module,  hence  $ \, \overline{H} \in \HA \, $)  and  $ \; H_F :=
F(R) \otimes_R \overline{H} = \Bbb{H} \, $.

\vskip3pt

   We are now ready to introduce the notion of ``quantum group''.

\vskip7pt

\proclaim {Definition 1.4} \!\!\! ({\sl ``Global quantum groups'' [or
``algebras'']})  Let  $ R $,  
   \hbox{$ \h \! \in \! R \! \setminus \! \{0\} $ 
\! be as in \S 1.3.}
                                          \par
   (a) \, We call  {\sl quantized restricted universal enveloping
algebra}   (in short,  {\sl QrUEA\/})  {\sl (at  $ \h $)}  any
$ \U_\h \in \HA \, $  such that  $ \, \U_\h{\big|}_{\h=0} \! :=
\U_\h \big/ \h \, \U_\h \, $  is (isomorphic to) an rUEA.
                                          \par
   We call  $ \, \QrUEA \, $  the full subcategory
of  $ \HA $  whose objects are all the QrUEAs (at  $ \h $).
                                          \par
   (b) \, We call  {\sl quantized function algebra}  (in short, 
{\sl QFA\/})  {\sl (at  $ \h $)}  any  $ \, F_\h \in \HA $
such that  $ \, {(F_\h)}_\infty = {I_{\!\scriptscriptstyle
F_\h}}^{\!\!\infty} \, $
            (notation of \S 1.3)\footnote{This requirement turns
                                          out to be a natural one,
                                          see Theorem 3.8.}
and  $ \, F_\h{\big|}_{\h=0} \! := F_\h \big/ \h \, F_\h \, $
is (isomorphic to) an FA.
                                          \par
   We call  $ \, \QFA \, $  the full subcategory
of  $ \HA $  whose objects are all the QFAs (at  $ \h $).
\endproclaim

\vskip5pt

  {\bf Remark 1.5:} \, If  $ \, \U_\h \, $  is a QrUEA (at  $ \h \, $,
that is w.r.t.~to  $ \h \, $)  then  $ \, \U_\h{\big|}_{\h=0} \, $  is
a co-Poisson Hopf algebra, w.r.t.~the Poisson cobracket  $ \delta $
defined as follows: if  $ \, x \in \U_\h{\big|}_{\h=0} \, $  and  $ \,
x' \in \U_\h \, $  gives  $ \, x = x' \mod \h \, \U_\h \, $,  \,
then  $ \, \delta(x) := \big( \h^{-1} \, \big( \Delta(x') -
\Delta^{\text{op}}(x') \big) \big) \mod \h \, \big( \U_\h \otimes
\U_\h \big) \, $.  If  $ \, \Bbbk := R \big/ \h \, R \, $  is a field,
then  $ \, \U_\h{\big|}_{\h=0} \cong \U(\gerg) \, $  for some Lie
algebra  $ \gerg \, $,  and by [Dr], \S 3, the restriction of 
$ \delta $  makes  $ \gerg $  into a  {\sl Lie bialgebra\/}  (the
isomorphism  $ \, \U_\h{\big|}_{\h=0} \cong \U(\gerg) \, $  being
one of  {\sl co-Poisson\/}  Hopf algebras); in this case we write 
$ \, \U_\h = \U_\h(\gerg) \, $.   
                                                  \par   
   Similarly, if  $ F_\h $  is a QFA at $ \h $,  then  $ \,
F_\h{\big|}_{\h=0} \, $  is a  {\sl Poisson\/}  Hopf algebra,
w.r.t.~the Poisson bracket  $ \{\,\ ,\ \} $  defined as follows:
if  $ \, x $,  $ y \in F_\h{\big|}_{\h=0} \, $  and  $ \, x' $, 
$ y' \in F_\h \, $  give  $ \, x = x' \mod \h \, F_\h \, $,  $ \,
y = y' \mod \h \, F_\h \, $,  \, then  $ \, \{x,y\} := \big( \h^{-1}
(x' \, y' - y' \, x') \big) \mod \h \, F_\h \, $.  Therefore , if 
$ \, \Bbbk := R \big/ \h \, R \, $  is a field, then  $ \, F_\h
{\big|}_{\h=0} \cong F[G] \, $  for some connected  {\sl Poisson\/} 
algebraic group  $ G $  (the isomorphism being one of  {\sl Poisson\/} 
Hopf algebras): in this case we write  $ \, F_\h = F_\h[G] \, $.   

\vskip7pt

\proclaim{Definition 1.6}
                                           \hfill\break
   \indent   (a) \, Let  $ R $  be any (integral) domain, and let
$ F $  be its field of fractions.  Given two  $ F $--modules
$ \Bbb{A} $,  $ \Bbb{B} $,  and an  $ F $--bilinear  pairing
$ \, \Bbb{A} \times \Bbb{B} \longrightarrow F \, $,  for any
$ R $--submodule  $ \, A \subseteq \Bbb{A} \, $  and  $ \, B
\subseteq \Bbb{B} \, $  we set  $ \; \displaystyle{ A^\bullet
\, := \Big\{\, b \in \Bbb{B} \;\Big\vert\; \big\langle A, \, b
\big\rangle \subseteq R \Big\} } \; $  and  $ \; \displaystyle{
B^\bullet \, := \Big\{\, a \in \Bbb{A} \;\Big\vert\; \big\langle a,
B \big\rangle \subseteq R \Big\} } \, $.
                                           \hfill\break
   \indent   (b) \, Let  $ R $  be a domain.  Given  $ \, H $,  $ K \in
\HA \, $,  \, we say that  {\sl  $ H $  and  $ K $  are dual to each
other}  if there exists a perfect Hopf pairing between them for which
$ \, H = K^\bullet \, $  and  $ \, K = H^\bullet \, $.
\endproclaim

\vskip1,9truecm

\centerline {\bf \S \; 2 \  The global quantum duality principle }

\vskip10pt

  {\bf 2.1 Drinfeld's functors.} \,  (Cf.~[Dr], \S 7) Let  $ R $,
$ \HA $  and  $ \, \h \in R \, $  be as in \S 1.3.  For any
$ \, H \in \HA \, $,  \, let  $ \, I = I_{\scriptscriptstyle
\! H} := \hbox{\sl Ker} \Big( H \,{\buildrel \epsilon \over
{\relbar\joinrel\twoheadrightarrow}}\, R \,{\buildrel {\h \mapsto 0}
\over \llongtwoheadrightarrow}\, R \big/ \h \, R = \Bbbk \Big) =
\hbox{\sl Ker} \Big( H \,{\buildrel {\h \mapsto 0} \over
\llongtwoheadrightarrow}\, H \big/ \h \, H \,{\buildrel
\bar{\epsilon} \over {\relbar\joinrel\twoheadrightarrow}}\,
\Bbbk \Big) \, $  (as in \S 1.3), a maximal Hopf ideal of  $ H $
(where  $ \bar{\epsilon} $  is the counit of  $ \, H{\big|}_{\h=0}
\, $,  \, and the two composed maps clearly coincide): we define
  $$  H^\vee \; := \; {\textstyle \sum_{n \geq 0}} \, \h^{-n} I^n \;
= \; {\textstyle \sum_{n \geq 0}} \, {\big( \h^{-1} I \, \big)}^n \;
= \; {\textstyle \bigcup_{n \geq 0}} \, {\big( \h^{-1} I \, \big)}^n
\quad  \big( \! \subseteq H_F \, \big) \; .  $$
If  $ \, J = J_{\scriptscriptstyle \! H} := \hbox{\sl Ker}\,
(\epsilon_{\scriptscriptstyle \! H}) \, $  then  $ \, I = J +
\h \cdot 1_{\scriptscriptstyle H} \, $,  \, so  $ \; H^\vee =
\sum_{n \geq 0} \h^{-n} J^n = \sum_{n \geq 0} {\big( \h^{-1}
J \, \big)}^n \, $  too.   
%
%
 \eject   
  Given any Hopf algebra  $ H $,  for every  $ \, n \in \N \, $  define
$ \; \Delta^n \colon H \longrightarrow H^{\otimes n} \; $  by  $ \,
\Delta^0 := \epsilon \, $,  $ \, \Delta^1 := \id_{\scriptscriptstyle
H} $,  \, and  $ \, \Delta^n := \big( \Delta \otimes
\id_{\scriptscriptstyle H}^{\otimes (n-2)} \big)
\circ \Delta^{n-1} \, $  if  $ \, n > 2 $.  For any
ordered subset  $ \, \Sigma = \{i_1, \dots, i_k\} \subseteq
\{1, \dots, n\} \, $  with  $ \, i_1 < \dots < i_k \, $,  \, define
the morphism  $ \; j_{\scriptscriptstyle \Sigma} : H^{\otimes k}
\longrightarrow H^{\otimes n} \; $  by  $ \; j_{\scriptscriptstyle
\Sigma} (a_1 \otimes \cdots \otimes a_k) := b_1 \otimes \cdots
\otimes b_n \; $  with  $ \, b_i := 1 \, $  if  $ \, i \notin
\Sigma \, $  and  $ \, b_{i_m} := a_m \, $  for  $ \, 1 \leq m
\leq k $~;  then set  $ \; \Delta_\Sigma := j_{\scriptscriptstyle
\Sigma} \circ \Delta^k \, $,  $ \, \Delta_\emptyset := \Delta^0
\, $,  and  $ \; \delta_\Sigma := \sum_{\Sigma' \subset \Sigma}
{(-1)}^{n- \left| \Sigma' \right|} \Delta_{\Sigma'} \, $,
$ \; \delta_\emptyset := \epsilon \, $.  By the inclusion-exclusion
principle, this definition admits the inverse formula  $ \; \Delta_\Sigma
= \sum_{\Psi \subseteq \Sigma} \delta_\Psi \, $.  We shall also use the
notation  $ \, \delta_0 := \delta_\emptyset \, $,  $ \, \delta_n :=
\delta_{\{1, 2, \dots, n\}} \, $,  and the useful formula  $ \; \delta_n
= {(\id_{\scriptscriptstyle H} - \epsilon)}^{\otimes n} \circ \Delta^n
\, $,  \, for all  $ \, n \in \N_+ \, $.
                                              \par
   Now consider any  $ \, H \in \HA \, $  and  $ \, \h \in R \, $  as
in \S 1.3: we define
  $$  H' \; := \; \big\{\, a \in H \,\big\vert\, \delta_n(a) \in
\h^n H^{\otimes n} ,  \; \forall \,\, n \in \N \, \big\}  \quad
\big( \! \subseteq H \, \big) \, .  $$

\vskip9pt

\proclaim {Theorem 2.2} \! ({\sl ``The Global Quantum Duality
Principle''})  Assume  $ \, \Bbbk \! := \! R \big/ \h R \, $ 
is a field.   
                                        \hfill\break
  \indent   (a) \, The assignment  $ \, H \mapsto H^\vee \, $,
resp.~$ \, H \mapsto H' \, $,  defines a functor  $ \; {(\ )}^\vee
\colon \, \HA \, \longrightarrow \, \HA \, $,  \, resp.~$ \; {(\ )}'
\colon \, \HA \, \longrightarrow \, \HA \, $,  \, whose image lies
in  $ \QrUEA \, $,  resp.~in  $ \QFA \, $.  Moreover, for all  $ \,
H \in \HA \, $  we have  $ \, H \subseteq {\big( H^\vee \big)}' \, $ 
and  $ \, H \supseteq {\big( H' \big)}^{\!\vee} \! $,  hence also 
$ \, H^\vee = \big( \big(H^\vee\big)' \,\big)^{\!\vee} $  and  $ \,
H' = \big( \big(H'\big)^{\!\vee} \big)' $.  Finally, if  $ \, H \in
\HA \, $  is flat, then  $ H^\vee $  and  $ H' $  are flat as well.
                                        \hfill\break
   \indent   (b) \, Assume that  $ \, \hbox{\it Char}\,(\Bbbk) = 0 \, $. 
Then for any  $ \, H \in \HA \, $   
 \vskip1pt
   \centerline{ $ \displaystyle{ H = {\big(H^\vee\big)}'
\,\Longleftrightarrow\, H \in \QFA  \qquad  \hbox{and}
\qquad  H = {\big( H' \big)}^{\!\vee} \,\Longleftrightarrow\,
H \in \QrUEA \quad ; } $ }
 \vskip1pt
\noindent   thus we have two induced equivalences, namely
$ \; {(\ )}^\vee \colon \, \QFA \, \llongrightarrow
\, \QrUEA \, $,  $ \, H \mapsto H^\vee \, $,  \; and
$ \; {(\ )}' \colon \, \QrUEA \, \llongrightarrow \, \QFA \, $,
$ \, H \mapsto H' \, $,  \; which are inverse to each other.
                                        \hfill\break
  \indent   (c) \, (``Quantum Duality Principle'')  Assume that 
$ \, \hbox{\it Char}\,(\Bbbk) = 0 \, $.  Then
 \vskip5pt
   \centerline{ $ {F_\h[G]}^\vee{\Big|}_{\h=0} := {F_\h[G]}^\vee \Big/
\h \, {F_\h[G]}^\vee \, = \, U(\gerg^\times) \; ,  \quad
\displaystyle{ {U_\h(\gerg)}'{\Big|}_{\h=0} := {U_\h(\gerg)}' \Big/
\h \, {U_\h(\gerg)}' \, = \, F\big[G^\star\big] } $ }
 \vskip5pt
\noindent   (with  $ G $,  $ \gerg $,  $ \gerg^\times $,
$ \gerg^\star $  and  $ G^\star $  as in \S 1.1, and  $ U_\h(\gerg) $
has the obvious meaning, cf.~\S 1.5) where the choice of the group
$ G^\star $   --- among all the connected Poisson algebraic groups with
tangent Lie bialgebra  $ \gerg^\star $  ---   depends on the choice of
the QrUEA  $ U_\h(\gerg) $.  In other words,  $ \, {F_\h[G]}^\vee \, $
is a QrUEA for the Lie bialgebra  $ \gerg^\times $,  and  $ \,
{U_\h(\gerg)}' \, $  is a QFA for the Poisson group  $ G^\star $.
                                        \hfill\break
  \indent   (d) \, Let  $ \, \text{\it Char}\,(\Bbbk) = 0 \, $.  Let
$ \, F_\h \in \QFA \, $,  $ \, U_\h \in \QrUEA \, $  be dual to each
other (with respect to some pairing).  Then  $ {F_\h}^{\!\vee} $  and
$ {U_\h}' $  are dual to each other (w.r.t.~the  {\sl same}  pairing).
                                        \hfill\break
  \indent   (e) \, Let  $ \, \hbox{\it Char}\,(\Bbbk) = 0 \, $.
Then for all  $ \, \H \in \HA_F \, $  the following are equivalent:
                                        \hfill\break
   \indent \indent  \;  $ \H $  has an  $ R $--integer  form
$ H_{(f)} $  which is a QFA at  $ \h \, $;
                                        \hfill\break
   \indent \indent  \;  $ \H $  has an  $ R $--integer  form
$ H_{(u)} $  which is a QrUEA at  $ \h \, $.
\endproclaim

\vskip9pt

{\bf Remarks 2.3:}  after stating our main theorem, some
comments are in order.
                                                 \par
   {\it (a)} \, {\sl The Global Quantum Duality Principle as
a ``Galois correspondence'' type theorem.}
                                                 \par
\noindent   Let  $ \, L \subseteq E \, $  be a Galois (not
necessarily finite) field extension, and let  $ \, G := \hbox{\it
Gal}\,\big(E/L\big) \, $  be its Galois group.  Let  $ \, \Cal{F}
\, $  be the set of intermediate extensions (i.e.~all fields
$ F $  such that  $ \, L \subseteq F \subseteq E \, $),  \, let
$ \, \Cal{S} \, $  be the set of all subgroups of  $ G $  and let
$ \, \Cal{S}^c \, $  be the set of all subgroups of  $ G $  which
are  {\sl closed\/}  w.r.t.~the Krull topology of  $ G $.  Note
that  $ \Cal{F} $,  $ \Cal{S} $  and  $ \Cal{S}^s $  can all
be seen as lattices w.r.t.~set-theoretical inclusion   ---
$ \Cal{S}^c $  being a sublattice of  $ \Cal{S} $  ---   hence as
categories too.  The celebrated Galois Theorem yields two maps,
namely  $ \; \varPhi \, \colon \, \Cal{F} \loongrightarrow \Cal{S}
\, $,  $ \, F \mapsto \hbox{\it Gal}\,\big(E/F\big) := \big\{\, \gamma
\in G \;\big\vert\;\, \gamma{\big\vert}_F = \hbox{id}_F \,\big\} \, $,
and  $ \; \varPsi \, \colon \, \Cal{S} \loongrightarrow \Cal{F} \, $,
$ \, H \mapsto E^H := \big\{\, e \in E \;\big\vert\; \eta(e) = e \;\;
\forall \; \eta \in H \,\big\} \, $,
%
%
%
%
such that:
                                      \par
   {\it --- 1)} \;  $ \varPhi $  and  $ \varPsi $  are contravariant
functors (that is, they are order-reversing maps of lattices, i.e.{}
lattice antimorphisms); moreover, the image of  $ \varPhi $  lies in
the subcategory  $ \, \Cal{S}^c \, $;
                                      \par
   {\it --- 2)} \; for  $ \, H \in \Cal{S} \, $  one has  $ \,
\varPhi\big(\varPsi(H)\big) = \overline{H} \, $,  \, the  {\sl
closure\/}  of  $ H $  w.r.t.~the Krull topology: thus  $ \, H
\subseteq \varPhi \big( \varPsi(H) \big) \, $,  \, and  $ \,
\varPhi \circ \varPsi \, $  is a  {\sl closure operator},  so
that   $ \, H \in \Cal{S}^c \, $  iff  $ \, H = \varPhi \big(
\varPsi(H) \big) \, $;
                                      \par
   {\it --- 3)} \; for  $ \, F \in \Cal{F} \, $  one has  $ \,
\varPsi\big(\varPhi(F)\big) = F \, $;
                                      \par
   {\it --- 4)} \;  $ \varPhi $  and  $ \varPsi $  restrict to
antiequivalences  $ \, \varPhi : \Cal{F} \rightarrow \Cal{S}^c \, $
and  $ \, \varPsi : \Cal{S}^c \rightarrow \Cal{F} \, $  which are
inverse to each other.
                                      \par
   Then one can see that Theorem 2.2 establishes a strikingly
similar result, which in addition is much more symmetric:
$ \HA $  plays the role of both  $ \Cal{F} $  and  $ \Cal{S} $,
whereas  $ {(\ )}' $  stands for  $ \varPsi $  and  $ {(\ )}^\vee $
stands for  $ \varPhi $.  $ \QFA $  plays the role of the
distinguished subcategory  $ \Cal{S}^c $,  and symmetrically we
have the distinguished subcategory  $ \QrUEA $.  The composed
operator  $ \, {\big({(\ )}^\vee\big)}' = {(\ )}' \circ {(\ )}^\vee
\, $  plays the role of a ``closure operator'', and symmetrically
$ \, {\big( {(\ )}' \big)}^\vee = {(\ )}^\vee \circ {(\ )}' \, $
plays the role of a ``taking-the-interior operator'': in other
words, QFAs may be thought of as ``closed sets''  and QrUEAs
as ``open sets'' in  $ \HA \, $.
                                                 \par
   {\it (b)} \, {\sl Duality between Drinfeld's functors}.
For any  $ \, n \in \N \, $  let  $ \; \mu_n \, \colon
\, {J_{\scriptscriptstyle H}}^{\!\otimes n}
\lhook\joinrel\longrightarrow H^{\otimes n} \,{\buildrel
{\; m^n} \over \loongrightarrow}\, H \, $  be the composition of the
natural embedding of  $ {J_{\scriptscriptstyle H}}^{\!\otimes n} $
into  $ H^{\otimes n} $  with the  $ n $--fold  multiplication (in
$ H \, $):  then  $ \mu_n $  is the ``Hopf dual'' to  $ \delta_n \, $.
By construction we have  $ \, H^\vee = \sum_{n \in \N} \mu_n\big(
\h^{-n} {J_{\scriptscriptstyle H}}^{\!\otimes n}\big) \, $  and
$ \, H' = \bigcap_{n \in \N} {\delta_n}^{\!-1}\big(\h^{+n}
{J_{\scriptscriptstyle H}}^{\!\otimes n} \big) \, $:  \, this
shows that the two functors are built up as ``dual'' to each
other (cf.~also part  {\it (d)\/}  of Theorem 2.2).
                                                 \par
   {\it (c)} \, {\sl Ambivalence \;
       \hbox{QrUEA  $ \leftrightarrow $
QFA  \; in  $ \HA_F \, $.} \; Part  {\it (e)} of  Theorem 2.2
means that some}
   Hopf algebras  {\sl over\/}  $ F(R) $  might be thought
of  {\sl both\/}  as ``quantum function algebras''  {\sl and\/}
as ``quantum enveloping algebras'': examples are  $ U_F $  and
$ F_F $  for  $ \, U \in \QrUEA \, $  and  $ \, F \in \QFA \, $.
                                                 \par
   {\it (d)} \, {\sl Drinfeld's functors for algebras, coalgebras and
bialgebras}.  The definition of either of Drinfeld functors requires
only ``half of'' the notion of Hopf algebra.  In fact, one can define
$ (\ )^\vee $  for all ``augmented algebras'' (that is, roughly
speaking, ``algebras with a counit'') and  $ (\ )' $  for all
``coaugmented coalgebras'' (roughly, ``coalgebras with a unit''),
and in particular for bialgebras: this yields again nice functors,
and neat results extending the global quantum duality principle
hold for them; we shall prove all this in the next section.

\vskip1,9truecm
 \eject   

\centerline {\bf \S \; 3 \  General properties of Drinfeld's functors }

\vskip10pt

{\bf 3.1 Augmented algebras, coaugmented coalgebras and Drinfeld's
functors for them.}  Let  $ R $  be a commutative ring with 1,
$ \calM $  the category of torsion-free  $ R $--modules.
                                            \par
  We call  {\sl augmented algebra\/}  the datum of a unital associative
algebra  $ \, A \in \calM \, $  with a distinguished unital algebra
morphism  $ \, \undepsilon : \, A \longrightarrow R \, $  (so the unit
map  $ \, u : \, R \longrightarrow A \, $  is a section of  $ \undepsilon
\, $):  \, these form a category in the obvious way.  We call  {\sl
indecomposable elements\/}  of an augmented algebra  $ A $  the
elements of the set  $ \, Q(A) := J_{\scriptscriptstyle A} \big/
{J_{\scriptscriptstyle A}}^{\hskip-3pt 2} \, $  with
          \hbox{$ \, J_{\scriptscriptstyle A} := \text{\sl Ker}\,
\big(\, \undepsilon \, \colon A \longrightarrow R \,) \, $.}
We denote  $ \calA^+ $
the category of all augmented algebras in  $ \calM \, $.
                                            \par
   We call  {\sl coaugmented coalgebra\/}  any counital coassociative
coalgebra  $ C $    with a distinguished counital coalgebra morphism
$ \, \underline{u} : \, R \longrightarrow C \, $  (so  $ \underline{u} $ 
is a section of the counit map  $ \, \epsilon : \, C \longrightarrow R
\, $),  and let  $ \, \und1 := \underline{u}(1) \, $,  \, a group-like
element in  $ C \, $:  \, these form a category in the obvious way. 
For such a  $ C $  we said  {\sl primitive\/}  the elements of the
set  $ \, P(C) := \{\, c \in C \,|\; \Delta(c) \!=\! c \otimes \und1
+ \und1 \otimes c \,\} \, $.  We denote  $ \calC^+ $
 \hbox{the category of all coaugmented coalgebras in  $ \calM \, $.}
                                            \par
   We denote  $ \calB $  the category of all bialgebras in
$ \calM \, $;  clearly each bialgebra  $ B $  can be seen both as
an augmented algebra, w.r.t.~$ \, \undepsilon = \epsilon \equiv
\epsilon_{\! \scriptscriptstyle B} \, $  (the counit of  $ B \, $)
and as a coaugmented coalgebra, w.r.t.~$ \, \underline{u} = u \equiv
u_{\! \scriptscriptstyle B} \, $  (the unit map of  $ B \, $),  so
that  $ \, \und1 = 1 = 1_{\! \scriptscriptstyle B} \, $:  then
$ Q(B) $  is naturally a Lie coalgebra and  $ P(B) $  a Lie algebra
over  $ R \, $.  In the following we'll do such an interpretation
throughout, looking at objects of  $ \calB $  as objects
of  $ \calA^+ $  and of  $ \calC^+ $.
%
%
                                              \par
   Now let  $ R $  be a domain, and fix  $ \, \h \in R \setminus
\{0\} \, $  as in \S 1.3.  Let  $ \, A \in \calA^+ \, $,  \, and
        $ \, I = I_{\scriptscriptstyle \! A} := $\break
  $ := \hbox{\sl Ker} \, \Big( A
\,{\buildrel \undepsilon \over {\relbar\joinrel\twoheadrightarrow}}\,
R \,{\buildrel {\h \mapsto 0} \over \llongtwoheadrightarrow}\, R \big/
\h \, R = \Bbbk \Big) = \hbox{\sl Ker} \, \Big( A \,{\buildrel {\h
\mapsto 0} \over \llongtwoheadrightarrow}\, A \big/ \h \, A \,{\buildrel
\undepsilon{\,|}_{\h=0} \over {\relbar\joinrel\llongtwoheadrightarrow}}\,
\Bbbk \Big) \, $  as in \S 1.3, a maximal Hopf ideal of  $ A $
(where  $ \undepsilon{\,\big|}_{\h=0} $  is the counit of  $ \,
A{\big|}_{\h=0} \, $,  \, and the two composed maps do
coincide): like in \S 2.1, we define
  $$  A^\vee \; := \; {\textstyle \sum_{n \geq 0}} \, \h^{-n}
I^n \; = \; {\textstyle \sum_{n \geq 0}} \, {\big( \h^{-1} I \,
\big)}^n \; = \; {\textstyle \bigcup_{n \geq 0}} \, {\big( \h^{-1}
I \, \big)}^n  \quad  \big( \! \subseteq A_F \, \big) \; .  $$
If  $ \, J = J_{\scriptscriptstyle \! A} := \hbox{\sl Ker}\,
(\epsilon_{\scriptscriptstyle \! A}) \, $  then  $ \, I = J +
\h \cdot 1_{\scriptscriptstyle A} \, $,  \, thus  $ \; A^\vee =
\sum_{n \geq 0} \h^{-n} J^n = \sum_{n \geq 0} {\big( \h^{-1} J \,
\big)}^n \, $.
                                              \par
  Given any coalgebra  $ C $,  for every  $ \, n \!\in\! \N \, $
define  $ \; \Delta^n \colon C \longrightarrow C^{\otimes n}
\; $  by  $ \, \Delta^0 := \epsilon \, $,  $ \, \Delta^1 :=
\id_{\scriptscriptstyle C} $,  \, and  $ \, \Delta^n := \big(
\Delta \otimes \id_{\scriptscriptstyle C}^{\otimes (n-2)} \big) \circ
\Delta^{n-1} \, $  if  $ \, n > 2 \, $.  If  $ C $  is  {\sl coaugmented},
for any ordered subset  $ \, \Sigma = \{i_1, \dots, i_k\} \subseteq
\{1, \dots, n\} \, $  with  $ \, i_1 < \dots < i_k \, $,  \, define
the morphism  $ \; j_{\scriptscriptstyle \Sigma} : C^{\otimes k}
\longrightarrow C^{\otimes n} \; $  by  $ \; j_{\scriptscriptstyle
\Sigma} (a_1 \otimes \cdots \otimes a_k) := b_1 \otimes \cdots
\otimes b_n \; $  with  $ \, b_i := \und1 \, $  if  $ \, i \notin
\Sigma \, $  and  $ \, b_{i_m} := a_m \, $  for  $ \, 1 \leq m
\leq k $~;  then set  $ \; \Delta_\Sigma := j_{\scriptscriptstyle
\Sigma} \circ \Delta^k \, $,  $ \, \Delta_\emptyset := \Delta^0
\, $,  and  $ \; \delta_\Sigma := \sum_{\Sigma' \subset \Sigma}
{(-1)}^{n- \left| \Sigma' \right|} \Delta_{\Sigma'} \, $,  $ \;
\delta_\emptyset := \epsilon \, $.  Like in \S 2.1, the inverse
formula  $ \; \Delta_\Sigma = \sum_{\Psi \subseteq \Sigma} \delta_\Psi
\, $  holds.  We'll also use notation  $ \, \delta_0 := \delta_\emptyset
\, $,  $ \, \delta_n := \delta_{\{1, 2, \dots, n\}} \, $,  and the
useful formula  $ \; \delta_n = {(\id_{\scriptscriptstyle C} -
\epsilon \cdot \und1 \,)}^{\otimes n} \circ \Delta^n \, $,
\, for all  $ \, n \in \N_+ \, $.
                                              \par
   Now consider any  $ \, C \in \calC^+ \, $  and  $ \, \h \in R \, $
as in \S 1.3.  We define
  $$  C' \; := \; \big\{\, c \in C \,\big\vert\, \delta_n(c) \in
\h^n C^{\otimes n} ,  \; \forall \,\, n \in \N \, \big\}  \quad
\big( \! \subseteq C \, \big) \, .  $$
                                              \par
   In particular, according to our general remark above for any
$ \, B \in \calB \, $  (and any prime element  $ \, \h \in R \, $
as above)  $ \, B^\vee \, $  is defined w.r.t.~$ \, \undepsilon
= \epsilon_{\scriptscriptstyle \! B} \, $  and  $ \, B' \, $  is
defined w.r.t.~$ \, \und1 = 1_{\scriptscriptstyle \! B} \, $.

%
%
%
%
%

\vskip7pt

\proclaim{Lemma 3.2} \, Let  $ \, H \in \HA \, $,  \, and set  $ \,
\overline{H} := H \big/ H_\infty \, $  (notation of \S 1.3).  Then:
                                         \hfill\break
   \indent   {\it (a)}  $ \, H_\infty = {(H')}_\infty \, $,  $ \,
H_\infty \subseteq {\big(H^\vee\big)}_\infty \, $,  $ \; H_\infty $ 
is a Hopf ideal of  $ H \, $,  \, and  $ \, {\big( \overline{H}
\,\big)}_\infty \! = \{0\} \, $.  Moreover, there are natural
isomorphisms  $ \; {\big( \overline{H} \,\big)}^{\!\vee} =
H^\vee \! \Big/ H_\infty \, $,  $ \; {\big( \overline{H}
\,\big)}' = H' \! \Big/ H_\infty \, $.  
                                         \hfill\break
   \indent   {\it (b)}  $ \; \overline{H} \in \HA \, $,  \, and
$ \, {\overline{H}}{\big|}_{\h=0} = H{\big|}_{\h=0} \, $.  In
particular, if  $ H{\big|}_{\h=0} $  has no zero-divisors the
same holds for  $ H $,  and if  $ H $  is a QFA, resp.~a QrUEA,
then  $ \overline{H} $  is a QFA, resp.~a QrUEA.
                                         \hfill\break
   \indent   {\it (c)}  Analogous statements hold for any  $ \; A \in
\calA^+ $,  any  $ \, C \in \calC^+ $  and any  $ \; B \in \calB \, $.
\endproclaim

\demo{Proof} Trivial from definitions.   \qed
\enddemo

\vskip7pt

\proclaim{Proposition 3.3}  Let  $ \, A \!\in\! \calA^+ \! $,
$ B \!\in\! \calB $,  $ H \!\in\! \HA \, $.  Then  $ \, A^\vee
\!\!\in\! \calA^+ \! $,  $ B^\vee \!\!\in\! \calB \, $,  and
$ \, H^\vee \!\!\in\! \HA \, $.  If in addition  $ A $, 
resp.~$ B $,  is flat, then  $ A^\vee $,  resp.~$ B^\vee $, 
is flat as well.   
\endproclaim

\demo{Proof}  First, we have  $ \, A^\vee $,  $ B^\vee $,
$ H^\vee \in \calM \, $,  \, for they are clearly torsion-free. 
In addition,  $ A^\vee $  is obtained from  $ A $  in two steps:
localisation   --- namely,  $ \, A \succ\joinrel\rightsquigarrow
A\big[\h^{-1}\big] \, $  ---   and restriction, i.e.~taking
a submodule   --- namely,  $ \, A\big[\h^{-1}\big]
\succ\joinrel\rightsquigarrow A^\vee \big( \subseteq
A\big[\h^{-1}\big] \,\big) \, $.  Both these steps preserve
flatness, hence if  $ A $  is flat then  $ A^\vee $  is flat
too, and the same for  $ B $  and  $ B^\vee $. 
                                            \par
   Second,  $ \, A^\vee := \sum_{n=0}^\infty \h^n J^n \, $  where
$ \, J := \hbox{\sl Ker}\, (\epsilon_{\scriptscriptstyle A}) \, $,
so  $ A^\vee $  is clearly an  $ R $--subalgebra  of  $ \, A_F \, $,
\, hence  $ \, A^\vee \in \calA^+ $;  \, similarly holds for  $ B $
and  $ H $  of course.  Moreover,  $ J_{\scriptscriptstyle \! B} $
is bi-ideal of  $ B $,  so  $ \, \Delta(J_{\scriptscriptstyle \! B})
\subseteq B \otimes J_{\scriptscriptstyle B} + J_{\scriptscriptstyle
\! B} \otimes B \, $,  hence  $ \, \Delta \big( {J_{\scriptscriptstyle
\! B}}^{\!n} \big) \subseteq \! \sum\limits_{r+s=n} \!
{J_{\scriptscriptstyle \! B}}^{\!r} \otimes {J_{\scriptscriptstyle
\! B}}^{\!s} \, $  for all  $ \, n \in \N \, $,  \, thus  $ \, \Delta
\big( \h^{-n} {J_{\scriptscriptstyle \! B}}^{\!n} \big) \subseteq \h^{-n}
\! \sum\limits_{r+s=n} {J_{\scriptscriptstyle \! B}}^{\!r} \otimes
{J_{\scriptscriptstyle \! B}}^{\!s} \, = \sum\limits_{r+s=n} \! \big(
\h^{-r} {J_{\scriptscriptstyle \! B}}^{\!r} \big) \otimes \big( \h^{-s}
{J_{\scriptscriptstyle \! B}}^{\!s} \big) \, \subseteq \, B^\vee \otimes
B^\vee \, $  for all  $ \, n \in \N \, $,  \, whence  $ \, \Delta \big(
B^\vee \big) \subseteq B^\vee \otimes B^\vee \, $  which means  $ \,
B^\vee \in \calB \, $.  Finally, for  $ H $  we have in addition
$ \, S \big( {J_{\scriptscriptstyle \! H}}^{\!n} \big) =
{J_{\scriptscriptstyle \! H}}^{\!n} \, $  (for all  $ \, n
\in \N \, $)  because  $ J_{\scriptscriptstyle \! H} $  is a
Hopf ideal, therefore  $ \, S \big( \h^{-n} {J_{\scriptscriptstyle
\! H}}^{\!n} \big) = \h^{-n} {J_{\scriptscriptstyle \! H}}^{\!n} \, $
(for all  $ \, n \in \N \, $),  thus  $ \, S \big( H^\vee \big) =
H^\vee \, $  and so  $ \, H^\vee \in \HA \, $.   \qed
\enddemo

\vskip7pt

\proclaim{Lemma 3.4} Let  $ B $  be any bialgebra.  Let  $ \, a $,
$ b \in B $,  and let  $ \, \Phi \! \subseteq \! \N $,  with
$ \Phi $  finite.  Then\break
   \indent  (a) \quad \;\;  $ \delta_\Phi(ab) = \sum\limits_{\Lambda
\cup Y = \Phi} \delta_\Lambda(a) \, \delta_Y(b)  \, $;
                                     \hfill\break
   \indent  (b) \quad  if  $ \, \Phi \not= \emptyset \, $,  \;
then  \quad  $ \delta_\Phi(ab - ba) =
\sum\limits_{\Sb  \Lambda \cup Y = \Phi  \\
\Lambda \cap Y \not= \emptyset  \endSb}
\big( \delta_\Lambda(a) \, \delta_Y(b) - \delta_Y(b) \,
\delta_\Lambda(a) \big) \, $;
                                     \hfill\break
   \indent  (c) \quad  if the ground ring of  $ \, B $  is a field,
and if  $ \, D_n := \hbox{\sl Ker}\,(\delta_{n+1}) \, $  (for all
$ \, n \in \N \, $),  then
  $$  \underline{D} \;\;\; \colon \quad \{0\} =: D_{-1} \subseteq
D_0 \subseteq D_1 \subseteq \cdots D_n \subseteq \cdots  \quad
(\, \subseteq B \, )  $$
is a bialgebra filtration of  $ \, B \, $  with  $ \; \big[D_m,D_n\big]
\subseteq D_{m+n-1} \; (\,\forall \, m, n \in \N \,) \, $,  \, hence
the associated graded bialgebra is commutative.  If
%
%
$ \, B = H \, $  is a Hopf algebra, then  $ \underline{D} $  is a
Hopf algebra filtration, so the associated graded bialgebra is a
commutative graded Hopf algebra.
\endproclaim

\demo{Proof}  {\it (a)} \, (cf.~[KT], Lemma 3.2) First, notice that
the inversion formula  $ \; \Delta_\Phi = \sum_{\Psi \subseteq \Phi}
\delta_\Psi \; $  (see \S 2.1) gives  $ \; \sum_{\Psi \subseteq \Phi}
\delta_\Psi(ab) = \Delta_\Phi(ab) = \Delta_\Phi(a) \, \Delta_\Phi(b)
= \sum_{\Lambda, Y \subseteq \Phi} \delta_\Lambda(a) \, \delta_Y(b)
\, $;  \; this can be rewritten as
  $$  {\textstyle \sum\nolimits_{\Psi \subseteq \Phi}} \delta_\Psi(ab)
\; = \; {\textstyle \sum\nolimits_{\Psi \subseteq \Phi}} \,\,
{\textstyle \sum\nolimits_{\Lambda \cup Y = \Psi}}
\delta_\Lambda(a) \, \delta_Y(b) \; .   \eqno (3.1)  $$
   \indent   We prove the claim by induction on the cardinality
$ \vert \Phi \vert $  of  $ \Phi $.  If  $ \, \Phi = \emptyset \, $
then  $ \, \delta_\Phi = j_{\scriptscriptstyle \emptyset} \circ
\epsilon \, $,  which is a morphism of algebras, so the claim does
hold.  Now assume it holds for all sets of cardinality less than
$ \vert \Phi \vert $,  hence also for all proper subsets of
$ \Phi \, $:  then the right-hand-side of (3.1) equals  $ \;
\sum\limits_{\Psi \subsetneqq \Phi} \delta_\Psi(ab) \, +
\sum\limits_{\Lambda \cup Y = \Phi} \hskip-5pt \delta_\Lambda(a)
\, \delta_Y(b) \, $.  Then the claim follows by subtracting from
both sides of (3.1) the summands corresponding to the proper
subsets  $ \Psi $  of  $ \Phi $.
                                                    \par
   {\it (b)} \, (cf.~[KT], Lemma 3.2) The very definitions give  $ \,
\delta_\Lambda(a) \, \delta_Y(b) = \delta_Y(b) \, \delta_\Lambda(a)
\, $  when  $ \, \Lambda \cap Y = \emptyset \, $,  \, so the claim
follows from this and from  {\it (a)}.
                                                    \par
   {\it (c)} \, Let  $ \, a \in D_m \, $,  $ \, b \in D_n \, $:  \,
then  $ \, a \, b \in D_{m+n} \, $  because part  {\it(a)\/}  gives
$ \; \delta_{m+n+1}(a{}b) \, = \, \sum_{\Lambda \cup Y = \{1, \dots,
m+n+1\}} \delta_\Lambda(a) \, \delta_Y(b) \, = \, 0 \; $  since in the
sum one has  $ \, |\Lambda| > m \, $  or  $ \, |Y| > n \, $  which
forces  $ \, \delta_\Lambda(a) = 0 \, $  or  $ \, \delta_Y(b) = 0
\, $.  Similarly,  $ \, [a,b] \in D_{m+n-1} \leq m+n-1 \, $  because
part  {\it (b)\/}  yields  $ \; \delta_{m+n}\big([a,b]\big) \, = \,
\sum_{\Sb  \Lambda \cup Y = \{1,\dots,m+n\}  \\
           \Lambda \cap Y \not= \emptyset  \endSb}
\delta_\Lambda(a) \, \delta_Y(b) \, = \, 0 \; $.
%
%
%
%
                                                    \par
   Second, we prove that  $ \, \Delta\big(D_n) \subseteq \sum_{r+s=n}
D_r \otimes D_s \, $,  for all  $ \, n \in \N \, $.  Let  $ \, \eta
\in D_n \setminus D_{n-1} \, $.  Then  $ \, \Delta(\eta) =
\epsilon(\eta) \cdot 1 \otimes 1 + \eta \otimes 1 + 1 \otimes \eta
+ \delta_2(\eta) \, $;  \, since  $ \, D_0 := \hbox{\sl Ker}\,(\delta_1)
= \langle 1 \rangle = \Bbbk \cdot 1 \, $  we need only to show that
$ \, \delta_2(\eta) \in \sum_{r+s=n} D_r \otimes D_s \, $.  We can
write  $ \, \delta_2(\eta) = \sum_j u_j \otimes v_j \, $  with
$ \, u_j $,  $ v_j \in J := \hbox{\sl Ker}\,(\epsilon) \, $
---  so that  $ \, \delta_1(u_j) = u_j \, $  for all  $ j $  ---
and the  $ u_j $'s  linearly independent among themselves.  By
coassociativity of  $ \Delta $  one has  $ \, (\delta_r \otimes
\delta_s) \circ \delta_2 = \delta_{r+s} \, $  (for all  $ \, r $,
$ s \in \N \, $);  \, therefore,  $ \, 0 = \delta_{n+1}(\eta) =
\sum_j \delta_1(u_j) \otimes \delta_n(v_j) = \sum_j u_j \otimes
\delta_n(v_j) \, $:  \, since the  $ u_j $'s  are linearly independent,
this yields  $ \, \delta_n(v_j) = 0 \, $,  \, that is  $ \, v_j \in
\hbox{\sl Ker}\,(\delta_n) =: D_{n-1} \, $,  \, for all  $ j \, $.
                                                     \par
   Now, set  $ \, D^J_n := D_n \cap J \, $  for  $ \, n \in \N \, $.
Splitting  $ \, J \, $  as  $ \, J = D^J_1 \oplus W_1 \, $   --- for
some subspace  $ W_1 $  of  $ J $  ---   we can rewrite  $ \delta_2
(\eta) $  as  $ \, \delta_2(\eta) = \sum_i u^{(1)}_i \otimes v^{(n-1)}_i
+ \sum_h u^+_h \otimes v^+_h \, $,  \, where  $ \, u^{(1)}_i \in D_1^J
\, $,  $ \, u^+_h \in W_1 \, $,  $ \, v^{(n-1)}_i $,  $ v^+_h \in
D^J_{n-1} \, $  (for all  $ i $,  $ h \, $)  and the  $ u^+_h $'s  are
linearly independent.  Then also the  $ \delta_2\big(u^+_h\big) $'s
are linearly independent: indeed, if  $ \, \sum_h c_h \, \delta_2
\big(u^+_h\big) = 0 \, $  for some scalars  $ \, c_h \, $  then  $ \,
\sum_h c_h \, u^+_h \in \hbox{\sl Ker}\,(\delta_2) =: D_1 \, $,  \,
forcing  $ \, c_h = 0 \, $  for all  $ h \, $.  Then again by
coassociativity  $ \, 0 = \delta_{n+1}(\eta) = (\delta_2 \otimes
\delta_{n-1}) \big( \delta_2(\eta) \big) = \sum_i \delta_2 \big(
u^+_h \big) \otimes \delta_{n-1}\big(v^+_h\big) \, $,  \, which
--- as the  $ \delta_2\big(u^+_h\big) $'s  are linearly independent
---   yields  $ \, \delta_{n-1}\big(v^+_h\big) = 0 \, $,  \, i.e.~$ \,
v^+_h \in D_{n-2} \, $,  \, for all  $ h \, $.
                                            \par
   Now we repeat the argument.  Splitting  $ \, J \, $  as  $ \, J =
D_2 \oplus W_2 \, $   --- for some subspace  $ W_2 $  of  $ J $  ---
we can rewrite  $ \delta_2(\eta) $  as  $ \, \delta_2(\eta) = \sum_i
u^{(1)}_i \otimes v^{(n-1)}_i + \sum_j u^{(2)}_j \otimes v^{(n-2)}_j
+ \sum_k u^*_k \otimes v^*_k \, $,  \, where  $ \, u^{(2)}_j \in D_2
\bigcap J \, $,  $ \, v^{(n-2)}_j $,  $ v^*_k \in D_{n-2} \bigcap J
\, $,  $ \, u^*_k \in W_2 \, $  (for all  $ j $,  $ k \, $)  and
the  $ u^*_k $'s  are linearly independent.  Then also the
$ \delta_3\big(u^*_k\big) $'s  are linearly independent (as
above), and by coassociativity we get  $ \, 0 = \delta_{n+1}(\eta)
= (\delta_3 \otimes \delta_{n-2}) \big( \delta_2(\eta) \big) =
\sum_k \delta_3\big(u^*_k\big) \otimes \delta_{n-2}\big(v^*_k\big)
\, $,  \, which gives  $ \, \delta_{n-2}\big(v^*_k\big) = 0 \, $,
\, i.e.~$ \, v^*_k \in D_{n-3} \, $,  \, for all  $ k \, $.
Iterating this argument, we eventually stop getting  $ \;
\delta_2(\eta) = \sum_i u^{(1)}_i \otimes v^{(n-1)}_i +
\sum_j u^{(2)}_j \otimes v^{(n-2)}_j + \cdots + \sum_\ell
u^{(n-1)}_\ell \otimes v^{(1)}_\ell = \sum_{s=1}^{n-1} \sum_t
u^{(s)}_{s,t} \otimes v^{(n-s)}_{s,t} \; $  with  $ \, u^{(s)}_{s,t}
\in D_s \, $,  $ \, v^{(n-s)}_{s,t} \in D_{n-s} \, $  for all
$ \, s $,  $ t \, $,  \, so  $ \, \delta_2(\eta) \in \sum_{a+b=n}
D_a \otimes D_b \, $,  \, q.e.d.
                                               \par
   Finally, if  $ \, B = H \, $  is a Hopf algebra then  $ \, \Delta
\circ S = S^{\otimes 2} \circ \Delta \, $,  hence  $ \, \Delta^n \circ
S = S^{\otimes n} \circ \Delta^n \, $  ($ n \in \N \, $),  and  $ \,
\epsilon \circ S = S \circ \epsilon \, $,  \, thus  $ \, \delta^n
\circ S = S^{\otimes n} \circ \delta^n \, $  (for all  $ \, n \in
\N \, $)  follows, which yields  $ \, S\big(D_n\big) \subseteq
D_n \, $  for all  $ \, n \in \N \, $.  Thus  $ \underline{D} $
is a Hopf algebra filtration, and the rest follows.   \qed
\enddemo

\vskip7pt

\proclaim{Proposition 3.5}  Assume that  $ \, \Bbbk := R \big/ \h R \, $ 
is a field.  Let  $ \, C \in \calC^+ \! $,  $ B \in \calB $,  $ H \in
\HA \, $.  Then  $ \, C' \! \in \calC^+ \! $,  $ B' \! \in \calB \, $, 
and  $ \, H' \! \in \HA \, $.  Moreover, if  $ C $,  $ B $,  is flat,
then  $ C' $,  $ B' $,  is flat too.   
\endproclaim

\demo{Proof}  First, by definition  $ C' $  is an  $ R $--submodule
of  $ C $,  because the maps  $ \delta_n $  ($ n \in \N \, $)  are
$ R $--linear;  since  $ C $  is torsion-free, its submodule
$ C' $  is torsion-free too, i.e.~$ \, C' \! \in \calM \, $.  In
addition,  $ C' $  is an  $ R $--submodule  of  $ C $,  and taking
a submodule preserve flatness: hence if  $ C $  is flat then  $ C' $ 
is flat too.  The same holds for  $ B $  and  $ B' $  as well.   
                                        \par
   We must show that  $ C' $  is a subcoalgebra.  Due to
Lemma 3.2{\it (c)},  we can reduce to prove it for  $ \,
\big(\overline{C}\,\big)' $,  \, that is we can assume from
scratch that  $ \, C_\infty = \{0\} \, $.
                                        \par
   Let  $ \, \Rhat \, $  be the  $ \h $--adic  completion of  $ R \, $.
Let also  $ \, \widehat{C} \, $  be the  $ \h $--adic  completion of
$ C \, $:  this is a separated complete topological  $ \Rhat $--module,
hence it is topologically free (i.e.~of type  $ \Rhat^{\,Y} $  for
some set  $ Y \, $);  moreover, it is a topological Hopf algebra,
whose coproduct takes values into the  $ \h $--adic  completion
$ \, C \, \widehat{\otimes} \, C \, $  of  $ \, C \otimes C \, $.
Since  $ \, C_\infty = \{0\} \, $,  \, the natural map  $ \, C
\longrightarrow \widehat{C} \, $  is a  {\sl monomorphism\/}  of
(topological) Hopf  $ R $--algebras,  so  $ C $  identifies with
a Hopf  $ R $--subalgebra  of  $ \widehat{C} $.  Further, we have
$ \, \widehat{C} \Big/ \h^n \widehat{C} = C \Big/ \h^n C \, $  for
all  $ \, n \in \N \, $.  Finally, we set  $ \, {\widehat{C}}^*
:= {\text{\sl Hom}}_{\Rhat} \Big( \widehat{C} \, , \, \Rhat \,
\Big) \, $  for the dual of  $ \widehat{C} $.
                                        \par
   Pick  $ \, a \in C' \, $;  \, first we prove that  $ \, \Delta(a)
\in C' \otimes C \, $:  \, to this end, since  $ \widehat{C} $  is
topologically free it is enough to show that  $ \; (\hbox{id} \otimes f)
\big( \Delta(a) \big) \in C' \otimes_R \Rhat \cdot \und1 \; $  for all
$ \, f \in \widehat{C}^* $,  \, which amounts to show that  $ \; \big(
(\delta_n \otimes f) \circ \Delta \big) (a) \in \h^n C^{\otimes n}
\otimes_R \Rhat \cdot \und1 \; $  for all  $ \, n \in \N_+ \, $,
$ \, f \in \widehat{C}^* \, $.  Now, we can rewrite the latter term as
  $$  \big( (\delta_n \otimes f) \circ \Delta \big)(a) =
\big( \big( {(\hbox{id} - \epsilon \cdot \und1 \,)}^{\otimes n}
\otimes f \, \big) \circ \Delta^{n+1} \big)(a) = \delta_n(a) \otimes
f(\und1 \,) \cdot \und1 + \big( \hbox{id}^{\otimes n} \otimes f \,
\big) \big( \delta_{n+1}(a) \big)  $$
and the right-hand-side term does lie in  $ \, \h^n C^{\otimes n}
\otimes_R \Rhat \cdot \und1 \, $,  \; for  $ \, a \in C' $,  \, q.e.d.
                                              \par
   Definitions imply  $ \, \Delta(x) = - \epsilon(x) \cdot \und1
\otimes \und1 + x \otimes \und1 + \und1 \otimes x + \delta_2(x) \, $
for all  $ \, x  \in C \, $.  Due to the previous analysis, we argue
that  $ \, \delta_2(a) \in C' \otimes C \, $  for all  $ \, a \in C'
\, $,  \, and we only need to show that  $ \, \delta_2(a) \in C' \otimes
C' \, $:  \, this will imply  $ \, \Delta(a) \in C' \otimes C' \, $
since  $ \, \und1 \in C' \, $  (as it is group-like).
                                              \par
   Let  $ \, \widehat{C'} \, $  be the  $ \h $--adic  completion of
$ C' $:  again, this is a topologically free  $ \Rhat $--module,
and since  $ \, {(C')}_\infty = C_\infty = \{0\} \, $  (by  Lemma
3.2{\it (a)}  and our assumptions)  the natural map  $ \, C'
\longrightarrow \widehat{C'} \, $  is in fact an embedding, so
$ C' $  identifies with an  $ R $--submodule  of  $ \widehat{C'} \, $.
If  $ \, \big\{\, \beta_j \,\big\vert\, j \in \Cal{J} \,\big\} \, $
is a subset of  $ C' $  whose image in  $ \, C'{\big|}_{\h=0} \, $
is a basis of the latter  $ \Bbbk $--vector  space, then it is easy
to see that  $ \, C' = \bigoplus\limits_{j \in \Cal{J}} \Rhat \,
\beta_j \; $:  \, fixing a section  $ \, \nu \, \colon \, \Bbbk
\lhook\joinrel\longrightarrow R \, $  of the projection map  $ \, R
\relbar\joinrel\twoheadrightarrow R \big/ \h \, R =: \Bbbk \, $,  \,
this implies that each element  $ \, a \in C' \, $  has a unique
expansion as a series  $ \, a = \sum_{n \in \N} \sum_{j \in \Cal{J}}
\nu(\kappa_{j,n}) \, \h^n \, \beta_j \, $  for some  $ \, \kappa_{j,n}
\in \Bbbk \, $  which, for fixed  $ n $,  are almost all zero.
Finally  $ \; \widehat{C'} \Big/ \h \, \widehat{C'} = C'
\Big/ \h \, C' = B \Big/ \h \, B \, $,  \, with  $ \, B
:= \bigoplus\limits_{j \in \Cal{J}} R \, \beta_j \, $.
                                              \par
   We shall also consider  $ \, \big( \widehat{C'} \, {\big)}^* :=
{\text{\sl Hom}}_{\Rhat} \Big( \widehat{C'} \, , \, \Rhat \,\Big)
\, $,  \, the dual of  $ \, \widehat{C'} \, $.
                                              \par
   Now, we have  $ \, \delta_2(a) \in C' \otimes C \subseteq
\widehat{C'} \otimes C \, $,  \, so we can expand  $ \, \delta_2(a)
\, $  inside  $ \, \widehat{C'} \otimes C \, $  as  $ \; \delta_2(a)
= \sum_{i \in \Cal{I}} \Big( \sum_{n \in \N} \sum_{j \in \Cal{J}}
\nu\big(\kappa^i_{j,n}\big) \, \h^n \, \beta_j \Big) \otimes c_i \; $
for some  $ \, \kappa^i_{j,n} \in \Bbbk \, $  as above and  $ \, c_i
\in C \, $  ($ \Cal{I} $  being some finite set).  Then we can
rewrite  $ \delta_2(a) $  as
  $$  \delta_2(a) \, = \, {\textstyle \sum\limits_{j \in \Cal{J}}}
\, \beta_j \, \otimes \bigg( \, {\textstyle \sum\limits_{n \in \N}} \,
{\textstyle \sum\limits_{i \in \Cal{I}}} \, \nu\big(\kappa^i_{j,n}\big)
\, \h^n \cdot c_i \bigg) \, = \, {\textstyle \sum\limits_{j \in \Cal{J}}}
\, \beta_j \otimes \gamma_j \; \in \, C' \,\widetilde{\otimes}\,
\widehat{C}  $$
where  $ \, \gamma_j := \sum_{n \in \N} \sum_{i \in \Cal{I}}
\nu\big(\kappa^i_{j,n}\big) \, \h^n \cdot c_i \in \widehat{C} \, $
for all  $ j $,  and  $ \, C' \,\widetilde{\otimes}\, \widehat{C}
\, $  is the completion of  $ \, C' \otimes \widehat{C} \, $
w.r.t.~the weak topology.  We contend that all the  $ \gamma_j $'s
belong to  $ \big( \widehat{C} \,{\big)}' $.
                                              \par
   In fact, assume this is false: then there is  $ \, s \in \N_+
\, $  such that  $ \, \delta_s(\gamma_i) \not\in \h^s \widehat{C}^{\,
\widehat{\otimes}\, s} \, $  for some  $ \, i \in \Cal{J} \, $;  \,
we can choose such an  $ i $  so that  $ s $  be minimal, thus  $ \,
\delta_{s'}(\gamma_j) \in \h^{s'} \widehat{C}^{\,\widehat{\otimes}\,
s'} \, $  for all  $ \, j \in \Cal{J} \, $  and  $ \, s'<s \, $;
\, since  $ \, \delta_s = (\delta_2 \otimes \text{id}) \circ
\delta_{s-1} \, $   --- by coassociativity ---   we have also  $ \,
\delta_s(\gamma_j) \in \h^{s-1} \widehat{C}^{\,\widehat{\otimes}\, s}
\, $  for all  $ \, j \in \Cal{J} $.  Now consider the element
  $$  A \, := \, {\textstyle \sum_{j \in \Cal{J}}} \, \overline{\beta_j}
\otimes \overline{\delta_s(\gamma_j)} \,\; \in \; \Big( \widehat{C'}
\Big/ \h \, \widehat{C'} \,\Big) \,\widetilde{\otimes}_\Bbbk
\Big( \h^{s-1} \widehat{C}^{\,\widehat{\otimes}\, s} \Big/
\h^s \widehat{C}^{\,\widehat{\otimes}\, s} \Big)  $$
the right-hand-side space being equal to  $ \, \big( C' \big/ \h \, C'
\big) \,\widetilde{\otimes}_\Bbbk \big( \h^{s-1} C^{\otimes s} \big/
\h^s C^{\otimes s} \big) \, $;  hereafter, such notation as
$ \overline{x} $  will always denote the coset of  $ x $
in the proper quotient space.  By construction, the
$ \overline{\beta_j} $'s  are linearly independent and
some of the  $ \overline{\delta_s(\gamma_j)} $'s  are non
zero: therefore  $ A $  is non zero, and we can write it as
$ \; A = \sum_{\ell \in \Cal{L}} \overline{\lambda_\ell} \otimes
\overline{\chi_\ell} \;\; (\not= 0) \; $  where  $ \Cal{L} $
is a suitable non-empty index set,  $ \lambda_\ell $  (for all
$ \ell $)  belongs to the completion  $ \, \widetilde{C'} \, $
of  $ \widehat{C'} $  w.r.t.~the weak topology,  $ \, \chi_\ell
\in \h^{s-1} \widehat{C}^{\,\widehat{\otimes}\, s} \, $,  the
$ \overline{\lambda_\ell} $'s  are linearly  independent in
the  $ \Bbbk $--vector  space  $ \, \widetilde{C'} \big/ \h \,
\widetilde{C'} \, $  (which is just the completion of  $ \,
C'{\big|}_{\h=0} := C' \big/ \h \, C' \, $  w.r.t.~the weak topology),
and the  $ \overline{\chi_\ell} \, $'s  are linearly independent in
the  $ \Bbbk $--vector  space  $ \, \h^{s-1} \widehat{C}^{\,
\widehat{\otimes}\, s} \Big/ \h^s \widehat{C}^{\,\widehat{\otimes}
\, s} \, $.  In particular  $ \, \lambda_\ell \not\in \h \,
\widetilde{C'} \, $  for all  $ \ell \, $:  \, so there is
$ \, r \in \N_+ \, $  such that  $ \, \delta_r(\lambda_\ell)
\in \h^r \widetilde{C}^{\,\widetilde{\otimes}\, r} \setminus \h^{r+1}
\widetilde{C}^{\,\widetilde{\otimes}\, r} \, $  for all  $ \, \ell
\in \Cal{L} \, $  (hereafter,  $ \, K^{\,\widetilde{\otimes}\, m}
\, $  denotes the completion of  $ \, K^{\otimes m} \, $  w.r.t.~the
weak topology),  \, hence  $ \, \overline{\delta_r(\lambda_\ell)}
\not= 0 \in \h^r \widetilde{C}^{\,\widetilde{\otimes}\, r} \Big/
\h^{r+1} \widetilde{C}^{\,\widetilde{\otimes}\, r} \, $.  Now
write  $ \, \overline{\delta_n} \, $  for the composition
of  $ \delta_n $  with a projection map (such as  $ \, X
\relbar\joinrel\twoheadrightarrow X \big/ \h \, X \, $,
\, say):  then the outcome of this analysis is that
\vskip3pt
   \centerline{ $ \displaystyle{
  \big( \overline{\delta_r} \otimes \overline{\delta_s}
\,\big) \big(\delta_2(a)\big) \, = \,  \big(\, \overline{\delta_r}
\otimes \text{id} \big) \Big( {\textstyle \sum\nolimits_{j \in J}}
\, \overline{\beta_j} \otimes \overline{\delta_s(\gamma_j)} \,
\Big) \, = \, {\textstyle \sum\nolimits_{\ell \in \Cal{L}}} \,
\overline{\delta_r (\lambda_\ell)} \otimes \overline{\chi_\ell}
\, \not= \, 0 } $ }
\vskip2pt
\noindent
in the  $ \Bbbk $--vector  space
$ \;
%
%
 \Big( \h^r C^{\otimes r} \Big/ \h^{r+1} C^{\otimes r} \Big)
\,\widetilde{\otimes}_\Bbbk \Big( \h^{s-1} C^{\otimes s} \Big/ \h^s
C^{\otimes s} \Big) \, $.
                                              \par
   On the other hand, coassociativity yields  $ \, (\delta_r \otimes
\delta_s) \big( \delta_2(a) \big) = \delta_{r+s}(a) \, $.  Therefore,
since  $ \, a \in C' \, $  we have  $ \; \delta_{r+s}(a) \in \h^{r+s}
C^{\otimes (r+s)} \, $,  \, hence  $ \; \overline{\delta_{r+s}(a)} =
0 \; $  in the  $ \Bbbk $--vector  space  $ \, \h^{r+s-1} C^{\otimes
(r+s)} \Big/ \h^{r+s} C^{\otimes (r+s)} \, $.  Now, there are standard
isomorphisms
\vskip-11pt
  $$  \displaylines{
   \h^\ell C^{\otimes \ell} \Big/ \h^{\ell+1} C^{\otimes r} \; \cong
\; \big( \h^\ell C^{\otimes \ell} \, \big) \otimes_R \Bbbk  \qquad
\hbox{for}  \quad  \ell \in \{r, s-1, r+s-1\}  \cr
   \big( \h^{r+s-1} C^{\otimes s} \, \big) \otimes_R \Bbbk \, \cong
\, \Big( \! \big( \h^r C^{\otimes r} \, \big) \otimes_R \Bbbk \Big)
\otimes_\Bbbk \Big( \! \big( \h^{s-1} C^{\otimes s} \, \big)
\otimes_R \Bbbk \Big)  \cr
   \quad  \h^{r+s-1} C^{\otimes (r+s)} \Big/ \h^{(r+s)} C^{\otimes
(r+s)} \; \cong \; \Big( \h^r C^{\otimes r} \Big/ \h^{r+1}
C^{\otimes r} \Big) \otimes_\Bbbk \Big( \h^{s-1} C^{\otimes s}
\Big/ \h^s C^{\otimes s} \Big) \, .  \quad  \cr }  $$
Moreover,  $ \, \Big( \h^r C^{\otimes r} \Big/ \h^{r+1}
C^{\otimes r} \Big) \otimes_\Bbbk \Big( \h^{s-1} C^{\otimes s}
\Big/ \h^s C^{\otimes s} \Big) \, $  naturally embeds, as a dense
subset, into  $ \, \Big( \h^r C^{\otimes r} \Big/ \h^{r+1} C^{\otimes r}
\Big) \,\widetilde{\otimes}_\Bbbk\, \Big( \h^{s-1} C^{\otimes s} \Big/
\h^s C^{\otimes s} \Big) \, $,  \, so via the last isomorphism above
we get
 \eject   
  $$  \h^{r+s-1} C^{\otimes (r+s)} \Big/ \h^{(r+s)} C^{\otimes (r+s)} \;
\lhook\joinrel\relbar\joinrel\longrightarrow\; \Big( \h^r C^{\otimes r}
\Big/ \h^{r+1} C^{\otimes r} \Big) \,\widetilde{\otimes}_\Bbbk \Big(
\h^{s-1} C^{\otimes s} \Big/ \h^s C^{\otimes s} \Big) \; .  $$
This last monomorphism maps  $ \, \overline{ (\delta_r \otimes
\delta_s) \big( \delta_2(a) \big)} = \overline{\delta_{r+s}(a)}
= 0 \, $  onto  $ \, \big(\, \overline{\delta_r} \otimes
\overline{\delta_s} \,\big) \big( \delta_2(a) \big) \not= 0 \, $,
\, a contradiction!  Therefore we must have  $ \, \gamma_j \in
\big( \widehat{C} \,{\big)}' \, $  for all  $ \, j \in \Cal{J}
\, $,  \, as contended.
                                               \par
   The outcome is  $ \, \delta_2(a) = \sum\nolimits_{j \in \Cal{J}}
\beta_j \otimes \gamma_j \, \in C' \widetilde{\otimes} \big(
\widehat{C} \,{\big)}' $,  \, so  $ \, \delta_2(a) \in \big(
C' \otimes C \big) \bigcap \big( C' \widetilde{\otimes} \big(
\widehat{C} \,{\big)}' \big) $.
                                               \par
   For all  $ n \in \! \N $,  the result above yields  $ \, \big(
\text{id} \otimes \delta_n \big) \! \big( \delta_2(a) \big)
\in \big( C' \otimes C^{\otimes n} \big) \bigcap \Big( C'
\widetilde{\otimes}\, \delta_n
   \big(\! \big( \widehat{C} \,\big)' \big) \!\Big) \subseteq $\break
$ \subseteq \big( C' \otimes C^{\otimes n} \big) \bigcap \Big( C'
\,\widetilde{\otimes}\, \h^n \widehat{C}^{\,\widehat{\otimes}\, n}
\Big) = C' \otimes \h^n C^{\otimes n} \, $,  \; because  $ \,
\widehat{C} \Big/ \h^n \widehat{C} = C \Big/ \h^n C \, $  (see
above) implies  $ \, C \bigcap \h^n \widehat{C} = \h^n C \, $.  So
we found  $ \, \big( \text{id} \otimes \delta_n \big) \! \big(
\delta_2(a) \big) \in C' \otimes \h^n C^{\otimes n} \, $  for all
$ \, n \in \N \, $.  Acting like in the first part of the proof,
we'll show that this implies  $ \, \delta_2(a) \in C' \otimes C'
\, $.  To this end, it is enough to show that  $ \; (f \otimes
\hbox{id}) \big( \delta_2(a) \big) \in \Rhat \cdot \und1 \otimes_R
C' \; $  for all  $ \, f \in \big( \widehat{C'} \,{\big)}^* $,  \,
which amounts to show that  $ \; (f \otimes \delta_n) \big( \delta_2(a)
\big) \in \Rhat \cdot \und1 \otimes_R \h^n C^{\otimes n} \; $  for all
$ \, n \in \N_+ \, $,  $ \, f \in \big( \widehat{C'} \,{\big)}^* \, $.
But this is true because  $ \; (f \otimes \delta_n) \big( \delta_2(a)
\big) = (f \otimes \hbox{id}) \big( (\hbox{id} \otimes \delta_n) \big(
\delta_2(a) \big) \big) \in (f \otimes \hbox{id}) \big( C' \otimes
\h^n C^{\otimes n} \big) \subseteq \Rhat \cdot \und1 \otimes_R \h^n
C^{\otimes n} \, $.  We conclude that  $ \, C' \in \calC^+ $,
\, q.e.d.
                                                 \par
   Now look at  $ \, B \in \calB \, $.  By the previous part we have
$ \, B' \in \calC^+ $.  Moreover,  $ B' $  is multiplicatively closed,
thanks to Lemma 3.4{\it (a)},  and  $ \, 1 \in B' \, $  by the very
definitions.  Thus  $ B' $  is an  $ R $--sub-bialgebra  of  $ B $,
so  $ \, B' \in \calB \, $.
                                                 \par
   Finally, for  $ \, H \in \HA \, $  one has in addition  $ \, \Delta
\circ S = S^{\otimes 2} \circ \Delta \, $,  \, which implies $ \,
\Delta^n \circ S = S^{\otimes n} \circ \Delta^n \, $  hence  $ \,
\delta_n \circ S = S^{\otimes n} \circ \delta_n \, $,  for all
$ \, n \in \N \, $.  This clearly yields  $ \, S(H') = H' \, $,
\, whence  $ H' $  is a Hopf subalgebra of  $ H $,  thus  $ \,
H' \in \HA \, $,  \, q.e.d.   \qed
\enddemo

\vskip3pt

   {\sl $ \underline{\hbox{\it Remark}} $:}  \; The ``hard step'' in
the previous proof   --- i.e.~proving that  $ \, \Delta(C') \subseteq
C' \otimes C' \, $  ---   is much simpler when, after the reduction
step to  $ \, C_\infty = \{0\} \, $,  one has that  $ C $  is free,
as an  $ R $--module  (note also that for  $ C $  free one has
automatically  $ \, C_\infty = \{0\} \, $).  In fact, in this
case   ---  i.e.~if  $ C $  is free ---   we don't need to use
completions.  The argument to prove that  $ \, \delta_2(a) \in C'
\otimes C \, $  goes through untouched, just using  $ C $  instead
of  $ \widehat{C} $,  the freeness of  $ C $  taking the role of the
topological freeness of  $ \widehat{C} \, $;  \, similarly, later
on if  $ C' $  also is free (for instance, when  $ R $  is a PID,
for  $ C' $  is an  $ R $--submodule  of the free  $ R $--module
$ C \, $)  we can directly use it instead of the topologically
free module  $ \widehat{C'} $,  just taking  $ \, \big\{\, \beta_j
\,\big\vert\, j \in \Cal{J} \,\big\} \, $  to be an  $ R $--basis
of  $ C' \, $:  \, then we can write  $ \, \delta_2(a) = \sum_{j \in
\Cal{J}} \beta_j \otimes \gamma_j \, \in \, C' \otimes C \, $  for
some  $ \, \gamma_j \in C \, $, \, and the argument we used applies
      \hbox{again to show that now  $ \, \gamma_j \in C' \, $  for
all  $ j \, $,  so that  $ \, \delta_2(a) \in C' \otimes C' \, $,
\, q.e.d.}

\vskip7pt

\proclaim{Theorem 3.6}  Assume that  $ \, \Bbbk := R \big/ \h R \, $ 
is a field.   
                                           \hfill\break
   \hbox{\hskip5,3pt   (a)  $ X \mapsto \! X^\vee \! $  gives
well-defined functors from  $ \calA^+ \! $  to  $ \calA^+ \! $,
from  $ \calB $  to  $ \calB $,  from  $ \HA $  to  $ \HA \, $.}
                                           \hfill\break
   \hbox{\hskip5,3pt   (b)  $ X \mapsto \! X' \! $  gives
well-defined functors from  $ \, \calC^+ \! $  to  $ \,
\calC^+ $,  from  $ \calB $ to  $ \calB $,  from  $ \HA $ 
to  $ \HA \, $.}
                                           \hfill\break
   {\phantom{.}} \hskip-1,5pt   (c) For any  $ \, B \in \calB \, $
we have  $ \, B \! \subseteq \! {\big(B^\vee\big)}' $,  $ \, B \!
\supseteq \! {\big(B'\big)}^{\!\vee} \! $,  hence  $ \, B^\vee \!
= \! \Big( \! {\big(B^\vee\big)}' \Big)^{\!\!\vee} \! $,  $ \, B'
\! = \! \Big( \! {\big(B'\big)}^{\!\vee} \Big)' \! $.
\endproclaim

\demo{Proof}  In force of Propositions 3.3--5, to define the functors
we only have to set them on morphisms. So let  $ \, \varphi \in
\hbox{\sl Mor}_{\scriptscriptstyle \calA^+}(A,E) \, $  be a morphism
in  $ \, \calA^+ \, $:   by scalar extension it gives a morphism  $ \,
A_F \longrightarrow E_F \, $  of Hopf  $ F(R) $--algebras,  which
maps  $ \, \h^{-1} J_{\scriptscriptstyle A} \, $  into  $ \, \h^{-1}
J_{\scriptscriptstyle E} \, $,  hence  $ A^\vee $  into  $ E^\vee
\, $:  \, this yields the morphism  $ \; \varphi^\vee \in \hbox{\sl
Mor}_{\scriptscriptstyle \calA^+} \big( A^\vee, \, E^\vee \big) \; $
we were looking for.  On the other hand, if  $ \, \psi \in \hbox{\sl
Mor}_{\scriptscriptstyle \calC^+}(C,\varGamma) \, $  is a morphism in
$ \, \calC^+ \, $  then $ \, \delta_n \circ \psi = \psi^{\otimes n}
\circ \delta_n \, $  for all  $ n \, $,  so  $ \, \psi \big( C'
\big) \subseteq \varGamma' \, $:  thus as  $ \, \psi' \in \hbox{\sl
Mor}_{\scriptscriptstyle \calC^+} \big( C', \varGamma' \big) \, $
we take the restriction  $ \, \psi{\big\vert}_{C'} \, $  of  $ \,
\psi \, $  to  $ C' $.
                                                 \par
   Now consider  $ \, B \in \calB \, $.  For any  $ \, n \in \N \, $
we have  $ \, \delta_n(B) \subseteq {J_{\scriptscriptstyle B}}^{\!
\otimes n} \, $  (see \S 2.1); this can be read as  $ \, \delta_n(B)
\subseteq {J_{\scriptscriptstyle B}}^{\! \otimes n} = \h^n {\big(
\h^{-1} J_{\scriptscriptstyle B} \big)}^{\otimes n} \subseteq \h^n
{\big( B^\vee \big)}^{\otimes n} \, $,  which gives  $ \, B \subseteq
{\big( B^\vee \big)}' \, $,  q.e.d.  On the other hand, let  $ \, I'
:= \hbox{\sl Ker} \Big( B' \,{\buildrel \epsilon \over
{\relbar\joinrel\twoheadrightarrow}}\, R \,{\buildrel {\h \mapsto 0}
\over {\llongtwoheadrightarrow}}\, \Bbbk \Big) \, $;  since  $ \,
{\big( B' \big)}^{\!\vee} := \bigcup_{n=0}^\infty {\big( \h^{-1} I'
\big)}^n \, $,  \, in order to show that  $ \, B \supseteq {\big(
B' \big)}^{\!\vee} \, $  it is enough to check that  $ \, B \supseteq
\h^{-1} I' \, $.  So let  $ \, x' \in I' \, $:  then  $ \, \delta_1(x')
\in \h \, B \, $,  \, hence  $ \, x' = \delta_1(x') + \epsilon(x') \in
\h \, B \, $.  Therefore  $ \, \h^{-1} x' \in B \, $,  q.e.d.  Finally,
the last two identities follow easily from the two inclusions we've
just proved.   \qed
\enddemo

\vskip7pt

\proclaim {Theorem 3.7} \, Assume that  $ \, \Bbbk := R \big/ \h R \, $ 
is a field.   
                                       \par   
   Let  $ \, B \in \calB \, $.  Then  $ \, B^\vee{\Big|}_{\h=0} $   
 is an rUEA (see \S 1.1),   
%
%
 generated by  $ \;\; \h^{-1} J \mod\, \h \, B^\vee \, $. 
                                       \par   
   In particular, if  $ \, H \in \HA \, $  then  $ \, H^\vee \in
\QrUEA \, $.
\endproclaim

\demo{Proof} A famous characterization theorem in Hopf algebra theory
claims the following (cf.~for instance [Ab], Theorem 2.5.3, or [Mo],
Theorem 5.6.5, and references therein, noting also that in the
cocommutative case  {\sl connectedness\/}  and  {\sl irreducibility\/}
coincide):
                                       \par   
   {\sl A Hopf algebra  $ H $  over a ground field  $ \Bbbk $  is
the restricted universal enveloping algebra of a restricted Lie
algebra  $ \gerg $  if and only if  $ H $  is generated by
$ P(H) $  (the set of primitive elements of  $ H $)  and it is
cocommutative and connected.  In that case,  $ \, \gerg = P(H) $.}
                                       \par
   Thus we must prove that the bialgebra  $ \, B^\vee{\big|}_{\h=0}
\, $  is in fact a Hopf algebra, it is generated by its primitive
part  $ \, P \big( B^\vee{\big|}_{\h=0} \big) \, $  and it is
cocommutative and connected, for then  $ \, B^\vee{\big|}_{\h=0}
= \U(\gerg) \, $  with  $ \, \gerg = P \big( B^\vee{\big|}_{\h=0}
\big) \, $  being a restricted Lie  {\sl bialgebra\/}  (by Remark 1.5).
                                       \par
  Since  $ \, B^\vee = \sum_{n \geq 0} {\big( \h^{-1} J \, \big)}^n
\, $,  \, it is generated, as a unital algebra, by  $ \, J^\vee :=
\h^{-1} J \, $.  Consider  $ \, j^\vee \in J^\vee \, $,  \, and
$ \, j := \h \, j^\vee \in J \, $;  \, then
  $$  \Delta(j) = \delta_2(j) + j \otimes 1 + 1 \otimes j -
\epsilon(j) \cdot 1 \otimes 1 \in j \otimes 1 + 1 \otimes j
+ J \otimes J  $$
for  $ \, \Delta = \delta_2 + \hbox{id} \otimes 1 + 1 \otimes \hbox{id}
- \epsilon \cdot 1 \otimes 1 \, $  and  $ \, \hbox{\sl Im}\,(\delta_2)
\subseteq J \otimes J \, $  by construction.  Therefore
  $$  \displaylines{
   \quad  \Delta \big( j^\vee \big) = \delta_2 \big( j^\vee \big) +
j^\vee \otimes 1 + 1 \otimes j^\vee - \epsilon\big( j^\vee \big)
\cdot 1 \otimes 1 = \delta_2 \big( j^\vee \big) +
j^\vee \otimes 1 + 1 \otimes j^\vee \, \in   \hfill  \cr
   \hfill   \in \, j^\vee \otimes 1 + 1 \otimes j^\vee +
\h^{-1} J \otimes J = j^\vee \otimes 1 + 1 \otimes j^\vee
+ \h^{+1} J^\vee \otimes J^\vee  \cr }  $$
whence
  $$  \hskip25pt   \Delta \big( j^\vee \big) \equiv j^\vee \otimes 1
+ 1 \otimes j^\vee \mod \h \, B^\vee   \hskip25pt  (\, \forall \;
j^\vee \in J^\vee) \, .   \eqno (3.2)  $$
This proves that  $ \, J^\vee{\big|}_{\h=0} \subseteq P \big(
B^\vee{\big|}_{\h=0} \big) \, $,  \, and since  $ J^\vee{\big|}_{\h=0} $
generates  $ B^\vee{\big|}_{\h=0} $  (for  $ J^\vee $  generates
$ B^\vee $),  {\it a fortiori\/}  $ B^\vee{\big|}_{\h=0} $  is
generated by  $ P \big( B^\vee{\big|}_{\h=0} \big) $,  \, hence
$ B^\vee{\big|}_{\h=0} $  is cocommutative too.  In addition, (3.2)
enables us to apply Lemma 5.5.1 in [Mo]   --- which is stated there
for Hopf algebras, but  {\sl holds indeed for bialgebras as well\/}
---   to the bialgebra  $ B^\vee{\big|}_{\h=0} \, $,  with  $ \,
A_0 = \Bbbk \cdot 1 \, $  and  $ \, A_1 = J^\vee \Big/ \big( J^\vee
\cap \, \h \, B^\vee \big) \, $:  \; then that lemma proves that  $ \,
B^\vee{\big|}_{\h=0} \, $  is connected.  Another classical result
(cf.~[Ab], Theorem 2.4.24) then ensures that  $ \, B^\vee{\big|}_{\h=0}
\, $  is indeed a Hopf algebra; as it is also connected, cocommutative
and generated by its primitive part, we can apply the characterization
theorem and get the claim.   \qed
\enddemo

\vskip7pt

\proclaim {Theorem 3.8} \, Assume that  $ \, \Bbbk := R \big/ \h R \, $ 
is a field.   
                                       \par   
   Let  $ \, B \in \calB \, $.  Then
$ \, {(B')}_\infty = {I_{\scriptscriptstyle \! B'}}^{\!\infty} \, $
and  $ \, B'{\Big|}_{\h=0} \, $  is commutative and has no non-trivial
idempotents.  In addition, when  $ \, p:= \Char(\Bbbk) > 0 \, $  each
non-zero element of  $ \, J_{B'{|}_{\h=0}} \, $  has nilpotency
order  $ p \, $,  \, that is  $ \, \overline{\eta}^{\,p} = 0 \, $
for all  $ \, \overline{\eta} \in J_{B'{|}_{\h=0}} \, $.
                                         \hfill\break
   \indent   In particular, if  $ \; H \in \HA \, $  then  $ \, H'
\in \QFA \, $.   
%
%
\endproclaim

\demo{Proof} The second part of the claim (about  $ \, H \in \HA \, $)
is simply a straightforward reformulation of the first part (about
$ \, B \in \calB \, $),  \, so in fact it is enough to prove the
latter.
                                       \par
   First we must show that  $ \, B'{\big|}_{\h=0} \, $
is commutative,  $ \, {(B')}_\infty = {I_{\scriptscriptstyle \!
B'}}^{\!\infty} \, $  and  $ \, B'{\big|}_{\h=0} \, $  has
no non-trivial idempotents (cf.~\S 1.3--4).  For later use,
set
       \hbox{$ \, I \! := I_{\scriptscriptstyle \! B} \, $,
$ J \! := \! J_{\scriptscriptstyle \! B} \, $,  $ J' \! :=
\! J_{\scriptscriptstyle \! B'} \, $,  $ I' \! := \!
I_{\scriptscriptstyle \! B'} \, $.}
                                       \par
   As for commutativity, we have to show that  $ \; a b - b a \in
\h \, B' \; $  for all  $ \, a $,  $ b \in B' \, $.  First, by the
inverse formula for  $ \, \Delta^n \, $  (see \S 2.1) we have  $ \,
\id_{\scriptscriptstyle B} = \Delta^1 = \delta_1 + \delta_0 =
\delta_1 + \epsilon \, $;  \, so  $ \, x = \delta_1(x) + \epsilon(x)
\, $  for all  $ \, x \in B \, $.  If  $ \, x \in B' \, $  we have
$ \, \delta_1(x) \in \h \, B \, $,  hence there exists  $ \, x_1 \in
B \, $  such that  $ \, \delta_1(x) = \h \, x_1 \, $.  Now take  $ \,
a $,  $ b \in B' \, $:  \, then  $ \; a = \h \, a_1 + \epsilon(a)
\, $,  $ \, b = \h \, b_1 + \epsilon(b) \, $,  \; whence  $ \; a b
- b a = \h \, c \; $  with  $ \; c = \h \, (a_1 b_1 - b_1 a_1) \, $;
\; therefore we are left to show that  $ \, c \in B' $.  To this
end, we have to check that  $ \, \delta_\Phi(c) \, $  is divisible
by  $ \, \h^{\vert \Phi \vert} \, $  for any nonempty finite subset
$ \Phi $  of  $ \N_+ \, $:  as multiplication by  $ \h $  is
injective (for  $ B $  is torsion-free!), it is enough to
show that  $ \, \delta_\Phi(a b - b a) \, $  is divisible
by  $ \, \h^{|\Phi| + 1} $.
                                            \par
   Let  $ \Lambda $  and  $ Y $  be subsets of  $ \Phi $  such that
$ \, \Lambda \cup Y = \Phi \, $  and  $\, \Lambda \cap Y \not=
\emptyset \, $:  then  $ \, \vert \Lambda \vert + \vert Y \vert
\geq |\Phi| + 1 \, $.  Now,  $ \, \delta_\Lambda(a) \, $  is
divisible by  $ \, \h^{|\Lambda|} \, $  and  $ \, \delta_Y(b) \, $
is divisible by  $ \, \h^{|Y|} $.  From this and from  Lemma
3.4{\it (b)}  it follows that  $ \, \delta_\Phi(a b - b a) \, $
is divisible by  $ \, \h^{|\Phi| + 1} $,  \, q.e.d.
                                            \par
   Second, we show that  $ \, {(B')}_\infty = {(I')}^\infty \, $.
By definition  $ \, \h \, B' \! \subseteq \! I' \, $,  whence  $ \,
B'_\infty := \bigcap_{n=0}^{+\infty} \h^n B' \! \subseteq \!
\bigcap_{n=0}^{+\infty} {\big( I' \big)}^n =: {\big( I'
\big)}^\infty $,  i.e.~$ {\big( B' \big)}_\infty \! \subseteq
\! {\big( I' \big)}^\infty $.  Conversely,  $ \, I' = \h \, B'
+ J' \, $  with  $ \, \h \, B' \subseteq \h \, B \, $  and  $ \,
J' = \delta_1(J') \subseteq \h \, B \, $:  \, thus  $ \, I' \!
\subseteq \! \h \, B \, $,  \, hence  $ \, {\big( I' \big)}^\infty
\! \subseteq \! \bigcap_{n=0}^{+\infty} \h^n B =: B_\infty \, $.
Now definitions give  $ \, B_\infty \! \subseteq \! B' \, $  and
$ \, \h^\ell B_\infty = B_\infty \, $  for all  $ \, \ell \in \Z
\, $,  \, so  $ \, \h^{-n} {\big( I' \big)}^\infty \! \subseteq \!
\h^{-n} B_\infty = B_\infty \! \subseteq \! B' \, $  hence  $ \,
{\big( I' \big)}^{\!\infty} \! \subseteq \! \h^n B' \, $  for all
$ \, n \in \N $,  thus finally  $ \, {\big( I' \big)}^{\!\infty}
\! \subseteq \! {\big( B' \big)}_\infty $.
                                            \par
   Third, we prove that  $ \, B'{\big|}_{\h=0} \, $  has no
non-trivial idempotents.
                                            \par
   Let  $ \, a \in B' $,  and suppose that  $ \; \overline{a} := a
\, \mod \h \, B' \in B'{\big|}_{\h=0} \; $  is idempotent, i.e.~$ \,
{\overline{a}}^{\,2} = \overline{a} \, $.  Then  $ \, a^2 = a +
\h \, c \, $  for some  $ \, c \in \h \, B' \, $.  Set  $ \, a_0
:= \epsilon(a) \, $,  $ \, a_1 := \delta_1(a) \, $,  and  $ \,
c_0 := \epsilon(c) \, $,  $ \, c_1 := \delta_1(c) \, $;  \,
since  $ \, a $,  $ c \in B' \, $  we have  $ \, a_1 $,
$ c_1 \in \h \, B \cap J = \h \, J \, $.
                                            \par
   First, applying  $ \, \delta_n \, $  to the identity  $ \; a^2
= a + \h \, c \; $  and using  Lemma 3.4{\it (a)}  we get
  $$  {\textstyle \sum\limits_{\Lambda \cup Y = \{1,\dots,n\}}}
\hskip-5pt  \delta_\Lambda(a) \, \delta_Y(a) \; = \; \delta_n\big(a^2\big)
\, = \, \delta_n(a) + \h \, \delta_n(c)   \eqno \forall \;\;\;
n \in \N_+ \;\; .  \qquad  (3.3)  $$
Since  $ \, a $,  $ c \in B' \, $  we have  $ \, \delta_n(a) $,
$ \delta_n(c) \in \h^n B^{\otimes n} \, $  for all  $ \, n \in
\N \, $.  Therefore (3.3) yields
  $$  \delta_n(a) \; \equiv  \hskip-1pt  {\textstyle \sum\limits_{\Lambda
\cup Y = \{1,\dots,n\}}}  \hskip-5pt  \delta_\Lambda(a) \, \delta_Y(a)
\; = \; 2 \, \delta_0(a) \, \delta_n(a) \; +  \hskip-5pt  {\textstyle
\sum\limits_{\Sb \Lambda \cup Y = \{1,\dots,n\} \\  \Lambda, Y \not=
\emptyset \endSb}}  \hskip-5pt  \delta_\Lambda(a) \, \delta_Y(a)
\hskip6pt  \mod \, \h^{n+1} B^{\otimes n}  $$
for all  $ \, n \in \N_+ \, $,  \, which, recalling that  $ \, a_0
:= \delta_0(a) \, $,  \, gives (for all  $ \, n \in \N_+ $)
  $$  \big( 1 - 2 \, a_0 \big) \; \delta_n(a) \hskip3pt \equiv
\hskip-1pt {\textstyle \sum\limits_{\Sb \Lambda \cup Y = \{1,\dots,n\} \\
\Lambda, Y \not= \emptyset \endSb}}  \hskip-7pt  \delta_\Lambda(a) \,
\delta_Y(a)  \hskip11pt  \mod \, \h^{n+1} B^{\otimes n} \;\; .  
\eqno (3.4)  $$
   \indent   Now, applying  $ \, \epsilon \, $  to the identity
$ \; a^2 = a + \h \, c \; $  gives  $ \; {a_0}^{\!2} = a_0 + \h \,
c_0 \, $.  This implies  $ \, \big( 1 - 2 \, a_0 \big) \not\in \h \,
B \, $:  \, this is trivial if  $ \, \text{\it Char}\,(R) = 2 \, $
or  $ \, a_0 = 0 \, $;  \, otherwise, if  $ \, \big( 1 - 2 \, a_0
\big) \in \h \, B \, $  then  $ \, a_0 = {1/2} + \h \, \alpha \, $
for some  $ \, \alpha \in B \, $,  \, and so  $ \; {a_0}^{\!2} =
{\big( {1/2} + \h \, \alpha \big)}^2 = {1/4} + \h \, \alpha + \h^2
\, \alpha^2 \not= {1/2} + \h \, \alpha + \h \, c_0 = a_0 + \h \, c_0
\, $,  \; thus contradicting the identity  $ \, {a_0}^{\!2} = a_0
+ \h \, c_0 \, $.  Now using  $ \, \big( 1 - 2 \, a_0 \big) \not\in
\h \, B \, $  and formulas (3.4)   --- for all  $ \, n \in \N_+ $
---   an easy induction argument gives  $ \, \delta_n(a) \in \h^{n+1}
B \, $,  \, for all  $ \, n \in \N_+ \, $.  Now consider  $ \, a_1 =
a - a_0 = \h \, \alpha \, $  for some  $ \, \alpha \in \h \, J \, $:
\, we have  $ \, \delta_0(\alpha) = \epsilon(\alpha) = 0 \, $  and
$ \, \delta_n(\alpha) = \h^{-1} \delta_n(a) \in \h^n B \, $,  \, for
all  $ \, n \in \N_+ \, $,  \, which mean  $ \; \alpha \in B' \, $.
Thus  $ \;\; a = a_0 + \h \, \alpha \equiv a_0 \, \mod \h \, B' \, $,
\, whence  $ \; \overline{a} = \overline{a_0} \in B'{\big|}_{\h=0}
\, $;  \, then  $ \; {\overline{a_0}}^{\,2} = \overline{a_0} \in \Bbbk
\, $  gives us  $ \, \overline{a_0} \in \{0,1\} \, $,  \, hence
$ \, \overline{a} = \overline{a_0} \in \{0,1\} \, $,  \, q.e.d.
                                            \par
   Finally, assume that  $ \, p := \hbox{\it Char}\,(\Bbbk) > 0 \, $;
\, then we have to show that  $ \, \overline{\eta}^{\,p} = 0 \, $  for
each  $ \, \overline{\eta} \in J_{B'{|}_{\h=0}} \, $,  \, or simply
$ \, \eta^p \in \h \, J_{\scriptscriptstyle B'} \, $  for each  $ \,
\eta \in J_{\scriptscriptstyle B'} \, $.  Indeed, for any  $ \, n
\in \N \, $  from the multiplicativity of  $ \Delta^n $  and from
$ \, \Delta^n(\eta) = \sum_{\Lambda \subseteq \{1,\dots,n\}}
\delta_\Lambda(\eta) \, $  (cf.~\S 2.1) we have
  $$  \displaylines{
 \textstyle
   \Delta^n(\eta^p) = {\big( \Delta^n(\eta) \big)}^p = \Big(
\sum\nolimits_{\Lambda \subseteq \{1,\dots,n\}} \delta_\Lambda(\eta)
\Big)^p  \; \in \;  \sum\nolimits_{\Lambda \subseteq \{1,\dots,n\}}
{\delta_\Lambda(\eta)}^p  \; +   \hfill  \cr
 \textstyle
   \quad   +  \; \sum\limits_{\Sb e_1, \dots, e_p < p \\   e_1 +
\cdots e_p = p \endSb}  \hskip-7pt  {p \choose e_1, \dots, e_p}
\hskip-3pt  \sum\limits_{\Lambda_1, \dots,\Lambda_p \subseteq
\{1,\dots,n\}}  \hskip-3pt  \prod_{k=1}^p {\delta_{\Lambda_k}
(\eta)}^{e_k} \; + \; \h \, \cdot \sum\limits_{k=0}^{n-1}
\sum\limits_{\Sb \Psi \subseteq \{1,\dots,n\}  \\  |\Psi|=k  \endSb}
\hskip-9pt  j_\Psi \big( {J_{\!{}_{B'}}}^{\!\!\otimes k} \big) \;
+ \; \h \, \cdot {J_{\!{}_{B'}}}^{\!\!\otimes n}  \cr }  $$
because  $ \, \delta_\Lambda(\eta) \in j_\Lambda \Big(
{J_{\!{}_{B'}}}^{\!\!\otimes |\Lambda|} \Big) \, $  (for
all  $ \, \Lambda \subseteq \{1,\dots,n\} \, $)  and
$ \, \big[J_{\!{}_{B'}},J_{\!{}_{B'}}\big] \subseteq
\h \, J_{\!{}_{B'}} \, $.  Then
  $$  \displaylines{
 \textstyle
   \delta^n(\eta^p) = {(\id_{\scriptscriptstyle B} \! -
\epsilon)}^{\otimes n} \big( \Delta^n(\eta^p) \big)
\, \in \,  {\delta_n(\eta)}^p \, + \hskip-11pt
\sum\limits_{\Sb e_1, \dots, e_p < p \\   e_1 + \cdots e_p = p \endSb}
\hskip-9pt  {p \choose {e_1, \dots, e_p}}  \hskip-3pt  \sum\limits_{\cup_k
\Lambda_k = \{1,\dots,n\}}  \hskip-13pt  \prod_{k=1}^p
{\delta_{\Lambda_k}(\eta)}^{e_k}  + \, \h \,
{J_{\!{}_{B'}}}^{\!\!\otimes n} .  \cr }  $$
Now,  $ \, {\delta^n(\eta)}^p \in {\big( \h^n B^{\otimes n} \big)}^p
\subseteq \h^{n+1} B^{\otimes n} \, $  because  $ \, \eta \in B' \, $,
\, and similarly we have also  $ \, \prod\nolimits_{k=1}^p
{\delta_{\Lambda_k}(\eta)}^{e_k} \in \h^{\sum_k |\Lambda_k| \, e_k}
B^{\otimes n} \subseteq \h^n B^{\otimes n} \, $  whenever  $ \,
\bigcup_{k=1}^n \Lambda_k = \{1,\dots,n\} \, $;  \, in addition, the
multinomial coefficient  $ \, {p \choose e_1, \dots, e_p} \, $  (with
$ \, e_1 $,  $ \dots $,  $ e_p < p \, $)  is a multiple of  $ p \, $,
hence it is zero in  $ \, R \big/ \h \, R = \Bbbk \, $,  that is
$ \, {p \choose e_1, \dots, e_p} \in \h \, R \, $:  \, therefore
  $$  \textstyle \sum\limits_{\Sb e_1, \dots, e_p < p \\   e_1 + \cdots
+ e_p = p \endSb}  {p \choose e_1, \dots, e_p}  \sum\limits_{\cup_{k=1}^p
\Lambda_k = \{1,\dots,n\}}  \, \prod\limits_{k=1}^p {\delta_{\Lambda_k}
(\eta)}^{e_k} \; \in \; \h^{n+1} B^{\otimes n} \;\; .  $$
Finally, since  $ \, J_{\!{}_{B'}} \subseteq \h \, J_{\!{}_B} \, $  we
have also  $ \, \h \, {J_{\!{}_{B'}}}^{\!\!\otimes n} \subseteq \h^{n+1}
B^{\otimes n} \, $.  The outcome then is that  $ \, \delta_n(\eta^p) \in
\h^{n+1} B^{\otimes n} \, $  for all  $ \, n \in \N \, $,  \, thus  $ \,
\eta \in \h \, B' \, $  as expected.   \qed
\enddemo

%
%
 \eject

\centerline {\bf \S \; 4 \  Drinfeld's functors on
quantum groups }

\vskip10pt

   {\sl From now on, we  {\bf assume that  $ \, \Bbbk :=
R \big/ \h R \, $  is a field}.}

\vskip7pt

\proclaim{Lemma 4.1} \, Let  $ \, F_\h \in \QFA \, $,  \, and assume
that  $ \, F_\h{\big|}_{\h=0} \, $  is reduced.  Let  $ \, I :=
I_{\scriptscriptstyle F_\h} \, $,  \, let  $ \widehat{F_\h} $
be the  $ I $--adic  completion of  $ F_\h \, $,  \, and
$ \widehat{I{\phantom|\!}^n} $  the  $ I $--adic  closure
of  $ I^n $  in  $ \widehat{F_\h} \, $,  for all
$ \, n \! \in \! \N \, $.
                                        \break
   \indent   (a) \;  $ \, \widehat{F_\h} \, $  is isomorphic as
an  $ \Rhat $--module  (where  $ \Rhat $  is the  $ \h $--adic
completion of  $ R \, $)  to a formal power series algebra of
type  $ \, \Rhat \big[ \big[ {\{Y_b\}}_{b \in \Cal{S}} \big]
\big] \, $  (where  $ \Cal{S} $  stands for some index set).
                                        \hfill\break
   \indent   (b) \; Letting  $ \, \nu \, \colon \, \Bbbk
\,\lhook\joinrel\relbar\joinrel\longrightarrow\, R \; $
be a section of the quotient map  $ \, R \longtwoheadrightarrow
R \big/ \h \, R =: \Bbbk \, $,  \, use it to identify
(set-theoretically)  $ \, \widehat{F_\h} \cong \Rhat \big[\big[
{\{Y_b\}}_{b \in \Cal{S}} \big]\big] \, $  with  $ \, \nu(\Bbbk)
\big[\big[ \{Y_0\} \cup {\{Y_b\}}_{b \in \Cal{S}} \big]\big] \, $
(with  $ \, \h \cong Y_0 \, $).  Then via such an identification
both  $ \, \widehat{\big( \widehat{I} \;{\big)}^n} \, $  and
$ \, \widehat{I{\phantom|\!}^n} \, $  coincide with the set of
all formal series of (least) degree  $ n $  (in the  $ Y_i $'s,
with  $ \, i \in \{0\} \cup \Cal{S} \, $),  for all  $ \, n
\in \N \, $.
                                        \hfill\break
   \indent   (c) \; There exist  $ \Bbbk $--module  isomorphisms
$ \; G_I \big( F_\h \big) \cong \Bbbk \big[ Y_0, {\{Y_b\}}_{b \in
\Cal{S}} \big] \cong G_{\widehat{I}} \big(\, \widehat{F} \,\big)
\, $  for the graded rings  associated to  $ \, F_\h \, $  and
$ \, \widehat{F_\h} \, $  with the  $ I $--adic  and the
$ \widehat{I} $--adic  filtration.
                                        \hfill\break
   \indent   (d) \; Let  $ \, \mu : F_\h \loongrightarrow
\widehat{F_\h} \, $  be the natural map.  Then  $ \; \mu(F_\h)
\bigcap \widehat{I{\phantom|\!}^n} = \mu(I^n) \, $  for all
$ \, n \in \N \, $.
\endproclaim

\demo{Proof} Let  $ \, F[G] \equiv F_\h{\big|}_{\h=0} := F_\h \big/
\h \, F_\h \, $,  \, and let  $ \, \widehat{F[G]} = F[[G]] \, $  be
the  $ \germ $--adic  completion of  $ F[G] $,  where  $ \, \germ
= \text{\sl Ker}\, \big( \epsilon_{\scriptscriptstyle F[G]} \big)
\, $  is the maximal ideal of  $ \, F[G] \, $  at the unit element
of  $ G $.  Then  $ \, I = \pi^{-1}(\germ) \, $,  \, the preimage
of  $ \germ $  under the specialization map  $ \, \pi \, \colon
F_\h \longtwoheadrightarrow F_\h \big/ \h \, F_\h = F[G] \, $.
Therefore  $ \pi $  induces a continuous epimorphism
$ \, \widehat{\pi} \, \colon \, \widehat{F_\h}
\longtwoheadrightarrow \widehat{F[G]} = F[[G]] \, $,  \,
which again is nothing but specialization at  $ \, \h = 0 \, $.
Note also that the ground ring of  $ \widehat{F_\h} $  is
$ \Rhat $,  because the ground ring of the  $ I $--adic
completion of a unital  $ R $--algebra  is the
$ (R \cap I) $--adic  completion of  $ R $,  and
the  $ \, R \cap I = \h \, R \, $.  Then of course
$ \widehat{F_\h} $  is also a topological  $ \Rhat $--module.
Moreover, by construction we have  $ \, {\big( {\widehat{F_\h}}
\, \big)}_\infty = \{0\} \, $.
                                                \par
   Now, let  $ \, {\{y_b\}}_{b \in \Cal{S}} \, $  be a
$ \Bbbk $--basis  of  $ \, \germ \big/ \germ^2 = Q\big(F[G]\big)
\, $;  \, by hypothesis  $ F[G] $  is reduced, thus  $ F[[G]] $
is just the formal power series algebra in the  $ y_b $'s,
i.e.~$ \, F[[G]] \cong \Bbbk \big[\big[ {\{Y_b\}}_{b \in
\Cal{S}} \big]\big] \, $.  For any  $ \, b \in \Cal{S} $,
\, pick a  $ \, j_b \in \pi^{-1} (y_b) \bigcap J \, $  (with
$ \, J := \text{\sl Ker}\,(\epsilon_{\scriptscriptstyle F_\h})
\, $),  \, and fix also a section  $ \; \nu \, \colon \, \Bbbk
\,\lhook\joinrel\relbar\joinrel\longrightarrow\, R \; $
of the quotient map  $ \, R \longtwoheadrightarrow R \big/
\h \, R = \Bbbk \, $  as in  {\it (b)}.  Using these, we can
define a continuous morphism of  $ \Rhat $--modules  $ \, \Psi
\, \colon \, \Rhat \big[\big[ {\{Y_b\}}_{b \in \Cal{S}} \big]\big]
\llongrightarrow \widehat{F_\h} \, $  mapping
%
%
%
%
   $ \, Y^{\,\underline{e}} := \prod_{b \in \Cal{S}}
Y_b^{\,\underline{e}(b)} \, $  to  $ \, j^{\,\underline{e}}
:= \prod_{b \in \Cal{S}} j_b^{\,\underline{e}(b)} \, $
for all  $ \, \underline{e} \in \N^{\Cal{S}}_f := \big\{\, \sigma
\in \N^{\Cal{S}} \,\big\vert\, \sigma(b) = 0 \;\; \text{\sl for
almost all} \;\; b \in \Cal{S} \,\big\} \, $  (hereafter, monomials
like the previous ones are  {\sl ordered\/}  w.r.t.~any fixed order
of the index set  $ \Cal{S} \, $).  In addition, using  $ \nu $  one
can identify (set-theoretically)  $ \, \Rhat \cong \nu(\Bbbk)[[Y_0]]
\, $  (with  $ \, \h \cong Y_0 \, $),  whence a bijection  $ \,
\nu(\Bbbk) \big[\big[ Y_0 \cup {\{Y_b\}}_{b \in \Cal{S}} \big]
\big] \cong \Rhat \big[\big[ {\{Y_b\}}_{b \in \Cal{S}} \big]\big]
\, $  arises.
                                          \par
   We claim that  $ \Psi $  is surjective.  Indeed, since  $ \,
{\big( \widehat{F_\h} \, \big)}_\infty = \{0\} \, $,  \, for any
$ \, f \in \widehat{F_\h} \, $  there is a unique  $ \, v_\h(f) \in
\N \, $  such that  $ \, f \in \h^{v_\h(f)} \widehat{F_\h} \setminus
\h^{v_\h(f)+1} \widehat{F_\h} \, $,  \, so  $ \; \widehat{\pi} \Big(
\h^{-v_\h(f)} f \Big) = \sum_{\underline{e} \in \N^{\Cal{S}}_f }
c_{\underline{e}} \cdot y^{\underline{e}} \; $  for some  $ \,
c_{\underline{e}} \in \Bbbk \, $  not all zero.  Then for  $ \,
f_1 := f - \h^{v_\h(f)} \cdot \sum_{\underline{e} \in \N^{\Cal{S}}_f }
\nu(c_{\underline{e}}) \cdot j^{\,\underline{e}} \, $  we have
$ \, v_\h(f_1) > v_\h(f) \, $.  Iterating, we eventually find for
$ f $  a formal power series expression of the type  $ \; f =
\sum_{n \in \N} \h^n \cdot \sum_{\underline{e} \in \N^{\Cal{S}}_f }
\nu(c_{\underline{e},n}) \cdot j^{\,\underline{e}} =
\sum_{(e_0,\underline{e}) \in \N \times \N^{\Cal{S}}_f }
\nu(\kappa_{\underline{e}}) \cdot \h^{e_0} j^{\,\underline{e}} \, $,
\, so  $ \, f \in \text{\it Im}\,(\Psi) \, $,  \, q.e.d.  Thus in
order to prove  {\it (a)\/}  we are left to show that  $ \Psi $
is injective too.
                                                \par
   Consider the graded ring associated to the  $ \h $--adic
filtration of  $ \widehat{F_\h} $,  that is  $ \, G_\h \big( \,
\widehat{F_\h} \, \big) := \oplus_{n=0}^{+\infty} \Big( \h^n
\widehat{F_\h} \big/ \h^{n+1} \widehat{F_\h} \Big) \, $:  \, this
is  {\sl commutative\/},  because  $ \, \widehat{F_\h} \Big/ \h \,
\widehat{F_\h} = \widehat{\pi} \big(\, \widehat{F_\h} \,\big) =
F[[G]] \, $  is commutative, and more precisely  $ \;  G_\h \big( \,
\widehat{F_\h} \, \big) \; \cong \; \Bbbk \big[ \,\overline{\h} \;\big]
\otimes_\Bbbk \Big( \widehat{F_\h} \Big/ \h \, \widehat{F_\h} \Big)
\; \cong \; \Bbb\Bbbk[Y_0] \otimes_\Bbbk
                       F[[G]] \; \cong $\break
$ \cong \big(F[[G]]\big)[Y_0] \; $  as graded  $ \Bbbk $--algebras.
In addition, the epimorphism (of  {\sl  $ \Rhat $--modules\/})
$ \; \Psi \, \colon \, \Rhat \otimes_\Bbbk \! F[[G]]
\relbar\joinrel\relbar\joinrel\twoheadrightarrow \!
\widehat{F_\h} \; $  induces an epimorphism  \hbox{$ \; G_{Y_0}(\Psi)
\, \colon \, G_{Y_0} \big( \Rhat \otimes_\Bbbk \! F[[G]] \big)
\relbar\joinrel\relbar\joinrel\twoheadrightarrow G_\h \big(
\widehat{F_\h} \,\big) $}  \; of  {\sl graded\/  $ \Bbbk $--algebras},
and by the very construction  $ \, G_{Y_0} (\Psi) \, $  is clearly
an isomorphism yielding  $ \; \big(F[[G]]\big)[Y_0] \cong G_{Y_0}
\big( \Rhat \otimes_\Bbbk F[[G]] \big) \cong G_\h \big(\,
\widehat{F_\h} \,\big) $:  \, then by a standard argument
(cf.~[Bo], Ch.~III, \S 2.8, Corollary 1) we conclude that
$ \, \Psi \, $  is an isomorphism as well, q.e.d.
                                                \par
   As for part  {\it (b)},  we start by noting that  $ \,
\widehat{I} = \widehat{\pi}^{-1}\big(\text{\sl Ker}\,
(\epsilon_{\scriptscriptstyle F[[G]]})\big) = \text{\sl Ker}\,
(\epsilon_{\scriptscriptstyle \widehat{F_\h}}) + \h \, \widehat{F_\h}
\, $,  \, so each element of  $ \widehat{I} $  is expressed   ---
via the isomorphism  $ \Psi $  ---   by a series of degree at least
1; moreover, for all  $ \, b $,  $ d \in \Cal{S} \, $  we have  $ \,
j_b \, j_d - j_d \, j_b = \h \, j_+ \, $  for some  $ \, j_+ \in
\text{\sl Ker}\, (\epsilon_{\scriptscriptstyle \widehat{F_\h}})
\, $.  This implies that when multiplying  $ n $  factors from
$ \widehat{I} $  expressed by  $ n $  series of positive degree,
we can reorder the unordered monomials in the  $ y_b $'s  occurring
in the multiplication process and eventually get a formal series
--- with  {\sl ordered\/}  monomials ---   of degree at least
$ n \, $.  This proves the claim for both  $ \, \widehat{I^n}
\, $  and  $ \, \widehat{\big( \widehat{I} \;{\big)}^n} \, $.
                                                \par
   For part  {\it (c)},  the analysis above shows that the natural
map  $ \, \mu : F_\h \longrightarrow \widehat{F_\h} \, $  induces
$ \Bbbk $--module  isomorphisms  $ \, {\big( \widehat{I} \;\big)}^n
\Big/ {\big( \widehat{I} \;\big)}^{n+1} \cong  \widehat{\big(
\widehat{I} \;{\big)}^n} \Big/  \widehat{\big( \widehat{I} \;
{\big)}^{n+1}} \cong \widehat{I^n} \Big/ \widehat{I{\phantom{|}
\!}^{n+1}} \, $  (for all  $ \, n \in \N \, $),  so  $ \, G_I
\big( F_\h \big) := \bigoplus\limits_{n=0}^{+\infty} I^n \Big/
I^{n+1} \cong \bigoplus\limits_{n=0}^{+\infty} {\big( \widehat{I}
\;\big)}^n \Big/ {\big( \widehat{I} \;\big)}^{n+1} =: G_{\widehat{I}}
\big(\, \widehat{F_\h} \,\big) \cong \bigoplus\limits_{n=0}^{+\infty}
\widehat{I{\phantom{|}\!}^n} \Big/ \widehat{I{\phantom{|}\!}^{n+1}}
\, $;  \, moreover, the given description of the  $ \widehat{I^n} $'s
implies  $ \, G_{\widehat{I}} \big(\, \widehat{F_\h} \,\big) :=
\bigoplus\limits_{n=0}^{+\infty} \widehat{I{\phantom{|}\!}^n} \Big/
\widehat{I{\phantom{|}\!}^{n+1}} \cong \Bbbk \big[ Y_0, {\{Y_b\}}_{b \in
\Cal{S}} \big] \, $  as  $ \Bbbk $--modules,  and the like for  $ \, G_I
\big( F_\h \big) $,  \, thus  {\it (c)\/}  is proved.
                                                \par
   Finally,  {\it (d)\/}  is a direct consequence of  {\it (c)\/}:  for
the latter yields  $ \Bbbk $--module  isomorphisms  $ \, F_\h \big/ I^n
\cong G_I \big( F_\h \big) \Big/ G_I \big( I^n \big) \cong G_I \big(\,
\widehat{F_\h} \,\big) \Big/ G_I \big(\, \widehat{I{\phantom{|}\!}^n}
\,\big) \cong \widehat{F_\h} \Big/ \widehat{I{\phantom{|}\!}^n} $,
\, thus  $ \, \mu(I^n) = \widehat{I{\phantom{|}\!}^n} \bigcap
\mu(F_\h) \, $.   \qed
\enddemo

\vskip7pt

   {\sl $ \underline{\hbox{\it Remark}} $:}  \; the previous description
of the ``formal quantum group''  $ \widehat{F_\h} $  shows that the
latter looks exactly like expected.  In particular, in the finite
dimensional case we can say it is a local ring which is also
``regular'', in the sense that the four numbers
                                                \par
   --- dimension of the ``cotangent space''  $ \,
I_{\scriptscriptstyle F_\h} \Big/ {I_{\scriptscriptstyle F_\h}}^{\!2}
\, $,
                                                \par
   --- least number of generators of the maximal ideal
$ I_{\scriptscriptstyle F_\h} \, $,
                                                \par
   --- Hilbert dimension ( = degree of the Hilbert
polynomial of the graded ring  $ \, G_{\widehat{I}}
\big(\, \widehat{F_\h} \,\big) \, $),
                                                \par
   --- Krull dimension of the associated graded ring
$ \, G_{\widehat{I}} \big(\, \widehat{F_\h} \,\big) \, $,
                                                \par
\noindent   are  {\sl all equal\/}.  Another way to say it is to
note that, if  $ \, \big\{ j_1, \dots, j_d \big\} \, $  is a lift
in  $ J_{\scriptscriptstyle F_\h} $  of any system of parameters
of  $ \, G = {\text{\it Spec}}\,\big( F_\h{\big|}_{\h=0} \big) \, $
around the identity (with  $ \, d = \dim(G) \, $),  \, then the set
$ \, \big\{ j_0 := \h, j_1, \dots, j_d \big\} \, $  is a ``system of
parameters'' for  $ \, F_\h \, $  (or, more precisely, for the local
ring  $ \widehat{F_\h} \, $).  A suggestive way to interpret all this
is to think at quantization as ``adding one dimension (or deforming)
in the direction of the quantization parameter  $ \h \, $'':  and
here we stress the fact that this is to be done ``in a  {\sl
regular\/}  way''.

\vskip7pt

\proclaim{Lemma 4.2}  Let  $ \, F_\h \in \, \QFA \, $,  \, and assume
that  $ F_\h{\big|}_{\h=0} $  is reduced.  Then:
                           \hfill\break
   \indent (a) \, if  $ \, \varphi \in F_\h \, $  and
$ \, \h^s \, \varphi \in {I_{\scriptscriptstyle F_\h}
\phantom{)}}^{\!\!\!\!\!n} \, $  ($ s, n \in \N $),  \,
then  $ \; \varphi \in {I_{\scriptscriptstyle F_\h}
\phantom{)}}^{\!\!\!\!\!n-s} $;
                        \hfill\break
   \indent (b) \, if  $ \, y \in I_{\scriptscriptstyle F_\h}
\setminus {I_{\scriptscriptstyle F_\h} \phantom{)}}^{\!\!\!\!\!2} $,
\, then  $ \; \h^{-1} y \not\in \h \, {F_\h}^{\!\vee} \, $;
                        \hfill\break
   \indent (c) \,  $ \; {\big(F_\h^{\,\vee}\big)}_\infty =
{\big( F_\h \big)}_\infty \; \big( = {I_{\scriptscriptstyle
F_\h}}^{\!\!\!\infty} \,\big) \, $.
                        \hfill\break
   \indent (d) \, Let  $ \, \text{\it Char}\,(\Bbbk) = 0 \, $,  and
let  $ \, U_\h \in \QrUEA \, $.  Let  $ \, x' \in {U_\h}' $,  and let
$ \, x \in U_\h \setminus \h \, U_\h \, $,  $ \, n \in \N \, $,  be
such that  $ \, x' = \h^n x \, $.  Set  $ \, \bar{x} := x \,\mod \h
\, U_\h \, $.  Then  $ \, \partial(\bar{x}) \leq n \, $  (hereafter
$ \, \partial(\bar{x}) \, $  is the degree of  $ \, \bar{x} \, $
w.r.t.~the standard filtration of the universal enveloping algebra
$ \, U_\h{\big|}_{\h=0} \, $).
\endproclaim

\demo{Proof}  {\it (a)} \, Set  $ \, I := I_{\scriptscriptstyle F_\h}
\, $.  Consider  $ I^\infty $  (cf.~Definition 1.4{\it (b)\/})
and the quotient Hopf algebra  $ \, \overline{F}_\h := F_\h \Big/
I^\infty \, $:  \, then  $ \, \bar{I} := I_{\scriptscriptstyle
\overline{F}_\h} = I \Big/ I^\infty \, $.  By  Lemma 3.2{\it (a)},
$ \overline{F}_\h $  is again a QFA, having the same specialization
at  $ \, \h=0 \, $  than  $ \, F_\h \, $,  and such that  $ \,
\bar{I}^\infty := {I_{\scriptscriptstyle \overline{F}_\h}
\phantom{)}}^{\hskip-7pt \infty} = \{0\} \, $.  Now,  $ \;
\phi \in I^\ell \,\Longleftrightarrow\, \overline{\phi} \in
\bar{I}^\ell \; $  for all  $ \, \phi \in F_\h \, $,  $ \, \ell
\in \N \, $,  \, with  $ \, \overline{\phi} := \phi + I^\infty \in
\overline{F}_\h \, $:  \, thus it is enough to make the proof for
$ \overline{F}_\h \, $,  \, i.~e.~we can assume from scratch that
$ \, I^\infty = \{0\} \, $.  In particular the natural map from
$ F_\h $  to its  $ I $--adic  completion  $ \widehat{F_\h} \, $
is injective, as its kernel is  $ I^\infty \, $.
                                                \par
   Consider the embedding  $ \, F_\h \hookrightarrow \widehat{F_\h}
\, $:  from the proof of Lemma 4.1 one easily sees that  $ \,
\widehat{I^\ell} \bigcap F_\h = I^\ell $,  \, for all  $ \, \ell
\, $  (because  $ \, F_\h \Big/ I^\ell \cong \widehat{F_\h} \Big/
\widehat{I^\ell} \; $):  \, then, using the description of
$ \widehat{I^\ell} $  in Lemma 4.1,  $ \, \varphi \in F_\h \, $
and  $ \, \h^s \varphi \in \! I^n \, $  give at once  $ \, \varphi
\in \widehat{I^{n-s}} \bigcap F_\h = I^{n-s} $,  \; q.e.d.
                                                \par
   {\it (b)} \, Let  $ \, y \in I_{\scriptscriptstyle F_\h} \setminus
{I_{\scriptscriptstyle F_\h}}^{\!\!2} \, $.  Assume  $ \, \h^{-1} y =
\h \, \eta \, $  for some  $ \, \eta \in {F_\h}^{\!\vee} \setminus
\{0\} \, $.  Since  $ \, {F_\h}^{\!\vee} :=
%
%
%
%
\bigcup_{N \geq 0} \h^{-N}
{I_{\scriptscriptstyle F_\h}}^{\!\!n} \, $  we have  $ \, \eta =
\h^{-N} i_N \, $  for some  $ \, N \in \N_+ \, $,  $ \, i_N \in
{I_{\scriptscriptstyle F_\h}}^{\!\!N} \, $.  Then we have  $ \,
\h^{-1} y = \h \, \eta = \h^{1-N} i_N \, $,  \, whence  $ \,
\h^{N-1} y = \h \, i_N \, $:  \, but the right-hand-side belongs
to  $ {I_{\scriptscriptstyle F_\h}}^{\!\!{N+1}} $,  whilst
the left-hand-side cannot belong to  $ {I_{\scriptscriptstyle
F_\h}}^{\!\!{N+1}} $,  due to  {\it (a)\/},  because  $ \, y
\not\in {I_{\scriptscriptstyle F_\h}}^{\!\!2} \, $,  \, a
contradiction.
                                                \par
   {\it (c)} \, Clearly  $ \, F_\h \subset {F_\h}^{\!\vee} \, $
implies  $ \, {\big( F_\h \big)}_\infty \! := \bigcap_{n=0}^{+\infty}
\h^n F_\h \! \subseteq \! \bigcap_{n=0}^{+\infty} \h^n {F_\h}^{\!\vee}
=: \! {\big( {F_\h}^{\!\vee} \big)}_\infty $.  For the converse
inclusion, note that by definitions  $ {\big( F_\h \big)}_\infty $
is a two-sided ideal both inside  $ F_\h $  and inside
$ {F_\h}^{\!\vee} $,  and  $ \, {\overline{F}_\h}^{\!\vee} \equiv
{\left( F_\h \Big/ {\big( F_\h \big)}_\infty \right)}^{\!\vee}
= {F_\h}^{\!\vee} \Big/ {\big( F_\h \big)}_\infty \, $,  \, so we
have also  $ \, {\big( {F_\h}^{\!\vee} \big)}_\infty \!\mod\!
{\big( F_\h \big)}_\infty \subseteq \! {\left( {F_\h}^{\!\vee}
\Big/ \! {\big( F_\h \big)}_\infty \right)}_\infty \!\!\! = \!
{\left( {\overline{F}_\h}^{\!\vee} \right)}_\infty \, $,  \,
with  $ \, \overline{F}_\h := F_\h \Big/ {I_{\scriptscriptstyle
F_\h}}^{\!\!\!\infty} = F_\h \Big/ {\big( F_\h \big)}_\infty \, $
(a QFA, by  Lemma 3.2{\it (a)}).  So, we prove that  $ \, {\left(
{\overline{F}_\h}^{\!\vee} \right)}_\infty \!\! = \! \{0\} \, $
for then  $ \, {\big( {F_\h}^{\!\vee} \big)}_\infty \subseteq
{\big( F_\h \big)}_\infty \, $  will follow.
                                                \par
   Let  $ \, \mu \, \colon \, F_\h \,\longrightarrow\,
\widehat{F_\h} \, $  be the natural map from  $ F_\h $  to its
$ I_{\scriptscriptstyle F_\h} $--adic  completion, whose kernel
is  $ \, {I_{\scriptscriptstyle F_\h}}^{\!\!\!\infty} = {\big(
F_\h \big)}_\infty \, $:  \, this makes  $ \overline{F}_\h $
embed into  $ \, \widehat{F_\h} \, $,  \, and gives  $ \,
{\overline{F}_\h}^{\!\vee} \subseteq {\widehat{F_\h}}^{\,\vee}
:= \bigcup_{n \geq 0} \h^{-n} \widehat{I^n} \, $  (notation of
Lemma 4.1), whence  $ \, {\left( {\overline{F}_\h}^{\!\vee}
\right)}_\infty \subseteq {\left( {\widehat{F_\h}}^{\,\vee}
\right)}_\infty \, $.  Now, the description of
$ \widehat{F_\h} $  and  $ \widehat{I^n} $  in
Lemma 4.1 yields that  $ {\widehat{F_\h}}^{\,\vee} $
is contained in the  $ \h $--adic  completion of the
$ R $--subalgebra  of  $ \, F(R) \otimes_R \widehat{F_\h} \, $
generated by  $ \, {\big\{ \h^{-1} j_b \big\}}_{b \in \Cal{S}} \, $
(as in the proof of Lemma 4.1), which is a  {\sl polynomial\/}
algebra.  But then  $ {\widehat{F_\h}}^{\,\vee} $  is separated
in the  $ \h $--adic  topology,  i.e.~$ \, {\left(
{\widehat{F_\h}}^{\, \vee} \right)}_\infty = \{0\} \, $.
                                                \par
  {\it (d)} \, (cf.~[EK], Lemma 4.12)  By hypothesis  $ \,
\delta_{n+1}(x') \in \h^{n+1} {U_\h}^{\otimes (n+1)} $,  \,
whence  $ \, \delta_{n+1}(x) \in \h \, {U_\h}^{\otimes (n+1)} $,
\, so  $ \, \delta_{n+1}(\bar{x}) = 0 \, $,  i.e.~$ \, \bar{x}
\in \hbox{\sl Ker}\, \Big( \delta_{n+1} \colon \, U(\gerg)
\longrightarrow {U(\gerg)}^{\otimes (n+1)} \Big) $,  where
$ \gerg $  is the Lie bialgebra such that  $ \, U_\h{\big|}_{\h=0}
:= U_\h \Big/ \h \, U_\h = U(\gerg) \, $.  But since  $ \,
\text{\it Char}\,(\Bbbk) = 0 \, $,  \, the latter kernel
equals  $ \, {U(\gerg)}_n := \big\{\, \bar{y} \in U(\gerg)
\,\big\vert\, \partial(\bar{y}) \leq n \,\big\} \, $
(cf.~[KT], \S 3.8), whence the claim.   \qed
\enddemo

\vskip7pt

\proclaim {Proposition 4.3}  Let  $ \, \text{\it Char}\,(\Bbbk)
= 0 \, $.  Let  $ \, F_\h \in \QFA \, $.  Then  $ \, {\big(
{F_\h}^{\!\vee} \big)}' \! = F_\h \, $.
\endproclaim

\demo{Proof}  Theorem 3.6 gives  $ \, F_\h \subseteq {\big(
{F_\h}^{\!\vee} \big)}' $,  \, so we have to prove only the
converse.  Let  $ \, \overline{F}_\h := F_\h \Big/ {(F_\h)}_\infty
\, $;  \, by  Lemma 4.2{\it (c)\/}  we have  $ \, {\big( {F_\h}^{\!\vee}
\big)}_\infty = {(F_\h)}_\infty \, $;  \, by Lemma 3.2{\it (a)\/}
we have  $ \, {\big( \overline{F}_\h \big)}^{\!\vee} =
{F_\h}^{\!\vee} \Big/ {(F_\h)}_\infty = {F_\h}^{\!\vee} \Big/ {\big(
{F_\h}^{\!\vee} \big)}_\infty \, $,  \, whence again by Lemma
3.2{\it (a)\/}  we get  $ \, {\Big(\! {\big( \overline{F}_\h
\big)}^{\!\vee} \Big)}' = {\Big( {F_\h}^{\!\vee} \Big/ {\big(
{F_\h}^{\!\vee} \big)}_\infty \Big)}' = {\big( {F_\h}^{\!\vee}
\big)}' \Big/ {\big( {F_\h}^{\!\vee} \big)}_\infty = {\big(
{F_\h}^{\!\vee} \big)}' \Big/ {(F_\h)}_\infty \, $.  Thus, if the
claim is true for  $ \, \overline{F}_\h \, $  then  $ \, F_\h \Big/
{(F_\h)}_\infty =: \overline{F}_\h = {\Big(\! {\big( \overline{F}_\h
\big)}^{\!\vee} \Big)}' = {\big( {F_\h}^{\!\vee} \big)}' \Big/
{(F_\h)}_\infty \, $,  \, whence clearly  $ \, {\big( {F_\h}^{\!\vee}
\big)}' = F_\h\, $.  Therefore it is enough to prove the claim
for  $ \, \overline{F}_\h \, $:  \, in other words, we can assume
$ \, {I_{\scriptscriptstyle F_\h}}^{\hskip-5pt \infty} =
{(F_\h)}_\infty = {\big( {F_\h}^{\!\vee} \big)}_\infty
= \{0\} \, $  (see  Lemma 4.2{\it (c)\/}).  In the
sequel, set  $ \, I := I_{\scriptscriptstyle F_\h} \, $.
                                               \par
   Let  $ \, x' \in {\big( {F_\h}^{\!\vee} \big)}' \, $  be given;
since  $ \, {(F_\h)}_\infty = \{0\} \, $  there are  $ \, n \in
\N \, $  and  $ \, x^\vee \in {F_\h}^{\!\vee} \setminus \h \,
{F_\h}^{\!\vee} \, $  such that  $ \, x' = \h^n x^\vee \, $.
By Theorem 3.7,  $ {F_\h}^{\!\vee} $  is a QrUEA, with
semiclassical limit  $ U(\gerg) $  where the Lie bialgebra
$ \gerg $  is  $ \, \gerg = I^\vee \Big/ \big( \h \, {F_\h}^{\!\vee}
\bigcap I^\vee \big) \, $,  \, with  $ \, I^\vee := \h^{-1} I \, $.
                                             \par
   Fix an ordered basis  $ \, {\{b_\lambda\}}_{\lambda \in \Lambda}
\, $  of  $ \gerg $  over  $ \Bbbk $,  and fix also a subset  $ \,
{\big\{ x^\vee_\lambda \big\}}_{\lambda \in \Lambda} \, $  of
$ {I_{\scriptscriptstyle F_\h}}^{\hskip-5pt \vee} $  such that
$ \, x^\vee_\lambda  \hskip-1,4pt  \mod \h \, {F_\h}^{\!\vee} =
b_\lambda \, $  for all  $ \, \lambda \in \Lambda \, $:  \, so
$ \, x^\vee_\lambda = \h^{-1} x_\lambda \, $  for some  $ \,
x_\lambda \in J $,  \, for all  $ \lambda \, $.
                                             \par
   Lemma 4.2{\it (d)}  gives  $ \, d:= \partial(\bar{x}) \leq n
\, $,  \, so we can write  $ \overline{x^\vee} $  as a polynomial
$ \, P \big( {\{b_\lambda\}}_{\lambda \in \Lambda} \big) \, $  in
the  $ b_\lambda $'s  of degree  $ \, d \leq n \, $;  \, hence
$ \, x^\vee \equiv P \big( {\big\{ x^\vee_\lambda \big\}}_{\lambda
\in \Lambda} \big) \mod\, \h \, {F_\h}^{\!\vee} \, $,  \, so  $ \,
x^\vee = P \big( {\big\{ x^\vee_\lambda \big\}}_{\lambda \in
\Lambda} \big) + \h \, x^\vee_{[1]} \, $  for some  $ \, x^\vee_{[1]}
\in {F_\h}^{\!\vee} $.  Now  $ \, x' = \h^n x^\vee = \h^n P \big(
{\big\{ x^\vee_\lambda \big\}}_{\lambda \in \Lambda} \big) +
\h^{n+1} x^\vee_{[1]} \, $  with
  $$  \h^n P \big( {\big\{ x^\vee_\lambda \big\}}_{\lambda \in
\Lambda} \big) = \h^n P \Big( {\big\{ \h^{-1} \, x_\lambda
\big\}}_{\lambda \in \Lambda} \Big) \, \in \, F_\h  $$
because  $ P $  has degree  $ \, d \leq n \, $;  \, thus since
$ \, F_\h \subseteq {\big( {F_\h}^{\!\vee} \big)}' \, $  (by
Theorem 3.6) we get
  $$  x'_1 := x' - \h^n P \big( {\big\{ x^\vee_\lambda
\big\}}_{\lambda \in \Lambda} \big) \in {\big( {F_\h}^{\!\vee}
\big)}'  \qquad  \hbox{and}  \qquad  x'_1 = \h^{n+1} x^\vee_{[1]}
= \h^{n_1} x^\vee_1  $$
for some  $ \, n_1 \in \N $,  $ \, n_1 > n \, $,  \, and some
$ \, x^\vee_1 \in {F_\h}^{\!\vee} \setminus \h \, {F_\h}^{\!\vee} $.
Therefore, we can repeat this construction with  $ x'_1 $  instead
of  $ x' $,  $ \, n_1 $  instead of  $ n $,  and  $ x^\vee_1 $
instead of  $ x^\vee $,  and so on.  Iterating, we eventually get
an increasing sequence  $ \, {\big\{ n_s \big\}}_{s \in \N} \, $
of natural numbers and a  sequence  $ \, {\Big\{ P_s \big( {\{
X_\lambda \}}_{\lambda \in \Lambda} \big) \Big\}}_{s \in \N} \, $
of polynomials such that the degree of  $ P_s \big( {\{ X_\lambda
\}}_{\lambda \in \Lambda} \big) $  is at most  $ n_s $,  for all
$ \, s \in \N $,  \, and   $ \, x' = \sum_{s \in \N} \h^{n_s} P_s
\big( {\big\{ x^\vee_\lambda \big\}}_{\lambda \in \Lambda}
\big) \, $.
                                             \par
   How should we look at the latter formal series?  By
construction, each one of the summands  $ \, \h^{n_s} P_s
\big( {\big\{ x^\vee_\lambda \big\}}_{\lambda \in \Lambda}
\big) \, $  belongs to  $ F_\h \, $:  \, more precisely,
$ \; \h^{n_s} P_s \big( {\big\{ x^\vee_\lambda \big\}}_{\lambda
\in \Lambda} \big) \in {I_{\scriptscriptstyle F_\h}}^{\hskip-3pt
n_s} \; $  for all  $ \, s \in \N \, $;  \; this means that
$ \, \sum_{s \in \N} \h^{n_s} P_s \big( {\big\{ x^\vee_\lambda
\big\}}_{\lambda \in \Lambda} \big) \, $  is a well-defined element
of  $ \widehat{F_\h} $,  the  $ I_{\scriptscriptstyle F_\h} $--adic
completion of $ F_\h $,  and the formal expression  $ \, x' =
\sum_{s \in \N} \h^{n_s} P_s \big( {\big\{ x^\vee_\lambda
\big\}}_{\lambda \in \Lambda} \big) \, $  is an identity
in  $ \widehat{F_\h} $.  So we find  $ \, x' \in {\big(
{F_\h}^{\!\vee} \big)}' \bigcap \widehat{F_\h} \, $.
%
%
   Now, consider the embedding  $ \, \mu \, \colon \, F_\h
\hookrightarrow \widehat{F_\h} \, $  and the specialization
map  $ \, \pi \, \colon \, F_\h \relbar\joinrel\twoheadrightarrow
F_\h{\big|}_{\h=0} = F[G] \, $:  \; like in the proof of Lemma 4.1,
$ \pi $  extends by continuity to  $ \, \widehat{\pi} \, \colon \,
\widehat{F_\h} \longtwoheadrightarrow \widehat{F[G]} = F[[G]] \, $:
\, then one easily checks that the map  $ \, \mu{\big|}_{\h=0}
\, \colon \, F[G] = F_\h{\big|}_{\h=0} \relbar\joinrel\llongrightarrow
\widehat{F_\h}{\big|}_{\h=0} = F[[G]] \, $  is injective too.  Since
$ \, \hbox{\sl Ker}\,(\pi) = \h \, F_\h \, $  and  $ \, \hbox{\sl Ker}
\,(\widehat{\pi}) = \h \, \widehat{F_\h} \, $,  \, this implies  $ \,
F_\h \bigcap \, \h \, \widehat{F_\h} = \h \, F_\h \, $,  \, whence $ \,
F_\h \bigcap \h^\ell \, \widehat{F_\h} = \h^\ell F_\h \, $  for all
$ \, \ell \in \N \, $.  Getting back to our  $ \, x' \in {\big(
{F_\h}^{\!\vee} \big)}' \bigcap \widehat{F_\h} \, $,  \, we have
$ \, x' = \h^{-n} y \, $  for some  $ \, n \in \N \, $  and  $ \,
y \in F_\h \, $;  \, thus, we conclude that  $ \, y = \h^n
x' \in F_\h \bigcap \h^n \, \widehat{F_\h} = \h^n F_\h\, $,
\, so that  $ \, x' \in F_\h \, $,  \; q.e.d.   \qed
\enddemo

\vskip7pt

\proclaim{Proposition 4.4} Let  $ \, H, K \! \in \! \HA \, $,
\, and  $ \; \langle \,\ , \,\ \rangle \, \colon H \times K
\loongrightarrow R \, $  a Hopf pairing.  Then
                                      \hfill\break
   \indent   (a)  $ \; H^\vee \subseteq {\big( K' \big)}^\bullet \; $
and  $ \; K' \subseteq {\big( H^\vee \big)}^\bullet \; $  (and
viceversa).  Therefore, the above pairing induces a Hopf pairing
$ \; \displaystyle{ \langle \,\ ,\,\ \rangle \colon \, H^\vee
\! \times K' \loongrightarrow R } \, $.
                                            \hfill\break
   \indent   (b) \, If in addition the pairing  $ \, H \times K
\longrightarrow R \, $  and its specialization
$  \, H{\big|}_{\h=0} \times K{\big|}_{\h=0}
                 \longrightarrow \Bbbk \, $\break
       at  $ \h = 0 $  are both perfect, and  $ \, K = {H}^\bullet \, $,
\, then we have also  $ \; K' = {\big( H^\vee \big)}^\bullet \, $.
                                            \hfill\break
   \indent   (c) \, Similar results hold for  $ \, B, \varOmega \in
\calB \, $  and  $ \; \langle \,\ , \,\ \rangle \, \colon B \times
\varOmega \loongrightarrow R \, $  a bialgebra pairing (i.e.~a pairing
with the properties of  Definition 1.2{\it(a)\/}  but the one about
the antipode).
\endproclaim

\demo{Proof}  {\it (a)} \, The Hopf pairing  $ \, H_F \times K_F
\longrightarrow F(R) \, $  given by scalar extension clearly
restricts to a similar Hopf pairing  $ \, H^\vee \! \times K'
\longrightarrow F(R) \, $:  \, we must prove this takes values
in  $ \, R \, $.
                                                \par
   Let  $ \, I = I_{\scriptscriptstyle H} \, $,  \, so  $ \,
H^\vee \! = \bigcup_{n=0}^\infty \h^{-n} I^n \, $  (cf.~\S 2.1).
Pick  $ \, c_1 $,  $ \dots $,  $ c_n \in I \, $,  $ \, y \in
K' \, $:  \, then
  $$  \displaylines{
   \Big\langle {\textstyle \prod_{i=1}^n} c_i \, , \, y
\Big\rangle = \Big\langle \otimes_{i=1}^n c_i \, , \, \Delta^n(y)
\Big\rangle = \left\langle \otimes_{i=1}^n c_i \, , \, {\textstyle
\sum_{\Psi \subseteq \{1,\dots,n\}}} \hskip-1pt \delta_\Psi(y)
\right\rangle =   \hfill  \cr
   = {\textstyle \sum_{\Psi \subseteq \{1,\dots,n\}}} \Big\langle
\otimes_{i=1}^n c_i \, , \, \delta_\Psi(y) \Big\rangle =
{\textstyle \sum_{\Psi \subseteq \{1,\dots,n\}}} \Big\langle
\otimes_{i \in \Psi} c_i \, , \, \delta_{\vert\Psi\vert}(y)
\Big\rangle \cdot {\textstyle \prod_{j \not\in \Psi}}
\langle c_j \, , 1 \rangle \; \in  \cr
   \hfill   \in \, {\textstyle \sum_{\Psi \subseteq \{1,\dots,n\}}}
\h^{n - |\Psi|} R \cdot \h^{\vert\Psi\vert} R = \h^n R \; .  \cr }  $$
  The outcome is  $ \; \big\langle I^n, K' \big\rangle \subseteq
\h^n R \, $,  \, whence  $ \; \big\langle \h^n I^n, K' \big\rangle
\subseteq R \, $,  \, for all  $ \, n \in \N \, $;  \, since
$ \, H^\vee = \bigcup_{n=0}^\infty \h^{-n} I^n \, $,  \, we get
$ \, H^\vee \subseteq {\big( K' \big)}^\bullet $  and  $ \, K'
\subseteq {\big( H^\vee \big)}^\bullet \, $:  \, then it follows
also that the restricted pairing  $ \, H^\vee \! \times K'
\loongrightarrow F(R) \, $  does take values in  $ R $,
as claimed.
                                             \par
   {\it (b)} \, We revert the previous argument to show that
$ \, {\big( H^\vee \big)}^\bullet \subseteq K' \, $.
                                             \par
   Let  $ \, \psi \in {\big( H^\vee \big)}^\bullet \, $:  then  $ \,
\big\langle \h^{-s} I^s \!, \, \psi \big\rangle \in R \, $  so  $ \,
\big\langle I^s \! , \, \psi \big\rangle \in \h^s R $,  for all
$ s $.  For  $ s \! = \! 0 $  we get  $ \, \big\langle H, \,
\psi \big\rangle \in R \, $,  thus  $ \, \psi \in H^\bullet
\! = K \, $  \, 
 and so  $ \, \delta_n(\psi) \in K^{\otimes n} \, $  for all  $ n \, $. 
If  $ \, n \in \N $,  $ \, i_1, \dots, i_n \in I $,  
  $$  \displaylines{
   {} \hskip7pt   \big\langle \otimes_{k=1}^n i_k \, , \,
\delta_n(\psi) \big\rangle \; = \; {\textstyle \sum_{\Psi \subseteq
\{1,\dots,n\}}} {(-1)}^{n - |\Psi|} \cdot \big\langle {\textstyle
\prod_{k \in \Psi}} i_k \, , \psi \big\rangle \cdot {\textstyle
\prod_{k \not\in \Psi}} \epsilon(i_k) \; \in  \qquad \qquad  \cr
   {} \hfill   \in \, {\textstyle \sum_{\Psi \subseteq  
\{1,\dots,n\}}} \! \big\langle I^{|\Psi|}, \psi \big\rangle
\cdot \h^{n - |\Psi|} \, R \; \subseteq \; {\textstyle \sum_{s=0}^n}
\h^s \cdot \h^{n-s} R = \h^n R \hskip7pt  \cr }  $$
therefore  $ \, \big\langle I^{\otimes n}, \, \delta_n(\psi)
\big\rangle \subseteq \h^n R \, $.  Now,  $ H $  splits as
$ \, H = R \cdot 1_{\scriptscriptstyle H} \oplus
J_{\scriptscriptstyle H} \, $,  with  $ \, J_{\scriptscriptstyle
H} \! := \hbox{\sl Ker}\,(\epsilon_{\scriptscriptstyle H}) \, $;
then  $ \, H^{\otimes n} \, $  splits into direct sum of  $ \,
{J_{\scriptscriptstyle H}}^{\!\otimes n} \, $  plus other direct
summands which are tensor products with at least one tensor
factor  $ \, R \cdot 1_{\scriptscriptstyle H} \, $.  As
$ \, J_{\scriptscriptstyle K} \! := \hbox{\sl Ker}\,
(\epsilon_{\scriptscriptstyle K}) = \big\{\, y \in \! K
\,\big\vert\, \langle 1_{\scriptscriptstyle H}, y \rangle =
0 \,\big\} $,  \, we have  $ \, \big\langle H^{\otimes n}, \,
j^\otimes \big\rangle = \big\langle {J_{\scriptscriptstyle H}}^{\!
\otimes n}, \, j^\otimes \big\rangle \, $  for  $ \, j^\otimes
\in {J_{\scriptscriptstyle K}}^{\!\otimes n} \, $.  Now  $ \,
\delta_n(\psi) \in {J_{\scriptscriptstyle K}}^{\!\otimes n} \, $:
\, this and the previous analysis together give  $ \, \big\langle
H^{\otimes n}, \, \delta_n(\psi) \big\rangle \subseteq \big\langle
{I_{\scriptscriptstyle H}}^{\otimes n}, \, \delta_n(\psi) \big\rangle
\subseteq \h^n \, R \, $,  for all  $ \, n \in \N \, $. 
                                             \par
   Now,  $ \, H^\bullet = K \, $  implies  $ \, {\big( H^{\otimes n}
\big)}^\bullet = K^{\otimes n} \, $  for the induced pairing  $ \;
H^{\otimes n} \times K^{\otimes n} \longrightarrow R \, $.  On the
other hand,  $ \, \Big\langle H^{\otimes n}, \, \delta_n(\psi)
\Big\rangle \subseteq \h^n R \, $  (for all  $ n \, $)  gives
$ \, \h^{-n} \delta_n(\psi) \in {\big( H^{\otimes n} \big)}^\bullet
= K^{\otimes n} \, $,  \, that is  $ \, \delta_n(\psi) \in \h^n
K^{\otimes n} \, $  for all  $ \, n \in \N \, $,  whence
finally  $ \, \psi \in K' \, $,  q.e.d.
                                             \par
   {\it (c)} \, We don't need antipode to prove  {\it (a)\/}  and
{\it (b)\/}:  the like arguments prove  {\it (c)\/}  too.   \qed
\enddemo

\vskip7pt

\proclaim{Proposition 4.5} \, Let  $ \, \text{\it Char}\,(\Bbbk) = 0
\, $.  Let  $ \, U_\h \in \QrUEA \, $.  Then  $ \, {\big( {U_\h}' \,
\big)}^{\!\vee} = U_\h \, $.
\endproclaim

\demo{Proof}  First, let  $ \, \overline{U}_\h := U_\h \Big/ {\big(
U_\h \big)}_\infty \, $,  \, and assume the claim holds for  $ \,
\overline{U}_\h \, $:  \, then repeated applications of  Lemma
3.2{\it (a)}  give  $ \, U_\h \Big/ {\big(U_\h\big)}_\infty =
\overline{U}_\h = \Big( \! \big( \overline{U}_\h \,\big)'
\Big)^{\!\vee} = {\big( {U_\h}' \big)}^{\!\vee} \Big/ {\big(
U_\h \big)}_\infty \, $,  \, whence  $ \, {\big( {U_\h}' \,
\big)}^{\!\vee} = U_\h \, $  follows at once; therefore we
are left to prove the claim for  $ \, \overline{U}_\h \, $,  \,
which means we may assume  $ \, {\big( U_\h \big)}_\infty = \{0\}
\, $.  In order to simplify notation, we set  $ \, H := U_\h \, $.
                                             \par
   Our purpose now is essentially to resort to a similar result
which holds for quantum groups ``\`a la Drinfeld'': so we mimic
the procedure followed in [Ga4] (in particular Proposition 3.7
therein), noting in addition that in the present case we can get
rid of the hypotheses  $ \, \dim(\gerg) < +\infty \, $  (with
$ \, U(\gerg) = U_\h \big/ \h \, U_\h \, $),  as one can check
getting through the entire procedure developed in [Ga4] in
light of  [{\it loc.~cit.}],  \S 3.9.
                                             \par
   Let  $ \, \widehat{H} \, $  be the  $ \h $--adic  completion
of  $ H $:  this is a separated complete topological
$ \Rhat $--module,  $ \Rhat $  being the  $ \h $--adic 
completion of  $ R \, $,  and a topological Hopf algebra, whose
coproduct takes values into  $ \, \widehat{H} \, \widehat{\otimes}
\, \widehat{H} := H \, \widehat{\otimes} \, H \, $,  \, the
$ \h $--adic  completion of  $ \, H \otimes H \, $  (indeed,
$ \widehat{H} $  is a  {\sl quantized universal enveloping
algebra\/}  in the sense of Drinfeld).  As  $ \, H_\infty =
\{0\} \, $,  \, the natural map  $ \, H \longrightarrow
\widehat{H} \, $  embeds  $ H $  as a (topological) Hopf
$ R $--subalgebra  of  $ \widehat{H} $.  Then we set also
$ \; \widehat{H}' := \big\{\, \eta \in \widehat{H} \,\big\vert\,
\delta_n (\eta) \in \h^n \widehat{H}^{\,\widehat{\otimes}\, n}
\, \big\} \; $  and  $ \; \big( \widehat{H}' \big)^{\!\times}
\! := \bigcup_{n \geq 0} \h^{-n} {I_{\scriptscriptstyle \!
\widehat{H}'}}^{\hskip-3pt n} \; \Big( \! \subseteq Q \big(
\Rhat \, \big) \otimes_{\Rhat} \widehat{H} \, \Big) $,  \,
where  $ \, I_{\scriptscriptstyle \! \widehat{H}'} := \hbox{\sl
Ker}\,(\epsilon_{\scriptscriptstyle \! \widehat{H}'}) + \h \cdot
\widehat{H}' \, $  (as in \S 1.3), and we let  $ \, \big( \widehat{H}'
\big)^{\!\vee} \, $
    \hbox{be the  $ \h $--adic  completion of
$ \big( \widehat{H}' \big)^{\!\times} \, $.}
                                        \par
   Now consider  $ \, \widehat{K} := {\widehat{H}}^* \equiv
{\hbox{\sl Hom}}_{\Rhat} \Big( \widehat{H} \, , \, \Rhat \,\Big)
\, $,  \, the dual of  $ \widehat{H} \, $:  \, this is a topological
Hopf  $ \Rhat $--algebra,  w.r.t.~the weak topology, in natural
perfect Hopf pairing with  $ \widehat{H} \, $:  \, in Drinfeld's
terminology, it is a  {\sl quantized formal series Hopf algebra}.
We define  $ \, \widehat{K}^\times := \sum_{n \geq 0} \h^{-n}
{J_{\scriptscriptstyle \! \widehat{K}}}^{\hskip-3pt n} \; \Big(
\! \subseteq Q \big( \Rhat \, \big) \otimes_{\Rhat} \widehat{K}
\, \Big) $,  \, where  $ \, J_{\scriptscriptstyle \! \widehat{K}}
:= \hbox{\sl Ker}\,(\epsilon_{\scriptscriptstyle \! \widehat{K}})
\, $  (as in \S 1.3)  and we let  $ \widehat{K}^\vee $  be the
$ \h $--adic  completion of  $ \widehat{K}^\times $,  \, and we
define  $ \, {\big( \widehat{K}^\vee \big)}' \, $  in the obvious
way.  With much the same arguments used for Proposition 4.3,
one proves that  $ \, {\big( \widehat{K}^\vee \big)}' =
\widehat{K} \, $.  Like in [Ga4], one proves   --- with much
the same arguments as for Proposition 4.4 ---   that  $ \,
\widehat{H}' = {\big( {\widehat{K}}^\vee \big)}^{\!\bullet} \, $
and  $ \, \widehat{K}^\vee \subseteq {\big( {\widehat{H}}'
\big)}^{\!\bullet} \, $;  \, moreover, one has also  $ \, \widehat{H}
= \widehat{K}^* $,  \, whence one argues  $ \, \widehat{K}^\vee
= {\big( {\widehat{H}}' \big)}^{\!\bullet} \, $.  Using this and the
equality  $ \, {\big( \widehat{K}^\vee \big)}' = \widehat{K} \, $
one proves  $ \, \big( \widehat{H}' \big)^{\!\vee} = \widehat{H}
\, $  as well (see [Ga4] for details).
                                        \par
   Now, definitions imply  $ \, \widehat{H} \Big/ \h^n \widehat{H} =
H \Big/ \h^n H \, $  for all  $ \, n \in \N \, $:  \, thus  $ \, \h^n
\widehat{H} \bigcap H = \h^n H \, $,  \, and similarly  $ \, \h^n
\widehat{H}^{\widehat{\otimes} \ell} \bigcap H^{\otimes \ell} = \h^n
H^{\otimes \ell} $,  \, for all  $ \, n $,  $ \ell \in \N \, $,  \,
whence  $ \, \widehat{H}' \bigcap H = H' \, $  follows at once; this
easily implies  $ \, I_{\widehat{H}'} \bigcap H = I_{H'} \, $  as well.
By construction  $ H $  is dense inside  $ \widehat{H} $  w.r.t.~the
$ \h $--adic  topology; then  $ H' $  is dense inside  $ \widehat{H}' $
w.r.t.~the topology induced on the latter by the  $ \h $--adic  topology
of  $ \widehat{H} $.  Now, the description of  $ \widehat{H}' $  in
[Ga4], \S 3.5 (which can be given also when  $ \, \dim(\gerg) =
+\infty \, $),  tells us that  $ \, \widehat{H}' \bigcap \h^n H =
I_{\widehat{H}'}^n \, $;  then we argue that  $ H' $  is dense within
$ \widehat{H}' $  w.r.t.~the  $ I_{\widehat{H}'} $--adic  topology
of  $ \widehat{H}' $.  This together with  $ \, I_{\widehat{H}'}
\bigcap H = I_{H'} \, $  implies, by a standard argument, that
$ \, I{\phantom{|}\!}^n_{\!\!\widehat{H}'} \bigcap H =
I{\phantom{|}\!}^n_{\!\!H'} \, $  for all  $ \, n \in \N \, $.
                                        \par
   Finally, take  $ \, \eta \in H \setminus \h \, H \, $.
We can show that there exists an  $ \, \eta' \in I^{\,
\partial(\overline{\eta})}_{\widehat{H}'} \, $
(notation of Lemma 4.2{\it (d)\/})  such that  $ \,
\eta' = \h^{\partial(\overline{\eta})} \eta + \eta'_+
\, $  for some  $ \, \eta'_+ \in I^{\,\partial(\overline{\eta})
+ 1}_{\widehat{H}'} \, $,  just proceeding like in [Ga4], \S 3.5
(noting again that we can drop the assumption  $ \, \dim(\gerg)
< +\infty \, $):  roughly, we consider any basis of  $ \,
H{\big\vert}_{\h=0} = \widehat{H}{\big\vert}_{\h=0} \, $
containing  $ \overline{\eta} $,  we look at the dual basis
inside  $ \, \widehat{K}{\big\vert}_{\h=0} \, $  and lift it
to a topological basis of  $ K $,  then rescale the latter
--- dividing out each element by the proper power of  $ \h \, $
---   to sort a topological basis of  $ K^\vee $:  the dual basis
of  $ \widehat{H}' $  will contain an element  $ \eta' $  as
required.  But then  $ \, \h^{\partial(\overline{\eta})} \eta
= \eta' - \eta'_+ \in I^{\,\partial(\overline{\eta})}_{\widehat{H}'}
\bigcap H = I^{\,\partial(\overline{\eta})} \, $,  \, thanks to the
previous analysis: therefore  $ \, \eta = \h^{-\partial(\overline{\eta})}
\cdot \h^{\partial(\overline{\eta})} \eta \in \h^{-\partial
(\overline{\eta})} I^{\,\partial(\overline{\eta})} \subseteq
{\big( H' \big)}^{\!\vee} \, $.  The outcome is  $ \, H \subseteq
{\big( H' \big)}^{\!\vee} \, $,  \, whilst the reverse inclusion
follows from Theorem 3.6.   \qed
\enddemo

\vskip7pt

\proclaim {Corollary 4.6} \, Let  $ \, \text{\it Char}\,
(\Bbbk) = 0 \, $.  Let  $ \, U_\h \in \QrUEA \, $.  Then
$ \, {\big( {U_\h}' \big)}_F = {(U_\h)}_F \, $.
\endproclaim

\demo{Proof}  Definitions give  $ \, {H^\vee}_{\!\!F} = H_F \, $  for
all  $ \, H \in \HA \, $.  Therefore, since  $ \, U_\h = {\big( {U_\h}'
\big)}^\vee \, $  by Proposition 4.5, we have  $ \, {\big( {U_\h}'
\big)}_F = {\Big( \! {\big( {U_\h}' \big)}^{\!\vee} \Big)}_{\!F}
= {(U_\h)}_F \, $,  \, q.e.d.   \qed
\enddemo

\vskip7pt

   {\sl $ \underline{\hbox{\it Remark}} $:}  \; it is worth noticing
that, while  $ \, {H^\vee}_{\!\!F} = H_F \, $  for all  $ \, H \in \HA
\, $,  \,  {\sl we have not in general}  $ \, {H'}_{\!F} = H_F \, $;
\, in particular, example exist of non-trivial  $ \, H \in \HA \, $
such that  $ \, H' = R \cdot 1_{\scriptscriptstyle H} \, $,  \, so that
$ \, {H'}_F = F(R) \cdot 1_{\scriptscriptstyle H} \subsetneqq H \, $.
These cases also yield counterexamples to Proposition 4.5, namely
some  $ \, H \in \HA $  for which  $ \, {\big( H' \, \big)}^{\!\vee}
\subsetneqq H \, $.

\vskip7pt

\proclaim{Theorem 4.7} \, Let  $ \, F_\h[G] \in \QFA \, $
(notation of Remark 1.5) such that  $ F_\h[G]{\big|}_{\h=0} $
is reduced.  Then  $ \, {F_\h[G]}^\vee{\big|}_{\h=0} \, $  is
a universal enveloping algebra, namely
  $$  {F_\h[G]}^\vee{\Big|}_{\h=0} := {F_\h[G]}^\vee \Big/
\h \, {F_\h[G]}^\vee = U(\gerg^\times)  $$
where  $ \, \gerg^\times $  is the cotangent Lie bialgebra
of  $ \, G \, $  (cf.~\S 1.1).
\endproclaim

\demo{Proof}  Set for simplicity  $ \, F_\h := F_\h[G] \, $,  $ \,
F_0 := F_\h \big/ \h \, F_\h = F[G] \, $,  \, and  $ \, {F_\h}^{\!
\vee} := {F_\h[G]}^\vee $,  $ \, {F_0}^{\!\vee} := {F_\h}^{\!\vee}
\big/ \h \, {F_\h}^{\!\vee} $.  By Theorem 3.7,  $ \, {F_\h}^{\!\vee}
\, $  is a QrUEA, so  $ \, {F_0}^{\!\vee} = {F_\h}^{\!\vee}{\big|}_{\h=0}
\, $  is the restricted universal enveloping algebra  $ \u(\gerk) $  of
some restricted Lie bialgebra  $ \gerk $.  Our purpose is to prove that:
first,  $ \, {F_0}^{\!\vee} = \u(\gerk) = U(\gerh) \, $  for some Lie
bialgebra  $ \gerh $;  second,  $ \, \gerh \cong \gerg^\times $.
                                             \par
   Once again we can reduce to the case when  $ F_\h $  is
separated w.r.t.~the  $ \h $--adic  topology.  Indeed, we
have  $ \, {\big( F_\h \big)}_\infty = {\big( {F_\h}^{\!\vee}
\big)}_\infty \, $  by  Lemma 4.2{\it (c)\/};  then  $ \,
\overline{{F_\h}^{\!\vee}} := {F_\h}^{\!\vee} \Big/ {\big(
{F_\h}^{\!\vee} \big)}_\infty = {F_\h}^{\!\vee} \big/ {\big(
F_\h \big)}_\infty = {\big( \overline{F_\h} \, \big)}^{\!\vee}
\, $  by  Lemma 3.2{\it (a)\/}  (taking notation from
there), and so  $ \, {F_\h}^{\!\vee}{\Big|}_{\h=0} =
\overline{{F_\h}^{\!\vee}}{\Big|}_{\h=0} = {\big( \overline{F_\h}
\, \big)}^{\!\vee}{\Big|}_{\h=0} \, $,  \, where the first identity
follows from  Lemma 3.2{\it (b)}.  Therefore it is enough to prove
the claim for  $ \, \overline{F_\h} \, $,  \, which means that we
can assume that  $ \, {\big( F_\h \big)}_\infty = \{0\} \, $.
                                             \par
   Like in Lemma 4.1, let  $ \, I := I_{\scriptscriptstyle F_\h}
\, $,  \, and let  $ \widehat{F_\h} $  be the  $ I $--adic
completion of  $ F_\h \, $.  By assumption  $ \, I^\infty =
{\big( F_\h \big)}_\infty = \{0\} \, $,  \, hence the natural map
$ \, F_\h \loongrightarrow \widehat{F}_\h \, $  is a monomorphism.
                                             \par
   Consider  $ \, J := \text{\sl Ker}\, \big(\epsilon \, \colon \,
F_\h \longrightarrow R \,\big) \, $,  \, and let  $ \, J^\vee :=
\h^{-1} J \subset {F_\h}^{\!\vee} \, $.  As in the proof of Lemma 4.1,
let  $ \, {\{y_b\}}_{b \in \Cal{S}} \, $  be a  $ \Bbbk $--basis
of  $ \, J_0 \Big/ {J_0}^{\!2} = Q\big(F[G]\big) \, $,  \, where
$ \, J_0 := \text{\sl Ker}\,(\epsilon_{\scriptscriptstyle F[G]})
= \germ \, $,  \,  and pull it back to a subset  $ \,
{\{j_b\}}_{b \in \Cal{S}} \, $  of  $ J \, $.
%
%
%
%
Using notation of Lemma 4.1, we have  $ \, I^n \big/ I^{n+1} \cong
\widehat{I^n} \big/ \widehat{I^{n+1}} \, $  for all  $ \, n \in \N $;
\, then from the description of the various  $ \widehat{I^\ell} $
($ \ell \in \N $)  given there we see that  $ \, I^n \big/ I^{n+1}
\, $  is a  $ \Bbbk $--vector  space with basis the set of (cosets
of) ordered monomials  $ \, \big\{\, \h^{e_0} j^{\,\underline{e}}
\mod I^{n+1} \;\big|\; e_0 \in \N, \underline{e} \in
\N^{\Cal{S}}_f \, , \; e_0 + |\underline{e}| = n \,\big\}
\, $  where  $ \, |\underline{e}| := \sum_{b \in \Cal{S}}
\underline{e}(b) \, $.  As a consequence, and noting that
$ \, \h^{-n} I^{n+1} = \h \cdot \h^{-(n+1)} I^{n+1} \equiv
0 \mod \h \, {F_\h}^{\!\vee} \, $,  \, we argue that  $ \,
\h^{-n} I^n \mod \h \, {F_\h}^{\!\vee} \, $  is spanned
over  $ \Bbbk $  by  $ \, \big\{\, \h^{-|\underline{e}|}
j^{\,\underline{e}} \mod \h \, {F_\h}^{\!\vee} \;\big|\; \underline{e}
\in \N^{\Cal{S}}_f \, , \; |\underline{e}| \leq n \,\big\} \, $:  \,
we claim this set is in fact a basis of  $ \, \h^{-n} I^n \mod \h \,
{F_\h}^{\!\vee} \, $.  Indeed, if not we find a non-trivial linear
combination of the elements of this set which is zero: multiplying
by  $ \h^n $  this gives an element  $ \, \gamma_n \in I^n \setminus
I^{n+1} \, $  such that  $ \, \h^{-n} \gamma_n \equiv 0 \mod \h \,
{F_\h}^{\!\vee} \, $;  \, then there is  $ \, \ell \in \N \, $  such
that  $ \, \h^{-n} \gamma_n \in \h \cdot \h^{-\ell} I^\ell \, $,
\, so  $ \, \h^\ell \gamma_n = \h^{1+n} I^\ell \subseteq I^{1 + n +
\ell} \, $:  \, but then  Lemma 4.2{\it (a)}  yields  $ \, \gamma_n
\in I^{n+1} $,  \, a contradiction!   \hbox{The outcome is that
$ \, \big\{\, \h^{-|\underline{e}|} j^{\,\underline{e}} \mod \h
\, {F_\h}^{\!\vee} \;\big|\; \underline{e} \in \N^{\Cal{S}}_f
\,\big\} \, $  is a $ \Bbbk $--basis  of  $ {F_0}^{\!\vee} $.}
                                             \par
   Now let  $ \, j_\beta^\vee := \h^{-1} j_\beta \, $  for all
$ \, \beta \in \Cal{S} \, $.  Since  $ \, j_\mu \, j_\nu - j_\nu
\, j_\mu \in \h \, J \, $,  \, for any  $ \, \mu, \nu \in \Cal{S}
\, $,  \, we can write  $ \; j_\mu \, j_\nu - j_\nu \, j_\mu = \h \sum_{\beta
\in \Cal{S}} c_\beta \, j_\beta + \h^2 \gamma_1 +
\h \, \gamma_2 \; $  for some  $ \, c_\beta \in R \, $,  $ \gamma_1
\in J \, $  and  $ \, \gamma_2 \in J^2 $,  \, whence  $ \; \big[
j_\mu^\vee, j_\nu^\vee \big] := j_\mu^\vee \, j_\nu^\vee -
j_\nu^\vee \, \gamma_\mu^\vee = \sum_{\beta \in \Cal{S}} c_\beta
\, j_\beta^\vee + \gamma_1 + \h^{-1} \gamma_2 \equiv \sum_{\beta
\in \Cal{S}} c_\beta \, j_\beta^\vee \; $  mod  $ \, J + J^\vee
J \, $:  \, but  $ \, J + J^\vee J = \h \, \big( J^\vee + J^\vee
J^\vee \big) \subseteq \h \, {F_\h}^{\!\vee} \, $,  \, so  $ \;
\big[ j_\mu^\vee, j_\nu^\vee \big] \equiv \sum_{\beta \in \Cal{S}}
c_\beta \, j_\beta^\vee \mod \h \, {F_\h}^{\!\vee} \, $  which shows
that  $ \, \gerh := J^\vee \mod \h \, {F_\h}^{\!\vee} \, $  is a Lie
subalgebra of  $ {F_0}^{\!\vee} $.  By the previous analysis
$ {F_0}^{\!\vee} $  has the  $ \Bbbk $--basis  $ \,
\big\{ {\big( j^\vee \big)}^{\underline{e}} \mod \h \,
{F_\h}^{\!\vee} \;\big|\; \underline{e} \in \N^{\Cal{S}}_f
\,\big\} \, $,  \, hence the Poincar\'e-Birkhoff-Witt theorem
tells us that  $ \, {F_0}^{\!\vee} = U(\gerh) \, $  as associative
algebras.  On the other hand, we saw in the proof of Theorem 3.7
that  $ \; \Delta \big( j^\vee \big) \equiv j^\vee \otimes 1 + 1
\otimes j^\vee \mod \, \h \, {\big( {F_\h}^{\!\vee} \big)}^{\otimes 2}
\; $  for all  $ \, j^\vee \in J^\vee $,  \, which gives  $ \;
\Delta(\text{j}) = \text{j} \otimes 1 + 1 \otimes \text{j} \; $
for all  $ \, \text{j} \in \gerh \, $,  \, whence  $ \, {F_0}^{\!\vee}
= U(\gerh) \, $  as  {\sl Hopf\/}  algebras too.
                                             \par
   Now for the second step.  The specialization map  $ \; \pi^\vee
\colon \, {F_\h}^{\!\vee} \relbar\joinrel\twoheadrightarrow
{F_0}^{\!\vee} = U(\gerh) \; $  restricts to  $ \; \eta \,
\colon \, J^\vee \relbar\joinrel\twoheadrightarrow \;
\gerh \, := J^\vee \!\!\! \mod \h \, {F_\h}^{\!\vee} = J^\vee \Big/
J^\vee \cap {\big( \h \, {F_\h}^{\!\vee} \big)} = J^\vee \Big/ \big(
J + J^\vee J_\h \big) \, $,  \, for  $ \, J^\vee \cap {\big( \h \,
{F_\h}^{\!\vee} \big)} = J^\vee \cap \h^{-1} {I_{\scriptscriptstyle
F_\h}}^{\hskip-3pt 2} = J_\h + J^\vee J_\h \, $  by  Lemma 4.2{\it
(b)\/}.  Moreover, multiplication by  $ \h^{-1} $  yields an
$ R $--module  isomorphism  $ \, \mu \, \colon \, J
\, {\buildrel \cong \over
{\lhook\joinrel\relbar\joinrel\relbar\joinrel\twoheadrightarrow}}
\, J^\vee $.  Let  $ \; \rho \, \colon \, J_0
\relbar\joinrel\twoheadrightarrow J_0 \big/ {J_0}^{\!2} =:
\gerg^\times \, $  be the natural projection map, and  $ \; \nu
\, \colon \, \gerg^\times \lhook\joinrel\longrightarrow J_0 \, $
a section of  $ \rho \, $.  The specialization map  $ \; \pi \,
\colon \, F_\h \relbar\joinrel\twoheadrightarrow
F_0 \; $  restricts to  $ \; \pi' \colon \, J
\relbar\joinrel\twoheadrightarrow J \big/ (J \cap \h \, F_\h) = J_\h
\big/ \h \, J_\h = J_0 \, $:  \, we fix a section  $ \, \gamma
\, \colon \, J_0 \lhook\joinrel\relbar\joinrel\rightarrow J_\h
\, $  of  $ \pi' $.
                                             \par
   Consider the composition map  $ \, \sigma := \eta \circ \mu \circ
\gamma \circ \nu \, \colon \, \gerg^\times \! \loongrightarrow \gerh
\, $.  This is well-defined, i.e.~it is independent of the choice
of  $ \nu $  and  $ \gamma \, $.  Indeed, if  $ \, \nu, \nu' \colon
\, \gerg^\times \lhook\joinrel\relbar\joinrel\rightarrow J_0 \, $
are two sections of  $ \rho $,  and  $ \sigma $,  $ \sigma' $  are
defined correspondingly (with the same fixed  $ \gamma $  for both),
then  $ \, \text{\it Im}\,(\nu-\nu') \subseteq \hbox{\sl Ker}\,
(\rho) = {J_0}^{\!2} \subseteq \hbox{\sl Ker}\, (\eta \circ \mu
\circ \gamma) \, $,  \, so that  $ \, \sigma = \eta \circ \mu
\circ \gamma \circ \nu = \eta \circ \mu \circ \gamma \circ \nu'
= \sigma' \, $.  Similarly, if  $ \, \gamma, \gamma' \colon \,
J_0 \lhook\joinrel\relbar\joinrel\rightarrow J_\h \, $  are two
sections of  $ \pi' $,  and  $ \sigma $,  $ \sigma' $  are defined
correspondingly (with the same  $ \nu $  for both), we have  $ \,
\text{\it Im}\,(\gamma-\gamma') \subseteq \hbox{\sl Ker}\,(\pi')
= \h \, J = \hbox{\sl Ker}\,(\eta \circ \mu) \, $,  thus  $ \,
\sigma = \eta \circ \mu \circ \gamma \circ \nu = \eta \circ \mu
\circ \gamma' \circ \nu = \sigma' \, $,  \, q.e.d.  In a nutshell,
$ \sigma $  is the composition map
  $$  \gerg^\times \;{\buildrel \bar{\nu} \over
{\lhook\joinrel\relbar\joinrel\relbar\joinrel\relbar\joinrel\twoheadrightarrow}}
\;
J_0 \Big/ {J_0}^{\! 2} \;{\buildrel \bar{\gamma} \over
{\lhook\joinrel\relbar\joinrel\relbar\joinrel\relbar\joinrel\rightarrow}}\;
J \Big/ \! \big( J^2 + \h \, J \,\big) \;{\buildrel \bar{\mu} \over
{\lhook\joinrel\relbar\joinrel\relbar\joinrel\relbar\joinrel\twoheadrightarrow}}
\;
J^\vee \Big/ \big( J + J^\vee J \,\big)
\;{\buildrel \bar{\eta} \over
{\lhook\joinrel\relbar\joinrel\relbar\joinrel\relbar\joinrel\twoheadrightarrow}}
\;
\gerh  $$
where the maps  $ \bar{\nu} $,  $ \bar{\gamma} $,  $ \bar{\mu} $,
$ \bar{\eta} $,  are the ones canonically induced by  $ \nu $,
$ \gamma $,  $ \mu $,  $ \eta $,  and  $ \bar{\nu} $,  {\sl
resp.~$ \bar{\gamma} $,  does not depend on the choice of
$ \nu $,  resp.~$ \gamma $},  as it is the inverse of
the isomorphism  $ \, \bar{\rho} \, \colon \, J_0
\big/ {J_0}^{\! 2} \,{\buildrel \cong \over
{\lhook\joinrel\relbar\joinrel\twoheadrightarrow}}\,
\gerg^\times \, $,  resp.~$ \, \overline{\pi'} \, \colon
\, J \Big/ \big(J^2 + \h \, J \big) \,{\buildrel \cong \over
{\lhook\joinrel\relbar\joinrel\twoheadrightarrow}}\,
J_0 \big/ {J_0}^{\!2} \, $,  induced by  $ \, \rho \, $,
resp.~by  $ \pi' $.  We use this remark to show that  $ \sigma $
is also an isomorphism of the Lie bialgebra structure.
                                             \par
   Using the vector space isomorphism  $ \, \sigma \, \colon \,
\gerg^\times \,{\buildrel \cong \over \longrightarrow}\, \gerh \, $
we pull-back the Lie bialgebra structure of  $ \gerh $  onto
$ \gerg^\times $,  and denote it by  $ \big( \gerg^\times, {[\,\ ,
\ ]}_\bullet, \delta_\bullet \big) $;  \, on the other hand,
$ \gerg^\times $  also carries its natural structure of Lie
bialgebra, dual to that of  $ \gerg $  (e.g., the Lie bracket
is induced by restriction of the Poisson bracket of  $ F[G] \, $),
denoted by  $ \big( \gerg^\times, {[\,\ ,\ ]}_\times, \delta_\times
\big) $:  we must prove that these two structures coincide.
                                             \par
   First,  {\it for all  $ \, x_1 $,  $ x_2 \in \gerg^\times \, $
we have}  $ \; {\big[x_1,x_2\big]}_\bullet = {\big[ x_1, x_2
\big]}_\times \, $.
                                             \par
   Indeed, let  $ \, f_i := \nu(x_i) \, $,  $ \, \varphi_i :=
\gamma(f_i) \, $,  $ \, \varphi^\vee_i := \mu(\varphi_i) \, $,
$ \, y_i := \eta \big( \varphi^\vee_i \big) \, $  \, ($ \, i =
1, 2 $).  Then
  $$  \displaylines{
   {\big[ x_1, x_2 \big]}_\bullet := \sigma^{-1} \Big( {\big[
\sigma(x_1), \sigma(x_2) \big]}_\gerh \Big) = \sigma^{-1} \big(
[y_1,y_2] \big)  = \big( \rho \circ \pi' \circ \mu^{-1} \big) \Big(
\big[ \varphi^\vee_1, \varphi^\vee_2 \big] \Big) =   \hfill   \cr
    \hfill   = \big( \rho \circ \pi' \big) \big( \h^{-1} [\varphi_1,
\varphi_2] \big) = \rho \big( \{f_1,f_2\} \big) =: {\big[ x_1, x_2
\big]}_\times \; ,  \qquad  \text{q.e.d.}  \cr }  $$
   \indent   The case of Lie cobrackets can be treated similarly;
but since they take values in tensor squares, we make use of
suitable maps  $ \, \nu_\otimes := \nu^{\otimes 2} $,  $ \,
\gamma_\otimes := \gamma^{\otimes 2} $,  etc.; we also make
use of notation  $ \, \chi_\otimes \! := \eta_\otimes \! \circ
\mu_\otimes = {(\eta \circ \mu)}^{\otimes 2} \, $  and  $ \,
\nabla \! := \Delta - \Delta^{\text{op}} $.
                                                 \par
   Now,  {\it for all  $ \, x \in \gerg^\times \, $  we have}
$ \; \delta_\bullet(x) = \delta_\times(x) \, $.
                                                 \par
   Indeed, let  $ \, f := \nu(x) \, $,  $ \, \varphi := \gamma(f)
\, $,  $ \, \varphi^\vee := \mu(\varphi) \, $,  $ \, y := \eta
\big( \varphi^\vee \big) \, $.  Then we have
  $$  \displaylines{
   \delta_\bullet(x) := {\sigma_\otimes}^{\!-1} \big( \delta_\gerh
(\sigma(x)) \big) = {\sigma_\otimes}^{\!\!-1} \big( \delta_\gerh
\big( \eta \big( \varphi^\vee \big) \big) \big) =
{\sigma_\otimes}^{\!\!-1} \Big( \eta_\otimes \big( \h^{-1}
\nabla \big( \varphi^\vee \big) \big) \Big) =   \hfill  \cr
   = {\sigma_\otimes}^{\!\!-1} \Big( \! {\big( \eta \circ \mu
\big)}_\otimes \big( \nabla(\varphi) \! \big) \! \Big) =
{\sigma_\otimes}^{\!\!-1} \Big( \! {\big( \eta \circ \mu \circ
\gamma \big)}_\otimes \big( \nabla(f) \! \big) \! \Big) =
\rho_\otimes \big( \nabla(f) \! \big) = \rho_\otimes \big(
\nabla(\nu(x)) \! \big) = \delta_\times(x)  \cr }  $$
where the last equality holds because  $ \delta_\times(x) $  is
uniquely defined as the unique element in  $ \, \gerg^\times
\otimes \gerg^\times \, $  such that  $ \; \big\langle u_1
\otimes u_2 \, , \, \delta_\times(x) \big\rangle = \big\langle
[u_1,u_2] \, , \, x \big\rangle \; $  for all  $ \, u_1, u_2
\! \in \! \gerg \, $,  \, and we have
  $$  \big\langle [u_1,u_2] \, , \, x \big\rangle = \big\langle
[u_1,u_2] \, , \, \rho(f) \big\rangle = \big\langle u_1 \otimes u_2
\, , \nabla(f) \big\rangle = \big\langle u_1 \otimes u_2 \, ,
\rho_\otimes \big(\nabla(\nu(x))\big) \big\rangle \; .
\quad  \square  $$
\enddemo

\vskip7pt

{\bf 4.8 Interlude: quantizations of pointed Poisson manifolds and
of their linear approximation.} \, The proof of Theorem 4.7 in
fact leads to a more general result; to mention it, we need some
more terminology.
                                             \par   
  Namely, among algebraic  $ \Bbbk $--varieties  let us
consider the  {\sl pointed Poisson varieties},  defined to be
pairs  $ (M,\bar{m}) $  where  $ M $  is a Poisson variety and 
$ \, \bar{m} \in M \, $  is a point of  $ M $  where the rank of
the Poisson bivector is zero: in other words,  $ \{\bar{m}\} $  is
a symplectic leaf of  $ M \, $.  A  {\sl morphism of pointed Poisson
varieties  $ (M, \bar{m}) $  and  $ (N,\bar{n}) $\/}  is any Poisson
map  $ \, \varphi : M \longrightarrow N \, $  such that  $ \,
\varphi(\bar{m}) = \bar{n} \, $.  Clearly this defines a subcategory
of the category of all Poisson varieties, whose morphisms are those
morphisms of Poisson varieties which map distinguished points into
distinguished points.  In terms of their function algebras, any
pointed Poisson variety  $ M $  is given by the datum  $ \big(
F[M], \germ_{\bar{m}}\big) $  where  $ \germ_{\bar{m}} $  is
the defining ideal of  $ \, \bar{m} \in M \, $.  
                                               \par    
   By assumptions, the Poisson bracket of  $ F[M] $  restricts
to a Lie bracket onto  $ \germ_{\bar{m}} \, $:  \, from this the
quotient space  $ \, \L_M := \germ_{\bar{m}} \big/ \germ_{\bar{m}}^{\,2}
\, $  (the cotangent space to  $ M $  at  $ \bar{m} $)  inherits a Lie
algebra structure too, the so-called ``linear approximation of  $ M $ 
at  $ \bar{m} \, $''  (see e.g.~[We], \S 4).  In the following we
also call it the cotangent Lie algebra of  $ (M,\bar{m}) $,  or
simply of  $ M \, $.  
                                         \par    
   Natural examples of pointed Poisson varieties are the coisotropic
Poisson homogeneous spaces, also called  {\sl Poisson quotients}, 
i.e.~those Poisson homogeneous spaces of the form  $ \, G \big/ H \, $, 
\, where  $ G $  is a Poisson group and  $ H $  is a (closed) coisotropic
subgroup, where coisotropic means that the ideal  $ I(H) $  in  $ F[G] $ 
of all functions vanishing on  $ H $  is a Poisson subalgebra of 
$ F[G] \, $.  The distinguished point is the coset  $ eH $  of the unit
element  $ \, e \in G \, $.    
                                         \par    
   Another special class is given by the category of Poisson monoids
(\,=\,unital Poisson semigroups): each one of them is naturally pointed
by its unit element.  If  $ \, (M,\bar{m}) = (\varLambda,e) \, $  is
any Poisson monoid, then  $ F[\varLambda] $  is a bialgebra (and
conversely), and  $ \L_\varLambda \, $  has a natural structure of 
{\sl Lie bialgebra}   --- the  {\sl cotangent Lie bialgebra\/}  of 
$ \varLambda \, $  ---   the Lie cobracket being induced by the
coproduct of  $ F[\varLambda] $,  hence (dually) by the multiplication
in  $ \varLambda \, $.  It follows then that  $ U(\L_\varLambda) $  is
a co-Poisson Hopf algebra.  When in particular the monoid  $ \varLambda $ 
is a Poisson group  $ G $  we have  $ \, \varLambda = \gerg^\times \, $.   
                                    \par   
   We call  {\sl quantization of a pointed Poisson variety 
$ (M,\bar{m}) $}  any  $ \, A \in \calA^+ \, $  such that  $ \,
A{\big|}_{\h=0} \! \cong F[M] \, $  as Poisson  $ \Bbbk $--algebras, 
and if  $ \, \pi : A \relbar\joinrel\twoheadrightarrow A{\big|}_{\h=0}
\! \cong F[M] \, $  is the specialisation map  ($ \, \h \mapsto 0 \, $), 
then  $ \, \text{\sl Ker}\,(\pi \circ \undepsilon_{{}_M}) = \germ_{\bar{m}}
\, $;  \, in this case we write  $ \, A = F_\h[M] \, $.  For any such
object we set  $ \, J_M := \text{\sl Ker}\,(\undepsilon_{{}_M}) \, $ 
and  $ \, I_M := J_M + \h \, A \, $.  A  {\sl morphism of quantizations
of Poisson varieties\/}  is any morphism  $ \, \phi : F_\h[M]
\longrightarrow F_\h[N] \; $  in  $ \calA^+ \, $  such that 
$ \, \phi(J_M) \subseteq J_N \, $.  Quantizations of pointed
Poisson varieties and their morphisms form a subcategory
of  $ \calA^+ \, $.  
                                    \par   
   A quick check throughout the proof of Theorem 4.7 (and of
Theorem 3.7 for the last part of the claim) then shows that the
same arguments also prove the following:   

\vskip7pt

\proclaim{Theorem 4.9} \, Let  $ \, F_\h[M] \in \calA^+ \, $  be
a quantization of a pointed Poisson manifold  $ (M,\bar{m}) $ 
(as above) such that  $ F_\h[M]{\big|}_{\h=0} $  is reduced. 
Then  $ \, {F_\h[M]}^\vee{\big|}_{\h=0} \, $  is a universal
enveloping algebra, namely
  $$  {F_\h[M]}^\vee{\Big|}_{\h=0} \, := \; {F_\h[M]}^\vee \Big/
\h \, {F_\h[M]}^\vee \, = \; U(\L_M)  $$
(notation of \S 4.8).  If in addition  $ M $  is a Poisson monoid
and  $ F_\h[M] $  is a quantization of  $ F[M] $  in  $ \calB \, $,
then the last identification above is one of Hopf algebras.   $ \square $
\endproclaim

\vskip7pt

\proclaim{Theorem 4.10}  Let  $ \, \text{\it Char}\,(\Bbbk) = 0 \, $.
Let  $ \, U_\h(\gerg) \in \QrUEA \, $  (notation of Remark 1.5).  Then
  $$  {U_\h(\gerg)}'{\Big|}_{\h=0} := {U_q(\gerg)}' \Big/
\h \, {U_\h(\gerg)}' = F\big[G^\star\big]  $$
where  $ \, G^\star $  is a connected  algebraic
Poisson group dual to  $ G $  (as in \S 1.1).
\endproclaim

\demo{Proof}  Due to Theorem 3.8,  $ \, {U_\h(\gerg)}' \, $  is a
QFA,  with  $ \; {U_\h(\gerg)}' {\buildrel \,{\h \rightarrow 0}\,
\over \llongrightarrow}\, F[H] \; $  for some connected algebraic
Poisson group  $ H \, $;  \, in addition we know by assumption that
$ F[H] $  is reduced: we have to show that  $ H $  is a group of
type  $ G^\star $  as in the claim.
                                            \par
   Applying Theorem 4.7 to the QFA  $ \, {U_\h(\gerg)}' \, $
yields  $ \; {\big( {U_\h(\gerg)}' \big)}^{\!\vee} {\buildrel
\,{\h \rightarrow 0}\, \over \llongrightarrow}\; U(\gerh^\times)
\, $.  Since Proposition 4.5 gives  $ \, U_\h(\gerg) = {\big(
{U_\h(\gerg)}' \big)}^{\!\vee} $, \,  we have then
  $$  U(\gerg) \,{\buildrel \,{0 \leftarrow \h}\, \over
{\leftarrow\joinrel\relbar\joinrel\relbar\joinrel\relbar\joinrel\relbar}}\,
U_\h(\gerg) = {\big({U_\h(\gerg)}'\big)}^{\!\vee}
{\buildrel \,{\h \rightarrow 0}\, \over
{\relbar\joinrel\llongrightarrow}}\; U(\gerh^\times)  $$
so that  $ \, \gerh^\times = \gerg \, $:  \, thus  $ \, \gerh
:= {\big( \gerh^\times \big)}^\star = \gerg^\star \, $,  \,
whence  $ \, H = G^\star \, $,  \, q.e.d.   \qed
\enddemo

\vskip7pt

\proclaim {Theorem 4.11} \, Let  $ \, \text{\it Char}\,(\Bbbk)
= 0 \, $.  Consider  $ \, F_\h \in \QFA \, $,  $ \, U_\h \in
\QrUEA \, $,  \, and a perfect Hopf pairing  $ \; \langle \,\ ,
\, \ \rangle \, \colon F_\h \times U_\h \loongrightarrow R \, $
such that  $ \, F_\h = {U_\h}^\bullet \, $  and  $ \, U_\h =
{F_\h}^\bullet $.  Then
  $$  {U_\h}' = {\big( {F_\h}^{\!\vee} \big)}^\bullet  \qquad
\text{and}  \qquad  {F_\h}^{\!\vee} = {\big( {U_\h}' \big)}^\bullet
\; .  $$
\endproclaim

\demo{Proof} First of all notice that the assumptions imply that
the specialized Hopf pairing  $ \, F_\h{\big|}_{\h=0} \times
U_\h{\big|}_{\h=0} \!\longrightarrow \Bbbk \, $  is perfect as
well: then  Proposition 4.4{\it (b)\/}  gives  $ \, {U_\h}'
= {\big( {F_\h}^{\!\vee} \big)}^\bullet \, $.  In addition  $ \,
{F_\h}^{\!\vee} \subseteq {\big( {U_\h}' \big)}^\bullet \, $  by
Proposition 4.4{\it (a)\/},  and we have to prove the reverse
inclusion.
                                              \par
   Let   $ \, \varphi \in {\big( {U_\h}' \big)}^\bullet \, $;
\, in particular, we can choose  $ \varphi $  such that  $ \,
\big\langle \varphi \, , {U_\h}' \big\rangle = R \, $.  Since  $ \,
{\big( {U_\h}' \big)}^\bullet \subseteq F(R) \otimes_R F_\h = F(R)
\otimes_R {F_\h}^{\!\vee} \, $,  \, there exists  $ \, c \in R
\setminus \{0\} \, $  such that  $ \, \varphi_+ := c \, \varphi
\in {F_\h}^{\!\vee} \setminus \h \, {F_\h}^{\!\vee} \, $:  \, it
follows that  $ \, \big\langle \varphi_+ \, , {U_\h}' \big\rangle
= c \, R \, $.  If  $ \, F_\h = F_\h[G] \, $,  $ \, U_\h = U_\h(\gerg)
\, $,  then Theorems 4.7--8 give  $ \, {F_\h}^{\!\vee}{\big|}_{\h=0}
= U(\gerg^\times) \, $  and  $ \, {U_\h}'{\big|}_{\h=0} = F[G^\star]
\, $.  Thus there is  $ \, \overline{\eta} \in F[G^\star] \, $  such
that  $ \, \big\langle \varphi_+{\big|}_{\h=0} \, , \overline{\eta}
\, \big\rangle = 1 \, $,  \, hence there is  $ \, \eta \in {U_\h}'
\, $  (a lift of  $ \overline{\eta} \, $)  such that  $ \, \big\langle
\varphi_+ , \eta \big\rangle = 1 + \h \, \kappa \, $  for some  $ \,
\kappa \in R \, $;  \, but  $ \, \big\langle \varphi_+ , \eta
\big\rangle \in c \, R \, $  by construction, hence  $ c \, $
divides  $ (1 + \h \, \kappa) $  in  $ R \, $.
                                              \par
   As  $ \, \varphi_+ \in {F_\h}^{\!\vee} := \bigcup_{n \in \N} \h^{-n}
I{\phantom{|}\!}^n_{\!\!F_\h} \, $  we have  $ \, \varphi_+ = \h^{-n}
\varphi_0 \, $  for some  $ \, n \in \N \, $  and  $ \, \varphi_0
\in I{\phantom{|}\!}^n_{\!\!F_\h} \, $;  \, therefore  $ \, \big\langle
\varphi_0 , {U_\h}' \big\rangle = c \, \h^n R \, $.  On the other hand,
since  $ \, U_\h = \big({U_\h}'\big)^\vee \, $  (by Proposition 4.5)
each  $ \, y \in U_\h \, $  can be written as  $ \, y = \h^{-\ell}
y{}' \, $  for some  $ \, \ell \in \N \, $  and  $ \, y{}' \in {U_\h}'
\, $;  then  $ \, \big\langle \varphi_0 , y \big\rangle = c \,
\h^{n-\ell} \big\langle \varphi , y{}' \big\rangle \in R \bigcap
c \, \h^{n-\ell} R \, $  because  $ \, \big\langle \varphi_0 ,
y \big\rangle \in R \, $  and  $ \, \big\langle \varphi , y{}'
\big\rangle \in R \, $.  Now, if  $ \, \h^{n-\ell} \big\langle
\varphi , y{}' \big\rangle \not\in R \, $  then  $ \, n - \ell
< 0 \, $  and so  $ \h \, $  divides  $ c \, $.  Since  $ c \, $
divides  $ (1 + \h \, \kappa) $  we get an absurd, unless  $ c \, $
is invertible in  $ R \, $:  hence  $ \, \varphi = c^{-1} \varphi_+
\in {F_\h}^{\!\vee} \, $,  \, q.e.d.  Otherwise, we have always
$ \, \h^{n-\ell} \big\langle \varphi , y{}' \big\rangle \in R \, $,
which means  $ \, \big\langle \varphi_0 , y \big\rangle \in c \, R
\, $  for all  $ \, y \in U_\h \, $;  \, thus  $ \, c^{-1} \varphi_0
\in {U_\h}^\bullet = F_\h \, $.  Now consider the  $ I_{F_\h} $--adic
completion  $ \widehat{F_\h} $  of  $ F_\h \, $:  the kernel of the
natural map  $ \, \mu : F_\h \loongrightarrow \widehat{F_\h} \, $
is  $ \, I{\phantom{|}\!}^\infty_{\!\!F_\h} = {(F_\h)}_\infty \, $
(because  $ \, F_\h \in \QFA \, $),  and the latter is zero because
it is contained in the trivial left radical of the perfect pairing
between  $ F_\h $  and  $ U_\h $;  therefore  $ F_\h $  embeds into
$ \widehat{F_\h} $  via  $ \mu \, $.  We have  $ \, c^{-1} \varphi_0
\in F_\h \subseteq \widehat{F_\h} \, $  and  $ \, \varphi_0 \in
I{\phantom{|}\!}^n_{\!\!F_\h} \subseteq \widehat{I{\phantom{|}
\!}^n_{\!\!F_\h}} \, $:  then from  Lemma 4.1{\it (a)--(b)\/}
we argue that  $ \, c^{-1} \varphi_0 \in \widehat{I{\phantom{|}
\!}^n_{\!\!F_\h}} \, $,  \, hence finally  $ \, c^{-1} \varphi_0
\in \widehat{I{\phantom{|}\!}^n_{\!\!F_\h}} \bigcap F_\h =
I{\phantom{|}\!}^n_{\!\!F_\h} \, $,  \, thanks to Lemma 4.1{\it
(d)}.  The outcome is  $ \, \varphi = c^{-1} \h^{-n} \varphi_0
\in h^{-n} I{\phantom{|}\!}^n_{\!\!F_\h} \subseteq {F_\h}^{\!\vee}
\, $,  \, q.e.d.   \qed
\enddemo
%
%
 \eject   

   At last, we can gather our partial results to prove the
main Theorem:

\vskip7pt

\demo{$ \underline{\hbox{\it Proof of Theorem 2.2}} $}  Part {\it
(a)\/}  is proved by patching together Proposition 3.3, Proposition
3.5  and Theorems 3.6--8.  Now recall that if  $ \, \Char(\Bbbk)
= 0 \, $  every commutative Hopf  $ \Bbbk $--algebra  is reduced;
then part  {\it (b)}  follows from Theorem 3.6 and Propositions 4.3
and 4.5.  Part  {\it (c)\/}  is proved by Theorems 4.7--8, whereas
Theorem 4.11 proves part  {\it (d)}.  Finally, assume  $ \, \Char(\Bbbk)
= 0 \, $  and consider  $ \, \H \in \HA_F \, $:  if  $ \, H_{(f)}
\in \QFA \, $  (w.r.t.~a fixed prime  $ \, \h \in R \, $)  is an
$ R $--integer  form  of  $ \H $,  then  $ H_{(f)}^\vee $  is an
integer form too   --- by the very definitions ---   and it is a
QrUEA (at  $ \h \, $),  by Proposition 3.3; conversely, if  $ \,
H_{(u)} \in \QrUEA \, $  (w.r.t.~a fixed prime  $ \, \h \in R \, $)
is an  $ R $--integer  form of  $ \H $,  then also  $ H_{(u)}' $
is an integer form   --- by Corollary 4.6 ---   and it is a
QFA (again at  $ \h \, $),  by Proposition 3.5; this proves
part  {\it (e)}.   \qed
\enddemo

\vskip1,3truecm

\centerline {\bf \S \; 5 \  Application to trivial deformations:
the Crystal Duality Principle }

\vskip10pt

  {\bf 5.1 Trivial deformations and GQDP.} \, In this section,
we apply the GQDP to the framework of trivial deformations of
Hopf algebras over a field.  In particular, we consider more
closely some key examples: function algebras over algebraic groups,
universal enveloping algebras of Lie algebras, and group algebras of
abstract groups.  The outcome seems to be of special interest in its
own, as a chapter of classical   --- rather than ``quantum'' ---  Hopf
algebra theory, and we propose it as a new tool for specialists in that
matter.
                                            \par
   To be short we perform our analysis for Hopf algebras only: however,
as Drinfeld's functors are defined not only for Hopf algebras but for
augmented algebras and coaugmented coalgebras too, we might do the
same study for them as well (indeed, we do it in [Ga5]).
                                      \par
   Let us now be more precise.  Let  $ \HA_\Bbbk $  be the category of
all Hopf algebras over the field  $ \Bbbk \, $.  For all  $ \, n \in
\N \, $,  \, let  $ \; J^n := {\big( \text{\sl Ker}\,(\epsilon \colon
\, H \longrightarrow \Bbbk) \big)}^n \, $  and  $ \; D_n := \text{\sl
Ker}\, \big( \delta_{n+1} \colon \, H \longrightarrow H^{\otimes n}
\big) \, $,  \; and set  $ \; \underline{J} := {\big\{ J^n \big\}}_{n
\in \N} \, $,  $ \; \underline{D} := {\big\{ D_n \big\}}_{n \in \N}
\, $.  Of course  $ \underline{J} $  is a decreasing filtration of
$ H \, $  (maybe with  $ \, \bigcap_{n \geq 0} J^n \supsetneqq \{0\}
\, $),  and  $ \underline{D} $  is an increasing filtration of  $ H
\, $  (maybe with  $ \, \bigcup_{n \geq 0} D_n \subsetneqq H \, $),
by coassociativity of the  $ \delta_n $'s.
                                      \par
   Let  $ \, R := \Bbbk[\h\,] \, $  be the polynomial ring in the
indeterminate  $ \h \, $:  then  $ R $  is a PID (= \, principal
ideal domain), and  $ \h $  is a non-zero prime in  $ R $  such
that  $ R \big/ \h \, R $  is the field  $ \Bbbk \, $.  Let  $ \,
H_\h := H[\h] = R \otimes_{\Bbbk} H \, $,  the scalar extension of 
$ H \, $:  \, this is the  {\sl trivial deformation\/}  of  $ H \, $. 
Clearly  $ H_\h $  is a torsion free Hopf algebra over  $ R $,  hence
one can apply Drinfeld's functors to it;  in this section we do it
with respect to the element  $ \h $  itself.  We shall see that the
outcome is quite neat, and can be expressed purely in terms of Hopf
algebras in  $ \HA_\Bbbk \, $:  because of the special relation between
some features of  $ H $   --- namely, the filtrations  $ \underline{J} $ 
and  $ \underline{D} \, $  ---   and some properties of Drinfeld's
functors, we call this result ``Crystal Duality Principle'', in
that it is obtained through sort of a ``crystallization'' process.
Here we bear in mind, in a sense, Kashiwara's motivation for the
terminology ``crystal bases'' in the context of quantum groups:
see [CP], \S 14.1, and references therein.  Indeed, this theorem
can also be proved almost entirely using only classical Hopf
algebraic methods within  $ \HA_\Bbbk $,  i.e.~without resorting
to deformations: this is accomplished in [Ga5].
                                            \par
  Note that the same analysis and results (with only a bit more
annoying details to take care of) still hold if we take as  $ R $
any domain which is also a  $ \Bbbk $--algebra  and as  $ \h $  any
element in  $ \, R \setminus \{0\} \, $  such that  $ \, R \big/ \h
\, R = \Bbbk \, $;  \, for instance, one can take  $ \, R = \Bbbk[[h]]
\, $  and  $ \, \h := h \, $,  or  $ \, R = \Bbbk \big[ q, q^{-1} \big]
\, $  and  $ \h := q-1 \, $.
%
%
%
%
   Finally, in the sequel to be short we perform our analysis for
Hopf algebras only: however, as Drinfeld's functors are defined not
only for Hopf algebras but for augmented algebras and coaugmented
coalgebras too, we might do the same study for them as well.  In
particular, the final result (the Crystal Duality Principle) has
a stronger version which concerns these more general objects too
(see [Ga5]).
                                            \par
   We first discuss the general situation (\S\S 5.2--4), second we
look at the case of function algebras and enveloping algebras (\S\S
5.6--7), then we state and prove the theorem of Crystal Duality
Principle and eventually (\S\S 5.12--13) we dwell upon two other
interesting applications: hyperalgebras, and group algebras and
their duals.

\vskip7pt

\proclaim{Lemma 5.2}
 \vskip-19pt   
  $$  \eqalignno{
   {H_\h}^{\hskip-3pt\vee} \;  &  = \; {\textstyle \sum_{n \geq 0}} \,
R \cdot \h^{-n} J^n \, = \; R \cdot J^0 + R \cdot \h^{-1} J^1 + \cdots
+ R \cdot \h^{-n} J^n + \cdots  \;  &   (5.1)  \cr
   {H_\h}' \;  &  = \; {\textstyle \sum_{n \geq 0}} \, R \cdot \h^n
D_n \, = \; R \cdot D_0 + R \cdot \h \, D_1 + \cdots + R \cdot \h^n
D_n + \cdots  &   (5.2)  \cr }  $$
 \vskip-9pt   
\endproclaim

\demo{Proof} As for (5.1), we have  $ \, J_{H_\h} = R \cdot J \, $,
\, whence  $ \; {H_\h}^{\hskip-3pt\vee} := \, \sum_{n \geq 0} \h^{-n}
{J_\h}^{\!n} = \sum_{n \geq 0} \h^{-n} J^n \, $.
                                              \par
   On the other hand, one has trivially  $ \; {H_\h}' \supseteq
\sum_{n=0}^{+\infty} R \cdot \h^n D_n \, $.  Conversely, let  $ \, \eta
\in {H_\h}' \, $:  \, then  $ \, \eta = \sum_k c_k \, \eta_k \, $  for
some  $ \, c_k \in R \, $  and  $ \, \eta_k \in H \, $;  \, in addition,
we can assume the  $ \eta_k $'s  enjoy the following:  $ \, \eta_1 $,
$ \ldots $,  $ \eta_{k_1} \in D_{\ell_1} \setminus D_{\ell_1-1} \, $,
$ \, \eta_{k_1+1} $,  $ \ldots $,  $ \eta_{k_2} \in D_{\ell_2} \setminus
D_{\ell_2-1} \, $,  $ \ldots $,  $ \, \eta_1 $,  $ \ldots $,  $ \eta_{k_t}
\in D_{\ell_t} \setminus D_{\ell_t-1} \; $  for some  $ \, k_i $,
$ \ell_j $,  $ t \in \N \, $  with  $ \, k_1 < k_2 < \cdots < k_t \, $,
\, they are linearly independent over  $ \Bbbk $,  and moreover  $ \,
\eta_{k_i+1} $,  $ \dots $,  $ \eta_{k_{i+1}} \, $  belong to a vector
subspace  $ W_i $  of  $ H $  such that  $ \, W_i \bigcap D_{\ell_i} \!
= \{0\} \, $,  \, for all  $ i \, $.  By the assumptions on  $ R \, $,
for any  $ \, k \, $  there is a unique  $ \, v_k \in \N \, $  such
that  $ \, c_k \in \h^{v_k} R \setminus \h^{v_k+1} R \, $;  \, then
for all  $ \, n \in \N \, $  we have  $ \, c_k \equiv 0 \mod \h^n R
\, $  when  $ \, v_k \geq n \, $  and  $ \, c_k \equiv \sum_{s =
v_k}^{n-1} a^{(k)}_s \, \h^s \mod \h^n R \, $,  \, for some  $ \,
a^{(k)}_s \in \Bbbk \, $  with  $ \, a^{(k)}_{v_k} \not= 0 \, $,
\, when  $ \, v_k < n \, $.  Then  $ \, \sum_{v_k < n} \sum_{s =
v_k}^{n-1} a^{(k)}_s \h^s \, \delta_n(\eta_k) \equiv \delta_n(\eta)
\equiv 0 \mod \h^n \, $  and  $ \, \delta_n(\eta_k) \in H^{\otimes n}
\subset {H_\h}^{\!\otimes n} \setminus \h \, {H_\h}^{\!\otimes n} \, $
imply   --- since  $ \, {H_\h}^{\!\otimes n} \Big/ \h^n {H_\h}^{\!
\otimes n} \cong \Big( R \big/ (\h^n R) \Big) \otimes_{\Bbbk}
H^{\otimes n} \cong \Big( \Bbbk[x] \big/ (x^n) \Big) \otimes_{\Bbbk}
H^{\otimes n} \, $  ---   that  $ \; \sum_{n > v_k = v_-} a^{(k)}_{v_-}
\delta_n(\eta_k) = 0 \, $,  \; where  $ \, v_- := \underset{k}\to{\min}
\,\{v_k\} \, $,  \, hence  $ \, \sum_{n > v_k = v_-} a^{(k)}_{v_-}
\eta_k \in \hbox{\sl Ker}\,(\delta_n) =: D_{n-1} \, $:  \, since all
coefficients  $ a^{(k)}_{v_-} $  in this sum are non-zero, by our
assumptions on the  $ \eta_k $'s  this forces  $ \, \eta_k \in
\hbox{\sl Ker}\,(\delta_n) =: D_{n-1} \, $,  for all  $ k \, $
such that  $ \, v_k = v_- \, $.  The outcome is:  $ \, v_k < n
\, \Longrightarrow \, \eta_k \in D_{n-1} \, $  (for all  $ k $,
$ n \, $),  \, whence we get  $ \, \eta_k \in D_{v_k} \, $  for
all  $ k \, $,  so that  $ \; \eta = \sum_k c_k \, \eta_k \in
\sum_{s=0}^{+\infty} R \cdot \h^s D_s \, $,  \; q.e.d.   \qed
\enddemo

\vskip7pt

  {\bf 5.3 Rees Hopf algebras and their specializations.} \,
We need some more terminology.  Let  $ M $  be a module over
a commutative unitary ring  $ R $,  and let
 \eject   
  $$  \underline{M} \; := \; {\{M_z\}}_{z \in \Z} \; = \; \Big(
\cdots \subseteq M_{-m} \subseteq \cdots \subseteq M_{-1} \subseteq M_0
\subseteq M_1 \subseteq \cdots \subseteq M_n \subseteq \cdots \Big)  $$
be a bi-infinite filtration of  $ M $  by submodules  $ M_z $
($ z \in \Z $).  In particular, we consider increasing filtrations
(i.e., those with  $ \, M_z = \{0\} \, $  for all  $ \, z < 0 \, $)
and decreasing filtrations (those with  $ \, M_z = \{0\} \, $  for
all  $ \, z > 0 \, $)  as special cases of bi-infinite filtrations.
First we define the associated  {\sl blowing module\/}  to be the
$ R $--submodule  $ \Cal{B}_{\underline{M}}(M) $  of  $ M \big[ t,
t^{-1} \big] $  (where  $ t $  is any indeterminate) given by  $ \,
\Cal{B}_{\underline{M}}(M) := \sum_{z \in \Z} t^z M_z \, $;  \,
this is isomorphic to the
         {\sl first graded module\footnote{Hereafter, I pick such
terminology from Serge Lang's textbook  {\it ``Algebra''}.}
associated to  $ M $},  namely  $ \, \bigoplus_{z \in \Z} M_z \, $.
Second, we define the associated  {\sl Rees module\/}  to be the
$ R[t] $--submodule  $ \Cal{R}^t_{\underline{M}}(M) $  of  $ M \big[
t, t^{-1} \big] $  generated by  $ \Cal{B}_{\underline{M}}(M) $;  \,
straight\-forward computations then give  $ R $--module  isomorphisms
  $$  \Cal{R}^t_{\underline{M}}(M) \Big/ (t-1) \,
\Cal{R}^t_{\underline{M}}(M) \; \cong \; {\textstyle \bigcup\limits_{z
\in \Z}} M_z \; ,  \qquad  \Cal{R}^t_{\underline{M}}(M) \Big/ t \,
\Cal{R}^t_{\underline{M}}(M) \; \cong \; G_{\underline{M}}(M)  $$
where  $ \, G_{\underline{M}}(M) := \bigoplus_{z \in Z} M_z \big/
M_{z-1} \, $  is the  {\sl second graded module associated to
$ M \, $}.  In other words,  $ \Cal{R}^t_{\underline{M}}(M) $
is an  $ R[t] $--module  which specializes to  $ \, \bigcup_{z
\in \Z} M_z \, $  for  $ \, t = 1 \, $  and specializes to  $ \,
G_{\underline{M}}(M) \, $  for  $ \, t = 0 \, $;  \, therefore
the  $ R $--modules  $ \, \bigcup_{z \in \Z} M_z \, $  and  $ \,
G_{\underline{M}}(M) \, $  can be seen as 1-parameter (polynomial)
deformations of each other via the 1-parameter family of  $ R $--modules
given by  $ \Cal{R}^t_{\underline{M}}(M) $.
                                        \par
   We can repeat this construction within the category of algebras,
coalgebras, bialgebras or Hopf algebras over  $ R $  with a filtration
in the proper sense (by subalgebras, subcoalgebras, etc.):  then we'll
end up with corresponding objects  $ \Cal{B}_{\underline{M}}(M) $,
$ \Cal{R}^t_{\underline{M}}(M) $,  etc.{} of the like type (algebras,
coalgebras, etc.).  In particular we'll cope with Rees Hopf algebras.

\vskip7pt

  {\bf 5.4 Drinfeld's functors on  $ \boldkey{H}_\h \, $  and
filtrations on  $ \boldkey{H} \, $.} \, Lemma 5.2 sets a link
between properties of  $ {H_\h}' $,  resp.~of  $ {H_\h}^{\!
\vee} $,  and properties of the filtration  $ \underline{D}
\, $,  resp.~$ \underline{J} \, $,  \, of  $ H \, $.
                                              \par
   First, formula (5.1) together with the fact that  $ \, {H_\h}^{\!\vee}
\in \HA \, $  implies that  $ \, \underline{J} \, $  is a Hopf algebra
filtration of  $ H \, $;  \, conversely, if one proves that  $ \,
\underline{J} \, $  is a Hopf algebra filtration of  $ H \, $  (which
is straightforward)  then from (5.1) we get a one-line direct proof
that  $ \, {H_\h}^{\!\vee} \in \HA \, $.  Second, we can look at
$ \underline{J} $  as a bi-infinite filtration by reversing the index
notation and then extending it trivially on the positive indices, namely
  $$  \underline{J} \; = \; \Big( \cdots \subseteq J^n \subseteq \cdots
J^2 \subseteq J \subseteq J^0 \big( = H \big) \subseteq H \subseteq
\cdots \subseteq H \subseteq \cdots \Big) \, ;  $$
then the Rees Hopf algebra  $ \Cal{R}^\h_{\underline{J}}(H) $  is
defined (see \S 5.3).  Now (5.1) give  $ \, {H_\h}^{\!\vee} =
\Cal{R}^\h_{\underline{J}}(H) \, $,  \, so  $ \; {H_\h}^{\!\vee}
\Big/ \h \, {H_\h}^{\!\vee} \cong \Cal{R}^\h_{\underline{J}}(H) \Big/
\h \, \Cal{R}^\h_{\underline{J}}(H) \cong G_{\underline{J}}(H) \, $.
Thus  $ \, G_{\underline{J}}(H) \, $  is cocommutative because  $ \,
{H_\h}^{\!\vee} \Big/ \h \, {H_\h}^{\!\vee} \, $  is; conversely, we
get an easy proof of the cocommutativity of  $ \, {H_\h}^{\!\vee} \Big/
\h \, {H_\h}^{\!\vee} \, $  once we prove that  $ \, G_{\underline{J}}(H)
\, $  is cocommutative, which is straightforward.  Finally,  $ \,
G_{\underline{J}}(H) \, $  is generated by  $ \, Q(H) = J \big/
J^{\,2} \, $  whose elements are primitive, so  {\sl a fortiori\/}
$ \, G_{\underline{J}}(H) \, $  is generated by its primitive
elements; then the latter holds for  $ \, {H_\h}^{\!\vee} \Big/ \h
\, {H_\h}^{\!\vee} \, $  as well.  To sum up, as  $ {H_\h}^{\!\vee}
\in \QrUEA \, $  we argue that  $ \, G_{\underline{J}}(H) = \U(\gerg)
\, $  for some restricted Lie bialgebra  $ \gerg \, $;  conversely, we
can get  $ {H_\h}^{\!\vee} \in \QrUEA \, $  directly from the properties
of the filtration  $ \underline{J} \, $  of  $ H $.  Moreover, since
$ \, G_{\underline{J}}(H) = \U(\gerg) \, $  is graded,  $ \gerg $
{\sl as a restricted Lie algebra is graded too}.
                                             \par
   On the other hand, it is straightforward to see that (5.2) together
with the fact that  $ \, {H_\h}' \in \HA \, $  implies that  $ \,
\underline{D} \, $  is a Hopf algebra filtration of  $ H \, $;  \,
conversely, if one proves that  $ \, \underline{D} \, $  is a Hopf
algebra filtration of  $ H \, $  (as we did in  Lemma 3.4{\it (c)})
then from (5.2) we get an easy direct proof that  $ \, {H_\h}' \in
\HA \, $.  Second, we can look at  $ \underline{D} $  as a bi-infinite
filtration by extending it trivially on the negative indices, namely
  $$  \underline{D} \; = \; \Big( \cdots \subseteq \{0\} \subseteq \cdots
\{0\} \subseteq \big( \{0\} = \big) D_0 \subseteq D_1 \subseteq \cdots
\subseteq D_n \subseteq \cdots \Big) \quad ;  $$
then the Rees Hopf algebra  $ \Cal{R}^\h_{\underline{D}}(H) $
is defined (see \S 5.3).  Now (5.2) gives  $ \, {H_\h}' =
\Cal{R}^\h_{\underline{D}}(H) \, $;  \, but then  $ \; {H_\h}'
\Big/ \h \, {H_\h}' \cong \Cal{R}^\h_{\underline{D}}(H) \Big/ \h \,
\Cal{R}^\h_{\underline{D}}(H) \cong G_{\underline{D}}(H) \, $.
Thus  $ \, G_{\underline{D}}(H) \, $  is commutative because  $ \,
{H_\h}' \Big/ \h \, {H_\h}' \, $  is; or, conversely, we get an easy
proof of the commutativity of  $ \, {H_\h}' \Big/ \h \, {H_\h}' \, $
once we prove that  $ \, G_{\underline{D}}(H) \, $  is commutative, as
we did in  Lemma 3.4{\it (c)}.  Finally,  $ \, G_{\underline{D}}(H) \, $
is graded with  $ 1 $-dimensional  $ 0 $-component   --- by construction
---   hence it has no non-trivial idempotents; therefore the latter is
true for  $ \, {H_\h}' \Big/ \h \, {H_\h}' \, $  as well.  Note also that
$ \, {I_{\!{H_\h}'}\phantom{\big|}}^{\hskip-11pt\infty} = \{0\} \, $  by
construction (because  $ H_\h $  is free over  $ R \, $).  To sum up,
since  $ {H_\h}' \in \QFA \, $  we get that  $ \, G_{\underline{D}}(H)
= F[G] \, $  for some connected algebraic Poisson group  $ G \, $;
conversely, we can argue that  $ \, {H_\h}' \in \QFA \, $  directly
from the properties of the filtration  $ \underline{D} \, $.
                                             \par
   In addition, since  $ \, G_{\underline{D}}(H) = F[G] \, $  is
graded, when  $ \, \Char(\Bbbk) = 0 \, $  the (pro)affine variety
$ \, G_{(\text{\it cl\/})} \, $  of closed points of  $ G $  is a
          (pro)affine space\footnote{For it is a  {\sl cone}
--- since  $ H $  is graded ---   without vertex   --- since
$ G_{(\text{\it cl\/})} $,  being a group, is smooth.},
that is  $ \, G_{(\text{\it cl\/})} \cong \Bbb{A}_{\,\Bbbk}^{\times
\Cal{I}} = \Bbbk^{\Cal{I}} \, $  for some index set  $ \Cal{I} \, $,
\, and so  $ \, F[G] = \Bbbk \big[ {\{x_i\}}_{i \in \Cal{I}} \big] \, $
is a polynomial algebra.
                                             \par
   Finally, when  $ \, p := \Char(\Bbbk) > 0 \, $  the group  $ G $
has dimension 0 and height 1: indeed, we can see this as a consequence
of the last part of Theorem 3.8 via the identity  $ \; H'_\h \Big/ \h
\, H'_\h = G_{\underline{D}}(H) \, $,  \, or conversely we can prove
that part of Theorem 3.8 via this identity by observing that  $ G $
has those properties.  In fact, we must show that  $ \, \bar{\eta}^p
= 0 \, $  for each  $ \, \eta \in \widetilde{H} := G_{\underline{D}}(G)
\, $:  \, letting  $ \, \eta \in H'_\h \, $  be any lift of
$ \bar{\eta} $  in  $ H'_\h $,  we have  $ \, \eta \in D_\ell \, $
for some  $ \, \ell \in \N \, $,  hence  $ \, \delta_{\ell+1}(x) =
0 \, $.  From  $ \, \Delta^{\ell+1}(\eta) = \sum_{\Lambda \subseteq
\{1,\dots,\ell+1\}} \delta_\Lambda(\eta) \, $  (cf.~\S 2.1) and
the multiplicativity of  $ \Delta^{\ell+1} $  we have
 \eject   
  $$  \displaylines{
   \quad   \Delta^{\ell+1}(\eta^p) = {\big( \Delta^{\ell+1}(\eta)
\big)}^p = {\Big( {\textstyle \sum_{\Lambda \subseteq \{1, \dots,
\ell+1\}}} \delta_\Lambda(\eta) \Big)}^p  \; \in \;  {\textstyle
\sum_{\Lambda \subseteq \{1,\dots,\ell+1\}}} \, {\delta_\Lambda(\eta)}^p
\; +   \hfill  \cr
   +  \; {\textstyle \sum_{\Sb e_1, \dots, e_p < p \\   e_1 + \cdots e_p
= p \endSb}} \Big(\! {\textstyle {p \atop {e_1, \dots, e_p}}} \!\Big)
{\textstyle \sum\limits_{\Lambda_1, \dots,\Lambda_p \subseteq \{1,
\dots, \ell + 1\}}} {\textstyle \prod_{k=1}^p} \, {\delta_{\Lambda_k}
(\eta)}^{e_k}  \; +  \cr
   {} \hfill   +  \; {\textstyle \sum_{k=0}^{\ell}} \,
{\textstyle \sum_{\Sb \Psi \subseteq \{1,\dots,\ell+1\}  \\
|\Psi|=k  \endSb}}  j_\Psi \big( {J_{\!{}_{H'}}}^{\!\!\otimes k} \big)
\;  +  \; {\big( \hbox{ad}_{[\ ,\ ]}(D_{(\ell)})\big)}^{p-1} \big(
D_{(\ell)} \big)  \cr }  $$
(since  $ \, \delta_\Lambda(\eta) \in j_\Lambda \Big(
{J_{\!{}_{H'}}}^{\!\!\otimes |\Lambda|} \Big) \, $  for all
$ \, \Lambda \subseteq \{1,\dots,\ell+1\} \, $)  where  $ \; D_{(\ell)}
:= \sum\limits_{\sum_k s_k = \ell} \otimes_{k=1}^{\ell+1} D_{s_k} \; $
and  $ \; {\big(\hbox{ad}_{[\ ,\ ]}(D_{(\ell)})\big)}^{p-1} \big(
D_{(\ell)} \big) := \big[ \undersetbrace{p}\to{D_{(\ell)}, \big[
D_{(\ell)}, \dots, \big[ D_{(\ell)}, \big[ D_{(\ell)}, D_{(\ell)}}
\big] \big] \cdots \big] \big] \, $.  Then
  $$  \displaylines{
   \quad   \delta^{\ell+1}(\eta^p) = {(\id_{\scriptscriptstyle H}
- \epsilon)}^{\otimes (\ell+1)} \big( \Delta^{\ell+1}(\eta^p) \big)
\; \in \;  {\delta_{\ell+1}(\eta)}^p  \; + \;  {\textstyle \sum_{\Sb
e_1, \dots, e_p < p \\   e_1 + \cdots e_p = p \endSb}}  \Big(\!
{\textstyle {p \atop {e_1, \dots, e_p}}} \!\Big) \; \times  \hfill  \cr
   \hfill   \times \;  {\textstyle \sum_{\cup_{k=1}^p \Lambda_k =
\{1,\dots,\ell+1\}}} \, {\textstyle \prod_{k=1}^p} {\delta_{\Lambda_k}
(\eta)}^{e_k} \;  +  \; {(\id_{\scriptscriptstyle H} - \epsilon)}^{\otimes
(\ell+1)} \Big( \! {\big( \hbox{ad}_{[\ ,\ ]}(D_{(\ell)}) \big)}^{p-1}
\big( D_{(\ell)} \big) \Big) \; .  \cr }  $$
Now,  $ \, {\delta^{\ell+1}(\eta)}^p = 0 \, $  by construction, and
$ \, \Big(\! {\textstyle {p \atop {e_1, \dots, e_p}}} \!\Big) \, $
(with  $ \, e_1 $,  $ \dots $,  $ e_p < p \, $)  is a multiple of
$ p \, $,  hence it is zero because  $ \, p = \Char(\Bbbk) \, $;
\, therefore we end up with
  $$  \delta_{\ell+1}(\eta)  \; \in \;  {(\id_{\scriptscriptstyle H}
- \epsilon)}^{\otimes (\ell+1)} \Big( \! {\big( \hbox{ad}_{[\ ,\ ]}
(D_{(\ell)}) \big)}^{p-1} \big( D_{(\ell)} \big) \Big) \; .  $$
Now, by Lemma 3.4 we have  $ \, D_{s_i} \cdot D_{s_j} \subseteq D_{s_i
+ s_j} \, $  and  $ \, \big[ D_{s_i}, D_{s_j} \big] \subseteq D_{(s_i
+ s_j) - 1} \, $;  \, these together with Leibniz' rule imply that
$ \; {\big( \hbox{ad}_{[\ ,\ ]}(D_{(\ell)}) \big)}^{p-1} \big(
D_{(\ell)} \big) \subseteq \sum\limits_{\sum_t r_t = p\,\ell+1-p}
\hskip-5pt  \otimes_{t=1}^{\ell+1} D_{r_t} \, $;  \, moreover,
since  $ \, D_0 = \hbox{\sl Ker}\,(\delta_1) = \hbox{\sl Ker}\,
(\hbox{id} - \epsilon) \, $  we have
  $$  {(\id_{\scriptscriptstyle H} - \epsilon)}^{\otimes (\ell+1)}
\Big( \! {\big( \hbox{ad}_{[\ ,\ ]} (D_{(\ell)}) \big)}^{p-1}
\big( D_{(\ell)} \big) \Big)  \hskip3pt  \subseteq  \hskip3pt
{\textstyle \sum_{\Sb  \sum_t r_t = p\,\ell+1-p  \\  r_1, \dots,
r_{\ell+1} > 0  \endSb}}  \hskip-5pt  \otimes_{t=1}^{\ell+1}
D_{r_t} \;\; .  $$
In particular, in the last term above we have  $ \, D_{r_1} \subseteq
D_{(p-1)\ell+1-p} := \hbox{\sl Ker}\,(\delta_{(p-1)\ell+2-p}) \subseteq
\hbox{\sl Ker}\,(\delta_{(p-1)\ell}) \, $:  \, therefore, using the
coassociativity of the maps  $ \delta_n $'s,  we get
  $$  \delta_{p\,\ell}(\eta) = \big( \big( \delta_{(p-1)\ell} \otimes
\id^\ell \big) \circ \delta_{\ell+1} \big)(\eta)  \hskip3pt  \subseteq
{\textstyle \sum\limits_{\Sb  \sum_t r_t = p\,\ell-1  \\  r_1, \dots,
r_{\ell+1} > 0  \endSb}}  \hskip-7pt  \delta_{(p-1)\ell}(D_{r_1})
\otimes D_{r_2} \otimes \cdots \otimes D_{r_{\ell+1}}  \hskip3pt
=  \hskip3pt  0  $$
i.e.~$ \, \delta_{p\,\ell}(\eta) = 0 \, $.  This means  $ \, \eta \in
D_{p\,\ell-1} \, $,  \, whereas, on the other hand,  $ \, \eta^p \in
D_\ell^{\;p} \subseteq D_{p\,\ell} \, $:  \, then  $ \, \bar{\eta}^p
:= \overline{\eta^p} = \bar{0} \in D_{p\,\ell} \Big/ D_{p\,\ell - 1}
\subseteq G_{\underline{D}}(H) \, $,  \, by the definition of the
product in  $ \, G_{\underline{D}}(H) \, $.  Finally, by general
theory since  $ G $  has dimension 0 and height 1 the function
algebra  $ \, F[G] = G_{\underline{D}}(H) = H'_\h \Big/ \h \,
H'_\h \, $  is  {\sl truncated polynomial},  namely  $ \, F[G]
= \Bbbk \big[ {\{x_i\}}_{i \in \Cal{I}} \big] \Big/ \big(
{\{x_i^{\,p}\}}_{i \in \Cal{I}} \big) \, $.

\vskip7pt

  {\bf 5.5 Special fibers of  $ {{\boldkey{H}}_\h}^{\!\boldsymbol\vee} $
and  $ {\boldkey{H}_\h}^{\!\boldsymbol\prime} $  and deformations.}  \,
Given  $ \, H \in \HA_\Bbbk \, $,  \, consider  $ \, H_\h \, $:  \, our
goal is to study  $ \, {H_\h}^{\!\vee} \, $  and  $ \, {H_\h}' \, $.
                                        \par
   As for  $ {H_\h}^{\!\vee} $,  the natural map from  $ H $
to  $ \, \widehat{H} := G_{\underline{J}}(H) = {H_\h}^{\!\vee}
\Big/ \h \, {H_\h}^{\!\vee} =: {H_\h}^{\!\vee}{\Big|}_{\h=0} \, $
sends  $ \, {J\phantom{|}}^{\hskip-1pt \infty} := \bigcap_{n \geq 0}
{J\phantom{|}}^{\hskip-1pt n} \, $  to zero, by definition; also,
letting  $ \, H^\vee := H \Big/ {J\phantom{|}}^{\hskip-1pt \infty} \, $
(which is a Hopf algebra quotient of  $ H $  because  $ \underline{J} $
is a Hopf algebra filtration), we have  $ \, \widehat{H} =
\widehat{H^\vee} \, $.  Thus  $ \, {{(H^\vee)}_\h}^{\!\!\vee}
{\Big|}_{\h=0} = \widehat{H^\vee} = \widehat{H} = \U(\gerg_-) \, $
for some graded restricted Lie bialgebra  $ \gerg_- \, $.  On the
other hand,  $ \, {{(H^\vee)}_\h}^{\!\!\vee}{\Big|}_{\h=1} :=
{{(H^\vee)}_\h}^{\!\!\vee} \Big/ (\h-1) \, {{(H^\vee)}_\h}^{\!\!\vee}
= \sum_{n \geq 0} \overline{J}^{\,n} = H^\vee \, $  (see \S 5.3).  Thus
we can see  $ \, {{(H^\vee)}_\h}^{\!\!\vee} = \Cal{R}^\h_{\underline{J}}
(H^\vee) \, $  as a 1-parameter family inside  $ \HA_\Bbbk $  with
{\sl regular\/}  fibers   --- that is, they are isomorphic to each
other as  $ \Bbbk $--vector  spaces (indeed, we switch from  $ H $  to
$ H^\vee $  just in order to achieve this regularity) ---   which links
$ \widehat{H^\vee} $  and  $ H^\vee $  as (polynomial) deformations
of each other, namely
  $$  \U(\gerg_-) = \widehat{H^\vee} = {{(H^\vee)}_\h}^{\!\!\vee}
{\Big|}_{\h=0}  \hskip5pt  \underset{\;{{(H^\vee)}_\h}^{\!\!\vee}}
\to  {\overset{0 \,\leftarrow\, \h \,\rightarrow\, 1}
\to{\longleftarrow\joinrel\relbar\joinrel%
\relbar\joinrel\relbar\joinrel\relbar\joinrel\llongrightarrow}}
\hskip6pt  {{(H^\vee)}_\h}^{\!\!\vee}{\Big|}_{\h=1} = H^\vee \; .  $$
   \indent   Now look at  $ \, {\big( {{(H^\vee)}_\h}^{\!\!\vee}
\big)}' \, $:  by construction, we have  $ \, {\big( {{( H^\vee
)}_\h}^{\!\!\vee} \big)}'{\Big|}_{\h=1} \! = {{( H^\vee )}_\h}^{\!
\!\vee}{\Big|}_{\h=1} \! = H^\vee \, $,  \, whereas  $ \, {\big(
{{(H^\vee)}_\h}^{\!\!\vee} \big)}'{\Big|}_{\h=1} \! = F[K_-] \, $
for some connected algebraic Poisson group  $ K_- \, $:  in addition,
if  $ \, \text{\it Char}\,(\Bbbk) = 0 \, $  then  $ \, K_- = G_-^\star
\, $  by  Theorem 2.2{\it (c)}.  So  $ \, {\big( {{(H^\vee)}_\h}^{\!
\!\vee} \big)}' \, $  can be thought of as a 1-parameter family inside
$ \HA_\Bbbk \, $,  with regular fibers, linking  $ H^\vee $  and
$ F[G_-^\star] $  as (polynomial) deformations of each other, namely
  $$  H^\vee = {\big( {{(H^\vee)}_\h}^{\!\!\vee} \big)}'{\Big|}_{\h=1}
\hskip-4pt  \underset{{({{(H^\vee)}_\h}^{\!\!\vee})}'}  \to
{\overset{1 \,\leftarrow\, \h \,\rightarrow\, 0}
\to{\longleftarrow\joinrel\relbar\joinrel%
\relbar\joinrel\relbar\joinrel\relbar\joinrel\llongrightarrow}}
\hskip0pt  {\big( {{(H^\vee)}_\h}^{\!\!\vee} \big)}'{\Big|}_{\h=0}
\hskip-5pt  = F[K_-]
\hskip4pt  \Big( = F[G_-^\star]  \,\text{\ if \ }  \text{\it Char}\,
(\Bbbk) = 0 \Big) \, .  $$
Therefore  $ H^\vee $  is both a deformation of an enveloping algebra
and a deformation of a function algebra, via two different 1-parameter
families (with regular fibers) in  $ \HA_\Bbbk $  which match at the
value  $ \, \h = 1 \, $,  \, corresponding to the common element
$ \, H^\vee \, $.  At a glance,
  $$  \U(\gerg_-)  \hskip2pt
\underset{{{(H^\vee)}_\h}^{\!\!\vee}}  \to
{\overset{0 \,\leftarrow\, \h \,\rightarrow\, 1}
\to{\longleftarrow\joinrel\relbar\joinrel%
\relbar\joinrel\relbar\joinrel\relbar\joinrel\llongrightarrow}}
\hskip2pt  H^\vee  \hskip1pt
\underset{\;{({{(H^\vee)}_\h}^{\!\!\vee})}'}  \to
{\overset{1 \,\leftarrow\, \h \,\rightarrow\, 0}
\to{\longleftarrow\joinrel\relbar\joinrel%
\relbar\joinrel\relbar\joinrel\relbar\joinrel\llongrightarrow}}
\hskip2pt  F[K_-]
\hskip4pt  \Big( = F[G_-^\star]  \,\text{\ if \ }  \!\text{\it Char}\,
(\Bbbk) = 0 \Big) \, .   \eqno (5.3)  $$
   \indent   Now consider  $ {H_\h}' $.  We have  $ \,
{H_\h}'{\Big|}_{\h=0} := {H_\h}' \Big/ \h \, {H_\h}' =
G_{\underline{D}}(H) =: \widetilde{H} \, $,  \, and  $ \,
\widetilde{H} = F[G_+] \, $  for some connected algebraic
Poisson group  $ G_+ \, $.  On the other hand, we have also
$ \, {H_\h}'{\Big|}_{\h=1} := {H_\h}' \Big/ (\h-1) \, {H_\h}' =
\sum_{n \geq 0} D_n =: H' \, $;  \, note that the latter is a
Hopf subalgebra of  $ H $,  because  $ \underline{D} $  is a
Hopf algebra filtration; moreover we have  $ \, \widetilde{H}
= \widetilde{H'} \, $,  by the very definitions.  Therefore we
can think at  $ \, {H_\h}' = \Cal{R}^\h_{\underline{D}}(H') \, $
as a 1-parameter family inside  $ \HA_\Bbbk $  with regular fibers
which links  $ \widetilde{H} $  and  $ H' $  as (polynomial)
deformations of each other, namely
  $$  F[G_+] = \widetilde{H} = {H_\h}'{\Big|}_{\h=0}  \hskip5pt
\underset{\;{H_\h}^{\!\prime}}  \to
{\overset{0 \,\leftarrow\, \h \,\rightarrow\, 1}
\to{\longleftarrow\joinrel\relbar\joinrel%
\relbar\joinrel\relbar\joinrel\relbar\joinrel\llongrightarrow}}
\hskip6pt  {H_\h}'{\Big|}_{\h=1} = H' \; .  $$
   \indent   Consider also  $ \, {\big({H_\h}'\big)}^\vee \, $:  by
construction, we have  $ \, {\big({H_\h}'\big)}^\vee{\Big|}_{\h=1}
\! = {H_\h}'{\Big|}_{\h=1} \! = H' \, $,  \, whereas  $ \, {\big(
{H_\h}' \big)}^\vee{\Big|}_{\h=0} \! = \U(\gerk_+) \, $  for some
restricted Lie bialgebra  $ \gerk_+ \, $:  in addition, if  $ \,
\Char(\Bbbk) = 0 \, $  then  $ \, \gerk_+ = \gerg_+^{\,\times} \, $
thanks to  Theorem 2.2{\it (c)}.  Thus  $ \, {\big({H_\h}'\big)}^\vee
\, $  can be thought of as a 1-parameter family inside  $ \HA_\Bbbk $
with regular fibers which links  $ \U(\gerk_+) $  and  $ H' $  as
(polynomial) deformations of each other, namely
  $$  H' = {\big({H_\h}'\big)}^\vee{\Big|}_{\h=1}
\hskip3pt  \underset{\;{({H_\h}^{\!\prime})}^\vee}  \to
{\overset{1 \,\leftarrow\, \h \,\rightarrow\, 0}
\to{\longleftarrow\joinrel\relbar\joinrel%
\relbar\joinrel\relbar\joinrel\relbar\joinrel\llongrightarrow}}
\hskip4pt  {\big({H_\h}'\big)}^\vee{\Big|}_{\h=0} = \U(\gerk_+)
\hskip7pt  \Big( = U(\gerg_+^{\,\times})  \,\text{\ if \ }
\text{\it Char}\,(\Bbbk) = 0 \Big) \, .  $$
Therefore,  $ H' $  is at the same time a deformation of a function
algebra and a deformation of an enveloping algebra, via two different
1-parameter families inside  $ \HA_\Bbbk $  (with regular fibers)
which match at the value  $ \, \h = 1 \, $,  \, corresponding (in
both families) to  $ \, H' \, $.  In short,
  $$  F[G_+]  \hskip3pt
\underset{\;{H_\h}^{\!\prime}}  \to
{\overset{0 \,\leftarrow\, \h \,\rightarrow\, 1}
\to{\longleftarrow\joinrel\relbar\joinrel%
\relbar\joinrel\relbar\joinrel\relbar\joinrel\llongrightarrow}}
\hskip2pt  H'  \hskip1pt
\underset{\;{({H_\h}^{\!\prime})}^\vee}  \to
{\overset{1 \,\leftarrow\, \h \,\rightarrow\, 0}
\to{\longleftarrow\joinrel\relbar\joinrel%
\relbar\joinrel\relbar\joinrel\relbar\joinrel\llongrightarrow}}
\hskip3pt  \U(\gerk_+)
\hskip7pt  \Big( = U(\gerg_+^{\,\times})  \,\text{\ if \ }
\text{\it Char}\,(\Bbbk) = 0 \Big) \, .   \eqno (5.4)  $$
   \indent   Finally, it is worth noticing that in the special
case  $ \, H' = H = H^\vee \, $  we can splice together (5.3)
and (5.4) to get
  $$  \eqalignno{
   \hskip-17pt   F[G_+]  \hskip3pt
\underset{\;{H_\h}^{\!\prime}}  \to
{\overset{0 \,\leftarrow\, \h \,\rightarrow\, 1}
\to{\longleftarrow\joinrel\relbar\joinrel%
\relbar\joinrel\relbar\joinrel\relbar\joinrel\llongrightarrow}}
\hskip2pt {}  &  H'  \hskip1pt
\underset{\;{({H_\h}^{\!\prime})}^\vee}  \to
{\overset{1 \,\leftarrow\, \h \,\rightarrow\, 0}
\to{\longleftarrow\joinrel\relbar\joinrel%
\relbar\joinrel\relbar\joinrel\relbar\joinrel\llongrightarrow}}
\hskip3pt  \U(\gerk_+)
\hskip7pt  \Big( = U(\gerg_+^{\,\times})  \,\text{\ if \ }
\text{\it Char}\,(\Bbbk) = 0 \Big)   &  \cr
   &  \hskip1pt   ||   &  \cr
   &  H_{\phantom{\displaystyle I}}   & (5.5)  \cr
   &  \hskip1pt   ||   &  \cr
   \hskip-17pt   \U(\gerg_-)  \hskip2pt
\underset{{{(H^\vee)}_\h}^{\!\!\vee}}  \to
{\overset{0 \,\leftarrow\, \h \,\rightarrow\, 1}
\to{\longleftarrow\joinrel\relbar\joinrel%
\relbar\joinrel\relbar\joinrel\relbar\joinrel\llongrightarrow}}
\hskip2pt {}  &  H^\vee  \hskip1pt
\underset{\;{({{(H^\vee)}_\h}^{\!\!\vee})}'}  \to
{\overset{1 \,\leftarrow\, \h \,\rightarrow\, 0}
\to{\longleftarrow\joinrel\relbar\joinrel%
\relbar\joinrel\relbar\joinrel\relbar\joinrel\llongrightarrow}}
\hskip2pt  F[K_-]
\hskip4pt  \Big( = F[G_-^\star]  \,\text{\ if \ }  \!\text{\it Char}\,
(\Bbbk) = 0 \Big)   &  \cr}  $$
which gives  {\sl four\/}  different regular 1-parameter deformations
from  $ H $  to Hopf algebras encoding geometrical objects of Poisson
type (i.e.~Lie bialgebras or Poisson algebraic groups).

\vskip7pt

  {\bf 5.6 The function algebra case.} \, Let  $ G $  be any algebraic
group over the field  $ \Bbbk $.  Let  $ \, R := \Bbbk[\h] \, $  be
as in \S 5.1, and set  $ \, F_\h[G] := {\big( F[G] \big)}_\h = R
\otimes_\Bbbk F[G] \, $:  \, this is trivially a QFA at  $ \h \, $,  for
$ \, F_\h[G] \big/ \h \, F_\h[G] = F[G] \, $,  inducing on  $ G $  the
trivial Poisson structure, so that its cotangent Lie bialgebra is
simply  $ \gerg^\times $  with trivial Lie bracket and Lie cobracket
dual to the Lie bracket of  $ \gerg \, $.  In the sequel we identify
\hbox{$ F[G] $  with  $ \, 1 \! \otimes \! F[G] \subset F_\h[G] \, $.}
                                          \par
   We begin by computing  $ \, {F_\h[G]}^\vee \, $  (w.r.t.~the
non-zero element $ \h \, $).
                                          \par
   Let  $ \, J := J_{F[G]} \equiv \text{\sl Ker}\, \big( \epsilon_{F[G]}
\big) \, $,
%
%
%
%
let
$ \, {\{y_b\}}_{b \in \Cal{S}} \, $  be a  $ \Bbbk $--basis  of  $ \,
Q\big(F[G]\big) = J \Big/ J^2 = \gerg^\times \, $,  \, and pull it
back to a subset $ \, {\{j_b\}}_{b \in \Cal{S}} \, $  of  $ J \, $.
Then we see that  $ \, J^n \big/ J^{n+1} \, $  is a  $ \Bbbk $--vector
space spanned by the set of (cosets of) ordered monomials (using
multiindices and all the notation introduced in the proof of
Theorem 4.7)  $ \, \big\{\, j^{\,\underline{e}} \mod J^{n+1}
\;\big|\; \underline{e} \in \N^{\Cal{S}}_f \, , \; |\underline{e}|
= n \,\big\} \, $  where  $ \, |\underline{e}| := \sum_{b \in \Cal{S}}
\underline{e}(b) \, $;  \, therefore  $ \, I^n \big/ I^{n+1} \, $  as
a  $ \Bbbk $--vector  space is spanned by  $ \, \big\{\, \h^{e_0}
j^{\,\underline{e}} \mod I^{n+1} \;\big|\; e_0 \in \N, \underline{e}
\in \N^{\Cal{S}}_f \, , \; e_0 + |\underline{e}| = n \,\big\} \, $.
Noting that  $ \, \h^{-n} I^{n+1} = \h \cdot \h^{-(n+1)} I^{n+1}
\equiv 0 \mod \h \, {F_\h[G]}^\vee \, $,  \, we can argue that
$ \, \h^{-n} I^n \mod \h \, {F_\h[G]}^\vee \, $  is spanned over
$ \Bbbk $  by  $ \, \big\{\, \h^{-|\underline{e}|} j^{\,\underline{e}}
\mod \h \, {F_\h[G]}^\vee \;\big|\; \underline{e} \in \N^{\Cal{S}}_f
\, , \; |\underline{e}| \leq n \,\big\} \, $.  Now we have to
distinguish the various cases.
                                          \par
  First of all, let  $ \, y \in F[G] \, $  be idempotent: switching
if necessary to  $ \, y_+ := y - \epsilon(y) \, $  we can assume that
$ \, y \in J \, $.  Then  $ \, y = y^2 = \cdots = y^n \in J^n \subset
I^n \, $  for all  $ \, n \in \N \, $,  so that  $ \, y \equiv 0 \mod
\h \, F_\h[G] \, $.  Thus in order to compute  $ \, F_\h[G] \Big/ \h \,
F_\h[G] \, $  the idempotents of  $ F[G] $  are irrelevant: this means
$ \, F_\h[G] \Big/ \h \, F_\h[G] = F_\h\big[G^0\big] \Big/ \h \,
F_\h\big[G^0\big] \, $,  where  $ G^0 $  is the connected component
of  $ G \, $;  \, thus we can assume from scratch that  $ G $  be
connected.
                                          \par
   First assume  $ \, G \, $  is  {\sl smooth},  i.e.~$ F[G] $
is  {\sl reduced},  which is always the case if  $ \hbox{\it Char}
\,(\Bbbk) = 0 \, $.  Then the set above is a basis of  $ \, \h^{-n}
I^n \mod \h \, {F_\h[G]}^\vee \, $:  \, for if we have a non-trivial
linear combination of the elements of this set which is zero,
multiplying by  $ \h^n $  gives an element  $ \, \gamma_n \in
I^n \setminus I^{n+1} \, $  such that  $ \, \h^{-n} \gamma_n
\equiv 0 \mod \h \, {F_\h[G]}^\vee \, $;  \, then there is
$ \, \ell \in \N \, $  such that  $ \, \h^{-n} \gamma_n
\in \h \cdot \h^{-\ell} I^\ell \, $,  \, so  $ \, \h^\ell
\gamma_n = \h^{1+n} I^\ell \subseteq I^{1+n+\ell} \, $,
\, whence clearly  $ \, \gamma_n \in I^{n+1} $,  \, a
contradiction!  The outcome is that  $ \, \big\{\,
\h^{-|\underline{e}|} j^{\,\underline{e}} \mod \h \,
{F_\h[G]}^\vee \;\big|\; \underline{e} \in \N^{\,\Cal{S}}_f
\,\big\} \, $  is a  $ \Bbbk $--basis  of  $ \,
{F_\h[G]}^\vee{\Big|}_{\h=0} := {F_\h[G]}^\vee \big/
\h \, {F_\h[G]}^\vee \, $.  Now let  $ \, j_\beta^\vee :=
\h^{-1} j_\beta \, $  for all  $ \, \beta \in \Cal{S} \, $.
By the previous analysis  $ \, {F_\h[G]}^\vee{\Big|}_{\h=0}
\, $  has  $ \Bbbk $--basis  $ \, \big\{ {\big( j^\vee
\big)}^{\underline{e}} \mod \h \, {F_\h[G]}^{\!\vee} \;\big|\;
\underline{e} \in \N^{\Cal{S}}_f \,\big\} \, $,  \, hence
the Poincar\'e-Birkhoff-Witt theorem tells us that  $ \,
{F_\h[G]}^\vee{\Big|}_{\h=0} = U(\gerh) \equiv S(\gerh) \, $  as
associative algebras, where  $ \, \gerh \, $  is the Lie algebra
spanned by  $ \, {\big\{ j^\vee_b \mod \h \, {F_\h[G]}^\vee
\big\}}_{b \in \Cal{S}} \, $  (as in the proof of Theorem 4.7),
whose Lie bracket is trivial for it is given by  $ \, \big[ j^\vee_b,
j^\vee_\beta \big] := \h^{-2} (j_b \, j_\beta - j_\beta \, j_b)
\mod \h \, {F_\h[G]}^\vee \equiv 0 \, $.  Further, we have also
$ \; \Delta \big( j^\vee \big) \equiv j^\vee \otimes 1 + 1 \otimes
j^\vee \mod \, \h \, {\big( F_\h[G]^\vee \big)}^{\otimes 2} \; $
for all  $ \, j^\vee \in J^\vee := \h \, J $  (cf.~the proof of
Theorem 3.7), whence  $ \; \Delta(\text{j}) = \text{j} \otimes 1
+ 1 \otimes \text{j} \; $  for all  $ \, \text{j} \in \gerh \, $,
\, so  $ \, {F_\h[G]}^\vee{\Big|}_{\h=0} \!\! = U(\gerh) \equiv
S(\gerh) \, $  as  {\sl Hopf\/}  algebras too.  Now, consider the
linear map  $ \; \sigma \, \colon \, \gerg^\times = J \Big/ J^2
\longrightarrow \gerh \; \big( \! \subset \! U(\gerh) \big) \; $
given by  $ \, y_b \mapsto j^\vee_b \, $  ($ \, b \in \Cal{S} \, $).
By construction this is clearly a vector space isomorphism, and also
a  {\sl Lie algebra\/}  isomorphism, since the Lie bracket is trivial
on both sides  ($ G $  has the trivial Poisson structure!).  In
addition, one has
 \vskip-15pt
  $$  \displaylines{
   \; \Big\langle u_1 \otimes u_2 \, , \, \delta_\gerh\big(
\sigma(y_b) \big) \Big\rangle = \Big\langle u_1 \otimes u_2
\, , \, \h^{-1} \big( \Delta - \Delta^{\text{op}} \big) \big(
\sigma(y_b) \big) \mod \h \Big\rangle =   \hfill  \cr
   \hfill   = \Big\langle u_1 \otimes u_2 \, , \, \h^{-2}
\big( \Delta - \Delta^{\text{op}} \big) (j_b) \mod \h \Big\rangle
= \Big\langle \big( u_1 \cdot u_2 - u_2 \cdot u_1 \big) \, , \,
\h^{-2} j_b \mod \h \Big\rangle =  \cr
   \;  = \Big\langle \hskip-1pt [u_1,u_2] \, , \, \h^{-2} j_b
\mod \h \Big\rangle = \Big\langle u_1 \otimes u_2 \, , \, \h^{-2}
\delta_{\gerg^\times}(y_b) \mod \h \Big\rangle = \Big\langle
\hskip-1pt u_1 \otimes u_2 \, , (\sigma \otimes \sigma)
\big( \delta_{\gerg^\times}(y_b) \big) \hskip-2pt
\Big\rangle  \cr }  $$
for all  $ \, u_1, u_2 \in \gerg \, $  (with  $ \, (u_1 \cdot u_2
- u_2 \cdot u_1) \in U(\gerg) \, $)  and  $ \, b \in \Cal{S} \, $,
\, which is enough to prove that  $ \, \delta_\gerh \circ \sigma
= (\sigma \otimes \sigma) \circ \delta_{\gerg^\times} \, $,
\, i.e.~$ \sigma $  is a Lie  {\sl bi\/}algebra  morphism as
well.  Therefore the outcome is  $ \; {F_\h[G]}^\vee{\Big|}_{\h=0}
\! = \, U(\gerg^\times) \, \equiv \, S(\gerg^\times) \; $  as
{\sl co-Poisson Hopf algebras}.
                                             \par
   Another extreme case is when  $ G $  is a  {\sl finite
connected group scheme\/}:  then, assuming for simplicity that
$ \, \Bbbk \, $  be perfect, we have  $ \, F[G] = \Bbbk[x_1,
\dots, x_n] \Big/ \big( x_1^{p^{e_1}}, \dots, x_n^{p^{e_n}}
\big) \, $  for some  $ \, n , e_1, \dots, e_n \in \N \, $.
The previous analysis, with minor changes, then shows that  $ \,
{F_\h[G]}^\vee{\Big|}_{\h=0} $  is now a quotient of  $ \,
U(\gerg^\times) \equiv S(\gerg^\times) \, $:  \, namely, the
$ x_i $'s  take place of the  $ j_b $'s,  so the cosets of the
$ x^\vee_i $'s  ($ i = 1, \dots, n $)  modulo  $ \, \h \,
{F_\h[G]}^\vee \, $  generate  $ \gerg^\times $,  and we find
$ \, {F_\h[G]}^\vee{\Big|}_{\h=0} = U(\gerg^\times) \Big/ \big(
{(x^\vee_1)}^{p^{e_1}}, \dots, {(x^\vee_n)}^{p^{e_n}} \big) \, $.
Now, recall that for any Lie algebra  $ \gerh $  we can consider
$ \, \gerh^{{[p\hskip0,7pt]}^\infty} := \Big\{\, x^{{[p\hskip0,7pt]}^n}
\! := x^{p^n} \,\Big\vert\, x \in \gerh \, , n \in \N \,\Big\} \, $,
\, the  {\sl restricted Lie algebra generated by  $ \gerh $}  inside
$ U(\gerh) $,  with the  $ p $--operation  given by  $ \, x^{[p
\hskip0,7pt]} := x^p \, $;  \, then one always has  $ \, U(\gerh)
= \u\big( \gerh^{{[p\hskip0.7pt]}^\infty} \big) \, $.  In the present
case, the subset  $ \, \Big\{ {(x^\vee_1)}^{p^{e_1}}, \, \ldots, \,
{(x^\vee_n)}^{p^{e_n}} \Big\} \, $  generates a  $ p $--ideal
$ \Cal{I} $  of  $ (\gerg^\times)^{{[p\hskip0,7pt]}^\infty} $,  so
$ \, \gerg^\times_{\text{\it res}} := \gerg^{{[p\hskip0,7pt]}^\infty}
\! \Big/ \Cal{I} \, $  is a restricted Lie algebra too, with  $ \,
\Big\{ {(x^\vee_1)}^{p^{a_1}}, \, \ldots, \, {(x^\vee_n)}^{p^{a_n}}
\,\Big|\; a_1 < e_1, \dots, a_n < e_n \Big\} \, $  as a  $ \Bbbk $--basis.  Then
the previous
analysis proves  $ \; {F_\h[G]}^\vee{\Big|}_{\h=0} = \u\left(
\gerg^\times_{\text{\it res}} \right) \equiv S(\gerg^\times) \Big/
\Big( \Big\{ {(x^\vee_1)}^{p^{e_1}}, \, \ldots, \, {(x^\vee_n)}^{p^{e_n}}
\Big\} \Big) \; $  as co-Poisson Hopf algebras.
                                             \par
   The general case is intermediate; we get it via the relation  $ \,
F_\h[G] \Big/ \h \, F_\h[G] = G_{\underline{J}}\,\big(F[G]\big) \, $.
                                             \par
Assume again for simplicity that
$ \Bbbk $  be perfect.  Let  $ \, F[[G]] \, $  be the  $ J $--adic
completion of  $ \, H = F[G] \, $.  By standard results on algebraic
groups (cf.~[DG]) there is a subset  $ \, {\{x_i\}}_{i \in \Cal{I}}
\, $  of  $ \, J \, $  such that  $ \, {\big\{\, \overline{x}_i :=
x_i \! \mod J^2 \,\big\}}_{i \in \Cal{I}} \, $  is a basis of  $ \,
\gerg^\times = J \big/ J^2 \, $  and  $ \; F[[G]] \, \cong \,
\Bbbk\big[\big[ {\{x_i\}}_{i \in \Cal{I}} \big]\big] \Big/ \Big(
\Big\{ x_i^{\,p^{n(x_i)}} \Big\}_{i \in \Cal{I}_0} \Big) \; $  (the
algebra of truncated formal power series), for some subset  $ \,
\Cal{I}_0 \subset \Cal{I} \, $  and some  $ \, {\big( n(x_i)
\big)}_{i \in \Cal{I}_0} \in \N^{\,\Cal{I}_0} $.  Since  $ \,
G_{\underline{J}}\,\big(F[G]\big) = G_{\underline{J}}\,\big(
F[[G]] \big) \, $,  \, we argue that  $ \; G_{\underline{J}}\,\big(
F[G] \big) \, \cong \, \Bbbk\big[ {\{ \overline{x}_i \}}_{i \in
\Cal{I}} \big] \Big/ \Big( \Big\{ \overline{x}_i^{\,p^{n(x_i)}}
\Big\}_{i \in \Cal{I}_0} \,\Big) \; $;  \; finally, since  $ \;
\Bbbk\big[ {\{ \overline{x}_i \}}_{i \in \Cal{I}} \big] \cong
S(\gerg^\times) \; $  we get  $ \, G_{\underline{J}}\,\big(F[G]\big)
\cong S(\gerg^\times) \bigg/ \! \Big( \Big\{\, \overline{x}^{\,
p^{n(x)}} \Big\}_{x \in \Cal{N}(F[G])} \,\Big) \, $  as
$ \Bbbk $--algebras,  where  $ \Cal{N} \big( F[G] \big) $
is the nilradical of  $ F[G] $  and  $ \, p^{n(x)} $  is the
order of nilpotency of  $ \, x \in \Cal{N}\big(F[G]\big) $.
                                        \par
   Now, let  $ \, 0 \not= \overline{\eta} \in J \big/ J^2 \, $,  \,
and let  $ \, \eta \in J \, $  be a lift of  $ \, \overline{\eta}
\, $:  \, then  $ \, \Delta(\eta) = \epsilon(\eta) \cdot 1 \otimes 1
+ \delta_1(\eta) \otimes 1 + 1 \otimes \delta_1(\eta) + \delta_2(\eta)
= \eta \otimes 1 + 1 \otimes \eta + \delta_2(\eta) \, $,  \, by the
very definitions.  But  $ \, \delta_2(\eta) \in J \otimes J \, $,
\, hence  $ \, \overline{\delta_2(\eta)} = \overline{0} \in
%
%
%
%
G_{\underline{J}}\,\big(F[G]\big) \otimes G_{\underline{J}}\,
\big(F[G]\big)
%
%
%
%
\, $:  \,
so  $ \, \Delta(\overline{\eta}\,) := \overline{\Delta(\eta)}
%
%
%
%
= \overline{\eta}
\otimes 1 + 1 \otimes \overline{\eta} \, $.  Therefore, all
elements of  $ \, J \big/ J^2 = \gerg^\times \, $  are primitive:
this implies that the previous isomorphism respects also the Hopf
structure.  As for the Lie cobracket of  $ \, G_{\underline{J}}\,
\big(F[G]\big) \, $,  by construction it is given by  $ \,
\delta_{G_{\underline{J}}\,(F[G])}(\overline{x}\,) :=
\overline{\nabla(x)} = \overline{\Delta(x) - \Delta^{\text{op}}(x)}
\, $.  Now, in the natural pairing between  $ F[G] $  and  $ U(\gerg) $,
for all  $ \, x \in J \, $  we have  $ \, \big\langle \nabla(x), Y \otimes
Z \big\rangle = \big\langle \Delta(x) - \Delta^{\text{op}}(x), Y \otimes Z
\big\rangle = \big\langle x, Y \, Z - Z \, Y \big\rangle = \big\langle x,
[Y,Z] \big\rangle \, $,  hence  $ \, \big\langle \delta_{G_{\underline{J}}
\,(F[G])}(\overline{x}\,), Y \otimes Z \big\rangle = \big\langle x, [Y,Z]
\big\rangle \, $  for all  $ \, Y $,  $ Z \in \gerg \, $;  similarly
$ \, \big\langle \delta_{\gerg^\times}(\overline{x}\,), y \otimes z
\big\rangle = \big\langle x, [Y,Z] \big\rangle \, $  for all  $ \, Y $,
$ Z \in \gerg \, $.  We then argue that  $ \, \delta_{G_{\underline{J}}
\,(F[G])}(\overline{x}\,) = \delta_{\gerg^\times}(\overline{x}\,) \, $
for all  $ \, x \in \gerg^\times \, $,  \, whence the two Lie cobrackets
do correspond to one another in the isomorphism above.  Since  $ \,
{F_\h[G]}^\vee{\Big|}_{\h=0} = G_{\underline{J}}\,\big(F[G]\big) \, $,
\, the outcome is that  $ \; {F_\h[G]}^\vee{\Big|}_{\h=0} \, \cong \,
S(\gerg^\times) \bigg/ \! \Big( \Big\{\, \overline{x}^{\,p^{n(x)}}
\,\Big|\, x \in \Cal{N}\big(F[G]\big) \Big\} \Big) \; $  as
co-Poisson Hopf algebras.
                                             \par
   Note also that the description of  $ \, {F_\h[G]}^\vee \, $  in
the general case is exactly like the one we gave for the smooth
case: one simply has to mod out the ideal generated by  $ \,
{\Cal{N}\big(F[G]\big)}^\vee := \h^{-1} \Cal{N}\big(F[G]\big)
\, $,  i.e.~(roughly) to set  $ \; {(x_i^\vee)}^{p^{n(x_i)}}
\hskip-3pt = 0 \; $  (with  $ \, x_i^\vee := \h^{-1} x_i \, $)
for all  $ \, i \in \Cal{I} \, $.  By the way, to have this
description we do not need  $ \Bbbk $  to be perfect.  As for
$ \, F[G]^\vee := F[G] \big/ J^\infty \, $,  it is known (cf.~[Ab],
Lemma 4.6.4) that  $ \, F[G]^\vee =F[G] \, $  whenever  $ G $  is
finite dimensional and there exists no  $ \, f \in F[G] \setminus
\Bbbk \, $  which is separable algebraic over  $ \Bbbk \, $.
                                             \par
   It is also interesting to consider  $ {\big( {F_\h[G]}^\vee
\big)}' $.  If  $ \, \text{\it Char}\,(\Bbbk) = 0 \, $,  then the
proof of Proposition 4.3 does work, with no special simplifications,
giving  $ \, {\big( {F_\h[G]}^\vee \big)}' = F_\h[G] \, $.  If instead
$ \, \text{\it Char}\,(\Bbbk) = p > 0 \, $,  then the situation might
change dramatically.  Indeed, if the group  $ G $  has height 1   ---
i.e., if  $ \, F[G] = \Bbbk\big[ {\{x_i\}}_{i \in \Cal{I}} \big] \big]
\Big/ \big( \big\{ x_i^p \,\big|\, i \in \Cal{I} \big\} \big) \, $ 
as a  $ \Bbbk $--algebra ---   then the same analysis as in the
characteristic zero case may be applied, with a few minor changes,
whence one gets again  $ \, {\big( {F_\h[G]}^\vee \big)}' = F_\h[G] \, $.
Otherwise, let  $ \, y \in J \setminus \{0\} \, $  be primitive and such
that  $ \, y^p \not= 0 \, $:  for instance, this occurs for  $ \, F[G]
= \Bbbk[x] \, $,  i.e.~$ \, G \cong \Bbb{G}_a \, $,  and  $ \, y =
x \, $.  Then  $ \, y^p \, $  is primitive as well, hence  $ \,
\delta_n(y^p) = 0 \, $  for each  $ \, n > 1 \, $.  It follows
that  $ \, 0 \not= \h \, {(y^\vee)}^p \in {\big( {F_\h[G]}^\vee
\big)}' \, $,  \, whereas  $ \, \h \, {(y^\vee)}^p \not\in F_\h[G]
\, $,  as follows from our previous description of  $ {F_\h[G]}^\vee $.
Thus  $ \, {\big( {F_\h[G]}^\vee \big)}' \supsetneqq {F_\h[G]}^\vee
\, $,  \, a  {\sl counterexample\/}  to Proposition 4.3.
                                             \par
   What for  $ \, {F_\h[G]}' \, $  and  $ \, \widetilde{F[G]} \, $?
Again, this depends on the group  $ G $  under consideration.  We
give two simple examples, both ``extreme'', in a sense, and
opposite to each other.
                                             \par
   Let  $ \, G := \Bbb{G}_a = \hbox{\it Spec}\big(\Bbbk[x]\big)
\, $,  \, so  $ \, F[G] = F[\Bbb{G}_a] = \Bbbk[x] \, $  and  $ \,
F_\h[\Bbb{G}_a] := R \otimes_\Bbbk \Bbbk[x] = R[x] \, $.  Then
since  $ \, \Delta(x) := x \otimes 1 + 1 \otimes x \, $  and  $ \,
\epsilon(x) = 0 \, $  we find  $ \, {F_\h[\Bbb{G}_a]}' = R[\h{}x]
\, $  (like in \S 5.7 below: indeed, this is just a special instance,
for  $ \, F[\Bbb{G}_a] = U(\gerg) \, $  where  $ \gerg $  is the
1-dimensional Lie algebra).  Moreover, iterating one gets easily
$ \, {\big( {F_\h[\Bbb{G}_a]}' \big)}' = R\big[\h^2{}x\big] \, $,
$ \, {\Big( {\big( {F_\h [\Bbb{G}_a]}' \big)}' \Big)}' = R \big[
\h^3{}x \big] \, $,  \, and in general  $ \; \Big( \Big( \big(
F_\h [\Bbb{G}_a \underbrace{]'{\big)}'{\Big)}' \cdots
{\Big)}'}_{n} = R\big[\h^n{}x\big] \cong R[x] = F_\h
[\Bbb{G}_a] \; $  for all  $ \, n \in \N \, $.
                                             \par
   Second, let  $ \, G := \Bbb{G}_m = \hbox{\it Spec} \, \Big(
\Bbbk \big[ z^{+1}, z^{-1} \big] \Big) \, $,  \, that is  $ \,
F[G] = F[\Bbb{G}_m] = \Bbbk \big[ z^{+1}, z^{-1} \big] \, $  so
that  $ \, F_\h[\Bbb{G}_m] := R \otimes_\Bbbk \Bbbk \big[ z^{+1},
z^{-1} \big] = R \big[ z^{+1}, z^{-1} \big] \, $.  Then since
$ \, \Delta \big( z^{\pm 1} \big) := z^{\pm 1} \otimes z^{\pm 1}
\, $  and  $ \, \epsilon \big( z^{\pm 1} \big) = 1 \, $  we
find  $ \, \Delta^n \big( z^{\pm 1} \big) = {\big( z^{\pm 1}
\big)}^{\otimes n} \, $  and  $ \, \delta_n \big( z^{\pm 1}
\big) = {\big( z^{\pm 1} - 1 \big)}^{\otimes n} \, $  for
all  $ \, n \in \N \, $.  From that it follows easily
$ \, {F_\h[\Bbb{G}_m]}' = R \cdot 1 \, $,  \, the
trivial possibility (see also \S 5.13 later on).

\vskip7pt

  {\bf 5.7 The enveloping algebra case.} \, Let  $ \gerg $  be any
Lie algebra over the field  $ \Bbbk $,  and  $ U(\gerg) $  its
universal enveloping algebra with its standard Hopf structure.
Assume  $ \, \Char(\Bbbk) = 0 \, $,  and let  $ \, R = \Bbbk[\h] \, $
as in \S 5.1,  and set  $ \, U_\h(\gerg) := R \otimes_\Bbbk U(\gerg) =
{\big( U(\gerg) \big)}_\h \, $.  Then  $ U_\h(\gerg) $  is trivially
a QrUEA at  $ \h $,  for  $ \, U_\h(\gerg) \big/ \h \, U_\h(\gerg) =
U(\gerg) \, $,  inducing on  $ \gerg $  the trivial Lie cobracket.
Therefore the dual Poisson group is nothing but  $ \gerg^\star $
(the topological dual of  $ \gerg $  w.r.t.~the weak topology), an
Abelian group w.r.t.~addition, with  $ \gerg $  as cotangent Lie
bialgebra and function algebra  $ \, F[\gerg^\star] = S(\gerg) \, $:
\, the Hopf structure is the standard one, given by  $ \, \Delta(x)
= x \otimes 1 + 1 \otimes x \, $  (for all  $ \, x \in \gerg \, $),
and the Poisson structure is the one induced by  $ \, \{x,y\} :=
[x,y] \, $  for all  $ \, x $,  $ y \in \gerg \, $.  This is
the so-called Kostant-Kirillov structure on  $ \gerg^\star $.
                                          \par
   Similarly, if  $ \, \text{\it Char}\,(\Bbbk) = p > 0 \, $  and
$ \gerg $  is any restricted Lie algebra over  $ \Bbbk $,  let
$ \, \u(\gerg) \, $  be its restricted universal enveloping algebra,
with its standard Hopf structure.  Then if  $ \, R = \Bbbk[\h] \, $
the Hopf  $ R $--algebra  $ \, U_\h(\gerg) := R \otimes_\Bbbk \u(\gerg)
= {\big(\u(\gerg)\big)}_\h \, $  is a QrUEA at  $ \h $,  because  $ \,
\u_\h(\gerg) \big/ \h \, \u_\h(\gerg) = \u(\gerg) \, $,  inducing on
$ \gerg $  the trivial Lie cobracket: then the dual Poisson group is
again  $ \gerg^\star $, with cotangent Lie bialgebra  $ \gerg $  and
function algebra  $ \, F[\gerg^\star] = S(\gerg) \, $  (the Poisson
Hopf structure being as above).  Recall also that  $ \, U(\gerg) =
\u\big( \gerg^{{[p\hskip0.7pt]}^\infty} \big) \, $  (cf.~\S 5.6).
                                          \par
   First we compute  $ \, {\u_\h(\gerg)}' $  (w.r.t.~the prime
$ \h \, $)  using (5.2), i.e.~computing
              \hbox{the filtration  $ \underline{D} \, $.}
                                          \par
   By the PBW theorem, once an ordered basis  $ B $  of  $ \gerg $  is
fixed  $ \u(\gerg) $  admits as basis the set of ordered monomials
in the elements of  $ B $  whose degree (w.r.t.~each element of  $ B $)
is less than  $ p \, $;  this yields a Hopf algebra filtration of
$ \u(\gerg) $  by the total degree, which we refer to as  {\sl the
standard filtration}.  Then from the very definitions a straightforward
calculation shows that  $ \underline{D} $  coincides with the standard
filtration.  This together with (5.2) immediately implies  $ \,
{\u_\h(\gerg)}' = \langle \tilde{\gerg} \rangle = \langle \h \,
\gerg \rangle \, $:  \, hereafter  $ \, \tilde{\gerg} := \h \,
\gerg \, $,  \, and similarly we set  $ \, \tilde{x} := \h \, x \, $
for all  $ \, x \in \gerg \, $.  Then the relations  $ \, x \, y -
y \, x = [x,y] \, $  and  $ \, z^p = z^{[p\hskip0,7pt]} \, $  in
$ \u(\gerg) $  yield  $ \, \tilde{x} \, \tilde{y} - \tilde{y} \,
\tilde{x} = \h \, \widetilde{[x,y]} \equiv 0 \mod \h \, {\u_\h(\gerg)}'
\, $  and  $ \, \tilde{z}^p = \h^{p-1} \widetilde{z^{[p\hskip0,7pt]}}
\equiv 0 \mod \h \, {\u_\h(\gerg)}' \, $;  therefore from the
           presentation\footnote{Hereafter,  $ \, T_A(M) \, $,
resp.~$ \, S_A(M) \, $,  \, is the tensor, resp.~symmetric
algebra of an  $ A $--module  $ M $.}
  $ \; \u_\h(\gerg) = T_R(\gerg) \Big/ \big( \big\{\,
x \, y - y \, x - [x,y] \, , \, z^p - z^{[p\hskip0,7pt]} \;\big\vert\;
x, y, z \in \gerg \,\big\} \big) \; $  we get
 \vskip-10pt
  $$  \displaylines{
   {\u_\h(\gerg)}' = \langle \tilde{\gerg} \rangle
\,\;{\buildrel {\h \rightarrow 0} \over
{\relbar\joinrel\relbar\joinrel\llongrightarrow}}\;\,
\widetilde{\u(\gerg)} \, = \, T_\Bbbk(\tilde{\gerg})
\bigg/ \Big( \Big\{\, \tilde{x} \, \tilde{y} -
\tilde{y} \, \tilde{x} \, , \, \tilde{z}^p \;\Big\vert\;
\tilde{x}, \tilde{y}, \tilde{z} \in \tilde{\gerg} \,\Big\} \Big)
=   \hfill  \cr
   = T_\Bbbk (\gerg) \! \Big/ \! \big( \big\{\, x \, y - y \,
x \, , \, z^p \;\big\vert\; x, y, z \in \gerg \,\big\} \big)
= S_\Bbbk (\gerg) \! \Big/ \! \big( \big\{\, z^p \,\big\vert\;
z \in \gerg \,\big\} \big) = F[\gerg^\star] \! \Big/ \! \big(
\big\{\, z^p \,\big\vert\; z \in \gerg \,\big\} \big)  \cr }  $$
that is  $ \; \widetilde{\u(\gerg)} := G_{\underline{D}} \big(
\u(\gerg) \big) = {\u_\h(\gerg)}' \big/ \h \, {\u_\h(\gerg)}' \cong
F[\gerg^\star] \Big/ \big( \big\{\, z^p \,\big\vert\; z \in \gerg
\,\big\} \big) \; $  {\sl as Poisson Hopf algebras}.  In particular,
{\it this means that  $ \widetilde{\u(\gerg)} $  is the function
algebra of, and  $ {\u_\h(\gerg)}' $  is a QFA (at $ \h \, $)  for,
a non-reduced algebraic Poisson group of dimension 0 and height 1,
whose cotangent Lie bialgebra is  $ \gerg \, $,  hence which is dual
to  $ \, \gerg \, $};  \, thus, in a sense, part  {\it (c)\/}  of
Theorem 2.2 is still valid in this (positive characteristic) case.

\vskip5pt

   {\sl $ \underline{\hbox{\it Remark}} $:}  \, Note that this last
result reminds the classical formulation of the analogue of Lie's
Third Theorem in the context of group-schemes:  {\it Given a
restricted Lie algebra  $ \gerg $,  there exists a group-scheme
$ G $  {\sl of dimension 0 and height 1\/}  whose tangent Lie algebra
is  $ \, \gerg \, $}  (see e.g.~[DG]).  Here we have just given sort
of a ``dual Poisson-theoretic version'' of this fact, in that our
result sounds as follows:  {\it Given a restricted Lie algebra
$ \gerg $,  there exists a Poisson group-scheme  $ G $  {\sl of
dimension 0 and height 1\/}  whose  {\sl cotangent}  Lie algebra
is  $ \, \gerg \, $}.

\vskip5pt

   As a byproduct, since  $ \, U_\h(\gerg) = \u_\h\big( \gerg^{{[p
\hskip0,7pt]}^\infty} \big) \, $  we have also  $ \; {U_\h(\gerg)}' =
{\u_\h\big( \gerg^{{[p\hskip0,7pt]}^\infty} \big)}' \, $,  \, whence
  $$  {U_\h(\gerg)}' = {\u_\h\big( \gerg^{{[p\hskip0,7pt]}^\infty}
\big)}' \,{\buildrel {\h \rightarrow 0} \over
{\relbar\joinrel\relbar\joinrel\loongrightarrow}}\; S_\Bbbk \Big(
\gerg^{{[p\hskip0,7pt]}^\infty} \Big) \hskip-2pt \bigg/ \hskip-2pt \Big(
{\big\{\, z^p \,\big\}}_{z \in \gerg^{{[p\hskip0,7pt]}^\infty}} \Big)
\, = \, F \! \left[ {\big( \gerg^{{[p\hskip0,7pt]}^\infty} \big)}^{\!
\star} \right] \hskip-2pt \bigg/ \hskip-2pt \Big( {\big\{\, z^p \,
\big\}}_{z \in \gerg^{{[p\hskip0,7pt]}^\infty}} \Big) \, .  $$
   \indent   Furthermore,  $ \, {\u_\h(\gerg)}' = \langle \tilde{\gerg}
\rangle \, $  implies that  $ \, I_{{\u_\h(\gerg)}'} \, $  is generated
(as a two-sided ideal) by  $ \, \h \, R \cdot 1_{\u_\h(\gerg)} + R \,
\tilde{\gerg} \, $,  \, hence  $ \, \h^{-1} I_{{\u_\h(\gerg)}'} \, $
is generated by  $ \, R \cdot 1 + R \, \gerg \, $,  \, thus  $ \,
{\big( {\u_\h(\gerg)}' \big)}^\vee := \bigcup\nolimits_{n \geq 0}
{\big( \h^{-1} I_{{\u_\h(\gerg)}'} \big)}^n = \bigcup\nolimits_{n
\geq 0} {\big( R \cdot 1 + R \, \gerg \big)}^n = \u_\h(\gerg) \, $;
\, this means that also part  {\it (b)}  of Theorem 2.2 is still
valid, though now  $ \, \text{\it Char}\,(\Bbbk) > 0 \, $.
                                          \par
   When  $ \, \text{\it Char}\,(\Bbbk) = 0 \, $  and we look at
$ U(\gerg) $,  the like argument applies:  $ \underline{D} $
coincides with the standard filtration of  $ U(\gerg) $  given
by the total degree, via the PBW theorem.  This and (5.2) immediately
imply  $ \, {U(\gerg)}' = \langle \tilde{\gerg} \rangle =
\langle \h \, \gerg \rangle \, $,  so that from the presentation
$ \; U_\h(\gerg) = T_R(\gerg) \Big/ \! \big( {\big\{\, x \, y -
y \, x - [x,y] \,\big\}}_{x, y, z \in \gerg } \big) \, $  we get
$ \, {U_\h(\gerg)}' \! = T_{\!R}(\tilde{\gerg}) \Big/ \! \Big( {\Big\{
\tilde{x} \, \tilde{y} - \tilde{y} \, \tilde{x} - \h \cdot
\widetilde{[x,y]} \Big\}}_{\tilde{x}, \tilde{y} \in \tilde{\gerg}}
\Big) $,  \, whence we get at once
  $$  \displaylines{
   {U_\h(\gerg)}' = T_{\!R}(\tilde{\gerg})
\Big/ \Big( \Big\{ \tilde{x} \, \tilde{y} - \tilde{y} \, \tilde{x}
- \h \cdot \widetilde{[x,y]} \;\Big\vert\; \tilde{x}, \tilde{y} \in
\tilde{\gerg} \Big\} \Big) \,{\buildrel {\h \rightarrow 0} \over
\llongrightarrow}\, T_\Bbbk(\tilde{\gerg}) \Big/ \big( \big\{\,
\tilde{x} \, \tilde{y} - \tilde{y} \, \tilde{x} \;\big\vert\;
\tilde{x}, \tilde{y} \in \tilde{\gerg} \,\big\} \big) =  \cr
   {} \hfill   = T_\Bbbk (\gerg) \Big/ \big( \big\{\, x \, y -
y \, x \;\big\vert\; x, y \in \gerg \,\big\} \big) = S_\Bbbk(\gerg)
= F[\gerg^\star]  \cr }  $$
i.e.~$ \; \widetilde{U(\gerg)} := G_{\underline{D}} \big( U(\gerg)
\big) = {U_\h(\gerg)}' \big/ \h \, {U_\h(\gerg)}' \cong F[\gerg^\star]
\; $  {\sl as Poisson Hopf algebras},  as predicted by  Theorem
2.2{\it (c)}.  Moreover,
      \hbox{$ \, {U_\h(\gerg)}' = \langle \tilde{\gerg}
\rangle = T(\tilde{\gerg}) \Big/ \big( \big\{\, \tilde{x} \, \tilde{y}
- \tilde{y} \, \tilde{x} = \h \cdot \widetilde{[x,y]} \;\big\vert\;
\tilde{x}, \tilde{y} \in \tilde{\gerg} \,\big\} \big) \, $}
  implies
that  $ \, I_{{U_\h(\gerg)}'} \, $  is generated by  $ \, \h \, R
\cdot 1_{U_\h(\gerg)} + R \, \tilde{\gerg} \, $:  \, therefore  $ \,
\h^{-1} I_{{U_\h(\gerg)}'} \, $  is generated by  $ \, R \cdot
1_{U_\h(\gerg)} + R \, \gerg \, $,  \, whence  $ \, {\big( {U_\h
(\gerg)}' \big)}^\vee := \bigcup\nolimits_{n \geq 0} {\big( \h^{-1}
I_{{U_\h(\gerg)}'} \big)}^n = \bigcup\nolimits_{n \geq 0} {\big(
R \cdot 1_{U_\h(\gerg)} + R \, \gerg \big)}^n = U_\h(\gerg) \, $,
\, which agrees with  Theorem 2.2{\it (b)}.
                                          \par
   What for the functor  $ \, {(\ )}^\vee \, $?  This heavily
depends on the  $ \gerg $  we start from!
                                          \par
   First assume  $ \, \hbox{\it Char}\,(\Bbbk) = 0 \, $.  Let
$ \, \gerg_{(1)} := \gerg \, $,  $ \, \gerg_{(k)} := \big[ \gerg,
\gerg_{(k-1)} \big] \, $  ($ k \in \N_+ $),  be the  {\sl lower
central series\/}  of  $ \gerg \, $.  Pick subsets  $ \, B_1 \, $,
$ B_2 \, $,  $ \dots \, $,  $ B_k \, $,  $ \dots \, $  ($ \subseteq
\gerg \, $)  such that  $ \, B_k \! \mod \gerg_{(k+1)} \, $  be a
$ \Bbbk $--basis  of  $ \, \gerg_{(k)} \big/ \gerg_{(k+1)} \, $
(for all  $ \, k \in \N_+ \, $),  \, pick also a  $ \Bbbk $--basis
$ B_{\infty} $  of  $ \, \gerg_{(\infty)} := \bigcap_{\, k \in \N_+}
\, $,  \, and set  $ \, \partial(b) := k \, $  for any  $ \, b \in B_k
\, $  and each  $ \, k \in \N_+ \cup \{\infty\} \, $.  Then  $ \, B :=
\left( \bigcup_{k \in \N_+} B_k \right) \cup B_{\infty} \, $  is a
$ \Bbbk $--basis  of  $ \gerg \, $;  \, we fix a total order on it.
Applying the PBW theorem to this ordered basis of  $ \gerg $  we
get that  $ \, J^n \, $  has basis the set of ordered monomials
$ \, \big\{\, b_1^{e_1} b_2^{e_2} \cdots b_s^{e_s} \;\big|\; s \in
\N_+ \, , b_r \in B \, , \, \sum_{r=1}^s b_r \, \partial(b_r)
\geq n \,\big\} \, $.  Then one easily finds that  $ \,
{U_\h(\gerg)}^\vee \, $  is generated by  $ \, \big\{\,
\h^{-1} b \;\big|\; b \in B_1 \setminus B_2 \,\big\} \, $
(as a unital  $ R $--algebra)  and it is the direct sum
 \vskip2pt
   \centerline{ $  {U_\h(\gerg)}^\vee = \, \bigg( \hskip-3pt \oplus_{\Sb
\hskip-15pt  s \in \N_+  \\   b_r \in B \setminus B_\infty  \endSb}
\hskip-15pt  R \, \big( \h^{-\partial(b_1)} b_1 \big)^{e_1} \cdots
\big( \h^{-\partial(b_s)} b_s \big)^{e_s} \hskip-2pt \bigg) \bigoplus
\bigg( \hskip-3pt \oplus_{\Sb  s \in \N_+, \, b_r \in B  \\   \exists
\, \bar{r} : b_{\bar{r}} \in B_\infty \endSb} R\big[ \h^{-1} \big] \,
b_1^{e_1} \cdots b_s^{e_s} \hskip-2pt \bigg) $ }
From this it follows at once that  $ \; {U_\h(\gerg)}^\vee \Big/ \h \,
{U_\h(\gerg)}^\vee \cong \, U\left( \gerg \big/ \gerg_{(\infty)} \right)
\; $  via an isomorphism which maps  $ \, \h^{-\partial(b)} b \mod
\h \, {U_\h(\gerg)}^\vee \, $  to  $ \; b \! \mod \gerg_{(\infty)}
\in \gerg \big/ \gerg_{(\infty)} \subset U \left( \gerg \big/
\gerg_{(\infty)} \right) \, $  for all  $ \, b \in B \setminus
B_{\infty} \, $  and maps  $ \, \h^{-n} b \! \mod \h \,
{U_\h(\gerg)}^\vee \, $  to  $ \, 0 \, $  for all  $ \, b
\in B \setminus B_{\infty} \, $  and all  $ \, n \in \N \, $.
                                          \par
   Now assume  $ \, \hbox{\it Char}\,(\Bbbk) = p > 0 \, $.  Then
in addition to the previous considerations one has to take into
account the filtration of  $ \, \u(\gerg) \, $  induced by both
the lower central series of  $ \gerg $  {\sl and\/}  the
$ p $--filtration  of  $ \gerg $,  that is  $ \, \gerg \supseteq
\gerg^{[p\hskip0,7pt]} \supseteq \gerg^{{[p\hskip0,7pt]}^2} \supseteq
\cdots \supseteq \gerg^{{[p\hskip0,7pt]}^n} \supseteq \cdots \, $,
\, where  $ \, \gerg^{{[p\hskip0,7pt]}^n} \, $  is the restricted
Lie subalgebra generated by  $ \, \big\{\, x^{{[p\hskip0,7pt]}^n}
\,\big|\, x \in \gerg \,\big\} \, $  and  $ \, x \mapsto
x^{[p\hskip0,5pt]} \, $  is the  $ p \, $--operation  in
$ \gerg \, $:  these encode the  $ J $--filtration  of
$ U(\gerg) $,  hence of  $ \, H = \u_\h(\gerg) \, $,
       \hbox{so permit to describe  $ H^\vee $.}
                                          \par
   In detail, for any restricted Lie algebra  $ \gerh $,
let  $ \, \gerh_n := \Big\langle \bigcup_{(m \, p^k \geq n}
{(\gerh_{(m)})}^{[p^k]} \Big\rangle \, $  for all  $ \, n
\in \N_+ \, $ (where  $ \langle X \rangle $  denotes the Lie
subalgebra of  $ \gerh $  generated by  $ X \, $)  and  $ \,
\gerh_\infty := \bigcap_{\,n \in \N_+} \gerh_n \, $:  \, we
call  $ \, {\big\{ \gerh_n \big\}}_{n \in \N_+} \, $  {\sl the
$ p $--lower  central series of\/}  $ \gerh \, $.  It is a  {\sl
strongly central series\/}  of  $ \gerh $,  i.e.~it is a central
series  ($ = $  decreasing filtration of ideals, each one centralizing
the previous one) of  $ \gerh $  such that  $ \, [\gerh_m,\gerh_n] \leq
\gerh_{m+n} \, $  for all  $ m $,  $ n \, $;  in addition, it verifies
$ \, \gerh_n^{\,[p\,]} \leq \gerh_{n+1} \, $.  When  $ \gerh $  is
Abelian  $ \, {\big\{ \gerh_n \big\}}_{n \in \N_+} \, $  coincides
(after index rescaling) with the  {\sl $ p $--power  series\/}
$ \, {\Big\{ \gerh^{[p^n]} \Big\}}_{n \in \N} \, $.
                                          \par
   Applying these tools to  $ \, \gerg \subseteq \u(\gerg) \, $
the very definitions give  $ \, \gerg_n \subseteq J^n \, $  (for
all  $ \, n \in \N \, $)  where  $ \, J := J_{\u(\gerg)} \, $:
more precisely, if  $ B $  is an ordered basis of  $ \gerg $
then the (restricted) PBW theorem for  $ \, \u(\gerg) \, $
implies that  $ \, J^n \big/ J^{n+1} \, $  admits as
$ \Bbbk $--basis  the set of ordered monomials of the
form  $ \, x_{i_1}^{e_1} x_{i_2}^{e_2} \cdots x_{i_s}^{e_s}
\, $  such that  $ \, \sum_{r=1}^s e_r \partial(x_{i_r}) = n \, $
where  $ \, \partial(x_{i_r}) \in \N \, $  is uniquely determined
by the condition  $ \, x_{i_r} \in \gerg_{\partial(x_{i_r})}
\setminus \gerg_{\partial(x_{i_r})+1} \, $  and each  $ \,
x_{i_k} \, $  is a fixed lift in  $ \gerg $  of an element of
a fixed ordered basis of  $ \, \gerg_{\partial(x_{i_k})} \Big/
\gerg_{\partial(x_{i_k})+1} \, $.  This yields an explicit
description of  $ \underline{J} \, $,  hence of  $ {\u(\gerg)}^\vee $
and  $ {\u_\h(\gerg)}^\vee $,  like before: in particular  $ \;
{\u_\h(\gerg)}^\vee \Big/ \h \, {\u_\h(\gerg)}^\vee \cong \,
\u\left( \gerg \big/ \gerg_{\infty} \right) \, $.

\vskip7pt

\proclaim{Definition 5.8} \, For any\/  $ \Bbbk $--coalgebra  $ C $,
define  $ \, X \bigwedge Y := \Delta^{-1} \big( X \otimes C + C
\otimes Y \big) \, $  for all subspaces  $ X $,  $ Y $  of  $ C \, $.
Set also  $ \, \bigwedge^1 \! X := X \, $  and  $ \, \bigwedge^{n+1}
X := \big( \bigwedge^n X \big) \bigwedge X  \, $  for all  $ \,
n \in \N_+ \, $,  and also  $ \, \bigwedge^0 \! X := \Bbbk
\cdot 1 \, $  if  $ \, C $  is a  $ \Bbbk $--bialgebra.
\endproclaim

\vskip5pt

\proclaim{Lemma 5.9} \, Let  $ H $  be a Hopf  $ \, \Bbbk $--algebra.
Then  $ \; D_n = \bigwedge^{n+1} \! (\Bbbk \cdot 1) \; $  for all
$ \, n \in \N \, $.
\endproclaim

\demo{Proof}  Definitions give  $ \; D_0 := \text{\sl Ker}\,(\delta_1)
= \Bbbk \cdot 1 = \bigwedge^{\!1} \! (\Bbbk \cdot 1) \, $.  By
coassociativity we have  $ \, D_n := \text{\sl Ker}\,(\delta_{n+1})
= \text{\sl Ker}\,\big( (\delta_n \otimes \delta_1) \circ \delta_2
\big) = \text{\sl Ker}\,\big( (\delta_n \otimes \delta_1) \circ
\Delta \big) = \Delta^{-1} \big( \text{\sl Ker}\,(\delta_n \otimes
\delta_1) \big) = \Delta^{-1} \big( \text{\sl Ker}\,(\delta_n)
\otimes H + H \otimes \text{\sl Ker}\,(\delta_1) \big) = \Delta^{-1}
\big( D_{n-1} \otimes H + H \otimes D_0 \big) = D_{n-1} \! \bigwedge
D_0 = D_{n-1} \! \bigwedge (\Bbbk \cdot 1) \, $  for all  $ \, n \!
\in \! \N_+ \, $;  \,
           \hbox{so by induction
$ \, D_n = D_{n-1} \! \bigwedge (\Bbbk \! \cdot \! 1) =
\big( \bigwedge^{\!n} (\Bbbk \! \cdot \! 1) \big) \! \bigwedge
(\Bbbk \! \cdot \! 1) = \bigwedge^{\! n+1} (\Bbbk \! \cdot \! 1)
\, . \; \square $ }
\enddemo

\vskip5pt

\proclaim{Definition 5.10}
                                          \hfill\break
  \indent  (a) \, We call  {\sl pre-restricted universal enveloping
algebra\/}  (in short,  {\sl PrUEA\/})  any  $ \, H \in \HA_\Bbbk
\, $  which is down-filtered by  $ \underline{J} $  (that is,  $ \,
\bigcap_{n \in \N} J^n = \{0\} \, $).  We call  $ \, \PrUEA \, $
the full subcategory of  $ \HA_\Bbbk $  of all the PrUEAs.
                                         \hfill\break
  \indent  (b) \, We call  {\sl pre-function algebra\/}  (in short,
{\sl PFA\/})  any  $ \, H \in \HA_\Bbbk \, $  which is up-filtered
by  $ \underline{D} $  (that is,  $ \, \bigcup_{n \in \N} D_n =
H \, $).  We call  $ \, \PFA \, $  the full subcategory of
$ \HA_\Bbbk $  of all the PFAs.
\endproclaim

\vskip6pt

   The content of the notions of PrUEA and of PFA is revealed
by parts  {\it (a)\/}  and  {\it (b)\/}  of next theorem, which
collects the main results of this section.

\vskip7pt

\proclaim{Theorem 5.11} \, ({\sl ``The Crystal Duality Principle''})
                                       \hfill\break
  \indent   (a) \, The assignment  $ \, H \mapsto H^\vee := H \big/
J_{{}_H}^{\, \infty} \, $,  \, resp.~$ \, H \mapsto H' := \bigcup_{n
\in \N} D_n \, $,  \, defines a functor  $ \; {(\ )}^\vee \colon \,
\HA_\Bbbk \longrightarrow \HA_\Bbbk \, $,  \, resp.~$ \; {(\ )}'
\colon \, \HA_\Bbbk \longrightarrow \HA_\Bbbk \, $,  \, whose image
is  $ \PrUEA $,  \, resp.~$ \PFA $.  More in general, the assignment
$ \, A \mapsto A^\vee := A \big/ J_{{}_A}^{\, \infty} \, $,  \,
resp.~$ \, C \mapsto C' := \bigcup_{n \in \N} D_n(C) \, $,  \,
defines a surjective functor from augmented\/  $ \Bbbk $--algebras,
resp.~coaugmented  $ \Bbbk $--coalgebras,  to augmented\/
$ \Bbbk $--algebras  which are down-filtered by  $ \underline{J} $,
resp.~coaugmented\/  $ \Bbbk $--coalgebras  which are up-filtered by
$ \underline{D} \, $;  \, and similarly for  $ \Bbbk $--bialgebras.
                                          \hfill\break
  \indent   (b) \, Let  $ \, H \in \HA_\Bbbk \, $.  Then  $ \,
\widehat{H} := G_{\underline{J}}(H) \cong \U(\gerg) \, $,  as graded
co-Poisson Hopf algebras, for some restricted Lie bialgebra  $ \gerg $
which is graded as a Lie algebra.  In particular, if  $ \, \Char(\Bbbk)
\!=\! 0 \, $  and  $ \, \dim(H) \!\in \N \, $  then  $ \, \widehat{H}
= \Bbbk \!\cdot\! 1 \, $  and  $ \, \gerg = \{0\} \, $.
                                        \hfill\break
   \indent   More in general, the same holds if  $ \, H = B \, $
is a\/  $ \Bbbk $--bialgebra.
                                        \hfill\break
  \indent   (c) \, Let  $ \, H \in \HA_\Bbbk \, $.  Then  $ \,
\widetilde{H} := G_{\underline{D}}(H) \cong F[G] \, $,  as graded
Poisson Hopf algebras, for some connected algebraic Poisson group
$ G $  whose variety of closed points form a (pro)affine space.  If
$ \, \Char(\Bbbk) = 0 \, $  then  $ \, F[G] = \widetilde{H} \, $  is
a polynomial algebra, i.e.~$ \, F[G] = \Bbbk \big[ {\{x_i \}}_{i \in
\Cal{I}} \big] \, $  (for some set  $ \Cal{I} $);  in particular, if
$ \, \dim(H) \in \N \, $  then  $ \, \widetilde{H} = \Bbbk \cdot 1 \, $
and  $ \, G = \{1\} \, $.  If  $ \, p := \Char(\Bbbk) > 0 \, $  then
$ G $  has dimension 0 and height 1, and if\/  $ \Bbbk $  is perfect
then  $ \, F[G] = \widetilde{H} \, $  is a truncated polynomial
algebra, i.e.~$ \, F[G] = \Bbbk \big[ {\{x_i\}}_{i \in \Cal{I}}
\big] \Big/ \big( {\{x_i^{\,p}\}}_{i \in \Cal{I}} \big) \, $
(for some set  $ \Cal{I} $).
                                          \hfill\break
   \indent   More in general, the same holds if  $ \, H = B \, $
is a\/  $ \Bbbk $--bialgebra.
                                          \hfill\break
   \indent   (d) \, For every  $ \, H \in \HA_\Bbbk \, $,  \, there
exist two 1-parameter families  $ \, {{(H^\vee)}_\h}^{\!\!\vee}
= \Cal{R}^\h_{\underline{J}}(H^\vee) \, $  and  $ \, {\big(
{{(H^\vee)}_\h}^{\!\!\vee} \big)}' \, $  in  $ \, \HA_\Bbbk \, $
giving deformations of  $ H^\vee $  with regular fibers
 \vskip-9pt
  $$  \left.  \hbox{$ \matrix
 \text{\ if \ }  \!\text{\it Char}\,(\Bbbk) = 0 \, ,
\quad\;  U(\gerg_-)  \\
 \text{\ if \ }  \!\text{\it Char}\,(\Bbbk) > 0 \, ,
\quad\;  \u(\gerg_-)
                      \endmatrix $}  \right\}
 = \widehat{H}  \hskip3pt
\underset{{{(H^\vee)}_\h}^{\!\!\vee}}  \to
{\overset{0 \,\leftarrow\, \h \,\rightarrow\, 1}
\to{\longleftarrow\joinrel\relbar\joinrel%
\relbar\joinrel\relbar\joinrel\llongrightarrow}}
\hskip2pt  H^\vee  \hskip0pt
\underset{\;{({{(H^\vee)}_\h}^{\!\!\vee})}'}  \to
{\overset{1 \,\leftarrow\, \h \,\rightarrow\, 0}
\to{\longleftarrow\joinrel\relbar\joinrel%
\relbar\joinrel\relbar\joinrel\llongrightarrow}}
\hskip3pt  \hbox{$ \cases
         F[K_-] = F[G_-^\star]  \\
         F[K_-]
                   \endcases $}  $$
 \vskip-3pt
 \noindent
and two 1-parameter families  $ \, {H_\h}^{\!\prime}
= \Cal{R}^\h_{\underline{D}}(H') \, $  and  $ \,
{({H_\h}^{\!\prime})}^\vee \, $  in  $ \, \HA_\Bbbk \, $
giving deformations of  $ H' $  with regular fibers
 \vskip-12pt
  $$  F[G_+] = \widetilde{H}  \hskip3pt
\underset{\;{H_\h}^{\!\prime}}  \to
{\overset{0 \,\leftarrow\, \h \,\rightarrow\, 1}
\to{\longleftarrow\joinrel\relbar\joinrel%
\relbar\joinrel\relbar\joinrel\llongrightarrow}}
\hskip2pt  H'  \hskip0pt
\underset{\;{({H_\h}^{\!\prime})}^\vee}  \to
{\overset{1 \,\leftarrow\, \h \,\rightarrow\, 0}
\to{\longleftarrow\joinrel\relbar\joinrel%
\relbar\joinrel\relbar\joinrel\llongrightarrow}}
\hskip3pt  \hbox{$ \cases
      U(\gerk_+) = U(\gerg_+^{\,\times})
&   \quad  \text{\ if \ }  \text{\it Char}\,(\Bbbk) = 0  \\
      \u(\gerk_+)
&   \quad  \text{\ if \ }  \text{\it Char}\,(\Bbbk) > 0  \\
                   \endcases $}  $$
where  $ G_+ $  is like  $ G $  in (c),  $ K_- $  is a connected
algebraic Poisson group,  $ \gerg_- $  is like  $ \gerg $  in (b),
$ \gerk_+ $  is a (restricted, if  $ \, \text{\it Char}\,(\Bbbk) >
0 \, $)  Lie bialgebra,  $ \gerg_+^{\,\times} $  is the cotangent
Lie bialgebra to  $ G_+ $  and  $ G_-^\star $  is a connected
algebraic Poisson group whose cotangent Lie bialgebra is
$ \gerg_- \, $.
                                        \hfill\break
   \indent   (e) \, If  $ \, H = F[G] \, $  is the function algebra
of an algebraic Poisson group  $ G $,  then  $ \widehat{F[G]} $  is
a bi-Poisson Hopf algebra (see [KT], \S 1), namely
  $$  \widehat{F[G]} \;\, \cong \,\; S(\gerg^\times) \bigg/ \!
\Big( \Big\{\, \overline{x}^{\,p^{n(x)}} \,\Big|\, x \in
\Cal{N}_{F[G]} \Big\} \Big) \;\, \cong \,\;
U(\gerg^\times) \bigg/ \! \Big( \Big\{\, \overline{x}^{\,p^{n(x)}}
\,\Big|\, x \in \Cal{N}_{F[G]} \Big\} \Big)  $$
where  $ \Cal{N}_{F[G]} $  is the nilradical of  $ F[G] $,
$ \, p^{n(x)} \, $  is the order of nilpotency of  $ \, x \in
\Cal{N}_{F[G]} $  and the bi-Poisson Hopf structure of  $ \;
S(\gerg^\times) \bigg/ \! \Big( \Big\{\, \overline{x}^{\,p^{n(x)}}
\,\Big|\, x \in \Cal{N}_{F[G]} \Big\} \Big) \; $  is the quotient
one from  $ \, S(\gerg^\times) \, $;  \, in particular, if the
group  $ G $  is reduced then  $ \; \widehat{F[G]} \, \cong \,
S(\gerg^\times) \, \cong \, U(\gerg^\times) \; $.
                                        \hfill\break
   \indent   (f)  If  $ \, \text{Char}\,(\Bbbk) = 0 \, $  and  $ \,
H = U(\gerg) \, $  is the universal enveloping algebra of some Lie
bialgebra  $ \gerg $,  then  $ \, \widetilde{U(\gerg)} \, $  is
a bi-Poisson Hopf algebra, namely
  $$  \widetilde{U(\gerg)} \;\, \cong \,\; S(\gerg)
\;\, = \,\; F[\gerg^\star] \;  $$
where the bi-Poisson Hopf structure on  $ S(\gerg) $  is the
canonical one.
                                      \hfill\break
   \indent   If  $ \, \text{Char}\,(\Bbbk) = p > 0 \, $  and  $ \, H =
\u(\gerg) \, $  is the restricted universal enveloping algebra of some
restricted Lie bialgebra  $ \gerg $,  then  $ \, \widetilde{\u(\gerg)}
\, $  is a bi-Poisson Hopf algebra, namely
  $$  \widetilde{\u(\gerg)} \;\, \cong \,\; S(\gerg) \Big/
\big( \big\{ x^p \,\big|\, x \in \gerg \big\} \big) \;\,
= \,\; F[G^\star] \;  $$
where the bi-Poisson Hopf structure on  $ \, S(\gerg) \Big/
\big( \big\{ x^p \,\big|\, x \in \gerg \big\} \big) \, $  is induced
by the canonical one on  $ S(\gerg) $,  and  $ G^\star $  is a
connected algebraic Poisson group of dimension 0 and height 1
whose cotangent Lie bialgebra is  $ \gerg \, $.
                                      \hfill\break
  \indent   (g) \, Let  $ \, H $,  $ K \in \HA_\Bbbk \, $  and let
$ \, \pi \colon \, H \times K \loongrightarrow \Bbbk \, $  be a Hopf
pairing.  Then  $ \pi $  induce a  {\sl filtered}  Hopf pairing  $ \,
\pi_f \, \colon \, H^\vee \times K' \loongrightarrow \Bbbk \, $,  \,
a  {\sl graded}  Hopf pairing  $ \, \pi_{{}_G} \, \colon \, \widehat{H}
\times \widetilde{K} \loongrightarrow \Bbbk \, $,  both perfect on the
right, and Hopf pairings over  $ \Bbbk[\h\,] $  (notation of\/ \S 5.1)
$ \, H_\h \times K_\h \loongrightarrow \Bbbk[\h\,] \, $  and  $ \,
{H_\h}^{\!\vee} \times {K_\h}^{\!\prime} \loongrightarrow \Bbbk[\h\,]
\, $,  \, the latter being perfect on the right.  If in addition the
pairing  $ \, \pi_f \, \colon \, H^\vee \times K' \loongrightarrow
\Bbbk \, $  is perfect, then  $ \pi_{{}_G} $  is perfect as well,
and  $ {H_\h}^{\!\vee} $  and  $ {K_\h}^{\!\prime} $  are dual to
each other.  The left-right symmetrical results hold too.
\endproclaim

\demo{Proof} \, Parts  {\it (a)\/}  through  {\it (c)\/}  of the
statement are proved by the analysis in \S 5.4, but for the naturality
of  $ \, H \mapsto H^\vee \, $  and  $ \, H \mapsto H' \, $,  \, which
is however clear because,  $ \, \varphi\big(J_{{}_H}^{\, \infty}\big)
\subseteq J_{{}_K}^{\, \infty} \, $  and  $ \, \varphi\big(D_n(H)\big)
\subseteq D_n(K) \, $  for any morphism  $ \, \varphi \, \colon \,
H \longrightarrow K \, $  within  $ \HA_\Bbbk \, $.  In addition, for
part  {\it (b)\/}  when  $ \, \Char(\Bbbk) = 0 \, $  and  $ \, \dim(H)
\in \N \, $  we have to notice that  $ \, \widehat{H} = U(\gerg) \, $
is finite dimensional too, hence  $ \, \widehat{H} = U(\gerg) = \Bbbk
\cdot 1 \, $  and  $ \, \gerg = \{0\} \, $;  \, similarly for  {\it
(c)\/}  these assumptions imply that  $ \, \widetilde{H} = F[G] \, $
is finite dimensional too, so  $ \, \widetilde{H} = F[G] = \Bbbk
\cdot 1 \, $  and  $ \, G = \{1\} \, $,  \, q.e.d.  Finally, if
$ \, H = B \, $  is just a  $ \Bbbk $--bialgebra  then both  $ \,
\widehat{B} := G_{\underline{J}}(B) \, $  and  $ \, \widetilde{B}
:= G_{\underline{D}}(B) \, $  are  {\sl irreducible\/}  graded
$ \Bbbk $--bialgebras:  then by [Ab], Theorem 2.4.24, they are
also graded Hopf algebras, whence we conclude as if  $ B $
were a Hopf algebra.
                                       \par
   Part  {\it (d)\/}  is proved by \S 5.5.
                                       \par
   As for part  {\it (e)},  it is almost entirely proved by the
analysis in \S 5.6, noting also that in the case of  $ \, H = F[G]
\, $  one has  $ \, S(\gerg^\times) = U(\gerg^\times) \, $
because $ \gerg^\times $  is Abelian.  What is left to check is
whatever refers to  {\sl bi-Poisson\/}  structures.  Indeed, the
Lie bracket of  $ \gerg^\times $  extends to a Poisson bracket
which makes $ S(\gerg^\times) $  into a bi-Poisson Hopf algebra
(see \S 5.1); then  $ \; \Big( \Big\{\, \overline{x}^{\,p^{n(x)}}
\Big\}_{x \in \Cal{N}_{F[G]}} \,\Big) \; $  is a bi-Poisson Hopf
ideal, thus $ \; S(\gerg^\times) \bigg/ \! \Big( {\Big\{\,
\overline{x}^{\,p^{n(x)}} \Big\}}_{x \in \Cal{N}_{F[G]}} \,\Big)
\; $  is a bi-Poisson Hopf algebra as well.  But  $ \widehat{F[G]} $
also inherits a Poisson bracket from  $ F[G] $  which makes it
into a bi-Poisson Hopf algebra too: it is then clear that the
isomorphism  $ \; S(\gerg^\times) \bigg/ \! \Big( {\Big\{\,
\overline{x}^{\,p^{n(x)}} \Big\}}_{x \in \Cal{N}_{F[G]}} \,\Big)
\, \cong \, \widehat{F[G]} \; $  is one of {\sl bi-Poisson\/}
Hopf algebras.
                                       \par
   Similarly, part  {\it (f)\/}  is proved by the analysis in \S
5.7, noting also that both  $ \, \widetilde{U(\gerg)} \, $  and
$ \, S(\gerg) \, = \, F[G^\star] \, $  are naturally bi-Poisson
Hopf algebras, isomorphic to each other via the previously
considered isomorphism.  In addition, the same holds also for
$ \, \widetilde{\u(\gerg)} \, $  and  $ \, S(\gerg) \Big/ \big(
\big\{ x^p \,\big|\, x \in \gerg \big\} \big) \, = \, F[G^\star]
\, $,  \, because  $ \, \big( \big\{ x^p \,\big|\, x \in \gerg
\big\} \big) \, $  is a bi-Poisson Hopf ideal of  $ \, S(\gerg)
\, $.
                                       \par
   Finally, we go for part  {\it (g)}.  Let  $ \, \pi \colon \,
H \times K \loongrightarrow \Bbbk \, $  be the Hopf pairing under
study.  Consider the filtrations  $ \, \underline{J} = {\big\{
{J_{\!{}_H}}^{\! n} \big\}}_{n \in \N} \, $  and  $ \, \underline{D}
= {\big\{ D^{\scriptscriptstyle K}_n \big\}}_{n \in \N} \, $.  The
key fact is that
  $$  D^{\scriptscriptstyle K}_n = \big( {J_{\!{}_H}}^{\! n+1}
\big)^\perp  \qquad  \hbox{and}  \qquad  {J_{\!{}_H}}^{\! n+1}
\subseteq \big( D^{\scriptscriptstyle K}_n \big)^\perp  \qquad
\qquad  \hbox{for all} \;\; n \in \N \, .   \eqno (5.6)  $$
   \indent   Indeed, if  $ X $  is a subspace of a coalgebra  $ C $
and  $ C $  is in perfect ``Hopf-like'' pairing with an algebra
$ A \, $,  \, one has  $ \, \bigwedge^n X = {\big( {(X^\perp)}^n
\big)}^\perp \, $  (cf.~Definition 5.8) for all  $ \, n \in \N \, $,
\, where the superscript  $ \perp $  means ``orthogonal subspace''
(either in  $ A $  or in  $ C $)  w.r.t.~the pairing under exam
(cf.~[Ab] or [Mo]).  Now, Lemma 5.9 gives  $ \, D^{\scriptscriptstyle
K}_n = \bigwedge^{n+1} (\Bbbk \!\cdot\! 1_{\!{}_K}) \, $,  \, thus
$ \, D^{\scriptscriptstyle K}_n = \bigwedge^{n+1} (\Bbbk \!\cdot\!
1_{\!{}_K}) = {\Big( {\big( {(\Bbbk \!\cdot\! 1_{\!{}_K})}^\perp
\big)}^{n+1} \Big)}^\perp = {\big( {J_{\!{}_H}}^{\! n+1} \big)}^\perp
\, $  because  $ \, {(\Bbbk \!\cdot\! 1_{\!{}_K})}^\perp =
J_{\!{}_H} \, $  (w.r.t.~the pairing  $ \pi $  above).  Therefore
$ \, D^{\scriptscriptstyle K}_n = \big( {J_{\!{}_H}}^{\! n+1}
\big)^\perp \, $,  \, and this also implies  $ \, {J_{\!
{}_H}}^{\! n+1} \subseteq \big( D^{\scriptscriptstyle K}_n
\big)^\perp $.
                                       \par
   Now  $ \, K' := \bigcup_{n \in \N} D^{\scriptscriptstyle K}_n =
\bigcup_{n \in \N} \left( {J_{\!{}_H}}^{\! n+1} \right)^\perp \! =
\left( \, \bigcap_{n \in \N} {J_{\!{}_H}}^{\! n+1} \right)^\perp \!
= \big( {J_{\!{}_H}}^{\!\infty} \big)^\perp \, $.  Thus  $ \pi $
induces a Hopf pairing  $ \, \pi_f \, \colon \, H^\vee \times K'
\loongrightarrow \Bbbk \, $  as required, and by (5.6) this respects
the filtrations on either side.  Then by general theory  $ \pi_f $
induces a graded Hopf pairing  $ \pi_{{}_G} $  as required: in
particular  $ \pi_{{}_G} $  is well-defined because  $ \,
D^{\scriptscriptstyle K}_n \subseteq \big( {J_{\!{}_H}}^{\! n+1}
\big)^\perp \, $  and  $ \, {J_{\!{}_H}}^{\! n+1} \subseteq \big(
D^{\scriptscriptstyle K}_n \big)^\perp \, $  (for all  $ \, n \in
\N_+ \, $)  by (5.6), and both  $ \pi_f $  and  $ \pi_{{}_G} $
are perfect on the right because all the inclusions  $ \,
D^{\scriptscriptstyle K}_n \subseteq \big( {J_{\!{}_H}}^{\! n+1}
\big)^\perp \, $  happen to be identities.  Clearly by scalar
extension  $ \pi $  defines also a Hopf pairing  $ \, H_\h \times
K_\h \longrightarrow \Bbbk[\h\,] \, $;  \, then (5.6) and the
description of  $ {H_\h}^{\!\prime} $  and  $ {K_\h}^{\!\vee} $
in Lemma 5.2 directly imply that this yields another Hopf pairing
$ \, {H_\h}^{\!\vee} \times {K_\h}^{\!\prime} \longrightarrow
\Bbbk[\h\,] \, $  as claimed.
                                       \par
   Finally when  $ \pi_f $  is perfect it is easy to see that
$ \pi_{{}_G} $  is perfect as well; note that this improves
(5.6), for we have  $ \, {J_{\!{}_H}}^{\! n+1} = \big(
D^{\scriptscriptstyle K}_n \big)^\perp \, $  for all  $ \,
n \in \N \, $.  It is also clear that the pairing  $ \,
{H_\h}^{\!\vee} \times {K_\h}^{\!\prime} \longrightarrow
\Bbbk[\h\,] \, $  is perfect as well, and that  $ {H_\h}^{\!
\vee} $  and  $ {K_\h}^{\!\prime} $  are dual to each
other.   \qed
\enddemo

\vskip4pt

   {\sl $ \underline{\hbox{\it Remarks}} $:}  \;  {\it (a)} \, It is
worth noticing that, though usually introduced in a different way,
$ H' $  is an object which is pretty familiar to Hopf algebra theorists:
indeed, it is the  {\sl connected component\/}  of  $ H \, $ (cf.~[Ga5]
for a proof); in particular,  $ H $  is a PFA if and only if it is
connected.  Nevertheless, surprisingly enough the pretty remarkable
property of its associated graded Hopf algebra  $ \, \widetilde{H} =
G_{\underline{D}}(H) \, $  expressed by  Theorem 5.11{\it (c)\/}  seems
to have been unknown so far (at least, to the author's knowledge)!
Similarly, the ``dual'' construction of  $ H^\vee $  and the important
property of its associated graded Hopf algebra  $ \, \widehat{H}
= G_{\underline{J}}(H) \, $  stated in  Theorem 5.11{\it (b)\/}
seem to have escaped the specialists' attention.
                                    \par
   {\it (b)} \, Part  {\it (d)}  of Theorem 5.11 is quite interesting
for applications in physics.  In fact, let  $ H $  be a Hopf algebra
which describes the symmetries of some physically meaningful system,
but has no geometrical meaning (typically, when it is not commutative
nor cocommutative, as it usually happens in quantum physics), and
assume also  $ \, H' = H = H^\vee \, $.  Then  Theorem 5.11{\it (d)}
yields a recipe to deform  $ H $  to four different Hopf algebras
bearing a geometrical meaning, which means having two Poisson groups
and two Lie bialgebras attached to  $ H $,  hence a rich ``geometrical
symmetry'' (of Poisson type) underlying the physical system; if the
ground field has characteristic zero (as usual) we simply have two
pairs of mutually dual Poisson groups together with their tangent Lie
bialgebras.  In \S 10 we'll give a nice application of this kind
with the two pairs of groups strictly related, yet different.

\vskip7pt

  {\bf 5.12 The hyperalgebra case.} \, Let  $ G $  be an algebraic
group, which for simplicity we assume to be finite-dimensional.  By
$ \hyp(G) $  we mean the hyperalgebra associated to  $ G $,  defined
as  $ \, \hyp(G) := {\big( {F[G]}^\bullet \big)}_\epsilon = \big\{\,
\phi \in {F[G]}^\circ \,\big|\, \phi({\germ_e}^{\!n}) = 0 \, , \, \forall
\; n \gg 0 \,\big\} \, $,  \, that is the irreducible component of the
{\sl dual\/}  Hopf algebra  $ \, {F[G]}^\circ \, $  containing  $ \,
\epsilon = \epsilon_{\!{}_{F[G]}} \, $,  \, which is a Hopf subalgebra
of  $ {F[G]}^\circ $;  in particular,  $ \hyp(G) $  is connected
cocommutative.  Recall that there's a natural Hopf algebra morphism
$ \, \Phi : U(\gerg) \longrightarrow \hyp(G) \, $;  \, if  $ \,
\text{\it Char}\,(\Bbbk) = 0 \, $  then  $ \Phi $  is an isomorphism,
so  $ \hyp(G) $  identifies to  $ U(\gerg) $;  \, if  $ \, \text{\it
Char}\,(\Bbbk) > 0 \, $  then  $ \Phi $  factors through  $ \u(\gerg) $
and the induced morphism  $ \, \overline{\Phi} : \u(\gerg)
\longrightarrow \hyp(G) \, $  is injective, so that  $ \u(\gerg) $
identifies with a Hopf subalgebra of  $ \hyp(G) $.  Now we study
$ {\hyp(G)}' $,  $ {\hyp(G)}^\vee $,  $ \widetilde{\hyp(G)} $,
$ \widehat{\hyp(G)} $.
                                               \par
   As  $ \hyp(G) $  is connected, letting  $ \, C_0 := \text{\it
Corad}\,\big(\hyp(G)\big) \, $  be its coradical we have  $ \;
\hyp(G) = \bigcup_{n \in \N} \bigwedge^{n+1} C_0 = \bigcup_{n \in \N}
\bigwedge^{n+1} (\Bbbk \cdot 1) = \bigcup_{n \in \N} D_{n+1}\big(
\hyp(G)\big) =: {\hyp(G)}' \; $.  Now,  Theorem 5.11{\it (c)\/}
gives  $ \, \widetilde{\hyp(G)} := G_{\underline{D}}\big(\hyp(G)\big)
= F[\varGamma\,] \, $  for some connected algebraic Poisson group
$ \varGamma \, $;  Theorem 5.11{\it (e)}  yields  $ \, \widehat{F[G]}
\, \cong \, S(\gerg^*) \bigg/ \! \Big( \Big\{\, \overline{x}^{\,
p^{n(x)}} \Big\}_{x \in \Cal{N}_{F[G]}} \Big) \, = \, \u \bigg(
P \bigg( S(\gerg^*) \bigg/ \! \Big( \Big\{\, \overline{x}^{\,p^{n(x)}}
\Big\}_{x \in \Cal{N}_{F[G]}} \Big) \bigg) \bigg) \, = \, \u \Big(
\big( \gerg^* \big)^{p^\infty} \Big) \, $,  \, with  $ \, \big(
\gerg^* \big)^{p^\infty} := \text{\sl Span}\, \Big( \Big\{\, x^{p^n}
\,\Big\vert\; x \in \gerg^* \, , n \in \N \,\Big\} \Big) \subseteq
\widehat{F[G]} \, $,  \, and noting that  $ \, \gerg^\times =
\gerg^* \, $.  On the other hand, exactly like for  $ U(\gerg) $
and  $ \u(\gerg) $  respectively in case  $ \, \Char(\Bbbk) = 0 \, $
and  $ \, \Char(\Bbbk) > 0 \, $,  \, the filtration  $ \underline{D} $
of  $ \hyp(G) $  is nothing but the natural filtration given by the
order of differential operators: this implies immediately  $ \;
{\hyp(G)_\h}' := {\big( \Bbbk[\h\,] \otimes_\Bbbk \hyp(G) \big)}'
\, = \, \big\langle \big\{\, \h^n x^{(n)} \,\big|\; x \in \gerg \, ,
n \in \N \,\big\} \big\rangle \, $,  \, where hereafter notation like
$ x^{(n)} $  denotes the  $ n $--th  divided power of  $ \, x \in \gerg
\, $  (recall that  $ \hyp(G) $  is generated as an algebra by all the
$ x^{(n)} $'s,  some of which might be zero).  It is then immediate to
check that the graded Hopf pairing between  $ \, {\hyp(G)_\h}' \Big/
\h \, {\hyp(G)_\h}' = \widetilde{\hyp(G)} = F[\varGamma] \, $  and
$ \, \widehat{F[G]} \, $  given by  Theorem 5.11{\it (g)\/}  is
perfect.  From this we easily argue that the cotangent Lie bialgebra
of  $ \varGamma $  is isomorphic to  $ \, \Big( \big( \gerg^*
\big)^{p^\infty} \Big)^* \, $.
                                               \par
   As for  $ {\hyp(G)}^\vee $  and  $ \widehat{\hyp(G)} $,  the
situation is much like for  $ U(\gerg) $  and  $ \u(\gerg) $,
in that it strongly depends on the algebraic nature of  $ G $
(cf.~\S 5.7).

\vskip7pt

  {\bf 5.13 The CDP on group algebras and their duals.} \, In
this section,  $ G $  is any abstract group.  We divide the
subsequent material in several subsections.

\vskip5pt

\noindent   {\it  $ \underline{\hbox{\sl Group-related algebras}} $.}
\, For any commutative unital ring  $ \A \, $,  by  $ \, \A[G] \, $
we mean the group algebra of  $ G $  over  $ \A \, $  and, when  $ G $
is  {\sl finite},  we denote by  $ \, A_\A(G) := {\A[G]}^* \, $  (the
linear dual of  $ \A[G] \, $)  the function algebra of  $ G $  over
$ \A \, $.  Our purpose is to apply the Crystal Duality Principle to
$ \Bbbk[G] $  and  $ A_\Bbbk(G) \, $  with their standard Hopf algebra
structure: hereafter  $ \Bbbk $  is a field and  $ \, R := \Bbbk[\h]
\, $  as in \S 5.1, and we set  $ \, p := \text{\it Char}\,(\Bbbk) \, $.
                                          \par
   Recall that  $ \, H := \A[G] \, $  admits  $ G $  itself
as a distinguished basis, with Hopf algebra structure given by
$ \, g \cdot_{{}_H} \gamma := g \cdot_{{}_G} \gamma \, $,  $ \,
1_{{}_H} := 1_{{}_G} \, $,  $ \, \Delta(g) := g \otimes g \, $,
$ \, \epsilon(g) := 1 \, $,  $ \, S(g) := g^{-1} \, $,  \, for all
$ \, g, \gamma \in G \, $.  Dually,  $ \, H := A_\A(G) \, $  has
basis  $ \, \big\{ \varphi_g \,\big|\, g \! \in \! G \big\} \, $
dual to the basis  $ G $  of  $ \A[G] \, $,  \, with  $ \, \varphi_g
(\gamma) := \delta_{g,\gamma} \, $  for all  $ \, g, \gamma \in G
\, $;  its Hopf algebra structure is given by  $ \, \varphi_g \cdot
\varphi_\gamma := \delta_{g,\gamma} \varphi_g \, $,  $ \, 1_{{}_H}
:= \sum_{g \in G} \varphi_g \, $,  $ \, \Delta(\varphi_g) :=
\sum_{\gamma \cdot \ell = g} \varphi_\gamma \otimes \varphi_\ell
\, $,  $ \, \epsilon(\varphi_g) := \delta_{g,1_G} \, $,  $ \,
S(\varphi_g) := \varphi_{g^{-1}} \, $,  \, for all  $ \, g, \gamma
\in G \, $.  In particular,  $ \, R[G] = R \otimes_\Bbbk \Bbbk[G]
\, $  and  $ \, A_R[G] = R \otimes_\Bbbk A_\Bbbk[G] \, $.  Our
first result is

\vskip5pt

\noindent   {\it  $ \underline{\hbox{\sl Theorem A}} $:}  $ \,
{{\big(\Bbbk[G]\big)}_\h}^{\!\prime} = R \cdot 1 \, $,  $ \;
{\Bbbk[G]}' = \Bbbk \cdot 1 \; $  {\it and}  $ \; \widetilde{\Bbbk[G]}
= \Bbbk \cdot 1 = F\big[\{*\}\big] \, $.

\demo{Proof} The claim follows easily from the formula  $ \, \delta_n(g)
= {(g \! - \! 1)}^{\otimes n} $,  for  $ \, g \in G $,  $ n \in \N \, $.
\qed
\enddemo

\vskip5pt

\noindent   {\it  $ \underline{\hbox{\sl  $ {R[G]}^\vee $,  $ \,
{\Bbbk[G]}^\vee $,  $ \widehat{\Bbbk[G]} $  \ and the dimension
subgroup problem}} $.} \, In contrast with the triviality result in
Theorem A above, things are more interesting for  $ \, {R[G]}^\vee \!
= {{\big( \Bbbk[G] \big)}_\h}^{\!\!\vee} \, $,  $ \, {\Bbbk[G]}^\vee $
and  $ \, \widehat{\Bbbk[G]} \, $.  Note however that since  $ \,
\Bbbk[G] \, $  is cocommutative the induced Poisson cobracket on
$ \, \widehat{\Bbbk[G]} \, $  is trivial, and the Lie cobracket
of  $ \, \gerk_G := P \Big( \widehat{\Bbbk[G]} \Big) \, $  is
trivial as well.
                                          \par
   Studying  $ {\Bbbk[G]}^\vee $  and  $ \widehat{\Bbbk[G]} $
amounts to study the filtration  $ {\big\{ J^n \big\}}_{n \in \N}
\, $,  \, with  $ J := \text{\sl Ker}\,(\epsilon_{{}_{\Bbbk[G]}}) $,
\, which is a classical topic.  Indeed, for  $ \, n \! \in \! \N \, $
let  $ \, \Cal{D}_n(G) := \big\{\, g \in G \,\big|\, (g \! - \! 1) \in
J^n \,\big\} \, $: \, this is a characteristic subgroup of  $ G $,  called
{\sl the  $ n^{\text{th}} $  dimension subgroup of  $ G \, $}.  All these
form a filtration inside  $ G \, $:  \, characterizing it in terms of
$ G $  is the  {\sl dimension subgroup problem},  which (for group
algebras over fields) is completely solved (see [Pa], Ch.~11, \S 1,
and [HB], and references therein); this also gives a description of
$ \, \big\{ J^n \big\}_{n \in \N_+} \, $.  Thus we find ourselves
within the domain of classical group theory: now we use the results
which solve the dimension subgroup problem to argue a description of
$ {\Bbbk[G]}^\vee $,  $ \widehat{\Bbbk[G]} $  and  $ {R[G]}^\vee $,
and later on we'll get from this a description
      \hbox{of  $ \big( {R[G]}^\vee \big)' $  and its semiclassical
limit too.}
                                          \par
   By construction,  $ J $  has  $ \Bbbk $--basis  $ \, \big\{ \eta_g
\,\big|\; g \! \in \! G \setminus \{1_{{}_G}\} \big\} \, $,  \, where
$ \, \eta_g := (g\!-\!1) \, $.  Then  $ \, {\Bbbk[G]}^\vee \, $  is
generated by  $ \, \big\{\, \eta_g \! \mod J^\infty \,\big|\; g \in
G \setminus \{1_{{}_G}\} \big\} \, $,  \, and  $ \, \widehat{\Bbbk[G]}
\, $  by  $ \, \big\{\, \overline{\,\eta_g} \;\big|\; g \! \in \! G
\setminus \{1_{{}_G}\} \big\} \, $:  \, hereafter  $ \, \overline{x}
:= x \mod J^{n+1} \, $  for all  $ \, x \in J^n \, $,  \, that is
$ \overline{x} $  is the element in  $ \widehat{\Bbbk[G]} $  which
corresponds to  $ \, x \in \Bbbk[G] \, $.  Moreover,  $ \, \overline{g}
= \overline{\,1 + \eta_g} = \overline{1} \, $  for all  $ \, g \in G \, $;
\, also,  $ \; \Delta \big(\overline{\,\eta_g}\big) = \overline{\,\eta_g}
\otimes \overline{g} + 1 \otimes \overline{\,\eta_g} = \overline{\,\eta_g}
\otimes 1 + 1 \otimes \overline{\,\eta_g} \; $:  \; thus  $ \overline{\,
\eta_g} $  is primitive, so  $ \, \big\{\,\overline{\,\eta_g} \;\big|\;
g \! \in \! G \setminus \{1_{{}_G}\} \big\} \, $  generates  $ \,
\gerk_G := P \Big( \widehat{\Bbbk[G]} \Big) \, $.

\vskip5pt

\noindent   {\it  $ \underline{\hbox{\sl The Jennings-Hall theorem}} $.}
\, The description of  $ \Cal{D}_n(G) $  is given by the Jennings-Hall
theorem, which we now recall.  The construction involved strongly depends
on whether  $ \, p := \text{\it Char}\,(\Bbbk) \, $  is zero or not, so
we shall distinguish these two cases.
                                          \par
   First assume  $ \, p = 0 \, $.  Let  $ \, G_{(1)} := G \, $,  $ \,
G_{(k)} := (G,G_{(k-1)}) \, $  ($ k \in \N_+ $),  form the  {\sl lower
central series\/}  of  $ G \, $;  hereafter  $ (X,Y) $  is the
commutator subgroup of  $ G $  generated by the set of commutators
$ \, \big\{ (x,y) := x \, y \, x^{-1} y^{-1} \,\big|\, x \in X, y
\in Y \big\} \, $:  \, this is a  {\sl strongly central series\/}
in  $ G $,  which means a central series  $ \, {\big\{ G_k \big\}}_{k
\in \N_+} \, $  (= decreasing filtration of normal subgroups, each
one centralizing the previous one) of  $ G $  such that  $ \, (G_m,
G_n) \leq G_{m+n} \, $  for all  $ m \, $,  $ n \, $.  Then let  $ \,
\sqrt{G_{(n)}} := \big\{ x \in G \,\big|\, \exists \, s \in \N_+ : x^s
\in G_{(n)} \big\} \, $  for all  $ n \in \N_+ \, $:  these form a
descending series of characteristic subgroups in  $ G $,  such that
each composition factor  $ \, A^G_{(n)} := \sqrt{G_{(n)}} \Big/ \!
\sqrt{G_{{(n+1)}}} \, $  is torsion-free Abelian.  Therefore  $ \,
\L_0(G) := \bigoplus_{n \in \N_+} A^G_{(n)} \, $  is a graded Lie ring,
with Lie bracket  $ \, \big[ \overline{g}, \overline{\ell} \,\big] :=
\overline{(g,\ell\,)} \, $  for all  {\sl homogeneous}  $ \overline{g} $,
$ \overline{\ell} \in \L_0(G) \, $,  \, with obvious notation.  It is
easy to see that the map  $ \; \Bbbk \otimes_\Z \L_0(G) \longrightarrow
\gerk_G \, $,  $ \, \overline{g} \mapsto \overline{\eta_g} \, $,  \;
is an epimorphism  {\sl of graded Lie rings\/}:  \, therefore  {\sl
the Lie algebra  $ \, \gerk_G \, $  is a quotient of  $ \; \Bbbk
\otimes_\Z \L_0(G) \, $};  in fact, the above is an isomorphism,
see below.  We shall use notation  $ \, \partial(g) := n \, $  for
all  $ \, g \in \sqrt{G_{(n)}} \, \setminus \sqrt{G_{(n+1)}} \; $.
                                     \par
   For each $ \, k \in \N_+ \, $  pick in  $ A^G_{(k)} $  a subset
$ \overline{B}_k $  which is a  $ \Bbb{Q} $--basis  of  $ \, \Bbb{Q}
\otimes_\Z A^G_{(k)} \, $;  \, for each  $ \, \overline{b} \in
\overline{B}_k \, $,  \, choose a fixed  $ \, b \in \sqrt{G_{(k)}}
\, $  such that its coset in  $ A^G_{(k)} $  be  $ \overline{b} $,  \,
and denote by  $ \, B_k \, $  the set of all such elements  $ b \, $.
Let  $ \, B := \bigcup_{k \in \N_+} B_k \, $:  \, we call such a set
{\sl t.f.l.c.s.-net\/}  (\,=\,``torsion-free-lower-central-series-net'')
on  $ G $.  Clearly  $ \, B_k = \Big( B \cap \sqrt{G_{(k)}} \,\Big)
\setminus \Big( B \cap \sqrt{G_{(k+1)}} \,\Big) \, $  for all
$ k \, $.  By an  {\sl ordered t.f.l.c.s.-net\/}  is meant a
t.f.l.c.s.-net  $ B $  which is totally ordered in such a way that:
{\it (i)\/}  if  $ \, a \in B_m \, $,  $ \, b \in B_n \, $,  $ \, m
< n \, $,  \, then  $ \, a \preceq b \, $;  \, {\it (ii)\/}  for each
$ k $,  every non-empty subset of  $ B_k $  has a greatest element.
An ordered t.f.l.c.s.-net always exists.
                                     \par
   Now assume instead  $ \, p > 0 \, $.  The situation is similar, but
we must also consider the  $ p $--power  operation in the group  $ G $
and in the restricted Lie algebra  $ \gerk_G \, $.  Starting from the
lower central series  $ \, {\big\{G_{(k)}\big\}}_{k \in \N_+} $,  define
$ \, G_{[n]} := \prod_{k p^\ell \geq n} {(G_{(k)})}^{p^\ell} \; $  for
all  $ \, n \in \N_+ \, $  (hereafter, for any group  $ \varGamma $  we
denote  $ \varGamma^{p^e} $  the subgroup generated by  $ \, \big\{
\gamma^{p^e} \,\big|\, \gamma \! \in \! \varGamma \,\big\} \, $):  \,
this gives another strongly central series  $ \, {\big\{G_{[n]}\big\}}_{n
\in \N_+} $  in  $ G $,  \, with the additional property that  $ \,
{(G_{[n]})}^p \leq G_{[n+1]} \, $  for all  $ n \, $,  \, called
{\sl the  $ p $--lower  central series of\/}  $ G \, $.  Then  $ \,
\Cal{L}_p(G) := \oplus_{n \in \N_+} G_{[n]} \big/ G_{[n+1]} \, $  is
a graded restricted Lie algebra over  $ \, \Z_p := \Z \big/ p \, \Z \, $,
\, with operations  $ \, \overline{g} + \overline{\ell} := \overline{g
\cdot \ell} \, $,  $ \, \big[\overline{g},\overline{\ell}\,\big] :=
\overline{(g,\ell\,)} \, $,  $ \, \overline{g}^{\,[p\,]} := \overline{g^p}
\, $,  \, for all  $ \, g $,  $ \ell \in G \, $  (cf.~[HB], Ch.~VIII,
\S 9).  Like before, we consider the map  $ \; \Bbbk \otimes_{\Z_p}
\Cal{L}_p(G) \longrightarrow \gerk_G \, $,  $ \, \overline{g} \mapsto
\overline{\eta_g} \, $,  \; which now is an epimorphism  {\sl of graded
restricted Lie  $ \Z_p $--algebras},  whose image spans  $ \gerk_G $
over  $ \Bbbk \, $:  \, therefore  {\sl  $ \, \gerk_G \, $  is a quotient
of  $ \; \Bbbk \otimes_{\Z_p} \Cal{L}_p(G) \, $};  \, in fact, the above
is an isomorphism, see below.  Finally, we introduce also the notation
$ \, d(g) := n \, $  for all  $ \, g \in G_{[n]} \setminus G_{[n+1]} \, $.
                                      \par
   For each  $ \, k \in \N_+ \, $  choose a  $ \Z_p $--basis
$ \overline{B}_k $  of the  $ \Z_p $--vector  space  $ \, G_{[k]} \big/
G_{[k+1]} \, $;  \, for each  $ \, \overline{b} \in \overline{B}_k \, $,
\, fix  $ \, b \in G_{[k]} \, $  such that  $ \, \overline{b} = b \,
G_{[k+1]} \, $,  \, and let  $ \, B_k \, $  be the set of all such
elements  $ b \, $.  Let  $ \, B := \bigcup_{k \in \N_+} B_k \, $:
\, such a set will be called a  {\sl  $ p $-l.c.s.-net\/}  (=
``$ p $-lower-central-series-net''; the terminology in [HB] is
``$ \kappa $-net'') on  $ G $.  Of course  $ \, B_k = \big( B \cap
G_{[k]} \big) \setminus \big( B \cap G_{[k+1]} \big) \, $  for all
$ k \, $.  By an  {\sl ordered  $ p $-l.c.s.-net\/}  we mean a
$ p $-l.c.s.-net  $ B $  which is totally ordered in such a way
that:  {\it (i)\/}  if  $ \, a \in B_m \, $,  $ \, b \in B_n \, $,
$ \, m < n \, $,  \, then  $ \, a \preceq b \, $;  \, {\it (ii)\/}
for each  $ k $,  every non-empty subset of  $ B_k $  has a greatest
element (like for  $ \, p = 0 \, $).  Again,  $ p $-l.c.s.-nets  do
exist.
                                      \par
   We can now describe each  $ \Cal{D}_n(G) $,  hence also each graded
summand  $ \, J^n \big/ J^{n+1} \, $  of  $ \, \widehat{\Bbbk[G]} \, $, 
in terms of the lower central series or the  $ p $--lower  central series
of  $ G \, $,  more precisely in terms of a fixed ordered t.f.l.c.s.-net
or  $ p $-l.c.s.-net.  To unify notations, set  $ \, G_n := G_{(n)}
\, $,  $ \, \theta(g) := \partial(g) \, $  if  $ \, p \! = \! 0 \, $,
\, and  $ \, G_n := G_{[n]} \, $,  $ \, \theta(g) := d(g) \, $  if
$ \, p \! > \! 0 \, $,  \, set  $ \, G_\infty := \bigcap_{n \in N_+}
\!\! G_n \, $,  \, let  $ \, B := \bigcup_{k \in \N_+} B_k \, $  be
an ordered t.f.l.c.s.-net or  $ p $-l.c.s.-net according to whether
$ \, p \! = \! 0 \, $  or  $ \, p \! > \! 0 \, $,  \, and set  $ \,
\ell(0) := + \infty \, $  and  $ \, \ell(p) := p \, $  for  $ \, p
> 0 \, $.  The key result we need is

\vskip3pt

\noindent   {\it  $ \underline{\hbox{\sl Jennings-Hall theorem}} $
(cf.~[HB], [Pa] and references therein).  Let  $ \, p:= \text{\it
Char}\,(\Bbbk) \, $.
                                       \par
   (a) \, For all  $ \, g \in G \, $,  $ \; \eta_g
\in J^n \Longleftrightarrow g \in \! G_n \, $.  Therefore
$ \, \Cal{D}_n(G) = G_n \; $  for all  $ \, n \in \N_+ \, $.
                                       \par
   (b) \, For any  $ \, n \in \N_+ \, $,  the set of ordered monomials
 \vskip-19pt
  $$  \Bbb{B}_n \, := \, \Big\{\, {\overline{\,\eta_{b_1}}}^{\;e_1}
\cdots {\overline{\,\eta_{b_r}}}^{\;e_r} \;\Big|\; b_i \in B_{d_i}
\, , \; e_i \in \N_+ \, , \; e_i < \ell(p) \, , \; b_1 \precneqq
\cdots \precneqq b_r \, , \; {\textstyle \sum}_{i=1}^r e_i \, d_i
= n \,\Big\}  $$
 \vskip-9pt
\noindent   is a\/  $ \Bbbk $--basis  of  $ \, J^n \big/ J^{n+1} \, $,
\, and  $ \; \Bbb{B} \, := \, \{1\} \cup \bigcup_{n \in \N} \Bbb{B}_n
\; $  is a\/  $ \Bbbk $--basis  of  $ \; \widehat{\Bbbk[G]} \, $.
                                          \par
   (c) \;  $ \big\{\, \overline{\,\eta_b} \,\;\big|\;
b \in B_n \big\} \; $  is a  $ \Bbbk $--basis  of the  $ n $--th
graded summand  $ \, \gerk_G \cap \big( J^n \big/ J^{n+1} \big) \, $
of the graded restricted Lie algebra\/  $ \gerk_G \, $,  \, and
$ \, \big\{\, \overline{\,\eta_b} \,\;\big|\; b \in B \,\big\} \; $
is a\/  $ \Bbbk $--basis  of\/  $ \gerk_G \, $.
                                          \par
   (d) \;  $ \big\{\, \overline{\,\eta_b} \,\;\big|\;
b \in B_1 \big\} \; $  is a minimal set of generators of
the (restricted) Lie algebra\/  $ \gerk_G \, $.
                                          \par
   (e) \; The map  $ \; \Bbbk \otimes_\Z \L_p(G) \longrightarrow
\gerk_G \, $,  $ \, \overline{g} \mapsto \overline{\,\eta_g} \, $,
\, is an isomorphism of graded restricted Lie algebras.  Therefore
$ \; \widehat{\Bbbk[G]} \, \cong \, \U \big( \Bbbk \otimes_\Z \L_p(G)
\big) \; $  as Hopf algebras (notation of \S 1.1).
                                          \par
   (f) \;  $ J^\infty \, = \, \hbox{\sl Span} \big( \big\{\,
\eta_g \,\big|\, g \in G_\infty \,\big\} \big) \, $,  \, whence\/
$ \; {\Bbbk[G]}^\vee \cong \, \bigoplus_{\overline{g} \in G/G_\infty}
\! \Bbbk \cdot \overline{g} \; \cong \, \Bbbk \big[ G \big/ G_\infty
\big] \; $.   \qed}

\vskip5pt

   Recall that  $ A\big[x,x^{-1}\big] $  (for any  $ A $)  has
$ A $--basis  $ \, \big\{ {(x \!-\! 1)}^n x^{-[n/2]} \,\big|\, n \in
\N \big\} \, $,  \, where  $ [q] $  is the integer part of  $ \, q \in
\Bbb{Q} \, $.  Then from Jennings-Hall theorem and (5.2) we argue

\vskip5pt

\noindent   {\it  $ \underline{\hbox{\sl Proposition B}} $.
\, Let  $ \; \chi_g := \h^{-\theta(g)} \eta_g \, $,  \, for
all  $ \, g \in \{G\} \setminus \{1\} \, $.  Then
  $$  \displaylines{
   {R[G]}^\vee  = \,  \Big( {\textstyle
\bigoplus_{\Sb  b_i \in B, \; 0 < e_i < \ell(p)  \\
                r \in \N, \; b_1 \precneqq \cdots \precneqq b_r  \endSb}}
R \cdot \chi_{b_1}^{\;e_1} \, b_1^{\,-[e_1\!/2]} \cdots \chi_{b_r}^{\;
e_r} \, b_r^{\,-[e_r/2]} \Big) \,{\textstyle
\bigoplus}\, R\big[\h^{-1}\big] \cdot J^\infty  \; =   \hfill  \cr
   {} \;  = \,  \Big( {\textstyle
\bigoplus_{\Sb  b_i \in B, \; 0 < e_i < \ell(p)  \\
                r \in \N, \; b_1 \precneqq \cdots \precneqq b_r  \endSb}}
R \cdot \chi_{b_1}^{\,e_1} \, b_1^{\,-[e_1\!/2]} \cdots \chi_{b_r}^{\,
e_r} \, b_r^{\,-[e_r/2]} \Big) \,{\textstyle \bigoplus}\,
\Big( \hskip1pt {\textstyle \sum_{\gamma \in G_\infty}}
\hskip-0pt R\big[\h^{-1}\big] \cdot \eta_\gamma \Big) \; ;
\hskip-2,5pt  \cr }  $$
If  $ \, J^\infty \! = \! J^n $  for some  $ \, n \! \in \! \N $
(iff  $ \, G_\infty \! = G_n $)  we can drop the factors  $ \,
b_1^{-[e_1\!/2]}, \dots, b_r^{-[e_r/2]} \, . \, \square $}

\vskip5pt

\noindent   {\it  $ \underline{\hbox{\sl Poisson groups from
$ \Bbbk[G] $}} $.} \, The previous discussion attached to the
abstract group  $ G $  the (maybe restricted) Lie algebra  $ \gerk_G $
which, by Jennings-Hall theorem, is just the scalar extension of the Lie
ring  $ \L_{\text{\it Char}(\Bbbk)} $  associated to  $ G $  via the
central series of the  $ G_n $'s;  in particular the functor  $ \, G
\mapsto \gerk_G \, $  is one considered since long in group theory.
                                            \par
   Now, by  Theorem 5.8{\it (d)\/}  we know that  $ \, {\big( {R[G]}^\vee
\big)}' \, $  is a QFA, with  $ \, {\big({R[G]}^\vee \big)}'{\Big|}_{\h=0}
= F\big[\varGamma_G\big] \, $  for some connected Poisson group
$ \varGamma_G \, $.  This defines a functor  $ \, G \mapsto \varGamma_G
\, $  from abstract groups to connected Poisson groups, of dimension
zero and height 1 if  $ \, p > 0 \, $;  \, in particular, this
$ \varGamma_G $  {\sl is a new invariant for abstract groups}.
                                            \par
   The description of  $ {R[G]}^\vee $  in Proposition B above leads
us to an explicit description of  $ {\big({R[G]}^\vee \big)}' $,  \,
hence of  $ \, {\big({R[G]}^\vee \big)}'{\Big|}_{\h=0} \! = F \big[
\varGamma_G\big] \, $  and  $ \varGamma_G $  too.  Indeed, direct
inspection gives  $ \, \delta_n\big(\chi_g\big) = \h^{(n-1) \theta(g)}
\chi_g^{\;\otimes n} \, $,  \, so  $ \, \psi_g := \h \, \chi_g = \h^{1
- \theta(g)} \eta_g \in {\big( {R[G]}^\vee \big)}' \setminus \h \,
{\big( {R[G]}^\vee \big)}' \, $  for each  $ \, g \in G \setminus
G_\infty \, $,  \, whilst for  $ \, \gamma \in G_\infty $  we have
$ \, \eta_\gamma \in J^\infty \, $  which implies also  $ \, \eta_\gamma
\in {\big({R[G]}^\vee \big)}' \, $,  \, and even  $ \, \eta_\gamma \in
\bigcap_{n \in \N} \h^n {\big({R[G]}^\vee \big)}' \, $.  Therefore
$ {\big({R[G]}^\vee\big)}' $  is generated by  $ \, \big\{\, \psi_g
\;\big|\; g \in G \setminus \{1\} \big\} \cup \big\{\, \eta_\gamma
\,\big|\, \gamma \in G_\infty \big\} \, $.  Moreover,  $ \, g = 1 +
\h^{\theta(g)-1} \psi_g \in {\big( {R[G]}^\vee \big)}' \, $  for every
$ \, g \in G \setminus G_\infty \, $,  \, and  $ \, \gamma = 1 + (\gamma
- 1) \in 1 + J^\infty \subseteq {\big({R[G]}^\vee \big)}' \, $  for  $ \,
\gamma \in G_\infty \, $.  This and the previous analysis along with
Proposition B prove next result, which in turn is the basis for
Theorem D below.

\vskip5pt

\noindent   {\it  $ \underline{\hbox{\sl Proposition C}} $.
  $$  \displaylines{
   {\big({R[G]}^\vee\big)}' \; = \; \Big( {\textstyle
\bigoplus_{\Sb  b_i \in B, \; 0 < e_i < \ell(p)  \\
                r \in \N, \; b_1 \precneqq \cdots \precneqq b_r  \endSb}}
R \cdot \psi_{b_1}^{\;e_1} \, b_1^{\,-[e_1\!/2]} \cdots \psi_{b_r}^{\;
e_r} \, b_r^{\,-[e_r/2]} \Big) \,{\textstyle \bigoplus}\,
R\big[\h^{-1}\big] \cdot J^\infty  \; =   \hfill  \cr
   {} \hskip9pt   \; = \;  \Big( {\textstyle
\bigoplus_{\Sb  b_i \in B, \; 0 < e_i < \ell(p)  \\
                r \in \N, \; b_1 \precneqq \cdots \precneqq b_r  \endSb}}
R \cdot \psi_{b_1}^{\,e_1} \, b_1^{\,-[e_1\!/2]} \cdots \psi_{b_r}^{\,
e_r} \, b_r^{\,-[e_r/2]} \Big) \,{\textstyle \bigoplus}\,
\Big( \hskip1pt {\textstyle \sum_{\gamma \in G_\infty}}
\hskip-0pt R\big[\h^{-1}\big] \cdot \eta_\gamma \Big) \; .  \cr }  $$
In particular,  $ \; {\big({R[G]}^\vee\big)}' = R[G] \; $  {\sl if
and only if  $ \; G_2 = \{1\} = G_\infty \; $.}  If in addition  $ \,
J^\infty \! = \! J^n \, $  for some  $ \, n \! \in \! \N \, $  (iff
$ \, G_\infty = G_n $)  then we can drop the factors  $ \, b_1^{\,
-[e_1\!/2]}, \dots, b_r^{\,-[e_r/2]} \, $.   \qed}

\vskip5pt

\noindent   {\it  $ \underline{\hbox{\sl Theorem D}} $.  \, Let
$ \; x_g := \psi_g \mod \h \; {\big( {R[G]}^\vee \big)}' \, $,
$ \, z_g := g \mod \h \; {\big( {R[G]}^\vee \big)}' \, $  for
all  $ \, g \not= 1 \, $,  \, and  $ \, B_1 := \big\{\, b \in
\! B \,\big\vert\, \theta(b) = 1 \big\} \, $,  $ \, B_> :=
\big\{\, b \in \! B \,\big\vert\, \theta(b) > 1 \big\} \, $.
                                        \hfill\break
   \indent   (a) \, If  $ \, p = 0 \, $,  \, then  $ \, F \big[
\varGamma_G \big] = {\big( {R[G]}^\vee \big)}'{\Big|}_{\h=0} $  is
{\sl a polynomial/Laurent polynomial algebra}, namely  $ \, F \big[
\varGamma_G \big] = \Bbbk \big[ {\big\{ {z_b}^{\!\pm 1} \big\}}_{b
\in B_1} \!\! \cup {\{x_b\}}_{b \in B_>} \big] \, $,  \; the 
$ z_b $'s  being group-like and the  $ x_b $'s  being primitive. 
In particular  $ \, \varGamma_G \cong \big( \Bbb{G}_m^{\,\times B_1}
\big) \times \big( \Bbb{G}_a^{\, \times B_>} \big) \, $ as algebraic
groups, that is  $ \varGamma_G $  is a torus times a (pro)affine space.
                                        \hfill\break
   \indent   (b) \, If  $ \, p > 0 \, $,  \, then  $ \, F \big[
\varGamma_G \big] = {\big( {R[G]}^\vee \big)}'{\Big|}_{\h=0} $  is
{\sl a truncated polynomial/Laurent polynomial algebra}, namely  $ \,
F \big[ \varGamma_G \big] = \, \Bbbk \big[ {\big\{ {z_b}^{\!\pm 1}
\big\}}_{b \in B_1} \!\! \cup {\{x_b\}}_{b \in B_>} \big] \Big/ \! \big(
{\{ z_b^{\,p} \! - \! 1 \}}_{b \in B_1} \!\! \cup \big\{ x_b^{\,p}
\big\}_{b \in B_>} \big) \, $,  \, the  $ z_b $'s  being group-like
and the  $ x_b $'s  being primitive.  In particular  $ \, \varGamma_G
\cong \big( {{\boldsymbol\mu}_p}^{\! \times B_1} \big) \times \big(
{{\boldsymbol\alpha}_p}^{\!\times B_>} \big) \, $  as
algebraic groups of dimension zero and height 1.
                                        \hfill\break
   \indent   (c) \, The Poisson group  $ \varGamma_G $  has cotangent
Lie bialgebra  $ \gerk_G \, $,  that is  $ \, \text{\sl coLie}\,
(\varGamma_G) = \gerk_G \, $.}

\demo{Proof}  {\it (a)} \, Definitions give  $ \, \partial(g\,\ell\,)
\geq \partial(g) + \partial(\ell\,) \, $  for all  $ \, g, \ell \in G
\, $,  \, so that  $ \; [\psi_g,\psi_\ell\big] =
%
%
\h^{1 - \partial(g)
- \partial(\ell) + \partial((g,\ell))} \, \psi_{(g,\ell)} \, g \, \ell
\in \h \cdot {\big({R[G]}^\vee\big)}' \, $,  \; which proves (directly)
that  $ \, {\big({R[G]}^\vee\big)}'{\Big|}_{\h=0} \, $  is commutative!
Moreover, the relation  $ \, 1 = g^{-1} \, g =
g^{-1} \, \big(1 + \h^{\partial(g)-1} \psi_g \big) \, $  (for any
$ \, g \in G \, $)  yields  $ \, z_{g^{-1}} = {z_g}^{\!-1} \, $  iff
$ \, \partial(g) = 1 \, $  and  $ \, z_{g^{-1}} = 1 \, $  iff  $ \,
\partial(g) > 1 \, $.  Noting also that  $ \, J^\infty \equiv 0 \!\!
\mod \h \, {\big({R[G]}^\vee\big)}' \, $  and  $ \, g = 1 + \h^{\partial(g)
- 1} \psi_g \equiv 1 \mod \h \, {\big({R[G]}^\vee \big)}' \, $  for  $ \,
g \in G \setminus G_\infty \, $,  \, and also  $ \, \gamma = 1 +
(\gamma - 1) \in 1 + J^\infty \equiv 1 \mod \h \, {\big({R[G]}^\vee
\big)}' \, $  for  $ \, \gamma \in G_\infty \, $,  \, Proposition C
gives
  $$  F\big[\varGamma_G\big]  \; = \;
{\big({R[G]}^\vee\big)}'{\Big|}_{\h=0}
= \,  \Big( {\textstyle
\bigoplus_{\hskip-3pt   \Sb  b_i \in B_1, \; a_i \in \Z  \\
                s \in \N, \; b_1 \precneqq \cdots \precneqq b_s  \endSb}}
\hskip-2pt  \Bbbk \cdot z_{b_1}^{\,a_1} \cdots z_{b_s}^{\,a_s} \Big)
\,{\textstyle \bigotimes}\,  \Big( {\textstyle
\bigoplus_{\hskip-1pt   \Sb  b_i \in B_>, \; e_i \in \N_+  \\
                r \in \N, \; b_1 \precneqq \cdots \precneqq b_r  \endSb}}
\hskip-2pt  \Bbbk \cdot x_{b_1}^{\,e_1} \cdots x_{b_r}^{\,e_r} \Big)  $$  
which means that  $ F\big[\varGamma_G\big] $  is a polynomial-Laurent
polynomial algebra as claimed.  Again definitions imply   $ \, \Delta(z_g)
= z_g \otimes z_g \, $  for all  $ \, g \in G \, $  and  $ \, \Delta(x_g)
= x_g \otimes 1 + 1 \otimes x_g \, $  iff  $ \, \partial(g) > 1 \, $;  \,
thus the  $ z_b $'s  are group-like and the  $ x_b $'s  are primitive as
claimed.
                                          \par
   {\it (b)} \, The definition of  $ d $  implies  $ \, d(g\,\ell\,)
\geq d(g) + d(\ell\,) \, $  ($ g, \ell \in G $),  \, whence we get
$ \; [\psi_g,\psi_\ell] \, = \, \h^{1 - d(g) - d(\ell) + d((g,\ell))}
\, \psi_{(g,\ell)} \, g \, \ell \, \in \, \h \cdot {\big( {R[G]}^\vee
\big)}' \, $,  \, proving that  $ \, {\big( {R[G]}^\vee \big)}'
{\Big|}_{\h=0} \, $  is commutative.  In addition  $ \; d(g^p)
\geq p \; d(g) \, $,  \, so  $ \; \psi_g^{\;p} = \h^{\, p \,
(1 - d(g))} \, \eta_g^{\;p} =
%
%
\h^{\, p - 1 + d(g^p) - p\,d(g)} \,
\psi_{g^p} \in \h \cdot {\big({R[G]}^\vee\big)}' \, $,  \, whence
$ \, {\big( \psi_g^{\;p}{\big|}_{\h=0} \big)}^p = 0 \, $  inside  $ \,
{\big( {R[G]}^\vee \big)}'{\Big|}_{\h=0} \! = F \big[ \varGamma_G \big]
\, $,  \, which proves that  $ \varGamma_G $  has dimension 0 and height
1.  Finally  $ \; b^p = {(1 + \psi_b)}^p = 1 + {\psi_b}^{\!p} \equiv 1
\mod \h \, {\big( {R[G]}^\vee \big)}' \, $  for all  $ \, b \in B_1 \, $,
\, so  $ \, b^{-1} \equiv b^{p-1} \!\! \mod \h \, {\big( {R[G]}^\vee
\big)}' \, $.
   \hbox{Thus letting  $ \, x_g := \psi_g \!\! \mod \h \; {\big(
{R[G]}^\vee \big)}' $  (for  $ \, g \! \not= \! 1 $)  we get}
  $$  F\big[\varGamma_G\big] = {\big({R[G]}^\vee\big)}'{\Big|}_{\h=0} \! =
\Big( {\textstyle \bigoplus_{\Sb  b_i \in B_1, \; 0 < e_i < p  \\
                s \in \N, \; b_1 \precneqq \cdots \precneqq b_s  \endSb}}
\Bbbk \cdot z_{b_1}^{\,e_1} \cdots z_{b_s}^{\,e_s} \Big)
 {\textstyle \bigotimes}
\Big( {\textstyle \bigoplus_{\Sb  b_i \in B_>, \; 0 < e_i < p  \\
                r \in \N, \; b_1 \precneqq \cdots \precneqq b_r  \endSb}}
\Bbbk \cdot x_{b_1}^{\,e_1} \cdots x_{b_r}^{\,e_r} \Big)  $$
just like for  {\it (a)\/}  and also taking care that  $ \, z_b = x_b
+ 1 \, $  and  $ \, z_b^{\,p} = 1 \, $  for  $ \, b \in B_1 \, $.
Therefore  $ \, {\big({R[G]}^\vee\big)}'{\Big|}_{\h=0} \, $  is a
truncated polynomial/Laurent polynomial algebra as claimed.  The
properties of the  $ x_b $'s  and the  $ z_b $'s  w.r.t.~the Hopf
     \hbox{structure are then proved like for \!  {\it (a)\/}  again.}
                                                 \par
   {\it (c)} \, The augmentation ideal  $ \, \germ_e \, $  of  $ \,
{\big( {R[G]}^\vee \big)}'{\Big|}_{\h=0} = F\big[\varGamma_G\big] \, $
is generated by  $ \, {\{ x_b \}}_{b \in B} \, $;  \, then  $ \; \h^{-1}
\, [\psi_g,\psi_\ell\big] \, =
%
%
\, \h^{\,
\theta((g,\ell)) - \theta(g) - \theta(\ell)} \, \psi_{(g,\ell)}
\, \big( 1 + \h^{\, \theta(g) - 1} \psi_g \,\big) \, \big( 1 + \h^{\,
\theta(\ell) - 1} \psi_\ell \,\big) \, $
by the
previous computation, whence at  $ \, \h = 0 \, $
one has  $ \; \big\{ x_g \, , x_\ell \big\} \, \equiv \,
x_{(g,\ell)} \mod \, \germ_e^{\,2} \, $  if  $ \, \theta
\big( (g,\ell\,) \big) = \theta(g) + \theta(\ell\,) \, $,  and
$ \; \big\{ x_g \, , x_\ell \big\} \, \equiv \, 0 \mod \, \germ_e^{\,2}
\, $  if  $ \, \theta \big( (g,\ell\,) \big) > \theta(g) + \theta(\ell\,)
\, $.  This means that the cotangent Lie bialgebra  $ \, \germ_e \Big/
\germ_e^{\,2} \, $  of  $ \varGamma_G \, $  is isomorphic to  $ \gerk_G
\, $,  \, as claimed.   \qed
\enddemo

\vskip3pt

\noindent   {\it  $ \underline{\hbox{\sl Remarks}} $:}  {\it (a)}
\, Theorem D claims that the connected Poisson group  $ \, K_G^\star
:= \varGamma_G \, $  is  {\sl dual\/}  to  $ \gerk_G $  in the sense
of \S 1.1.  Since  $ \; {R[G]}^\vee{\Big|}_{\h=0} \! = \U(\gerk_G)
\; $  and  $ \; {\big( {R[G]}^\vee \big)}'{\Big|}_{\h=0} \! = F \big[
K_G^\star \big] \, $,  \,  {\sl this gives a close analogue, in positive
characteristic, of the second half of Theorem 2.2{\it (c)}.}
                                                \par
  {\it (b)} \, Theorem D gives functorial recipes to attach to
each abstract group  $ G $  and each field  $ \Bbbk $  a connected
Abelian algebraic Poisson group over\/  $ \Bbbk $,  namely  $ \; G
\mapsto \varGamma_G = K_G^\star \; $,  \; explicitly described as
algebraic group and such that  $ \, \text{\sl coLie}\,(K_G^\star) =
\gerk_G \, $.  {\sl Every such  $ \varGamma_G $  (for given  $ \Bbbk $)
is then an invariant of  $ G \, $,  a new one to the author's
knowledge.  Indeed, it is perfectly equivalent to the well-known
invariant\/  $ \gerk_G $}  (over the same  $ \Bbbk $),  because
clearly  $ \, \varGamma_{G_1} \! \cong \varGamma_{G_2} \, $ 
implies  $ \, \gerk_{G_1} \! \cong \gerk_{G_2} \, $,  \, whereas 
$ \, \gerk_{G_1} \! \cong \gerk_{G_2} \, $  implies that 
$ \varGamma_{G_1} $  and  $ \varGamma_{G_2} $  are isomorphic
as algebraic groups   --- by  Theorem D{\it (a--b)}  ---   and
bear isomorphic Poisson structures   --- by  part  {\it (c)} 
of Theorem D  ---   whence  $ \, \varGamma_{G_1} \! \cong
\varGamma_{G_2} \, $  as algebraic Poisson groups.  

\vskip5pt

\noindent   {\it  $ \underline{\hbox{\sl The case of  $ A_\Bbbk(G)
\, $}} $.} \, Let's now dwell upon  $ \, H := A_\Bbbk(G) \, $,  \,
for a  {\sl finite\/}  group  $ G \, $.
                                          \par
   Let  $ \A $  be any commutative unital ring, and let  $ \Bbbk $,
$ \, R := \Bbbk[\h] \, $  be as before.  By definition  $ \, A_\A(G)
= {\A[G]}^* \, $,  \, hence  $ \, \A[G] = {A_\A(G)}^* \, $,  \, and
we have a natural perfect Hopf pairing  $ \, A_\A(G) \times \A[G]
\longrightarrow \A \, $.
%
%
 Our first result is one of triviality:

\vskip5pt

\noindent   {\it  $ \underline{\hbox{\sl Theorem E}} $.  $ \,
{A_R(G)}^\vee \! = R \cdot 1 \oplus R \big[ \h^{-1} \big] \,
J = {\big({A_R(G)}^\vee\big)}' \, $,  $ \, {A_\Bbbk(G)}^\vee
= \Bbbk \!\cdot\! 1 \, $,  $ \, \widehat{A_\Bbbk(G)} =
{A_R(G)}^\vee{\Big|}_{\h=0} \! = \, \Bbbk \cdot 1 = \,
\U(\boldkey{0}) \; $  and  $ \; {\big( {A_R(G)}^\vee
\big)}'{\Big|}_{\h=0} \! = \, \Bbbk \cdot 1 =
F\big[\{*\}\big] \; $.}

\demo{Proof} By construction  $ \, J := \text{\sl Ker}\,(\epsilon_{\!
{}_{A_\Bbbk(G\,)}}) \, $  has  $ \Bbbk $--basis  $ \, \big\{ \varphi_g
\big\}_{g \in G \setminus \{1_{{}_G}\}} \cup \big\{ \varphi_{1_G} \! -
1_{\!{}_{A_\Bbbk(G\,)}} \big\} \, $,  \, and since  $ \, \varphi_g =
{\varphi_g}^{\!2} \, $  for all  $ g $  and  $ \, {(\varphi_{1_G} \!
- \! 1)}^2 = -(\varphi_{1_G} \!-\! 1) \, $  we have  $ \, J = J^\infty
\, $,  \, so  $ \, {A_\Bbbk(G)}^\vee = \Bbbk \!\cdot\! 1 \, $  and
$ \, \widehat{A_\Bbbk(G)} = \Bbbk \!\cdot\! 1 \, $.  Similarly,  $ \,
{A_R(G)}^\vee \, $  is generated by  $ \, \big\{ \h^{-1} \varphi_g
\big\}_{g \in G \setminus \{1_{{}_G}\}} \cup \big\{ \h^{-1} (\varphi_{1_G}
\! - \! 1_{\!{}_{A_R(G\,)}}) \big\} \, $;  \, moreover,  $ \, J = J^\infty
\, $  implies  $ \, \h^{-n} J \subseteq {A_R(G)}^\vee \, $  for all  $ \,
n \, $,  \, whence  $ \, {A_R(G)}^\vee = R \, 1 \oplus R[\h^{-1}]
J \, $.  Then  $ \, J_{{A_R(G)}^\vee} = R\big[\h^{-1}\big] J \subseteq
\h \, {A_R(G)}^\vee \, $,  \, which implies  $ \, {\big( {A_R(G)}^\vee
\big)}' = {A_R(G)}^\vee \, $:  \, in particular,  $ \, {\big( {A_R(G)}^\vee
\big)}'{\Big|}_{\h=0} = {A_R(G)}^\vee{\Big|}_{\h=0} \! = \Bbbk \cdot 1
\, $,  \, as claimed.   \qed
\enddemo

%
%
 \vskip1pt

\noindent   {\it  $ \underline{\hbox{\sl Poisson groups from
$ A_\Bbbk(G) \, $}} $.} \, Now we look at  $ {A_R(G)}' $,
$ {A_\Bbbk(G)}' $  and  $ \widetilde{A_\Bbbk(G)} \, $.  By
construction  $ A_R(G) $  and  $ R[G] $  are in perfect Hopf pairing,
and are free  $ R $--modules  of  {\sl finite rank}.  In this case, since
Proposition 4.4 yields  $ \, {A_R(G)}' = {\big({R[G]}^\vee\big)}^\bullet
\, $  we have in fact  $ \, {A_R(G)}' = {\big({R[G]}^\vee\big)}^\bullet
= {\big( {R[G]}^\vee \big)}^* \, $:  thus  $ \, {A_R(G)}' \, $  is the
{\sl dual\/}  Hopf algebra to  $ {R[G]}^\vee $;  then from Proposition
B we can argue an explicit description of  $ {A_R(G)}' $,  whence
also of  $ \big({A_R(G)}'\big)^{\!\vee} $.  By Theorem 5.10{\it (g)\/}
and its proof, namely that  $ \, {A_\Bbbk(G)}' = \big( J_{\Bbbk[G]}^{\;
\infty} \big)^\perp $,  there is a perfect filtered Hopf pairing  $ \,
{\Bbbk[G]}^\vee \times {A_\Bbbk(G)}' \longrightarrow \Bbbk \, $  and
a perfect graded Hopf pairing  $ \, \widehat{\Bbbk[G]} \times
\widetilde{A_\Bbbk(G)} \! \longrightarrow \Bbbk \, $:  \, thus  $ \,
{A_\Bbbk(G)}' \! \cong \! {\big( {\Bbbk[G]}^\vee \big)}^* \, $  as
filtered Hopf algebras and  $ \, \widetilde{A_\Bbbk(G)} \cong {\big(
\widehat{\Bbbk[G]}\big)}^* \, $  as graded Hopf algebras.  If  $ p = 0
\, $  then  $ J = J^\infty $,  \, as each  $ \, g \in G \, $  has finite
order and  $ \, g^n = 1 \, $  implies  $ \, g \in G_\infty \, $:  \, then
$ \, {\Bbbk[G]}^\vee \! = \Bbbk \cdot 1 = \widehat{\Bbbk[G]} \, $,
\, so  $ {A_\Bbbk(G)}' = \Bbbk \cdot 1 = \widetilde{A_\Bbbk(G)}
\, $.  If  $ \, p > 0 \, $  instead, this analysis gives  $ \,
\widetilde{A_\Bbbk(G)} = {\big( \widehat{\Bbbk[G]} \big)}^* = {\big(
\u(\gerk_G) \big)}^* = F[K_G] \, $,  \, where  $ \, K_G  \, $  is a
connected Poisson group of dimension 0, height 1 and tangent Lie
bialgebra  $ \gerk_G \, $.  Thus

\vskip7pt

\noindent   {\it  $ \underline{\hbox{\sl Theorem F}} $.
                                                \par
   (a) \, There is a second functorial recipe to attach to each
finite abstract group a connected algebraic Poisson group of
dimension zero and height 1 over any field\/  $ \Bbbk $  with
$ \, \text{\it Char}\,(\Bbbk) > 0 \, $,  namely  $ \; G \mapsto
K_G := \text{\it Spec}\,\Big( \widetilde{A_\Bbbk(G)} \Big) $.
This  $ K_G $  is Poisson dual to  $ \varGamma_G $  of Theorem
D in the sense of \S 1.1, in that  $ \; \text{\sl Lie}\,(K_G)
= \, \gerk_G = \text{\sl coLie}\,(\varGamma_G) \; $.
                                                \par
   (b) \, If  $ \, p:= \text{\it Char}\,(\Bbbk) > 0 \, $,  then
$ \; \big({A_R(G)}'\big)^{\!\vee}\Big|_{\h=0} \! = \, \u \big(
\gerk_G^{\,\times} \big) = S\big(\gerk_G^{\,\times}\big) \Big/ \big(
\big\{ x^p \,\big|\, x \in \gerk_G^{\,\times} \big\} \big) \; $.}

\demo{Proof} Claim  {\it (a)\/}  is the outcome of the discussion
above.  Part  {\it (b)\/}  instead requires an explicit description
of  $ \, \big({A_R(G)}'\big)^{\!\vee} $.  Since  $ \, {A_R(G)}'
\cong \big({R[G]}^\vee\big)^* \, $,  \, from Proposition B we get
$ \; {A_R(G)}' = \, \Big(\! \bigoplus_{\Sb  b_i \in B, \; 0 < e_i < p  \\
                 r \in \N, \; b_1 \precneqq \cdots \precneqq b_r  \endSb}
R \cdot \rho_{b_1, \dots, b_r}^{e_1, \dots, e_r} \Big) \; $
where each  $ \, \rho_{b_1, \dots, b_r}^{e_1, \dots, e_r} \, $  is
defined by
                      \hfill\break
 \vskip2pt
   \centerline{$ \Big\langle \, \rho_{b_1, \dots, b_r}^{e_1, \dots,
e_r} \; , \; \chi_{\beta_1}^{\;\varepsilon_1} \, \beta_1^{\,
-[\varepsilon_1\!/2]} \cdots \chi_{\beta_s}^{\;\varepsilon_s} \,
\beta_s^{\,-[\varepsilon_s/2]} \,\Big\rangle \; = \; \delta_{r,s}
\, \prod_{i=1}^r \delta_{b_i,\beta_i} \delta_{e_i, \varepsilon_i} $}
 \vskip3pt
\noindent   (for all  $ \, b_i, \beta_j \in B \, $  and  $ \, 0 <
e_i, \varepsilon_j < p \, $).  Now, using notation of \S 1.3,  $ \,
K_\infty \subseteq K' \, $  for any  $ \, K \in \HA \, $,  \, whence
$ \, K' = \pi^{-1} \big(\, \overline{K}^{\;\prime} \,\big) \, $  where
$ \; \pi \, \colon \, K \relbar\joinrel\twoheadrightarrow K \big/
K_\infty =: \overline{K} \; $  is the canonical projection.  So let
$ \, K := {R[G]}^\vee \, $,  $ \, H := {A_R(G)}' \, $;  \, Proposition
B gives  $ \, K_\infty = R\big[\h^{-1}\big] \cdot J^\infty \, $  and
provides at once a description of  $ \, \overline{K} \, $;  \, from this
and the previous description of  $ H $  one sees also that in the present
case  $ K_\infty $  is exactly the right kernel of the natural pairing
$ \, H \times K \longrightarrow R \, $,  \, which is perfect on the left,
so that the induced pairing  $ \, H \times \overline{K} \longrightarrow
R \, $  is perfect.  By construction its specialization at  $ \, \h =
0 \, $  is the natural pairing  $ \, F[K_G] \times \u(\gerk_G)
\longrightarrow \Bbbk \, $,  \, which is perfect too.  Then we can
apply  Proposition 4.4{\it (c)\/}  (with  $ \overline{K} $  playing
the r\^{o}le of  $ K $  therein)  which yields  $ \, \overline{K}^{\;
\prime} \! = \big(H^\vee\big)^\bullet = \Big(\!\! \big({A_R(G)}'\big)^{\!
\vee} \Big)^{\!\bullet} \, $.  By construction,  $ \, \overline{K}^{\;
\prime} \! = \big({R[G]}^\vee\big)' \! \Big/ \! \big( R\big[\h^{-1}\big]
\cdot J^\infty \big) \, $,  \, and Proposition C describes the latter as
$ \; \overline{K}^{\;\prime} = \, \Big(\! \bigoplus_{\Sb  b_i \in B,
                                                     \; 0 < e_i < p  \\
               r \in \N, \; b_1 \precneqq \cdots \precneqq b_r  \endSb}
\hskip-3pt  R \cdot \overline{\psi}_{b_1}^{\;e_1} \cdots
\overline{\psi}_{b_r}^{\;e_r} \Big) \, $,
\; where  $ \, \overline{\psi}_{b_i} := \psi_{b_i}
\mod R\big[\h^{-1}\big] \cdot J^\infty \, $  for all  $ i \, $;  \,
since  $ \, \overline{K}^{\;\prime} \! = \Big(\!\! \big( {A_R(G)}'
\big)^{\!\vee} \Big)^{\!\bullet}  \, $  and  $ \, \psi_g = \h^{+1}
\chi_g \, $,  \, this analysis yields
$ \; \big({A_R(G)}'\big)^{\!\vee} = \, \Big(\!
             \bigoplus_{\Sb  b_i \in B, \; 0 < e_i < p  \\
      r \in \N, \; b_1 \precneqq \cdots \precneqq b_r  \endSb}
R \cdot \h^{- \sum_i e_i} \rho_{b_1, \dots, b_r}^{e_1,\dots,e_r} \Big)
\cong \big(\overline{K}^{\;\prime}\big)^* \, $,
\; whence we get  $ \; \big({A_R(G)}'\big)^{\!\vee}\Big|_{\h=0} \! \cong
\big(\overline{K}^{\;\prime}\big)^*\Big|_{\h=0} \! = \big(K'\big|_{\h=0}
\big)^* \! = \Big(\! \big({R[G]}^\vee\big)'\big|_{\h=0}\Big)^* \! \cong
{F\big[\varGamma_G\big]}^* = \u \big( \gerk_G^{\,\times} \big) =
S\big(\gerk_G^{\,\times}\big) \Big/ \big( \big\{ x^p \,\big|\, x \in
\gerk_G^{\,\times} \big\} \big) \; $  as claimed, the latter identity
being trivial (as  $ \gerk_G^{\,\times} $  is Abelian).   \qed
\enddemo

\vskip5pt

\noindent   {\it  $ \underline{\hbox{\sl Remarks}} $:}  {\it (a)} \,
this  $ K_G $  is another invariant for  $ G $,  but again equivalent
to  $ \gerk_G \, $.
                                                \par
   {\it (b)} \, {\sl Theorem F{\it (b)\/}  is a positive
characteristic analogue for  $ \, F_\h[G] = {A_R(G)}' \, $
of the first half of  Theorem 2.2{\it (c)}.}

\vskip11pt

\noindent   {\it  $ \underline{\hbox{\sl Examples}} $:}
 \vskip5pt
   {\sl (1) \, Finite Abelian  $ p \, $--groups.} \, Let  $ p \, $
be a prime number and  $ \, G := \Z_{p^{e_1}} \times \Z_{p^{e_2}}
\times \cdots \times \Z_{p^{e_k}} \, $  ($ k, e_1, \ldots, e_k \in
\N \, $),  \, with  $ \, e_1 \geq e_2 \geq \cdots \geq e_k \, $.
Let  $ \Bbbk $  be a field with  $ \, \Char(\Bbbk) = p > 0 \, $,
\, and  $ \, R := \Bbbk[\h] \, $  as above, so that  $ \,
\Bbbk[G]_\h = R[G] \, $.
                                              \par
   First,  $ \, \gerk_G \, $  is Abelian, because  $ G $  is.
Let  $ g_i $  be a generator of  $ \, \Z_{p^{e_i}} \, $  (for
all  $ i \, $),  identified with its image in  $ G \, $.  Since
$ G $  is Abelian we have  $ \, G_{[n]} = G^{p^n} \, $  (for all
$ n \, $),  and an ordered  $ p $-l.c.s.-net  is  $ \, B :=
\bigcup_{r \in \N_+} B_r \, $  with  $ \, B_r := \Big\{\, g_1^{\,p^r},
\, g_2^{\,p^r}, \, \ldots , \, g_{j_r}^{\,p^r} \Big\} \, $  where
$ j_r $  is uniquely defined by  $ \; e_{j_r} > r \, $,  $ \; e_{j_r
+ 1} \leq r \, $.  Then  $ \, \gerk_G \, $  has  $ \Bbbk $--basis
$ \, {\big\{\, \overline{\,\eta_{g_i^{p^{s_i}}}} \,\big\}}_{1 \leq i
\leq k; \; 0 \leq s_i < e_i} \, $,  \, and minimal set of generators
(as a restricted Lie algebra)  $ \, \big\{\,\overline{\,\eta_{g_1}}
\, , \, \overline{\,\eta_{g_2}} \, , \, \ldots, \, \overline{\,
\eta_{g_k}} \,\big\} \, $,  \, for the  $ p $--operation  of
$ \gerk_G $  is  $ \, {\big(\overline{\,\eta_{g_i^{p^s}}}
\,\big)}^{[p\hskip0,5pt]} = \overline{\,\eta_{g_i^{p^{s+1}}}} \, $,
\, and the order of nilpotency of each  $ \, \overline{\,\eta_{g_i}}
\, $  is exactly  $ p^{e_i} $,  i.e.~the order of  $ g_i \, $.  In
addition  $ \, J^\infty = \{0\} \, $  so  $ \, {\Bbbk[G]}^\vee \!
= \Bbbk[G] \, $.  The outcome is  $ \; {\Bbbk[G]}^\vee \! = \,
\Bbbk[G] \; $  and
 \vskip3pt
   \centerline{ $  \widehat{\Bbbk[G]}  \, = \,  \u(\gerk_G)  \; = \;
U(\gerk_G) \bigg/ \Big( {\Big\{ {\big(\overline{\,\eta_{g_i^{p^s}}}
\,\big)}^p - \overline{\,\eta_{g_i^{p^{s+1}}}} \,\Big\}}_{1 \leq i
\leq k}^{0 \leq s < e_i} \bigcup \; {\Big\{ {\big(\overline{\,
\eta_{g_i^{p^{e_i-1}}}}\,\big)}^p \,\Big\}}_{1 \leq i \leq k}
\, \Big) $ }
 \vskip3pt
\noindent
whence  $ \; \widehat{\Bbbk[G]} \, \cong \, \Bbbk[x_1,\dots,x_k] \bigg/
\! \Big( \Big\{\, x_i^{p^{e_i}} \;\Big|\; 1 \leq i \leq k \,\Big\} \Big)
\, $,  \; via  $ \; \overline{\,\eta_{g_i^{p^s}}} \mapsto x_i^{\,p^s}
\, $  (for all  $ i $,  $ s \, $).
                                              \par
   As for  $ {\Bbbk[G]}_\h^{\,\vee} $,  for all  $ \, r < e_i \, $  we
have  $ \, d\big(g_i^{p^r}\big) = p^r \, $  and so  $ \, \chi_{g_i^{p^r}}
= \h^{-p^r} \big( g_i^{p^r} \!-\! 1 \big) \, $  and  $ \, \psi_{g_i^{p^r}}
= \h^{1-p^r} \big( g_i^{p^r} \!-\! 1 \big) \, $;  \, since  $ \,
G_{[\infty]} = \{1\} \, $  (or, equivalently,  $ \, J^\infty = \{0\}
\, $)  and everything is Abelian, from the general theory we conlude
that both  $ {\Bbbk[G]}_\h^{\,\vee} $  and  $ \Big( {\Bbbk[G]}_\h^{\,
\vee} \Big)' $  are truncated-polynomial algebras, in the
$ \chi_{g_i^{p^r}} $'s  and in the  $ \psi_{g_i^{p^r}} $'s
respectively, namely
 \vskip3pt
   \centerline{ $ \eqalign{
   {\Bbbk[G]}_\h^{\,\vee}  &  = \; \Bbbk[\h] \Big[ \big\{\,
\chi_{g_i^{p^s}} \big\}_{1 \leq i \leq k \, ; \; 0 \leq s < e_i} \Big]
\; \cong \;  \Bbbk[\h]\big[\,y_1,\dots,y_k\big] \bigg/ \! \Big(
\Big\{\, y_i^{\,p^{e_i}} \;\Big|\; 1 \leq i \leq k \,\Big\} \Big)  \cr
   \Big( {\Bbbk[G]}_\h^{\,\vee} \Big)'  &  = \, \Bbbk[\h] \Big[
\big\{ \, \psi_{g_i^{p^s}} \big\}_{1 \leq i \leq k \, ; \; 0 \leq
s < e_i} \Big] \, \cong \, \Bbbk[\h] \Big[ \big\{\, z_{i,s} \big\}_{1
\leq i \leq k \, ; \; 0 \leq s < e_i} \Big] \bigg/ \! \Big( \Big\{\,
{z_{i,s}}^{\!p} \;\Big|\; 1 \leq i \leq k \,\Big\} \Big)  \cr } $ }
 \vskip3pt
\noindent
via the isomorphisms given by  $ \, \overline{\,\chi_{g_i^{p^s}}}
\mapsto y_i^{\,p^s} \, $  and  $ \, \overline{\,\psi_{g_i^{p^s}}}
\mapsto z_{i,s} \, $  (for all  $ i $,  $ s \, $).  When  $ \, e_1
> 1 \, $  this implies  $ \, {\big( {\Bbbk[G]}_\h^{\,\vee\,} \big)}'
\supsetneqq {\Bbbk[G]}_\h \, $,  \, that is a counterexample to
Theorem 2.2{\it (b)}.  Setting  $ \; \overline{\psi_{g_i^{p^s}}}
:= \psi_{g_i^{p^s}} \mod \h \; {\big( {\Bbbk[G]}_\h^{\,\vee\,} \big)}'
\; $  (for all  $ \, 1 \leq i \leq k \, $,  $ \, 0 \leq s < e_i \, $)
we have
                              \hfill\break
   \centerline{ $ F\big[\varGamma_G\big] \, = \, {\big(
{\Bbbk[G]}_\h^{\,\vee\,} \big)}'{\Big|}_{\h=0} = \,
\Bbbk \Big[ \big\{ \overline{\psi_{g_i^{p^s}}}
\,\big\}_{1 \leq i \leq k}^{0 \leq s < e_i} \Big] \, \cong
\, \Bbbk\,\Big[ \big\{\, w_{i,s} \big\}_{1 \leq i \leq k}^{0
\leq s < e_i} \Big] \bigg/ \! \Big( \Big\{\, w_{i,s}^{\,p}
\Big|\; 1 \!\leq\! i \!\leq\! k \,\Big\} \Big) $ }
(via  $ \, \overline{\psi_{g_i^{p^s}}} \mapsto w_{i,s} \, $)
as a  $ \Bbbk $--algebra.  The Poisson bracket trivial, and
the  $ w_{i,s} $'s  are primitive for  $ \, s > 1 \, $  and
$ \, \Delta(w_{i,1}) = w_{i,1} \otimes 1 + 1 \otimes w_{i,1}
+ w_{i,1} \otimes w_{i,1} \, $  for all  $ 1 \leq i \leq k \, $.
%
%
   If instead  $ \; e_1 = \cdots = e_k = 1 \, $,  \, then  $ \,
{\big( {\Bbbk[G]}_\h^{\,\vee\,} \big)}' = {\Bbbk[G]}_\h \, $. 
This is an analogue of  Theorem 2.2{\it (b)},  although now  $ \,
\Char(\Bbbk) > 0 \, $,  \, in that in this case  $ \, {\Bbbk[G]}_\h
\, $  is a QFA, with  $ \, {\Bbbk[G]}_\h{\Big|}_{\h=0} \! = \Bbbk[G]
= F\big[\widehat{G} \,\big] \, $  where  $ \, \widehat{G} \, $  is
the  {\sl group of characters\/}  of  $ G \, $.  But then  $ \, F
\big[\widehat{G}\,\big] = \Bbbk[G] = {\Bbbk[G]}_\h{\Big|}_{\h=0}
\! = {\big( {\Bbbk[G]}_\h^{\,\vee\,} \big)}'{\Big|}_{\h=0} \! =
F\big[\varGamma_G\big] \, $  (by general analysis) which means
that  $ \widehat{G} $  can be realized as a finite, connected,
    \hbox{Poisson group-scheme of dimension 0 and height 1 dual
to  $ \gerk_G \, $,  i.e.~$ \, \varGamma_G = K_G^\star \, $.}
                                                    \par
   Finally, a direct easy calculation shows that   --- letting  $ \,
\chi^*_g := \h^{\,d(g)} \, (\varphi_g - \varphi_1) \in {A_\Bbbk(G)}'_\h
\, $  and  $ \, \psi^*_g := \h^{\,d(g)-1} \, (\varphi_g - \varphi_1)
\in {\big({A_\Bbbk(G)}'\big)}^\vee_\h \, $  (for all  $ \, g \in G
\setminus \{1\} \, $)  ---   we have also
                              \hfill\break
   \centerline{ $  \eqalign{
   {A_\Bbbk(G)}_\h^{\,\prime} \,  &  = \; \Bbbk[\h] \Big[ \big\{\,
\chi^*_{g_i^{p^s}} \big\}_{1 \leq i \leq k}^{0 \leq s < e_i}
\Big] \, \cong \; \Bbbk[\h] \Big[ \big\{ Y_{i,j} \big\}_{1
\leq i \leq k}^{0 \leq s < e_i} \Big] \bigg/ \! \Big( \big\{
Y_{i,j}^{\;p} \big\}_{1 \leq i \leq k}^{0 \leq s < e_i} \Big)  \cr
   {\big({A_\Bbbk(G)}_\h^{\,\prime}\,\big)}^{\!\vee}  &  = \; \Bbbk[\h]
\Big[ \big\{\, \psi^*_{g_i^{p^s}} \big\}_{1 \leq i \leq k}^{0 \leq
s < e_i} \Big] \, \cong \; \Bbbk[\h] \Big[ \big\{ Z_{i,s} \big\}_{1
\leq i \leq k}^{0 \leq s < e_i} \Big] \bigg/ \! \Big( \big\{
Z_{i,s}^{\;p} - Z_{i,s} \big\}_{1 \leq i \leq k}^{0 \leq s
< e_i} \Big)  \cr } $ }
 \vskip3pt
\noindent
via the isomorphisms given by  $ \, \chi^*_{g_i^{p^s}} \mapsto
Y_{i,s} \, $  and  $ \, \psi^*_{g_i^{p^s}} \mapsto Z_{i,s} \, $,
\, from which one also gets the analogous descriptions of  $ \,
{A_\Bbbk(G)}_\h^{\,\prime}{\Big|}_{\h=0} \! = \widetilde{A_\Bbbk(G)}
= F[K_G] \, $  and of  $ \, {\big({A_\Bbbk(G)}_\h^{\,\prime}\,
\big)}^{\!\vee}{\Big|}_{\h=0} \! = \u(\gerk_G^\times) \; $.
 \vskip5pt
   {\sl (2) \, A non-Abelian  $ p \, $--group.} \, Let  $ p \, $
be a prime number,  $ \Bbbk $  be a field with  $ \, \Char(\Bbbk)
= p > 0 \, $,  \, and  $ \, R := \Bbbk[\h] \, $  as above, so
that  $ \, \Bbbk[G]_\h = R[G] \, $.
                                              \par
   Let  $ \, G := \Z_p \ltimes \Z_{p^{\,2}} \, $,  \, that is the group
with generators  $ \, \nu $,  $ \tau \, $  and relations  $ \, \nu^p
= 1 \, $,  $ \, \tau^{p^2} = 1 \, $,  $ \, \nu \, \tau \, \nu^{-1}
= \tau^{1+p} \, $.  In this case,  $ \, G_{[2]} = \cdots = G_{[p\,]}
= \big\{ 1, \tau^p \,\big\} \, $,  $ \, G_{[p+1]} = \{1\} \, $,  \,
so we can take  $ \, B_1 = \{\nu \, , \tau \,\} \, $  and  $ \, B_p
= \big\{ \tau^p \,\big\} \, $  to form an ordered  $ p $-l.c.s.-net
$ \, B := B_1 \cup B_p \, $  w.r.t.~the ordering  $ \, \nu \preceq
\tau \preceq \tau^p \, $.  Noting also that  $ \, J^\infty = \{0\}
\, $  (for  $ \, G_{[\infty]} = \{1\} \, $),  we have
 \vskip1pt
   \centerline{ $ {{\Bbbk[G]}_\h}^{\!\vee} \; = \; {\textstyle
\bigoplus_{a,b,c=0}^{p-1}} \, \Bbbk[\h] \cdot \chi_\nu^{\,a}
\, \chi_\tau^{\,b} \, \chi_{\tau^p}^{\,c} \; = \; {\textstyle
\bigoplus_{a,b,c=0}^{p-1}} \, \Bbbk[\h] \, \h^{- a - b - c \, p}
\cdot {(\nu-1)}^a \, {(\tau-1)}^b \, {\big( \tau^p - 1 \big)}^c $ }
 \vskip1pt
\noindent
as  $ \Bbbk[\h] $--modules,  since  $ \, d(\nu) = 1 = d(\tau) \, $
and  $ \, d\big(\tau^p\big)) = p \, $,  \, with  $ \, \Delta(\chi_g)
= \chi_g \otimes 1 + 1 \otimes \chi_g + \h^{d(g)} \, \chi_g \otimes
\chi_g \, $  for all  $ \, g \in B \, $.  As a direct consequence
we have also
 \vskip1pt
   \centerline{ $ {\textstyle \bigoplus_{a,b,c=0}^{p-1}} \, \Bbbk
\cdot \overline{\chi_\nu}^{\;a} \, \overline{\chi_\tau}^{\;b} \,
\overline{\chi_{\tau^p}}^{\;c} \; = \; {{\Bbbk[G]}_\h}^{\!\vee}
{\Big|}_{\h=0} \; \cong \; \widehat{\Bbbk[G]} \; = \; {\textstyle
\bigoplus_{a,b,c=0}^{p-1}} \, \Bbbk \cdot \overline{\eta_\nu}^{\;a} \,
\overline{\eta_\tau}^{\;b} \, \overline{\eta_{\tau^p}}^{\;c} \, . $ }
 \vskip1pt
   The two relations  $ \, \nu^p = 1 \, $  and  $ \, \tau^{p^2} =
1 \, $  within  $ G $  yield trivial relations inside  $ \Bbbk[G] $
and  $ {\Bbbk[G]}_\h \, $;  instead, the relation  $ \, \nu \, \tau
\, \nu^{-1} = \tau^{1+p} \, $  turns into  $ \, [\eta_\nu,\eta_\tau]
= \eta_{\tau^p} \cdot \tau \, \nu \, $,  \, which gives  $ \, [\chi_\nu,
\chi_\tau] = \h^{p-2} \, \chi_{\tau^p} \cdot \tau \, \nu \, $  in
$ {{\Bbbk[G]}_\h}^{\!\vee} $.  Therefore  $ \; [\, \overline{\chi_\nu}
\, , \, \overline{\chi_\tau}\,] = \delta_{p,2} \,
\overline{\chi_{\tau^p}} \, $.  Since  $ \, [\,
\overline{\chi_\tau} \, , \, \overline{\chi_{\tau^p}}\,] = 0 =
[\,\overline{\chi_\nu} \, , \, \overline{\chi_{\tau^p}} \,] \, $
(because  $ \, \nu \, \tau^p \, \nu^{-1} = {\big( \tau^{1+p}
\big)}^p = \tau^{p + p^2} = \tau^p \, $)  and  $ \, \{
\overline{\chi_\nu} \, , \, \overline{\chi_\tau} \, , \,
\overline{\chi_{\tau^p}} \,\} \, $  is a  $ \Bbbk $--basis  of
$ \, \gerk_G = \Cal{L}_p(G) \, $,  \, we conclude that the latter
has trivial or non-trivial Lie bracket according to whether  $ \,
p \not= 2 \, $  or  $ \, p = 2 \, $.  In addition, we have the
relations  $ \, \chi_\nu^{\;p} = 0 \, $,  $ \, \chi_{\tau^p}^{\;p}
= 0 \, $  and  $ \, \chi_\tau^{\;p} = \chi_{\tau^p} \, $:  \, these
give analogous relations in  $ \, {{\Bbbk[G]}_\h}^{\!\vee}{\Big|}_{\h=0}
\, $,  \, which read as formulas for the  $ p $--operation  of
$ \gerk_G $,  namely  $ \; \overline{\chi_\nu}^{\;[p\,]} = 0
\, $,  $ \; \overline{\chi_{\tau^p}}^{\;[p\,]} = 0 \, $,
$ \; \overline{\chi_\tau}^{\;[p\,]} = \chi_{\tau^p} \, $.
                                     \par
   To sum up, we have a complete presentation for
$ {{\Bbbk[G]}_\h}^{\!\vee} $  by generators and relations, i.e.
 \vskip3pt
  \centerline{ $ {{\Bbbk[G]}_\h}^{\!\vee}  \;\; \cong \;\;
\Bbbk[\h] \, \big\langle v_1, v_2, v_3 \big\rangle \bigg/
\! \bigg( \hskip-9pt
%
%
\hbox{ $ \matrix
 &  v_1 \, v_2 - v_2 \, v_1 - \h^{p-2} \, v_3 \,
(1 + \h \, v_\tau) \, (1 + \h \, v_\nu)  \\
 &  v_1 \, v_3 - v_3 \, v_1 \, ,  \quad  v_1^{\,p} \, ,
\quad  v_2^{\,p} - v_3 \, , \quad  v_3^{\,p} \, , \quad
v_2 \, v_3 - v_3 \, v_2
       \endmatrix $ }  \bigg) $ }
 \vskip3pt
\noindent
via  $ \, \chi_\nu \mapsto v_1 \, $,  $ \, \chi_\tau \mapsto v_2 \, $,
$ \, \chi_{\tau^p} \mapsto v_3 \, $.  Similarly (as a consequence) we
have the presentation
 \vskip3pt
  \centerline{ $ \widehat{\Bbbk[G]}  \, = \,
{{\Bbbk[G]}_\h}^{\!\vee}{\Big|}_{\h=0} \; \cong
\;\;  \Bbbk \, \big\langle y_1,y_2,y_3 \big\rangle
\bigg/ \!
%
%
\bigg( \hbox{ $ \matrix
     y_1 \, y_2 - y_2 \, y_1 - \delta_{p,2} \, y_3 \, ,
\qquad  y_2^{\,p} - y_3  \\
     y_1 \, y_3 - y_3 \, y_1 \, ,  \quad  y_1^{\,p} \, ,
\quad  y_3^{\,p} \, ,  \quad  y_2 \, y_3 - y_3 \, y_2
       \endmatrix $ }  \bigg) $ }
 \vskip3pt
\noindent
via  $ \, \overline{\chi_\nu} \mapsto y_1 \, $,  $ \,
\overline{\chi_\tau} \mapsto y_2 \, $,  $ \, \overline{\chi_{\tau^p}}
\mapsto y_3 \, $,  \, with  $ p $--operation  as above and the
           $ y_i $'s  being primitive\break
{\it  $ \underline{\text{Remark}} $:}  \, if  $ \, p \not= 2 \, $
exactly the same result holds for  $ \, G = \Z_p \times \Z_{p^2} \, $,
\, i.e.~$ \; \gerk_{\, \Z_p \hskip-0,5pt \ltimes \Z_{p^2}} =
\gerk_{\, \Z_p \hskip-0,5pt \times \Z_{p^2}} \; $:  \; this
shows that the restricted Lie bialgebra  $ \gerk_G $  may be
     \hbox{not enough to recover the group  $ G \, $.}
                                          \par

   As for  $ \, {\big( {{\Bbbk[G]}_\h}^{\!\vee} \big)}' $,  \, it is
generated by  $ \, \psi_\nu = \nu - 1 $,  $ \, \psi_\tau = \tau - 1 $,
$ \, \psi_{\tau^p} = \h^{1-p} \big(\tau^p - 1 \big) \, $,  \, with
relations  $ \; \psi_\nu^{\;p} = 0 \, $,  $ \; \psi_\tau^{\;p}
= \h^{p-1} \psi_{\tau^p} \, $,  $ \; \psi_{\tau^p}^{\;\;p} =
0 \, $,  $ \; \psi_\nu \, \psi_\tau - \psi_\tau \, \psi_\nu =
\h^{\,p-1} \psi_{\tau^p} \, (1 + \psi_\tau) \, (1 + \psi_\nu) \, $,
$ \; \psi_\tau \, \psi_{\tau^p} - \psi_{\tau^p} \, \psi_\tau = 0 \, $,
and  $ \; \psi_\nu \, \psi_{\tau^p} - \psi_{\tau^p} \, \psi_\nu = 0
\, $.  In particular  $ \; {\big( {{\Bbbk[G]}_\h}^{\!\vee} \big)}'
\supsetneqq {\Bbbk[G]}_\h \, $,  \, and
 \vskip3pt
  \centerline{ $ {\big( {{\Bbbk[G]}_\h}^{\!\vee} \big)}'  \;\; \cong
\;\;  \Bbbk[\h] \, \big\langle u_1, u_2, u_3 \big\rangle \bigg/ \!
\bigg( \hskip-13pt
\hbox{ $ \matrix
   u_1 \, u_3 - u_3 \, u_1 \, ,  \quad  u_2^{\;p} - \h^{p-1} u_3
\, ,  \quad  u_2 \, u_3 - u_3 \, u_2  \\
   \quad  u_1^{\;p} \, ,  \;\;\,  u_1 \, u_2 - u_2 \,
u_1 - \h^{\,p-1} u_3 \, (1 + u_2) \, (1 + u_1) \, , \;\;  u_3^{\;p}
       \endmatrix $ }  \hskip-4pt \bigg) $ }
 \vskip7pt
\noindent
via  $ \, \psi_\nu \mapsto u_1 \, $,  $ \, \psi_\tau \mapsto u_2
\, $,  $ \, \psi_{\tau^p} \mapsto u_3 \, $.  Letting  $ \; z_1 :=
\psi_\nu{\big|}_{\h=0} \! + 1 \, $,  $ \; z_2 := \psi_\tau{\big|}_{\h=0}
\! + 1 \; $  and  $ \; x_3 := \psi_{\tau^p}{\big|}_{\h=0} \; $  this gives
$ \, {\big( {{\Bbbk[G]}_\h}^{\!\vee} \big)}'{\Big|}_{\h=0} \!\! = \Bbbk
\big[z_1,z_2,x_3\big] \Big/ \big( z_1^{\,p} \! - \! 1, z_2^{\,p} \! -
\! 1, x_3^{\,p} \,\big) \, $  as a  $ \Bbbk $--algebra,  with the
$ z_i $'s  group-like,  $ x_3 $  primitive  (cf.~Theorem D{\it (b)\/}),
and Poisson bracket given by  $ \, \big\{ z_1, z_2 \big\} = \delta_{p,2}
\, z_1 \, z_2 \, x_3 \, $,  $ \, \big\{ z_2, x_3 \big\} = 0 \, $  and
$ \, \big\{z_1, x_3\big\} = 0 \, $.  Thus  $ \, {\big(
{{\Bbbk[G]}_\h}^{\!\vee} \big)}'{\Big|}_{\h=0} \! = F[\varGamma_G]
\, $  with  $ \, \varGamma_G \cong {\boldsymbol\mu}_p \times
{\boldsymbol\mu}_p \times {\boldsymbol\alpha}_p \, $  as algebraic
groups, with Poisson structure such that  $ \,\text{\sl coLie}\,
(\varGamma_G) \cong \gerk_G \, $.
                                         \par
   Since  $ \, G_\infty = \{1\} \, $  the general theory ensures that
$ \, {A_\Bbbk(G)}' = A_\Bbbk(G) \, $.  We leave to the interested
reader the task of computing the filtration  $ \underline{D} $  of
$ A_\Bbbk(G) $,  and consequently describe  $ \, {A_R(G)}' \, $,
$ \, \big({A_R(G)}'\big)^{\!\vee} \, $,  $ \, \widetilde{A_\Bbbk(G)}
\, $  and the connected Poisson group  $ \, K_G := \text{\it Spec}\,
\big( \widetilde{A_\Bbbk(G)} \big) \, $.
 \vskip4pt
   {\sl (3) \, An Abelian infinite group.} \,  Let  $ \, G = \Z^n
\, $  (written multiplicatively, with generators  $ \, e_1, \dots,
e_n \, $),  then  $ \, \Bbbk[G] = \Bbbk[\Z^n] = \Bbbk \big[ e_1^{\pm 1},
\dots, e_n^{\pm 1} \big] \, $  (the ring of Laurent polynomials).
This is the function algebra of the algebraic group
$ {\Bbb{G}_m}^{\hskip-3pt n} $,  i.e.~the  $ n $--dimensional
torus on  $ \Bbbk $  (which is exactly the character group of
$ \Z^n $),  thus we get back to the function algebra case.

 \vskip1,7truecm

\centerline {\bf \S \; 6 \  First example: the Kostant-Kirillov
structure }

\vskip10pt

  {\bf 6.1 Classical and quantum setting.} \, We study now another
quantization of the Kostant-Kirillov structure.  Let  $ \gerg $
and  $ \gerg^\star $  be as in \S 5.7, consider  $ \gerg $  as
             a Lie bialgebra with  
\noindent   trivial Lie cobracket and look at  $ \gerg^\star $
as its dual Poisson group, hence its Poisson structure is exactly
the Kostant-Kirillov one.
                                          \par
   Take as ground ring  $ \, R := \Bbbk[\nu] \, $  (a PID): we shall
consider the primes  $ \, \h = \nu \, $  and  $ \, \h = \nu - 1 \, $, 
\, and we'll find quantum groups at either of them for both  $ \gerg $ 
and  $ \gerg^\star $.
                                          \par
   To begin with, we assume  $ \, \hbox{\it Char}\,(\Bbbk)
= 0 \, $,  \, and postpone to \S 6.4 the case  $ \, \hbox{\it
Char}\,(\Bbbk) > 0 \, $.
                                  \par
   Let  $ \, \gerg_\nu := \gerg[\nu] = \Bbbk[\nu] \otimes_\Bbbk
\gerg \, $,  \, endow it with the unique  $ \Bbbk[\nu] $--linear
Lie bracket  $ \, {[\ ,\ ]}_\nu \, $  given by  $ \, {[x,y]}_\nu
:= \nu \, [x,y] \, $  for all  $ \, x $,  $ y \in \gerg \, $,
\, and define
  $$  H := U_{\Bbbk[\nu]}(\gerg_\nu) = T_{\Bbbk[\nu]}(\gerg_\nu)
\Big/ \big(\big\{\, x \cdot y - y \cdot x - \nu \, [x,y] \;\big|\;
x, y \in \gerg \,\big\} \big) \, ,  $$
the universal enveloping algebra of the Lie  $ \Bbbk[\nu] $--algebra
$ \, \gerg_\nu \, $,  \, endowed with its natural structure of Hopf
algebra.  Then  $ H $  is a free  $ \Bbbk[\nu] $--algebra,
so that  $ \, H \in \HA \, $  and  $ \, H_F := \Bbbk(\nu)
\otimes_{\Bbbk[\nu]} H \in \HA_F \, $  (see \S 1.3); its
specializations at  $ \, \nu = 1 \, $  and at  $ \, \nu
= 0 \, $  are
  $$  \eqalign{
   H \Big/ (\nu\!-\!1) \, H \, = \, U(\gerg)  \qquad  &
\text{as a  {\sl co-Poisson}  Hopf algebra} \, ,  \cr
   H \Big/ \nu \, H \, = \, S(\gerg) \, = \, F[\gerg^\star]
\qquad  &  \text{as a  {\sl Poisson}  Hopf algebra} \, ;
\cr }  $$
in a more suggesting way, we can also express this with notation
like  $ \; H \,{\buildrel {\, \nu \rightarrow 1 \,} \over
\llongrightarrow}\, U(\gerg) \, $,  $ \; H \,{\buildrel {\,
\nu \rightarrow 0 \,} \over \llongrightarrow}\, F[\gerg^\star]
\, $.  Therefore,  {\sl  $ H $  is a QrUEA at  $ \, \h := (\nu
\! - \! 1) \, $  and a QFA at  $ \, \h := \nu \, $};  thus now
we go and consider Drinfeld's functors for  $ H $  at
$ (\nu\!-\!1) $  and at  $ (\nu) $.

\vskip7pt

  {\bf 6.2 Drinfeld's functors at  $ (\nu) $.} \, Let  $ \;
{(\ )}^{\vee_{\!(\nu)}} \, \colon \, \HA \longrightarrow \HA \; $
and  $ \; {(\ )}^{\prime_{(\nu)}} \, \colon \, \HA \longrightarrow
\HA \; $  be the Drinfeld's functors at  $ \, (\nu) \, \big( \in
\text{\it Spec}\big(\Bbbk[\nu]\big) \, \big) \, $.  By definitions
$ \, J := \text{\sl Ker} \, \big( \epsilon \, \colon \, H
\longrightarrow \Bbbk[\nu] \big) \, $  is nothing but the 2-sided
ideal of  $ \, H := U(\gerg_\nu) \, $  generated by  $ \gerg_\nu $
itself; thus  $ \, H^{\vee_{\!(\nu)}} $,  which by definition is
the unital  $ \Bbbk[\nu] $--subalgebra  of  $ H_F $  generated by
$ \, J^{\vee_{\!(\nu)}} := \nu^{-1} J \, $,  \, is just the unital
$ \Bbbk[\nu] $--subalgebra  of  $ H_F $  generated by  $ \;
{\gerg_\nu}^{\!\vee_{\!(\nu)}} := \nu^{-1} \, \gerg_\nu \, $.
Now consider the  $ \Bbbk[\nu] $--module  isomorphism  $ \;
{(\ )}^{\!\vee_{\!(\nu)}} \, \colon \, \gerg_\nu \,{\buildrel
\cong \over \longrightarrow}\, {\gerg_\nu}^{\!\vee_{\!(\nu)}}
:= \nu^{-1} \, \gerg_\nu \; $  given by  $ \, z \mapsto z^\vee
:= \nu^{-1} z \in {\gerg_\nu}^{\!\vee_{\!(\nu)}} \, $  for all
$ z \in \gerg_\nu \, $;  consider on  $ \, \gerg_\nu := \Bbbk[\nu]
\otimes_\Bbbk \gerg \, $  the natural Lie algebra structure (with
trivial Lie cobracket), given by scalar extension from  $ \gerg
\, $,  and push it over  $ {\gerg_\nu}^{\!\vee_{\!(\nu)}} $ 
via  $ {(\ )}^{\!\vee_{\!(\nu)}} $,  so that  $ {\gerg_\nu}^{\!
\vee_{\!(\nu)}} $  is isomorphic to  $ \, \gerg_\nu^{\text{\it
nat}} \, $  (i.e.~$ \gerg_\nu $  carrying the natural Lie bialgebra
structure) as a Lie bialgebra.  Consider  $ \, x^\vee $,  $ y^\vee
\in {\gerg_\nu}^{\!\vee_{\!(\nu)}} $  (with  $ \, x $,  $ y \in
\gerg_\nu \, $):  \, then  $ \; H^{\vee_{\!(\nu)}} \ni \big(
x^\vee \, y^\vee - y^\vee \, x^\vee \big) = \nu^{-2} \big( x
\, y - y \, x \big) = \nu^{-2} \, {[x,y]}_\nu = \nu^{-2} \, \nu
\, [x,y] = \nu^{-1} \, [x,y] = {[x,y]}^\vee =: \big[ x^\vee,
y^\vee \big] \in {\gerg_\nu}^{\!\vee_{\!(\nu)}} \, $.  Therefore
we can conclude at once that  $ \; H^{\vee_{\!(\nu)}} = U
\big( {\gerg_\nu}^{\!\vee_{\!(\nu)}} \big) \cong U \big(
\gerg_\nu^{\text{\it nat}} \big) \; $.
                                        \par
   As a first consequence,  $ \, {\big( H^{\vee_{\!(\nu)}}
\big)}{\Big|}_{\nu=0} \cong U \big( \gerg_\nu^{\text{\it nat}}
\big) \Big/ \nu \, U \big( \gerg_\nu^{\text{\it nat}} \big)
= U \Big( \gerg_\nu^{\text{\it nat}} \big/ \nu \,
\gerg_\nu^{\text{\it nat}} \Big) = U(\gerg) \, $,  \; that is
$ \; H^{\vee_{\!(\nu)}} \,{\buildrel {\, \nu \rightarrow 0 \,}
\over \llongrightarrow}\, U(\gerg) \, $,  \, thus agreeing with
the second half of  Theorem 2.2{\it (c)}.
                                        \par
   Second, look at  $ \, {\big( H^{\vee_{\!(\nu)}}
\big)}^{\prime_{(\nu)}} $.  Since  $ \, H^{\vee_{\!(\nu)}}
= U \big( {\gerg_\nu}^{\!\vee_{\!(\nu)}} \big) \, $,  \, and
$ \, \delta_n(\eta) = 0 \, $  for all  $ \, \eta \in U \big(
{\gerg_\nu}^{\!\vee_{\!(\nu)}} \big) \, $  such that  $ \,
\partial(\eta) < n \, $  (cf.~the proof of  Lemma 4.2{\it
(d)\/}),  it is easy to see that
  $$  {\big( H^{\vee_{\!(\nu)}} \big)}^{\prime_{(\nu)}} =
\big\langle \nu \, {\gerg_\nu}^{\!\vee_{\!(\nu)}} \big\rangle
= \big\langle \nu \, \nu^{-1} \gerg_\nu \big\rangle =
U(\gerg_\nu) = H  $$
(hereafter  $ \langle S \, \rangle $  is the subalgebra
generated by  $ S \, $),  so  $ \, {\big( H^{\vee_{\!(\nu)}}
\big)}^{\prime_{(\nu)}} = H \, $,  \, which agrees with
Theorem 2.2{\it (b)}.
%
%
   Finally, proceeding as in \S 5.7 we see that  $ \,
H^{\prime_{(\nu)}} = U(\nu\,\gerg_\nu) \, $,  \, whence  $ \,
{\big( H^{\prime_{(\nu)}} \big)}{\Big|}_{\nu=0} \!\! = {\big( U
(\nu \, \gerg_\nu) \big)}{\Big|}_{\nu=0} \!\! \cong S({\gerg}_{ab})
= F \big[ \gerg^\star_{\delta-{\text{\it ab}}} \big] \, $  where
$ \, {\gerg}_{ab} \, $,  resp.~$ \, \gerg^\star_{\delta-{\text{\it
ab}}} \, $,  is simply  $ \gerg $,  resp.~$ \gerg^\star $,  endowed
with the trivial Lie bracket, resp.~cobracket,  so that  $ \,
{\big( H^{\prime_{(\nu)}} \big)} {\Big|}_{\nu=0} \! \cong
S({\gerg}_{ab}) = F \big[ \gerg^\star_{\delta-{\text{\it ab}}}
\big] \, $  has trivial Poisson bracket.  Similarly, we can
iterate this procedure and find that all further images  $ \,
\Big( \cdots {\big( {(H)}^{\prime_{(\nu)}} \big)}^{\prime_{(\nu)}}
\cdots \Big)^{\prime_{(\nu)}} \, $  of the functor
$ {(\ )}^{\prime_{(\nu)}} $  applied many times
to  $ H $  are all isomorphic, hence they all have
the same specialization at  $ (\nu) $,  namely
  $$  {\bigg( \! {\Big( \cdots {\big( {(H)}^{\prime_{(\nu)}}
\big)}^{\prime_{(\nu)}} \cdots \Big)}^{\prime_{(\nu)}}
\bigg)}{\bigg|}_{\nu=0} \, \cong \;\, S({\gerg}_{ab}) =
F \big[ \gerg^\star_{\delta-{\text{\it ab}}} \big] \; .  $$

\vskip7pt

  {\bf 6.3 Drinfeld's functors at  $ (\nu\!-\!1) $.} \, Now
consider  $ \, (\nu\!-\!1) \, \big( \! \in \! \text{\it
Spec}\big(\Bbbk[\nu]\big) \, \big) \, $,  \, and let  $ \;
{(\ )}^{\vee_{\!(\nu\!-\!1)}} \, \colon \, \HA \longrightarrow
\HA \; $  and  $ \; {(\ )}^{\prime_{(\nu\!-\!1)}} \, \colon \,
\HA \longrightarrow \HA \; $  be the corresponding Drinfeld's
functors.  Set  $ \, {\gerg_\nu}^{\prime_{(\nu\!-\!1)}} :=
(\nu\!-\!1) \, \gerg_\nu \, $,  \, let  $ \; \, \colon \,
{\gerg_\nu} \,{\buildrel \cong \over \longrightarrow}\,
{\gerg_\nu}^{\!\prime_{(\nu\!-\!1)}} := (\nu\!-\!1)
\, \gerg_\nu \; $  be the  $ \Bbbk[\nu] $--module  isomorphism
given by  $ \, z \mapsto z' := (\nu\!-\!1) \, z \in
{\gerg_\nu}^{\!\prime_{(\nu\!-\!1)}} \, $  for all
$ \, z \in \gerg_\nu \, $,  \, and push over via it the
Lie bialgebra structure of  $ \gerg_\nu $  to an isomorphic
Lie bialgebra structure on  $ \, {\gerg_\nu}^{\!\prime_{(\nu\!-\!1)}} \, $,  \,
whose Lie bracket will be denoted by
$ \, {[\ ,\ ]}_\ast \, $.  Notice then that we have Lie bialgebra
isomorphisms  $ \; \gerg \, \cong \, \gerg_\nu \big/ (\nu\!-\!1)
\, \gerg_\nu \, \cong \, {\gerg_\nu}^{\!\prime_{(\nu\!-\!1)}}
\big/ (\nu\!-\!1) \, {\gerg_\nu}^{\!\prime_{(\nu\!-\!1)}} \, $.
                                              \par
   Since  $ \, H := U(\gerg_\nu) \, $  it is easy to see by direct
computation that
  $$  H^{\prime_{(\nu\!-\!1)}} = \big\langle (\nu\!-\!1) \,
\gerg_\nu \big\rangle = U \big( {\gerg_\nu}^{\!\prime_{(\nu
\! - \! 1)}} \big) \, ,   \eqno (6.1)  $$
where  $ {\gerg_\nu}^{\!\prime_{(\nu\!-\!1)}} $  is considered as
a Lie  $ \Bbbk[\nu] $--subalgebra  of  $ \, \gerg_\nu \, $.  Now,
if  $ \, x' $,  $ y' \in {\gerg_\nu}^{\prime_{(\nu\!-\!1)}} \, $
(with  $ \, x $,  $ y \in \gerg_\nu \, $),  \, we have
  $$  x' \, y' - y' \, x' = {(\nu\!-\!1)}^2 \big( x \, y -
y \, x \big) = {(\nu\!-\!1)}^2 \, {[x,y]}_\nu = (\nu\!-\!1) \,
{{[x,y]}_\nu}^{\!\prime} = (\nu\!-\!1) \, {\big[x',y'\big]}_\ast
\, .   \eqno (6.2)  $$
   \indent   This and (6.1) show at once that  $ \, {\big(
H^{\prime_{(\nu\!-\!1)}} \big)}{\Big|}_{(\nu\!-\!1)=0} \! = S
\Big( {\gerg_\nu}^{\prime_{(\nu\!-\!1)}} \big/ (\nu\!-\!1) \,
{\gerg_\nu}^{\prime_{(\nu\!-\!1)}} \Big) \, $  as Hopf algebras,
and also as  {\sl Poisson}  algebras: indeed, the latter holds
because the Poisson bracket  $ \, \{\ ,\ \} \, $  of  $ \, S
\Big( {\gerg_\nu}^{\prime_{(\nu\!-\!1)}} \big/ (\nu\!-\!1) \,
{\gerg_\nu}^{\prime_{(\nu\!-\!1)}} \Big) \, $  inherited from
$ H^{\prime_{(\nu\!-\!1)}} $  (by specialization) is uniquely
determined by its restriction to  $ \, {\gerg_\nu}^{\prime_{(\nu
\! - \! 1)}} \big/ (\nu\!-\!1) \, {\gerg_\nu}^{\prime_{(\nu \! -
\! 1)}} \, $,  \, and on the latter space we have  $ \; \{\ ,\ \}
= {[\ ,\ ]}_\ast \, \mod (\nu\!-\!1) \, {\gerg_\nu}^{\prime_{(\nu
\! - \! 1)}} \; $  (by (6.2)).  Finally, since  $ \,
{\gerg_\nu}^{\prime_{(\nu\!-\!1)}} \big/ (\nu\!-\!1) \,
{\gerg_\nu}^{\prime_{(\nu\!-\!1)}} \cong \gerg \, $  as Lie
algebras we have  $ \, {\big( H^{\prime_{(\nu\!-\!1)}} \big)}
{\Big|}_{(\nu\!-\!1)=0} \! = S(\gerg) = F[\gerg^\star] \, $  as
Poisson Hopf algebras, or, in short,  $ \; H^{\prime_{(\nu
\! - \! 1)}} \,{\buildrel {\, \nu \rightarrow 1 \,} \over
\llongrightarrow}\, F[\gerg^\star] \, $,  \, as prescribed by
the ``first half\/'' of  Theorem 2.2{\it (c)}.
                                        \par
   Second, look at  $ \, {\big( H^{\prime_{(\nu\!-\!1)}}
\big)}^{\vee_{\!(\nu\!-\!1)}} $.  Since  $ \, H^{\prime_{(\nu
\! - \! 1)}} = U \big( {\gerg_\nu}^{\!\prime_{(\nu\!-\!1)}} \big)
\, $,  \, we have that  $ \, J^{\prime_{(\nu\!-\!1)}} := \text{\sl
Ker} \, \big( \epsilon \, \colon \, H^{\prime_{(\nu\!-\!1)}}
\longrightarrow \Bbbk[\nu] \big) \, $  is nothing but the
2-sided ideal of  $ \, H^{\prime_{(\nu\!-\!1)}} = U \big(
{\gerg_\nu}^{\!\prime_{(\nu\!-\!1)}} \big) \, $  generated by
$ {\gerg_\nu}^{\!\prime_{(\nu\!-\!1)}} \, $;  thus  $ \, {\big(
H^{\prime_{(\nu\!-\!1)}} \big)}^{\vee_{\!(\nu\!-\!1)}} $,
\, generated by  $ \, {\big( J^{\prime_{(\nu\!-\!1)}}
\big)}^{\vee_{\!(\nu\!-\!1)}} := {(\nu\!-\!1)}^{-1}
J^{\prime_{(\nu\!-\!1)}} \, $ as a unital
$ \Bbbk[\nu] $--subalgebra  of  $ \, {\big(
H^{\prime_{(\nu\!-\!1)}} \big)}_F = H_F \, $,  \,
is just the unital  $ \Bbbk[\nu] $--subalgebra
of  $ H_F $  generated by  $ \; {(\nu\!-\!1)}^{-1}
{\gerg_\nu}^{\!\prime_{(\nu\!-\!1)}} = {(\nu\!-\!1)}^{-1}
(\nu\!-\!1) \, \gerg_\nu = \gerg_\nu \, $,  \, that is to say
$ \, {\big( H^{\prime_{(\nu\!-\!1)}} \big)}^{\vee_{\!(\nu\!-\!1)}}
= U(\gerg_\nu) = H \, $,  \, according to  Theorem 2.2{\it (b)}.
                                        \par
   Finally, for  $ H^{\vee_{(\nu\!-\!1)}} $  one has essentially
the same feature as in \S 5.7, and the analysis therein can be
applied again; the final result then will depend on the nature
of  $ \gerg $,  in particular on its lower central series.
%
%
 \eject   
  {\bf 6.4 The case of positive characteristic.} \, Let us
consider now a field  $ \Bbbk $  such that  $ \, \hbox{\it
Char}\,(\Bbbk) = p > 0 \, $.  Starting from  $ \gerg $  and
$ \, R := \Bbbk[\nu] \, $  as in \S 6.1, define  $ \gerg_\nu $
like therein, and consider  $ \, H := U_{\Bbbk[\nu]}(\gerg_\nu)
= U_R(\gerg_\nu) \, $.  Then we have again
  $$  \eqalign{
   H \Big/ (\nu\!-\!1) \, H \, = \, U(\gerg) \, = \, \u \Big(
\gerg^{{[p\hskip0,7pt]}^\infty} \Big)  \qquad  &  \text{as a  {\sl
co-Poisson}  Hopf algebra} \, ,  \cr
   H \Big/ \nu \, H \, = \, S(\gerg) \, = \, F[\gerg^\star]
\qquad  &  \text{as a  {\sl Poisson}  Hopf algebra}  \cr }  $$
so that  $ H $  is a QrUEA at  $ \, \h := (\nu\!-\!1) \, $  (for
the restricted universal enveloping algebra  $  \u \Big(
\gerg^{{[p\hskip0,7pt]}^\infty} \Big) \, $)  and is a QFA at  $ \,
\h := \nu \, $  (for the function algebra  $ F[\gerg^\star]
\, $);  so now we study Drinfeld's functors for  $ H $  at
$ (\nu\!-\!1) $  and at  $ (\nu) $.
                                          \par
   Exactly the same procedure as before shows again that
$ \; H^{\vee_{\!(\nu)}} = U \big( {\gerg_\nu}^{\!\vee_{\!(\nu)}}
\big) \, $,  \, from which it follows that  $ \, {\big(
H^{\vee_{\!(\nu)}} \big)}{\Big|}_{\nu=0} \! \cong U(\gerg)
\, $,  \, i.e.~in short  $ \; H^{\vee_{\!(\nu)}} \,{\buildrel
{\, \nu \rightarrow 0 \,} \over \llongrightarrow}\, U(\gerg)
\, $,  \, which is a result quite ``parallel'' to the second
half of  Theorem 2.2{\it (c)}.
                                        \par
   Changes occur when looking at  $ \, {\big( H^{\vee_{\!(\nu)}}
\big)}^{\prime_{(\nu)}} $:  since  $ \, H^{\vee_{\!(\nu)}} =
U \big( {\gerg_\nu}^{\!\vee_{\!(\nu)}} \big) = \u \Big( {\big(
{\gerg_\nu}^{\!\vee_{\!(\nu)}} \big)}^{{[p\hskip0,7pt]}^\infty} \Big) \, $
we have  $ \, \delta_n(\eta) = 0 \, $  for all  $ \, \eta \in \u
\Big( \! {\big( {\gerg_\nu}^{\!\vee_{\!(\nu)}} \big)}^{{[p\hskip0,7pt]}^\infty}
\Big) \, $  such that  $ \, \partial(\eta) < n \, $  w.r.t.~the
standard filtration of  $ \, \u \Big( {\big( {\gerg_\nu}^{\!
\vee_{\!(\nu)}} \big)}^{{[p\hskip0,7pt]}^\infty} \Big) \, $  (cf.~the
proof of  Lemma 4.2{\it (d)\/},  which clearly adapts to
the present situation): this implies
  $$  {\big( H^{\vee_{\!(\nu)}} \big)}^{\prime_{(\nu)}}
= \Big\langle \nu \cdot {\big( {\gerg_\nu}^{\! \vee_{\!(\nu)}}
\big)}^{{[p\hskip0,7pt]}^\infty} \Big\rangle  \qquad  \Big( \subset
\u \Big( \nu \cdot {\big( {\gerg_\nu}^{\!\vee_{\!(\nu)}}
\big)}^{{[p\hskip0,7pt]}^\infty} \Big) \; \Big)  $$
which is strictly  {\sl bigger\/}  than  $ H $,  because
$ \; \Big\langle \nu \cdot {\big( {\gerg_\nu}^{\! \vee_{\!(\nu)}}
\big)}^{{[p\hskip0,7pt]}^\infty} \Big\rangle = \Big\langle \sum\limits_{n
\geq 0} \nu \cdot {\big( {\gerg_\nu}^{\! \vee_{\!(\nu)}}
\big)}^{{[p\hskip0,7pt]}^n}
               \Big\rangle = $\break
        $ = \Big\langle \gerg_\nu + \nu^{1-p} \, \big\{\, x^p
\,\big|\, x \in \gerg_\nu \big\} + \nu^{1-p^2} \Big\{\, x^{p^2}
\,\Big|\, x \in \gerg_\nu \Big\} + \cdots \Big\rangle \,
\supsetneqq \, U(\gerg_\nu) = H \, $.
                                        \par
   Finally, proceeding as above it is easy to see that  $ \,
H^{\prime_{(\nu)}} = \Big\langle \nu \, P \big( U(\gerg_\nu)
\big) \Big\rangle = \Big\langle \nu \, \gerg^{{[p\hskip0,7pt]}^\infty}
\Big\rangle \, $  whence, letting  $ \, \tilde{\gerg} :=
\nu \, \gerg \, $  and  $ \, \tilde{x} := \nu \, x \, $
for all  $ \, x \in \gerg \, $,  we have
  $$  H^{\prime_{(\nu)}} = T_R (\tilde{\gerg}) \bigg/
\Big( \Big\{\, \tilde{x} \, \tilde{y} - \tilde{y} \, \tilde{x}
- \nu^2 \, \widetilde{[x,y]} \, , \, \tilde{z}^p - \nu^{p-1}
\widetilde{z^{[p\hskip0,7pt]}} \;\Big\vert\;\, x, y, z \in \gerg
\,\Big\} \Big)  $$
so that  $ \; H^{\prime_{(\nu)}}
\;{\buildrel {\nu \rightarrow 0} \over
{\relbar\joinrel\llongrightarrow}}\;
T_\Bbbk (\tilde{\gerg}) \bigg/ \Big( \Big\{\, \tilde{x} \,
\tilde{y} - \tilde{y} \, \tilde{x} \, , \, \tilde{z}^p
\;\Big\vert\; \tilde{x}, \tilde{y}, \tilde{z} \in
\tilde{\gerg} \,\Big\} \Big) =
S_\Bbbk (\gerg_{\text{\it ab}}) \Big/ \big( \big\{\,
z^p \,\big\vert\; z \in \gerg
              \,\big\} \big) = $\break
      $ = F[\gerg^\star_{\delta-\text{\it ab}}] \Big/ \big(
\big\{\, z^p \,\big\vert\; z \in \gerg \,\big\} \big) \, $,
\, that is  $ \; H^{\prime_{(\nu)}}{\Big|}_{\nu=0} \!\! \cong
F[\gerg^\star_{\delta-\text{ab}}] \Big/ \big( \big\{\, z^p
\,\big\vert\; z \in \gerg \,\big\} \big) \; $  {\sl as Poisson
Hopf algebras},  where  $ \gerg_{\text{\it ab}} $  and
$ \gerg^\star_{\delta-\text{\it ab}} $  are as above.
Therefore  $ H^{\prime_{(\nu)}} $  is a QFA (at  $ \h = \nu
\, $)  for a non-reduced algebraic Poisson group of height
1, whose cotangent Lie bialgebra is the vector space
$ \gerg $  with trivial Lie bialgebra structure: this
again yields somehow an analogue of part  {\it (c)\/}  of
Theorem 2.2 for the present case.  If we iterate, we
find that all further images  $ \, \Big( \cdots {\big(
{(H)}^{\prime_{(\nu)}} \big)}^{\prime_{(\nu)}} \cdots
\Big)^{\prime_{(\nu)}} \, $  of the functor
$ {(\ )}^{\prime_{(\nu)}} $  applied to  $ H $  are
all pairwise isomorphic, so
  $$  {{\Big( \cdots {\big( {(H)}^{\prime_{(\nu)}}
\big)}^{\prime_{(\nu)}} \cdots \Big)}^{\prime_{(\nu)}}}
{\bigg|}_{\nu=0} \, \cong \;\, S({\gerg}_{ab}) \Big/
\big( \big\{\, z^p \,\big\vert\; z \in \gerg \,\big\}
\big) = F \big[ \gerg^\star_{\delta-{\text{\it ab}}}
\big] \Big/ \big( \big\{\, z^p \,\big\vert\; z \in
\gerg \,\big\} \big) \; .  $$
                                        \par
   Now for Drinfeld's functors at  $ (\nu\!-\!1) $.  Up to
minor changes, with the same procedure and notations as in
\S 6.3 we get analogous results.  First, an
                     analogue of (6.1), namely\break
   $ \; H^{\prime_{(\nu\!-\!1)}} = \Big\langle (\nu\!-\!1) \cdot
P\big( U(\gerg_\nu) \big) \Big\rangle = \Big\langle (\nu\!-\!1)
\, {( \gerg_\nu )}^{{[p\hskip0,7pt]}^\infty} \Big\rangle = \bigg\langle
{\Big( {(\gerg_\nu)}^{{[p\hskip0,7pt]}^\infty} \Big)}^{\!\prime_{(\nu
\! - \! 1)}} \bigg\rangle \, $,  \, holds and yields
  $$  \displaylines{
   H^{\prime_{(\nu\!-\!1)}} = T_R \left( \! {\Big(
{(\gerg_\nu)}^{{[p\hskip0,7pt]}^\infty} \Big)}^{\!\prime_{(\nu\!-\!1)}}
\right) \Bigg/ \bigg( \Big\{\, x' \, y' - y' \, x' - (\nu-1)
\, {\big[ x', y']}_* \, ,  \, {(x')}^p - {(\nu\!-\!1)}^{p-1}
{\big( x^{[p\hskip0,7pt]} \big)}' \;\Big|\;  \cr
   \hfill   \;\Big|\; x, y \in {(\gerg_\nu)}^{{[p\hskip0,7pt]}^\infty}
\,\Big\} \bigg)  \cr }  $$
and consequently  $ \; H^{\prime_{(\nu\!-\!1)}}{\Big|}_{(\nu \!
- \! 1) = 0} \! \cong S_\Bbbk (\gerg) \Big/ \big( \big\{\, x^p
\;\big|\; x \in \gerg \,\big\} \big) = F[\gerg^\star] \Big/ \big(
\big\{\, x^p \;\big|\; x \in \gerg \,\big\} \big) \, $  as Poisson
Hopf algebras: in short,  $ \; H^{\prime_{(\nu\!-\!1)}} \,{\buildrel
{\, \nu \rightarrow 1 \,} \over \llongrightarrow}\, F[\gerg^\star]
\Big/ \big( \big\{\, x^p \;\big|\; x \in \gerg \,\big\} \big)
\, $.
                                        \par
   Iterating, one finds again that all  $ \, \Big( \cdots {\big(
{(H)}^{\prime_{(\nu)}} \big)}^{\prime_{(\nu\!-\!1)}} \cdots
\Big)^{\prime_{(\nu)}} \, $   \hbox{are pairwise isomorphic, so}
  $$  {{\Big( \! \cdots {\big( {(H)}^{\prime_{(\nu\!-\!1)}}
\big)}^{\prime_{((\nu\!-\!1)}} \! \cdots \Big)}^{\prime_{(\nu
\! - \! 1)}}}{\bigg|}_{(\nu\!-\!1)=0}  \hskip-5pt  \cong \;
S({\gerg}_{ab}) \Big/ \big( \big\{ z^p \,\big\vert\; z \! \in \!
\gerg \big\} \big) = F \big[ \gerg^\star_{\delta-{\text{\it ab}}}
\big] \Big/ \big( \big\{ z^p \,\big\vert\; z \! \in \! \gerg
\,\big\} \big) \; .  $$
                                        \par
   Further on, one has  $ \, {\big( H^{\prime_{(\nu\!-\!1)}}
\big)}^{\vee_{(\nu\!-\!1)}} = {\Big\langle (\nu\!-\!1) \, {(
\gerg_\nu )}^{{[p\hskip0,7pt]}^\infty} \Big\rangle}^{\vee_{(\nu \! - \! 1)}}
= \big\langle {(\nu\!-\!1)}^{-1} \cdot (\nu\!-\!1)
                        \, \gerg_\nu \big\rangle = $\break
$ = \big\langle \gerg_\nu \big\rangle = U_R(\gerg_\nu) =:
H \, $,  \, which perfectly agrees with  Theorem 2.2{\it (b)}.
                                        \par
   Finally, as for  $ H^{\vee_{(\nu\!-\!1)}} $  one has again
the same feature as in \S 5.7: one has to apply the analysis
therein, however, the  $ p $--filtration  in this case is
``harmless'', since it is essentially encoded in the standard
filtration of  $ U(\gerg) $.  In any case the final result
will depend on the properties of the lower central series
of  $ \gerg $.
                                        \par
   Second, we assume in addition that  $ \gerg $  be a
{\sl restricted\/}  Lie algebra, and we consider  $ \,
H := \u_{\Bbbk[\nu]}(\gerg_\nu) = \u_R(\gerg_\nu) \, $.
Then we have
  $$  \eqalign{
   H \Big/ (\nu\!-\!1) \, H \,  &  = \, \u(\gerg)   \hskip80pt
\hfill   \text{as a  {\sl co-Poisson}  Hopf algebra} \, ,  \cr
   H \Big/ \nu \, H \, = \, S(\gerg) \Big/ \big( \big\{\,
z^p \,\big\vert\; z \in \gerg \,\big\} \big) \,  &  = \,
F[\gerg^\star] \Big/ \big( \big\{\, z^p \,\big\vert\;
z \in \gerg \,\big\} \big)   \qquad  \hfill   \text{as
a  {\sl Poisson}  Hopf algebra}  \cr }  $$
which means that  $ H $  is a QrUEA at  $ \, \h := (\nu\!-\!1) \, $
(for  $  \u(\gerg) \, $)  and is a QFA at  $ \, \h := \nu \, $
(for  $ F[\gerg^\star] \Big/ \big( \big\{\, z^p \,\big\vert\;
z \in \gerg \,\big\} \big) \, $);  we go and study Drinfeld's
functors for  $ H $  at  $ (\nu\!-\!1) $  and at  $ (\nu) $.
                                          \par
   As for  $ \, H^{\vee_{\!(\nu)}} \, $,  it depends again
on the  $ p $--operation  of  $ \gerg $,  in short because
the  $ I $--filtration  of  $ \u_\nu(\gerg) $  depends on the
$ p $--filtration   of  $ \gerg $.  In the previous case
--- i.e.~when  $ \, \gerg = {\gerh}^{{[p\hskip0,7pt]}^\infty} \, $  for
some Lie algebra  $ \gerh $  ---  the solution was a plain
one, because the  $ p $--filtration  of  $ \gerg $  is
``encoded'' in the standard filtration of  $ U(\gerh) $;
but the general case will be more complicated, and in
consequence also the situation for  $ \, {\big( H^{\vee_{\!
(\nu)}} \big)}^{\prime_{(\nu)}} $,  since  $ \, H^{\vee_{\!
(\nu)}} \, $  will be different according to the nature of
$ \gerg $.  Instead, proceeding exactly like before one sees
that  $ \, H^{\prime_{(\nu)}} = \Big\langle \nu \, P \big(
u(\gerg_\nu) \big) \Big\rangle = \big\langle \nu \, \gerg
\big\rangle \, $,  whence, letting  $ \, \tilde{\gerg} :=
\nu \, \gerg \, $  and  $ \, \tilde{x} := \nu \, x \, $
for all  $ \, x \in \gerg \, $,  we have
  $$  H^{\prime_{(\nu)}} = T_{\Bbbk[\nu]} (\tilde{\gerg}) \bigg/
\Big( \Big\{\, \tilde{x} \, \tilde{y} - \tilde{y} \, \tilde{x}
- \nu^2 \, \widetilde{[x,y]} \, , \, \tilde{z}^p - \nu^{p-1}
\widetilde{z^{[p\hskip0,7pt]}} \;\Big\vert\;\, x, y, z \in \gerg
\,\Big\} \Big)  $$
so that  $ \; H^{\prime_{(\nu)}}
\;{\buildrel {\nu \rightarrow 0} \over
{\relbar\joinrel\llongrightarrow}}\;
T_\Bbbk (\tilde{\gerg}) \bigg/ \Big( \Big\{\, \tilde{x} \,
\tilde{y} - \tilde{y} \, \tilde{x} \, , \, \tilde{z}^p
\;\Big\vert\; \tilde{x}, \tilde{y}, \tilde{z} \in
\tilde{\gerg} \,\Big\} \Big) =
S_\Bbbk (\gerg_{\text{\it ab}}) \Big/ \big( \big\{\,
z^p \,\big\vert\; z \in \gerg
              \,\big\} \big) = $\break
      $ =F[\gerg^\star_{\delta-\text{\it ab}}] \Big/ \big(
\big\{\, z^p \,\big\vert\; z \in \gerg \,\big\} \big) \, $,
\, that is  $ \; H^{\prime_{(\nu)}}{\Big|}_{\nu=0} \!\! \cong
F[\gerg^\star_{\delta-\text{ab}}] \Big/ \big( \big\{\, z^p
\,\big\vert\; z \in \gerg \,\big\} \big) \; $  {\sl as Poisson
Hopf algebras\/}  (using notation as before).  Thus
$ H^{\prime_{(\nu)}} $  is a QFA (at  $ \h = \nu \, $)
for a non-reduced algebraic Poisson group of height 1,
whose cotangent Lie bialgebra is  $ \gerg $  with the
trivial Lie bialgebra structure: so again we get an
analogue of part of  Theorem 2.2{\it (c)\/}.  Moreover,
iterating again one finds that all  $ \, \Big( \cdots
{\big( {(H)}^{\prime_{(\nu)}} \big)}^{\prime_{( \nu \!
- \! 1 )}} \cdots \Big)^{\prime_{(\nu \! - \! 1)}} \, $
are pairwise isomorphic, so
  $$  {{\Big( \! \cdots {\big( {(H)}^{\prime_{(\nu\!-\!1)}}
\big)}^{\prime_{((\nu\!-\!1)}} \! \cdots \Big)}^{\prime_{(\nu
\! - \! 1)}}}{\bigg|}_{(\nu\!-\!1)=0}  \hskip-5pt  \cong \;
S({\gerg}_{ab}) \Big/ \big( \big\{ z^p \,\big\vert\; z \! \in \!
\gerg \big\} \big) = F \big[ \gerg^\star_{\delta-{\text{\it ab}}}
\big] \Big/ \big( \big\{ z^p \,\big\vert\; z \! \in \! \gerg
\,\big\} \big) \; .  $$
                                        \par
   As for Drinfeld's functors at  $ (\nu\!-\!1) $,  the
situation is more similar to the previous case of  $ \,
H = U_R(\gerg_\nu) \, $.  First  $ \; H^{\prime_{(\nu\!-\!1)}}
= \Big\langle (\nu\!-\!1) \cdot P \big( \u(\gerg_\nu) \big)
\Big\rangle = \big\langle (\nu\!-\!1) \, \gerg_\nu \big\rangle
=: \big\langle {\gerg_\nu}^{\!\prime_{(\nu\!-\!1)}}
\big\rangle \, $,  \, hence
  $$  \displaylines{
   H^{\prime_{(\nu\!-\!1)}} = T_R \Big( {\gerg_\nu}^{\!\prime_{(
\nu \! - \! 1 )}} \Big) \bigg/ \Big( \Big\{\, x' \, y' - y' \, x'
- (\nu-1) \, {\big[ x', y']}_* \, ,  \, {(x')}^p - {(\nu\!-\!1)}^{p-1}
{\big( x^{[p\hskip0,7pt]} \big)}' \;\Big|\;  \cr
   \hfill   \;\Big|\; x, y \in \gerg_\nu \,\Big\} \Big)  \cr }  $$
thus again  $ \; H^{\prime_{(\nu\!-\!1)}}{\Big|}_{(\nu\!-\!1) = 0}
\! \cong S_\Bbbk (\gerg) \Big/ \big( \big\{\, x^p \;\big|\; x \in
\gerg \,\big\} \big) = F[\gerg^\star] \Big/ \big( \big\{\, x^p
\;\big|\; x \in \gerg \,\big\} \big) \, $  as Poisson Hopf
algebras, that is  $ \; H^{\prime_{(\nu\!-\!1)}} \,{\buildrel
{\, \nu \rightarrow 1 \,} \over \llongrightarrow}\, F[\gerg^\star]
\Big/ \big( \big\{\, x^p \;\big|\; x \in \gerg \,\big\} \big) \, $.
Iteration then shows that all  $ \, \Big( \cdots {\big(
{(H)}^{\prime_{(\nu)}} \big)}^{\prime_{(\nu\!-\!1)}} \cdots
\Big)^{\prime_{(\nu)}} \, $  are pairwise isomorphic, so that again
  $$  {{\Big( \! \cdots {\big( {(H)}^{\prime_{(\nu\!-\!1)}}
\big)}^{\prime_{((\nu\!-\!1)}} \! \cdots \Big)}^{\prime_{(\nu
\! - \! 1)}}}{\bigg|}_{(\nu\!-\!1)=0}  \hskip-5pt  \cong \;
S({\gerg}_{ab}) \Big/ \big( \big\{ z^p \,\big\vert\; z \! \in \!
\gerg \big\} \big) = F \big[ \gerg^\star_{\delta-{\text{\it ab}}}
\big] \Big/ \big( \big\{ z^p \,\big\vert\; z \! \in \! \gerg
\,\big\} \big) \; .  $$
                                        \par
   Further, we have  $ \, {\big( H^{\prime_{(\nu\!-\!1)}}
\big)}^{\vee_{(\nu\!-\!1)}} = {\big\langle (\nu\!-\!1) \,
\gerg_\nu \big\rangle}^{\vee_{(\nu\!-\!1)}} = \big\langle
{(\nu\!-\!1)}^{-1} \cdot (\nu\!-\!1) \, \gerg_\nu \big\rangle
= \big\langle \gerg_\nu
                 \big\rangle = $\break
$ = \u_R(\gerg_\nu) =: H \, $,  \, which agrees at all with
Theorem 2.2{\it (b)}.
%
%
   Finally,  $ H^{\vee_{(\nu\!-\!1)}} $  again has the same
feature as in \S 5.7: in particular, in this case the final
result will strongly depend  {\sl both\/}  on the properties
of the lower central series  {\sl and\/}  of the
$ p $--filtration  of  $ \gerg \, $.

 \vskip7pt

  {\bf 6.5 The hyperalgebra case.} \, Let  $ \Bbbk $  be again a field
with  $ \, \hbox{\it Char}\,(\Bbbk) = p > 0 \, $.  Like in \S 5.12, let
$ G $  be an algebraic group (finite-dimensional, for simplicity), and
let  $ \, \hyp(G) := {\big( {F[G]}^\circ \big)}_\epsilon = \big\{\, \phi
\in {F[G]}^\circ \,\big|\, \phi({\germ_e}^{\!n}) = 0 \, , \, \forall \;
n \gg 0 \,\big\} \, $  be the hyperalgebra associated to  $ G \, $.
                                               \par
   For each  $ \, \nu \in \Bbbk \, $,  \, let  $ \, \gerg_\nu := \big(
\gerg, {[\ ,\ ]}_\nu \big) \, $  be the Lie algebra given by  $ \, \gerg
\, $  endowed with the rescaled Lie bracket  $ \, {[\,\ ,\ ]}_\nu := \nu
\, {[\,\ ,\ ]}_\gerg \, $.  By general theory, the algebraic group  $ G $
is uniquely determined by a neighborhood of the identity together with
the formal group law uniquely determined by  $ {[\,\ ,\ ]}_\gerg \, $:
similarly, a neighborhood of the identity of  $ G $  together with
$ {[\,\ ,\ ]}_\nu $  uniquely determines a new connected algebraic
group  $ G_\nu \, $,  whose hyperalgebra  $ \hyp(G_\nu) $  is an
algebraic deformation of  $ \hyp(G) $;  moreover,  $ G_\nu $  is
birationally equivalent to  $ G $,  and for  $ \, \nu \not= 0 \, $
they are also isomorphic as algebraic groups, via an isomorphism
induced by  $ \, \gerg \cong \gerg_\nu \, $,  $ \, x \mapsto \nu^{-1}
x \, $  (however, this may not be the case when  $ \, \nu = 0 \, $).
Note that  $ \hyp(G_0) $  is clearly commutative, because  $ G_0 $  is
Abelian: indeed,
%
%
%
%
we have
  $$  \hyp(G_0)  \; = \;  S_\Bbbk \Big( \gerg^{{(p\hskip0,7pt)^\infty}}
\Big) \hskip-2pt \bigg/ \hskip-2pt \Big( {\big\{\, x^p \,\big\}}_{x \in
\gerg^{{(p\hskip0,7pt)}^\infty}} \Big)  \, = \,  F \Big[ \Big( \gerg^{{(p
\hskip0,7pt)^\infty}} \Big)^{\!\star} \Big] \bigg/ \hskip-2pt \Big(
{\big\{\, y^p \,\big\}}_{y \in \gerg^{{(p\hskip0,7pt)}^\infty}}
\Big)  $$
where  $ \; \gerg^{{(p\hskip0,7pt)}^\infty} \! := \text{\sl Span}\, \Big(
\Big\{\, x^{{(p^n\hskip0pt)}} \,\Big\vert\, x \in \gerg \, , n \in \N
\,\Big\} \Big) \, $;  \, here as usual  $ x^{(n)} $  denotes the
$ n $--th  divided power of  $ \, x \in \gerg \, $  (recall that
$ \hyp(G) $,  hence also  $ \hyp(G_\nu) $,  is generated as an
algebra by all the  $ x^{(n)} $'s,  some of which might be zero).
So  $ \, \hyp(G_0) = F[\varGamma] \, $  where  $ \varGamma $  is a
connected algebraic group of dimension zero and height 1: moreover,
$ \varGamma $  is a Poisson group, with cotangent Lie bialgebra
$ \gerg^{{(p\hskip0,7pt)^\infty}} $  and Poisson bracket induced
by the Lie bracket of  $ \gerg \, $.
                                  \par
%
%
%
%
%
   Now think at  $ \nu $  as a parameter in  $ \, R := \Bbbk[\nu] \, $
(as in \S 6.1), and set  $ \, H := \Bbbk[\nu] \otimes_\Bbbk \hyp(G_\nu)
\, $.  Then we find a situation much similar to that of \S 6.1,
which we shall shortly describe.
                                            \par
   Namely,  $ H $  is a free  $ \Bbbk[\nu] $--algebra,  thus  $ \,
H \in \HA \, $  and  $ \, H_F := \Bbbk(\nu) \otimes_{\Bbbk[\nu]} H
\in \HA_F \, $  (see \S 1.3); its specialization at  $ \, \nu = 1
\, $  is  $ \; H \Big/ (\nu\!-\!1) \, H \, = \, \hyp(G_1) \, = \,
\hyp(G) \, $,  \, and at  $ \, \nu = 0 \, $  is  $ \; H \Big/ \nu
\, H \, = \, \hyp(G_0) \, = \, F[\varGamma] \, $  (as a  {\sl Poisson}
Hopf algebra),  or  $ \; H \,{\buildrel {\, \nu \rightarrow 1 \,} \over
\llongrightarrow}\, \hyp(G) \, $  and  $ \; H \,{\buildrel {\, \nu
\rightarrow 0 \,} \over \llongrightarrow}\, F[\varGamma] \, $,
i.e.~$ H $  is a ``quantum hyperalgebra'' at  $ \, \h := (\nu-1) \, $
and a QFA at  $ \, \h := \nu \, $.  Now we study Drinfeld's functors
for  $ H $  at  $ \, \h = (\nu\!-\!1) \, $  and at  $ \h = \nu \, $.
                                        \par
   First, a straightforward analysis like in \S 6.2 yields  $ \,
H^{\vee_{\!(\nu)}} \cong \Bbbk[\nu] \otimes_\Bbbk \hyp(G) \, $
(induced by  $ \, \gerg \cong \gerg_\nu \, $,  $ \, x \mapsto \nu^{-1}
x \, $)  whence in particular  $ \, {\big( H^{\vee_{\!(\nu)}}
\big)}{\Big|}_{\nu=0} \cong \hyp(G) \, $,  \; that is  $ \;
H^{\vee_{\!(\nu)}} \,{\buildrel {\, \nu \rightarrow 0 \,}
\over \llongrightarrow}\, \hyp(G) \, $.  Second, one can also
see (essentially,  {\sl mutatis mutandis},  like in \S 6.2) that
$ \, {\big( H^{\vee_{\!(\nu)}} \big)}^{\prime_{(\nu)}} = H \, $,
\, whence  $ \; {\big( H^{\vee_{\!(\nu)}} \big)}^{\prime_{(\nu)}}
{\Big|}_{\nu=0} = \, H{\Big|}_{\nu=0} = \, \hyp(G_0) \, = \,
F[\varGamma] \, $.
                                        \par
   At  $ \, \h = (\nu-1) \, $,  \, we can see by direct computation
that  $ \; H^{\prime_{(\nu-1)}} = \Big\langle \big( \gerg^{{(p
\hskip0,7pt)}^\infty} \big)^{\!\prime_{(\nu-1)}} \Big\rangle \; $
where  $ \, \big( \gerg^{{(p\hskip0,7pt)}^\infty} \big)^{\!\prime_{(\nu
- 1)}} := \text{\sl Span}\, \Big( \Big\{\, (\nu-1)^{p^n} x^{{(p^n)}}
\,\Big\vert\, x \in \gerg \, , n \in \N \,\Big\} \Big) \, $.  Indeed the
structure of  $ H^{\prime_{(\nu-1)}} $  depends only on the coproduct
of  $ H $,  in which  $ \nu $  plays no role; therefore we can do the same
analysis as in the trivial deformation case (see \S 5.12): the filtration
$ \underline{D} $  of  $ \hyp(G_\nu) $  is just the natural filtration
given by the order (of divided powers), and this yields the previous
description of  $ H^{\prime_{(\nu-1)}} \, $.  When specializing at 
$ \, \nu = 1 \, $  we find   
 \vskip-15pt  
  $$  H^{\prime_{(\nu-1)}} \Big/ (\nu-1) \, H^{\prime_{(\nu-1)}}
\; \cong \;  S_\Bbbk \Big( \gerg^{{(p\hskip0,7pt)^\infty}} \Big)
\hskip-2pt \bigg/ \hskip-2pt \Big( {\big\{\, x^p \,\big\}}_{x \in
\gerg^{{(p\hskip0,7pt)}^\infty}} \Big)  \; = \;  \hyp(G_0)  \; = \;
F[\varGamma]  $$
 \vskip-5pt  
\noindent   
 as Poisson Hopf algebras: in a nutshell,  $ H^{\prime_{(\nu-1)}} $  is
a QFA, at  $ \, \h = \nu-1 \, $,  for the Poisson group  $ \varGamma $.
Similarly  $ \, H^{\prime_{(\nu)}} = \Big\langle \big( \gerg^{{(p
\hskip0,7pt)}^\infty} \big)^{\!\prime_{(\nu)}} \Big\rangle \, $  with
$ \, \big( \gerg^{{(p\hskip0,7pt)}^\infty} \big)^{\!\prime_{(\nu)}} :=
\text{\sl Span}\, \Big( \Big\{\, \nu^{p^n} x^{{(p^n)}} \,\Big\vert\,
x \in \! \gerg \, , \, n \in \N \,\Big\} \Big) \, $;  \, thus on
the upshot we have
 \vskip-5pt  
  $$  H^{\prime_{(\nu)}} \Big/ \nu \, H^{\prime_{(\nu)}}  \; \cong \;
S_\Bbbk \Big( \gerg_{\text{\it ab}}^{{(p\hskip0,7pt)^\infty}} \Big)
\hskip-2pt \bigg/ \hskip-2pt \Big( {\big\{\, x^p \,\big\}}_{x \in
\gerg^{{(p\hskip0,7pt)}^\infty}} \Big)  \; = \;  F \big[
\varGamma_{\text{\it ab\,}} \big]  $$  
 \vskip-5pt  
\noindent    
 where  $ \gerg_{\text{\it ab}} $  is simply  $ \gerg $  with trivialized
Lie bracket and  $ \varGamma_{\text{\it ab}} $  is the  {\sl same
algebraic group\/}  as  $ \varGamma $  but with  {\sl trivial\/}
Poisson bracket: this comes essentially like in \S 6.2, roughly
because  $ \, \big\{ \overline{\nu\,x} \, , \overline{\nu\,y} \big\}
:= \big( \nu^{-1} [\nu\,x \, , \nu\,y \,] \big){\Big|}_{\nu=0} =
\big( \nu^{-1} \cdot \nu^3 [x,y]_\gerg \big){\Big|}_{\nu=0} =
\big( \nu \cdot \nu \, [x,y]_\gerg \big){\Big|}_{\nu=0} = 0 \; $
($ x, y \in \gerg $).
                                       \par
   Finally, we have  $ \, {\big( H^{\prime_{(\nu-1)}} \big)}^{\vee_{\!
(\nu-1)}} = \Big\langle \Big\{\, (\nu-1)^{p^n - 1\,} x^{{(p^n)}}
\,\Big\vert\, x \in \gerg \, , \, n \in \N \,\Big\} \Big\rangle \,
\varsubsetneqq H \, $  and  $ \, {\big( H^{\prime_{(\nu)}}
\big)}^{\vee_{\!(\nu)}} = \, \Big\langle \Big\{\, \nu^{p^n - 1\,}
x^{{(p^n)}} \,\Big\vert\, x \in \gerg \, , \, n \in \N \,\Big\}
\Big\rangle \, \varsubsetneqq H \, $,  \, by direct computation.
For  $ H^{\vee_{(\nu-1)}} $  we have the same features as in \S 5.7:
the analysis therein can be repeated, the upshot depending on the
nature of  $ G $
      \hbox{(or of  $ \gerg $,  essentially, in particular
on its  $ p $--lower  central series).}

 \vskip1,7truecm

\centerline {\bf \S \, 7 \;  Second example: quantum  $ {\boldkey{S}
\boldkey{L}}_\boldkey{2} $,  $ {\boldkey{S}\boldkey{L}}_\boldkey{n} $,
finite and affine Kac-Moody groups }

\vskip10pt

  {\bf 7.1 The classical setting.} \, Let  $ \Bbbk $  be any field of
characteristic  $ \, p \geq 0 \, $.  Let  $ \, G := {SL}_2(\Bbbk) \equiv
{SL}_2 \, $;  \, its tangent Lie algebra  $ \, \gerg = \gersl_2 \, $
is generated by  $ \, f $,  $ h $,  $ e \, $  (the  {\it Chevalley
generators\/})  with relations  $ \, [h,e] = 2 \, e $,  $ [h,f] =
-2 f $,  $ [e,f] = h \, $.  The formulas  $ \, \delta(f) = h \otimes f
- f \otimes h \, $,  $ \, \delta(h) = 0 \, $,  $ \, \delta(e) = h \otimes
e - e \otimes h \, $,  define a Lie cobracket on  $ \gerg $  which gives
it a structure of Lie bialgebra, corresponding to a structure of Poisson
group on  $ G $.  These formulas give also a presentation of the
co-Poisson Hopf algebra  $ U(\gerg) $  (with the standard Hopf
structure).  If  $ \, p > 0 \, $,  the  $ p $--operation
in  $ \gersl_2 $  is given by  $ \, e^{[p\,]} = 0 \, $,
$ \, f^{[p\,]} = 0 \, $,  $ \, h^{[p\,]} = h \, $.
                                             \par
  On the other hand,  $ F[{SL}_2] $  is the unital associative
commutative  $ \Bbbk $--algebra  with generators  $ \, a $,
$ b $,  $ c $,  $ d \, $  and the relation  $ \, a d - b c
= 1 \, $,  and Poisson Hopf structure given by
  $$  \displaylines{
  \Delta(a) = a \otimes a + b \otimes c \, ,  \;\,
\Delta(b) = a \otimes b + b \otimes d \, ,  \;\,
\Delta(c) = c \otimes a + d \otimes c \, ,  \;\,
\Delta(d) = c \otimes b + d \otimes d  \cr
  \epsilon(a) = 1  \, ,  \hskip5pt  \epsilon(b) = 0 \, ,  \hskip5pt
\epsilon(c) = 0 \, ,  \hskip5pt  \epsilon(d) = 1 \, ,  \hskip13pt
S(a) = d \, ,  \hskip5pt  S(b) = -b \, ,  \hskip5pt  S(c) = - c \, ,
\hskip5pt  S(d) = a  \cr
  \{a,b\} = b \, a \, , \quad  \{a,c\} = c \, a \, ,  \quad  \{b,c\}
= 0 \, ,  \quad  \{d,b\} = - b \, d \, ,  \quad  \{d,c\} = - c \, d
\, ,  \quad  \{a,d\} = 2 \, b \, c \, .  \cr }  $$
                                             \par
  The dual Lie bialgebra  $ \, \gerg^* = {\gersl_2}^{\!*} \, $
is the Lie algebra with generators  $ \, \text{f} $,  $ \text{h} $,
$ \text{e} \, $,  and relations  $ \, [\text{h},\text{e}] = \text{e} $,
$ [\text{h},\text{f}\,] = \text{f} $,  $ [\text{e},\text{f}\,] = 0 \, $,
with Lie cobracket given by  $ \, \delta(\text{f}\,) = 2 (\text{f}
\otimes \text{h} - \text{h} \otimes \text{f}\,) $,  $ \, \delta(\text{h}) =
\text{e} \otimes \text{f} - \text{f} \otimes \text{e} $,  $ \,
\delta(\text{e}) = 2 (\text{h} \otimes \text{e} - \text{e} \otimes
\text{h}) \, $  (we choose as generators  $ \, \text{f} := f^* \, $,
$ \, \text{h} := h^* \, $,  $ \, \text{e} := e^* \, $,  where
$ \, \big\{ f^*, h^*, e^* \big\} \, $  is the basis of
$ {\gersl_2}^{\! *} $  which is the dual of the basis  $ \{ f, h, e \} $
of  $ \gersl_2 \, $).  This again yields also a presentation of  $ \,
U \left( {\gersl_2}^{\!*} \right) \, $.  If  $ \, p > 0 \, $,  the
$ p $--operation  in  $ {\gersl_2}^{\!*} $  is given by  $ \,
\text{e}^{[p\,]} = 0 \, $,  $ \, \text{f}^{\,[p\,]} = 0 \, $,  $ \,
\text{h}^{[p\,]} = \text{h} \, $.  The simply connected algebraic
Poisson group whose tangent Lie bialgebra is $ {\gersl_2}^{\!*} $
can be realized as the group of pairs of matrices (the left subscript
$ s $  meaning ``simply connected'')
  $$  {}_s{{SL}_2}^{\!*} = \Bigg\{\,
\bigg( \bigg( \matrix  z^{-1} & 0 \\  y & z  \endmatrix \bigg) \, ,
\bigg( \matrix  z & x \\  0 & z^{-1}  \endmatrix \bigg) \bigg)
\,\Bigg\vert\, x, y \in \Bbbk, z \in \Bbbk \setminus \{0\} \,\Bigg\}
\; \leq \; {SL}_2 \times {SL}_2 \; .  $$
This group has centre  $ \, Z := \big\{ (I,I), (-I,-I) \big\} \, $,
so there is only one other (Poisson) group sharing the same Lie
(bi)algebra, namely the quotient  $ \, {}_a{{SL}_2}^{\!*} :=
{}_s{SL_2}^* \Big/ Z \, $  (the adjoint of  $ \, {}_s{{SL}_2}^{\!*}
\, $,  as the left subscript  $ a $  means).  Therefore  $ F\big[
{}_s{{SL}_2}^{\!*} \big] $  is the unital associative commutative
$ \Bbbk $--algebra  with generators  $ \, x $,  $ z^{\pm 1} $,
$ y $,  with Poisson Hopf structure given by
  $$  \displaylines{
   \Delta(x) = x \otimes z^{-1} + z \otimes x \, ,  \hskip21pt
\Delta\big(z^{\pm 1}\big) = z^{\pm 1} \otimes z^{\pm 1} \, ,
\hskip21pt  \Delta(y) = y \otimes z^{-1} + z \otimes y  \cr
  \epsilon(x) = 0  \, ,  \hskip10pt  \epsilon\big(z^{\pm 1}\big) = 1
\, ,  \hskip10pt  \epsilon(y) = 0 \, ,  \hskip31pt
  S(x) = -x \, ,  \hskip10pt  S\big(z^{\pm 1}\big) = z^{\mp 1} \, ,
\hskip10pt  S(y) = -y  \cr
  \{x,y\} = \big( z^2 - z^{-2} \big) \big/ 2 \, ,  \hskip27pt
\big\{z^{\pm 1},x\big\} = \pm \, x \, z^{\pm 1} \, ,  \hskip27pt
\big\{z^{\pm 1},y\big\} = \mp \, z^{\pm 1} y  \cr }  $$
(N.B.: with respect to this presentation, we have  $ \, \text{f}
= {\partial_y}{\big\vert}_e \, $,  $ \, \text{h} =
%
%
z \, {\partial_z}{\big\vert}_e \, $,  $ \, \text{e} =
{\partial_x}{\big\vert}_e \, $,  where  $ e $  is the identity
element of  $ {}_s{SL_2}^* \, $).  Moreover,  $ F \big[
{}_a{{SL}_2}^{\!*} \big] $ can be identified with the Poisson
Hopf subalgebra of  $ F\big[{}_s{SL_2}^*\big] $  spanned by
products of an even number of generators   --- i.e. monomials
of even degree: this is generated, as a unital subalgebra, by
$ x z $,  $ z^{\pm 2} $,  and  $ z^{-1} y $.
                                             \par
   In general, we shall consider  $ \, \gerg = \gerg^\tau \, $  a
semisimple Lie algebra, endowed with the Lie cobracket   ---
depending on the parameter  $ \tau $  ---   given in [Ga1], \S 1.3;
in the following we shall also retain from  [{\it loc.~cit.}]  all
the notation we need: in particular, we denote by  $ Q $,
resp.~$ P $,  the root lattice, resp.~the weight lattice, of
$ \gerg \, $,  \, and by  $ r $  the rank of  $ \gerg \, $.

 \vskip7pt

  {\bf 7.2 The\footnote{In \S\S 7--9 we should use notation
$ \, U_{q-1}(\gerg) \, $  and  $ \, F_{q-1}[G] \, $,  after
Remark 1.5 (for  $ \, \h = q-1 \, $);  instead, we write
$ \, U_q(\gerg) \, $  and  $ \, F_q[G] \, $  to be consistent
with the standard notation in use for these quantum algebras.}
           QrUEAs  $ \, U_q(\gerg) \, $.} \, We turn now to quantum
groups, starting with the  $ \gersl_2 $  case.  Let  $ R $   be any
domain,  $ \, \h \in R \setminus \{0\} \, $  an element such that  $ \,
R \big/ \h \, R = \Bbbk \, $;  \, moreover, letting  $ \, q := \h + 1
\, $  we assume that  $ q $  be invertible in  $ R $,  i.e.~there
exists  $ \, q^{-1} = {(\h + 1)}^{-1} \in R  \, $.  E.g., one can
pick  $ \, R := \Bbbk \! \left[ q, q^{-1} \right] \, $  for an
indeterminate  $ q $  and  $ \, \h := q-1 \, $,  then  $ \,
F(R) = \Bbbk(q) \, $.
                                            \par
   Let  $ \, \Bbb{U}_q(\gerg) = \Bbb{U}_q(\gersl_2) \, $  be the
associative unital  $ F(R) $--algebra  with (Chevalley-like)
generators  $ \, F $,  $ K^{\pm 1} $,  $ E $,  and relations
  $$  K K^{-1} = 1 = K^{-1} K \, ,  \;  K^{\pm 1} F = q^{\mp 2} F
K^{\pm 1} \, ,  \;  K^{\pm 1} E = q^{\pm 2} E K^{\pm 1} \, ,  \;
E F - F E = {{\, K - K^{-1} \, \over \, q - q^{-1} \,}} \, .  $$
This is a Hopf algebra, with Hopf structure given by
  $$  \displaylines{
   \Delta(F) = F \otimes K^{-1} + 1 \otimes F \, ,  \hskip21pt
\Delta \big( K^{\pm 1} \big) = K^{\pm 1} \otimes K^{\pm 1} \, ,
\hskip21pt  \Delta(E) = E \otimes 1 + K \otimes E  \cr
   \epsilon(F) = 0  \, ,  \hskip3pt  \epsilon \big( K^{\pm 1} \big)
= 1 \, ,  \hskip3pt  \epsilon(E) = 0 \, ,  \hskip15pt  S(F) = - F K ,
\hskip3pt  S \big( K^{\pm 1} \big) = K^{\mp 1} ,  \hskip3pt
S(E) = - K^{-1} E \, .  \cr }  $$
Then let  $ \, U_q(\gerg) \, $  be the  $ R $--subalgebra  of
$ \Bbb{U}_q(\gerg) $  generated by  $ \, F $,
         $ H := \displaystyle{\, K - 1 \, \over \, q - 1 \,} $,\break
$ \varGamma := \displaystyle{\, K - K^{-1} \, \over \, q - q^{-1} \,} $,
$ K^{\pm 1} $,  $ E $.  From the definition of  $ \Bbb{U}_q(\gerg) $
one gets a presentation of  $ U_q(\gerg) $  as the associative
unital algebra with generators  $ \, F $,  $ H $,  $ \varGamma $,
$ K^{\pm 1} $,  $ E $  and relations
  $$  \displaylines{
   K K^{-1} = 1 = K^{-1} K \, ,  \;\quad  K^{\pm 1} H = H K^{\pm 1} \, ,
\;\quad  K^{\pm 1} \varGamma = \varGamma K^{\pm 1} \, ,  \;\quad  H
\varGamma = \varGamma H  \cr
   ( q - 1 ) H = K - 1 \, ,  \quad  \big( q - q^{-1} \big) \varGamma = K -
K^{-1} \, ,  \quad  H \big( 1 + K^{-1} \big) = \big( 1 + q^{-1} \big)
\varGamma \, ,  \quad  E F - F E = \varGamma  \cr
   K^{\pm 1} F = q^{\mp 2} F K^{\pm 1} \, ,  \,\quad  H F = q^{-2} F H -
(q+1) F \, ,  \,\quad  \varGamma F = q^{-2} F \varGamma - \big( q + q^{-1}
\big) F  \cr
   K^{\pm 1} E = q^{\pm 2} E K^{\pm 1} \, ,  \,\quad  H E = q^{+2} E H +
(q+1) E \, ,  \,\quad  \varGamma E = q^{+2} E \varGamma + \big( q + q^{-1}
\big) E  \cr }  $$
and with a Hopf structure given by the same formulas as above
for  $ \, F $,  $ K^{\pm 1} $,  and  $ E \, $  plus
  $$  \matrix
      \Delta(\varGamma) = \varGamma \otimes K + K^{-1} \otimes \varGamma
\, , &  \hskip10pt  \epsilon(\varGamma) = 0  \, ,  &  \hskip10pt
S(\varGamma) = - \varGamma  \\
      \Delta(H) = H \otimes 1 + K \otimes H \, ,  &  \hskip10pt
\epsilon(H) = 0 \, ,  &  \hskip10pt  S(H) = - K^{-1} H \, .  \\
     \endmatrix  $$
Note also that  $ \, K = 1 + (q-1) H \, $  and  $ \, K^{-1} = K -
\big( q - q^{-1} \big) \varGamma = 1 + (q-1) H - \big( q - q^{-1}
\big) \varGamma \, $,  hence  $ \, U_q(\gerg) \, $  is generated
even by  $ \, F $,  $ H $, $ \varGamma $  and $ E $  alone.  Further,
notice that
  $$  \eqalignno{
   \Bbb{U}_q(\gerg) \;  &  = \; \hbox{free  $ F(R) $--module  over} \;\;
\Big\{\, F^a K^z E^d \,\Big\vert\, a, d \in \N, z \in \Z \,\Big\}
&   (7.1)  \cr
   U_q(\gerg) \;  &  = \; \hbox{$ R $--span of} \;\; \Big\{\,
F^a H^b \varGamma^c E^d \,\Big\vert\, a, b, c, d \in \N \,\Big\}
\;\; \hbox{inside} \;\; \Bbb{U}_q(\gerg)  &   (7.2)  \cr }  $$
which implies that  $ \, F(R) \otimes_R U_q(\gerg) = \Bbb{U}_q(\gerg)
\, $.  Moreover, definitions imply at once that  $ \, U_q(\gerg) \, $
is torsion-free, and also that it is a Hopf  $ R $--subalgebra  of
$ \Bbb{U}_q(\gerg) \, $.  Therefore  $ \, U_q(\gerg) \in \HA \, $,  \,
and in fact  $ \, U_q(\gerg) \, $  is even a QrUEA, whose semiclassical
limit is  $ \, U(\gerg) = U(\gersl_2) \, $,  \, with the generators
$ \, F $,  $ K^{\pm 1} $,  $ H $,  $ \varGamma $,  $ E \, $
respectively mapping to  $ \, f $,  $ 1 $,  $ h $,  $ h $,
$ e \in U(\gersl_2) \, $.
                                       \par
   It is also possible to define a ``simply connected'' version of
$ \Bbb{U}_q(\gerg) $  and  $ U_q(\gerg) $,  obtained from the previous
ones   --- referred to as the ``adjoint (type) ones'' ---   as follows.
For  $ \Bbb{U}_q(\gerg) $,  one simply adds a square root of  $ K^{\pm 1} $,
call it  $ L^{\pm 1} $,  as new generator; for  $ U_q(\gerg) $  one
adds the new generators  $ L^{\pm 1} $  and also  $ \, D :=
\displaystyle{ {\, L - 1 \,} \over {\, q - 1 \,} } \, $.  Then
the same analysis as before shows that  $ U_q(\gerg) $  is another
quantization (containing the ``adjoint'' one) of  $ U(\gerg) \, $.
                                               \par
   In the general case of semisimple  $ \gerg \, $,  let  $ \Bbb{U}_q
(\gerg) $  be the Lusztig-like quantum group   --- over  $ R $  ---  
associated to  $ \, \gerg = \gerg^\tau \, $  as in [Ga1], namely  $ \,
\Bbb{U}_q(\gerg) := U_{q,\varphi}^{\scriptscriptstyle M}(\gerg) \, $ 
with respect to the notation in  [{\it loc.~cit.}],  where  $ M $  is
any intermediate lattice such that  $ \, Q \leq M \leq P \, $  (this
is just a matter of choice, of the type mentioned in the statement of 
Theorem 2.2{\it (c)}):  this is a Hopf algebra over  $ F(R) $,  generated
by elements  $ \, F_i \, $,  $ M_i $,  $ E_i $  for  $ \, i=1, \dots, r
=: \hbox{\it rank}\,(\gerg) \, $.  Then let  $ \, U_q(\gerg) \, $  be
the unital  $ R $--subalgebra  of  $ \Bbb{U}_q(\gerg) $  generated by
the elements  $ \, F_i $,  $ \, H_i := \displaystyle{\, M_i - 1 \, \over
\, q - 1 \,} $,  $ \, \varGamma_i := \displaystyle{\, K_i - K_i^{-1} \,
\over \, q - q^{-1} \,} \, $,  $ M_i^{\pm 1} $,  $ E_i \, $,  where the 
$ \, K_i = M_{\alpha_i} \, $  are suitable product of  $ M_j $'s,  defined
as in [Ga1], \S 2.2 (whence  $ \, K_i $,  $ K_i^{-1} \in U_q(\gerg) \, $). 
From [Ga1], \S\S 2.5, 3.3, we have that  $ \Bbb{U}_q(\gerg) $  is the
free  $ F(R) $--module  with basis the set of monomials
  $$  \bigg\{\, \prod_{\alpha \in \Phi^+} F_\alpha^{f_\alpha}
\cdot \prod_{i=1}^n K_i^{z_i} \cdot \prod_{\alpha \in \Phi^+}
E_\alpha^{e_\alpha} \,\bigg\vert\, f_\alpha, e_\alpha \in \N \, ,
\, z_i \in \Z \, , \;\; \forall \; \alpha \in \Phi^+, \, i = 1,
\dots, n \,\bigg\}  $$
while  $ U_q(\gerg) $  is the  $ R $--span  inside  $ \Bbb{U}_q(\gerg) $
of the set of monomials
  $$  \bigg\{\, \prod_{\alpha \in \Phi^+} F_\alpha^{f_\alpha} \cdot
 \prod_{i=1}^n H_i^{t_i} \cdot \prod_{j=1}^n \varGamma_j^{c_j}
\cdot \prod_{\alpha \in \Phi^+} E_\alpha^{e_\alpha} \,\bigg\vert\,
f_\alpha, t_i, c_j, e_\alpha \in \N \;\; \forall \; \alpha \in
\Phi^+, \, i, j= 1, \dots, n \,\bigg\}  $$
(hereafter,  $ \Phi^+ $  is the set of positive roots of  $ \gerg $,
each  $ \, E_\alpha \, $,  resp.~$ \, F_\alpha \, $,  is a root vector
attached to  $ \, \alpha \in \Phi^+ $,  resp.~to  $ \, -\alpha \in
(-\Phi^+) $,  and the products of factors indexed by  $ \Phi^+ $
are ordered with respect to a fixed convex order of  $ \Phi^+ $,
see [Ga1]), whence (as for  $ \, n = 2 \, $)  $ U_q(\gerg) $  is
a free  $ R $--module.  In this case again  $ \, U_q(\gerg) \, $
is a QrUEA, with semiclassical limit  $ \, U(\gerg) \, $.

 \vskip7pt

  {\bf 7.3 Computation of  $ \, {U_q(\gerg)}' \, $  and specialization
$ \, {U_q(\gerg)}' \,{\buildrel {q \rightarrow 1} \over
\llongrightarrow}\, F[G^*] \, $.} \, We begin with the simplest
case  $ \, \gerg = \gersl_2 \, $.  From the definition of  $ \,
U_q(\gerg) = U_q(\gersl_2) \, $  we have  ($ \forall \, n \in \N $)
  $$  \displaylines{
   {} \quad  \delta_n(E) = {(\text{id} - \epsilon)}^{\otimes n} \big(
\Delta^n(E) \big) = {(\text{id} - \epsilon)}^{\otimes n} \left( \;
{\textstyle \sum_{s=1}^n} K^{\otimes (s-1)} \otimes E \otimes 1^{\otimes (n-s)}
\right) =   \hfill  \cr
   {} \hfill   = {(\text{id} - \epsilon)}^{\otimes n} \big( K^{\otimes
(n-1)} \otimes E \big) = {(K-1)}^{\otimes (n-1)} \otimes E =
{(q-1)}^{n-1} \cdot H^{\otimes (n-1)} \otimes E  \cr }  $$
from which  $ \; \delta_n\big( (q-1) E \big) \in {(q-1)}^n
U_q(\gerg) \setminus {(q-1)}^{n+1} U_q(\gerg) \, $,  \; whence  $ \,
(q-1) E \in {U_q(\gerg)}' $,  whereas  $ \, E \notin {U_q(\gerg)}' $.
Similarly,  $ \, (q-1) F \in {U_q(\gerg)}' $,  whilst
$ \, F \notin {U_q(\gerg)}' $.  As for generators  $ H $,
$ \varGamma $,  $ K^{\pm 1} $,  we have  $ \; \Delta^n(H) =
\sum_{s=1}^n K^{\otimes (s-1)} \otimes H \otimes 1^{\otimes (n-s)} $,
$ \, \Delta^n \big( K^{\pm 1} \big) = {\big( K^{\pm 1}
\big)}^{\otimes n} $,  $ \, \Delta^n(\varGamma) = \sum_{s=1}^n
K^{\otimes (s-1)} \otimes \varGamma \otimes {\big( K^{-1}
\big)}^{\otimes (n-s)} $,  \; hence for  $ \, \delta_n =
{(\text{id} - \epsilon)}^{\otimes n} \circ \Delta^n \, $  we have
  $$  \displaylines{
   \delta_n(H) = {(q-1})^{n-1} \cdot H^{\otimes n} \, ,  \qquad
\delta^n \big( K^{-1} \big) = {(q-1})^n \cdot {(- K^{-1}
H)}^{\otimes n}  \cr
   \delta^n(K) = {(q-1})^n \cdot H^{\otimes n} \, ,  \quad
\delta^n(\varGamma) = {(q-1})^{n-1} \cdot {\textstyle
\sum\limits_{s=1}^n} {(-1)}^{n-s} H^{\otimes (s-1)} \otimes
\varGamma \otimes {\big( H K^{-1} \big)}^{\otimes (n-s)}  \cr }  $$
for all  $ \, n \in \N \, $,  so that  $ \, (q-1) H $,  $ (q-1)
\varGamma $,  $ K^{\pm 1} \in {U_q(\gerg)}' \setminus (q-1)
{U_q(\gerg)}' \, $.  Therefore  $ {U_q(\gerg)}' $  contains
the subalgebra  $ U' $  generated by  $ \, (q-1) F $,  $ K $,
$ K^{-1} $,  $ (q-1) H $,  $ (q-1) \varGamma $,  $ (q-1) E \, $.  On
the other hand, using (7.2) a thorough   --- but straightforward ---
computation along the same lines as above shows that any element in
$ {U_q(\gerg)}' $  does necessarily lie in  $ U' $  (details are
left to the reader: everything follows from definitions and the
formulas above for  $ \, \Delta^n \, $).  Thus  $ {U_q(\gerg)}' $
is nothing but the subalgebra of  $ U_q(\gerg) $  generated by
$ \, \dot{F} := (q-1) F $,  $ K $,  $ K^{-1} $,  $ \dot{H} :=
(q-1) H $,  $ \dot{\varGamma} := (q-1) \varGamma $,  $ \dot{E} :=
(q-1) E \, $;  notice also that the generator  $ \, \dot{H} $  is
unnecessary, for  $ \, \dot{H} = K - 1 \, $.  As a consequence,
$ {U_q(\gerg)}' $  can be presented as the unital associative
$ R $--algebra  with generators  $ \, \dot{F} $,
$ \dot{\varGamma} $,  $ K^{\pm 1} $,  $ \dot{E} \, $
and relations
  $$  \displaylines{
   K K^{-1} = 1 = K^{-1} K ,  \,\;  K^{\pm 1} \dot{\varGamma} =
\dot{\varGamma} K^{\pm 1} ,  \,\;  \big( 1 + q^{-1} \big)
\dot{\varGamma} = K - K^{-1} ,  \,\;  \dot{E} \dot{F} - \dot{F}
\dot{E} = (q-1) \dot{\varGamma}  \cr
   K - K^{-1} = \big( 1 + q^{-1}\big) \dot{\varGamma} \, ,
\quad  K^{\pm 1} \dot{F} = q^{\mp 2} \dot{F} K^{\pm 1} \, ,
\quad  K^{\pm 1} \dot{E} = q^{\pm 2} \dot{E} K^{\pm 1}  \cr
   \dot{\varGamma} \dot{F} = q^{-2} \dot{F} \dot{\varGamma} -
(q-1) \big( q + q^{-1} \big) \dot{F} \, ,  \qquad  \dot{\varGamma}
\dot{E} = q^{+2} \dot{E} \dot{\varGamma} + (q-1) \big( q + q^{-1} \big)
\dot{E}  \cr }  $$
with Hopf structure given by
  $$  \matrix
      \Delta\big(\dot{F}\big) = \dot{F} \otimes K^{-1} + 1 \otimes
\dot{F} \, ,  &  \hskip25pt  \epsilon\big(\dot{F}\big) = 0 \, ,  &
\hskip25pt  S\big(\dot{F}\big) = - \dot{F} K \, \phantom{.}  \\
      \Delta\big(\dot{\varGamma}\big) = \dot{\varGamma} \otimes K +
K^{-1} \otimes \dot{\varGamma} \, ,  &  \hskip25pt  \epsilon\big(
\dot{\varGamma}\big) = 0 \, ,  &  \hskip25pt  S\big(\dot{\varGamma}\big)
= - \dot{\varGamma}  \\
      \Delta \big( K^{\pm 1} \big) = K^{\pm 1} \otimes K^{\pm 1} \, ,
&  \hskip25pt  \epsilon \big( K^{\pm 1} \big) = 1 \, ,  &  \hskip25pt
S \big( K^{\pm 1} \big) = K^{\mp 1} \, \phantom{.}  \\
      \Delta\big(\dot{E}\big) = \dot{E} \otimes 1 + K \otimes \dot{E} \, ,
&  \hskip25pt  \epsilon\big(\dot{E}\big) = 0 \, ,  &  \hskip25pt
S\big(\dot{E}\big) = - K^{-1} \dot{E} \, .  \\
     \endmatrix  $$
   \indent   When  $ \, q \rightarrow 1 \, $,  an easy direct
computation shows that this gives a presentation of the function
algebra  $ \, F \big[ {}_a{{SL}_2}^{\!*} \big] $,  and the Poisson
structure that $ \, F \big[ {}_a{{SL}_2}^{\!*} \big] \, $  inherits
from this quantization process is exactly the one coming from the
Poisson structure on  $ \, {}_a{{SL}_2}^{\!*} \, $:  in fact,
there is a Poisson Hopf algebra isomorphism
 $$  {U_q(\gerg)}' \Big/ (q-1) \, {U_q(\gerg)}'
\,{\buildrel \cong \over \llongrightarrow}\, F \big[ {}_a{{SL}_2}^{\!*} \big]
\quad \Big( \, \subseteq F \big[ {}_s{{SL}_2}^{\!*} \big] \, \Big)  $$
given by:  $ \;\, \dot{E} \mod (q-1) \mapsto x{}z \, $,
$ \, K^{\pm 1} \mod (q-1) \mapsto z^{\pm 2} \, $,
$ \, \dot{H} \mod (q-1)
                \mapsto z^2 - 1 \, $,\break
$ \, \dot{\varGamma} \mod (q-1) \mapsto \big( z^2 - z^{-2}
\big) \Big/ 2 \, $,  $ \dot{F} \mod (q-1) \mapsto z^{-1}{}y \, $.
In other words, $ \, {U_q(\gerg)}' \, $ specializes to  $ \, F \big[
{}_a{{SL}_2}^{\!*} \big] \, $  as a Poisson Hopf algebra.  {\sl Note
that this was predicted by Theorem 2.2{\it (c)\/}  when  $ \,
\Char(\Bbbk) = 0 \, $,  \, but our analysis now proved it also
for  $ \, \Char(\Bbbk) > 0 \, $.}
                                            \par
  Note that we got the  {\sl adjoint}  Poisson group dual of
$ \, G = {SL}_2 \, $,  that is  $ {}_a{{SL}_2}^{\!*} \, $; a
different choice of the initial QrUEA leads us to the {\sl simply
connected}  one, i.e.~$ \, {}_s{{SL}_2}^{\!*} $.  Indeed, if we
start from the ``simply connected'' version of  $ U_q(\gerg) $  (see
\S 7.2) the same analysis shows that  $ {U_q(\gerg)}' $  is like
above but for containing also the new generators  $ L^{\pm 1} $,
and similarly when specializing  $ q $  at  $ 1 $:  thus we get the
function algebra of a Poisson group which is a double covering of
$ {}_a{{SL}_2}^{\!*} $,  namely  $ {}_s{{SL}_2}^{\!*} $.  So changing
the QrUEA quantizing  $ \gerg $  we get two different QFAs, one for each
of the two connected Poisson algebraic groups dual of  $ {SL}_2 $,
i.e.~with tangent Lie bialgebra  $ \, {\gersl_2}^{\!*} \, $;  this
shows the dependence of the group  $ G^\star $  (here denoted
$ G^* $  since  $ \gerg^\times \! = \! \gerg^* $)  in  Theorem
2.2{\it (c)}  on the choice of the QrUEA (for a fixed  $ \gerg $).
                                                 \par
   With a bit more careful study, exploiting the analysis in [Ga1],
one can treat the general case too: we sketch briefly our arguments
--- restricting to the simply laced case, to simplify the exposition
---   leaving to the reader the (straightforward) task of filling in
details.
                                                 \par
   So now let  $ \, \gerg = \gerg^\tau $  be a semisimple Lie algebra,
as in \S 7.1, and let  $ \, U_q(\gerg) \, $  be the QrUEA introduced in
\S 7.2: our aim again is to compute the QFA  $ \, {U_q(\gerg)}' \, $.
                                               \par
   The same computations as for  $ \, \gerg = {\frak s}{\frak l}(2)
\, $  show that  $ \; \delta_n(H_i) = {(q-1})^{n-1} \cdot
H_i^{\otimes n} \, $  and  $ \, \delta^n(\varGamma_i) = {(q-1})^{n-1}
\cdot \sum_{s=1}^n {(-1)}^{n-s} H_i^{\otimes (s-1)} \otimes
\varGamma_i \otimes {\big( H_i K_i^{-1} \big)}^{\otimes (n-s)} $,
which gives
  $$  \dot{H}_i := (q-1) H_i \in {U_q(\gerg)}' \setminus (q-1) \,
{U_q(\gerg)}'  \quad   \hbox{and}  \quad  \dot{\varGamma}_i :=
(q-1) \, \varGamma_i \in {U_q(\gerg)}' \setminus (q-1) \,
{U_q(\gerg)}' \, .  $$
                                                 \par
   As for root vectors, let  $ \, \dot{E}_\gamma := (q-1) E_\gamma
\, $  and  $ \, \dot{F}_\gamma := (q-1) F_\gamma \, $  for all
$ \, \gamma \in \Phi^+ \, $:  using the same type of arguments
          as in [Ga1]\footnote{Note that in [Ga1] the assumption
$ \, \Char(\Bbbk) = 0 \, $  is made throughout: nevertheless,  {\sl
this hypothesis is not necesary\/}  for the analysis we are concerned
with right now!},
   \S 5.16, we can prove that  $ \, E_\alpha \not\in {U_q(\gerg)}'
\, $  but  $ \, \dot{E}_\alpha \in {U_q(\gerg)}' \setminus (q-1) \,
{U_q(\gerg)}' \, $.  In fact, let  $ \, \Bbb{U}_q({\frak b}_+) \, $  and
$ \, \Bbb{U}_q({\frak b}_-) \, $  be quantum Borel subalgebras, and
$ \, {\frak U}_{\varphi,\geq}^{\scriptscriptstyle M} \, $,
$ \, {\Cal U}_{\varphi,\geq}^{\scriptscriptstyle M} \, $,
$ \, {\frak U}_{\varphi,\leq}^{\scriptscriptstyle M} \, $,
$ \, {\Cal U}_{\varphi,\leq}^{\scriptscriptstyle M} \, $
their  $ R $--subalgebras  defined in [Ga1], \S 2: then both
$ \, \Bbb{U}_q({\frak b}_+) \, $  and  $ \, \Bbb{U}_q({\frak b}_-) \, $
are Hopf subalgebras  of  $ \Bbb{U}_q(\gerg) $;  in addition, letting
$ M' $  be the lattice between  $ Q $  and  $ P $  dual of  $ M $
(in the sense of [Ga1], \S 1.1, there exists an  $ F(R) $--valued
perfect Hopf pairing between  $ \Bbb{U}_q({\frak b}_\pm) $  and
$ \Bbb{U}_q({\frak b}_\mp) $   --- one built up on  $ M $  and
the other on  $ M' $  ---   such that  $ \,
{\frak U}_{\varphi,\geq}^{\scriptscriptstyle M} = {\Big(
{\Cal U}_{\varphi,\leq}^{\scriptscriptstyle M'} \Big)}^{\!\bullet}
\, $,  $ \, {\frak U}_{\varphi,\leq}^{\scriptscriptstyle M} = {\Big(
{\Cal U}_{\varphi,\geq}^{\scriptscriptstyle M'} \Big)}^{\!\bullet}
\, $,  $ \, {\Cal U}_{\varphi,\geq}^{\scriptscriptstyle M} = {\Big(
{\frak U}_{\varphi,\leq}^{\scriptscriptstyle M'} \Big)}^{\!\bullet}
\, $,  and  $ \, {\Cal U}_{\varphi,\leq}^{\scriptscriptstyle M} =
{\Big( {\frak U}_{\varphi,\geq}^{\scriptscriptstyle M'}
\Big)}^{\!\bullet} \, $.  Now,  $ \, \big( q - q^{-1} \big) E_\alpha
\in {\Cal U}_{\varphi,\geq}^{\scriptscriptstyle M} = {\Big(
{\frak U}_{\varphi,\leq}^{\scriptscriptstyle M'}
\Big)}^{\!\bullet} \, $,  hence   --- since
$ {\frak U}_{\varphi,\leq}^{\scriptscriptstyle M'} $
is an algebra ---   we have  $ \, \Delta \Big( \big(
q - q^{-1} \big) E_\alpha \Big) \in {\Big(
{\frak U}_{\varphi,\leq}^{\scriptscriptstyle M'} \otimes
{\frak U}_{\varphi,\leq}^{\scriptscriptstyle M'} \Big)}^{\!\bullet}
= {\Big( {\frak U}_{\varphi,\leq}^{\scriptscriptstyle M'}
\Big)}^{\!\bullet} \otimes
{\Big( {\frak U}_{\varphi,\leq}^{\scriptscriptstyle M'}
\Big)}^{\!\bullet} = {\Cal U}_{\varphi,\geq}^{\scriptscriptstyle M}
\otimes {\Cal U}_{\varphi,\geq}^{\scriptscriptstyle M} \, $.
Therefore, by definition of
$ {\Cal U}_{\varphi,\geq}^{\scriptscriptstyle M} $
and by the PBW theorem for it and for
$ {\frak U}_{\varphi,\leq}^{\scriptscriptstyle M'} $
(cf.~[Ga1], \S 2.5) we have that  $ \, \Delta \Big( \big(
q - q^{-1} \big) E_\alpha \Big) \, $  is an  $ R $--linear
combination like  $ \; \Delta \Big( \big( q - q^{-1} \big)
E_\alpha \Big) = \sum_r \, A^{(1)}_r \otimes A^{(2)}_r \; $
in which the  $ A^{(j)}_r $'s  are monomials in the  $ M_j $'s
and in the  $ \, \overline{E}_\gamma $'s,  where  $ \,
\overline{E}_\gamma := \big( q - q^{-1} \big) E_\gamma \, $
for all  $ \, \gamma \in \Phi^+ \, $:  iterating, we find that
$ \, \Delta^\ell \Big( \big( q - q^{-1} \big) E_\alpha \Big) \, $
is an  $ R $--linear  combination
  $$  \Delta^\ell \Big( \big( q - q^{-1} \big) E_\alpha \Big) \; =
\; {\textstyle \sum_r} \, A^{(1)}_r \otimes A^{(2)}_r \otimes \cdots
\otimes A^{(\ell)}_r   \eqno (7.3)  $$
in which the  $ A^{(j)}_r $'s  are again monomials in the  $ M_j $'s
and in the  $ \, \overline{E}_\gamma $'s.  Now, we distinguish
two cases: either  $ A^{(j)}_r $  does contain some  $ \,
\overline{E}_\gamma \, (\, \in \big( q - q^{-1} \big) \, U_q(\gerg)
\big) \, $,  thus  $ \, \epsilon \Big( A^{(j)}_r \Big) = A^{(j)}_r
\in (q-1) \, U_q(\gerg) \, $  whence  $ \, (\text{id} - \epsilon)
\Big( A^{(j)}_r \Big) = 0 \, $;  or  $ A^{(j)}_r $  does not contain
any  $ \, \overline{E}_\gamma \, $  and is only a monomial in the
$ M_t $'s,  say  $ \, A^{(j)}_r = \prod_{t=1}^n M_t^{m_t} \, $:  then
$ \, (\text{id} - \epsilon) \Big( A^{(j)}_r \Big) = \prod_{t=1}^n
M_t^{m_t} - 1 = \prod_{t=1}^n {\big( (q-1) \, H_t + 1 \big)}^{m_t} -
1 \in (q-1) \, U_q(\gerg) \, $.  In addition, for some  ``$ Q $--grading
reasons'' (as in [Ga1], \S 5.16), in each one of the summands in (7.3)
the sum of all the  $ \gamma $'s  such that the (rescaled) root
vectors  $ \overline{E}_\gamma $  occur in any of the factors
$ A^{(1)}_r $,  $ A^{(2)}_r $,  $ \dots $,  $ A^{(n)}_r $  must
be equal to  $ \alpha $:  therefore, in each of these summands at
least one factor  $ \overline{E}_\gamma $  does occur.  The conclusion
is that  $ \; \delta_\ell \big( \overline{E}_\alpha \big) \in \big(
1 + q^{-1} \big) (q-1)^\ell \, {U_q(\gerg)}^{\otimes \ell} \; $
(the factor  $ \, \big( 1 + q^{-1} \big) \, $  being there because at
least one rescaled root vector  $ \overline{E}_\gamma $  occurs in each
summand of  $ \, \delta_\ell \big( \overline{E}_\alpha \big) \, $,
thus providing a coefficient  $ \, \big( q - q^{-1} \big) \, $  the
term  $ \, \big( 1 + q^{-1} \big) \, $  is factored out of), whence
$ \; \delta_\ell \big( \dot{E}_\alpha \big) \in (q-1)^\ell \,
{U_q(\gerg)}^{\otimes \ell} \, $.  More precisely, we have also
$ \; \delta_\ell \big( \dot{E}_\alpha \big) \not\in (q-1)^{\ell+1} \,
{U_q(\gerg)}^{\otimes \ell} \, $,  \; for we can easily check that
$ \, \Delta^\ell \big( \dot{E}_\alpha \big) \, $  is the sum of  $ \,
M_\alpha \otimes M_\alpha \otimes \cdots \otimes M_\alpha \otimes
\dot{E}_\alpha \, $  plus other summands which are  $ R $--linearly
independent of this first term: but then  $ \, \delta_\ell \big(
\dot{E}_\alpha \big) \, $  is the sum of  $ \, {(q-1)}^{\ell-1}
H_\alpha \otimes H_\alpha \otimes \cdots \otimes H_\alpha \otimes
\dot{E}_\alpha \, $  (where  $ \, H_\alpha := {\, M_\alpha - 1 \,
\over \, q - 1 \,} \, $  is equal to an  $ R $--linear  combination
of products of  $ M_j $'s  and  $ H_t $'s)  plus other summands which
are  $ R $--linearly  independent of the first one, and since
$ \, H_\alpha \otimes H_\alpha \otimes \cdots \otimes H_\alpha \otimes
\dot{E}_\alpha \not\in (q-1)^2 \, {U_q(\gerg)}^{\otimes \ell} \, $
we can conclude as claimed.  Therefore  $ \; \delta_\ell \big(
\dot{E}_\alpha \big) \in (q-1)^\ell \, {U_q(\gerg)}^{\otimes \ell}
\setminus (q-1)^{\ell+1} \, {U_q(\gerg)}^{\otimes \ell} \, $,
\; whence we get  $ \; \displaystyle{ \dot{E}_\alpha := (q-1)
E_\alpha \in {U_q(\gerg)}' \setminus (q-1) \, {U_q(\gerg)}' \;\;
\forall \; \alpha \in \Phi^+ } \, $.  \; An entirely similar analysis
yields also  $ \; \displaystyle{ \dot{F}_\alpha := (q-1) F_\alpha
\in {U_q(\gerg)}' \setminus (q-1) \, {U_q(\gerg)}' \;\; \forall \;
\alpha \in \Phi^+ } \, $.
                                                 \par
   Summing up, we have found that  $ \, {U_q(\gerg)}' \, $  contains
for sure the subalgebra  $ U' $  generated by  $ \, \dot{F}_\alpha \, $,
$ \, \dot{H}_i \, $,  $ \, \dot{\varGamma}_i \, $,  $ \, \dot{E}_\alpha
\, $  for all  $ \, \alpha \in \Phi^+ \, $  and all  $ \, i = 1, \dots,
n \, $.  On the other hand, using (7.2) a thorough   --- but
straightforward ---   computation along the same lines as above
shows that any element in  $ {U_q(\gerg)}' $  must lie in  $ U' $
(details are left to the reader).  Thus finally  $ \, {U_q(\gerg)}'
= U' \, $,  so we have a concrete description of  $ {U_q(\gerg)}' $.
                                                \par
   Now compare  $ \, U' = {U_q(\gerg)}' \, $  with the algebra
$ \, {\Cal U}_\varphi^{\scriptscriptstyle M}(\gerg) \, $  in [Ga1],
\S 3.4 (for  $ \, \varphi = 0 \, $),  the latter being just the
$ R $--subalgebra  of  $ \Bbb{U}_q(\gerg) $  generated by the set
$ \, \big\{\, \overline{F}_\alpha, M_i, \overline{E}_\alpha
\,\big\vert\, \alpha \in \Phi^+, i=1, \dots, n \,\big\} \, $.
First of all, by definition, we have  $ \,
{\Cal U}_\varphi^{\scriptscriptstyle M}(\gerg)
\subseteq U' = {U_q(\gerg)}' \, $;  moreover,
 \eject   
  $$  \dot{F}_\alpha \equiv {\, 1 \, \over \, 2 \,}
\overline{F}_\alpha \, ,  \quad  \dot{E}_\alpha \equiv
{\, 1 \, \over \, 2 \,} \overline{E}_\alpha \, ,  \quad
\dot{\varGamma}_i \equiv {\, 1 \, \over \, 2 \,}
\big( K_i - K_i^{-1} \big) \, \mod\, (q-1) \,
{\Cal U}_\varphi^{\scriptscriptstyle M}(\gerg)
\qquad  \forall\;\, \alpha, \; \forall\;\, i \, .  $$
Then   
  $$  {\big(U_q(\gerg)'\big)}_1 := {U_q(\gerg)}' \Big/ (q-1) \,
{U_q(\gerg)}' \, = \, {\Cal U}_\varphi^{\scriptscriptstyle M}
(\gerg) \Big/ (q-1) \, {\Cal U}_\varphi^{\scriptscriptstyle M}
(\gerg) \, \cong \, F \big[ G^*_{\scriptscriptstyle M} \big]  $$
where  $ \, G^*_{\scriptscriptstyle M} \, $  is the Poisson
group dual of  $ \, G = G^\tau \, $  with centre  $ \,
Z(G^*_{\scriptscriptstyle M}) \cong M \big/ Q \, $  and
fundamental group  $ \, \pi_1(G^*_{\scriptscriptstyle M})
\cong P \big/ M \, $,  and the isomorphism (of Poisson Hopf
algebras) on the right is given by [Ga1], Theorem 7.4 (see
also references therein for the original statement and proof
of this result).  In other words,  $ \, {U_q(\gerg)}' \, $
specializes to  $ \, F \big[ G^*_{\scriptscriptstyle M} \big] \, $
{\sl as a Poisson Hopf algebra},  as prescribed by Theorem 2.2.  By
the way, notice that in the present case the dependence of the dual
group  $ \, G^\star = G^*_{\scriptscriptstyle M} \, $  on the choice
of the initial QrUEA (for fixed  $ \gerg $)   --- mentioned in the
last part of the statement of Theorem 2.2{\it (c)}  ---   is evident.
                                                 \par
   By the way, the previous discussion applies as well to the case
of  $ \gerg $  {\it an untwisted affine Kac-Moody algebra\/}:  one
just has to substitute any quotation from [Ga1]   --- referring to
some result about  {\sl finite}  Kac-Moody algebras ---   with a
similar quotation from [Ga3]   --- referring to the corresponding
analogous result about untwisted  {\sl affine}  Kac-Moody algebras.

\vskip7pt

  {\bf 7.4 The identity  $ \, {\big({U_q(\gerg)}'\big)}^\vee =
U_q(\gerg) \, $.} \,  In the present section we check that part of
Theorem 2.2{\it (b)}  claiming that, when  $ \, p = 0 \, $,  \, one
has  $ \;  H \in \QrUEA \,\Longrightarrow\, {\big( H' \big)}^\vee =
H \; $  for  $ \, H = U_q(\gerg) \, $  as above.  In addition, our
proof now will work for the case  $ \, p > 0 \, $  as well.  Of
course, we start once again from  $ \, \gerg = \gersl_2 \, $.
                                                     \par
   Since  $ \, \epsilon\big(\dot{F}\big) = \epsilon\big(\dot{H}\big)
= \epsilon\big(\dot{\varGamma}\big) = \epsilon\big(\dot{E}\big) = 0
\, $,  the ideal  $ \, J := \text{\sl Ker}\,\big( \epsilon \, \colon
\, {U_q(\gerg)}' \! \loongrightarrow R \,\big) \, $  is generated by
$ \, \dot{F} $,  $ \dot{H} $,  $ \dot{\varGamma} $,  and  $ \dot{E}
\, $.  This implies that  $ J $  is the  $ R $--span  of  $ \,
\Big\{\, {\dot{F}}^\varphi {\dot{H}}^\kappa {\dot{\varGamma}}^\gamma
{\dot{E}}^\eta \,\Big\vert\, $
$ (\varphi, \kappa, \gamma, \eta) \in \N^4
\setminus \{(0,0,0,0)\} \! \Big\} \, $.  Now  $ \,I := \text{\sl
Ker} \, \Big( {U_q(\gerg)}' \, {\buildrel \epsilon \over
{\relbar\joinrel\twoheadrightarrow}} \, R \,
{\buildrel {q \mapsto 1} \over \llongtwoheadrightarrow} \, \Bbbk
\Big) = (q-1) \cdot {U_q(\gerg)}' + J $,  \; therefore we get that
$ \; {\big({U_q(\gerg)}'\big)}^\vee := \sum_{n \geq 0} {\Big(
{(q-1)}^{-1} I \Big)}^n \; $  is generated, as a unital
$ R $-subalgebra  of  $ \Bbb{U}_q(\gerg) $,  by the elements
$ \, {(q-1)}^{-1} \dot{F} = F $,  $ {(q-1)}^{-1} \dot{H} = H $,
$ {(q-1)}^{-1} \dot{\varGamma} = \varGamma $,  $ {(q-1)}^{-1} \dot{E}
= E $,  hence it coincides with  $ U_q(\gerg) $,  q.e.d.
                                                 \par
  An entirely similar analysis works in the ``adjoint'' case as well;
and also,  {\it mutatis mutandis},  for the general semisimple or
affine Kac-Moody case.

\vskip7pt

  {\bf 7.5 The quantum hyperalgebra  $ \hyp_q(\gerg) $.} \, Let  $ G $
be a semisimple (affine) algebraic group, with Lie algebra  $ \gerg $,
and let  $ \Bbb{U}_q(\gerg) $  be the quantum group considered in the previous
sections.  Lusztig introduced (cf.~[Lu1-2]) a ``quantum hyperalgebra'',
i.e.~a Hopf subalgebra of  $ \Bbb{U}_q(\gerg) $  over  $ \Z \big[ q, q^{-1}
\big] $  whose specialization at  $ \, q = 1 \, $  is exactly the
Kostant's  $ \Z $--integer  form  $ U_\Z(\gerg) $  of  $ U(\gerg) $
from which one gets the hyperalgebra  $ \hyp(\gerg) $  over any field
$ \Bbbk $  of characteristic  $ \, p > 0 \, $  by scalar extension,
namely  $ \, \hyp(\gerg) = \Bbbk \otimes_\Z U_\Z(\gerg) \, $.  In
fact, to be precise one needs a suitable enlargement of the algebra
given by Lusztig, which is given in [DL], \S 3.4, and denoted
by  $ \, \varGamma(\gerg) $.  Now we study Drinfeld's functors (at
$ \, \h = q - 1 \, $)  on  $ \, \hyp_q(\gerg) := R \otimes_{\Z[q,q^{-1}]}
\varGamma(\gerg) \, $  (with  $ R $  like in \S 7.2), taking as sample
the case of  $ \, \gerg = \gersl_2 \, $.
                                             \par
   Let  $ \, \gerg = \gersl_2 \, $.  Let  $ \hyp^\Z_q(\gerg) $  be
the unital  $ \Z\big[q,q^{-1}\big] $--subalgebra  of  $ \Bbb{U}_q(\gerg) $
(say the one of ``adjoint type'' defined like above  {\sl but over\/}
$ \Z\big[q,q^{-1}\big] $)  generated the ``quantum divided powers''
$ \, \displaystyle{ F^{(n)} \! := F^n \! \Big/ {[n]}_q! } \; $,
$ \, \displaystyle{ \left( {{K \, ; \, c} \atop {n}} \right) \!
:= \prod_{s=1}^n {{\; q^{c+1-s} K - 1 \,} \over {\; q^s - 1 \;}} } \, $,
$ \, \displaystyle{ E^{(n)} \! := E^n \! \Big/ {[n]}_q! } \; $  (for all
$ \, n \in \N \, $,  $ \, c \in \Z \, $)  and by  $ K^{-1} $,  where
$ \, {[n]}_q! := \prod_{s=1}^n {[s]}_q \, $  and  $ \, {[s]}_q = \big(
q^s - q^{-s} \big) \Big/ \big( q - q^{-1} \big) \, $  for all  $ \, n $,
$ s \in \N \, $.  Then (cf.~[DL]) this is a Hopf subalgebra of  $ \Bbb{U}_q
(\gerg) $,  and  $ \, \hyp^\Z_q(\gerg){\Big|}_{q=1} \cong\, U_\Z(\gerg)
\, $;  \, therefore  $ \, \hyp_q(\gerg) := R \otimes_{\Z[q,q^{-1}]}
\hyp^\Z_q(\gerg) \, $  (for any  $ R $  like in \S 7.2, with  $ \,
\Bbbk := R \big/ \h \, R \, $  and  $ \, p := \Char(\Bbbk) \, $)
specializes at  $ \, q = 1 \, $  to the  $ \Bbbk $--hyperalgebra
$ \hyp(\gerg) $.  Moreover, among all the  $ \left( {{K ; \, c}
\atop {n}} \right) $'s  it is enough to take only those with
$ \, c = 0 \, $.  {\sl From now on we assume  $ \, p > 0 \, $.}
                                             \par
   Using formulas for the iterated coproduct in [DL], Corollary 3.3
(which uses the opposite coproduct than ours, but this doesn't matter),
and exploiting the PBW-like theorem for  $ \hyp_q(\gerg) $  (see [DL]
again) we see by direct inspection that  $ \, {\hyp_q(\gerg)}' \, $
is the unital  $ R $--subalgebra  of  $ \hyp_q(\gerg) $  generated by
$ K^{-1} $  and the ``rescaled quantum divided powers''  $ \, {(q-1)}^n
F^{(n)} \, $,  $ \, {(q-1)}^n \left( {{K ; \, 0} \atop {n}} \right)
\, $  and  $ \, {(q-1)}^n E^{(n)} \, $  for all  $ \, n \in \N \, $.
Since  $ \, {[n]}_q!{\Big|}_{q=1} = n! = 0 \, $  iff  $ \, p \,\Big|
n \, $,  we argue that  $ \, {\hyp_q(\gerg)}'{\Big|}_{q=1} \, $  is
generated by the corresponding specializations of  $ \, {(q-1)}^{p^s}
F^{(p^s)} \, $,  $ \, {(q-1)}^{p^s} \left( {{K ; \, 0} \atop {p^s}}
\right) \, $  and  $ \, {(q-1)}^{p^s} E^{(p^s)} \, $  for all  $ \, s
\! \in \! \N \, $:  in particular this shows that the spectrum of  $ \,
{\hyp_q(\gerg)}'{\Big|}_{q=1} \, $  has dimension 0 and height 1, and
its cotangent Lie algebra  $ \, J \Big/ J^{\,2} \, $   --- where  $ J $
is the augmentation ideal of  $ {\hyp_q(\gerg)}' {\Big|}_{q=1} $  ---
has basis  $ \, \Big\{\, {(q\!-\!1)}^{p^s} F^{(p^s)}, \,
{(q\!-\!1)}^{p^s} \! \left( {{K ; \, 0} \atop {p^s}} \right) ,
\, {(q\!-\!1)}^{p^s} E^{(p^s)} \; \mod \, (q\!-\!1) \,
{\hyp_q(\gerg)}' \, \mod J^{\,2} \;\Big|\; s \in \N \,\Big\} \, $.
Furthermore,  $ \, \big( {\hyp_q(\gerg)}' \big)^{\!\vee} \, $  is
generated by  $ \, {(q-1)}^{p^s-1} F^{(p^s)} \, $,  $ \, {(q-1)}^{p^s-1}
\left( {{K ; \, 0} \atop {p^s}} \right) \, $,  $ \, K^{-1} \, $  and
$ \, {(q-1)}^{p^s-1} E^{(p^s)} \, $  for all  $ \, s \in \N \, $:  \,
in particular we have that  $ \, \big( {\hyp_q(\gerg)}' \big)^{\!\vee}
\subsetneqq \hyp_q(\gerg) \, $,  \, and  $ \, \big( {\hyp_q(\gerg)}'
\big)^{\!\vee}{\Big|}_{q=1} \, $  is generated by the cosets modulo
$ (q-1) $  of the previous elements, which do form a basis of the
restricted Lie bialgebra  $ \gerk $  such that  $ \, \big(
{\hyp_q(\gerg)}' \big)^{\!\vee}{\Big|}_{q=1} = \, \u(\gerk) \, $.
                                             \par
   We performed the previous study using the ``adjoint'' version of
$ U_q(\gerg) $  as starting point: instead, we can use as well its
``simply connected'' version, thus obtaining a ``simply connected
version of  $ \hyp_q(\gerg) $''  which is defined exactly like before
but for using  $ L^{\pm 1} $  instead of  $ K^{\pm 1} $  throughout;
up to these changes, the analysis and its outcome will be exactly
the same.  Note that all quantum objects involved   --- namely,
$ \hyp_q(\gerg) $,  $ {\hyp_q(\gerg)}' $  and  $ \big( {\hyp_q(\gerg)}'
\big)^{\!\vee} $  ---   will strictly contain the corresponding
``adjoint'' quantum objects; on the other hand, the semiclassical
limit is the same in the case of  $ \hyp_q(\gerg) $  (giving
$ \hyp(\gerg) $,  in both cases)  and in the case of  $ \big(
{\hyp_q(\gerg)}' \big)^{\!\vee} $  (giving  $ \u(\gerk) $,  in both
cases), whereas the semiclassical limit of  $ {\hyp_q(\gerg)}' $  in
the ``simply connected'' case is a (countable) covering of that in
the ``adjoint'' case.
                                             \par
   The general case of semisimple or affine Kac-Moody  $ \gerg $
can be dealt with similarly, with analogous outcome.  Indeed,
$ \hyp^\Z_q(\gerg) $  is defined as the unital  $ \Z \big[ q,
q^{-1} \big] $--subalgebra  of  $ \Bbb{U}_q(\gerg) $  (defined like
before  {\sl but over\/}  $ \Z\big[q,q^{-1}\big] $) generated by
$ K_i^{-1} $  and the ``quantum divided powers'' (in the above
sense)  $ \, F_i^{(n)} \, $,  $ \, \left( {{K_i ; \, c} \atop {n}}
\right) \, $,  $ \, E_i^{(n)} \, $  for all  $ \, n \in \N \, $,
$ \, c \in \Z \, $  and  $ \, i =1, \dots, \text{\it rank}\,(\gerg)
\, $  (notation of \S 7.2, but now each divided power relative to
$ i $  is built upon  $ q_i $,  see [Ga1]).  Then (cf.~[DL]) this
is a Hopf subalgebra of  $ \Bbb{U}_q(\gerg) $  with  $ \, \hyp^\Z_q(\gerg)
{\Big|}_{q=1} \cong\, U_\Z(\gerg) \, $,  \, so  $ \, \hyp_q(\gerg)
:= R \otimes_{\Z[q,q^{-1}]} \hyp^\Z_q(\gerg) \, $  (for any
$ R $  like before)  specializes at  $ \, q = 1 \, $  to the
$ \Bbbk $--hyperalgebra  $ \hyp(\gerg) $;  \, and among the
$ \left( {{K_i ; \, c} \atop {n}} \right) $'s  it is enough
to take those with  $ \, c = 0 \, $.
                                             \par
   Again a PBW-like theorem holds for  $ \hyp_q(\gerg) $  (see [DL]),
where powers of root vectors are replaced by quantum divided powers
like  $ \, F_\alpha^{(n)} \, $,  $ \, \left( {{K_i ; \, c} \atop {n}}
\right) \cdot K_i^{-\text{\it Ent}(n/2)} \, $  and  $ \, E_\alpha^{(n)}
\, $,  \, for all positive roots  $ \alpha $  of  $ \gerg $  (each
divided power being relative to  $ q_\alpha $,  see [Ga1])  both in
the finite and in the affine case.  Using this and the same type of
arguments as in \S 7.3   --- i.e.~the perfect graded Hopf pairing
between quantum Borel subalgebras ---   we see by direct inspection
that  $ \, {\hyp_q(\gerg)}' \, $  is the unital  $ R $--subalgebra
of  $ \hyp_q(\gerg) $  generated by the  $ K_i^{-1} $'s  and the
``rescaled quantum divided powers''  $ \, {(q_\alpha - 1)}^n
F_\alpha^{(n)} \, $,  $ \, {(q_i - 1)}^n \left( {{K_i ; \, 0} \atop
{n}} \right) \, $  and  $ \, {(q_\alpha - 1)}^n E_\alpha^{(n)} \, $
for all  $ \, n \in \N \, $.  Since  $ \, {[n]}_{q_\alpha\!}!
{\,\Big|}_{q=1} = n! = 0 \, $  iff  $ \, p \,\Big| n \, $,  one
argues like before that  $ \, {\hyp_q(\gerg)}'{\Big|}_{q=1} \, $
is generated by the corresponding specializations of  $ \,
{(q_\alpha - 1)}^{p^s} F_\alpha^{(p^s)} \, $,  $ \, {(q_i-1)}^{p^s}
\left( {{K_i ; \, 0} \atop {p^s}} \right) \, $  and  $ \, {(q_\alpha
- 1)}^{p^s} E_\alpha^{(p^s)} \, $  for all  $ \, s \in \N \, $  and all
positive roots  $ \alpha \, $:  \, this shows that the spectrum of  $ \,
{\hyp_q(\gerg)}'{\Big|}_{q=1} \, $  has (dimension 0 and) height 1, and
its cotangent Lie algebra  $ \, J \Big/ J^{\,2} \, $  (where  $ J $
is the augmentation ideal of  $ {\hyp_q(\gerg)}' {\Big|}_{q=1} $)
has basis
   \hbox{$ \, \Big\{ {(q_\alpha \!\! - \! 1)}^{p^s} \! F_\alpha^{(p^s)},
\, {(q_i \! - \! 1)}^{p^s} \!\! \left( {{K_i ; \, 0} \atop {p^s}}
\right) , \, {(q_\alpha \!\! - \! 1)}^{p^s} E^{(p^s)} \mod
(q \! - \! 1) {\hyp_q(\gerg)}' \mod \! J^{\,2\!} \;\Big|\,
s \in \N \Big\} \, $.}
  \allowbreak
Moreover,  $ \, \big( {\hyp_q(\gerg)}' \big)^{\!\vee} \, $  is
generated by  $ \, {(q_\alpha - 1)}^{p^s-1} F_\alpha^{(p^s)} \, $,
$ \, {(q_i-1)}^{p^s-1} \! \left( {{K_i ; \, 0} \atop {p^s}} \right)
\, $,  $ K_i^{-1} $  and  $ \, {(q_\alpha - 1)}^{p^s-1} E_\alpha^{(p^s)}
\, $  for all  $ s \, $,  $ i \, $  and  $ \alpha \, $:  \, in particular
$ \, \big( {\hyp_q(\gerg)}' \big)^{\!\vee} \subsetneqq \hyp_q(\gerg)
\, $,  \, and  $ \, \big( {\hyp_q(\gerg)}' \big)^{\!\vee}{\Big|}_{q=1}
\, $  is generated by the cosets modulo  $ (q-1) $  of the previous
elements, which in fact form a basis of the restricted Lie bialgebra
$ \gerk $  such that  $ \, \big( {\hyp_q(\gerg)}' \big)^{\!\vee}
{\Big|}_{q=1} = \, \u(\gerk) \, $.

\vskip7pt

  {\bf 7.6 The QFA  $ F_q[G] \, $.} \, In this and the following
sections we pass to look at Theorem 2.2 the other way round: namely,
we start from QFAs and produce QrUEAs.
                                             \par
   We begin with  $ \, G = {SL}_n \, $,  with the standard Poisson
structure,  for which an especially explicit description of the QFA
is available.  Namely, let  $ \, F_q[{SL}_n] \, $  be the unital
associative  \hbox{$R$--alge}bra  generated by \,  $ \{\,
\rho_{ij} \mid i, j = 1, \ldots, n \,\} $  \, with relations
  $$  \eqalignno{
   \rho_{ij} \rho_{ik} = q \, \rho_{ik} \rho_{ij} \; ,  \quad \quad
\rho_{ik} \rho_{hk}  &  = q \, \rho_{hk} \rho_{ik}  &   \forall\, j<k,
\, i<h  \quad  \cr
   \rho_{il} \rho_{jk} = \rho_{jk} \rho_{il} \; ,  \quad \quad
\rho_{ik} \rho_{jl} - \rho_{jl} \rho_{ik}  &  = \left( q - q^{-1}
\right) \, \rho_{il} \rho_{jk}  \hskip90pt  &   \forall\, i<j, \, k<l
\quad  \cr
   {det}_q (\rho_{ij}) := {\textstyle \sum\limits_{\sigma \in S_n}}
{(-q)}^{l(\sigma)} \rho_{1,\sigma(1)}  &  \rho_{2,\sigma(2)} \cdots
\rho_{n,\sigma(n)} = 1 \, .  &  \cr }  $$
  \indent   This is a Hopf algebra, with comultiplication,
counit and antipode given by
  $$  \Delta (\rho_{ij}) = {\textstyle \sum\limits_{k=1}^n} \rho_{ik}
\otimes \rho_{kj} \, ,  \quad  \epsilon(\rho_{ij}) = \delta_{ij} \, ,
\quad  S(\rho_{ij}) = {(-q)}^{i-j} \, {det}_q \! \left(
{(\rho_{hk})}_{h \neq j}^{k \neq i} \right)  $$
for all  $ \, i, j = 1, \dots, n \, $.  Let  $ \, \F_q[{SL}_n] :=
F(R) \otimes_R F_q[{SL}_n] \, $.  The set of ordered monomials
  $$  M \; := \; \left\{\; {\textstyle \prod\limits_{i>j}}
\rho_{ij}^{N_{ij}} {\textstyle \prod\limits_{h=k}} \rho_{hk}^{N_{hk}}
{\textstyle \prod\limits_{l<m}} \rho_{lm}^{N_{lm}} \;\bigg\vert\;
N_{st} \in \N  \;\; \forall \; s, t \; ; \; \min \big\{ N_{1,1},
\dots, N_{n,n} \big\} = 0 \;\right\}   \eqno (7.4)  $$
is an  $ R $--basis  of  $ F_q[{SL}_n] $  and an  $ F(R) $--basis of
$ \F_q[{SL}_n] $  (cf.~[Ga2], Theorem 7.4, and [Ga7],  Theorem 2.1{\it
(c)\/}).  Moreover,  $ \, F_q[{SL}_n] \, $  is a QFA (at  $ \, \h =
q - 1 \, $),  with  $ \, F_q[{SL}_n] \, {\buildrel \, q \rightarrow
1 \, \over \loongrightarrow} \, F[{SL}_n] \, $.

\vskip7pt

  {\bf 7.7 Computation of  $ \, {F_q[G]}^\vee \, $  and specialization
$ \, {F_q[G]}^\vee \,{\buildrel {q \rightarrow 1} \over \llongrightarrow}\,
U(\gerg^\times) \, $.} \, In this section we compute  $ \, {F_q[G]}^\vee
\, $  and its semiclassical limit (= specialization at  $ \, q = 1 \, $).
                           Note that    \break
  $$  M' \, := \, \left\{\; {\textstyle \prod\limits_{i>j}}
\rho_{ij}^{N_{ij}} {\textstyle \prod\limits_{h=k}} {(\rho_{hk} -
1)}^{N_{hk}} {\textstyle \prod\limits_{l<m}} \rho_{lm}^{N_{lm}}
\;\bigg\vert\; N_{st} \in \N \;\, \forall \; s, t \; ; \; \min
\big\{ N_{1,1}, \dots, N_{n,n} \big\} = 0 \,\right\}  $$
is an  $ R $--basis  of  $ F_q[{SL}_n] $  and an  $ F(R) $--basis
of  $ \F_q[{SL}_n] $;  then, from the definition of the counit, it
follows that  $ \, M' \setminus \{1\} \, $  is an  $ R $--basis  of
$ \, \hbox{\sl Ker}\,\big( \epsilon: F_q[{SL}_n] \longrightarrow R
\,\big) \, $.  Now, by definition  $ \, I := \hbox{\sl Ker} \left(
F_q[{SL}_n] {\buildrel \epsilon \over \llongtwoheadrightarrow} R
{\buildrel {q \mapsto 1} \over \llongtwoheadrightarrow} \Bbbk \right)
\, $,  \, whence  $ \, I = \hbox{\sl Ker}\,(\epsilon) + (q-1) \cdot
F_q[{SL}_n] \, $;  \, therefore  $ \, \big( M' \setminus \{1\}
\big) \cup \big\{ (q-1) \cdot 1 \big\} \, $  is an  $ R $--basis
of  $ I $,  hence  $ \; {(q-1)}^{-1} I \; $  has  $ R $--basis
$ \; {(q-1)}^{-1} \cdot \big( M' \setminus \{1\} \big) \cup
\{1\} \, $.  The outcome is that  $ \, {F_q[{SL}_n]}^\vee :=
\sum_{n \geq 0} {\Big( {(q-1)}^{-1} I \Big)}^n \, $  is just
the unital  $ R $--subalgebra of  $ \F_q[{SL}_n] $  generated by
 $$  \left\{\, r_{i{}j} := {\, \rho_{i{}j} - \delta_{i{}j} \, \over
\, q-1 \,} \;\bigg\vert\; i, j= 1, \dots, n \,\right\} \, .  $$
  \indent   Then one can directly show that this is a Hopf algebra,
and that  $ \, {F_q[{SL}_n]}^\vee \! {\buildrel \, q \rightarrow 1
\, \over \llongrightarrow} \, U({\gersl_n}^{\!*}) $  as predicted by
Theorem 2.2.  Details can be found in [Ga2], \S\S \, 2, 4, looking at
the algebra  $ \widetilde{F}_q[{SL}_n] $  considered therein, up to
the following changes.  The algebra which is considered in
[{\it loc.~cit.}]  has generators  $ \, \displaystyle{
\big( 1 + q^{-1} \big)^{\delta_{i{}j}} {\, \rho_{i{}j} - \delta_{i{}j}
\, \over \, q - q^{-1} \,}} \, $  ($ \, i, j = 1, \dots, n \, $)  instead
of our  $ \, r_{i{}j} \, $'s  (they coincide iff  $ \, i = j \, $)  and
also generators  $ \, \rho_{i{}i} = 1 + (q-1) \, r_{i{}i} \, $  ($ \, i =
1, \dots, n \, $);  then the presentation in \S 2.8 of  [{\it loc.~cit.}]
must be changed accordingly; computing the specialization then goes
exactly the same, and gives the same result   --- specialized
generators are rescaled, though, compared with the standard
ones given in  [{\it loc.~cit.}],  \S 1.
%
%
 \eject   
  We sketch the case of  $ \, n = 2 \, $  (see also [FG]).  Using
notation $ \, \text{a} := \rho_{1,1} \, $,  $ \, \text{b} :=
\rho_{1,2} \, $,  $ \, \text{c} := \rho_{2,1} \, $,  $ \, \text{d}
:= \rho_{2,2} \, $,  we have the relations
 $$  \displaylines{
  \text{a} \, \text{b} = q \, \text{b} \, \text{a} \, ,  \;\hskip27pt
\text{a} \, \text{c} = q \, \text{c} \, \text{a} \, , \;\hskip27pt
\text{b} \, \text{d} = q \, \text{d} \, \text{b} \, ,  \;\hskip27pt
\text{c} \, \text{d} = q \, \text{d} \, \text{c} \, , \cr
  {}  \hskip9pt  \text{b} \, \text{c} = \text{c} \, \text{b} \, ,
\;\hskip41pt  \text{a} \, \text{d} - \text{d} \, \text{a} = \big( q -
q^{-1} \big) \text{b} \, \text{c} \, ,  \;\hskip31pt  \text{a} \, \text{d}
- q \, \text{b} \, \text{c} = 1 \cr }  $$
holding in  $ F_q[{SL}_2] $  and in  $ \F_q[{SL}_2] $,  with
  $$  \displaylines{
   \Delta(\text{a}) = \text{a} \otimes \text{a} + \text{b} \otimes \text{c} ,
\; \Delta(\text{b}) = \text{a} \otimes \text{b} + \text{b} \otimes
\text{d} , \;  \Delta(\text{c}) = \text{c} \otimes \text{a} + \text{d}
\otimes \text{c} , \;  \Delta(\text{d}) = \text{c} \otimes \text{b} +
\text{d} \otimes \text{d}  \cr
   \epsilon(\text{a}) = 1 ,  \hskip1pt  \epsilon(\text{b}) = 0 ,
\hskip1pt  \epsilon(\text{c}) = 0 ,  \hskip1pt  \epsilon(\text{d})
= 1 ,  \hskip9,7pt  S(\text{a}) = \text{d} ,  \hskip1pt  S(\text{b})
= - q^{-1} \text{b} ,  \hskip1pt  S(\text{c}) = - q^{+1} \text{c} ,
\hskip1pt  S(\text{d}) = \text{a} \, .  \cr }  $$
Then the elements  $ \, H_+ := r_{1,1} = \displaystyle{\, \text{a} - 1 \,
\over \, q - 1\,} \, $,  $ \, E := r_{1,2} = \displaystyle{\, \text{b} \,
\over \, q - 1\,} \, $, $ \, F := r_{2,1} = \displaystyle{\, \text{c} \,
\over \, q - 1\,} \, $ and $ \, H_- := r_{2,2} = \displaystyle{\, \text{d}
- 1 \, \over \, q - 1\,} \, $  generate  $ {F_q[{SL}_2]}^\vee $:
these generators have relations
  $$  \displaylines{
   H_+ E = q \, E H_+ + E \, ,  \;\;  H_+ F = q \, F H_+ + F \, ,  \;\;
E H_- = q \, H_- E + E \, ,  \;\; F H_- = q \, H_- F + F \, ,  \cr
   E F = F E \, ,  \;\;\;  H_+ H_- - H_- H_+ = \big( q - q^{-1} \big) E F
\, , \;\;\;  H_- + H_+ = (q-1) \big( q \, E F - H_+ H_- \big)  \cr }  $$
and Hopf operations given by
  $$  \displaylines{
   \Delta(H_+) = H_+ \otimes 1 + 1 \otimes H_+ + (q-1) \big( H_+ \otimes
H_+ + E \otimes F \big) \, ,  \quad  \epsilon(H_+) = 0 \, ,  \quad  S(H_+)
= H_-  \cr
   {} \, \Delta(E) = E \otimes 1 + 1 \otimes E + (q-1) \big( H_+
\otimes E + E \otimes H_- \big) \, ,  \,\;\;\quad  \epsilon(E) = 0 \, ,
\;\;\quad  S(E) = - q^{-1} E  \cr
   {} \, \Delta(F) = F \otimes 1 + 1 \otimes F + (q-1) \big( F \otimes
H_+ + H_- \otimes F \big) \, ,  \,\;\;\quad  \epsilon(F) = 0 \, ,
\;\;\quad  S(F) = - q^{+1} F  \cr
   \Delta(H_-) = H_- \otimes 1 + 1 \otimes H_- + (q-1) \big( H_- \otimes
H_- + F \otimes E \big) \, ,  \quad  \epsilon(H_-) = 0 \, ,  \quad
S(H_-) = H_+  \cr }  $$
from which one easily checks that  $ \, {F_q[{SL}_2]}^\vee
\,{\buildrel \, q \rightarrow 1 \, \over \llongrightarrow}\,
U({\gersl_2}^{\! *}) \, $  as co-Poisson Hopf algebras, for a
co-Poisson Hopf algebra isomorphism
  $$  {F_q[{SL}_2]}^\vee \Big/ (q-1) \, {F_q[{SL}_2]}^\vee
\;{\buildrel \cong \over \llongrightarrow}\; U({\gersl_2}^{\! *})  $$
exists, given by:  $ \;\, H_\pm \mod (q-1) \mapsto \pm \text{h} \, $,
$ \, E \mod (q-1) \mapsto \text{e} \, $,  $ \, F \mod (q-1) \mapsto
\text{f} \, $;  that is, $ \, {F_q[{SL}_2]}^\vee \, $  specializes
to  $ \, U({\gersl_2}^{\! *}) \, $  {\sl as a co-Poisson Hopf
algebra},  q.e.d.
                                             \par
   Finally, the general case of any semisimple group  $ \, G = G^\tau
\, $,  with the Poisson structure induced from the Lie bialgebra
structure of  $ \, \gerg = \gerg^\tau \, $,  can be treated in a
different way.  Following [Ga1], \S\S 5--6,  $ \F_q[G] $  can be
embedded into a (topological) Hopf algebra  $ \, \Bbb{U}_q(\gerg^*) =
{\Bbb U}_{q,\varphi}^{\scriptscriptstyle M}(\gerg^*) \, $,  so
that the image of the integer form  $ F_q[G] $
lies into a suitable (topological) integer form  $ \,
{\Cal U}_{q,\varphi}^{\scriptscriptstyle M}(\gerg^*) \, $
of  $ \, \Bbb{U}_q(\gerg^*) \, $.  Now, the analysis given in
[{\it loc.~cit.\/}],  when carefully read, shows that  $ \,
{F_q[G]}^\vee = \F_q[G] \cap {{\Cal U}_{q,\varphi}^{\scriptscriptstyle M}
(\gerg^*)}^\vee \, $;  moreover, the latter (intersection) algebra
``almost'' coincides   --- it is its closure in a suitable topology ---
with the integer form  $ \, {\Cal F}_q[G] \, $  considered in  [{\it
loc.~cit.\/}]:  in particular, they have the same specialization
at  $ \, q = 1 \, $.  Since in addition  $ \, {\Cal F}_q[G] \, $
does specialize to  $ U(\gerg^*) $,  the same is true for
$ {F_q[G]}^\vee $,  q.e.d.
                                             \par
   The last point to stress is that, once more, the whole analysis
above is valid for  $ \, p := \Char(\Bbbk) \geq 0 \, $, i.e.~also
for  $ \, p > 0 \, $,  \, which was not granted by Theorem 2.2.
%
%
 \eject   

  {\bf 7.8 The identity  $ \; {\Big({F_q[G]}^\vee\Big)}' = F_q[G] \, $.}
\, In this section we verify the validity of that part of Theorem
2.2{\it (b)}  claiming that  $ \; H \in \QFA \,\Longrightarrow\,
{\big(H^\vee\big)}' = H \; $  for  $ \, H = F_q[G] \, $  as above;
moreover we show that this holds for  $ \, p > 0 \, $  too.  We
begin with  $ \, G = {SL}_n \, $.
                                                     \par
   From  $ \; \Delta(\rho_{i{}j}) = \sum\limits_{k=1}^n \rho_{i,k}
\otimes \rho_{k,j} \, $,  \; we get  $ \; \Delta^{N} (\rho_{i{}j})
= \!\! \sum\limits_{k_1, \ldots, k_{N-1} = 1}^n \hskip-15pt \rho_{i,k_1}
\otimes \rho_{k_1,k_2} \otimes \cdots \otimes \rho_{k_{N-1},j} \, $,
\, by repeated iteration, whence a simple computation yields
  $$  \delta_N(r_{i{}j}) = \sum_{k_1, \ldots, k_{N-1} = 1}^n \hskip-9pt
{(q-1)}^{-1} \cdot \big( (q-1) \, r_{i,k_1} \otimes (q-1) \,
r_{k_1,k_2} \otimes \cdots \otimes (q-1) \, r_{k_{N-1},j} \big)
\qquad \;  \forall \hskip5pt  i, j  $$
%
%
%
%
so that
  $$  \delta_N \big( (q-1) r_{i{}j} \big) \in {(q-1)}^N {F_q[{SL}_n]}^\vee
\setminus {(q-1)}^{N+1} {F_q[{SL}_n]}^\vee  \qquad \hskip11pt \forall
\hskip5pt i, j \, .   \eqno (7.5)  $$
   \indent   Now, consider again the set  $ \; M' := \bigg\{\,
\prod\limits_{i>j} \rho_{ij}^{N_{ij}} \prod\limits_{h=k}
{(\rho_{hk} -1)}^{N_{hk}} \prod\limits_{l<m} \rho_{lm}^{N_{lm}}
\,\bigg\vert\, N_{st} \in \N \; \forall \, s, t \, ; $\break
$ \, \min \big\{\, N_{1,1}, \dots, N_{n,n} \,\big\} = 0 \,\bigg\} \, $:
\; since this is an  $ R $--basis  of $ F_q[{SL}_n] $,  we
have also that
 $$  M'' \, := \; \left\{\; {\textstyle \prod\limits_{i>j}}
r_{ij}^{N_{ij}} {\textstyle \prod\limits_{h=k}} r_{hk}^{N_{hk}}
{\textstyle \prod\limits_{l<m}} r_{lm}^{N_{lm}} \;\bigg\vert\;
N_{st} \in \N  \;\; \forall \; s, t \, ; \; \min \big\{\, N_{1,1},
\dots, N_{n,n} \,\big\} = 0 \;\right\}  $$
is an  $ R $--basis  of  $ {F_q[{SL}_n]}^\vee $.  This and  $ (7.5) $
above imply that  $ \, {\big( {F_q[{SL}_n]}^\vee \big)}' $  is the
unital  $ R $--subalgebra of  $ \, \F_q[{SL}_n] \, $  generated by
the set  $ \, \big\{\, (q-1) r_{i{}j} \,\big\vert\, i, j= 1, \dots,
            n \,\big\} \, $;  since\break
$ \, (q-1) \, r_{i{}j} = \rho_{i{}j} - \delta_{i{}j} \, $,  the latter
algebra does coincide with  $ \, F_q[{SL}_n] \, $,  as expected.
                                             \par
   For the general case of any semisimple group  $ \, G = G^\tau
\, $,  the result can be obtained again by looking at the immersions
$ \, \F_q[G] \subseteq \Bbb{U}_q(\gerg^*) \, $  and  $ \, F_q[G] \subseteq
{\Cal U}_{q,\varphi}^{\scriptscriptstyle M}(\gerg^*) \, $,
and at the identity  $ \, {F_q[G]}^\vee = \F_q[G] \cap
{{\Cal U}_{q,\varphi}^{\scriptscriptstyle M}(\gerg^*)}^\vee \, $
(cf.~\S 7.6);  if we try to compute  $ \, {\Big(
{{\Cal U}_{q,\varphi}^{\scriptscriptstyle M}(\gerg^*)}^\vee \Big)}' \, $
(noting that  $ \, {\big( {\Cal U}_{q,\varphi}^{\scriptscriptstyle M}
(\gerg^*) \big)}^\vee \, $  is a QrUEA), we have just to apply much the
like methods as for  $ \, {U_q(\gerg)}' \, $,  \, thus finding a similar
result; then from this and the identity  $ \, {F_q[G]}^\vee = \F_q[G]
\cap {{\Cal U}_{q,\varphi}^{\scriptscriptstyle M}(\gerg^*)}^\vee \, $
we eventually find  $ \, {\Big({F_q[G]}^\vee\Big)}' = F_q[G] \, $,
q.e.d.
                                             \par
   We'd better point out once more that the previous analysis is valid
for  $ \, p := \Char(\Bbbk) \geq 0 \, $, i.e.~also for  $ \, p > 0 \, $,
\, so the outcome is stronger than what ensured by Theorem 2.2.

\vskip7pt

   {\sl $ \underline{\hbox{\it Remark}} $:}  \; Formula (7.4)
gives an explicit  $ R $--basis  $ M $  of  $ F_q[{SL}_2] $.
By direct computation one sees that  $ \, \delta_n(\mu)
\in {F_q[{SL}_2]}^{\otimes n} \setminus (q-1) \,
{F_q[{SL}_2]}^{\otimes n} $  for all  $ \, \mu \in M \setminus
\{1\} \, $  and  $ \, n \in \N \, $,  whence  $ \, {F_q[{SL}_2]}'
= R \cdot 1 \, $,  which implies  $ \, {\big( {F_q[{SL}_2]}' \big)}_F
= F(R) \cdot 1 \subsetneqq \F_q[{SL}_2] \, $  (cf.~the Remark after
Corollary 4.6) and also  $ \, {\big( {F_q[{SL}_2]}' \big)}^\vee =
\, R \cdot 1 \, \subsetneqq \, F_q[{SL}_2] \, $.

\vskip7pt

  {\bf 7.9 Drinfeld's functors and  $ L $--operators  for
$ U_q(\gerg) $  when  $ \gerg $  is classical.} \, Let now
$ \Bbbk $  have characteristic zero, and let  $ \gerg $  be
a finite dimensional semisimple Lie algebra over  $ \Bbbk $
whose simple Lie subalgebra are all of classical type.  It is
known from [FRT2] that in this case  $ \Bbb{U}^P_q(\gerg) $
(where the subscript  $ P $  means that we are taking a
``simply-connected'' quantum group) admits an alternative
presentation, in which the generators are the so-called
$ L $--operators,  denoted  $ l_{i,j}^{(\varepsilon)} $  with
$ \, \varepsilon = \pm 1 \, $  and  $ i $,  $ j $  ranging
in a suitable set of indices (see [FRT2], \S 2).  Now, if we
consider instead the  $ R $--subalgebra  $ H $  generated by
the  $ L $--operators,  we get at once from the very description
of the relations between the  $ l^{(\varepsilon)}_{i,j} $'s
given in [FRT2] that  $ H $  is a  {\sl Hopf\/}  $ R $--subalgebra
of  $ \Bbb{U}^P_q(\gerg) $,  and more precisely it is a QFA for
the connected simply-connected dual Poisson group  $ G^\star \, $.
                                          \par
   When computing  $ H^\vee $,  it is generated by the elements
$ \, {(q-1)}^{-1} l^{(\varepsilon)}_{i,j} \, $;  \, even more,
the elements  $ \, {(q-1)}^{-1} l^{(+)}_{i,i+1} \, $  and  $ \,
{(q-1)}^{-1} l^{(-)}_{i+1,i} \, $  are enough to generate.  Now,
Theorem 12 in [FRT2] shows that these latter generators are simply
multiples of the Chevalley generators of  $ U^P_q(\gerg) $  (in
the sense of Jimbo, Drinfeld, etc.), by a coefficient  $ \, \pm
q^s \big( 1 + q^{-1} \big) \, $,  for some  $ \, s \in \Z \, $,
times a ``toral'' generator: this proves directly that  $ H^\vee $
is a QrUEA associated to  $ \gerg \, $,  that is the dual Lie
bialgebra of  $ G^* $,  as prescribed by Theorem 2.2.
Conversely, if we start from  $ U_q^P(\gerg) $,  again
Theorem 12 of [FRT2] shows that the  $ {\big( q - q^{-1}
\big)}^{-1} l^{(\varepsilon)}_{i,j} $'s  are quantum root vectors
in  $ U_q^P(\gerg) $.  Then when computing  $ {U_q^P(\gerg)}' $
we can shorten a lot the analysis in \S 5.3, because the explicit
expression of the coproduct on the  $ L $--operators  given in
[FRT2]   --- roughly,  $ \Delta $  is given on them by a standard
``matrix coproduct'' ---   tells us directly that all the
$ {\big( 1 + q^{-1} \big)}^{-1} l^{(\varepsilon)}_{i,j} $'s
do belong to  $ {U_q^P(\gerg)}' $,  and again by a PBW argument
we conclude that  $ {U_q^P(\gerg)}' $  is generated by these
rescaled  $ L $--operators,  i.e.~the  $ {\big( 1 + q^{-1}
\big)}^{-1} l^{(\varepsilon)}_{i,j} \, $.
                                          \par
   Therefore, we can say in short that shifting from  $ H $  to
$ H^\vee $  or from  $ U_q^P(\gerg) $  to  $ {U_q^P(\gerg)}' $
essentially amounts   --- up to rescaling by irrelevant factors
(in that they do not vanish at  $ \, q = 1 \, $)  ---   to
switching from the presentation of  $ \Bbb{U}^P_q(\gerg) $
via  $ L $--operators  (after [FRT2]) to the presentation of
Serre-Chevalley type (after Drinfeld and Jimbo), and conversely.
See also the analysis in [Ga7] for the cases  $ \, \gerg =
\frak{gl}_n \, $  and  $ \, \gerg = \frak{sl}_n \, $.

\vskip7pt

  {\bf 7.10 The cases  $ U_q(\frak{gl}_n) \, $,  $ F_q[{GL}_n] \, $
and  $ F_q[M_n] \, $.} \, In [Ga2], \S 5.2, a certain algebra
$ U_q(\frak{gl}_n) $  is considered as a quantization of
$ \frak{gl}_n \, $;  \, due to their strict relationship, from
the analysis we did for the case of  $ \gersl_n $  one can easily
deduce a complete description of  $ {U_q(\frak{gl}_n)}' $  and its
specialization at  $ \, q = 1 \, $,  and also verify that  $ \,
\big({U_q(\frak{gl}_n)}'\big)^{\!\vee} = \, U_q(\frak{gl}_n) \, $.
                                             \par
   Similarly, we can consider the unital associative $R$--algebra  $ \,
F_q[M_n] \, $  with generators  $ \, \rho_{ij} \, $  ($ i $,  $ j = 1,
\ldots, n \, $)  and relations  $ \; \rho_{ij} \rho_{ik} = q \, \rho_{ik}
\rho_{ij} \, $,  $ \, \rho_{ik} \rho_{hk} = q \, \rho_{hk} \rho_{ik} \, $
(for all  $ \, j < k \, $,  $ \, i < h \, $),  $ \; \rho_{il} \rho_{jk} =
\rho_{jk} \rho_{il} \, $,  $ \, \rho_{ik} \rho_{jl} - \rho_{jl} \rho_{ik}
= \left( q - q^{-1} \right) \, \rho_{il} \rho_{jk} \, $  (for all  $ \,
i < j \, $,  $ \, k < l \, $)   --- i.e.~like for  $ {SL}_n \, $,  \,
but for skipping the last relation.  This is the celebrated standard
quantization of  $ F[M_n] $,  the function algebra of the variety
$ M_n $  of  ($ n \times n $)--matrices  over  $ \Bbbk \, $:  \, it
is a  $ \Bbbk $--bialgebra,  whose structure is given by formulas
$ \; \Delta (\rho_{ij}) = \sum_{k=1}^n \rho_{ik} \otimes \rho_{kj} \, $,
$ \, \epsilon(\rho_{ij}) = \delta_{ij} \, $  (for all  $ \, i $,  $ j
= 1, \dots, n \, $)  again, but it is  {\sl not a Hopf algebra}.  The
quantum determinant  $ \; {det}_q (\rho_{ij}) := \sum_{\sigma \in S_n}
{(-q)}^{l(\sigma)} \rho_{1,\sigma(1)} \, \rho_{2,\sigma(2)} \cdots
\rho_{n,\sigma(n)} \; $  is central in  $ F_q[M_n] $,  so by standard
theory we can extend  $ F_q[M_n] $  by adding a formal inverse to
$ \, {det}_q (\rho_{ij}) \, $,  \, thus getting a larger algebra
$ \, F_q[{GL}_n] := F_q[M_n] \big[ {{det}_q(\rho_{ij})}^{-1} \big]
\, $:  \, this is now a Hopf algebra, with antipode  $ \; S(\rho_{ij})
= {(-q)}^{i-j} \, {det}_q \! \left( {(\rho_{hk})}_{h \neq j}^{k \neq i}
\right) \; $  (for all  $ \, i $,  $ j = 1, \dots, n \, $),  the
well-known standard quantization of  $ F[{GL}_n] $,  due to Manin
(see [Ma]).
                                             \par
   Applying Drinfeld's functor  $ (\ )^\vee $  w.r.t.~$ \, \h := (q - 1)
\, $  at  $ F_q[{GL}_n] $  we can repeat stepwise the analysis made for
$ F_q[{SL}_n] $:  \, then we have that  $ {F_q[{GL}_n]}^\vee $  is
generated by the  $ r_{i{}j} $'s  and  $ \, {(q-1)}^{-1} \big(
{det}_q(\rho_{ij}) - 1 \big) \, $,  the sole real difference being
the lack of the relation  $ \, {det}_q(\rho_{ij}) = 1 \, $,  \,
which implies one relation less among the  $ r_{i{}j} $'s  inside
$ {F_q[{GL}_n]}^\vee $,  hence also one relation less among their
cosets modulo  $ (q-1) $.  The outcome is pretty similar, in particular
$ \, {F_q[{GL}_n]}^\vee{\Big|}_{q=1} = \, U({\frak{gl}_n}^{\!\! *})
\, $  (cf.~[Ga2], \S 6.2).  Even more, we can do the same with
$ F_q[M_n] \, $:  \, things are even easier, because we have only
the  $ r_{i{}j} $'s  alone which generate  $ {F_q[M_n]}^\vee $,  with
no relation coming from the relation  $ \, {det}_q(\rho_{ij}) = 1 \, $;
\, nevertheless at  $ \, q = 1 \, $  the relations among the cosets
of the  $ r_{i{}j} $'s  are exactly the same as in the case of
$ {F_q[{GL}_n]}^\vee{\Big|}_{q=1} $,  whence we get  $ \,
{F_q[M_n]}^\vee{\Big|}_{q=1} = \, U({\frak{gl}_n}^{\!\! *}) \, $.
{\sl In particular, we get that  $ {F_q[M_n]}^\vee{\Big|}_{q=1} $  is
a Hopf algebra, although both  $ F_q[M_n] $  and  $ {F_q[M_n]}^\vee $
are only bialgebras,  {\sl not\/}  Hopf algebras: so this gives a
non-trivial explicit example of what claimed in the first part
of Theorem 3.7.}
                                             \par
   Finally, an analysis of the relationship between Drinfeld functors
and  $ L $--operators  about  $ \Bbb{U}^P_q(\frak{gl}_n) $  can be
done again, exactly like in \S 7.9, leading to entirely similar
results.

\vskip1,7truecm

\centerline {\bf \S \; 8 \  Third example: quantum three-dimensional
Euclidean group }

\vskip10pt

  {\bf 8.1  The classical setting.} \, Let  $ \Bbbk $  be any field of
characteristic  $ \, p \geq 0 \, $.  Let  $ \, G := {E}_2(\Bbbk) \equiv
{E}_2 \, $,  the three-dimensional Euclidean group; its tangent Lie
algebra  $ \, \gerg = \gere_2 \, $  is generated by  $ \, f $,
$ h $,  $ e \, $  with relations  $ \, [h,e] = 2 e $,  $ [h,f] =
-2 f $,  $ [e,f] = 0 \, $.  The formulas  $ \, \delta(f) = h \otimes
f - f \otimes h \, $,  $ \, \delta(h) = 0 \, $,  $ \, \delta(e) = h
\otimes e - e \otimes h \, $,  make  $ {\gere}_2 $  into a Lie bialgebra,
hence  $ {E}_2 $  into a Poisson group.  These also give a presentation
of the co-Poisson Hopf algebra  $ U({\gere}_2) $  (with standard Hopf
structure).  If  $ \, p > 0 \, $,  \, we consider on  $ \gere_2 $  the
$ p $--operation  given by  $ \, e^{[p\,]} = 0 \, $,  $ \, f^{[p\,]}
= 0 \, $,  $ \, h^{[p\,]} = h \, $.
                                             \par
  On the other hand, the function algebra  $ F[{E}_2] $  is the unital
associative commutative  $ \Bbbk $--algebra  with generators  $ \, b $,
$ a^{\pm 1} $,  $ c $,  with Poisson Hopf algebra structure given by
  $$  \displaylines{
   \Delta(b) = b \otimes a^{-1} + a \otimes b \, ,  \hskip29pt
\Delta\big(a^{\pm 1}\big) = a^{\pm 1} \otimes a^{\pm 1} \, ,
\hskip29pt  \Delta(c) = c \otimes a + a^{-1} \otimes c  \cr
   \epsilon(b) = 0  \, ,  \hskip7pt  \epsilon\big(a^{\pm 1}\big) = 1
\, ,  \hskip7pt  \epsilon(c) = 0 \, ,  \hskip39pt  S(b) = -b \, ,
\hskip7pt  S\big(a^{\pm 1}\big) = a^{\mp 1} \, ,  \hskip7pt
S(c) = -c  \cr
   \big\{a^{\pm 1},b\big\} = \pm a^{\pm 1} b \, ,  \hskip33pt
\big\{a^{\pm 1},c\big\} = \pm a^{\pm 1} c \, ,  \hskip33pt
\{b,c\} = 0  \cr }  $$
  \indent   We can realize  $ E_2 $  as  $ \, E_2 = \big\{\, (b,a,c)
\,\big\vert\, b, c \in k, a \in \Bbbk \setminus \{0\} \,\big\} \, $,
with group operation
  $$  (b_1,a_1,c_1) \cdot (b_2,a_2,c_2) = \big( b_1 a_2^{-1} + a_1 b_2
\, , \, a_1 a_2 \, , \, c_1 a_2 + a_1^{-1} c_2 \big) \, ;  $$
in particular the centre of  $ E_2 $  is simply  $ \, Z := \big\{
(0,1,0), (0,-1,0) \big\} \, $,  so there is only one other connected
Poisson group having  $ \gere_2 $  as Lie bialgebra, namely the adjoint
group  $ \, {}_aE_2 := E_2 \Big/ Z \, $  (the left subscript
$ a $  stands for ``adjoint'').  Then  $ F[{}_aE_2] $  coincides
with the Poisson Hopf subalgebra of  $ F[{}_aE_2] $  spanned
by products of an even number of generators, i.e. monomials of
even degree: as a unital subalgebra, this is generated by
$ b a $,  $ a^{\pm 2} $, and $ a^{-1} c $.
                                                \par
  The dual Lie bialgebra  $ \, \gerg^* = {\gere_2}^{\!*} \, $  is
the Lie algebra with generators  $ \, \text{f} $,  $ \text{h} $,
$ \text{e} \, $,  and relations  $ \, [\text{h},\text{e}] = 2 \text{e} $,
$ [\text{h},\text{f}\,] = 2 \text{f} $,  $ [\text{e},\text{f}\,] = 0 \, $,
with Lie cobracket given by  $ \, \delta(\text{f}\,) = \text{f} \otimes
\text{h} - \text{h} \otimes \text{f} $,  $ \, \delta(\text{h}) = 0 $,
$ \, \delta(\text{e}) = \text{h} \otimes \text{e} - \text{e} \otimes
\text{h} \, $  (we choose as generators  $ \, \text{f} := f^* \, $,
$ \, \text{h} := 2 h^* \, $,  $ \, \text{e} := e^* \, $,  where  $ \,
\big\{ f^*, h^*, e^* \big\} \, $  is the basis of  $ {\gere_2}^{\! *} $
which is the dual of the basis  $ \{ f, h, e \} $  of  $ \gere_2 \, $).
If  $ \, p > 0 \, $,  the  $ p $--operation  of  $ {\gere_2}^{\!*} $
is given by  $ \, \text{e}^{[p\,]} = 0 \, $,  $ \, \text{f}^{\,[p\,]}
= 0 \, $,  $ \, \text{h}^{[p\,]} = \text{h} \, $.  All this again gives
a presentation of  $ \, U \left( {{\gere}_2}^{\!*} \right) \, $  too.
The simply connected algebraic Poisson group with tangent Lie bialgebra
$ {{\gere}_2}^{\!*} $  can be realized as the group of pairs of matrices
 $$  {{}_s{E}_2}^{\! *} := \Bigg\{\,
\bigg( \bigg( \matrix  z^{-1} & 0 \\  y & z  \endmatrix \bigg) ,
\bigg( \matrix  z & x \\  0 & z^{-1}  \endmatrix \bigg) \bigg)
\,\Bigg\vert\, x, y \in k, z \in \Bbbk \setminus \{0\} \,\Bigg\} \; ;  $$
this group has centre  $ \, Z := \big\{ (I,I), (-I,-I) \big\} \, $,
so there is only one other (Poisson) group with Lie (bi)algebra
$ {{\gere}_2}^{\! *} \, $,  namely the adjoint group
$ \, {{}_a{E}_2}^{\! *} := {{}_s{E}_2}^{\! *} \Big/ Z \, $.
                                                  \par
  Therefore  $ F\big[{{}_s{E}_2}^{\! *}\big] $  is the unital
associative commutative  $ \Bbbk $--algebra  with generators  $ \, x $,
$ z^{\pm 1} $,  $ y $,  with Poisson Hopf structure given by
  $$  \displaylines{
   \Delta(x) = x \otimes z^{-1} + z \otimes x \, ,  \hskip29pt
\Delta\big(z^{\pm 1}\big) = z^{\pm 1} \otimes z^{\pm 1} \, ,
\hskip29pt  \Delta(y) = y \otimes z^{-1} + z \otimes y  \cr
   \epsilon(x) = 0  \, ,  \hskip7pt  \epsilon\big(z^{\pm 1}\big) = 1
\, ,  \hskip7pt  \epsilon(y) = 0 \, ,  \hskip39pt  S(x) = -x \, ,
\hskip7pt  S\big(z^{\pm 1}\big) = z^{\mp 1} \, ,  \hskip7pt
S(y) = -y  \cr
   \{x,y\} = 0 \, ,  \hskip33pt  \big\{z^{\pm 1},x\big\} =
\pm z^{\pm 1} x \, ,  \hskip33pt  \big\{z^{\pm 1},y\big\} =
\mp z^{\pm 1} y  \cr }  $$
(N.B.: with respect to this presentation, we have  $ \, \text{f} =
{\partial_y}{\big\vert}_e \, $,  $ \, \text{h} = z \,
{\partial_z}{\big\vert}_e \, $,  $ \, \text{e} =
{\partial_x}{\big\vert}_e \, $,  where  $ e $  is the
identity element of  $ \, {{}_s{E}_2}^{\! *} \, $).  Moreover,
$ F\big[{{}_a{E}_2}^{\! *}\big] $  can be identified with the
Poisson Hopf subalgebra of  $ F\big[{{}_s{E}_2}^{\! *}\big] $
spanned by products of an even number of generators, i.e.{}
monomials of even degree: this is generated, as a unital
subalgebra, by  $ x z $,  $ z^{\pm 2} $,  and  $ z^{-1} y $.

\vskip7pt

  {\bf 8.2 The QrUEAs  $ \, U_q^s(\gere_2) \, $  and  $ \,
U_q^a(\gere_2) \, $.} \, We turn now to quantizations: the
situation is much similar to the  $ \gersl_2 $  case, so we
follow the same pattern, but we stress a bit more the occurrence
of different groups sharing the same tangent Lie bialgebra.
                                                  \par
  Let  $ R $  be a domain and let  $ \, \h \in R \setminus \{0\} \, $
and  $ \, q := \h + 1 \in R \, $  be like in \S 7.2.
                                                  \par
   Let $ \, \Bbb{U}_q(\gerg) = \Bbb{U}_q^s({\gere}_2) \, $ 
(where the superscript  $ s $  stands for ``simply connected'')
be the associative unital  $ F(R) $--algebra  with generators
$ \, F $,  $ L^{\!\pm 1} $,  $ E $,  and relations
  $$  L L^{-1} = 1 = L^{-1} L \, ,  \;\quad  L^{\!\pm 1} F =
q^{\mp 1} F L^{\pm 1} \, ,  \;\quad  L^{\pm 1} E = q^{\pm 1}
E L^{\pm 1} \, , \;\quad  E F = F E \, .  $$
This is a Hopf algebra, with Hopf structure given by
  $$  \displaylines{
     \Delta(F) = F \otimes L^{-2} + 1 \otimes F \, ,  \hskip19pt
\Delta \big( L^{\pm 1} \big) = L^{\pm 1} \otimes L^{\pm 1} \, ,
\hskip19pt  \Delta(E) = E \otimes 1 + L^2 \otimes E  \cr
\epsilon(F) = 0 \, ,  \hskip5pt  \epsilon \big( L^{\pm 1} \big) = 1
\, ,  \hskip5pt  \epsilon(E) = 0 \, ,  \hskip15pt  S(F) = - F L^2 ,
\hskip5pt  S \big( L^{\pm 1} \big) = L^{\mp 1} ,  \hskip5pt
S(E) = - L^{-2} E \, .  \cr }  $$
Then let  $ \, U_q^s(\gere_2) \, $  be the  $ R $--subalgebra
of  $ \Bbb{U}_q^s(\gere_2) $  generated by  $ \, F $,  $ D_\pm :=
     \displaystyle{\, L^{\pm 1} - 1 \, \over \, q - 1 \,} $,\break
$ E $.  From the definition of  $ \Bbb{U}_q^s(\gere_2) $  one gets a
presentation of  $ U_q^s(\gere_2) $  as the associative unital
algebra with generators  $ \, F $,  $ D_\pm $,  $ E $  and relations
  $$  \displaylines{
   D_+ E = q E D_+ + E \, ,  \quad  F D_+ = q D_+ F + F \, , \quad
E D_- = q D_- E + E \, ,  \quad  D_- F = q F D_- + F  \cr
   E F = F E \, ,  \;\;\qquad  D_+ D_- = D_- D_+ \, ,  \;\;\qquad
D_+ + D_- + (q-1) D_+ D_- = 0  \cr }  $$  
with a Hopf structure given by
  $$  \displaylines{
   \Delta(E) = E \otimes 1 + 1 \otimes E + 2(q-1) D_+ \otimes E +
{(q-1)}^2 \cdot D_+^2 \otimes E  \cr
  \Delta(D_\pm) = D_\pm \otimes 1 + 1 \otimes D_\pm + (q-1) \cdot D_\pm
\otimes D_\pm  \cr
  \Delta(F) = F \otimes 1 + 1 \otimes F + 2(q-1) F \otimes D_- +
{(q-1)}^2 \cdot F \otimes D_-^2  \cr 
}  $$
 \vskip-15pt
  $$  \matrix
   \epsilon(E) = 0 \, ,  &  \quad \qquad  S(E) = -E - 2(q-1) D_- E
-{(q-1)}^2 D_-^2 E \, \phantom{.}  \\
   \epsilon(D_\pm) = 0 \, ,  &  \quad \qquad  S(D_\pm) = D_\mp \,
\phantom{.} \\
   \epsilon(F) = 0 \, ,  &  \quad \qquad  S(F) = -F - 2(q-1) F D_+ -
{(q-1)}^2 F D_+^2 \, .  \\
      \endmatrix  $$  
   \indent   The ``adjoint version'' of  $ \Bbb{U}_q^s(\gere_2) $  is the unital
subalgebra  $ \Bbb{U}_q^a(\gere_2) $  generated  by  $ F $,  $ K^{\pm 1} :=
L^{\pm 2} $,  $ E $,  which is clearly a Hopf subalgebra.  It also has
an  $ R $--integer  form  $ U_q^a(\gere_2) $,  the unital
$ R $--subalgebra  generated by  $ \, F $,  $ \, H_\pm :=
\displaystyle{\, K^{\pm 1} - 1 \, \over \, q - 1 \,} $,  $ \, E \, $:
this has relations
 $$  \displaylines{
  E F = F E \, ,  \,\;  H_+ E = q^2 E H_+ + (q+1) E \, ,  \,\;
F H_+ = q^2 H_+ F + (q+1) F \, ,  \,\; H_+ H_- = H_- H_+  \cr
  E H_- = q^2 H_- E + (q+1) E \, ,  \; H_- F = q^2 F H_- +
(q+1) F \, ,  \; H_+ + H_- + (q-1) H_+ H_- = 0  \cr }  $$  
and it is a Hopf subalgebra, with Hopf operations given by
 $$  \displaylines{
  \Delta(E) = E \otimes 1 + 1 \otimes E + (q-1) \cdot H_+
\otimes E \, ,  \hskip19pt  \epsilon(E) = 0 \, ,  \hskip19pt  S(E) = - E -
(q-1) H_- E \, \phantom{.}  \cr
  \Delta(H_\pm) = H_\pm \otimes 1 + 1 \otimes H_\pm + (q-1) \cdot H_\pm
\otimes H_\pm \, ,  \hskip17pt  \epsilon(H_\pm) = 0 \, ,  \hskip17pt
S(H_\pm) = H_\mp \,  \hskip25pt \phantom{.}  \cr
  \Delta(F) = F \otimes 1 + 1 \otimes F + (q-1) \cdot F \otimes H_-
\, ,  \hskip19pt  \epsilon(F) = 0 \, ,  \hskip19pt  S(F) = - F - (q-1) F
H_+ \, .  \cr }  $$  
  \indent   It is easy to check that  $ \, U_q^s(\gere_2) \, $  is a
QrUEA, whose semiclassical limit is $ \, U(\gere_2) \, $:  in fact,
mapping the generators  $ \, F \mod (q-1) $,  $ \, D_\pm \mod (q-1) $,
$ \, E \mod (q-1) \, $  respectively to $ \, f $, $ \pm h\big/2 $,
$ e \in U(\gere_2) $  gives an isomorphism  $ \, U_q^s(\gere_2) \Big/
(q-1) \, U_q^s(\gere_2) \,{\buildrel \cong \over \longrightarrow}\,
U(\gere_2) \, $  of co-Poisson Hopf algebras.  Similarly, $ \,
U_q^a(\gere_2) \, $  is a QrUEA too, with semiclassical limit
$ \, U(\gere_2) \, $  again: here a co-Poisson Hopf algebra
isomorphism  $ \, U_q^a(\gere_2) \Big/ (q\!-\!1) \, U_q^a(\gere_2)
\, \cong \, U(\gere_2) \, $  is given mapping  $ \, F \! \mod
(q\!-\!1) $,  $ \, H_\pm \! \mod (q\!-\!1) $,  $ \, E \! \mod
(q\!-\!1) \, $  respectively to  $ \, f $,  $ \pm h $,
$ e \in U(\gere_2) \, $.

\vskip7pt

  {\bf 8.3 Computation of  $ \, {U_q(\gere_2)}' \, $  and
specialization  $ \, {U_q(\gere_2)}' \,{\buildrel {q \rightarrow 1}
\over \llongrightarrow}\, F\big[{E_2}^{\!*}\big] \, $.} \, This
section is devoted to compute  $ {U_q^s(\gere_2)}' $  and
$ {U_q^a(\gere_2)}' $,  and their specialization at  $ \, q = 1
\, $:  everything goes on as in \S 7.3, so we can be more sketchy.
From definitions we have, for any  $ \, n \in \N \, $,  $ \, \Delta^n(E)
= \sum_{s=1}^n K^{\otimes (s-1)} \otimes E \otimes 1^{\otimes (n-s)} $,
\, so  $ \; \delta_n(E) = {(K-1)}^{\otimes (n-1)} \otimes E =
{(q-1)}^{n-1} \cdot H_+^{\otimes (n-1)} \otimes E \, $,  \; whence
$ \; \delta_n\big( (q-1) E \big) \in {(q-1)}^n \, U_q^a(\gere_2)
\setminus {(q-1)}^{n+1} \, U_q^a(\gere_2) \; $  thus  $ \,
(q-1) E \in {U_q^a(\gere_2)}' $,  whereas  $ \, E \notin
{U_q^a(\gere_2)}' $.  Similarly, we have  $ \, (q-1) F $,
$ (q-1) H_\pm \in {U_q^a(\gere_2)}' \setminus (q-1) \,
{U_q^a(\gere_2)}' \, $.  Therefore  $ {U_q^a(\gere_2)} $
contains the subalgebra  $ U' $  generated by
$ \, \dot{F} := (q-1) F $,  $ \, \dot{H}_\pm := (q-1)
H_\pm $,  $ \, \dot{E} := (q-1) E \, $.  On the other hand,
$ {U_q^a(\gere_2)}' $  is clearly the  $ R $--span  of
the set  $ \; \Big\{\, F^a H_+^b H_-^c E^d \,\Big\vert\, a,
b, c, d \in \N \,\Big\} \, $:  \; to be precise, the set
  $$  \Big\{\, F^a H_+^b K^{-[b/2]} E^d \,\Big\vert\, a, b,
d \in \N \,\Big\} = \Big\{\, F^a H_+^b {\big( 1 + (q-1) H_-
\big)}^{[b/2]} E^d \,\Big\vert\, a, b, d \in \N \,\Big\}  $$   
is an  $ R $--basis  of  $ {U_q^a(\gere_2)}' $;  therefore,
a straightforward computation shows that any element in
$ {U_q^a(\gere_2)}' $  does necessarily lie in  $ U' \, $,  thus
$ {U_q^a(\gere_2)}' $  coincides with  $ U ' \, $.  Moreover,
since  $ \, \dot{H}_\pm = K^{\pm 1} - 1 \, $,  the unital
algebra  $ {U_q^a(\gere_2)}' $  is generated by  $ \dot{F} $,
$ K^{\pm 1} $  and  $ \dot{E} $  as well.
                                       \par
  The previous analysis   --- {\it mutatis mutandis} ---   ensures
also that $ \, {U_q^s(\gere_2)}' \, $  coincides with the unital
$ R $--subalgebra  $ U'' $  of  $ \Bbb{U}_q^s(\gere_2) $  generated by
$ \, \dot{F} := (q-1) F $,  $ \, \dot{D}_\pm := (q-1) D_\pm $,
$ \dot{E} := (q-1) E \, $;  in particular,  $ \, {U_q^s(\gere_2)}'
\supset {U_q^a(\gere_2)}' \, $.  Moreover, as  $ \, \dot{D}_\pm =
L^{\pm 1} - 1 \, $,  the unital algebra  $ {U_q^s(\gere_2)}' $  is
generated by  $ \dot{F} $,  $ L^{\pm 1} $  and  $ \dot{E} $  as
well.  Thus  $ {U_q^s(\gere_2)}' $  is the unital associative
$ R $--algebra  with generators  $ \, {\Cal F} := L \dot{F} $,
$ {\Cal L}^{\pm 1} := L^{\pm 1} $,  $ {\Cal E} := \dot{E}
L^{-1} \, $  and relations
  $$  \displaylines{
 {\Cal L} {\Cal L}^{-1} = 1 = {\Cal L}^{-1} {\Cal L} \, ,  \qquad
{\Cal E} {\Cal F} = {\Cal F} {\Cal E} \, ,  \qquad  {\Cal L}^{\pm 1}
{\Cal F} = q^{\mp 1} {\Cal F} {\Cal L}^{\pm 1} \, ,  \qquad
{\Cal L}^{\pm 1} {\Cal E} = q^{\pm 1} {\Cal E} {\Cal L}^{\pm 1}
\cr }  $$
with Hopf structure given by
  $$  \displaylines{
   \Delta({\Cal F}) = {\Cal F} \otimes {\Cal L}^{-1} + {\Cal L}
\otimes {\Cal F} \, ,  \hskip19pt  \Delta \big( {\Cal L}^{\pm 1} \big)
= {\Cal L}^{\pm 1} \otimes {\Cal L}^{\pm 1} \, ,  \hskip19pt
\Delta({\Cal E}) = {\Cal E} \otimes {\Cal L}^{-1} +
{\Cal L} \otimes {\Cal E}  \cr
   \epsilon({\Cal F}) = 0  \, ,  \hskip7pt  \epsilon \big(
{\Cal L}^{\pm 1} \big) = 1 \, ,  \hskip7pt  \epsilon({\Cal E}) = 0
\, ,  \hskip29pt  S({\Cal F}) = - {\Cal F} \, ,  \hskip7pt
S \big( {\Cal L}^{\pm 1} \big) = {\Cal L}^{\mp 1} \, ,  \hskip7pt
S({\Cal E}) = - {\Cal E} \, .  \cr }  $$
   As  $ \, q \rightarrow 1 \, $,  this yields a presentation of the
function algebra  $ \, F \big[ {}_s{E_2}^{\!*} \big] $,  and the
Poisson bracket that  $ \, F \big[ {}_s{E_2}^{\!*} \big] \, $
earns from this quantization process coincides with the one
coming from the Poisson structure on  $ \, {}_s{E_2}^{\!*} \, $:
namely, there is a Poisson Hopf algebra isomorphism
 $$  {U_q^s(\gere_2)}' \Big/ (q-1) \, {U_q^s(\gere_2)}'
\,{\buildrel \cong \over \llongrightarrow}\, F \big[
{}_s{E_2}^{\!*} \big]  $$
given by  $ \;\, {\Cal E} \mod (q-1) \mapsto x \, $,
$ \, {\Cal L}^{\pm 1} \mod (q-1) \mapsto z^{\pm 1} \, $,
$ {\Cal F} \mod (q-1) \mapsto y \, $.  That is, $ \,
{U_q^s(\gere_2)}' \, $  specializes to  $ \, F \big[
{}_s{E_2}^{\!*} \big] \, $  {\sl as a Poisson Hopf
algebra},  as predicted by Theorem 2.2.
                                            \par
  In the ``adjoint case'', from the definition of  $ U' $  and
from  $ \, {U_q^a(\gere_2)}' = U' \, $  we find that
$ {U_q^a(\gere_2)}' $  is the unital associative
$ R $--algebra  with generators  $ \, \dot{F} $,
$ K^{\pm 1} $,  $ \dot{E} \, $  and relations
  $$  K K^{-1} = 1 = K^{-1} K \, ,  \;\quad  \dot{E} \dot{F} =
\dot{F} \dot{E} \, ,  \;\quad  K^{\pm 1} \dot{F} = q^{\mp 2} \dot{F}
K^{\pm 1} \, ,  \;\quad  K^{\pm 1} \dot{E} = q^{\pm 2} \dot{E}
K^{\pm 1}  $$
with Hopf structure given by
  $$  \displaylines{
   \Delta\big(\dot{F}\big) = \dot{F} \otimes K^{-1} + 1 \otimes
\dot{F} \, ,  \hskip19pt  \Delta \big( K^{\pm 1} \big) = K^{\pm 1}
\otimes K^{\pm 1} \, ,  \hskip19pt  \Delta\big(\dot{E}\big) = \dot{E}
\otimes 1 + K \otimes \dot{E}  \cr
   \epsilon\big(\dot{F}\big) = 0  \, ,  \hskip3pt
\epsilon \big( K^{\pm 1} \big) = 1 \, ,  \hskip3pt
\epsilon\big(\dot{E}\big) = 0 \, ,  \hskip11pt
S\big(\dot{F}\big) = - \dot{F} K \, ,  \hskip3pt
S \big( K^{\pm 1} \big) = K^{\mp 1} ,  \hskip3pt
S\big(\dot{E}\big) = - K^{-1} \dot{E} \, .  \cr }  $$
   \indent   The conclusion is that a Poisson Hopf algebra isomorphism
  $$  {U_q^a(\gere_2)}' \Big/ (q-1) \, {U_q^a(\gere_2)}' \,{\buildrel
\cong \over \llongrightarrow}\, F \big[ {}_a{E_2}^{\!*} \big]  \quad
\Big( \subset F \big[ {}_s{E_2}^{\!*} \big] \Big) $$
exists, given by  $ \;\, \dot{E} \mod (q-1) \mapsto x z \, $,
$ \, K^{\pm 1} \mod (q-1) \mapsto z^{\pm 2} \, $,  $ \dot{F} \mod
(q-1) \mapsto z^{-1} y \, $,  i.e.~$ \, {U_q^a(\gere_2)}' \, $
specializes to  $ \, F \big[ {}_a{E_2}^{\!*} \big] \, $  {\sl as
a Poisson Hopf algebra},  according to Theorem 2.2.
                                            \par
   To finish with, note that  {\sl all this analysis (and its
outcome) is entirely characteristic-free}.

\vskip7pt

  {\bf 8.4 The identity  $ \, {\big({U_q(\gere_2)}'\big)}^{\!\vee}
= U_q(\gere_2) \, $.} \,  The goal of this section is to check that
part of Theorem 2.2{\it (b)}  claiming that  $ \; H \in \QrUEA
\,\Longrightarrow\, {\big(H'\big)}^{\!\vee} = H \; $  both for
$ \, H = U_q^s(\gere_2) \, $  and  $ \, H = U_q^a(\gere_2) \, $.
{\sl In addition, the proof below will work for  $ \, \Char(\Bbbk)
= 0 \, $  and  $ \, \Char(\Bbbk) > 0 \, $  too, thus giving a
stronger result than predicted by Theorem 2.2{\it (b)}}.
                                            \par
   First,  $ \, {U_q^s(\gere_2)}' $  is clearly a free
$ R $--module,  with basis  $ \, \Big\{\, {\Cal F}^a
{\Cal L}^d {\Cal E}^c \,\Big\vert\, a, c \in \N, d \in \Z
\,\Big\} \, $,  hence the set  $ \, {\Bbb B} := \Big\{\, {\Cal F}^a
{({\Cal L}^{\pm 1} - 1)}^b {\Cal E}^c \,\Big\vert\, a, b, c \in \N
\,\Big\} \, $,  is an  $ R $--basis  as well.  Second, as  $ \,
\epsilon({\Cal F}) = \epsilon\big({\Cal L}^{\pm 1} - 1\big) =
\epsilon({\Cal E}) = 0 \, $,  the ideal  $ \, J := \hbox{\sl Ker}
\Big( \epsilon \, \colon \, {U_q^s(\gere_2)}' \longrightarrow R \!
\Big) \, $  is the span of  $ \, {\Bbb B} \setminus \{1\} $.   Now
$ \; I := \hbox{\sl Ker} \Big( {U_q^s(\gere_2)}' \, {\buildrel \epsilon
\over {\relbar\joinrel\relbar\joinrel\twoheadrightarrow}} \, R \,
{\buildrel {q \mapsto 1} \over\llongtwoheadrightarrow} \, \Bbbk \Big)
= J + (q-1) \cdot {U_q^s(\gere_2)}' \, $,  \; therefore  $ \;
{\big({U_q^s(\gere_2)}'\big)}^{\!\vee} := \sum_{n \geq 0} {\Big(
{(q-1)}^{-1} I \Big)}^n \; $  is generated   --- as a unital
$ R $--subalgebra  of  $ \Bbb{U}_q^s(\gere_2) $  ---   by  $ \,
{(q-1)}^{-1} {\Cal F} = L F $,  $ {(q-1)}^{-1} ({\Cal L} - 1) =
D_+ $,  $ {(q-1)}^{-1} \big( {\Cal L}^{-1} - 1 \big) = D_- $,
$ {(q-1)}^{-1} {\Cal E} = E L^{-1} $,  hence by  $ F $,  $ D_\pm $,
$ E $,  so it coincides with  $ U_q^s(\gere_2) $,  q.e.d.
                                                 \par
  The situation is entirely similar for the adjoint case: one simply
has to change  $ {\Cal F} $,  $ {\Cal L}^{\pm 1} $,  $ {\Cal E} $
respectively with  $ \dot{F} $,  $ K^{\pm 1} $,  $ \dot{E} $,  and
$ D_\pm $  with  $ H_\pm $,  then everything goes through as above.

\vskip7pt

  {\bf 8.5 The quantum hyperalgebra  $ \hyp_q(\gere_2) $.} \, Like
for semisimple groups, we can define ``quantum hyperalgebras''
attached to  $ \gere_2 $  mimicking what done in \S 7.5.  Namely,
we can first define a Hopf subalgebra of  $ \Bbb{U}_q^s(\gere_2) $  over
$ \Z \big[ q, q^{-1} \big] $  whose specialization at  $ \, q = 1
\, $  is exactly the Kostant-like  $ \Z $--integer  form  $ U_\Z
(\gere_2) $  of  $ U(\gere_2) $  (generated by divided powers, and
giving the hyperalgebra  $ \hyp(\gere_2) $  over any field  $ \Bbbk $
by scalar extension, namely  $ \, \hyp(\gere_2) = \Bbbk \otimes_\Z
U_\Z(\gere_2) \, $),  \, and then take its scalar extension over
$ R \, $.
                                             \par
   To be precise, let  $ \hyp^{s,\Z}_q(\gere_2) $  be the unital
$ \Z\big[q,q^{-1}\big] $--subalgebra  of  $ \Bbb{U}^s_q(\gere_2) $
(defined like above  {\sl but over\/}  $ \Z\big[q,q^{-1}\big] $)
generated by the ``quantum divided powers''  $ \, \displaystyle{
F^{(n)} \! := F^n \! \Big/ {[n]}_q! } \; $,  $ \, \displaystyle{
\left( {{L \, ; \, c} \atop {n}} \right) \! := \prod_{r=1}^n {{\;
q^{c+1-r} L - 1 \,} \over {\; q^r - 1 \;}} } \, $,  $ \, \displaystyle{
E^{(n)} \! := E^n \! \Big/ {[n]}_q! } \; $  (for all  $ \, n \in \N \, $
and  $ \, c \in \Z \, $,  with notation of \S 7.5)  and by  $ L^{-1} $.
Comparing with the case of  $ \gersl_2 $  one easily sees that this is
a Hopf subalgebra of  $ \Bbb{U}_q^s(\gere_2) $,  and  $ \, \hyp^{s,\Z}_q
(\gere_2){\Big|}_{q=1} \cong\, U_\Z(\gere_2) \, $;  \, thus  $ \,
\hyp^s_q(\gere_2) := R \otimes_{\Z[q,q^{-1}]} \hyp^{s,\Z}_q(\gere_2)
\, $  (for any  $ R $  like in \S 8.2, with  $ \, \Bbbk := R \big/
\h \, R \, $  and  $ \, p := \Char(\Bbbk) \, $)  specializes at  $ \,
q = 1 \, $  to the  $ \Bbbk $--hyperalgebra  $ \hyp(\gere_2) $.  In
addition, among all the  $ \left( {{L \, ; \, c} \atop {n}} \right) $'s
it is enough to take only those with  $ \, c = 0 \, $.  {\sl From now
on we assume  $ \, p > 0 \, $.}
                                             \par
   Again a strict comparison with the  $ \gersl_2 $  case   --- with
some shortcuts, since the defining relations of  $ \hyp^s_q(\gere_2) $
are simpler! ---   shows us that  $ \, {\hyp^s_q(\gere_2)}' \, $  is
the unital  $ R $--subalgebra  of  $ \hyp^s_q(\gere_2) $  generated by
$ L^{-1} $  and the ``rescaled quantum divided powers''  $ \, {(q \!
- \! 1)}^n F^{(n)} \, $,  $ \, {(q\!-\!1)}^n \! \left( {{L \, ; \, 0}
\atop {n}} \right) \, $  and  $ \, {(q\!-\!1)}^n E^{(n)} \, $  for all
$ \, n \in \N \, $.  It follows that  $ \, {\hyp^s_q(\gere_2)}'
{\Big|}_{q=1} \, $  is generated by the corresponding specializations
of  $ \, {(q-1)}^{p^r} F^{(p^r)} \, $,  $ \, {(q-1)}^{p^r} \! \left(
{{L \, ; \, 0} \atop {p^r}} \right) \, $  and  $ \, {(q-1)}^{p^r}
E^{(p^r)} \, $  for all  $ \, r \in \N \, $:  \, this
proves that the spectrum of  $ \, {\hyp^s_q(\gere_2)}'{\Big|}_{q=1}
\, $  has dimension 0 and height 1, and its cotangent Lie algebra
has basis  $ \, \Big\{\, {(q\!-\!1)}^{p^r} F^{(p^r)}, \,
{(q\!-\!1)}^{p^r} \! \left( {{L \, ; \, 0} \atop {p^r}} \right) \! ,
 \allowbreak
  \, {(q\!-\!1)}^{p^r} E^{(p^r)} \, \mod (q\!-\!1) \, {\hyp^s_q(\gerg)}'
\, \mod J^{\,2} \;\Big|\; r \!\in\! \N \,\Big\} \, $  (where  $ J $  is
the augmentation ideal of  $ {\hyp^s_q(\gere_2)}'{\Big|}_{q=1} \, $,  so
that  $ \, J \Big/ J^{\,2} \, $  is the aforementioned cotangent Lie
bialgebra).  Moreover,  $ \, \big( {\hyp^s_q(\gere_2)}' \big)^{\!\vee}
\, $  is generated by  $ \, {(q-1)}^{p^r-1} F^{(p^r)} \, $,  $ \,
{(q-1)}^{p^r-1} \left( {{L \, ; \, 0} \atop {p^r}} \right) \, $,
$ L^{-1} $  and  $ \, {(q-1)}^{p^r-1} E^{(p^r)} \, $  (for all  $ \,
r \in \N \, $):  \, in particular  $ \, \big( {\hyp^s_q(\gere_2)}'
\big)^{\!\vee} \subsetneqq \hyp^s_q(\gere_2) \, $,  \, and finally
$ \, \big( {\hyp^s_q(\gere_2)}' \big)^{\!\vee}{\Big|}_{q=1} \, $  is
generated by the cosets modulo  $ (q-1) $  of the elements above, which
in fact form a basis of the restricted Lie bialgebra  $ \gerk $  such
that  $ \, \big( {\hyp^s_q(\gere_2)}' \big)^{\!\vee} {\Big|}_{q=1} =
\, \u(\gerk) \, $.
                                             \par
   All this analysis was made starting from  $ \Bbb{U}_q^s(\gere_2) $,
which gave ``simply connected quantum objects''.  If we start instead
from  $ \Bbb{U}_q^a(\gere_2) $,  we get ``adjoint quantum objects'' following
the same pattern but for replacing everywhere  $ L^{\pm 1} $  by
$ K^{\pm 1} \, $:  apart from these changes, the analysis and its
outcome will be exactly the same.  Like for  $ \, \gersl_2 $  (cf.~\S
7.5), all the adjoint quantum objects   ---  i.e.~$ \hyp^a_q(\gere_2) $,
$ {\hyp^a_q(\gere_2)}' $  and  $ \big( {\hyp^a_q(\gere_2)}' \big)^{\!
\vee} $  ---   will be strictly contained in the corresponding simply
connected quantum objects; nevertheless, the semiclassical limits will
be the same in the case of  $ \hyp_q(\gere_2) $  (always yielding
$ \hyp(\gere_2) \, $)  and in the case of  $ \big( {\hyp_q(\gere_2)}'
\big)^{\!\vee} $  (giving  $ \u(\gerk) $,  in both cases), while the
semiclassical limit of  $ {\hyp_q(\gere_2)}' $  in the simply connected
case will be a (countable) covering of that in the adjoint case.

\vskip7pt

  {\bf 8.6 The QFAs  $ \, F_q[E_2] \, $  and  $ \, F_q[{}_aE_2] \, $.}
\, In this and the following sections we look at Theorem 2.2 starting
from QFAs, to get QrUEAs out of them.
                                             \par
   We begin by introducing a QFA for the Euclidean groups  $ E_2 $
and  $ {}_aE_2 \, $.  Let  $ \, F_q[E_2] \, $  be the unital associative
\hbox{$R$--alge}bra  with generators  $ \, \text{a}^{\pm 1} $,
$ \text{b} $,  $ \text{c} \, $  and relations
  $$  \text{a} \, \text{b} = q \, \text{b} \, \text{a} \, ,  \qquad\qquad
\text{a} \, \text{c} = q \, \text{c} \, \text{a} \, ,  \qquad\qquad
\text{b} \, \text{c} = \text{c} \, \text{b}  $$
endowed with the Hopf algebra structure given by
  $$  \displaylines{
   \Delta\big(\text{a}^{\pm 1}\big) = \text{a}^{\pm 1} \otimes
\text{a}^{\pm 1} \, ,  \hskip19pt    \Delta(\text{b}) = \text{b}
\otimes \text{a}^{-1} + \text{a} \otimes \text{b} \, ,  \hskip19pt
\Delta(\text{c}) = \text{c} \otimes \text{a} + \text{a}^{-1}
\otimes \text{c}  \cr
   \epsilon\big(\text{a}^{\pm 1}\big) = 1 \, ,  \hskip5pt
\epsilon(\text{b}) = 0 \, ,  \hskip5pt  \epsilon(\text{c}) = 0 \, ,
\hskip21pt  S\big(\text{a}^{\pm 1}\big) = \text{a}^{\mp 1} \, ,
\hskip5pt  S(\text{b}) = - q^{-1} \, \text{b} \, ,  \hskip5pt
S(\text{c}) = - q^{+1} \, \text{c} \, .  \cr }  $$
   \indent   Define  $ \, F_q[{}_aE_2] \, $  as the
$ R $--submodule  of  $ F_q[E_2] $  spanned by the
products of an even number of generators, i.e.~monomials
of even degree in  $ \text{a}^{\pm 1} $,  $ \text{b} $,
$ \text{c} \, $:  this is a unital subalgebra of  $ F_q[E_2] $,
generated by  $ \, \beta := \text{b} \, \text{a} $,  $ \alpha^{\pm 1}
:= \text{a}^{\pm 2} $,  and  $ \gamma := \text{a}^{-1} \text{c} \, $.
Let also  $ \, \F_q[E_2] := {\big(F_q[E_2]\big)}_F \, $  and  $ \,
\F_q[{}_aE_2] := {\big(F_q[{}_aE_2]\big)}_F \, $,  \, which have the
same presentation than  $ F_q[E_2] $  and  $ F_q[{}_aE_2] $  but over
$ F(R) $.  Essentially by definition, both $ \, F_q[E_2] \, $  and
$ \, F_q[{}_aE_2] \, $  are QFAs (at  $ \, \h = q - 1 \, $),  whose
semiclassical limit is  $ F[E_2] $  and  $ F[{}_aE_2] $  respectively.

\vskip7pt

  {\bf 8.7 Computation of  $ {F_q[E_2]}^\vee $  and
$ {F_q[{}_aE_2]}^\vee $  and specializations  $ {F_q[E_2]}^\vee
\!\!{\buildrel {q \rightarrow 1} \over \llongrightarrow}\,
U(\gerg^*) $  and  $ \, {F_q[{}_aE_2]}^\vee \!\!{\buildrel
{q \rightarrow 1} \over \llongrightarrow}\, U(\gerg^*) \, $.} \,
In this section we go and compute  $ \, {F_q[G]}^\vee \, $  and
its semiclassical limit (i.e.~its specialization at  $ \, q = 1
\, $)  for both  $ \, G = E_2 \, $  and  $ \, G = {}_aE_2 \, $.
                                             \par
   First,  $ F_q[E_2] $  is free over  $ R $,  with basis  $ \,
\Big\{\, \text{b}^b \text{a}^a \text{c}^c \,\Big\vert\, a \in \Z,
b, c \in \N \,\Big\} \, $,  so the set  $ \, {\Bbb B}_s := \Big\{\,
\text{b}^b {(\text{a}^{\pm 1} - 1)}^a \text{c}^c \,\Big\vert\, a,
b, c \in \N \,\Big\} \, $  is an  $ R $--basis  as well.  Second,
since  $ \, \epsilon(\text{b}) =   \hbox{$ \epsilon \big(
\text{a}^{\pm 1} - 1 \big) $}   = \epsilon(\text{c}) = 0 \, $,
the ideal  $ \, J := \text{\sl Ker} \, \Big( \epsilon \, \colon
\, F_q[E_2] \loongrightarrow R \Big) \, $  is the span of  $ \,
{\Bbb B}_s \setminus \{1\} \, $.  Now  $ \; I := \text{\sl Ker} \,
\Big( F_q[E_2] \, {\buildrel \epsilon \over \llongtwoheadrightarrow}
\, R \, {\buildrel {q \mapsto 1} \over \llongtwoheadrightarrow}
\, \Bbbk \Big) = J + (q-1) \cdot F_q[E_2] \, $,  \; thus  $ \;
{F_q[E_2]}^\vee := \sum_{n \geq 0} {\Big( \! {(q\!-\!1)}^{-1}
I \Big)}^n \, $  turns out to be the unital  $ R $--algebra
(subalgebra of  $ \F_q[E_2] $)  with generators  $ \, D_\pm :=
\displaystyle{\, \text{a}^{\pm 1} - 1 \, \over \, q - 1\,} \, $,
$ \, E := \displaystyle{\, \text{b} \, \over \, q - 1\,} \, $,
and  $ \, F := \displaystyle{\, \text{c} \, \over \, q - 1\,} \, $
and relations
  $$  \displaylines{
   D_+ E = q E D_+ + E \, ,  \quad  D_+ F = q F D_+ + F \, , \quad
E D_- = q D_- E + E \, ,  \quad  F D_- = q D_- F + F  \cr
   E F = F E \, ,  \;\qquad  D_+ D_- = D_- D_+ \, ,  \;\qquad
D_+ + D_- + (q-1) D_+ D_- = 0  \cr }  $$
with a Hopf structure given by
  $$  \matrix
   \Delta(E) = E \otimes 1 + 1 \otimes E + (q-1) \big( E \otimes D_-
+ D_+ \otimes E \big) \, ,  &  \;  \epsilon(E) = 0 \, ,  &  \;  S(E) =
- q^{-1} E \, \phantom{.}  \\
   \Delta(D_\pm) = D_\pm \otimes 1 + 1 \otimes D_\pm +
(q-1) \cdot D_\pm \otimes D_\pm \, ,  &  \;  \epsilon(D_\pm) = 0 \, ,
&  \;  S(D_\pm) = D_\mp \, \phantom{.}  \\
   \Delta(F) = F \otimes 1 + 1 \otimes F + (q-1) \big( F \otimes D_+ +
D_- \otimes F \big) \, ,  &  \;  \epsilon(F) = 0 \, ,  &  \;
S(F) = - q^{+1} F \, .  \\
     \endmatrix  $$
This implies that  $ \, {F_q[E_2]}^\vee \,{\buildrel \, q
\rightarrow 1 \, \over \llongrightarrow}\, U({\gere_2}^{\! *}) \, $  as
co-Poisson Hopf algebras, for a co-Poisson Hopf algebra isomorphism
  $$  {F_q[E_2]}^\vee \Big/ (q-1) \, {F_q[E_2]}^\vee
\,{\buildrel \cong \over \llongrightarrow}\, U({\gere_2}^{\! *})  $$
exists, given by  $ \;\, D_\pm \mod (q-1) \mapsto \pm \text{h} \big/ 2
\, $,  $ \, E \mod (q-1) \mapsto \text{e} \, $,  $ \, F \mod (q-1)
\mapsto \text{f} \, $;  thus  $ \, {F_q[E_2]}^\vee \, $  does
specialize to  $ \, U({\gere_2}^{\! *}) \, $ {\sl as a co-Poisson
Hopf algebra},  q.e.d.
                                             \par
  Similarly, if we consider  $ F_q[{}_aE_2] $  the same
analysis works again.  In fact,  $ F_q[{}_aE_2] $  is free
over $ R $,  with basis  $ \, {\Bbb B}_a :=  \Big\{\, \beta^b
{(\alpha^{\pm 1} - 1)}^a \gamma^c \,\Big\vert\, a, b, c \in \N
\,\Big\} \, $;  therefore, as above the ideal  $ \, J := \text{\sl
Ker} \, \Big( \epsilon \, \colon \, F_q[{}_aE_2] \rightarrow R \Big)
\, $  is the span of  $ \, {\Bbb B}_a \setminus \{1\} \, $.  Now, we
have  $ \; I := \text{\sl Ker} \, \Big( F_q[{}_aE_2] \, {\buildrel
\epsilon \over \llongtwoheadrightarrow} \, R \, {\buildrel {q
\mapsto 1} \over \llongtwoheadrightarrow} \, \Bbbk \Big) = J +
(q-1) \cdot F_q[{}_aE_2] \, $,  \; so  $ \; {F_q[{}_aE_2]}^\vee
:= \sum_{n \geq 0} \! {\Big( \! {(q\!-\!1)}^{-1} I \Big)}^n $
\; is nothing but the unital  $ R $--algebra  (subalgebra of
$ \F_q[{}_aE_2] \, $)  with generators  $ \, H_\pm :=
\displaystyle{\, \alpha^{\pm 1} - 1 \, \over \, q - 1\,}
\, $,  $ \, E' := \displaystyle{\, \beta \, \over \, q - 1\,}
\, $,  and  $ \, F' := \displaystyle{\, \gamma \, \over \,
q - 1\,} \, $  and relations
  $$  \displaylines{
   E' F' \! = \! q^{-2} F' E' \! ,  \, H_+ E' \! = \! q^2 E' H_+
+ (q+1) E' \! ,  \, H_+ F' \! = \! q^2 F' H_+ + (q+1) F' \! ,  \,
H_+ H_- \! = \! H_- H_+  \cr
   E' H_- = q^2 H_- E' + (q+1) E' \! ,  \, F' H_- = q^2 H_- F'
+ (q+1) F' \! ,  \, H_+ + H_- + (q-1) H_+ H_- = 0  \cr }  $$
with a Hopf structure given by
  $$  \displaylines{
   \Delta(E') = E' \otimes 1 + 1 \otimes E' + (q-1) \cdot H_+
\otimes E' \, ,  \hskip19pt  \epsilon(E') = 0 \, ,  \hskip19pt
S(E') = - E' - (q-1) H_- E' \, \phantom {.}  \cr
   \;\;  \Delta(H_\pm) = H_\pm \otimes 1 + 1 \otimes H_\pm +
(q-1) \cdot H_\pm \otimes H_\pm \, ,  \hskip17pt
\epsilon(H_\pm) = 0 \, ,  \hskip17pt S(H_\pm) = H_\mp \,
\hskip25pt  \phantom{.}  \cr
   \Delta(F') = F' \otimes 1 + 1 \otimes F' + (q-1) \cdot H_- \otimes F'
\, ,  \hskip19pt  \epsilon(F') = 0 \, ,  \hskip19pt  S(F') = - F' - (q-1)
H_+ F' \, .  \cr }  $$
This implies that  $ \, {F_q[{}_aE_2]}^\vee \,{\buildrel \, q
\rightarrow 1 \, \over \llongrightarrow}\, U({\gere_2}^{\! *}) \, $  as
co-Poisson Hopf algebras, for a co-Poisson Hopf algebra isomorphism
  $$  {F_q[{}_aE_2]}^\vee \Big/ (q-1) \, {F_q[{}_aE_2]}^\vee
\,{\buildrel \cong \over \llongrightarrow}\, U({\gere_2}^{\! *})  $$
is given by  $ \;\, H_\pm \mod (q-1) \mapsto \pm \text{h} \, $,  $ \, E'
\mod (q-1) \mapsto \text{e} \, $,  $ \, F' \mod (q-1) \mapsto \text{f}
\, $;  so  $ \, {F_q[{}_aE_2]}^\vee \, $  too specializes to
$ \, U({\gere_2}^{\! *}) \, $  {\sl as a co-Poisson Hopf algebra},
as expected.
                                            \par
   We finish noting that, once more,  {\sl this analysis (and its
outcome) is characteristic-free}.

\vskip7pt

  {\bf 8.8 The identities  $ \, {\big({F_q[E_2]}^\vee\big)}' =
F_q[E_2] \, $  and  $ \, {\big({F_q[{}_aE_2]}^\vee\big)}' =
F_q[{}_aE_2] \, $.} \, In this section we verify   for the QFAs
$ \, H = F_q[E_2] \, $  and  $ \, H = F_q[{}_aE_2] \, $  the
validity of the part of Theorem 2.2{\it (b)}  claiming that
$ \; H \in \QFA \,\Longrightarrow\, {\big(H^\vee\big)}' = H \, $.
Once more, our arguments will prove this result for  $ \, \Char(\Bbbk)
\geq 0 \, $,  \, thus going beyond what forecasted by Theorem 2.2.
                                           \par
   By induction we find formulas  $ \; \Delta^n(E) = \sum_{r+s+1=n}
\text{a}^{\otimes r} \otimes E \otimes {\big( \text{a}^{-1}
\big)}^{\otimes s} $,  $ \, \Delta^n(D_\pm) = \sum_{r+s+1=n}
{\big( \text{a}^{\pm 1} \big)}^{\otimes r} \otimes D_\pm \otimes
1^{\otimes s} $,  and  $ \, \Delta^n(F) = \sum_{r+s+1=n}
{\big( \text{a}^{-1} \big)}^{\otimes r}\otimes E \otimes
\text{a}^{\otimes s} $:  these imply
  $$  \eqalign{
   \delta_n(E)  &  = \! {\textstyle \sum\limits_{r+s+1=n}} \!
{(\text{a} - 1)}^{\otimes r} \otimes E \otimes {\big( \text{a}^{-1}
- 1 \big)}^{\otimes s} = {(q-1)}^{n-1} \! {\textstyle
\sum\limits_{r+s+1=n}} \! {D_+}^{\! \otimes r}
\otimes E \otimes {D_-}^{\! \otimes s}  \cr
  \delta_n  &  (D_\pm) = {\big( \text{a}^{\pm 1} - 1 \big)}^{\otimes
(n-1)} \otimes D_\pm = {(q-1)}^{n-1} {D_\pm}^{\!\otimes n}  \cr
  \delta_n(F)  &  = \! {\textstyle \sum\limits_{r+s+1=n}} \!
{\big( \text{a}^{-1} - 1 \big)}^{\otimes r} \otimes E \otimes
{(\text{a} - 1)}^{\otimes s} = {(q-1)}^{n-1} \! {\textstyle
\sum\limits_{r+s+1=n}} \! {D_-}^{\! \otimes r} \otimes
E \otimes {D_+}^{\! \otimes s}  \cr }  $$
which gives  $ \; \dot{E} := (q-1) E $,  $ \, \dot{D}_\pm
:= (q-1) D_\pm $,  $ \, \dot{F} := (q-1) F \in {\big( {F_q[E_2]}^\vee
\big)}' \setminus (q-1) \cdot $
$ \cdot {\big({F_q[E_2]}^\vee\big)}' \, $.  \; So  $ {\big(
{F_q[E_2]}^\vee \big)}' $  contains the unital  $ R $--subalgebra
$ A' $  generated (inside  $ \, \F_q[E_2] \, $)  by  $ \dot{E}$,
$ \dot{D}_\pm  $  and  $ \dot{F} \, $;  but  $ \, \dot{E} =
\text{b} $,  $ \, \dot{D}_\pm = \text{a}^{\pm 1} - 1 $,  and  $ \,
\dot{F} = \text{c} $,  thus  $ A' $  is just  $ F_q[E_2] $.  Since
$ {F_q[E_2]}^\vee $  is the  $ R $--span  of  $ \, \Big\{\, E^e
D_+^{d_+} D_-^{d_-} F^f \,\Big\vert\, e, d_+, d_-, f \in \N \,\Big\}
\, $,  one easily sees   --- using the previous formulas for
$ \Delta^n $  ---   that in fact  $ \, {\big({F_q[E_2]}^\vee\big)}' =
A' = F_q[E_2] \, $,  \, q.e.d.
                                            \par
   When dealing with the adjoint case, the previous arguments
go through again: in fact,  $ {\big({F_q[{}_aE_2]}^\vee\big)}' $
turns out to coincide with the unital  $ R $--subalgebra  $ A'' $
generated (inside  $ \, \F_q[{}_aE_2] \, $)  by  $ \, \dot{E}' :=
(q-1) E' = \beta \, $,  $ \, \dot{H}_\pm := (q-1) H_\pm = \alpha^{\pm 1}
- 1 \, $,  and  $ \, \dot{F}' := (q-1) F' = \gamma \, $;  but this is also
generated by  $ \beta $,  $ \alpha^{\pm 1} $  and  $ \gamma $,  thus it
coincides with  $ F_q[{}_aE_2] $,  q.e.d.
%
%
 \eject   

\centerline {\bf \S \; 9 \  Fourth example: quantum Heisenberg group }

\vskip10pt

  {\bf 9.1  The classical setting.} \, Let  $ \Bbbk $  be any field
of characteristic  $ \, p \geq 0 \, $.  Let  $ \, G := H_n(\Bbbk) =
H_n \, $,  the  $ (2 \, n + 1) $--dimensional  Heisenberg group; its
tangent Lie algebra  $ \, \gerg = \gerh_n \, $  is generated by
$ \, \{\, f_i, h, e_i \,\vert\, i = 1, \dots, n \,\} \, $  with
relations  $ \, [e_i,f_j] = \delta_{i{}j} h $,  $ [e_i,e_j] =
[f_i,f_j] = [h,e_i] = [h,f_j] = 0 \, $  ($ \forall \, i, j = 1,
\dots n \, $).  The formulas  $ \, \delta(f_i) = h \otimes f_i -
f_i \otimes h \, $,  $ \, \delta(h) = 0 \, $,  $ \, \delta(e_i) =
h \otimes e_i - e_i \otimes h \, $  ($ \forall \, i = 1, \dots n \, $)
make  $ \gerh_n $  into a Lie bialgebra, which yields  $ H_n $
with a structure of Poisson group; these same formulas give also a
presentation of the co-Poisson Hopf algebra  $ U({\gerh}_n) $  (with
the standard Hopf structure).  When  $ \, p > 0 \, $  we consider on
$ \gerh_n $  the  $ p $--operation  uniquely defined by  $ \, e_i^{\,
[p\,]} = 0 \, $,  $ \, f_i^{\,[p\,]} = 0 \, $,  $ \, h^{[p\,]} = h
\, $  (for all  $ \, i = 1, \dots, n \, $),  which makes it into a
restricted Lie bialgebra.  The group  $ H_n $  is usually realized
as the group of all square matrices  $ \, {\big( a_{i{}j} \big)}_{i,j
= 1, \dots, n+2;} \, $  such that  $ \, a_{i{}i} = 1 \; \forall \, i
\, $  and  $ \, a_{i{}j} = 0 \; \forall\, i, j \, $  such that either
$ \, i > j \, $  or  $ \, 1 \not= i < j \, $  or  $ \, i < j \not= n+2
\, $;  it can also be realized as  $ \, H_n = \Bbbk^n \times \Bbbk
\times \Bbbk^n \, $  with group operation given by  $ \; \big(
\underline{a}', c', \underline{b}' \big) \cdot \big( \underline{a}'',
c'', \underline{b}'' \big) = \big( \underline{a}' + \underline{a}'',
c' + c'' + \underline{a}' \ast \underline{b}'', \underline{b}' +
\underline{b}'' \big) \, $,  \;  where we use vector notation  $ \,
\underline{v} = (v_1, \dots, v_n) \in k^n \, $  and  $ \, \underline{a}'
\ast \underline{b}'' := \sum_{i=1}^n a'_i b''_i \, $  is the standard
scalar product in  $ \, k^n \, $;  in particular the identity of
$ H_n $  is  $ \, e = (\underline{0}, 0, \underline{0}) \, $  and the
inverse of a generic element is given by  $ \, {\big( \underline{a},
c, \underline{b} \big)}^{-1} = \big( \! - \! \underline{a} \, , - c
+ \underline{a} \ast \underline{b} \, , \! - \underline{b} \, \big)
\, $.  Therefore, the function algebra  $ F[{H}_n] $  is the unital
associative commutative  $ \Bbbk $--algebra  with generators  $ \, a_1 $,
$ \dots $,  $ a_n $,  $ c $,  $ b_1 $,  $ \dots $,  $ b_n $,  and
with Poisson Hopf algebra structure given by
  $$  \displaylines{
   \Delta(a_i) = a_i \otimes 1 + 1 \otimes a_i \, ,  \hskip5pt
\Delta(c) = c \otimes 1 + 1 \otimes c + {\textstyle \sum_{\ell=1}^{n}}
a_\ell \otimes b_\ell \, ,  \hskip5pt  \Delta(b_i) =
b_i \otimes 1 + 1 \otimes b_i  \cr
   \epsilon(a_i) = 0 \, ,  \hskip5pt  \epsilon(c) = 0 \, ,  \hskip5pt
\epsilon(b_i) = 0 \, ,  \hskip15pt  S(a_i) = - a_i \, ,  \hskip5pt
S(c) = - c + {\textstyle \sum_{\ell=1}^{n}} a_\ell b_\ell \, ,
\hskip5pt  S(b_i) = - b_i  \cr
   \{a_i,a_j\} = 0 \, ,  \hskip15pt  \{a_i,b_j\} = 0 \, ,  \hskip15pt
\{b_i,b_j\} = 0 \, , \hskip15pt  \{c \, , a_i\} = a_i \, ,  \hskip15pt
\{c \, , b_i\} = b_i  \cr }  $$
for all  $ \; i, j= 1, \dots, n \, $.  (N.B.: with respect to this
presentation, we have  $ \, f_i = {\partial_{b_i}}{\big\vert}_e \, $,
$ \, h = {\partial_c}{\big\vert}_e \, $,  $ \, e_i = {\partial_{a_i}}
{\big\vert}_e \, $,  where  $ e $  is the identity element of  $ H_n
\, $).  The dual Lie bialgebra  $ \, \gerg^* = {\gerh_n}^{\!*} \, $
is the Lie algebra with generators  $ \, \text{f}_i $,  $ \text{h} $,
$ \text{e}_i \, $,  and relations  $ \, [\text{h},\text{e}_i] =
\text{e}_i $,  $ [\text{h},\text{f}_i] = \text{f}_i $,  $ [\text{e}_i,
\text{e}_j] = [\text{e}_i,\text{f}_j] = [\text{f}_i,\text{f}_j] = 0
\, $,  with Lie cobracket given by  $ \, \delta(\text{f}_i) = 0 $,
$ \, \delta(\text{h}) = \sum_{j=1}^{n} (\text{e}_j \otimes \text{f}_j
- \text{f}_j \otimes \text{e}_j) $,  $ \, \delta(\text{e}_i) = 0 \, $
for all  $ \, i= 1, \dots, n \, $  (we take  $ \, \text{f}_i := f_i^*
\, $,  $ \, \text{h} := h^* \, $,  $ \, \text{e}_i := e_i^* \, $,
where  $ \, \big\{\, f_i^*, h^*, e_i^* \,\vert\, i= 1, \dots, n
\,\big\} \, $  is the basis of  $ {\gerh_n}^{\! *} $  which is
the dual of the basis  $ \{\, f_i, h, e_i \,\vert\, i= 1, \dots,
n \,\} $  of  $ \gerh_n \, $).  This again gives a presentation
of  $ \, U({\gerh_n}^{\!*}) \, $  too.  If  $ \, p > 0 \, $  then
$ {\gerh_n}^{\!*} $  is a restricted Lie bialgebra with respect to
the  $ p $--operation  given by  $ \, \text{e}_i^{\,[p\,]} = 0 \, $,
$ \, \text{f}_i^{\,\,[p\,]} = 0 \, $,  $ \, \text{h}^{[p\,]} = \text{h}
\, $  (for all  $ \, i = 1, \dots, n \, $).  The simply connected
algebraic Poisson group with tangent Lie bialgebra  $ {{\gerh}_n}^{\!*} $
can be realized (with  $ \, \Bbbk^\star := \Bbbk \setminus \{0\} \, $)
as  $ \; {{}_s{H}_n}^{\! *} = \Bbbk^n \times \Bbbk^\star \times \Bbbk^n
\, $,  with group operation  $ \; \big( \underline{\dot{\alpha}},
\underline{\dot{\gamma}}, \underline{\dot{\beta}} \,\big) \cdot
\big( \underline{\check{\alpha}}, \underline{\check{\gamma}},
\underline{\check{\beta}} \,\big) = \big( \check{\gamma}
\underline{\dot{\alpha}} + {\dot{\gamma}}^{-1}
\underline{\check{\alpha}}, \dot{\gamma} \check{\gamma},
\check{\gamma} \underline{\dot{\beta}} + {\dot{\gamma}}^{-1}
\underline{\check{\beta}} \,\big) \, $;  so the identity of
$ {{}_s{H}_n}^{\! *} $  is  $ \, e = (\underline{0}, 1, \underline{0})
\, $  and the inverse is given by  $ \, {\big( \underline{\alpha},
\gamma, \underline{\beta} \,\big)}^{-1} = \big( \! - \!
\underline{\alpha}, \gamma^{-1}, - \underline{\beta} \,\big) $.
Its centre is  $ \, Z \big( {{}_s{H}_n}^{\! *} \big) = \big\{
(\underline{0}, 1, \underline{0}), (\underline{0}, -1, \underline{0})
\big\} =: Z \, $,  so there is only one other (Poisson) group with
tangent Lie bialgebra  $ {{\gerh}_n}^{\! *} \, $,  that is the
adjoint group  $ \, {{}_a{H}_n}^{\! *} := {{}_s{H}_n}^{\! *}
\Big/ Z \, $.
                                                  \par
   It is clear that  $ F\big[{{}_s{H}_n}^{\! *}] $  is the unital
associative commutative  $ \Bbbk $--algebra  with generators
$ \, \alpha_1 $,  $ \dots $,  $ \alpha_n $,  $ \gamma^{\pm 1} $,
$ \beta_1 $,  $ \dots $,  $ \beta_n $,  and with Poisson Hopf
algebra structure given by
  $$  \displaylines{
   \Delta(\alpha_i) = \alpha_i \otimes \gamma + \gamma^{-1}
\otimes \alpha_i \, ,  \hskip17pt  \Delta\big(\gamma^{\pm 1}\big) =
\gamma^{\pm 1} \otimes \gamma^{\pm 1},  \hskip17pt  \Delta(\beta_i) =
\beta_i \otimes \gamma + \gamma^{-1} \otimes \beta_i  \cr
   \epsilon(\alpha_i) = 0 \, ,  \hskip5pt  \epsilon \big(
\gamma^{\pm 1} \big) = 1 \, ,  \hskip5pt  \epsilon(\beta_i) = 0 \, ,
\hskip29pt  S(\alpha_i) =  - \alpha_i \, ,  \hskip5pt  S \big(
\gamma^{\pm 1} \big) = \gamma^{\mp 1} ,  \hskip5pt  S(\beta_i)
= - \beta_i  \cr
   \{\alpha_i,\alpha_j\} = \{\alpha_i,\beta_j\} = \{\beta_i,
\beta_j\} = \{\alpha_i,\gamma\} = \{\beta_i,\gamma\} = 0 \, ,
\hskip13pt  \{\alpha_i,\beta_j\} = \delta_{i{}j} \big(
\gamma^2 - \gamma^{-2} \big) \big/ 2  \cr }  $$
for all  $ \; i, j= 1, \dots, n \, $  (N.B.: with respect to this
presentation, we have  $ \, \text{f}_i = {\partial_{\beta_i}}
{\big\vert}_e \, $,  $ \, \text{h} = {\,1\, \over \,2\,} \, \gamma \,
{\partial_\gamma}{\big\vert}_e \, $,  $ \, \text{e}_i =
{\partial_{\alpha_i}}{\big\vert}_e \, $,  where  $ e $
is the identity element of  $ {{}_s{H}_n}^{\! *} \, $),  and
$ F\big[{{}_a{H}_n}^{\! *}\big] $  can be identified   --- as in the
case of the Euclidean group ---   with the Poisson Hopf subalgebra of
$ F\big[{{}_a{H}_n}^{\! *}\big] $  which is spanned by products of an
even number of generators: this is generated, as a unital  subalgebra,
by  $ \alpha_i \gamma $,  $ \gamma^{\pm 2} $,  and  $ \gamma^{-1}
\beta_i \; $  ($ \, i= 1, \dots, n \, $).

\vskip7pt

  {\bf 9.2 The QrUEAs  $ \, U_q^s({\gerh}_n) \, $  and  $ \, U_q^a
({\gerh}_n) \, $.} \, We switch now to quantizations.  Once again,
let  $ R $  be a domain and let  $ \, \h \in R \setminus \{0\} \, $
and  $ \, q := 1 + \h \in R \, $  be like in \S 7.2.
                                                  \par
   Let $ \, \Bbb{U}_q(\gerg) = \Bbb{U}_q^s ({\gerh}_n) \, $  be the unital
associative  $ F(R) $--algebra  with generators  $ \, F_i $,
$ L^{\!\pm 1} $,  $ E_i $  ($ i = 1, \dots, n \, $)  and relations
  $$  L L^{-1} = 1 = L^{-1} L \, ,  \,\quad  L^{\!\pm 1} F = F L^{\pm 1}
\, ,  \,\quad  L^{\pm 1} E = E L^{\pm 1} \, ,  \,\quad  E_i F_j - F_j E_i
= \delta_{i{}j} {{\, L^2 - L^{-2} \, \over \, q - q^{-1} \,}}   $$
for all  $ \, i, j= 1, \dots, n \, $;  we give it a structure of Hopf
algebra, by setting  ($ \forall \; i, j = 1, \dots, n \, $)
  $$  \displaylines{
   \Delta(E_i) = E_i \otimes 1 + L^2 \otimes E_i \, ,  \hskip11pt
\Delta\big(L^{\pm 1}\big) = L^{\pm 1} \otimes L^{\pm 1} ,  \hskip11pt
\Delta(F_i) = F_i \otimes L^{-2} + 1 \otimes F_i  \cr
   \epsilon(E_i) = 0 \, ,  \hskip3pt  \epsilon\big(L^{\pm 1}\big) = 1 \, ,
\hskip3pt  \epsilon(F_i) = 0 \, ,  \hskip13pt  S(E_i) = - L^{-2} E_i ,
\hskip3pt  S\big(L^{\pm 1}\big) =  L^{\mp 1} ,  \hskip3pt
S(F_i) = - F_i L^2  \cr }  $$
Note that  $ \, \Big\{\, \prod_{i=1}^n F_i^{a_i} \! \cdot \! L^z \!
\cdot \! \prod_{i=1}^n E_i^{d_i} \,\Big\vert\, z \in \Z, \, a_i,
d_i \in \N, \, \forall\, i \,\Big\} \, $  is an  $ F(R) $--basis
of  $ \Bbb{U}_q^s(\gerh_n) $.
                                             \par
   Now, let  $ \, U_q^s(\gerh_n) \, $  be the unital
$ R $--subalgebra  of  $ \Bbb{U}_q^s(\gerh_n) $  generated by
                  $ \, F_1 $,  $ \dots $,  $ F_n $,\break
$ D := \displaystyle{\, L - 1 \, \over \, q - 1 \,} $,  $ \varGamma :=
\displaystyle{{\, L - L^{-2} \, \over \, q - q^{-1} \,}} $,  $ E_1 $,
$ \dots $,  $ E_n \, $.  Then  $ U_q^s(\gerh_n) $  can be presented as
the associative unital algebra with generators  $ \, F_1 $,   $ \dots $,
$ F_n $,  $ L^{\pm 1} $,  $ D $,  $ \varGamma $,  $ E_1 $,  $ \dots $,
$ E_n $ and relations
  $$  \displaylines{
   \hskip7pt  D X = X D \, ,  \hskip25pt  L^{\pm 1} X = X L^{\pm 1} \, ,
\hskip25pt  \varGamma X = X \varGamma \, ,  \hskip25pt  E_i F_j - F_j E_i =
\delta_{i{}j} \varGamma  \cr
  L = 1 + (q-1) D \, ,  \hskip15pt  L^2 - L^{-2} = \big( q -
q^{-1} \big) \varGamma \, ,  \hskip15pt  D (L + 1) \big(1 + L^{-2}\big)
= \big(1 + q^{-1}\big) \varGamma  \cr }  $$
for all  $ \; X \in {\big\{F_i, L^{\pm 1}, D, \varGamma, E_i
\big\}}_{i=1,\dots,n} \, $  and  $ \, i, j= 1, \dots, n \, $;
furthermore,  $ U_q^s(\gerh_n) $  is a Hopf subalgebra (over
$ R $),  with
  $$  \matrix
      \Delta(\varGamma) = \varGamma \otimes L^2 + L^{-2} \otimes \varGamma
\, , &  \hskip25pt  \epsilon(\varGamma) = 0  \, ,  &  \hskip25pt
S(\varGamma) = - \varGamma  \\
      \Delta(D) = D \otimes 1 + L \otimes D \, ,  &  \hskip25pt
\epsilon(D) = 0 \, ,  &  \hskip25pt  S(D) = - L^{-1} D \, .  \\
     \endmatrix  $$
Moreover, from relations  $ \, L = 1 + (q-1) D \, $  and
$ \, L^{-1} = L^3 - \big( q - q^{-1} \big) L \varGamma \, $  it follows
that
  $$  U_q^s(\gerh_n) \; = \; \hbox{$ R $--span  of}  \;\,
\bigg\{\, {\textstyle \prod\limits_{i=1}^n} F_i^{a_i} \! \cdot
\! D^b \varGamma^c \! \cdot \! {\textstyle \prod\limits_{i=1}^n}
E_i^{d_i} \,\bigg\vert\, a_i, b, c, d_i \in \N, \, \forall\,
i= 1, \dots, n \,\bigg\}   \eqno (9.1)  $$
                                                  \par
   The ``adjoint version'' of  $ \Bbb{U}_q^s(\gerh_n) $  is the unital
subalgebra  $ \Bbb{U}_q^a(\gerh_n) $  generated  by  $ \, F_i $,
$ K^{\pm 1} := L^{\pm 2} $,  $ E_i \; (i = 1,
\dots, n) $,  which is clearly a Hopf subalgebra.  It also has an
$ R $--integer  form  $ U_q^a(\gerh_n) $,  namely the unital
$ R $--subalgebra  generated by  $ \, F_1 \, $,  $ \dots $,  $ F_n
\, $, $ K^{\pm 1} $,  $ H := \displaystyle{\, K - 1 \, \over \, q - 1
\,} \, $,  $ \varGamma := \displaystyle{\, K - K^{-1} \, \over \, q -
q^{-1} \,} \, $,  $ \, E_1 $,  $ \dots $,  $ E_n  \, $:  this
has relations
  $$  \displaylines{
   \hskip7pt  H X = X H \, ,  \hskip25pt  K^{\pm 1} X = X K^{\pm 1}
\, ,  \hskip25pt  \varGamma X = X \varGamma \, ,  \hskip25pt
E_i F_j - F_j E_i = \delta_{i{}j} \varGamma  \cr
  K = 1 + (q-1) H \, ,  \hskip15pt  K - K^{-1} = \big( q -
q^{-1} \big) \varGamma \, ,  \hskip15pt  H \big(1 + K^{-1}\big) =
\big( 1 + q^{-1} \big) \varGamma  \cr }  $$  
for all  $ \; X \in {\big\{ F_i, K^{\pm 1}, H, \varGamma, E_i
\big\}}_{i=1,\dots,n} \, $  and  $ \, i, j= 1, \dots, n \, $,
and Hopf operations given by
  $$  \matrix
   \Delta(E_i) = E_i \otimes 1 + K \otimes E_i \, ,  &  \hskip21pt
\epsilon(E_i) = 0 \, ,  &  \hskip21pt  S(E_i) = - K^{-1} E_i  \\
   \Delta\big(K^{\pm 1}\big) = K^{\pm 1} \otimes K^{\pm 1} \, ,  &
\hskip21pt  \epsilon\big(K^{\pm 1}\big) = 1 \, ,  &  \hskip21pt
S\big(K^{\pm 1}\big) = K^{\mp 1}  \\
   \Delta(H) = H \otimes 1 + K \otimes H \, ,  &  \hskip21pt
\epsilon(H) = 0 \, ,  &  \hskip21pt  S(H) = - K^{-1} H  \\
   \Delta(\varGamma) = \varGamma \otimes K^{-1} + K \otimes \varGamma
\, ,   &  \hskip21pt  \epsilon(\varGamma) = 0 \, ,  &  \hskip21pt
S(\varGamma) = - \varGamma  \\
   \Delta(F_i) = F_i \otimes K^{-1} + 1 \otimes F_i \, ,  &
\hskip21pt  \epsilon(F_i) = 0 \, ,  &  \hskip21pt  S(F_i) =
- F_i K^{+1}  \\
      \endmatrix  $$   
for all  $ \, i= 1, \dots, n $.  One can easily check that  $ \, U_q^s
(\gerh_n) \, $  is a QrUEA, with  $ \, U(\gerh_n) \, $  as semiclassical
limit: in fact, mapping the generators  $ \, F_i \mod (q-1) $,  $ \,
L^{\pm 1} \mod (q-1) $,  $ \, D \mod (q-1) \, $,  $ \, \varGamma \mod
(q-1) \, $,  $ \, E_i \mod (q-1) \, $  respectively to  $ \, f_i $,
$ 1 $,  $ h\big/2 $,  $ h $,  $ e_i \in U(\gerh_n) \, $  yields a
co-Poisson Hopf algebra isomorphism between  $ \, U_q^s(\gerh_n)
\Big/ (q-1) \, U_q^s(\gerh_n) \, $  and  $ U(\gerh_n) $.  Similarly,
$ \,U_q^a(\gerh_n) \, $  is a QrUEA too, again with limit  $ \, U(\gerh_n)
\, $,  for a co-Poisson Hopf algebra isomorphism between  $ \,
U_q^a(\gerh_n) \Big/ (q-1) \, U_q^a(\gerh_n) \, $  and  $ \,
U(\gerh_n) \, $  is given by mapping the generators  $ \, F_i
\mod (q-1) $,  $ \, K^{\pm 1} \mod  (q-1) $,  $ \, H \mod (q-1) $,
$ \, \varGamma \mod (q-1) $,  $ E_i \mod (q-1) \, $  respectively
to  $ \, f_i $,  $ 1 $,  $ h $,  $ h $,  $ e_i \in U(\gerh_n) $.

\vskip7pt

  {\bf 9.3 Computation of  $ \, {U_q(\gerh_n)}' \, $  and
specialization  $ \, {U_q(\gerh_n)}' \,{\buildrel {q \rightarrow 1}
\over \llongrightarrow}\, F\big[{H_n}^{\!*}\big] \, $.} \, Here we
compute  $ {U_q^s(\gerh_n)}' $  and  $ {U_q^a(\gerh_n)}' $,  and
their semiclassical limits, along the pattern of \S 7.3.
                                               \par
   Definitions give, for any  $ \, n \in \N \, $,  $ \, \Delta^n(E_i)
= \sum_{s=1}^n {(L^2)}^{\otimes (s-1)} \otimes E_i \otimes 1^{\otimes
(n-s)} $,  \, hence  $ \; \delta_n(E_i) = {(q-1)}^{n-1} \cdot D^{\otimes
(n-1)} \otimes E_i \; $  so  $ \; \delta_n\big( (q-1) E \big) \in
{(q-1)}^n \, U_q^s(\gerh_n) \setminus {(q-1)}^{n+1} \, U_q^s(\gerh_n)
\; $  whence  $ \, \dot{E}_i := (q-1) \, E_i \in {U_q^s(\gerh_n)}' $,
whereas  $ \, E_i \notin {U_q^s(\gerh_n)}' $;  similarly, we have  $ \;
\dot{F}_i := (q-1) \, F_i $,  $ L^{\pm 1} $,  $ \dot{D} := (q-1) \, D
= L - 1 $,  $ \dot{\varGamma} := (q-1) \, \varGamma \in {U_q^s(\gerh_n)}'
\setminus (q-1) \, {U_q^s(\gerh_n)}' $,  for all  $ \, i = 1, \dots, n
\, $.  Therefore  $ {U_q^s(\gerh_n)}' $  contains the subalgebra  $ U' $
generated by  $ \, \dot{F}_i $,  $ L^{\pm 1} $,  $ \, \dot{D} $,
$ \dot{\varGamma} $,  $ \, \dot{E}_i \, $;  we conclude that in
fact  $ \, {U_q^s(\gerh_n)}' = U' \, $:  this is easily seen
--- like for  $ {SL}_2 $  and for  $ E_2 $  ---   using
the formulas above along with (9.1).  As a consequence,
$ {U_q^s(\gerh_n)}' $  is the unital  $ R $--algebra  with
generators  $ \, \dot{F}_1 $,  $ \dots $,  $ \dot{F}_n $,
$ L^{\pm 1} $,  $ \dot{D} $,  $ \dot{\varGamma} $,
$ \dot{E}_1 $,  $ \dots $,  $ \dot{E}_n \, $  and relations
  $$  \displaylines{
   \hskip7pt  \dot{D} \dot{X} = \dot{X} \dot{D} \, ,  \hskip25pt
L^{\pm 1} \dot{X} = \dot{X} L^{\pm 1} \, ,  \hskip25pt
\dot{\varGamma} \dot{X} = \dot{X} \dot{\varGamma} \, ,  \hskip25pt
\dot{E}_i \dot{F}_j - \dot{F}_j \dot{E}_i = \delta_{i{}j} (q-1)
\dot{\varGamma}  \cr
  L = 1 + \dot{D} \, ,  \hskip15pt  L^2 - L^{-2} = \big( 1 + q^{-1}
\big) \dot{\varGamma} \, ,  \hskip15pt  \dot{D} (L + 1) \big( 1 +
L^{-2}\big) = \big( 1 + q^{-1} \big) \dot{\varGamma}  \cr }  $$
for all  $ \; \dot{X} \in {\big\{ \dot{F}_i, L^{\pm 1}, \dot{D},
\dot{\varGamma}, \dot{E}_i \big\}}_{i=1,\dots,n} \, $  and  $ \, i,
j= 1, \dots, n \, $,  with Hopf structure given by
  $$  \matrix
   \Delta\big(\dot{E}_i\big) = \dot{E}_i \otimes 1 + L^2 \otimes \dot{E}_i
\, ,  &  \hskip7pt  \epsilon\big(\dot{E}_i\big) = 0 \, ,  &  \hskip7pt
S\big(\dot{E}_i\big) =  - L^{-2} \dot{E}_i  &  \hskip19pt  \forall
\; i = 1, \dots, n \, \phantom{.}  \\
   \Delta\big(L^{\pm 1}\big) = L^{\pm 1} \otimes L^{\pm 1} \, ,  &
\hskip7pt  \epsilon\big(L^{\pm 1}\big) = 1 \, ,  &  \hskip7pt
S\big(L^{\pm 1}\big) =  L^{\mp 1}  &  {}  \\
   \Delta\big(\dot{\varGamma}\big) = \dot{\varGamma} \otimes L^2 + L^{-2}
\otimes \dot{\varGamma} \, , &  \hskip7pt
\epsilon\big(\dot{\varGamma}\big) = 0  \, ,  &  \hskip7pt
S\big(\dot{\varGamma}\big) = - \varGamma  \\
   \Delta\big(\dot{D}\big) = \dot{D} \otimes 1 + L \otimes \dot{D} \, ,
&  \hskip7pt  \epsilon\big(\dot{D}\big) = 0 \, ,  &  \hskip7pt
S\big(\dot{D}\big) = - L^{-1} \dot{D}  \\
   \Delta\big(\dot{F}_i\big) = \dot{F}_i \otimes L^{-2} + 1 \otimes
\dot{F}_i \, ,  &  \hskip7pt  \epsilon\big(F_i\big) = 0 \, ,  &
\hskip7pt  S\big(\dot{F}_i\big) = - \dot{F}_i L^2  &  \hskip19pt
\forall \; i = 1, \dots, n \, .  \\
   \endmatrix  $$
   \indent   A similar analysis shows that  $ \, {U_q^a(\gerh_n)}' \, $
coincides with the unital  $ R $--subalgebra  $ U'' $  of
$ U_q^a(\gerh_n) $  generated by  $ \, \dot{F}_i $,  $ \, K^{\pm 1} $,
$ \, \dot{H} := (q-1) \, H $,  $ \, \dot{\varGamma} $,  $ \dot{E}_i \, $
($ i = 1, \dots, n $);  in particular,  $ \, {U_q^a(\gerh_n)}' \subset
{U_q^s(\gerh_n)}' \, $.  Therefore  $ {U_q^a(\gerh_n)}' $  can be
presented as the unital associative  $ R $--algebra  with
generators  $ \, \dot{F}_1 $,  $ \dots $,  $ \dot{F}_n $,
$ \dot{H} $,  $ K^{\pm 1} $,  $ \dot{\varGamma} $,  $ \dot{E}_1 $,
$ \dots $,  $ \dot{E}_n \, $  and relations
  $$  \displaylines{
   \hskip7pt  \dot{H} \dot{X} = \dot{X} \dot{H} \, ,  \hskip25pt
K^{\pm 1} \dot{X} = \dot{X} K^{\pm 1} \, ,  \hskip25pt  \dot{\varGamma}
\dot{X} = \dot{X} \dot{\varGamma} \, ,  \hskip25pt  \dot{E}_i \dot{F}_j
- \dot{F}_j \dot{E}_i = \delta_{i{}j} (q-1) \dot{\varGamma}  \cr
  K = 1 + \dot{H} \, ,  \hskip15pt  K - K^{-1} = \big( 1 + q^{-1} \big)
\dot{\varGamma} \, ,  \hskip15pt  \dot{H} \big( 1 + K^{-1} \big) =
\big( 1 + q^{-1} \big) \dot{\varGamma}  \cr }  $$
for all  $ \; \dot{X} \in {\big\{ \dot{F}_i, K^{\pm 1}, \dot{K},
\dot{\varGamma}, \dot{E}_i \big\}}_{i=1,\dots,n} \, $  and
$ \, i, j= 1, \dots, n \, $,  with Hopf structure given by
  $$  \matrix
   \Delta\big(\dot{E}_i\big) = \dot{E}_i \otimes 1 + K \otimes \dot{E}_i
\, ,  &  \hskip7pt  \epsilon\big(\dot{E}_i\big) = 0 \, ,  &  \hskip7pt
S\big(\dot{E}_i\big) =  - K^{-1} \dot{E}_i  &  \hskip19pt  \forall
\; i = 1, \dots, n \, \phantom{.}  \\
   \Delta\big(K^{\pm 1}\big) = K^{\pm 1} \otimes K^{\pm 1} \, ,  &
\hskip7pt  \epsilon\big(K^{\pm 1}\big) = 1 \, ,  &  \hskip7pt
S\big(K^{\pm 1}\big) = K^{\mp 1}  &  {}  \\
   \Delta\big(\dot{\varGamma}\big) = \dot{\varGamma} \otimes K +
K^{-1} \otimes \dot{\varGamma} \, , &  \hskip7pt
\epsilon\big(\dot{\varGamma}\big) = 0  \, ,  &  \hskip7pt
S\big(\dot{\varGamma}\big) = - \varGamma  \\
   \Delta\big(\dot{H}\big) = \dot{H} \otimes 1 + K \otimes \dot{H} \, ,
&  \hskip7pt  \epsilon\big(\dot{H}\big) = 0 \, ,  &  \hskip7pt
S\big(\dot{H}\big) = - K^{-1} \dot{H}  \\
   \Delta\big(\dot{F}_i\big) = \dot{F}_i \otimes K^{-1} + 1 \otimes
\dot{F}_i \, ,  &  \hskip7pt  \epsilon\big(F_i\big) = 0 \, ,  &
\hskip7pt S\big(\dot{F}_i\big) = - \dot{F}_i K  &  \hskip19pt
\forall \; i = 1, \dots, n \, .  \\
   \endmatrix  $$
   \indent   As  $ \, q \rightarrow 1 \, $,  the presentation above
yields an isomorphism of Poisson Hopf algebras
 \vskip-5pt
  $$  {U_q^s(\gerh_n)}' \Big/ (q-1) \, {U_q^s(\gerh_n)}'
\,{\buildrel \cong \over \llongrightarrow}\,
F \big[ {}_s{H_n}^{\!*} \big]  $$
 \vskip-2pt
\noindent   
given by  $ \;\, \dot{E}_i \mod (q-1) \mapsto \alpha_i \gamma^{+1}
\, $,  $ \, L^{\pm 1} \mod (q-1) \mapsto \gamma^{\pm 1} \, $,  $ \,
\dot{D} \mod (q-1) \mapsto \gamma - 1 \, $,  $ \, \dot{\varGamma}
\mod (q-1) \mapsto \big( \gamma^2 - \gamma^{-2} \big) \Big/ 2 \, $,
$ \dot{F}_i \mod (q-1) \mapsto \gamma^{-1} \beta_i \, $.  In other
words, the semiclassical limit of  $ \, {U_q^s(\gerh_n)}' \, $  is
$ \, F \big[ {}_s{H_n}^{\!*} \big] \, $,  as predicted by  Theorem
2.2{\it (c)\/}  for  $ \, p = 0 \, $.  Similarly, when considering
the ``adjoint case'', we find a Poisson Hopf algebra isomorphism
  $$  {U_q^a(\gerh_n)}' \Big/ (q-1) \, {U_q^a(\gerh_n)}'
\,{\buildrel \cong \over \llongrightarrow}\, F \big[ {}_a{H_n}^{\!*}
\big]  \quad  \Big( \subset F \big[ {}_s{H_n}^{\!*} \big] \Big)  $$
given by  $ \;\, \dot{E}_i \mod (q-1) \mapsto \alpha_i \gamma^{+1} \, $,
$ \, K^{\pm 1} \mod (q-1) \mapsto \gamma^{\pm 2} \, $,  $ \, \dot{H} \mod
(q-1) \mapsto \gamma^2 - 1 \, $,  $ \, \dot{\varGamma} \mod (q-1) \mapsto
\big( \gamma^2 - \gamma^{-2} \big) \Big/ 2 \, $,  $ \dot{F}_i \mod (q-1)
\mapsto \gamma^{-1} \beta_i \, $.  That is to say,  $ \, {U_q^a(\gerh_n)}'
\, $  has semiclassical limit  $ \, F \big[ {}_a{H_n}^{\!*} \big] \, $,
\, as predicted by  Theorem 2.2{\it (c)\/}  for  $ \, p = 0 \, $.
                                         \par
   We stress the fact that  {\sl this analysis is
characteristic-free},  so we get in fact that its outcome
does hold for  $ \, p > 0 \, $  as well, thus ``improving''
Theorem 2.2{\it (c)\/}  (like in \S\S 7--8).

\vskip7pt

  {\bf 9.4 The identity  $ \, {\big({U_q(\gerh_n)}'\big)}^{\!\vee} =
U_q(\gerh_n) \, $.} \,  In this section we verify the part of  Theorem
2.2{\it (b)}  claiming, for  $ \, p = 0 \, $,  that  $ \, H \in \QrUEA
\,\Longrightarrow\, {\big(H'\big)}^{\!\vee} = H \, $,  both for  $ \,
H = U_q^s(\gerh_n) \, $  and for  $ \, H = U_q^a(\gerh_n) \, $.  {\sl
In addition, the same arguments will prove such a result for  $ \,
p > 0 \, $  too.}
                                          \par
   To begin with, using (9.1) and the fact that  $ \, \dot{F}_i $,
$ \dot{D} $,  $ \dot{\varGamma} $,  $ \dot{E}_i \in \text{\sl
Ker} \Big( \epsilon \, \colon \, {U_q^s(\gerh_n)}' \,
\relbar\joinrel\twoheadrightarrow \, R \Big) \, $
we get that  $ \, J := \text{\sl Ker}\,(\epsilon) \, $
is the  $ R $--span  of  $ \, {\Bbb M} \setminus \{1\} \, $,
where  $ \, {\Bbb M} \, $  is the set in the right-hand-side of
(9.1).  Since  $ \, {\big({U_q^s(\gerh_n)}'\big)}^{\!\vee} :=
\sum_{n \geq 0} \! {\Big( \! {(q\!-\!1)}^{-1} I \Big)}^n \, $
with  $ \, I := \text{\sl Ker} \Big( {U_q^s(\gerh_n)}' \,
{\buildrel \epsilon \over {\relbar\joinrel\twoheadrightarrow}}
\, R \, {\buildrel {q \mapsto 1} \over \llongtwoheadrightarrow}
\, \Bbbk \Big) = J + (q-1) \cdot {U_q^s(\gerh_n)}' \; $  we
have that  $ \, {\big( {U_q^s (\gerh_n)}' \big)}^{\!\vee}
\, $  is generated   --- as a unital  $ R $--subalgebra  of
$ \Bbb{U}_q^s(\gerh_n) $  ---   by  $ \, {(q-1)}^{-1} \dot{F}_i
= F_i \, $,  $ {(q-1)}^{-1} \dot{D} = D \, $,  $ {(q-1)}^{-1}
\dot{\varGamma} = \varGamma \, $,  $ {(q-1)}^{-1} \dot{E}_i
= E_i \, $  ($ i = 1, \dots, n $),  so it coincides with
$ U_q^s(\gerh_n) $,  q.e.d.  In the adjoint case the
procedure is similar: one changes  $ L^{\pm 1} $,
resp.~$ \dot{D} $,  with  $ K^{\pm 1} $,  resp.~$ \dot{H} $,
and everything works as before.

\vskip7pt

  {\bf 9.5 The quantum hyperalgebra  $ \hyp_q(\gerh_n) $.} \, Like in
\S\S 7.5 and 8.5, we can define ``quantum hyperalgebras'' associated
to  $ \gerh_n \, $.  Namely, first we define a Hopf subalgebra of
$ \Bbb{U}_q^s(\gerh_n) $  over  $ \Z \big[ q, q^{-1} \big] $  whose
specialization at  $ \, q = 1 \, $  is the natural Kostant-like
$ \Z $--integer  form  $ U_\Z(\gerh_n) $  of  $ U(\gerh_n) $
(generated by divided powers, and giving the hyperalgebra
$ \hyp(\gerh_n) $  over any field  $ \Bbbk $  by scalar
  extension),
%
%
   and then take its scalar extension
over  $ R \, $.
                                             \par
   To be precise, let  $ \hyp^{s,\Z}_q(\gerh_n) $  be the unital
$ \Z\big[q,q^{-1}\big] $--subalgebra  of  $ \Bbb{U}^s_q(\gerh_n) $
(defined like above  {\sl but over\/}  $ \Z\big[q,q^{-1}\big] $)
generated by the ``quantum divided powers''  $ \, \displaystyle{
F_i^{(m)} \! := F_i^m \! \Big/ {[m]}_q! } \; $,  $ \, \displaystyle{
\left( {{L \, ; \, c} \atop {m}} \right) \! := \prod_{r=1}^n {{\;
q^{c+1-r} L - 1 \,} \over {\; q^r - 1 \;}} } \, $,  $ \, \displaystyle{
E_i^{(m)} \! := E_i^m \! \Big/ {[m]}_q! } \; $  (for all  $ \, m \in
\N \, $,  $ \, c \in \Z \, $  and  $ \, i = 1, \dots, n \, $,  with
notation of \S 7.5)  and by  $ L^{-1} $.  Comparing with the case of
$ \gersl_2 $   --- noting that for each  $ i $  the quadruple  $ (F_i,
L, L^{-1}, E_i) $  generates a copy of  $ \Bbb{U}^s_q(\gersl_2) $  ---   we
see at once that this is a Hopf subalgebra of  $ \Bbb{U}_q^s(\gerh_n) $,
and  $ \, \hyp^{s,\Z}_q(\gerh_n){\Big|}_{q=1} \cong\, U_\Z(\gerh_n)
\, $;  \, thus  $ \, \hyp^s_q(\gerh_n) := R \otimes_{\Z[q,q^{-1}]}
\hyp^{s,\Z}_q(\gerh_n) \, $  (for any  $ R $  like in \S 8.2, with
$ \, \Bbbk := R \big/ \h \, R \, $  and  $ \, p := \Char(\Bbbk) \, $)
specializes at  $ \, q = 1 \, $  to the  $ \Bbbk $--hyperalgebra
$ \hyp(\gerh_n) $.  Moreover, among all the  $ \left( {{L \, ; \, c}
\atop {n}} \right) $'s  it is enough to take only those with  $ \,
c = 0 \, $.  {\sl From now on we assume  $ \, p > 0 \, $.}
                                             \par
   Pushing forward the close comparison with the case of  $ \gersl_2 $
we also see that  $ \, {\hyp^s_q(\gerh_n)}' \, $  is the unital
$ R $--subalgebra  of  $ \hyp^s_q(\gerh_n) $  generated by  $ L^{-1} $
and the ``rescaled quantum divided powers''  $ \, {(q \! - \! 1)}^m
F_i^{(m)} \, $,  $ \, {(q\!-\!1)}^m \! \left( {{L \, ; \, 0} \atop {m}}
\right) \, $  and  $ \, {(q\!-\!1)}^m E_i^{(m)} \, $,  for all  $ \, m
\in \N \, $  and  $ \, i = 1, \dots, n \, $.  It follows that  $ \,
{\hyp^s_q(\gerh_n)}'{\Big|}_{q=1} \, $  is generated by the
specializations at  $ \, q = 1 \, $  of  $ \, {(q-1)}^{p^r}
F_i^{(p^r)} \, $,  $ \, {(q-1)}^{p^r} \! \left( {{L \, ; \, 0}
\atop {p^r}} \right) \, $  and  $ \, {(q-1)}^{p^r} E_i^{(p^r)} \, $,
for all  $ \, r \in \N \, $,  $ \, i = 1, \dots, n \, $:  \, this proves
directly that the spectrum of  $ \, {\hyp^s_q(\gerh_n)}'{\Big|}_{q=1}
\, $  has dimension 0 and height 1, and its cotangent Lie algebra has
basis  $ \, \Big\{\, {(q\!-\!1)}^{p^r} F_i^{(p^r)}, \, {(q\!-\!1)}^{p^r}
\! \left( {{L \, ; \, 0} \atop {p^r}} \right) \! ,  \, {(q\!-\!1)}^{p^r}
E_i^{(p^r)} \, \mod (q\!-\!1) \, {\hyp^s_q(\gerg)}' \, \mod J^{\,2}
\;\Big|\; r \!\in\! \N \, , \, i = 1, \dots, n \,\Big\} \, $  (with
$ J $  being the augmentation ideal of  $ {\hyp^s_q(\gerh_n)}'
{\Big|}_{q=1} \, $,  so that  $ \, J \Big/ J^{\,2} \, $  is the
aforementioned cotangent Lie bialgebra).  Finally,  $ \, \big(
{\hyp^s_q(\gerh_n)}' \big)^{\!\vee} \, $  is generated by  $ \,
{(q-1)}^{p^r-1} F_i^{(p^r)} \, $,  $ \, {(q-1)}^{p^r-1} \left(
{{L \, ; \, 0} \atop {p^r}} \right) \, $,  $ L^{-1} $  and  $ \,
{(q-1)}^{p^r-1} E_i^{(p^r)} \, $  (for  $ \, r \in \N \, $,  $ \, i =
1, \dots, n \, $):  \, in particular  $ \, \big( {\hyp^s_q(\gerh_n)}'
\big)^{\!\vee} \subsetneqq \hyp^s_q(\gerh_n) \, $,  \, and  $ \, \big(
{\hyp^s_q(\gerh_n)}' \big)^{\!\vee}{\Big|}_{q=1} \, $  is generated
by the cosets modulo  $ (q-1) $  of the elements above, which form
indeed a basis of the restricted Lie bialgebra  $ \gerk $  such
that  $ \, \big( {\hyp^s_q(\gerh_n)}' \big)^{\!\vee}{\Big|}_{q=1}
= \, \u(\gerk) \, $.
                                             \par
   The previous analysis stems from  $ \Bbb{U}_q^s(\gerh_n) $,  and so
gives ``simply connected quantum objects''.  Instead we can start
from  $ \Bbb{U}_q^a(\gerh_n) $,  thus getting ``adjoint quantum objects'',
moving along the same pattern but for replacing  $ L^{\pm 1} $  by
$ K^{\pm 1} $  throughout: apart from this, the analysis and its
outcome are exactly the same.  Like for  $ \, \gersl_2 $  (cf.~\S
7.5), all the adjoint quantum objects   ---  i.e.~$ \hyp^a_q(\gerh_n) $,
$ {\hyp^a_q(\gerh_n)}' $  and  $ \big( {\hyp^a_q(\gerh_n)}' \big)^{\!
\vee} $  ---   will be strictly contained in the corresponding simply
connected quantum objects; however, the semiclassical limits will be
the same in the case of  $ \hyp_q(\gerg) $  (giving  $ \hyp(\gerh_n) $,
in both cases) and in the case of  $ \big( {\hyp_q(\gerg)}' \big)^{\!
\vee} $  (always yielding  $ \u(\gerk) $), whereas the semiclassical
limit of  $ {\hyp_q(\gerg)}' $  in the simply connected case will
be a (countable) covering of the limit in the adjoint case.

\vskip7pt

  {\bf 9.6 The QFA  $ \, F_q[H_n] \, $.} \, Now we look at
Theorem 2.2 the other way round, i.e.~from QFAs to QrUEAs.
We begin by introducing a QFA for the Heisenberg group.
                                             \par
   Let  $ \, F_q[H_n] \, $  be the unital associative
\hbox{$R$--alge}bra  with generators  $ \, \text{a}_1 $,
$ \dots $,  $ \text{a}_n $,  $ \text{c} $,  $ \text{b}_1 $,
$ \dots $,  $ \text{b}_n $,  and relations  (for all  $ \, i,
j = 1, \dots, n \, $)
  $$  \text{a}_i \text{a}_j = \text{a}_j \text{a}_i \, ,  \hskip5pt
\text{a}_i \text{b}_j = \text{b}_j \text{a}_i \, ,  \hskip5pt  \text{b}_i
\text{b}_j = \text{b}_j \text{b}_i \, ,  \hskip5pt  \text{c} \, \text{a}_i
= \text{a}_i \, \text{c} + (q-1) \, \text{a}_i \, ,  \hskip5pt  \text{c}
\, \text{b}_j = \text{b}_j \, \text{c} + (q-1) \, \text{b}_j  $$
with a Hopf algebra structure given by (for all  $ \, i,
j = 1, \dots, n \, $)
  $$  \displaylines{
   \Delta(\text{a}_i) = \text{a}_i \otimes 1 + 1 \otimes \text{a}_i
\, ,  \hskip9pt  \Delta(\text{c}) = \text{c} \otimes 1 + 1 \otimes
\text{c} + {\textstyle \sum\limits_{j=1}^{n}} \text{a}_\ell \otimes
\text{b}_\ell \, ,  \hskip9pt  \Delta(\text{b}_i) = \text{b}_i \otimes 1
+ 1 \otimes \text{b}_i  \cr
   \epsilon(\text{a}_i) = 0 \, ,  \hskip5pt  \epsilon(\text{c})
= 0 \, ,  \hskip5pt \epsilon(\text{b}_i) = 0 \, ,  \hskip13pt
S(\text{a}_i) = - \text{a}_i \, ,  \hskip6pt  S(\text{c}) = - \text{c}
+ {\textstyle \sum\limits_{j=1}^{n}} \text{a}_\ell \text{b}_\ell \, ,
\hskip5pt  S(\text{b}_i) = - \text{b}_i  \cr }  $$
and let also  $ \, \F_q[H_n] \, $  be the  $ F(R) $--algebra  obtained
from  $ F_q[H_n] $  by scalar extension.  Then  $ \, \Bbb{B} :=
\Big\{ \prod_{i=1}^n \text{a}_i^{a_i} \cdot \text{c}^c \cdot
\prod_{j=1}^n \text{b}_j^{b_j} \,\Big\vert\, a_i, c, b_j \in
\N \; \forall \, i, j \,\Big\} \, $  is an  $ R $--basis  of
$ F_q[H_n] $,  hence an  $ F(R) $--basis  of  $ \F_q[H_n] $.
     \hbox{Moreover  $ F_q[H_n] $  is a QFA  (at  $ \, \h = q \!
- \! 1 $)  with semiclassical limit  $ F[H_n] \, $.}

\vskip11pt

  {\bf 9.7 Computation of  $ \, {F_q[H_n]}^\vee \, $  and
specialization  $ \, {F_q[H_n]}^\vee \,{\buildrel {q \rightarrow 1}
\over \llongrightarrow}\, U({\gerh_n}^{\!*}) \, $.} \, This section
is devoted to compute  $ \, {F_q[H_n]}^\vee \, $  and its semiclassical
limit (at  $ \, q = 1 \, $).
                                             \par
   Definitions imply that  $ \, \Bbb{B} \setminus \{1\} \, $  is an
$ R $--basis  of  $ \, J := \hbox{\sl Ker} \, \Big( \epsilon \, \colon
\, F_q[H_n] \! \twoheadrightarrow \! R \Big) $,  so  $ \, \big( \Bbb{B}
\setminus \{1\} \big) \cup \big\{ (q-1) \cdot 1 \big\} \, $  is an
$ R $--basis  of  $ \; I := \hbox{\sl Ker} \Big( F_q[H_n] \,
{\buildrel \epsilon \over {\relbar\joinrel\twoheadrightarrow}} \,
R \, {\buildrel {q \mapsto 1} \over {\llongtwoheadrightarrow}}
\, \Bbbk \Big) \, $,  for  $ \, I = J + (q-1) \cdot F_q[H_n] \, $.
Therefore  $ \; {F_q[H_n]}^\vee := \sum_{n \geq 0} {\Big(
{(q-1)}^{-1} I \Big)}^n \; $  is nothing but the unital
$ R $--algebra  (subalgebra of  $ \F_q[H_n] $)  with generators
$ \, E_i := \displaystyle{\, \text{a}_i \, \over \, q - 1 \,} \, $,
$ \, H := \displaystyle{\, \text{c} \, \over \, q - 1\,} \, $,  and
$ \, F_i := \displaystyle{\, \text{b}_i \, \over \, q - 1\,} \, $  ($ \,
i = 1, \dots, n \, $)  and relations (for all  $ \, i, j = 1, \dots,
n \, $)
  $$  E_i E_j = E_j E_i \, ,  \hskip7pt  E_i F_j = F_j E_i \, ,
\hskip7pt  F_i F_j = F_j F_i \, ,  \hskip7pt  H E_i = E_i H + E_i \, ,
\hskip7pt  H F_j = F_j H + F_j  $$
with Hopf algebra structure given by (for all  $ \, i, j = 1,
\dots, n \, $)
 \vskip-13pt
  $$  \displaylines{
  \Delta(E_i) = E_i \otimes 1 + 1 \otimes E_i \, ,  \hskip5pt
\Delta(H) = H \otimes 1 + 1 \otimes H + (q-1) {\textstyle
\sum\limits_{j=1}^{n}} E_j \otimes F_j \, ,  \hskip5pt
\Delta(F_i) = F_i \otimes 1 + 1 \otimes F_i  \cr
%
%
%
   \epsilon(E_i) = \epsilon(H) = \epsilon(F_i) = 0 \, ,  \hskip8pt
S(E_i) = \! - E_i \, ,  \hskip8pt  S(H) = \! - H + (q-1)
{\textstyle \sum\limits_{j=1}^{n}} E_j F_j \, ,  \hskip8pt
S(F_i) = \! - F_i \, .  \cr }  $$
 \vskip-5pt
   At  $ \, q = 1 \, $  this implies that  $ \, {F_q[H_n]}^\vee
\,{\buildrel \, q \rightarrow 1 \, \over \llongrightarrow}\,
U({\gerh_n}^{\! *}) \, $  as co-Poisson Hopf algebras, for a
co-Poisson Hopf algebra isomorphism
 \vskip-5pt
  $$  {F_q[H_n]}^\vee \Big/ (q-1) \, {F_q[H_n]}^\vee
\,{\buildrel \cong \over \llongrightarrow}\, U({\gerh_n}^{\! *})  $$
 \vskip-3pt
\noindent   exists, given by  $ \;\, E_i \mod (q-1) \mapsto \pm
\text{e}_i \, $,  $ \, H \mod (q-1) \mapsto \text{h} \, $,  $ \,
F_i \mod (q-1) \mapsto \text{f}_i \, $,  \; for all  $ \, i, j = 1,
\dots, n \, $;  so  $ \, {F_q[H_n]}^\vee \, $  specializes to  $ \,
U({\gerh_n}^{\! *}) \, $  {\sl as a co-Poisson Hopf algebra},  q.e.d.

\vskip7pt

  {\bf 9.8 The identity  $ \, {\big({F_q[H_n]}^\vee\big)}' = F_q[H_n]
\, $.} \, Finally, we check the validity of the part of  Theorem
2.2{\it (b)}  claiming, when  $ \, p = 0 \, $,  \, that  $ \; H
\in \QFA \,\Longrightarrow\, {\big(H^\vee\big)}' = H \; $  for the
QFA  $ \, H = F_q[H_n] \, $.  Once more the proof works for all
$ \, p \geq 0 \, $,  \, so we do improve  Theorem 2.2{\it (b)}.
                                         \par
 First of all, from definitions induction gives, for all
$ \, m \in \N \, $,
  $$  \Delta^m(E_i) = \hskip-4pt {\textstyle \sum\limits_{r+s=m-1}}
\hskip-4pt 1^{\otimes r} \otimes E_i \otimes 1^{\otimes s} \, ,
\hskip15pt  \Delta^m(F_i) = \hskip-4pt {\textstyle \sum\limits_{r+s
= m-1}}  \hskip-4pt 1^{\otimes r} \otimes F_i \otimes 1^{\otimes s}
\hskip25pt  \forall \; i = 1, \dots, n  $$
 \vskip-21pt
  $$  \Delta^m(H) = \hskip-4pt {\textstyle \sum\limits_{r+s=m-1}}
\hskip-4pt 1^{\otimes r} \otimes H \otimes 1^{\otimes s} +
{\textstyle \sum\limits_{i=1}^m} {\textstyle \sum\limits_{\Sb  j, k
= 1  \\   j < k   \endSb}^m} 1^{\otimes (j-1)} \otimes E_i \otimes
1^{\otimes (k-j-1)} \otimes F_i \otimes 1^{\otimes (m-k)}  $$
 \vskip-1pt
\noindent   so that  $ \; \delta_m(E_i) = \delta_\ell(H) = \delta_m(F_i)
= 0 \; $  for all  $ \, m > 1 $,  $ \ell > 2 \, $  and  $ \, i = 1,
\dots, n \, $;  moreover, for  $ \; \dot{E}_i := (q-1) E_i = \text{a}_i
\, $,  $ \dot{H} := (q-1) H = \text{c} \, $,  $ \dot{F}_i := (q-1) F_i
= \text{b}_i \, $  ($ i = 1, \dots, n $)  one has
  $$  \displaylines{
   \delta_1\big( \dot{E}_i \big) \! = (q-1) E_i \, ,  \; \delta_1\big(
\dot{H} \big) \! = (q-1) H ,  \;  \delta_1\big( \dot{F}_i \big) \!
= (q-1) F_i \in (q-1) \, {F_q[H_n]}^\vee \setminus {(q-1)}^2
{F_q[H_n]}^\vee  \cr
   \delta_2\big( \dot{H} \big) \, = \, {(q-1)}^2 \, {\textstyle
\sum_{i=1}^n} \, E_i \otimes F_i \in {(q-1)}^2 {\big( {F_q[H_n]}^\vee
\big)}^{\otimes 2} \setminus {(q-1)}^3 {\big( {F_q[H_n]}^\vee
\big)}^{\otimes 2} \, .  \cr }  $$
 \vskip-3pt
   The outcome is that  $ \, \dot{E}_i = \text{a}_i ,  \dot{H} =
\text{c} , \dot{F}_i = \text{b}_i \in {\big({F_q[H_n]}^\vee\big)}'
\, $,  so the latter algebra contains the one generated by these
elements, that is  $ F_q[H_n] $.  Even more,  $ {F_q[H_n]}^\vee $
is clearly the  $ R $--span  of the set  $ \, {\Bbb B}^\vee := \Big\{\,
\prod_{i=1}^n E_i^{a_i} \cdot H^c \cdot \prod_{j=1}^n F_j^{b_j}
\,\Big\vert\, a_i, c, b_j \in \N \; \forall \, i, j \,\Big\} \, $,
\, so from this and the previous formulas for $ \Delta^n $  one
gets that  $ \, {\big( {F_q[H_n]}^\vee \big)}' = F_q[H_n] \, $,
\, q.e.d.

\vskip1,9truecm

   \centerline {\bf \S \; 10 \  Fifth example:
%
%
non-commutative Hopf algebra of formal diffeomorphisms }
%

\vskip10pt

  {\bf 10.1 The goal: from ``quantum symmetries'' to ``classical
(geometrical) symmetries''.} \, The purpose of this section is to
give a highly significant example of how the global quantum duality
principle   --- more precisely, the crystal duality principle of
\S 5 ---   may be applied.  We consider a concrete sample, taken
from the theory of non-commutative renormalization of quantum
electro-dynamics (=QED) performed by Brouder and Frabetti in [BF2].
This is just one of several possible examples of the same type:
indeed, several cases of Hopf algebras built out of combinatorial
data have been introduced in last years both in (co)homological
theories (see for instance [LR] and [Fo1--3], and references therein)
and in renormalization studies (starting with [CK1]).  In most cases
these Hopf algebras are neither commutative nor cocommutative, and
our discussion apply almost  {\it verbatim\/}  to them, giving
analogous results.  So the present analysis of the ``toy model''
Hopf algebra of [BF2], can be taken more in general as a pattern
for all those cases.  See also [Ga6].   
                                           \par
   Note that the Hopf algebras under study are usually thought of as
``generalized symmetries'' (or ``quantum symmetries'', in physicists'
terminology); well, the crystal duality principle tells us how to get
out of them   --- via 1-parameter deformations! ---   ``classical
geometric symmetries'', i.e., Poisson groups and Lie bialgebras;
in other words, in a sense this method yields the classical
geometrical counterparts of a quantum symmetry object.

\vskip7pt

   {\bf 10.2 The classical data.} \, Let  $ \Bbbk $  be a fixed field
of characteristic zero.
                                                     \par
   Consider the set  $ \; \G^\dif := \big\{\, x + \sum_{n \geq 1} a_n
\, x^{n+1} \,\big|\; a_n \in \Bbbk \; \forall \; n \in \N_+ \,\big\}
\; $  of all formal series starting with $ x \, $:  endowed with the
composition product, this is a group, which can be seen as the group
of all ``formal diffeomorphisms''   \hbox{$ \, f \, \colon \, \Bbbk
\longrightarrow \Bbbk \, $}  such that  $ \, f(0) = 0 \, $  and  $ \,
f'(0) = 1 \, $  (i.e.~tangent to the identity), also known as the  {\sl
Nottingham group\/}  (see, e.g., [Ca] and references therein).  In fact,
$ \G^\dif $  is an infinite dimensional (pro)affine algebraic group, whose
function algebra  $ \, F \big[ \G^\dif \big] \, $  is generated by the
coordinate functions  $ \, a_n \, $  ($ \, n \in \N_+ $).  Giving to
     each  $ a_n $  the weight\footnote{We say  {\sl weight\/}
        instead of  {\sl degree\/}  because we save the latter term
        for the degree of polynomials.}
$ \, \partial(a_n) := n \, $,  we have that  $ F \big[ \G^\dif \big] $
is an  $ \N $--{\sl graded\/}  Hopf algebra, with polynomial structure
$ \, F \big[ \G^\dif \big] = \Bbbk [a_1,a_2,\dots,a_n,\dots] \, $  and
Hopf algebra structure given by
  $$  \eqalign{
   \Delta(a_n) = a_n \otimes 1 + 1 \otimes a_n + \sum\nolimits_{m=1}^{n-1}
a_m \otimes Q^m_{n-m}(a_*) \, ,  \qquad  \epsilon(a_n) = 0
\hskip15pt  \cr
   S(a_n) = - a_n - \sum\nolimits_{m=1}^{n-1} a_m \,
S\big(Q^m_{n-m}(a_*)\big) = - a_n - \sum\nolimits_{m=1}^{n-1}
S(a_m) \, Q^m_{n-m}(a_*) \cr }  $$
where  $ \; Q_t^\ell(a_*) := \sum_{k=1}^t {{\ell + 1} \choose k}
P_t^{(k)}(a_*) \; $  and  $ \; P_t^{(k)}(a_*) := \sum_{\hskip-4pt
\Sb  j_1, \dots, j_k > 0  \\   j_1 + \cdots + j_k = t  \endSb}
\hskip-5pt  a_{j_1} \cdots a_{j_k} \, $  (the symmetric monic
polynomial of weight  $ m $  and degree  $ k $  in the indeterminates
$ a_j $'s)  for all  $ \, m $,  $ k $,  $ \ell \in \N_+ \, $,  \, and
the formula for  $ S(a_n) $  gives the antipode by recursion.  From now
on, to simplify notation we shall use notation  $ \, \G := \G^\dif \, $
and  $ \, \G_\infty := \G = \G^\dif \, $.  Note also that the tangent
Lie algebra of  $ \G^\dif $  is just the Lie subalgebra  $ \, {W_1}^{\!
\geq 1} = \text{\sl Span}\, \big(\{\, d_n \,|\, n \in \N_+ \,\}\big)
\, $  of the one-sided Witt algebra  $ \, W_1 := \text{\sl Der} \big(
\Bbbk[t] \big) = \text{\sl Span}\, \big( \big\{\, d_n := t^{n+1} {{\,d\,}
\over {\,dt\,}} \,\big|\, n \in \N \cup \{-1\} \,\big\}\big) \, $.
                                                     \par
   In addition, for all  $ \, \nu \in \N_+ \, $  the subset
%
%
$ \, \G^\nu := \big\{\, f \in \G \,\big|\; a_n(f) = 0 , \;
\forall \; n \leq \nu \,\big\} \, $  is a normal subgroup of
$ \G \, $;  \, the corresponding quotient group  $ \, \G_\nu :=
\G \big/ \G^\nu \, $  is unipotent, with dimension  $ \nu $  and
function algebra  $ \, F \big[ \G_\nu \big] \, $  (isomorphic to)
the Hopf subalgebra of  $ F \big[ \G \big] $  generated by  $ \,
a_1 $,  $ \dots $,  $ a_\nu \, $.  In fact, the  $ \G^\nu $'s  form
exactly the lower central series of  $ \G $  (cf.~[Je2]).  Moreover,
$ \G $  is (isomorphic to) the inverse (or projective) limit of these
quotient groups  $ \G_\nu $  ($ \nu \in \N_+ $),  hence  $ \G $  is
pro-unipotent; conversely,  $ F[\G] $  is the direct (or inductive)
limit of the direct system of its graded Hopf subalgebras  $ F[\G_\nu]
\, $  ($ \nu \in \N_+ $).  Finally, the set  $ \; \G^{\text{odd}} :=
\big\{\, f \in \G^\dif \,\big|\; a_{2n+1}(f) = 0 \; \forall \;
n \in \N_+ \,\big\} \; $  is another normal subgroup of
$ \G^\dif $  (the group of  {\sl odd\/}
       formal diffeomorphisms\footnote{\hbox{The fixed-point set of
       the group homomorphism  $ \, \Phi : \G \rightarrow \G \, $,
       $ f \mapsto \Phi(f) \, \big( x \mapsto \big(\Phi(f)\big)(x)
       := -f(-x) \big) $}}
 after [CK3]), whose function algebra  $ F \big[ \G^{\text{odd}} \big] $
is (isomorphic to) the quotient Hopf algebra  $ \, F \big[ \G^\dif \big]
\Big/ \Big( \big\{ a_{2n-1} \big\}_{n \in \N_+} \Big) \, $.  The
latter has the following description: denoting again the cosets of the
$ a_{2n} $'s  with the like symbol, we have  $ \, F \big[ \G^{\text{odd}}
\big] = \Bbbk [a_2, a_4, \dots, a_{2n}, \dots] \, $  with Hopf algebra
structure
  $$  \eqalign{
   \Delta(a_{2n}) = a_{2n} \otimes 1 + 1 \otimes a_{2n} +
\sum\nolimits_{m=1}^{n-1} a_{2m} \otimes \bar{Q}^m_{n-m}(a_{2*})
\, ,  \qquad  \epsilon(a_{2n}) = 0   \hskip35pt  \cr
   S(a_{2n}) = - a_{2n} - \sum\nolimits_{m=1}^{n-1} a_{2m} \,
S\big(\bar{Q}^m_{n-m}(a_*)\big) = - a_{2n} - \sum\nolimits_{m=1}^{n-1}
S(a_{2m}) \, \bar{Q}^m_{n-m}(a_{2*}) \cr }  $$
where  $ \; \bar{Q}_t^\ell(a_{2*}) := \sum_{k=1}^t {{2\ell + 1}
\choose k} \bar{P}_t^{(k)}(a_{2*}) \; $  and  $ \; \bar{P}_t^{(k)}
(a_{2*}) := \sum_{\hskip-4pt  \Sb  j_1, \dots, j_k > 0  \\
j_1 + \cdots + j_k = t  \endSb}  \hskip-5pt  a_{2{}j_1} \cdots
a_{2{}j_k} \, $  for all  $ \, m $,  $ k $,  $ \ell \in \N_+ \, $.  For
each  $ \, \nu \in \N_+ \, $  we can consider also the normal subgroup
$ \, \G^\nu \cap \G^{\text{odd}} \, $  and the corresponding quotient
$ \, \G_\nu^{\text{odd}} := \G^{\text{odd}} \big/ \big( \G^\nu \cap
\G^{\text{odd}} \big) \, $:  \, then  $ F \big[ \G_\nu^{\text{odd}}
\big] $  is (isomorphic to) the quotient Hopf algebra  $ \, F \big[
\G^{\text{odd}} \big] \! \Big/ \! \Big( \! \big\{ a_{2n-1} \big\}_{(2n-1)
\in \N_\nu} \Big) \, $,  \, in particular it is the Hopf sub-\break
 \noindent   algebra of  $ F \big[ \G^{\text{odd}} \big] $  generated
by  $ \, a_2, \dots, a_{2\,[\nu/2]} \, $.  All the  $ F \big[
\G_\nu^{\text{odd}} \big] $'s  are graded Hopf (sub)al-\break
 \noindent   gebras forming a direct system with direct limit
$ F \big[ \G^{\text{odd}} \big] $;  \, conversely, the
$ \G_\nu^{\text{odd}} $'s  form an inverse system with inverse
limit  $ \G^{\text{odd}} $.  In the sequel we write  $ \, \G^+ :=
\G^{\text{odd}} \, $  and  $ \, \G^+_\nu := \G_\nu^{\text{odd}} \, $.
                                                     \par
   For each  $ \, \nu \in \N_+ \, $,  \, set  $ \, \N_\nu := \{1, \dots,
\nu\} \, $;  set also  $ \, \N_\infty := \N_+ \, $.  For each  $ \, \nu
\in \N_+ \cup \{\infty\} \, $,  \, let  $ \, \L_\nu = \L(\N_\nu) \, $
be the free Lie algebra over  $ \Bbbk $  generated by  $ {\{x_n\}}_{n
\in \N_\nu} $  and let  $ \, U_\nu = U(\L_\nu) \, $  be its universal
enveloping algebra; let also  $ \, V_\nu = V(\N_\nu) \, $  be the
$ \Bbbk $--vector  space with basis  $ \, {\{ x_n \}}_{n \in \N_\nu}
\, $,  \, and let  $ \, T_\nu = T(V_\nu) \, $  be its associated tensor
algebra.  Then there are canonical identifications  $ \, U(\L_\nu)
= T(V_\nu) = \Bbbk \big\langle \{\, x_n \,|\, n \in \N_\nu \,\}
\big\rangle \, $,  \, the latter being the unital  $ \Bbbk $--algebra
of non-commutative polynomials in the set of indeterminates  $ \,
{\{x_n\}}_{n \in \N_\nu} \, $,  \, and $ \L_\nu $  is just the Lie
subalgebra of  $ U_\nu = T_\nu $  generated by  $ \, {\{x_n\}}_{n \in
\N_\nu} \, $.  Moreover,  $ \L_\nu $  has a basis  $ B_\nu $  made of
Lie monomials in the  $ x_n $'s  ($ n \in \N_\nu $),  like  $ [x_{n_1},
x_{n_2}] $,  $ [[x_{n_1}, x_{n_2}], x_{n_3}] $,  $ [[[x_{n_1}, x_{n_2}],
x_{n_3}], x_{n_4}] $,  etc.: details can be found e.g.~in [Re], Ch.~4--5.
In the sequel I shall use these identifications with no further mention.
We consider on  $ \, U(\L_\nu) \, $  the standard Hopf algebra structure
given by  $ \; \Delta(x) = x \otimes 1 + 1 \otimes x \, $,  $ \,
\epsilon(x) = 0 \, $,  $ \, S(x) = -x \; $  for all  $ \, x \in
\L_\nu \, $,  \, which is also determined by the same formulas for
$ \, x \in {\{x_n\}}_{n \in \N_\nu} \, $  alone.  By construction
$ \, \nu \leq \mu \, $  implies  $ \, \L_\nu \subseteq \L_\mu \, $,
\, whence the  $ \L_\nu $'s  form a direct system (of Lie algebras)
whose direct limit is exactly  $ \L_\infty \, $;  \, similarly,
$ U(\L_\infty) $  is the direct limit of all the  $ U(\L_\nu) $'s.
Finally, with  $ \Bbb{B}_\nu $  we shall mean the obvious PBW-like
basis of  $ U(\L_\nu) $  w.r.t.~some fixed total order  $ \preceq $
of  $ B_\nu $,  namely  $ \, \Bbb{B}_\nu := \big\{\, x_{\underline{b}}
\;\big|\, \underline{b} = b_1 \cdots b_k \, ; \, b_1, \dots, b_k
\in B_\nu \, ; \, b_1 \preceq \cdots \preceq b_k \,\big\} \, $.
%
%
%
%
                                                     \par
   The same construction applies to define the corresponding ``odd''
objects, based on  $ {\{x_n\}}_{n \in \N^+_\nu} $,  with  $ \, \N^+_\nu
:= \N_\nu \cap 2\,\N $,  instead of  $ {\{x_n\}}_{n \in \N_\nu} $  (for
each  $ \, \nu \in \N \cup \{\infty\} \, $).  Thus we have  $ \, \L^+_\nu
= \L(\N^+_\nu) \, $,  $ \, U^+_\nu = U(\L^+_\nu) \, $,  \, $ \, V^+_\nu
= V(\N^+_\nu) \, $,  $ \, T^+_\nu = T(V^+_\nu) \, $,  \, with the obvious
canonical identifications  $ \, U(\L^+_\nu) = T(V^+_\nu) = \Bbbk
\big\langle \{\, x_n \,|\, n \in \N^+_\nu \,\} \big\rangle \, $;  \,
moreover,  $ \L^+_\nu $  has a basis  $ B^+_\nu $  made of Lie monomials
in the  $ x_n $'s  ($ n \in \N^+_\nu $),  etc.  The  $ \L^+_\nu $'s
form a direct system whose direct limit is  $ \L^+_\infty \, $,  \, and
$ U(\L^+_\infty) $  is the direct limit of all the  $ U(\L^+_\nu) $'s.

\vskip3pt

   {\it  $ \underline{\text{Warning}} \, $:}  \, in the sequel, we
shall often deal with subsets  $ {\{\hbox{\bf y}_b\}}_{b \in B_\nu} $
(of some algebra) in bijection with  $ B_\nu \, $,  the fixed basis
of  $ \L_\nu \, $.  Then we shall write things like  $ \, \hbox{\bf
y}_\lambda \, $  with  $ \, \lambda \in \L_\nu \, $:  \, this means
we extend the bijection  $ \, {\{\hbox{\bf y}_b\}}_{b \in B_\nu} \cong
B_\nu \, $  to  $ \, \text{\sl Span}\,\big( \{ \hbox{\bf y}_b \}_{b \in
B_\nu} \big) \cong \L_\nu \, $  by linearity, so that  $ \; \hbox{\bf
y}_\lambda \cong \sum_{b \in B_\nu} c_b \, b \; $  iff  $ \; \lambda
= \sum_{b \in B_\nu} c_b \, b \; $  ($ \, c_b \in \Bbbk \, $).  The
same kind of convention will be applied with  $ B^+_\nu $  instead
of  $ B_\nu $  and  $ \L^+_\nu $  instead of  $ \L_\nu \, $.

\vskip7pt

   {\bf 10.3 The noncommutative Hopf algebra of formal diffeomorphisms.}
\, For all  $ \, \nu \in \N_+ \cup \{\infty\} \, $,  \, let  $ \,
\calH_\nu \, $  be the Hopf  $ \Bbbk $--algebra  given as follows:
as a  $ \Bbbk $--algebra it is simply  $ \, \calH_\nu := \Bbbk
\big\langle \{\, \a_n \,|\, n \in \N_\nu \,\} \big\rangle \, $
(the  $ \Bbbk $--algebra  of non-commutative polynomials in the
set of indeterminates  $ \, {\{\a_n\}}_{n \in \N_\nu} \, $),  and its
Hopf algebra structure is given by (for all  $ \, n \in \N_\nu \, $)
 \vskip-15pt
  $$  \hbox{ $ \eqalign{
   \Delta(\a_n) = \a_n \otimes 1 + 1 \otimes \a_n +
\sum\nolimits_{m=1}^{n-1} \a_m \otimes Q^m_{n-m}(\a_*) \, ,
\qquad  \epsilon(\a_n) = 0  \hskip15pt  \cr
   S(\a_n) = -\a_n - \sum\nolimits_{m=1}^{n-1} \a_m \,
S\big(Q^m_{n-m}(\a_*)\big) = -\a_n - \sum\nolimits_{m=1}^{n-1}
S(\a_m) \, Q^m_{n-m}(\a_*)  \cr } $ }   \eqno  (10.1)  $$
 \vskip-2pt
\noindent   
 (notation like in \S 10.2) where the latter formula yields the
antipode by recursion.  Moreover,  $ \calH_\nu $  is in fact an  {\sl
$ \N $--graded  Hopf algebra},  once generators have been given degree
--- in the sequel called  {\sl weight\/}  ---   by the rule  $ \,
\partial(\a_n) := n \, $  (for all  $ \, n \in \N_\nu \, $).  By
construction the various  $ \calH_\nu $'s  (for all  $ \, \nu \in
\N_+ \, $)  form a direct system, whose direct limit is  $ \calH_\infty
\, $:  \, the latter was
    originally introduced\footnote{However, the formulas in [BF2]
    give the  {\sl opposite\/}  coproduct, hence change the antipode
    accordingly; we made the present choice to make these formulas
    ``fit well'' with those for  $ F \big[ \G^{\text{dif}} \big] $
    (see below).}
in [BF2], \S 5.1 (with  $ \, \Bbbk = \C \, $),  under the name
$ \, \calH^\dif \, $.
                                                     \par
   Similarly, for all  $ \, \nu \in \N_+ \cup \{\infty\} \, $  we
set  $ \, \calK_\nu := \Bbbk \big\langle \{\, \a_n \,|\, n \in
\N^+_\nu \,\} \big\rangle \, $  (where  $ \, \N^+_\nu := \N_\nu
\cap (2\,\N) \, $):  this bears a Hopf algebra structure given
by (for all  $ \, 2\,n \in \N^+_\nu \, $)
 \vskip-15pt
  $$  \hbox{ $ \eqalign{
   \Delta(\a_{2n}) = \a_{2n} \otimes 1 + 1 \otimes \a_{2n} +
\sum\nolimits_{m=1}^{n-1} \a_{2m} \otimes \bar{Q}^m_{n-m}(\a_{2*})
\, ,  \qquad  \epsilon(\a_{2n}) = 0  \hskip35pt  \cr
   S(\a_{2n}) = -\a_{2n} - \sum\nolimits_{m=1}^{n-1} \a_{2m}
\, S\big(\bar{Q}^m_{n-m}(\a_{2*})\big) = -\a_{2n} -
\sum\nolimits_{m=1}^{n-1} S(\a_{2m}) \, Q^m_{n-m}(\a_{2*}) 
\cr } $ }  $$
 \vskip-2pt
\noindent   
 (notation of \S 10.2).  Indeed, this is an  $ \N $--graded  Hopf algebra
where generators have degree   --- called  {\sl weight\/}  ---   given
by  $ \, \partial(\a_n) := n \, $  (for all  $ \, n \in \N^+_\nu \, $).
All the  $ \calK_\nu $'s  
   form a direct\break   
 \eject   
\noindent   
 system with direct limit
$ \calK_\infty \, $.  Finally, for each  $ \, \nu \in \N^+_\nu \, $
there is a graded Hopf algebra epimorphism  $ \, \calH_\nu
\longtwoheadrightarrow \calK_\nu \, $  given by  $ \, \a_{2n}
\mapsto \a_{2n} \, $,  $ \, \a_{2m+1} \mapsto 0 \, $  for
all  $ \, 2n, 2m+1 \in \N_\nu \, $.
                                             \par
   Definitions and \S 10.2 imply that   
  $$  {\big(\calH_\nu\big)}_{\text{\it ab}} \, := \,
\calH_\nu \Big/ \big( \big[ \calH_\nu, \calH_\nu \,\big] \big)
\; \cong \; F \big[ \G_\nu \big] \, ,   \qquad   \hbox{via}
\qquad   \a_n \mapsto a_n  \quad  \forall \; n \in \N_\nu  $$  
as  $ \N $--graded  Hopf algebras: in other words, the abelianization
of  $ \calH_\nu $  is nothing but  $ F \big[ \G_\nu \big] $.  Thus in
a sense one can think at  $ \calH_\nu $  as a  {\sl non-commutative
version}  (indeed, the ``coarsest'' one) of  $ F \big[ \G_\nu \big] $,
hence as a ``quantization'' of  $ \G_\nu $  itself: however, this is
{\sl not\/}  a quantization in the sense we mean in this paper, for
$ F \big[ \G_\nu \big] $  is attained through abelianization, not
through specialization (of some deformation parameter).  Similarly we
have also
  $$  {\big(\calK_\nu\big)}_{\text{\it ab}} \, := \,
\calK_\nu \Big/ \big( \big[ \calK_\nu, \calK_\nu \,\big] \big)
\; \cong \; F \big[ \G^+_\nu \big] \, ,   \qquad   \hbox{via}
\qquad   \a_{2n} \mapsto a_{2n}  \quad  \forall \; 2n \in \N^+_\nu  $$
as  $ \N $--graded  Hopf algebras: in other words, the abelianization
of  $ \calK_\nu $  is just  $ F \big[ \G^+_\nu \big] $.
                                             \par
   Note that  $ \, \big(\calH_\nu\big)^\vee = \calH_\nu = \big(
\calH_\nu \big)' \, $  (notation of \S 5.1) because  $ \calH_\nu $
is graded and connected: therefore applying the crystal duality
principle to  $ \calH_\nu $  we'll end up with (5.5), which means
we can deform  $ \calH_\nu $  in four different ways to Hopf algebras
bearing some (Poisson-type) geometrical content; and similarly for
$ \calK_\nu $.  In particular we'll describe the Poisson groups
$ G_+ $  and  $ G_-^\star $,  and their cotangent Lie bialgebras
$ \gerg_+^\times $  and  $ \gerg_- $,  attached to  $ \calH_\nu $
and to  $ \calK_\nu $  in this way.  We perform the analysis
explicitly for  $ \calH_\nu \, $;  \, the case  $ \calK_\nu $
is the like, and we leave to the reader the easy task to fill
in details.
                                                     \par
   We follow the recipe in \S\S 5.1--4.  Let's drop the subscript
$ \nu $  (which stands fixed) and write  $ \, \calH := \calH_\nu \, $.
Let  $ \, R := \Bbbk[\h] \, $,  \, and set  $ \; \calH_\h := \calH[\h]
\equiv \Bbbk[\h] \otimes_{\Bbbk} \calH \, $:  \, this is a Hopf algebra
over  $ \Bbbk[\h] $,  namely  $ \, \calH_\h = \Bbbk[\h] \big\langle
\{\, \a_n \,|\, n \in \N_\nu \,\} \big\rangle \, $  with Hopf structure
given by (10.1) again.  More precisely, we have  $ \, \calH[\h] \in \HA
\, $  w.r.t.~the ground ring  $ \, R := \Bbbk[\h] \, $  (a PID).  Then
   \hbox{$ \, F(R) = \Bbbk(\h) \, $,  \, and  $ \, {(\calH_\h)}_F
:= \Bbbk(\h) \otimes_{\Bbbk[\h]} \calH_\h = \Bbbk(\h) \otimes_{\Bbbk}
\calH = \calH(\h) = \Bbbk(\h) \big\langle \{\, \a_n \,|\, n \in \N_\nu
\,\} \big\rangle $.}

\vskip7pt

   {\bf 10.4 Drinfeld's algebra  $ \, {\calH_\h}^{\!\vee} :=
{\big( \calH[\h] \big)}^\vee \, $.} \, By the method in \S 5
leading to the Crystal Duality Principle, we can apply Drinfeld's
functors at the prime  $ \, \h \in \Bbbk[\h] \, $  to  $ \, \calH_\h
:= \calH[\h] \, $.  We begin with  $ \; {\calH_\h}^{\!\vee} :=
\sum_{n \geq 0} \h^{-n} J^n \; \big( \, \subseteq {(\calH_\h)}_F =
\calH(\h) \, \big) \, $,  \, where  $ \, J := \text{\sl Ker}\,\big(
\epsilon_{\scriptscriptstyle \calH_\h} \, \colon \calH_\h \longrightarrow
\Bbbk[\h] \, \big) \, $.  We'll describe  $ {\calH_\h}^{\!\vee} $
explicitly, thus checking that it is really a QrUEA, as predicted by
Theorem 2.2{\it (a)\/};  then we'll look at its specialization at  $ \,
\h = 1 \, $,  and finally we'll study  $ \big( {\calH_\h}^{\!\vee}
\big)' $  and its specializations at  $ \, \h = 0 \, $  and  $ \,
\h = 1 \, $.  The outcome will be an explicit description of the
diagram of deformations (5.3) for  $ \, H = \calH \, $  ($ \,
= \calH_\nu \, $).
                                      \par
   For all  $ \, n \in \N_\nu \, $,  \, set  $ \, \x_n := \h^{-1}
\a_n \, $.  Then clearly  $ \, {\calH_\h}^{\!\vee} \, $  is the
$ \Bbbk[\h] $--subalgebra of  $ \, \calH(\h) \, $  generated by
the set  $ \, {\{\x_n\}}_{n \in \N_\nu} \, $,  \, and thus  $ \;
\displaystyle{ {\calH_\h}^{\!\vee} \, = \, \Bbbk[\h] \big\langle
\{\, \x_n \,|\, n \in \N_\nu \,\} \big\rangle } \, $.  Moreover,
 \eject   
  $$  \hskip-7pt   \hbox{ $ \eqalign{
   \Delta(\hbox{\bf x}_n)  &  = \hbox{\bf x}_n \! \otimes \! 1
+ 1 \! \otimes \! \hbox{\bf x}_n + \sum\nolimits_{m=1}^{n-1}
\sum\nolimits_{k=1}^m \h^k {{n-m+1} \choose k} \, \hbox{\bf x}_{n-m}
\otimes P^{(k)}_m(\hbox{\bf x}_*) \, ,  \hskip6pt
\epsilon(\x_n) = 0  \cr
   S(\hbox{\bf x}_n)  &  = - \hbox{\bf x}_n -
\sum\nolimits_{m=1}^{n-1} \sum\nolimits_{k=1}^m
\h^k {{n-m+1} \choose k} \, \hbox{\bf x}_{n-m} \,
S \big( P^{(k)}_m(\hbox{\bf x}_*) \big) =  \cr
   &  \hskip61pt   = - \hbox{\bf x}_n - \sum\nolimits_{m=1}^{n-1}
\sum\nolimits_{k=1}^m \h^k {{n-m+1} \choose k} \,
S(\hbox{\bf x}_{n-m}) \, P^{(k)}_m(\hbox{\bf x}_*)  \cr } $ }
\hskip-18pt   (10.2)  $$   
for all  $ \, n \in \N_\nu \, $,  due to (10.1); from this one
sees by hands that the following holds:

\vskip7pt

\proclaim{Proposition 10.5} \, Formulas (10.2) make  $ \, {\calH_\h}^{\!
\vee} = \Bbbk[\h] \big\langle \{\, \x_n, \,|\, n \in \N_\nu \,\}
\big\rangle \, $  into a graded Hopf\/ $ \, \Bbbk[\h] $--algebra,
embedded into  $ \, \calH(\h) := \Bbbk(\h) \otimes_{\Bbbk} \calH \, $
as a graded Hopf subalgebra.  Moreover,  $ {\calH_\h}^{\!\vee} $  is
a deformation of  $ \, \calH $,
   \hbox{for its specialization at  $ \, \h = 1 \, $
is isomorphic to  $ \calH $,  i.e.}
  $$  {\calH_\h}^{\!\vee}{\Big|}_{\h=1} := \, {\calH_\h}^{\!\vee}
\Big/ (\h\!-\!1) \, {\calH_\h}^{\!\vee} \, \cong \, \calH  \quad
\text{via}  \quad  \x_n \!\!\mod (\h\!-\!1) \, {\calH_\h}^{\!\vee}
\, \mapsto \, \a_n  \quad  (\,\forall \;\;  n \in \N_\nu \,)  $$
as graded Hopf algebras over\/  $ \Bbbk \, $.   \qed
\endproclaim

\vskip7pt

   {\sl  $ \underline{\hbox{\it Remark}} $:}  \; The previous result
shows that  $ \calH_\h $  {\sl is a deformation of  $ \calH $,  which
is ``recovered'' as specialization limit (of  $ \calH_\h $)  at}  $ \,
\h = 1 \, $.  The next result instead shows that  $ \calH_\h $  {\sl
is also a deformation of  $ U(\L_\nu) $,  which is ``recovered'' as
specialization limit at}  $ \, \, \h = 0 \, $.  Altogether, this
gives the left-hand-side of (5.3) for  $ \, H = \calH := \calH_\nu
\, $,  \, with  $ \, \gerg_- = \L_\nu \, $.

\vskip7pt

\proclaim{Theorem 10.6} \,  $ {\calH_\h}^{\!\vee} $  is a QrUEA at  $ \,
\h = 0 \, $.  Namely, the specialization limit of  $ {\calH_\h}^{\!
\vee} $  at  $ \, \h \! = \! 0 \, $  is  $ \; {\calH_\h}^{\!\vee}
{\Big|}_{\h=0} := \, {\calH_\h}^{\!\vee} \Big/ \h \, {\calH_\h}^{\!\vee}
\, \cong \, U(\L_\nu) \; $  via  $ \; \x_n \mod \h \, {\calH_\h}^{\!\vee}
\mapsto x_n \; $  for all  $ \; n \in \N_\nu \, $,  \, thus inducing on
$ U(\L_\nu) $  the structure of co-Poisson Hopf algebra uniquely given
by the Lie bialgebra structure on  $ \L_\nu $  given by  $ \, \delta(x_n)
= \sum_{\ell=1}^{n-1} (\ell+1) \, x_\ell \wedge x_{n-\ell} \, $  (for
      all  $ \, n \in \N_\nu $)\,\footnote{Hereafter, I use notation
                   $ \, a \wedge b := a \otimes b - b \otimes a \, $.}.
In particular in the diagram (5.3) for  $ \, H = \calH \; (= \calH_\nu)
\, $  we have  $ \, \gerg_- = \L_\nu \, $.
                                   \hfill\break
   \indent   Finally, the grading  $ d $  given by  $ \, d(x_n) := 1
\;\, (n \in \N_+) \, $  makes  $ \, {\calH_\h}^{\!\vee}{\Big|}_{\h=0}
\! \cong U(\L_\nu) \, $  into a graded co-Poisson Hopf algebra;
similarly, the grading  $ \, \partial $  given by  $ \, \partial(x_n)
:= n \;\, (n \in \N_+) \, $  makes  $ \, {\calH_\h}^{\!\vee}
{\Big|}_{\h=0} \! \cong U(\L_\nu) \, $  into a graded Hopf algebra
and  $ \L_\nu $  into a graded Lie bialgebra.
\endproclaim

\demo{Proof}  First observe that since  $ \, {\calH_\h}^{\!\vee} =
\Bbbk[\h] \, \big\langle \{\, \x_n \,|\, n \in \N_\nu \,\} \big\rangle
\, $  and  $ \, U(\L_\nu) = T(V_\nu) = \Bbbk \, \big\langle \{\, x_n
\,|\, n \in \N_\nu \,\} \big\rangle \, $  mapping $ \; \x_n \mod \h \,
{\calH_\h}^{\!\vee} \mapsto x_n \; $  ($ \, \forall \; n \in \N_\nu
\, $)  does really define an isomorphism  {\sl of algebras}  $ \;
\Phi \, \colon \, {\calH_\h}^{\!\vee} \Big/ \h \, {\calH_\h}^{\!\vee}
\, \cong \, U(\L_\nu) \, $.  Second, formulas (10.2) give
  $$  \displaylines{
   \Delta(\x_n) \equiv \x_n \otimes 1 + 1 \otimes \x_n  \mod \h \,
\Big( {\calH_\h}^{\!\vee} \otimes {\calH_\h}^{\!\vee} \Big)  \cr
   \epsilon(\x_n) \equiv 0  \mod \h \, \Bbbk[\h] \, ,  \qquad
S(\x_n) \equiv - \x_n  \mod \h \, {\calH_\h}^{\!\vee}  \cr }  $$
for all  $ \, n \in \N_\nu \, $;  \, comparing with the standard Hopf
structure of  $ U(\L_\nu) $  this shows that  $ \Phi $  is in fact
an isomorphism  {\sl of Hopf algebras\/}  too.  Finally, as
$ {\calH_\h}^{\!\vee}{\Big|}_{\h=0} $  is cocommutative, a
Poisson co-bracket is defined on it by the standard recipe
in Remark 1.5: applying it yields
  $$  \eqalign{
   \delta(x_n) \;  &  := \; \big( \h^{-1} \, \big( \Delta(\x_n) -
\Delta^{\text{op}}(\x_n) \big) \big) \mod \h \, \Big( {\calH_\h}^{\!
\vee} \otimes {\calH_\h}^{\!\vee} \Big) \; =   \hskip75pt   \cr
                   &  \phantom{:}= \; {\textstyle \sum_{m=1}^{n-1}} \,
{\textstyle {{n-m+1} \choose 1}} \, x_{n-m} \wedge P_m^{(1)}(x_*) \;
= \; {\textstyle \sum_{\ell=1}^{n-1}} \, (\ell+1) \, x_\ell \wedge
x_{n-\ell}  \qquad  \forall \;\; n \in \N_\nu \; .   \quad
\hskip-3pt  \square  \cr }  $$
\enddemo

\vskip7pt

  {\bf 10.7 Drinfeld's algebra  $ \big( {\calH_\h}^{\!\vee} \big)' $.}
\, I look now at the other Drinfeld's functor (at  $ \h \, $),  and
consider  $ \; \big( {\calH_\h}^{\!\vee} \big)' := \Big\{\, \eta \in
{\calH_\h}^{\!\vee} \;\Big\vert\;\, \delta_n(\eta) \in \h^n \big(
{\calH_\h}^{\!\vee} \big)^{\otimes n} \; \forall\; n \in \N \,\Big\}
\; $  ($ \, \subseteq {\calH_\h}^{\!\vee} \, $).  Theorem 2.2 tells
us that  {\it  $ \big( {\calH_\h}^{\!\vee} \big)' $  is a Hopf\/
$ \Bbbk[\h] $--subalgebra  of  $ {\calH_\h}^{\!\vee} $,  and the
specialization of  $ \big( {\calH_\h}^{\!\vee} \big)' $  at  $ \,
\h = 0 \, $,  \, that is  $ \; \big( {\calH_\h}^{\!\vee} \big)'
{\Big|}_{\h=0} := \, \big( {\calH_\h}^{\!\vee} \big)' \Big/ \h
\, \big( {\calH_\h}^{\!\vee} \big)' \, $,  \; is the function
algebra of a connected algebraic Poisson group  $ {G_{\!\L_\nu}
\phantom{|}}^{\hskip-8pt \star} $  {\sl dual}  to}  $ \, G_{\!\L_\nu}
\, $,  the latter being the connected simply-connected Poisson algebraic
group with tangent Lie bialgebra  $ \L_\nu \, $.  In other words,
$ \, \big( {\calH_\h}^{\!\vee} \big)'{\Big|}_{\h=0} \, $  must be
isomorphic (as a Poisson Hopf algebra) to  $ F\big[{G_{\!\L_\nu}
\phantom{|}}^{\hskip-8pt \star}\big] $,  where  $ {G_{\!\L_\nu}
\phantom{|}}^{\hskip-8pt \star} $  is connected and has
cotangent Lie bialgebra  $ \, \text{\it Lie}\,\big({G_{\!\L_\nu}
\phantom{|}}^{\hskip-8pt \star}\big) = \L_\nu \, $.  Therefore we
must prove that  $ \big( {\calH_\h}^{\!\vee} \big)'{\Big|}_{\h=0} $
is a commutative Hopf $ \Bbbk $--algebra,  it has no non-trivial
idempotents, and  $ \, \Big( \text{\it co-Lie}\,\big({G_{\!\L_\nu}
\phantom{|}}^{\hskip-8pt \star}\big) := \Big) $
                $ \, J_0 \big/ J_0^{\,2}
\, \cong \, \L_\nu \, $  as Lie bialgebras,  where  $ \, J_0 :=
\text{\sl Ker}\, \Big( \epsilon \, \colon \big( {\calH_\h}^{\!\vee}
\big)'{\Big|}_{\h=0} \!\!\! \longrightarrow \Bbbk \Big) \, $.
We prove all this directly, via explicit description
of  $ \big( {\calH_\h}^{\!\vee} \big)' $  and
its specialization at  $ \, \h = 0 \, $.
                                              \par
   {\it Step I:} \,  {\it A direct check shows that  $ \, \tilde{\x}_n
:= \h \, \x_n = \a_n \in \big( {\calH_\h}^{\!\vee} \big)' \, $,  \, for
all  $ \, n \in \N_\nu \, $}.  Indeed, we have of course  $ \, \delta_0
(\tilde{\x}_n) = \epsilon(\tilde{\x}_n) \in \h^0 \, {\calH_\h}^{\!\vee}
\, $  and  $ \, \delta_1(\tilde{\x}_n) = \tilde{\x}_n - \epsilon
(\tilde{\x}_n) \in \h^1 \, {\calH_\h}^{\!\vee} \, $.  Moreover,
$ \, \delta_2(\tilde{\x}_n) = \sum_{m=1}^{n-1} \tilde{\x}_{n-m}
\otimes Q_m^{n-m} (\tilde{\x}_*) = \sum_{m=1}^{n-1} \sum_{k=1}^m
\h^{k+1} {{n-m+1} \choose k} \hbox{\bf x}_{n-m} \otimes P^{(k)}_m
(\hbox{\bf x}_*) \in \h^2 \Big( {\calH_\h}^{\!\vee} \otimes
{\calH_\h}^{\!\vee} \Big) \, $.  Since in general  $ \, \delta_\ell
= \big( \delta_{\ell-1} \otimes \id \big) \circ \delta_2 \, $  for
all  $ \, \ell \in \N_+ \, $,  \, we have
  $$  \delta_\ell(\tilde{\x}_n) = \big( \delta_{\ell-1} \otimes \id \big)
\big( \delta_2 (\tilde{\x}_n) \big) \, = \, \sum_{m=1}^{n-1} \sum_{k=1}^m
\, \h^k {{n-m+1} \choose k} \, \delta_{\ell-1}(\hbox{\bf x}_{n-m})
\otimes P^{(k)}_m(\hbox{\bf x}_*)  $$
whence induction gives  $ \, \delta_\ell(\tilde{\x}_n) \in \h^\ell
\, \big( {\calH_\h}^{\!\vee} \big)^{\otimes \ell} \, $  for all  $ \,
\ell \in \N \, $,  \, thus  $ \, \tilde{\x}_n \in \big( {\calH_\h}^{\!
\vee} \big)' \, $,  \, q.e.d.
                                              \par
   {\it Step II:} \, By  Theorem 2.2{\it (a)}  we have that  $ \big(
{\calH_\h}^{\!\vee} \big)'{\Big|}_{\h=0} $  is commutative: this means
$ \, [a,b] \equiv 0 \mod \h \, \big( {\calH_\h}^{\!\vee} \big)' \, $,
\, that is  $ \, [a,b] \in \h \, \big( {\calH_\h}^{\!\vee} \big)' \, $
hence also  $ \, \h^{-1} [a,b] \in \big( {\calH_\h}^{\!\vee} \big)' \, $,
\, for all  $ \, a $,  $ b \in \big( {\calH_\h}^{\!\vee} \big)' \, $.
{\it In particular,}  we get  $ \, \widetilde{[\x_n,\x_m]} := \h
\, [\x_n,\x_m] = \h^{-1} [\tilde{\x}_n,\tilde{\x}_m] \in \big(
{\calH_\h}^{\!\vee} \big)' \, $  for all  $ \, n $, $ m \in \N_\nu
\, $,  \, whence iterating (and recalling  $ \L_\nu $ is generated by
the  $ \x_n $'s)  {\it we get $ \, \tilde{\x} := \h \, \x \in \big(
{\calH_\h}^{\!\vee} \big)' \, $  for every  $ \, \x \in \L_\nu \, $}.
Hereafter we identify  $ \L_\nu $  with its image via the embedding
$ \, \L_\nu \hookrightarrow U(\L_\nu) = \Bbbk \big\langle {\{x_n\}}_{n
\in \N_\nu} \big\rangle \hookrightarrow \Bbbk[\h] \big\langle
{\{\x_n\}}_{n \in \N_\nu} \big\rangle = {\calH_\h}^{\!\vee} \, $
given by  $ \, x_n \mapsto \x_n \, $  ($ \, n \in \N_\nu \, $).
                                              \par
   {\it Step III:} \, The previous step showed that, if we embed  $ \,
\L_\nu \hookrightarrow U(\L_\nu) \hookrightarrow {\calH_\h}^{\!\vee}
\, $  via $ \, x_n \mapsto \x_n \, $  (for all  $ \, n \in \N_\nu \, $)
we find  $ \, \widetilde{\L_\nu} := \h \, \L_\nu \subseteq \big(
{\calH_\h}^{\!\vee} \big)' \, $.  {\it Let  $ \, \big\langle
\widetilde{\L_\nu} \,\big\rangle \, $  be the
$ \Bbbk[\h] $--subalgebra  of  $ \big( {\calH_\h}^{\!\vee} \big)' $
generated by  $ \widetilde{\L_\nu} \, $:  then  $ \, \big\langle
\widetilde{\L_\nu} \,\big\rangle \subseteq \big( {\calH_\h}^{\!\vee}
\big)' \, $},  \, because  $ \big( {\calH_\h}^{\!\vee} \big)' $  is
a subalgebra.  In particular, if  $ \, \b_b \in {\calH_\h}^{\!\vee}
\, $  is the image of any  $ \, b \in B_\nu \, $  (cf.~\S 10.2) we
have  $ \, \widetilde\b_b := \h \, \b_b \in \big( {\calH_\h}^{\!\vee}
\big)' \, $.
%
%
 \eject   
   {\it Step IV:} \,  {\it Conversely to  {\it Step III},  we have
$ \, \big\langle \widetilde{\L_\nu} \,\big\rangle \supseteq \big(
{\calH_\h}^{\!\vee} \big)' \, $}.  In fact, let  $ \, \eta \in \big(
{\calH_\h}^{\!\vee} \big)' \, $;  \, then there are unique  $ \, d
\in \N \, $,  $ \, \eta_+ \in {\calH_\h}^{\!\vee} \setminus \h \,
{\calH_\h}^{\!\vee} \, $  such that  $ \, \eta = \h^d \eta_+ \, $;
\, set also  $ \, \bar{y} := y \mod \h \, {H_\h}^{\!\vee} \in {H_\h}^{\!
\vee} \Big/ \h \, {H_\h}^{\!\vee} \, $  for all  $ \, y \in {H_\h}^{\!\vee}
\, $.  As  $ \, {\calH_\h}^{\!\vee} = \Bbbk[\h] \big\langle \{\, \x_n \,|\,
n \in \N_\nu \,\} \big\rangle \, $  there is a unique  $ \h $--adic
expansion  of  $ \eta_+ $,  namely  $ \; \eta_+ = \eta_0 + \h \, \eta_1
+ \cdots + \h^s \, \eta_s = \sum_{k=0}^s \h^k \, \eta_k \; $  with all  $ \;
\eta_k \in \Bbbk \big\langle \{\, \x_n \,|\, n \in \N_\nu \,\} \big\rangle
\; $  and  $ \; \eta_0 \not= 0 \, $.  Then  $ \, \bar\eta_+ = \bar\eta_0
:= \eta_0 \mod \h \, {\calH_\h}^{\!\vee} \, $; \, thus  Lemma 4.2{\it
(d)\/}  gives  $ \, \partial(\bar\eta_0) \leq d \, $,  \, where now
$ \partial(\bar\eta_0) $  denotes the degree of  $ \bar\eta_0 $  for
the standard filtration of  $ U(\L_\nu) $.  By the PBW theorem,
$ \partial(\bar\eta_0) $  is also the degree of  $ \bar\eta_0 $  as
a polynomial in the  $ \bar{\x}_b $'s,  hence also of  $ \eta_0 $
as a polynomial in the  $ \x_b $'s  ($ b \in B_\nu $):  \, then
$ \, \h^d \, \eta_0 \in \big\langle \widetilde{\L_\nu} \,\big\rangle
\subseteq \big( {\calH_\h}^{\!\vee} \big)' \, $  (using  {\it Step
III\/}),  hence we find
 \vskip-10pt
  $$  \eta_{(1)} \, := \, \h^{d+1} \, \big( \eta_1 + \h \, \eta_2
+ \cdots + \h^{s-1} \, \eta_s \big) = \eta - \h^d \, \eta_0 \in
\big( {\calH_\h}^{\!\vee} \big)' \, .  $$
 \vskip-1pt
\noindent   
 Thus we can apply our argument again, with  $ \eta_{(1)} $  instead
of  $ \eta $.  Iterating we find  $ \, \partial(\bar\eta_k) \leq
d+k \, $,  \, whence  $ \, \h^{d+k} \, \eta_k \in \big\langle
\widetilde{\L_\nu} \,\big\rangle \, $  $ \Big( \subseteq \big(
{\calH_\h}^{\!\vee} \big)' \Big) $  for all  $ k \, $,
\hbox{thus  $ \, \eta = \sum_{k=0}^s \h^{d+k}
   \, \eta_k \, \in \big\langle \widetilde{\L_\nu} \,\big\rangle \, $,
   \, q.e.d.}

\vskip3pt

   An entirely similar analysis clearly works with  $ \calK_\h $
taking the role of  $ \calH_\h $,  with similar results ({\sl
mutatis mutandis\/}).  On the upshot, we get the following
description:
%
%

\vskip11pt

\proclaim {Theorem 10.8} \, (a) \, With notation of  {\it Step II\/}
in \S 10.7 (and  $ \, [a,c\,] := a \, c - c \, a \, $),  we have
 \vskip-15pt
  $$  \big( {\calH_\h}^{\!\vee} \big)' \; = \; \Big\langle
\widetilde{\L_\nu} \;\Big\rangle \; = \; \Bbbk[\h] \, \Big\langle
{\big\{\, \widetilde\b_b \,\big\}}_{b \in B_\nu} \Big\rangle \Bigg/
\bigg( \Big\{\, \Big[ \widetilde\b_{b_1}, \widetilde\b_{b_2} \Big]
- \h \; \widetilde{\big[\,\b_{b_1},\b_{b_2}\big]} \;\Big|\;
\forall \; b_1, b_2 \in B_\nu \;\Big\} \bigg) \; .  $$
 \vskip-5pt
   (b) \,  $ \big( {\calH_\h}^{\!\vee} \big)' $  is a graded
Hopf\/  $ \, \Bbbk[\h] $--subalgebra  of  $ \, {\calH_\h}^{\!\vee}
\, $,  \, and  $ \, \calH $  is naturally embedded into  $ \, \big(
{\calH_\h}^{\!\vee} \big)' $  as a graded Hopf subalgebra via  $ \,\;
\calH \lhook\joinrel\loongrightarrow \big( {\calH_\h}^{\!\vee}
\big)' \, $,  $ \, \a_n \mapsto \tilde\x_n \; $  (for all
$ \, n \in \N_\nu $).
                                        \hfill\break
  \indent   (c) \,  $ \big( {\calH_\h}^{\!\vee} \big)'{\Big|}_{\h=0} :=
\big( {\calH_\h}^{\!\vee} \big)' \Big/ \h \, \big( {\calH_\h}^{\!\vee}
\big)' = F \big[ {G_{\!\L_\nu}\phantom{|}}^{\hskip-8pt \star} \big]
\, $,  \, where  $ \, {G_{\!\L_\nu}\phantom{|}}^{\hskip-8pt \star}
\, $  is an infinite dimensional connected Poisson algebraic group
with cotangent Lie bialgebra isomorphic to  $ \L_\nu $  (with the
graded Lie bialgebra structure of Theorem 10.6).  Indeed,  $ \big(
{\calH_\h}^{\!\vee} \big)'{\Big|}_{\h=0} $  is the free Poisson
(commutative) algebra over  $ \N_\nu \, $,  generated by all the\/
$ \tilde\x_n{\big|}_{\h=0} $  ($ \, n \in \N_\nu \, $)  with Hopf
structure given by (10.1) with\/  $ \tilde\x_* $  instead of\/
$ \a_* \, $.  Thus  $ \big( {\calH_\h}^{\!\vee} \big)'{\Big|}_{\h=0} $
is the polynomial algebra  $ \, \Bbbk \big[ {\{\, \beta_b \,\}}_{b \in
B_\nu} \big] \, $  generated by a set of indeterminates  $ \, {\{\,
\beta_b \,\}}_{b \in B_\nu} \, $  in bijection with the basis
$ B_\nu $  of\/  $ \L_\nu \, $,  so  $ \; {G_{\!\L_\nu}
\phantom{|}}^{\hskip-8pt \star} \cong \Bbb{A}_\Bbbk^{B_\nu}
\, $  (a (pro)affine\/  $ \Bbbk $--space)  as algebraic varieties.
Finally,  $ \, F \big[ {G_{\!\L_\nu}\phantom{|}}^{\hskip-8pt \star}
\big] = \big( {\calH_\h}^{\!\vee} \big)'{\Big|}_{\h=0} \cong \Bbbk
\big[ {\{\, \beta_b \,\}}_{b \in B_\nu} \big] \, $  bears the natural
{\sl algebra grading}  $ \, d $  of polynomial algebras and the  {\sl
Hopf algebra grading}  inherited from  $ \big( {\calH_\h}^{\!\vee}
\big)' $,  respectively given by  $ \, d\big(\widetilde\b_b\big) =
1 \, $  and  $ \, \partial\big(\widetilde\b_b\big) = \sum_{i=1}^k
n_i \, $  for all  $ \, b = [[\cdots[[x_{n_1},x_{n_2}],x_{n_3}],
\cdots],x_{n_k}] \in B_\nu \, $.
                                        \hfill\break
  \indent   (d) \,  $ F \big[ \G_\nu \big] $  is naturally embedded
into  $ \, \big( {\calH_\h}^{\!\vee} \big)'{\Big|}_{\h=0} = F \big[
{G_{\!\L_\nu}\phantom{|}}^{\hskip-8pt \star} \big] \, $  as a graded
Hopf subalgebra via  $ \; \mu \, \colon \, F \big[ \G_\nu \big]
\lhook\joinrel\loongrightarrow \big( {\calH_\h}^{\!\vee} \big)'
{\Big|}_{\h=0} = F \big[ {G_{\!\L_\nu} \phantom{|}}^{\hskip-8pt
\star} \big] \, $,  $ \, a_n \mapsto \Big( \, \tilde\x_n \! \mod
\h \, \big( {\calH_\h}^{\!\vee} \big)' \Big) \; $  (for all  $ \,
n \in \N_\nu \, $);  moreover,  $ F \big[ \G_\nu \big] $  freely
generates  $ F \big[ {G_{\!\L_\nu}\phantom{|}}^{\hskip-8pt \star}
\big] $  as a\/  {\rm Poisson}  algebra.   Thus there is an
algebraic group epimorphism  $ \, \mu_* \, \colon \,
{G_{\!\L_\nu} \phantom{|}}^{\hskip-8pt \star} \,
{\relbar\joinrel\relbar\joinrel\twoheadrightarrow}
\; \G_\nu \, $,  \, that is  $ \, {G_{\!\L_\nu}
\phantom{|}}^{\hskip-8pt \star} \, $  is
an extension of  $ \, \G_\nu \, $.
                                        \hfill\break
  \indent   (e) \,  Mapping  $ \; \Big( \, \tilde\x_n \! \mod
\h \, \big( {\calH_\h}^{\!\vee} \big)' \Big) \mapsto a_n \; $
(for all  $ \, n \in \N_\nu $)  gives a well-defined graded
Hopf algebra epimorphism  $ \; \pi \, \colon \, F \big[
{G_{\!\L_\nu}\phantom{|}}^{\hskip-8pt \star} \big] \,
{\relbar\joinrel\relbar\joinrel\twoheadrightarrow} \, F
\big[ \G_\nu \big] \, $.  Thus there is an algebraic
group monomorphism  $ \, \pi_* \, \colon \, \G_\nu
\lhook\joinrel\loongrightarrow\, {G_{\!\L_\nu}
\phantom{|}}^{\hskip-8pt \star} \, $,  \, that is
$ \, \G_\nu \, $  is an algebraic subgroup of  $ \,
{G_{\!\L_\nu} \phantom{|}}^{\hskip-8pt \star} \, $.
                                        \hfill\break
  \indent   (f) \,  The map  $ \mu $  is a section of  $ \pi $,
hence  $ \pi_* $  is a section of  $ \mu_* \, $.  Thus  $ \,
{G_{\!\L_\nu}\phantom{|}}^{\hskip-8pt \star} \, $  is a
semidirect product of algebraic groups, namely  $ \; {G_{\!\L_\nu}
\phantom{|}}^{\hskip-8pt \star} = \, \G_\nu \ltimes \Cal{N}_\nu \; $
where  $ \, \Cal{N}_\nu := \text{\sl Ker}\,(\mu_*) \trianglelefteq
{G_{\!\L_\nu}\phantom{|}}^{\hskip-8pt \star} \, $.
                                        \hfill\break
  \indent   (g) \,  The analogues of statements (a)--(f) hold
with  $ \calK $  instead of  $ \, \calH \, $,  \, with 
$ X^+ $  instead of  $ X $  for all  $ \, X = \L_\nu,
B_\nu, \N_\nu, \mu, \pi, \Cal{N}_\nu \, $,  and with 
$ {G_{\!\L_\nu^+}\phantom{|}}^{\hskip-8pt \star} $
instead of  $ {G_{\!\L_\nu}\phantom{|}}^{\hskip-8pt \star} $.
%
%
\endproclaim

\demo{Proof} {\it (a)} \, This part follows directly from
{\it Step III\/ {\rm and}  Step IV\/}  in \S 10.7.
                                              \par
   {\it (b)} \, To show that  $ \big( {\calH_\h}^{\!\vee} \big)' $
is a graded Hopf subalgebra we use its presentation in  {\it (a)}.
But first observe that by construction $ \, \a_n = \tilde\x_n \, $
(for all  $ \, n \in \N_\nu \, $),  so  $ \, \calH \, $  embeds into
$ \big( {\calH_\h}^{\!\vee} \big)' $  via an embedding
which is compatible with the Hopf operations: then this will be
a Hopf algebra monomorphism,  {\sl up to proving that  $ \big(
{\calH_\h}^{\!\vee} \big)' $  is a Hopf subalgebra\/}  (of
$ {\calH_\h}^{\!\vee} \, $).
                                              \par
   Now,  $ \epsilon_{{\calH_\h}^{\!\vee}} $  obviously restricts to
give a counit for  $ \big( {\calH_\h}^{\!\vee} \big)' $.  Second, we
show that  $ \, \Delta \Big( \big( {\calH_\h}^{\!\vee} \big)' \Big)
\subseteq \big( {\calH_\h}^{\!\vee} \big)' \otimes \big( {\calH_\h}^{\!
\vee} \big)' \, $,  \, so  $ \Delta $  restricts to a coproduct for
$ \big( {\calH_\h}^{\!\vee} \big)' $.  Indeed, each  $ \, b \in B_\nu
\, $  is a Lie monomial, say  $ \, b = [[[ \dots [x_{n_1}, x_{n_2}],
x_{n_3}], \dots], x_{n_k}] \, $  for some  $ \, k $,  $ n_1 $,
$ \dots $,  $ n_k \in \N_\nu \, $,  where  $ k $  is its Lie
degree: by induction on  $ k \, $  we'll prove  $ \, \Delta \big(
\widetilde\b_b \big) \in \big( {\calH_\h}^{\!\vee} \big)' \otimes
\big( {\calH_\h}^{\!\vee} \big)' \, $  (with  $ \, \widetilde\b_b
:= \h \, \b_b = \h \, [[[ \dots [\x_{n_1}, \x_{n_2}], \x_{n_3}],
\dots], \x_{n_k}] \, $).
                                              \par
   If  $ \, k = 1 \, $  then  $ \, b = x_n \, $  for some  $ \, n
\in \N_\nu \, $.  Then  $ \, \widetilde\b_b = \h \, \x_n = \a_n
\, $  and
  $$  \Delta\Big(\widetilde\b_b\Big) = \Delta(\a_n) = \a_n \otimes 1
+ 1 \otimes \a_n + \sum_{m=1}^{n-1} \a_{n-m} \otimes Q_m^{n-m}(\a_*)
\, \in \calH^{\text{dif}} \otimes \calH^{\text{dif}} \subseteq \big(
{\calH_\h}^{\!\vee} \big)' \otimes \big( {\calH_\h}^{\!\vee} \big)'
\, .  $$
   \indent   If  $ \, k > 1 \, $  then  $ \, b = [b^-,x_n] \, $  for
some  $ \, n \in \N_\nu \, $  and some  $ \, b^- \in B_\nu \, $
expressed by a Lie monomial of degree  $ \, k-1 \, $.  Then  $ \,
\widetilde\b_b = \h \, [\b^-,\x_n] = \Big[ \widetilde\b^-,\x_n \Big]
\, $  and
  $$  \displaylines{
   \Delta\Big(\widetilde\b_b\Big) \; = \; \Delta\Big( \Big[\widetilde\b^-,
\x_n \Big] \Big) \; = \; \Big[ \Delta\Big(\widetilde\b^-\Big),
\Delta(\x_n) \Big] \; = \; \h^{-1} \, \Big[ \Delta\Big(
\widetilde\b^-\Big), \Delta(\a_n) \Big] \; =   \hfill  \cr
   \hfill   = \; \h^{-1} \, \left[\; \sum\nolimits_{\left(
\widetilde\b^- \right)} \widetilde\b^-_{(1)} \otimes
\widetilde\b^-_{(2)} \; , \; \a_n \otimes 1 + 1 \otimes \a_n +
\sum\nolimits_{m=1}^{n-1} \a_{n-m} \otimes Q_m^{n-m}(\a_*)
\;\right] \; =  \cr
   = \; \sum\nolimits_{\left(\widetilde\b^-\right)} \h^{-1} \,
\left[\, \widetilde\b^-_{(1)} \, , \, \a_n \,\right] \otimes
\widetilde\b^-_{(2)} \; + \; \sum\nolimits_{\left(\widetilde\b^-\right)}
\widetilde\b^-_{(1)} \otimes \h^{-1} \, \left[\,
\widetilde\b^-_{(2)} \, , \, \a_n \,\right] \; +   \hfill  \cr
   \hfill   + \; \sum_{\left(\widetilde\b^-\right)} \sum_{m=1}^{n-1}
\bigg(\, \h^{-1} \left[\, \widetilde\b^-_{(1)} \, , \, \a_{n-m}
\,\right] \otimes \widetilde\b^-_{(2)} \, Q_m^{n-m}(\a_*) \; + \;
\widetilde\b^-_{(1)} \, \a_{n-m} \otimes \h^{-1} \left[\,
\widetilde\b^-_{(2)} \, , \, Q_m^{n-m}(\a_*) \,\right]
\hskip-1pt \bigg)  \cr }  $$
where we used the standard  $ \Sigma $--notation  for  $ \, \Delta
\Big(\widetilde\b^-\Big) = \sum\nolimits_{\left( \widetilde\b^- \right)}
\widetilde\b^-_{(1)} \otimes \widetilde\b^-_{(2)} \, $.  By inductive
hypothesis we have  $ \, \widetilde\b^-_{(1)} $,  $ \widetilde\b^-_{(2)}
\in \big( {\calH_\h}^{\!\vee} \big)' \, $;  \, then since also  $ \,
\a_\ell \in \big( {\calH_\h}^{\!\vee} \big)' \, $  for all  $ \ell $
and since  $ \big( {\calH_\h}^{\!\vee} \big)' $  is commutative
modulo $ \h $  we have
  $$  \h^{-1} \left[\, \widetilde\b^-_{(1)} \, , \, \a_n \,\right]
\, ,  \; \h^{-1} \left[\, \widetilde\b^-_{(2)} \, , \, \a_n
\,\right] \, ,  \; \h^{-1} \left[\, \widetilde\b^-_{(1)} \, ,
\, \a_{n-m} \,\right] \, ,  \; \h^{-1} \left[\, \widetilde\b^-_{(2)}
\, , \, Q_m^{n-m}(\a_*) \,\right] \, \in \, \big( {\calH_\h}^{\!\vee}
\big)'  $$
for all  $ n $  and  $ (n-m) $  above: so the previous formula
   \hbox{gives  $ \, \Delta\big(\,\widetilde\b_b\big)
\in \big( {\calH_\h}^{\!\vee} \big)' \! \otimes
\big( {\calH_\h}^{\!\vee} \big)' \, $,  \, q.e.d.}
                                              \par
   Finally, for the antipode we proceed as above.  Let  $ \, b \in
B_\nu \, $  be the Lie monomial  $ \, b = [[[ \dots [x_{n_1},x_{n_2}],
x_{n_3}], \dots], x_{n_k}] \, $,  \, so  $ \, \widetilde\b_b = \h
\, \b_b = \h \, [[[ \dots [\x_{n_1}, \x_{n_2}], \x_{n_3}], \dots],
\x_{n_k}] \, $.  We prove that  $ \, S \big( \, \widetilde\b_b \big)
\in \big( {\calH_\h}^{\!\vee} \big)' \, $  by induction on the
degree  $ k \, $.
                                              \par
   If  $ \, k = 1 \, $  then  $ \, b = x_n \, $  for some  $ n \, $,
so  $ \, \widetilde\b_b = \h \, \x_n = \a_n \, $  and
  $$  S\Big(\widetilde\b_b\Big) = S(\a_n) = - \a_n -
\sum\nolimits_{m=1}^{n-1} \a_{n-m} \, S\big(Q_m^{n-m}(\a_*)\big)
\; \in \; \calH^{\text{dif}} \subseteq \big( {\calH_\h}^{\!\vee}
\big)' \; ,  \quad  \text{q.e.d.}  $$
   \indent   If  $ \, k > 1 \, $  then  $ \, b = [b^-,x_n] \, $  for
some  $ \, n \in \N_\nu \, $  and some  $ \, b^- \in B_\nu \, $  which
is a Lie monomial of degree  $ \, k-1 \, $.  Then  $ \, \widetilde\b_b
= \h \, [\b^-,\x_n] = \Big[ \widetilde\b^-, \x_n \Big] = \h^{-1} \,
\Big[ \widetilde\b^-, \a_n \Big] \, $  and so
  $$  S\Big(\widetilde\b_b\Big) \, = \, S\Big( \Big[\widetilde\b^-,
\x_n \Big] \Big) \, = \, \h^{-1} \Big[ S(\a_n), S \big( \,
\widetilde\b^- \big) \Big] \; \in \; \h^{-1} \Big[ \big(
{\calH_\h}^{\!\vee} \big)', \big( {\calH_\h}^{\!\vee} \big)'
\Big] \; \subseteq \; \big( {\calH_\h}^{\!\vee} \big)'  $$
using the fact  $ \, S(\a_n) = S\big(\widetilde\x_n\big) =
S\big(\,\widetilde\b_{x_n}\big) \in \big( {\calH_\h}^{\!\vee} \big)'
\, $  (by the case $ \, k \! = \! 1 \, $)  along with the inductive
assumption  $ \, S\big(\,\widetilde\b^-\big) \in \big( {\calH_\h}^{\!
\vee} \big)' \, $  and the commutativity of  $ \big( {\calH_\h}^{\!\vee}
\big)' $  modulo  $ \h \, $.
                                              \par
   {\it (c)} \, As a consequence of  {\it (a)},  $ \,
\big( {\calH_\h}^{\!\vee} \big)'{\Big|}_{\h=0} $  is
a  {\sl polynomial\/}  $ \Bbbk $--algebra,  namely
  $$  \big( {\calH_\h}^{\!\vee} \big)'{\Big|}_{\h=0} \; = \;
\Bbbk \, \Big[ {\big\{\,\beta_b\,\big\}}_{b \in B} \Big]  \qquad
\text{with}  \qquad  \beta_b := \widetilde\b_b \mod \h \,
\big( {\calH_\h}^{\!\vee} \big)'  \quad  \text{for all}
\;\; b \in B_\nu \, .  $$
So  $ \big( {\calH_\h}^{\!\vee} \big)'{\Big|}_{\h=0} $  is the
algebra of regular functions  $ F[\varGamma] $ of some (affine)
algebraic variety  $ \varGamma \, $;  \, as $ \big( {\calH_\h}^{\!
\vee} \big)' $  is a Hopf algebra the same is true for  $ \, \big(
{\calH_\h}^{\!\vee} \big)'{\Big|}_{\h=0} = F[\varGamma] \, $,  \,
so  $ \varGamma $  is an (affine) algebraic group; and since  $ \,
F[\varGamma] = \big( {\calH_\h}^{\!\vee} \big)'{\Big|}_{\h=0} \, $
is a specialization limit of  $ \big( {\calH_\h}^{\!\vee} \big)' $,
it is endowed with a Poisson structure too, hence  $ \varGamma $
is a  {\sl Poisson\/}  (affine) algebraic group.
                                              \par
   We compute the cotangent Lie bialgebra of  $ \varGamma $.  First,
$ \, \germ_e := \text{\sl Ker}\, \big( \epsilon_{F[\varGamma]}
\big) = \Big( {\big\{\, \beta_b \,\big\}}_{b \in B_\nu} \Big) \, $
(the ideal generated by the $ \beta_b $'s)  by construction, so
$ \, {\germ_e}^{\hskip-3pt 2} = \Big( \! {\big\{\, \beta_{b_1}
\beta_{b_2} \,\big\}}_{b_1, b_2 \in B_\nu} \Big) $.  Therefore the
cotangent Lie bialgebra  $ \, Q\big(F[\varGamma]\big) := \germ_e \Big/
{\germ_e}^{\hskip-3pt 2} \, $  as a  $ \Bbbk $--vector space has basis
$ \, {\big\{\,\overline\beta_b \,\big\}}_{b \in B_\nu} \, $  where
$ \, \overline\beta_b := \beta_b \mod {{\germ_e}^{\hskip-3pt 2}}
\, $  for all  $ \, b \in B_\nu \, $.  For its Lie bracket
we have (cf.~Remark 1.5)
  $$  \displaylines{
   \big[ \, \overline\beta_{b_1}, \overline\beta_{b_2} \big] \;
:= \; \big\{ \beta_{b_1}, \beta_{b_2} \big\}  \hskip-3pt  \mod
{{\germ_e}}^{\hskip-3pt 2} \; = \; \Big( \h^{-1} \big[\,
\widetilde\b_{b_1}, \widetilde\b_{b_2} \big]  \hskip-3pt  \mod
\h \, \big( {\calH_\h}^{\!\vee} \big)' \Big)  \hskip-3pt  \mod
{{\germ_e}^{\hskip-3pt 2}} \; =   \hfill  \cr
   = \; \Big( \h^{-1} \h^2 \, \big[\,\b_{b_1},\b_{b_2}\big]
\hskip-3pt  \mod \! \h \, \big( {\calH_\h}^{\!\vee} \big)' \Big)
\hskip-3pt  \mod {{\germ_e}^{\hskip-3pt 2}} \; = \, \Big( \h \,
\b_{[b_1,b_2]} \hskip-3pt  \mod \! \h \, \big( {\calH_\h}^{\!\vee}
\big)' \Big)  \hskip-3pt  \mod  {{\germ_e}^{\hskip-3pt 2}} \; =  \cr
   \hfill   = \; \Big( \widetilde\b_{[b_1,b_2]}  \hskip-3pt  \mod
\h \, \big( {\calH_\h}^{\!\vee} \big)' \Big)  \hskip-3pt  \mod
{{\germ_e}^{\hskip-3pt 2}} \; = \; \beta_{[b_1,b_2]}  \hskip-3pt
\mod {{\germ_e}^{\hskip-3pt 2}} \; = \; \overline\beta_{[b_1,b_2]}
\; ,  \cr }  $$
thus the  $ \Bbbk $--linear  map  $ \; \Psi \, \colon \, \L_\nu
\longrightarrow \germ_e \big/ {\germ_e}^{\hskip-3pt 2} \; $
defined by  $ \, b \mapsto \overline\beta_b \, $  for all
$ \, b \in B_\nu \, $  is a  {\sl Lie algebra isomorphism}.  As for
the Lie cobracket, using the general identity  $ \, \delta = \Delta -
\Delta^{\text{op}} \mod \big( {{\germ_e}^{\hskip-3pt 2}} \otimes
F[\varGamma] + F[\varGamma] \otimes {{\germ_e}^{\hskip-3pt 2}} \big)
\, $  (written  $ \hskip-3pt \mod \widehat{{\germ_e}^{\hskip-3pt 2}}
\, $  for short) we get, for all  $ \, n \in \N_\nu \, $,
  $$  \displaylines{
   \delta\big(\,\overline\beta_{x_n}\big) = \big( \Delta \! - \!
\Delta^{\text{op}} \big) (\beta_{x_n})  \hskip-5pt  \mod  \hskip-2pt
\widehat{{\germ_e}^{\hskip-3pt 2}} =  \hskip-1pt  \Big(  \hskip-2pt
\big(  \hskip-1pt  \Delta \! - \! \Delta^{\text{op}}  \hskip-1pt
\big) (\tilde\x_n)  \hskip-5pt \mod \hskip-1,6pt  \h \Big( \! \big(
{\calH_\h}^{\!\vee} \big)' \! \otimes \big( {\calH_\h}^{\!\vee}
\big)' \Big)  \hskip-2pt  \Big)  \hskip-5pt  \mod  \hskip-2pt
\widehat{{\germ_e}^{\hskip-3pt 2}} =   \hfill  \cr
   = \; \left( \left( \a_n \wedge 1 \, + \, 1 \wedge \a_n \, + \,
\sum\nolimits_{m=1}^{n-1} \a_{n-m} \wedge Q_m^{n-m}(\a_*) \right)
\hskip-3pt  \mod \h \, \big( {\calH_\h}' \otimes {\calH_\h}'
\big) \right)  \hskip-3pt  \mod \widehat{{\germ_e}^{\hskip-3pt 2}}
\; =  \cr
   \hfill   = \; \left( \beta_{x_n} \wedge 1 \, + \, 1 \wedge
\beta_{x_n} \, + \, \sum\nolimits_{m=1}^{n-1} \beta_{x_{n-m}}
\wedge Q_m^{n-m}(\beta_{x_*}) \right)  \hskip-3pt  \mod
\widehat{{\germ_e}^{\hskip-3pt 2}} =  \cr
   \quad   = \; \left( \sum\nolimits_{m=1}^{n-1} \beta_{x_{n-m}}
\wedge Q_m^{n-m}(\beta_{x_*}) \right)  \hskip-3pt  \mod
\widehat{{\germ_e}^{\hskip-3pt 2}} \; =   \hfill  \cr
   \hfill   = \; \left( \sum\nolimits_{m=1}^{n-1}
\sum\nolimits_{k=1}^m {n-m+1 \choose k} \, \beta_{x_{n-m}}
\wedge P_m^{(k)}(\beta_{x_*}) \right)  \hskip-3pt  \mod
\widehat{{\germ_e}^{\hskip-3pt 2}} \; =  \cr
 }  $$   
  $$  \displaylines{
   \quad   = \; \left( \sum\nolimits_{m=1}^{n-1} {n-m+1 \choose 1} \,
\beta_{x_{n-m}} \wedge P_m^{(1)}(\beta_{x_*}) \right)  \hskip-3pt
\mod \widehat{{\germ_e}^{\hskip-3pt 2}} \; =   \hfill  \cr
   \hfill   = \; \sum\nolimits_{m=1}^{n-1} {n-m+1 \choose 1}
\; \overline\beta_{x_{n-m}} \wedge \overline\beta_{x_m} \; = \;
\sum\nolimits_{\ell=1}^{n-1} (\ell+1) \; \overline\beta_{x_\ell}
\wedge \overline\beta_{x_{n-\ell}}  \cr }  $$
because   --- among other things ---   one has  $ \; P_m^{(k)}
(\beta_{x_*}) \in {{\germ_e}^{\hskip-3pt 2}} \; $  for all
$ \, k > 1 \, $:  \, therefore
  $$  \delta\big(\,\overline\beta_{x_n}\big) \; = \;
\sum\nolimits_{\ell=1}^{n-1} (\ell+1) \; \overline\beta_{x_\ell}
\wedge \overline\beta_{x_{n-\ell}}  \qquad  \forall  \hskip7pt
n \in \N_\nu \; .   \eqno (10.3)  $$
Since  $ \L_\nu $  is generated (as a Lie algebra) by the  $ x_n $'s,
the last formula shows that the map  $ \; \Psi \, \colon \, \L_\nu
\longrightarrow \germ_e \big/ {\germ_e}^{\hskip-3pt 2} \; $  given
above is also an isomorphism  {\sl of Lie bialgebras},  \, q.e.d.
                                              \par
   Finally, the statements about gradings of  $ \, \big(
{\calH_\h}^{\!\vee} \big)'{\Big|}_{\h=0} $  should be trivially
clear.
                                              \par
   {\it (d)} \, The part about Hopf algebras is a direct consequence
of  {\it (a)\/}  and  {\it (b)},  noting that the  $ \tilde\x_n $'s
commute modulo  $ \, \h \, \big( {\calH_\h}^{\!\vee} \big)' \, $,
\, since  $ \big( {\calH_\h}^{\!\vee} \big)'{\Big|}_{\h=0} $  is
commutative.  Then, taking spectra (i.e.~sets of characters of each
Hopf algebras) we get (functorially) an algebraic group morphism
$ \, \mu_* \, \colon \, {G_{\!\L_\nu}\phantom{|}}^{\hskip-8pt \star}
\,{\relbar\joinrel\relbar\joinrel\rightarrow}\; \G_\nu \, $,  \, which
in fact happens to be  {\sl onto\/}  because, due to the special
polynomial form of these algebras, each character of  $ \, F \big[
\G_\nu \big] \, $  does extend to a character  of  $ \, F \big[
{G_{\!\L_\nu}\phantom{|}}^{\hskip-8pt \star} \big] \, $,  hence
the former does arise from restriction of the latter.
                                              \par
   {\it (e)} \, Due to the explicit description of  $ F \big[
{G_{\!\L_\nu} \phantom{|}}^{\hskip-8pt \star} \big] $  coming from
{\it (a)\/}  and  {\it (b)},  mapping  $ \; \Big( \, \tilde\x_n \! \mod
\h \, \big( {\calH_\h}^{\!\vee} \big)' \Big) \mapsto a_n \; $  (for all 
$ \, n \in \N_\nu $)  clearly yields a well-defined Hopf algebra epimorphism
$ \; \pi \, \colon \, F \big[ {G_{\!\L_\nu} \phantom{|}}^{\hskip-6pt
\star} \big] \, {\relbar\joinrel\relbar\joinrel\twoheadrightarrow}
\, F \big[ \G_\nu \big] \, $  (w.r.t.~the trivial Poisson bracket
on the right-hand-side)  is again a routine matter.  Then taking
spectra gives a monomorphism  $ \, \pi_* \, \colon \, \G_\nu
\,\lhook\joinrel\loongrightarrow\, {G_{\!\L_\nu}
\phantom{|}}^{\hskip-8pt \star} \, $  of algebraic
groups as required.
                                              \par
   {\it (f)} \, The map  $ \mu $  is a section of  $ \pi $  by
construction.  Then clearly  $ \pi_* $  is a section of  $ \mu_* \, $,
which implies  $ \; {G_{\!\L_\nu}\phantom{|}}^{\hskip-8pt \star}
= \, \G_\nu \ltimes \Cal{N}_\nu \; $  (with  $ \, \Cal{N}_\nu
:= \text{\sl Ker}\,(\mu_*) \trianglelefteq {G_{\!\L_\nu}
\phantom{|}}^{\hskip-8pt \star} \, $)  by general theory.
                                              \par
   {\it (g)} \, This ought to be clear from the whole discussion,
for all arguments apply again   --- {\sl mutatis mutandis\/} ---
when starting with  $ \calK $  instead of  $ \, \calH \, $;  \,
details are left to the reader.   \qed
\enddemo
%
%
 \eject   

   {\sl $ \underline{\hbox{\it Remark}} $:}  \; Roughly speaking,
we can say that the extension  $ \, F \big[ \G_\nu \big]
\,\lhook\joinrel\loongrightarrow\, F \big[ {G_{\!\L_\nu}
\phantom{|}}^{\hskip-8pt \star} \big] \, $  is performed simply
by adding to  $ F \big[ \G_\nu \big] $  a  {\sl free\/}  Poisson
structure, which happens to be compatible with the Hopf structure.
Then the Poisson bracket starting from the ``elementary'' coordinates
$ a_n $  (for  $ \, n \in \N_\nu \, $)  freely generates  {\sl new\/}
coordinates  $ \{a_{n_1},a_{n_2}\} $,  $ \big\{\! \{a_{n_1},
a_{n_2}\}, a_{n_3} \big\} $,  etc., thus enlarging  $ F \big[
\G_\nu \big] $  and generating  $ F \big[ {G_{\!\L_\nu}
\phantom{|}}^{\hskip-8pt \star} \big] $.  At the group level,
this means that  $ \G_\nu $  {\sl freely Poisson-generates\/}
the Poisson group  $ {G_{\!\L_\nu}\phantom{|}}^{\hskip-8pt \star}
\, $:  technically speaking, new 1-parameter subgroups, which are
build up in a  {\sl ``Poisson-free''\/}  manner from those attached
to the  $ a_n $'s,  are freely ``pasted'' to  $ \G_\nu \, $,  thus
expanding it and so building up  $ {G_{\!\L_\nu}
\phantom{|}}^{\hskip-8pt \star} \, $.  Then the algebraic group
epimorphism  $ \, {G_{\!\L_\nu}\phantom{|}}^{\hskip-8pt \star} \,
{\buildrel \mu_* \over {\relbar\joinrel\relbar\joinrel\twoheadrightarrow}}
\; \G_\nu \, $  is just a ``forgetful map'': it kills the new 1-parameter
subgroups and is injective (hence an isomorphism) on the subgroup
generated by the old ones.  On the other hand, definitions imply that
$ \, F \big[ {G_{\!\L_\nu}\phantom{|}}^{\hskip-8pt \star} \big] \Big/
\Big( \big\{ F \big[ {G_{\!\L_\nu}\phantom{|}}^{\hskip-8pt \star} \big],
F\big[{G_{\!\L_\nu}\phantom{|}}^{\hskip-8pt \star}\big] \big\} \Big)
\cong F \big[ \G_\nu \big] \, $,  \, and with this identification the
map  $ \, F \big[ {G_{\!\L_\nu}\phantom{|}}^{\hskip-8pt \star} \big] \,
{\buildrel \pi \over {\relbar\joinrel\relbar\joinrel\twoheadrightarrow}}
\; F \big[ \G_\nu \big] \, $  is just the canonical map, which ``mods
out'' all Poisson brakets  $ \{f_1,f_2\} $,  for  $ \, f_1, f_2 \in
F \big[ {G_{\!\L_\nu}\phantom{|}}^{\hskip-8pt \star} \big] \, $.

\vskip7pt

   {\bf 10.9 Specialization limits.} \, So far, we have already
pointed out (by Proposition 10.5, Theorem 10.6,  Theorem 10.8{\it
(c)\/})  the following specialization limits of  $ {\calH_\h}^{\!
\vee} $  and  $ \big( {\calH_\h}^{\!\vee} \big)' \, $:
  $$  {\calH_\h}^{\!\vee} \;{\buildrel \h \rightarrow 1 \over
\llongrightarrow}\; \calH \; ,  \qquad  {\calH_\h}^{\!\vee}
\;{\buildrel \h \rightarrow 0 \over \llongrightarrow}\; U(\L_\nu) \; ,
\qquad  \big( {\calH_\h}^{\!\vee} \big)' \;{\buildrel \h \rightarrow 0
\over \llongrightarrow}\; F \big[ {G_{\!\L_\nu} \phantom{|}}^{\hskip-8pt
\star} \big]  $$
as graded Hopf  $ \Bbbk $--algebras,  with some (co-)Poisson structures
in the last two cases.  As for the specialization limit of  $ \big(
{\calH_\h}^{\!\vee} \big)' $  at  $ \, \h = 1 \, $,  \, Theorem 10.8
implies that it is $ \calH $.  Indeed, by  Theorem 10.8{\it (b)\/}
$ \, \calH $  embeds into  $ \big( {\calH_\h}^{\!\vee} \big)' $  via
$ \; \a_n \mapsto \tilde\x_n \; $  (for all  $ \, n \in \N_\nu \, $):
\, then
  $$  [\a_n,\a_m] \, = \, \big[ \tilde\x_n,\tilde\x_m \big] \, = \,
\h \, \widetilde{[\x_n,\x_m]} \, \equiv \, \widetilde{[\x_n,\x_m]}
\mod (\h \! - \! 1) \, \big( {\calH_\h}^{\!\vee} \big)'   \eqno
\big(\, \forall \; n, m \in \N_\nu \big)  $$
whence, due to the presentation of  $ \big( {\calH_\h}^{\!\vee}
\big)' $  by generators and relations in  Theorem 10.8{\it (a)\/},
  $$  \big( {\calH_\h}^{\!\vee} \big)'{\Big|}_{\h=1} \; := \;
\big( {\calH_\h}^{\!\vee} \big)' \Big/ (\h \! - \! 1) \,
\big( {\calH_\h}^{\!\vee} \big)' \; = \; \Bbbk \big\langle
\overline{\tilde\x}_1, \overline{\tilde\x}_2, \dots,
\overline{\tilde\x}_n, \ldots \big\rangle \; = \;
\Bbbk \big\langle \overline\a_1, \overline\a_2, \dots,
\overline\a_n, \ldots \big\rangle  $$
(where  $ \, \overline{\bold{c}} := \bold{c} \mod (\h\!-\!1) \,
\big( {\calH_\h}^{\!\vee} \big)' \, $)  \, as  $ \Bbbk $--algebras,
and the Hopf structure is exactly the one of  $ \calH $  because it
is given by the like formulas on generators.  In a nutshell, we have
  $$  \big( {\calH_\h}^{\!\vee} \big)' \;{\buildrel {\h \rightarrow 1}
\over \llongrightarrow}\; \calH  $$
as Hopf  $ \Bbbk $--algebras.
   Therefore we got the bottom part of the diagram of deformations
(5.5), corresponding to (5.3), for  $ \, H = \calH \; (:= \calH_\nu)
\, $:  \, it is
  $$  U(\L_\nu) = {\calH_\h}^{\!\vee}{\Big|}_{\h=0}
\underset{{\calH_\h}^{\!\!\vee}} \to {\overset{0 \leftarrow \h
\rightarrow 1} \to{\longleftarrow\joinrel\llongrightarrow}}
{\calH_\h}^{\!\vee}{\Big|}_{\h=1} \! = \calH =
\big( {\calH_\h}^{\!\vee} \big)'{\Big|}_{\h=1}
\underset{({\calH_\h}^{\!\!\vee})'} \to {\overset{1 \leftarrow
\h \rightarrow 0} \to {\longleftarrow\joinrel\llongrightarrow}}
\big( {\calH_\h}^{\!\vee} \big)'{\Big|}_{\h=0} \! =
F \big[ {G_{\!\L_\nu}\phantom{|}}^{\hskip-8pt \star} \big]  $$   
or simply  $ \; U(\L_\nu) \underset{{\calH_\h}^{\!\!\vee}}
\to {\overset{0 \leftarrow \h \rightarrow 1} \to
{\longleftarrow\joinrel\relbar\joinrel\relbar\joinrel\llongrightarrow}}
\, \calH \underset{({\calH_\h}^{\!\!\vee})'}  \to
{\overset{1 \leftarrow \h \rightarrow 0} \to
{\longleftarrow\joinrel\relbar\joinrel\relbar\joinrel\llongrightarrow}}
F \big[ {G_{\!\L_\nu}\phantom{|}}^{\hskip-8pt \star} \big] \; $.
Therefore  $ \, \calH \, $   
   {\sl is intermediate}\break
 \eject   
\noindent   
 {\sl between the
(Poisson-type) ``geometrical symmetries''  $ U(\L_\nu) $  and
$ F \big[ {G_{\!\L_\nu} \phantom{|}}^{\hskip-8pt \star} \big] $},
hence the geometrical meaning of the latters should shed some light
on it; in turn, the physical meaning of  $ \calH $  should have
some reflect on the physical meaning of both  $ U(\L_\nu) $  and
$ F \big[ {G_{\!\L_\nu} \phantom{|}}^{\hskip-8pt \star} \big] $.

\vskip7pt

   {\bf 10.10 Drinfeld's algebra  $ \, {\calH_\h}^{\!\prime} :=
{\big( \calH[\h] \big)}^{\!\prime} \, $.} \, From now on we shall
deal with Drinfeld's functors in the opposite order: first  $ {(\ )}' $,
and then  $ {(\ )}^\vee $.  Like in \S 2.1, define  $ \; {\calH_\h}^{\!
\prime} \, := \, \big\{\, \eta \in \calH_\h \;\big\vert\;\, \delta_n(\eta)
\in \h^n {\calH_\h}^{\!\otimes n} \;\; \forall\; n \in \N \,\big\} \; $
$ \big( \, \subseteq \calH_\h \, \big) \, $.  We shall describe
$ {\calH_\h}^{\!\prime} $  explicitly, thus checking that it is
really a QFA, as predicted by  Theorem 2.2{\it (a)\/};  then we'll
look at its specialization at  $ \, \h = 0 \, $  and at  $ \, \h
= 1 \, $,  and finally we'll study  $ \big( {\calH_\h}^{\!\prime}
\big)^{\!\vee} $  and its specializations at  $ \, \h = 0 \, $  and
$ \, \h = 1 \, $.  The outcome will be an explicit description of
(5.4) for  $ \, H = \calH \, $  ($ \, = \calH_\nu \, $,  with
$ \, \nu \in \N \cup \{\infty\} \, $  fixed as before).
                                               \par
   Let  $ \, \underline{D} := \underline{D}(\calH) = {\big\{ D_n
\big\}}_{n \in \N} = {\big\{ \text{\sl Ker}\,\big(\delta_n \, \colon
\, \calH \longrightarrow (\calH^{\text{dif}})^{\otimes n} \big)
\big\}}_{n \in \N_\nu} \, $  be the Hopf algebra filtration of
$ \calH $  as considered in \S 5.1.  Then by Lemma 5.2, we have   
 \vskip-11pt   
  $$  {\calH_\h}^{\!\prime} \; = \; \Cal{R}^\h_{\underline{D}}
(\calH) \; := \; \Bbbk[\h] \cdot D_0 + \h \; \Bbbk[\h] \cdot D_1
+ \cdots + \h^n \, \Bbbk[\h] \cdot D_n + \cdots  $$
 \vskip-7pt   
\noindent   
 so we only need to compute the filtration  $ \underline{D} \, $.
The idea is to describe it in combinatorial terms, based on the
non-commutative polynomial nature of  $ \calH \, $.
                                               \par
   As before, we proceed in steps.

\vskip7pt

{\bf 10.11 Gradings and filtrations:} \, Let  $ \, \partial_- $
be the unique Lie algebra grading of  $ \L_\nu $  given by  $ \;
\partial_-(\boldalpha_n) := n - 1 + \delta_{n,1} \; $  (for all
$ \, n \in \N_\nu \, $).  Let also  $ d $  be the standard Lie algebra
grading associated with the central lower series of  $ \L_\nu \, $:
in down-to-earth terms,  $ d $  is defined by  $ \; d \big( [ \cdots
[[x_{s_1}, x_{s_2}], \dots x_{s_k}] \big) = k - 1 \; $  on any Lie
monomial of  $ \L_\nu \, $.  Since both  $ \partial_- $  and  $ d $
are Lie algebra gradings, their difference  $ \, (\partial_- - d) \, $
is a Lie algebra grading too.  Let
%
%
%
%
$ \, {\big\{ F_n \big\}}_{n \in \N} \, $  be the Lie algebra
filtration associated with the grading  $ (\partial_- - d \,) $;
then the down-shifted filtration  $ \, \underline{T} := {\big\{\,
T_n := F_{n-1} \,\big\}}_{n \in \N} \, $  is again a Lie algebra
filtration of  $ \L_\nu \, $.  There is a unique algebra filtration
on  $ U(\L_\nu) $  extending  $ \underline{T} $,  which we denote by
$ \, \underline{\varTheta} = {\big\{ \varTheta_n \big\}}_{n \in \N}
\, $;  \, as a matter of notation, we set also  $ \, \varTheta_{-1}
:= \{0\} \, $.  Finally, for each  $ \, y \in U(\L_\nu) \setminus
\{0\} \, $  there is a unique  $ \, \tau(y) \in \N \, $  with  $ \,
y \in \varTheta_{\tau(y)} \setminus \varTheta_{\tau(y)-1} \, $;  \,
in particular, we have  $ \; \tau(b) = \partial_-(b) - d(b) \, $,
$ \; \tau(b \, b') = \tau(b) + \tau(b') \; $  and  $ \; \tau \big(
[b,b'] \big) = \tau(b) + \tau(b') - 1 \; $  for all  $ \, b, b'
\in B_\nu \, $.
                                               \par
   We can explicitly describe  $ \underline{\varTheta} $.  Indeed, let
us fix any total order  $ \preceq $  on the basis  $ B_\nu $  of \S
10.2: then  $ \; \calX := \Big\{\, \underline{b} := b_1 \cdots b_k
\,\Big|\; k \in \N \, , \; b_1, \dots, b_k \in B_\nu \, , \; b_1
\preceq \cdots \preceq b_k \,\Big\} \; $  is a  $ \Bbbk $--basis  of
$ U(\L_\nu) $,  by the PBW theorem.  It follows that  $ \varTheta $
induces a set-theoretic filtration  $ \, \underline{\calX} = {\big\{
\calX_n \big\}}_{n \in \N} \, $  of  $ \calX $  with  $ \, \calX_n :=
\calX \cap \varTheta_n = \Big\{\, \underline{b} := b_1 \cdots b_k
\,\Big|\; k \in \N \, , \; b_1, \dots, b_k \in B_\nu \, , \; b_1
\preceq \cdots \preceq b_k \, , \; \tau(\underline{b}\,) = \tau(b_1)
+ \cdots + \tau(b_k) \leq n \,\Big\} \, $,  \; and also that  $ \;
\varTheta_n = \text{\sl Span}\,\big(\calX_n\big) \; $  for all
$ \, n \in \N \, $.
                                               \par
   Let us define  $ \, \boldalpha_1 := \a_1 \, $  and  $ \,
\boldalpha_n := \a_n - {\a_1}^{\!n} \, $  for all  $ \, n \in
\N_\nu \setminus \{1\} \, $.  This ``change of variables''
--- which switch from the  $ \a_n $'s  to their ``differentials'',
in a sense ---   will be the key to achieve a complete description
of the filtration  $ \underline{D} \, $;  in turn, this will pass
through a close comparison among  $ \calH $  and  $ U(\L_\nu) \, $.
%
%
 \eject   
   By definition  $ \, \calH = \calH_\nu \, $  is the free associative
algebra over  $ \{\a_n\}_{n \in \N_\nu} $,  hence   --- by definition
of the  $ \boldalpha $'s  ---   also over  $ \{\boldalpha_n\}_{n \in
\N_\nu} $;  so we have an algebra isomorphism  $ \; \Phi \, \colon
\, \calH \,{\buildrel \cong \over
{\lhook\joinrel\relbar\joinrel\twoheadrightarrow}}\,
U(\L_\nu) \; $  given by  $ \, \boldalpha_n \mapsto x_n \, $
($ \, \forall \; n \in \N_\nu \, $).  Via  $ \Phi $  we pull back
all data and results about gradings, filtrations, PBW bases and so
on mentioned above for  $ U(\L_\nu) \, $;  in particular we set  $ \,
\boldalpha_{\underline{b}} := \Phi(x_{\underline{b}}) = \boldalpha_{b_1}
\cdots \boldalpha_{b_k} \, $  (for all  $ \, b_1, \dots, b_k \in B_\nu
\, $),  $ \, \calA_n := \Phi(\calX_n) \, $  (for all  $ \, n \in \N
\, $)  and  $ \, \calA := \Phi(\calX) = \bigcup_{n \in \N} \calA_n
\, $.   For gradings on  $ \calH $  we stick to the like notation,
i.e.~$ \partial_- $,  $ d $  and  $ \tau \, $,  \, and similarly
for the filtration  $ \underline{\varTheta} \, $.
                                                    \par
   Finally, for all  $ \, a \in \calH \setminus \{0\} \, $  \,
we set also  $ \, \kappa\,(a) := k \, $  iff  $ \, a \in D_k
\setminus D_{k-1} \, $  (with  $ \, D_{-1} := \{0\} \, $).
                                                    \par
   {\sl Our goal is to prove an identity of filtrations, namely  $ \,
\underline{D} = \underline{\varTheta} \, $},  \, or equivalently
$ \, \kappa = \tau \, $.  In fact, this would give to the Hopf
filtration  $ \underline{D} $,  which is defined intrinsically
in Hopf algebraic terms, an explicit  {\sl combinatorial\/}
description, namely the one of  $ \varTheta $  explained above.

\vskip7pt

\proclaim{Lemma 10.12} For all  $ \, \ell, t \in \N \, $,
$ \, t \geq 1 \, $,  \, we have (notation of \S 10.11)
  $$  Z^\ell_t(\boldalpha_*) := \Big( Q^\ell_t(\a_*) -
{\textstyle {{\ell + t} \choose t}} \, {\a_1}^{\!t} \Big)
\in \varTheta_{t-1}  \quad  \text{and}  \quad  Q^\ell_t(\a_*)
\in \varTheta_t \setminus \varTheta_{t-1} \; .  $$
\endproclaim

\demo{Proof} \,  When  $ \, t = 1 \, $  definitions give
$ \, Q^\ell_1(\a_*) = (\ell+1) \, \a_1 \in \varTheta_1 \, $  and
so  $ \, Z^\ell_1(\boldalpha_*) = (\ell+1) \, \a_1 - {{\ell + 1}
\choose 1} \, \a_1 = 0 \in \varTheta_0 \, $,  \, for all  $ \,
\ell \in \N \, $.  Similarly, when  $ \, \ell = 0 \, $  we have
$ \, Q^0_t(\a_*) = \a_t \in \varTheta_t \, $  and so  $ \,
Z^0_t(\boldalpha_*) = \a_t - {1 \choose 1} \, {\a_1}^{\!t}
= \boldalpha_t \in \varTheta_{t-1} \, $  (by definition),
\, for all  $ \, t \in \N_+ \, $.
                                                  \par
   When  $ \, \ell > 0 \, $  and  $ \, t > 1 \, $,  \, we can prove
the claim using two independent methods.
                                                  \par
   {\it $ \underline{\text{First method}} $:} \,  The very
definitions imply that the following recurrence formula holds:
  $$  Q^\ell_t(\a_*) \, = \, Q^{\ell-1}_t(\a_*) \, + \,
\sum\nolimits_{s=1}^{t-1} Q^{\ell-1}_{t-s}(\a_*) \, \a_s +
\a_t  \qquad  \forall \quad \ell \geq 1 \, ,  \; t \geq 2 \, .  $$
From this formula we argue
  $$  \displaylines{
   Z^\ell_t(\boldalpha_*) \; := \; Q^\ell_t(\a_*) \, - \,
{\textstyle {{\ell+t} \choose t}} \, {\a_1}^{\!t} \; =
\; Q^{\ell-1}_t(\a_*) \, + \, {\textstyle \sum_{s=1}^{t-1}}
Q^{\ell-1}_{t-s}(\a_*) \, \a_s \, + \, \a_t \, - \, {\textstyle
{{\ell+t} \choose t}} \, {\a_1}^{\!t} \; =   \hfill  \cr
   {} \hfill   = \; Z^{\ell-1}_t(\a_*) \, + \, {\textstyle {{\ell-1+t}
\choose t}} \, {\a_1}^{\!t} \, + \, {\textstyle \sum_{s=1}^{t-1}}
\left( Z^{\ell-1}_{t-s}(\a_*) \, + \, {\textstyle {{\ell-1+t-s}
\choose {t-s}}} \, {\a_1}^{\!t-s} \right) \, \a_s + \a_t -
{\textstyle {{\ell+t} \choose t}} \, {\a_1}^{\!t} \; =  \cr
   = \; Z^{\ell-1}_t(\a_*) \, + \, {\textstyle {{\ell-1+t} \choose t}}
\, {\boldalpha_1}^{\!t} \, + \, {\textstyle \sum_{s=1}^{t-1}}
Z^{\ell-1}_{t-s}(\a_*) \, \big( \boldalpha_s + {\boldalpha_1}^{\!s}
\big) \, + \,   \hfill  \cr
   {} \hfill   + \, {\textstyle \sum_{s=1}^{t-1}} {\textstyle
{{\ell-1+t-s} \choose {t-s}}} \, {\boldalpha_1}^{\!t-s} \,
\big( \boldalpha_s + {\boldalpha_1}^{\!s} \big) \, + \,
\big( \boldalpha_t + {\boldalpha_1}^{\!t} \big) \, - \,
{\textstyle {{\ell+t} \choose t}} \, {\boldalpha_1}^{\!t} \; =  \cr
   = \; Z^{\ell-1}_t(\a_*) \, + \, {\textstyle \sum_{s=1}^{t-1}}
Z^{\ell-1}_{t-s}(\a_*) \, \big( \boldalpha_s + {\boldalpha_1}^{\!s}
\big) \, + \, {\textstyle \sum_{s=1}^{t-1}} {\textstyle
{{\ell-1+t-s} \choose {t-s}}} \, {\boldalpha_1}^{\!t-s} \,
\boldalpha_s \, + \, \boldalpha_t \, + \,   \hfill  \cr
   {} \hfill   + \, {\textstyle \sum_{s=1}^{t-1}} {\textstyle
{{\ell-1+t-s} \choose {t-s}}} \, {\boldalpha_1}^{\!t-s} \,
{\boldalpha_1}^{\!s} \, + \, {\boldalpha_1}^{\!t} \,
+ \, {\textstyle {{\ell-1+t} \choose t}}
\, {\boldalpha_1}^{\!t} \, - \,
{\textstyle {{\ell+t} \choose t}} \, {\boldalpha_1}^{\!t} \; =  \cr
   = \; Z^{\ell-1}_t(\a_*) \, + \, {\textstyle \sum_{s=1}^{t-1}}
Z^{\ell-1}_{t-s}(\a_*) \, \big( \boldalpha_s + {\boldalpha_1}^{\!s}
\big) \, + \, {\textstyle \sum_{s=1}^{t-1}} {\textstyle
{{\ell-1+t-s} \choose {t-s}}} \, {\boldalpha_1}^{\!t-s} \,
\boldalpha_s \, + \, \boldalpha_t \, + \,   \hfill  \cr
   {} \hfill   + \, \bigg( {\textstyle \sum_{r=0}^t} {\textstyle
{{\ell-1+r} \choose {\ell-1}}} \, - \, {\textstyle {{\ell+t}
\choose t}} \bigg) \, {\boldalpha_1}^{\!t} \; =  \cr
   {} \hfill   = \; Z^{\ell-1}_t(\a_*) \, + \, {\textstyle
\sum_{s=1}^{t-1}} Z^{\ell-1}_{t-s}(\a_*) \, \big( \boldalpha_s
+ {\boldalpha_1}^{\!s} \big) \, + \, {\textstyle \sum_{s=1}^{t-1}}
{\textstyle {{\ell-1+t-s} \choose {t-s}}} \, {\boldalpha_1}^{\!t-s}
\, \boldalpha_s \, + \, \boldalpha_t  \cr }  $$
because of the classical identity  $ \; {{\ell+t} \choose \ell} =
\sum_{r=0}^t {{\ell-1+r} \choose {\ell-1}} \; $.  Then induction
upon  $ \ell \, $  and the very definitions allow to conclude that
all summands in the final sum belong to  $ \varTheta_{t-1} $,
hence  $ \, Z^\ell_t(\boldalpha_*) \in \varTheta_{t-1} \, $  as
well.  Finally, this implies  $ \; Q^\ell_t(\a_*) \, = \, Z^\ell_t
(\boldalpha_*) \, + \, {{\ell+t} \choose t} \, {\boldalpha_1}^{\!t}
\in \varTheta_t \setminus \varTheta_{t-1} \, $.
                                                  \par
   {\it $ \underline{\text{Second method}} $:}  $ \; Q_t^\ell(\a_*)
:= \sum_{s=1}^t {{\ell + 1} \choose s} \, P_t^{(s)}(\a_*) =
\sum_{s=1}^t {{\ell + 1} \choose s} \sum_{\hskip-4pt  \Sb  j_1,
\dots, j_s > 0  \\   j_1 + \cdots + j_s \, = \, t  \endSb}
\hskip-5pt  \a_{j_1} \cdots \a_{j_s} \, $,  \;  by definition;
then expanding the  $ \a_j $'s  as  $ \, \a_1 = \boldalpha_1 \, $
and  $ \, \a_j = \boldalpha_j + {\boldalpha_1}^{\!j} \, $  (for
$ \, j > 1 \, $)  we find that  $ \, Q_t^\ell(\a_*) = Q_t^\ell \big(
\boldalpha_* + {\boldalpha_1}^{\!*} \big) \, $  is a linear combination
of monomials  $ \, \boldalpha_{(j_1)} \cdots \boldalpha_{(j_s)} \, $
with  $ \, j_1, \dots, j_s > 0 \, $,  $ \, j_1 + \cdots + j_s =
t \, $,  $ \, \boldalpha_{(j_r)} \in \big\{ \boldalpha_{j_r},
{\boldalpha_1}^{\!j_r} \big\} \, $  for all  $ r \, $.  Let
$ Q_- $  be the linear combination of those monomials such
that  $ \, (\boldalpha_{(j_1)}, \boldalpha_{(j_2)}, \dots,
\boldalpha_{(j_s)}\big) \not= \big( {\boldalpha_1}^{\!j_1},
{\boldalpha_1}^{\!j_2}, \dots, {\boldalpha_1}^{\!j_s}\big) \, $;
\, for the remaining monomials we have  $ \, \boldalpha_{j_1} \cdot
\boldalpha_{j_2} \cdots \boldalpha_{j_s} = {\boldalpha_1}^{\! j_1 +
\cdots + j_s} = {\boldalpha_1}^{\!t} \, $,  \, hence their linear
combination giving  $ \, Q_+ := Q^\ell_t(\a_*) - Q_- \, $  is a
multiple of  $ {\boldalpha_1}^{\!t} $,  say  $ \, Q_+ = N \,
{\boldalpha_1}^{\!t} \, $.
                                             \par
   Now we compute this coefficient  $ N \, $.  First, by construction
$ N $  is nothing but  $ \, N = Q^\ell_t(1_*) = Q^\ell_t(1, 1, \dots,
1, \dots) \, $  where the latter means the (positive integer) value
of the polynomial  $ Q^\ell_t $  when all its indeterminates are
set equal to  $ 1 $.  Thus we compute  $ Q^\ell_t(1_*) \, $.
                                             \par
   Recall that the  $ Q^\ell_t $'s  enter in the definition of the
coproduct of  $ F \big[ \G^\dif \big] $:  \, the latter is dual
to the (composition) product of series in  $ \G^\dif $,  thus
if  $ \{a_n\}_{n \in \N_+} $  and  $ \{b_n\}_{n \in \N_+} $
are two countable sets of commutative indeterminates then
  $$  \displaylines{
   \Big( \, x \, + {\textstyle \sum_{n=1}^{+\infty}} \, a_n \,
x^{n+1} \Big) \circ \Big( \, x \, + {\textstyle \sum_{m=1}^{+\infty}}
\, b_m \, x^{m+1} \Big) \; :=   \hfill  \cr
   {} \hfill   := \; \bigg( \! \Big( \, x \, +
{\textstyle \sum\limits_{m=1}^{+\infty}} \, b_m \, x^{m+1} \Big)
+ {\textstyle \sum\limits_{n=1}^{+\infty}} \, a_n \, \Big( \, x \, +
{\textstyle \sum\limits_{m=1}^{+\infty}} \, b_m \, x^{m+1} \Big)^{n+1}
\bigg) \; = \; x \, + {\textstyle \sum\limits_{k=0}^{+\infty}} \,
c_k \, x^{k+1}  \cr }  $$
with  $ \; c_k = Q^0_k(b_*) + \sum_{r=1}^k a_r \cdot Q^r_{k-r}(b_*)
\, $  (cf.~\S 10.2).  Specializing  $ \, a_\ell = 1 \, $  and  $ \,
a_r = 0 \, $  for all  $ \, r \not= \ell \, $  we get  $ \, c_{t+\ell}
= Q^0_{t+\ell}(b_*) + Q^\ell_t(b_*) = b_{t+\ell} + Q^\ell_t(b_*) \, $.
In particular setting  $ \, b_* = 1_* \, $  we have that  $ \, 1 +
Q^\ell_t(1_*) \, $  is the coefficient  $ \, c_{\ell+t} \, $  of
$ \, x^{\ell + t + 1} \, $  in the series
  $$  \displaylines{
   \big( \, x \, + x^{\ell+1} \big) \circ \Big( \, x \, +
{\textstyle \sum_{m=1}^{+\infty}} \, x^{m+1} \Big) \; =   \hfill  \cr
   {} \hfill   = \; \big( \, x \, + x^{\ell+1} \big) \circ
\big( x \cdot {(1-x)}^{-1} \big) \; = \; x \cdot {(1-x)}^{-1}
+ {\big( x \cdot {(1-x)}^{-1} \big)}^{\ell+1} = \;  \cr
   {} \quad   = \; {\textstyle \sum_{m=0}^{+\infty}} \, x^{m+1}
\, + \, x^{\ell+1} \Big( {\textstyle \sum_{m=0}^{+\infty}} \, x^m
\Big)^{\ell+1} \; = \; {\textstyle \sum_{m=0}^{+\infty}} \, x^{m+1}
\, + \, x^{\ell+1} \, {\textstyle \sum_{n=0}^{+\infty} {{\ell+n}
\choose \ell}} \, x^n \; =   \hfill  \cr
   {} \hfill   = \; {\textstyle \sum_{s=0}^{\ell-1}} \, x^{s+1} \, +
\, {\textstyle \sum_{s=\ell}^{+\infty} \left( 1 + {s \choose \ell}
\right)} \, x^{s+1}  \quad ;  \cr }  $$
therefore  $ \, 1 + Q^\ell_t(1_*) = c_{\ell+t} = 1 + {{\ell+t} \choose
\ell} \, $,  \, whence  $ \, Q^\ell_t(1_*) = {{\ell+t} \choose
\ell} \, $.  As an alternative approach, one can prove that
$ \, Q^\ell_t(1_*) = {{\ell+t} \choose \ell} \, $  by induction
using the recurrence formula  $ \; Q^\ell_t(\x_*) \, = \,
Q^{\ell-1}_t(\x_*) + \sum_{s=1}^{t-1} Q^{\ell-1}_{t-s}(\x_*)
\, \x_s + \x_t \; $  and the identity  $ \; {{\ell+t} \choose
\ell} = \sum_{s=0}^t {{\ell+t-1} \choose {\ell-1}} \; $.
                                               \par
   The outcome is  $ \, N = Q^\ell_t(1_*) = {{\ell+t} \choose \ell}
\, $  (for all  $ t, \ell \, $),  thus  $ \; Q^\ell_t(\a_*) - {{\ell+t}
\choose \ell} \, \a_t \, = \, Q_- + Q_+ - {{\ell+t} \choose \ell} \,
\a_t \, = \, Q_- + N \, \a_t - {{\ell+t} \choose \ell} \, \a_t \, = \,
Q_- \, $.  Now, by definition  $ \, \tau(\boldalpha_{j_r}) = j_r - 1
\, $  and  $ \, \tau\big({\boldalpha_1}^{\!j_r} \big) = j_r \, $.
Therefore if  $ \, \boldalpha_{(j_r)} \in \big\{ \boldalpha_{j_r},
{\boldalpha_1}^{\!j_r} \big\} \, $  (for all  $ \, r = 1, \dots, s
\, $)  and  $ \, (\boldalpha_{(j_1)}, \boldalpha_{(j_2)}, \dots,
\boldalpha_{(j_s)}) \not= \big( {\boldalpha_1}^{\!j_1},
{\boldalpha_1}^{\!j_2}, \dots, {\boldalpha_1}^{\!j_s} \big) \, $,
\, then  $ \, \tau \big( \boldalpha_{(j_1)} \cdots \boldalpha_{(j_s)}
\big) \leq j_1 + \cdots + j_s - 1 = t - 1 \, $.  Then by construction
$ \, \tau(Q_-) \leq t - 1 \, $,  \, whence, since  $ \,
Z^\ell_t(\boldalpha_*) := Q^\ell_t(\a_*) - {{\ell+t} \choose \ell}
\, \a_t = Q_- \, $,  \, we get also  $ \, \tau \big( Z^\ell_t
(\boldalpha_*) \big) \leq t - 1 \, $,  \, i.e.~$ \, Z^\ell_t
(\boldalpha_*) \in \varTheta_{t-1} \, $,  \, so  $ \;
Q^\ell_t(\a_*) \, = \, Z^\ell_t(\boldalpha_*) \, +
\, {{\ell+t} \choose t} \, {\boldalpha_1}^{\!t} \in
\varTheta_t \setminus \varTheta_{t-1} \, $.   \qed
\enddemo

\vskip7pt

\proclaim{Proposition 10.13}  $ \, \underline{\varTheta} $  is
a Hopf algebra filtration of  $ \, \calH \, $.
\endproclaim

\demo{Proof} \, By construction (cf.~\S 10.11)
$ \underline{\varTheta} $  is an algebra filtration; so to
check it is  {\sl Hopf\/}  too we are left only to show that
$ \; (\star) \, \Delta(\varTheta_n) \subseteq \sum_{r+s=n}
\varTheta_r \otimes \varTheta_s \; $  (for
all  $ \, n \in \N \, $),  \, for then  $ \, S(\varTheta_n) \subseteq
\varTheta_n \, $  (for all  $ \, n \, $)  will follow from that by
recurrence (and Hopf algebra axioms).
                                                  \par
   By definition  $ \, \varTheta_0 = \Bbbk \cdot 1_{\scriptscriptstyle
\calH} \, $;  then  $ \, \Delta(1_{\scriptscriptstyle \calH}) =
1_{\scriptscriptstyle \calH} \otimes 1_{\scriptscriptstyle \calH} \, $
proves  $ (\star) $  for  $ \, n = 0 \, $.  For  $ \, n = 1 \, $
definitions tell that  $ \varTheta_1 $  is nothing but the direct
sum of  $ \varTheta_0 $  with the (free) Lie (sub)algebra (of
$ \calH \, $)  generated by  $ \{\boldalpha_1,\boldalpha_2\} $.
Since  $ \, \Delta(\boldalpha_1) = \boldalpha_1 \otimes 1 + 1 \otimes
\boldalpha_1 \, $  and  $ \, \Delta(\boldalpha_2) = \boldalpha_2
\otimes 1 + 1 \otimes \boldalpha_2  \, $  (directly from definitions)
and since
  $$  \Delta\big([x,y]\big) \, = \, \big[\Delta(x),\Delta(y)\big]
\, = \, {\textstyle \sum_{(x),(y)}} \big( [x_{(1)},y_{(1)}] \otimes
x_{(2)} y_{(2)} + x_{(1)} y_{(1)} \otimes [x_{(2)},y_{(2)}] \big)  $$
(for all  $ \, x, y \in \calH \, $)  we argue that  $ (\star) $
holds for  $ \, n = 1 \, $  too.
                                                \par
   Further on, for every  $ \, n > 1 \, $  we have (setting
$ \, Q^n_0(\a_*) = 1 = \a_0 \, $  for short)
  $$  \displaylines{
   \Delta(\boldalpha_n) \; = \; \Delta(\a_n) - \Delta \big(
{\a_1}^{\!n} \big) \; = \; {\textstyle \sum_{k=0}^n} \, \a_k
\otimes Q^k_{n-k}(\a_*) - {\textstyle \sum_{k=0}^n \,
{n \choose k}} \, {\a_1}^k \otimes {\a_1}^{\!n-k} \;
= \;   \hfill  \cr
  {} \hfill   = \; {\textstyle \sum_{k=2}^n} \, \boldalpha_k
\otimes Q^k_{n-k}(\a_*) \, + \, {\textstyle \sum_{k=0}^{n-1}} \,
{\boldalpha_1}^{\!k} \otimes Z^k_{n-k}(\boldalpha_*)  \cr }  $$
hence  $ \, \Delta(\boldalpha_n) \in \sum_{r+s=n-1} \varTheta_r
\otimes \varTheta_s \, $  due to Lemma 10.12 (and to $ \,
\boldalpha_m \in \varTheta_{m-1} \, $  for  $ \, m > 1 \, $).
                                                   \par
   Finally, as  $ \, \Delta\big([x,y]\big) \! = \! \big[ \Delta(x),
\Delta(y) \big] \! = \! \sum_{(x),(y)} \! \big( [x_{(1)},y_{(1)}] \otimes
x_{(2)} y_{(2)} + x_{(1)} y_{(1)} \otimes [x_{(2)},y_{(2)}] \big) $
and similarly  $ \, \Delta(x\,y) = \Delta(x) \Delta(y) = \sum_{(x),(y)}
x_{(1)} y_{(1)} \otimes x_{(2)} y_{(2)} \, $  (for  $ \, x, y \in \calH
\, $),  we have that  $ \Delta $  does not increase  $ \, (\partial_-
- d\,) \, $:  \, as  $ \varTheta $  is exactly the (algebra) filtration
induced by  $ (\partial_- - d\,) \, $,  \, it is a Hopf algebra
filtration as well.   \qed
\enddemo

\vskip7pt

\proclaim{Lemma 10.14}  (notation of \S 10.11)
                                    \hfill\break
   \indent   (a)  $ \; \kappa\,(a) \leq \partial(a) \; $  for
every  $ \, a \in \calH \setminus \{0\} \, $  which is
$ \, \partial(a) $--homogeneous.
                                    \hfill\break
   \indent   (b)  $ \; \kappa\,(a\,a') \leq \kappa\,(a) + \kappa\,(a')
\; $  and  $ \; \kappa\,\big([a,a']\big) < \kappa\,(a) + \kappa\,(a')
\; $  for all  $ \, a, a' \in \calH \setminus \{0\} \, $.
                                    \hfill\break
   \indent   (c)  $ \; \kappa\,(\boldalpha_n) = \partial_-(\boldalpha_n)
= \tau(\boldalpha_n) \; $  for all  $ \, n \in \N_\nu \, $.
                                    \hfill\break
   \indent   (d)  $ \; \kappa \, \big( [\boldalpha_r,\boldalpha_s] \big)
= \partial_-(\boldalpha_r) + \partial_-(\boldalpha_s) - 1 = \tau \big(
[\boldalpha_r, \boldalpha_s] \big) \; $  for all  $ \, r $,  $ s \in
\N_\nu \, $  with  $ \, r \not= s \, $.
                                    \hfill\break
   \indent   (e)  $ \; \kappa\,(\boldalpha_b) = \partial_-(\boldalpha_b)
- d(\boldalpha_b) + 1 = \tau(\boldalpha_b) \; $  for every  $ \, b
%
%
 \in B_\nu \, $.
%
%
                                    \hfill\break
   \indent   (f)  $ \; \kappa\,(\boldalpha_{b_1} \boldalpha_{b_2}
\cdots \boldalpha_{b_\ell}) = \tau(\boldalpha_{b_1} \boldalpha_{b_2}
\cdots \boldalpha_{b_\ell}) \; $  for all  $ \, b_1, b_2, \dots,
b_\ell \in B_\nu \, $.
                                       \hfill\break
   \indent   (g)  $ \; \kappa \big([\boldalpha_{b_1},
\boldalpha_{b_2}]\big) = \kappa\,(\boldalpha_{b_1})
+ \kappa\,(\boldalpha_{b_2}) - 1 = \tau \big([\boldalpha_{b_1},
\boldalpha_{b_2}]\big) \, $,  \; for all  $ \, b_1, b_2 \in
B_\nu \, $.
\endproclaim

\demo{Proof} {\it (a)} \, Let  $ \, a \in \calH \setminus \{0\} \, $
be  $ \partial(a) $--homogeneous.  Since  $ \calH $  is graded, we
have  $ \, \partial \big( \delta_\ell(a) \big) = \partial(a) \, $  for
all  $ \ell \, $;  \, moreover,  $ \, \delta_\ell(a) \in \! J^{\otimes
\ell} \, $  (with  $ \, J := \hbox{\sl Ker}\,(\epsilon_{\calH}) \, $)
by definition, and  $ \, \partial(y) > 0 \, $  for each
$ \partial $--homogeneous  $ \, y \in J \setminus \{0\} \, $.
Then  $ \, \delta_\ell(a) = 0 \, $  for all  $ \, \ell >
\partial(a) \, $,  \, whence the claim.
                                             \par
   {\it (b)} \, This is just a reformulation of  Lemma 3.4{\it (c).}
                                                 \par
   {\it (c)} \, By part  {\it (a)\/}  we have  $ \, \kappa(\a_n)
\leq \partial(\a_n) = n \, $.  Moreover, by definition  $ \,
\delta_2(\a_n) = \sum_{k=1}^{n-1} \a_k \otimes Q^k_{n-k}(\a_*)
\, $,  thus  $ \, \delta_n(\a_n) = (\delta_{n-1} \otimes \delta_1)
\big( \delta_2(\a_n) \big) = \sum_{k=1}^{n-1} \delta_{n-1}(\a_k)
\otimes \delta_1 \big( Q^k_{n-k} (\a_*) \big) \, $  by coassociativity.
Since  $ \, \delta_\ell(\a_m) = 0 \, $  for  $ \, \ell > m \, $,  $ \,
Q^{n-1}_1(\a_*) = n \, \a_1 \, $  and  $ \, \delta_1(\a_1) = \a_1 \, $,
\, we have  $ \, \delta_n(\a_n) = \delta_{n-1}(\a_{n-1}) \otimes (n
\, \a_1) \, $,  \, thus by induction  $ \; \delta_n(\a_n) = n! \,
{\a_1}^{\!\otimes n} \, $  ($ \, \not= 0 \, $),  \, whence  $ \,
\kappa(\a_n) = n \, $.  But also  $ \, \delta_n({\a_1}^{\!n}) = n!
\, {\a_1}^{\!\otimes n} \, $.  Thus  $ \, \delta_n(\boldalpha_n) =
\delta_n(\a_n) - \delta_n({\a_1}^{\!n}) = 0 \, $  for  $ \, n > 1 \, $.
                                                  \par
   Clearly  $ \, \kappa(\boldalpha_1) = 1 \, $.  For the general case,
for all  $ \, \ell \geq 2 \, $  we have
  $$ \, \delta_{\ell-1}(\a_\ell) \; = \; (\delta_{\ell-2}
\otimes \delta_1) \big( \delta_2(\a_\ell) \big) \; = \;
\sum\nolimits_{k=1}^{\ell-1} \delta_{\ell-2}(\a_k) \otimes
\delta_1 \big( Q^k_{\ell-1-k}(\a_*) \big)  $$
which,
%
%
thanks to the previous analysis, gives
  $$  \displaylines{
   \quad  \delta_{\ell-1}(\a_\ell) \; = \; \delta_{\ell-2}(\a_{\ell-2})
\otimes \left( {(\ell-1)} \, \a_2 + {{\ell-1} \choose 2} \, {\a_1}^{\!2}
\right) \, + \, \delta_{\ell-2} (\a_{\ell-1}) \otimes \ell \, \a_1 \; =
\;   \hfill  \cr
   {} \hfill   = \; {(\ell-1)}! \cdot {\a_1}^{\!\otimes (\ell-2)}
\otimes \left( \a_2 + {{\;\ell-1\,} \over {\,2\,}} \cdot {\a_1}^{\!2}
\right) \, + \, \ell \cdot \delta_{\ell-2} (\a_{\ell-1}) \otimes \a_1
\; .  \quad  \cr }  $$
Iterating we get, for all  $ \, \ell \geq 2 \, $  (with
$ \; {{-1} \choose 2} := 0 \, $,  \, and changing indices)
  $$  \delta_{\ell-1}(\a_\ell) \; = \; \sum\nolimits_{m=1}^{\ell-1}
{{\,\ell\,!\,} \over {m+1}} \cdot {\a_1}^{\!\otimes (m-1)} \otimes
\left( \a_2 + {{\,m-1\,} \over {\,2\,}} \cdot {\a_1}^{\!2} \right)
\otimes {\a_1}^{\!\otimes (\ell-1-m)} \; .  $$   
   On the other hand, we have also
  $$  \delta_{\ell-1}\big({\a_1}^{\!\ell}\big) \; = \;
\sum\nolimits_{m=1}^{\ell-1} {{\,\ell\,!\,} \over {\,2\,}} \cdot
{\a_1}^{\!\otimes (m-1)} \otimes {\a_1}^{\!2} \otimes {\a_1}^{\!
\otimes (\ell-1-m)} \; .  $$
Therefore, for  $ \, \delta_{n-1}(\boldalpha_n) = \delta_{n-1}(\a_n)
- \delta_{n-1}({\a_1}^{\!n}) \, $  (for all  $ \, n \in \N_\nu \, $,
$ \, n \geq 2 \, $) the outcome is
  $$  \hbox{ $ \eqalign{
   \delta_{n-1}(\boldalpha_n) \;  &  = \; \sum\nolimits_{m=1}^{n-1}
{{\,n!\,} \over {\,m+1\,}} \cdot {\a_1}^{\!\otimes (m-1)} \otimes
\big( \a_2 - {\a_1}^{\!2} \big) \otimes {\a_1}^{\!\otimes (n-1-m)}
\; =  \cr
   {}  &  = \; \sum\nolimits_{m=1}^{n-1} {{\,n!\,} \over {\,m+1\,}}
\cdot {\boldalpha_1}^{\!\otimes (m-1)} \otimes \boldalpha_2
\otimes {\boldalpha_1}^{\!\otimes (n-1-m)}  \quad ;  \cr } $ }
\eqno (10.5)  $$
in particular  $ \, \delta_{n-1}(\boldalpha_n) \not= 0 \, $,  \,
whence  $ \, \boldalpha_n \not\in D_{n-2} \, $  and so
$ \, \kappa(\boldalpha_n) = n-1 \, $,  \, q.e.d.
                                                  \par
   {\it (d)} \, Let  $ \, r \not= 1 \not= s \, $.  From
{\it (b)--(c)\/} we get  $ \, \kappa \big([\boldalpha_r,
\boldalpha_s]\big) < \kappa(\boldalpha_r) + \kappa(\boldalpha_s)
= r+s-2 \, $.  In addition, we prove now that  $ \; \delta_{r+s-3}
\big([\boldalpha_r,\boldalpha_s]\big) \not= 0 \, $,  \, which
yields  {\it (d)}.  Lemma 3.4{\it(b)\/}  gives
  $$  \displaylines{
   \delta_{r+s-3}\big([\boldalpha_r,\boldalpha_s]\big) \; =
\hskip-9pt  \sum_{\Sb \Lambda \cup Y = \{1,\dots,r+s-3\}  \\
\Lambda \cap Y \not= \emptyset  \endSb}  \hskip-26pt
\big[ \delta_\Lambda(\boldalpha_r), \delta_Y(\boldalpha_s) \big]
\; =  \hskip-21pt
\sum_{\Sb \Lambda \cup Y = \{1,\dots,r+s-3\} \\
\Lambda \cap Y \not= \emptyset, \, |\Lambda| = r-1, \,
|Y| = s-1 \endSb}  \hskip-38pt  \big[ j_\Lambda \big( \delta_{r-1}
(\boldalpha_r) \big), j_Y \big( \delta_{s-1}(\boldalpha_s) \big)
\big] \, .  \cr }  $$
%
%
%
%
Using (10.5) in the form  $ \; \delta_{\ell-1}(\a_\ell) =
\sum_{m=1}^{\ell-1} {{\,\ell\,!\,} \over 2} \cdot \boldalpha_2
\otimes {\boldalpha_1}^{\!\otimes (\ell-2)} + \boldalpha_1
\otimes \eta_\ell \; $  (for some  $ \, \eta_\ell \in \calH \, $),
\, and counting how many  $ \Lambda $'s  and  $ Y $'s  exist with
$ \, 1 \in \Lambda \, $  and  $ \, \{1,2\} \subseteq Y \, $,  \,
and   --- conversely ---   how many of them exist with  $ \,
\{1,2\} \subseteq \Lambda \, $  and  $ \, 1 \in Y \, $,  \,
we argue
  $$  \delta_{r+s-3}\big([\boldalpha_r,\boldalpha_s]\big) \; =
\; c_{r,s} \cdot [\boldalpha_2,\boldalpha_1] \otimes \boldalpha_2
\otimes {\boldalpha_1}^{\! \otimes (r+s-5)} \, + \, \boldalpha_1
\otimes \varphi_1 \, + \, \boldalpha_2 \otimes \varphi_2 \, + \,
[\boldalpha_2,\boldalpha_1] \otimes \boldalpha_1 \otimes \psi  $$
for some  $ \, \varphi_1, \varphi_2 \in \calH^{\otimes (r+s-4)}
\, $,  $ \, \psi \in \calH^{\otimes (r+s-5)} \, $,  \, and with
  $$  c_{r,s} = {r! \over 2} \cdot {s! \over 3} \cdot {{r+s-5}
\choose {r-2}} - {s! \over 2} \cdot {r! \over 3} \cdot {{s+r-5}
\choose {s-2}}
%
%
 = {{\,2\,} \over {\,3\,}} \, {r \choose 2} {s \choose 2} (s-r) (r+s-5)!
\not= 0 \, .  $$
In particular  $ \; \delta_{r+s-3}\big([\boldalpha_r,\boldalpha_s]\big)
\, = \, c_{r,s} \cdot [\boldalpha_2,\boldalpha_1] \otimes \boldalpha_2
\otimes {\boldalpha_1}^{\! \otimes (r+s-5)} \, + \, \text{\sl l.i.t.}
\, $,  \; where  ``{\sl l.i.t.}''  stands for some further  {\sl
terms\/}  which are  {\sl linearly independent\/}  of  $ \,
[\boldalpha_2, \boldalpha_1] \otimes \boldalpha_2 \otimes
{\boldalpha_1}^{\! \otimes (r+s-5)} \, $  and  $ \, c_{r,s}
\not= 0 \, $.  Then  $ \; \delta_{r+s-3} \big( [\boldalpha_r,
\boldalpha_s] \big) \not= 0 \, $,  \; q.e.d.
                                              \par
   Finally, if  $ \, r > 1 = s \, $  (and similarly if  $ \, r
= 1 < s \, $)  things are simpler.  Indeed, again  {\it (b)\/}
and  {\it (c)\/}  together give  $ \, \kappa \big([\boldalpha_r,
\boldalpha_1]\big) < \kappa(\boldalpha_r) + \kappa(\boldalpha_1)
= (r-1) + 1 = r \, $,  \, and we prove that  $ \; \delta_{r-1}
\big( [\boldalpha_r, \boldalpha_1] \big) \not= 0 \, $.  Like
before,  Lemma 3.4{\it(b)\/}  gives (since  $ \,
\delta_1(\boldalpha_1) = \boldalpha_1 \, $)
  $$  \displaylines{
   \delta_{r-1}\big([\boldalpha_r,\boldalpha_1]\big) \; =
\hskip-15pt  \sum_{\Sb \Lambda \cup Y = \{1,2,\dots,r-1\} \\
\Lambda \cap Y \not= \emptyset, \, |\Lambda| = r-1, \, |Y| = 1
\endSb}  \hskip-35pt  \big[ \delta_\Lambda(\boldalpha_r),
\delta_Y(\boldalpha_1) \big]  =  \sum_{k=1}^{r-1}
\Big[ \delta_{r-1}(\boldalpha_r), 1^{\otimes (k-1)} \otimes
\boldalpha_1 \otimes 1^{\otimes (r-1-k)} \Big] =  \cr
   {} \hfill   = \; \sum_{m=1}^{r-1} {{\,r!\,} \over {\,m+1\,}}
\cdot {\boldalpha_1}^{\!\otimes (m-1)} \otimes [\boldalpha_2,
\boldalpha_1] \otimes {\boldalpha_1}^{\!\otimes (n-1-m)} \;
\not= \; 0 \; ,  \quad  \text{q.e.d.}  \cr }  $$
                                              \par
   {\it (e)} \, We perform induction upon  $ d(b) \, $:  \, the cases
$ \, d(b) = 0 \, $  and  $ \, d(b) = 1 \, $  are dealt with in parts
{\it (c)\/}  and  {\it (d)\/}  of the claim, thus we can assume  $ \,
d(b) \geq 2 \, $,  \, so that  $ \, b = \big[ b', x_\ell \big] \, $
for some  $ \, \ell \in \N_\nu \, $  and some other  $ \, b' \in
B_\nu \, $  with  $ \, d(b') = d(b) - 1 \, $;  \, then  $ \,
\tau(\boldalpha_b) = \tau\big([\boldalpha_{b'}, \boldalpha_\ell]\big)
= \tau(\boldalpha_{b'}) + \tau(\boldalpha_\ell) - 1 \, $,  \,
directly from definitions.  Moreover  $ \, \tau(\boldalpha_\ell)
= \kappa\,(\boldalpha_\ell) \, $  by part  {\it (c)},  \, and
$ \, \tau(\boldalpha_{b'}) = \kappa\,(\boldalpha_{b'}) \, $
by inductive assumption.
                                       \par
   From  {\it (b)\/}  we have  $ \, \kappa(\boldalpha_b) =
\kappa\big([\boldalpha_{b'},\boldalpha_\ell]\big) \leq
\kappa(\boldalpha_{b'}) + \kappa(\boldalpha_\ell) - 1  =
\tau(\boldalpha_{b'}) + \tau(\boldalpha_\ell) - 1 = \tau(\boldalpha_b)
\, $,  \, i.{} e.{}  $ \, \kappa(\boldalpha_b) \leq \tau(\boldalpha_b)
\, $;  \, we must prove the converse, for which it is enough to show
  $$  \delta_{\tau(\boldalpha_b)}(\boldalpha_b) \; = \; c_b \cdot
[\, \cdots [\, [ \undersetbrace{d(b)+1}\to{\boldalpha_1,\boldalpha_2],
\boldalpha_2], \dots, \boldalpha_2} \,] \otimes \boldalpha_2 \otimes
{\boldalpha_1}^{\!\otimes (\tau(\boldalpha_b)-2)} \; + \;
\text{\sl l.i.t.}   \eqno (10.6)  $$
for some  $ \, c_b \in \Bbbk \setminus \{0\} \, $,  \, where  ``{\sl
l.i.t.}''  means the same as before.
                                      \par
   Since  $ \, \tau(\boldalpha_b) = \tau \big( [\boldalpha_{b'},
\boldalpha_\ell] \big) = \tau(\boldalpha_{b'}) + \ell - 2 \, $,
\, computation via  Lemma 3.4{\it (b)\/}  gives
  $$  \displaylines{
   {} \quad   \delta_{\tau(\boldalpha_b)}(\boldalpha_b) \; = \;
\delta_{\tau(\boldalpha_b)} \big([\boldalpha_{b'},\boldalpha_\ell]\big)
\; = \; \sum_{\Sb  \Lambda \cup Y = \{1,\dots,\tau(\boldalpha_b)\}  \\
   \Lambda \cap Y \not= \emptyset  \endSb}
\big[ \delta_\Lambda(\boldalpha_{b'}), \delta_Y(\boldalpha_\ell) \big]
\; = \;   \hfill  \cr
   {} \hfill   = \; \sum_{\Sb  \Lambda \cup Y = \{1, \dots,
\tau(\boldalpha_b)\} \, , \; \Lambda \cap Y \not= \emptyset  \\
   |\Lambda| = \tau(\boldalpha_{b'}) \, , \, |Y| = \ell-1  \endSb}
\hskip-7pt
\big[ j_\Lambda\big(\delta_{\tau(\boldalpha_{b'})}(\boldalpha_{b'})\big),
j_Y\big(\delta_{\ell-1}(\boldalpha_\ell)\big) \big] \; =  \cr
%
%
   {} \hfill   =  \hskip-79pt \sum_{ \Sb  {} \hskip71pt \Lambda \cup Y
= \{1, \dots, \tau(\boldalpha_b)\} \, , \; \Lambda \cap Y \not=
\emptyset  \\
   {} \hskip71pt |\Lambda| = \tau(\boldalpha_{b'}) \, ,
\quad |Y| = \ell-1  \endSb}  \hskip-69pt
\Big[ j_\Lambda \big( c_{b'} \, [\, \cdots [
\undersetbrace{d(b')+1}\to{\boldalpha_1,\boldalpha_2],
\dots, \boldalpha_2} \,] \otimes \boldalpha_2 \otimes
{\boldalpha_1}^{\!\otimes (\tau(\boldalpha_{b'}\!)-2)} \big),
\, j_Y \big( \textstyle{{{\,\ell\,!\,} \over {\,2\,}}}
\, \boldalpha_2 \otimes {\boldalpha_1}^{\!\otimes (\ell-2)}
\big) \Big] + \, \text{\sl l.i.t.} =  \cr
   {} \hfill   = \; c_{b'} \cdot {{\,\ell\,!\,} \over {\,2\,}}
\cdot {{\tau(\boldalpha_b)-2} \choose {\ell-2}} \cdot [\,[\, \cdots
[[ \undersetbrace{d(b')+1+1 \, = \, d(b)+1} \to{\boldalpha_1,
\boldalpha_2], \boldalpha_2], \dots, \boldalpha_2], \boldalpha_2}
\,\big] \otimes \boldalpha_2 \otimes {\boldalpha_1}^{\!\otimes
(\tau(\boldalpha_b)-2)} \; + \; \text{\sl l.i.t.}  \cr }  $$
(using induction about  $ \boldalpha_{b'} $);  this proves (10.6)
with  $ \, c_b =  c_{b'} \cdot {{\,\ell\,!\,} \over {\,2\,}}
\cdot \Big( \! {{\tau(\boldalpha_b)-2} \atop {\ell-2}} \! \Big)
\not= 0 \, $.
                                                 \par
   Thus (10.6) holds, yielding  $ \, \delta_{\tau(\boldalpha_b)}
(\boldalpha_b) \not= 0 \, $,  \, hence  $ \, \kappa(\boldalpha_b)
\geq \tau(\boldalpha_b) \, $,  \, q.e.d.
                                                 \par
   {\it (f)} \, The case  $ \, \ell = 1 \, $  is proved by part
{\it (e)},  so we can assume  $ \, \ell > 1 \, $.  By part  {\it
(b)\/}  and the case  $ \, \ell = 1 \, $  we have  $ \, \kappa\,
(\boldalpha_{b_1} \boldalpha_{b_2} \cdots \boldalpha_{b_\ell})
\leq \sum_{i=1}^\ell \kappa\,(\boldalpha_{b_i}) = \sum_{i=1}^\ell
\tau(\boldalpha_{b_i}) = \tau(\boldalpha_{b_1} \boldalpha_{b_2}
\cdots \boldalpha_{b_\ell}) \, $;  \, so we must only prove the
converse inequality.  We begin with  $ \, \ell = 2 \, $  and
$ \, d(b_1) = d(b_2) = 0 \, $,  \, so  $ \, \boldalpha_{b_1}
= \boldalpha_r \, $,  $ \, \boldalpha_{b_2} = \boldalpha_s
\, $,  \, for some  $ \, r $,  $ s \in \N_\nu \, $.
                                                 \par
   If  $ \, r = s = 1 \, $  then  $ \, \kappa(\boldalpha_r) =
\kappa(\boldalpha_s) = \kappa(\boldalpha_1) = 1 \, $,  \, by
part  {\it (c)}.  Then
  $$  \delta_2(\boldalpha_1 \, \boldalpha_1) \; =
\delta_2(\a_1 \, \a_1) \; = \; {(\id - \epsilon)}^{\otimes 2}
\Delta\big({\a_1}^{\!2}\big) \; = \; 2 \cdot \a_1 \otimes \a_1 \;
= \; 2 \cdot \boldalpha_1 \otimes \boldalpha_1 \; \not= \; 0  $$
so that  $ \, \kappa(\boldalpha_1 \, \boldalpha_1) \geq 2
= \kappa(\boldalpha_1) + \kappa(\boldalpha_1) \, $,  \, hence
$ \, \kappa(\boldalpha_1 \, \boldalpha_1) = \kappa(\boldalpha_1)
+ \kappa(\boldalpha_1) \, $,  \, q.e.d.
                                                \par
   If  $ \, r > 1 = s \, $  (and similarly if  $ \, r = 1 < s
\, $)  then  $ \, \kappa(\boldalpha_r) = r-1 \, $,  $ \,
\kappa(\boldalpha_s) = \kappa(\boldalpha_1) = 1 \, $,  \,
by part  {\it (c)}.  Then  Lemma 3.4{\it(b)\/}  gives
  $$  \displaylines{
   \delta_r(\boldalpha_r \, \boldalpha_1) \; = \; \sum\nolimits_{\Sb
   \Lambda \cup Y = \{1,\dots,r\}  \\
   |\Lambda| = r-1 \, , \; |Y| = 1  \endSb}  \hskip-9pt
\delta_\Lambda(\boldalpha_r) \, \delta_Y(\boldalpha_1)
\; =   \hfill  \cr
   {} \quad   = \; \sum_{m=1}^r \sum_{k<m}
{{\,r!\,} \over {\,m+1\,}} \cdot \big( {\boldalpha_1}^{\!\otimes
(k-1)} \otimes 1 \otimes {\boldalpha_1}^{\!\otimes (m-1-k)} \otimes
\boldalpha_2 \otimes {\boldalpha_1}^{\!\otimes (r-1-m)} \big)
\times   \hfill  \cr
   {} \hfill   \times \big( 1^{\otimes (k-1)} \otimes \boldalpha_1
\otimes 1^{\otimes (r-k)} \big) \; +  \quad {}  \cr
   {} \quad   + \; \sum_{m=1}^r \sum_{k>m}
{{\,r!\,} \over {\,m+1\,}} \cdot \big( {\boldalpha_1}^{\!\otimes
(m-1)} \otimes \boldalpha_2 \otimes {\boldalpha_1}^{\!\otimes
(k-1-m)} \otimes 1 \otimes {\boldalpha_1}^{\!\otimes (r-1-k)}
\big) \times   \hfill  \cr
   {} \hfill   \times \big( 1^{\otimes (k-1)} \otimes \boldalpha_1
\otimes 1^{\otimes (r-k)} \big) \; =  \quad {}  \cr
   {} \hfill   = \; \sum\nolimits_{m=1}^r {{\,r!\,} \over {\,m+1\,}} \cdot
{\boldalpha_1}^{\!\otimes (m-1)} \otimes \boldalpha_2 \otimes
{\boldalpha_1}^{\!\otimes (r-1-m)} \; \not= \; 0  \cr }  $$
so that  $ \, \kappa(\boldalpha_r \, \boldalpha_1) \geq r
= \kappa(\boldalpha_r) + \kappa(\boldalpha_1) \, $  and so
$ \, \kappa(\boldalpha_r \, \boldalpha_1) = \kappa(\boldalpha_r)
+ \kappa(\boldalpha_1) \, $,  \, q.e.d.
                                                \par
   Finally let  $ \, r, s > 1 \, $  (and  $ \, r \not= s \, $).
Then  $ \, \kappa(\boldalpha_r) = r-1 \, $,  $ \,
\kappa(\boldalpha_s) = s-1 \, $,  \, by part  {\it (c)\/};
then  Lemma 3.4{\it(b)\/}  gives
  $$  \delta_{r+s-2}\big(\boldalpha_r \, \boldalpha_s \big) \;
=  \hskip-7pt  \sum_{\Sb \Lambda \cup Y = \{1,\dots,r+s-2\}  \\
|\Lambda| = r-1 \, , \; |Y| = s-1  \endSb}  \hskip-23pt
\delta_\Lambda(\boldalpha_r) \cdot \delta_Y(\boldalpha_s) \;
=  \hskip-13pt  \sum_{\Sb \Lambda \cup Y = \{1,\dots,r+s-2\} \\
|\Lambda| = r-1 \, , \; |Y| = s-1  \endSb}  \hskip-25pt
j_\Lambda\big( \delta_{r-1}(\boldalpha_r) \big) \cdot
j_Y \big( \delta_{s-1}(\boldalpha_s) \big) \, .  $$
Using (10.5) in the form  $ \; \delta_{t-1}(\a_t) =
\sum_{m=1}^{t-1} {{\,t\,!\,} \over 2} \cdot \boldalpha_2
\otimes {\boldalpha_1}^{\!\otimes (t-2)} + \boldalpha_1
\otimes \eta_t \; $  (for some  $ \, \eta_t \in \calH \, $
and  $ \, t \in \{r,s\} \, $)  and counting how many
$ \Lambda $'s  and  $ Y $'s  exist with  $ \, 1 \in \Lambda
\, $  and  $ \, 2 \in Y \, $  and viceversa   --- actually, it
is a matter of counting  $ (r-2,s-2) $-shuffles  ---   we argue
  $$  \delta_{r+s-2}\big(\boldalpha_r \, \boldalpha_s\big) \;
= \; e_{r,s} \cdot \boldalpha_2 \otimes \boldalpha_2 \otimes
{\boldalpha_1}^{\! \otimes (r+s-4)} \, + \, \boldalpha_1
\otimes \varphi  $$
for some  $ \, \varphi \in \calH^{\otimes (r+s-3)} \, $  and with
  $$  e_{r,s} \; = \; {{\,r!\,} \over {\,2\,}} \cdot {{\,s!\,}
\over {\,2\,}} \cdot \Bigg( {{r+s-4} \choose {r-2}} + {{s+r-4}
\choose {s-2}} \Bigg) \; = \; {{\, r! \, s! \,} \over {\,2\,}}
\cdot {{r+s-4} \choose {r-2}} \; \not= \; 0 \; .  $$
In particular  $ \; \delta_{r+s-2} \big( \boldalpha_r \,
\boldalpha_s \big) \, = \, e_{r,s} \cdot \boldalpha_2 \otimes
\boldalpha_2 \otimes {\boldalpha_1}^{\! \otimes (r+s-4)} \, +
\, \text{\sl l.i.t.} \, $,  \; where  ``{\sl l.i.t.}''  stands
again for some further  {\sl terms\/}  which are  {\sl linearly
independent\/}  of  $ \, \boldalpha_2 \otimes \boldalpha_2 \otimes
{\boldalpha_1}^{\! \otimes (r+s-4)} \, $  and  $ \, e_{r,s}
\not= 0 \, $.  Then  $ \; \delta_{r+s-2} \big( \boldalpha_r \,
\boldalpha_s \big) \not= 0 \, $,  \, so  $ \, \kappa(\boldalpha_r
\, \boldalpha_1) \geq r+s-2 = \kappa(\boldalpha_r) +
\kappa(\boldalpha_1) \, $,  \, q.e.d.
                                                \par
   Now let again  $ \, \ell = 2 \, $  but  $ \, d(b_1), d(b_2) >
0 \, $.  Set  $ \, \kappa_i := \kappa(\boldalpha_{b_i}) \, $  for
$ \, i = 1, 2 \, $.  Applying (10.6) to  $ \, b = b_1 \, $  and
$ \, b = b_2 \, $  (and reminding  $ \, \tau \equiv \kappa \, $)
gives   
  $$  \displaylines{
   \delta_{\kappa_1 + \kappa_2}(\boldalpha_{b_1} \, \boldalpha_{b_2})
\; =  \hskip-4pt  \sum_{\Lambda \cup Y = \{1, \dots, \kappa_1
+ \kappa_2\}}  \hskip-23pt
\delta_\Lambda(\boldalpha_{b_1}) \, \delta_Y(\boldalpha_{b_2})
\; =  \hskip-4pt  \sum_{\Sb
      \Lambda \cup Y = \{1, \dots, \kappa_1 + \kappa_2\}  \\
      |\Lambda| = \kappa_1, \, |Y| = \kappa_2  \endSb}  \hskip-23pt
j_\Lambda\big(\delta_{\kappa_1}(\boldalpha_{b_1})\big) \,
j_Y\big(\delta_{\kappa_2}(\boldalpha_{b_2})\big) \; =  \cr
   = \sum_{\Sb  \Lambda \cup Y = \{1, \dots, \kappa_1 + \kappa_2\}  \\
       |\Lambda| = \kappa_1, \, |Y| = \kappa_2  \endSb}  \hskip-9pt
j_\Lambda \Big( c_{b_1} \cdot [\, \cdots [\, [ \undersetbrace{d(b_1)+1}
\to{\boldalpha_1,\boldalpha_2],\boldalpha_2], \dots, \boldalpha_2} \,]
\otimes \boldalpha_2 \otimes {\boldalpha_1}^{\!\otimes (\kappa_1 - 2)}
\; + \; \hbox{\sl l.i.t.} \,\Big) \, \times   \hfill  \cr
   {} \hfill   \times \, j_Y \Big( c_{b_2} \cdot [\, \cdots [\, [
\undersetbrace{d(b_2)+1} \to{\boldalpha_1,\boldalpha_2],\boldalpha_2],
\dots, \boldalpha_2} \,] \otimes \boldalpha_2 \otimes {\boldalpha_1}^{\!
\otimes (\kappa_2 - 2)} \; + \; \hbox{\sl l.i.t.} \,\Big) \; =  \cr
   = \; c_{b_1} \, c_{b_2} \cdot 2 \, {{\kappa_1 + \kappa_2 - 4 \,}
\choose {\kappa_1 - 2}} \, \times   \hfill  \cr
   {} \hfill   \times \, [\, \cdots [\, [ \undersetbrace{d(b_1)+1} \to
{\boldalpha_1,\boldalpha_2], \boldalpha_2], \dots, \boldalpha_2} \,]
\otimes [\, \cdots [\, [ \undersetbrace{d(b_2)+1} \to {\boldalpha_2,
\boldalpha_1], \boldalpha_2], \dots, \boldalpha_2} \,] \otimes
\boldalpha_2 \otimes \boldalpha_2 \otimes {\boldalpha_1}^{\!
\otimes (\kappa_1 + \kappa_2 - 4)} \; + \, \hbox{\sl l.i.t.}
\cr }  $$
which proves the claim for  $ \, \ell = 2 \, $.  In addition, we can
take this last result as the basis of induction (on  $ \ell \, $)
to prove the following: for all  $ \, \underline{b} := (b_1, \dots,
b_\ell) \in {B_\nu}^{\!\ell} \, $,  \, one has
  $$  \delta_{|\underline{\kappa}|}\Bigg( \prod_{i=1}^\ell
\boldalpha_{b_i} \Bigg) = c_{\underline{b}} \, \Bigg(
\bigotimes_{i=1}^\ell \, [\, \cdots [\, [ \undersetbrace{d(b_i)+1}
\to {\boldalpha_1,\boldalpha_2], \boldalpha_2], \dots, \boldalpha_2}
\,] \Bigg) \otimes {\boldalpha_2}^{\!\otimes \ell} \otimes
{\boldalpha_1}^{\!\otimes (|\underline{\kappa}| - 2\,\ell)}
+  \hbox{\sl l.i.t.}   \hfill \hskip11pt (10.7)  $$
for some  $ \, c_{\underline{b}} \in \Bbbk \setminus \{0\} \, $,
\, with  $ \, |\underline{\kappa}| := \sum_{i=1}^\ell \kappa_i \, $
and  $ \, \kappa_i := \kappa(\boldalpha_{b_i}) \, $  ($ \, i = 1,
\dots, \ell \, $).  The induction step, from  $ \ell \, $  to
$ (\ell + 1) $,  amounts to compute  (with  $ \, \kappa_{\ell
+ 1} := \kappa (\boldalpha_{b_{\ell + 1}}) \, $)
  $$  \displaylines{
   \delta_{|\underline{\kappa}| + \kappa_{\ell + 1}}
(\boldalpha_{b_1} \cdots \boldalpha_{b_\ell} \cdot \boldalpha_{b_{\ell
+ 1}}) \;\; =  \sum_{\Lambda \cup Y = \{1,\dots,|\underline{\kappa}|
+ \kappa_{\ell + 1}\}}  \hskip-5pt
\delta_\Lambda(\boldalpha_{b_1} \cdots \boldalpha_{b_\ell})
\, \delta_Y(\boldalpha_{b_{\ell + 1}}) \; =   \hfill  \cr
   {} \hfill   = \sum_{\Sb
     \Lambda \cup Y = \{1, \dots, |\underline{\kappa}|
+ \kappa_{\ell + 1}\}  \\
     |\Lambda| = |\underline{\kappa}|, \,
|Y| = \kappa_{\ell + 1}  \endSb}  \hskip-5pt
j_\Lambda\big(\delta_{|\underline{\kappa}|}
(\boldalpha_{b_1} \cdots \boldalpha_{b_\ell})\big) \cdot
j_Y \big( \delta_{\kappa_{\ell + 1}} (\boldalpha_{b_{\ell
+ 1}}) \big) \; =  \cr   
   =  \hskip-7pt  \sum_{\Sb
       \Lambda \cup Y = \{1, \dots, |\underline{\kappa}|
+ \kappa_{\ell + 1}\}  \\
       |\Lambda| = |\underline{\kappa}|, \, |Y| =
\kappa_{\ell + 1} \endSb}  \hskip-11pt
j_\Lambda \Bigg( c_{\underline{b}} \, \cdot \Bigg(
\bigotimes_{i=1}^\ell \; [\, \cdots [\, [ \undersetbrace{d(b_i)+1}
\to {\boldalpha_1,\boldalpha_2], \boldalpha_2], \dots, \boldalpha_2}
\,] \Bigg) \otimes {\boldalpha_2}^{\!\otimes \ell} \otimes
{\boldalpha_1}^{\!\otimes (|\underline{\kappa}| - 2\,\ell)}
\; + \, \hbox{\sl l.{}i.{}t.}  \Bigg) \times   \hfill  \cr
   {} \hfill   \times j_Y \Big( c_{b_{\ell + 1}} \cdot [\, \cdots [\,
[ \undersetbrace{d(b_{\ell + 1})+1} \to{\boldalpha_1,\boldalpha_2],
\boldalpha_2], \dots, \boldalpha_2} \,] \otimes \boldalpha_2
\otimes {\boldalpha_1}^{\!\otimes (\kappa_{\ell + 1} - 2)}
\; + \; \hbox{\sl l.i.t.} \,\Big) \; =  \cr
 }  $$   
  $$  \displaylines{ 
   = \; c_{\underline{b}} \, c_{b_{\ell + 1}} \cdot (\ell + 1)
\, {{|\underline{\kappa}| + \kappa_{\ell + 1} - 2 \, (\ell + 1)}
\choose {|\underline{\kappa}| - 2 \, \ell}} \; \cdot
\Bigg( \bigotimes_{i=1}^\ell \; [\, \cdots [\, [
\undersetbrace{d(b_i)+1} \to {\boldalpha_1,\boldalpha_2],
\boldalpha_2], \dots, \boldalpha_2} \,] \Bigg) \otimes   \hfill  \cr
   {} \hfill   \otimes [\, \cdots [\, [ \undersetbrace{d(b_{\ell + 1})
+1} \to {\boldalpha_1, \boldalpha_2], \boldalpha_2], \dots,
\boldalpha_2} \,] \otimes {\boldalpha_2}^{\!\otimes (\ell + 1)} \otimes
{\boldalpha_1}^{\!\otimes (|\underline{\kappa}| + \kappa_{\ell + 1}
- 2 \, (\ell + 1))} \; + \, \hbox{\sl l.i.t.}  \cr }  $$
which proves (10.7) for  $ \, (\underline{b} \, , b_{\ell + 1}) \, $
with  $ \, c_{(\underline{b} \, , b_{\ell + 1})} = c_{\underline{b}}
\, c_{b_{\ell + 1}} \cdot (\ell + 1) \Big( {{|\underline{\kappa}| +
\kappa_{\ell + 1} - 2 \, (\ell + 1)} \atop {|\underline{\kappa}| -
2 \, \ell}} \Big) \, \not= \, 0 \, $.  Finally, (10.7) yields  $ \,
\delta_{|\underline{\kappa}|}(\boldalpha_{b_1} \hskip-0,7pt \cdots
\boldalpha_{b_\ell}) \not= 0 \, $,  \, so
       \hbox{$ \kappa(\boldalpha_{b_1} \cdots \boldalpha_{b_\ell})
\hskip-0,9pt \geq \hskip-1,5pt \kappa(\boldalpha_{b_1}) + \cdots
+ \hskip-0,3pt \kappa(\boldalpha_{b_\ell}) $,  q.e.d.}
                                         \par
   {\it (g)} \, Part  {\it (d)\/}  proves the claim for  $ \, d(b_1)
= d(b_2) = 0 \, $,  \, that is  $ \, b_1, b_2 \in {\{x_n\}}_{n \in \N}
\, $.  Moreover, when  $ \, b_2 = x_n \in {\{x_m\}}_{m \in \N_\nu}
\, $  we can replicate the proof of part  {\it (d)\/}  to show
that  $ \, \kappa \big( [\boldalpha_{b_1},\boldalpha_{b_2}]
\big) = \kappa \big( [\boldalpha_{b_1},\boldalpha_n] \big)
= \partial_-\big([\boldalpha_{b_1},\boldalpha_n]\big) -
d\big([\boldalpha_{b_1},\boldalpha_n]\big) \, $:  \, but
the latter is exactly  $ \, \tau \big([\boldalpha_{b_1},
\boldalpha_{b_2}]\big) \, $,  \, q.e.d.  Everything is
similar if  $ \, b_1 = x_n \in {\{x_m\}}_{m \in \N_\nu} \, $.
                                         \par
   Now let  $ \, b_1, b_2 \in B_\nu \setminus {\{x_n\}}_{n \in \N_\nu}
\, $.  Then  {\it (b)}  gives  $ \, \kappa\big([\boldalpha_{b_1},
\boldalpha_{b_2}]\big) \leq \kappa\,(\boldalpha_{b_1}) + \kappa\,
(\boldalpha_{b_2}) - 1 = \tau \big( [\boldalpha_{b_1},
\boldalpha_{b_2}] \big) $.  Applying (10.6) to  $ \, b =
b_1 \, $  and  $ \, b = b_2 \, $  we get, for  $ \, \kappa_i
:= \kappa(\boldalpha_{b_i}) $  ($ \, i = 1, 2 \, $)
  $$  \displaylines{
   \delta_{\kappa_1 + \kappa_2 - 1} \big( [\boldalpha_{b_1},
\boldalpha_{b_2}] \big) \; =  \hskip-4pt  \sum_{\Sb
\Lambda \cup Y = \{ 1, \dots, \kappa_1 + \kappa_2 - 1 \}  \\
           \Lambda \cap Y \not= \emptyset  \endSb}  \hskip-23pt
\big[ \delta_\Lambda(\boldalpha_{b_1}), \delta_Y(\boldalpha_{b_2})
\big] \; =   \hfill  \cr
   {} \hfill   =  \sum_{\Sb
      \Lambda \cup Y = \{1, \dots, \kappa_1 + \kappa_2 - 1\}  \\
      |\Lambda| = \kappa_1, \, |Y| = \kappa_2  \endSb}  \hskip-23pt
\Big[ j_\Lambda\big(\delta_{\kappa_1}(\boldalpha_{b_1})\big),
j_Y\big(\delta_{\kappa_2}(\boldalpha_{b_2})\big) \Big] \; =  \cr
   = \sum_{\Sb  \Lambda \cup Y = \{1, \dots, \kappa_1 + \kappa_2\}  \\
       |\Lambda| = \kappa_1, \, |Y| = \kappa_2  \endSb}  \hskip-5pt
\Big[ j_\Lambda \Big( c_{b_1} \cdot [\, \cdots [\, [
\undersetbrace{d(b_1)+1} \to{\boldalpha_1,\boldalpha_2],
\boldalpha_2], \dots, \boldalpha_2} \,] \otimes \boldalpha_2
\otimes {\boldalpha_1}^{\!\otimes (\kappa_1 - 2)} \; + \;
\hbox{\sl l.i.t.} \,\Big) \, \times   \hfill  \cr
   {} \hfill   \times \, j_Y \Big( c_{b_2} \cdot [\, \cdots [\, [
\undersetbrace{d(b_2)+1} \to{\boldalpha_1,\boldalpha_2],\boldalpha_2],
\dots, \boldalpha_2} \,] \otimes \boldalpha_2 \otimes {\boldalpha_1}^{\!
\otimes (\kappa_2 - 2)} \; + \; \hbox{\sl l.i.t.} \,\Big) \Big]
\; =  \cr
   = \; c_{b_1} \, c_{b_2} \cdot 2 \, {{\kappa_1 + \kappa_2 - 4 \,}
\choose {\kappa_1 - 2}} \, \times   \hfill  \cr
   {} \hfill   \times \, \big[ [\, \cdots [\, [ \undersetbrace{d(b_1)+1}
\to {\boldalpha_1,\boldalpha_2], \boldalpha_2], \dots, \boldalpha_2}
\,] ,  [\, \cdots [\, [ \undersetbrace{d(b_2)+1} \to {\boldalpha_1,
\boldalpha_2], \boldalpha_2], \dots, \boldalpha_2} \,] \big] \otimes
\boldalpha_2 \otimes \boldalpha_2 \otimes {\boldalpha_1}^{\!
\otimes (\kappa_1 + \kappa_2 - 4)} \; + \, \hbox{\sl l.i.t.}
\cr }  $$
(note that  $ \, d(b_i) \geq 1 \, $  because  $ b_i \not\in \big\{\,
x_n \,\big|\, n \in \N_\nu \,\big\} \, $  for  $ \, i = 1, 2 \, $).
In particular this means  $ \, \delta_{\kappa_1 + \kappa_2 - 1}
\big( [\boldalpha_{b_1}, \boldalpha_{b_2}] \big) \not= 0 \, $,
\, thus  $ \, \kappa\big([\boldalpha_{b_1},\boldalpha_{b_2}]\big)
\geq \kappa\,(\boldalpha_{b_1}) + \kappa\,(\boldalpha_{b_2}) - 1
= \tau\big([\boldalpha_{b_1},\boldalpha_{b_2}]\big) \, $.   \qed
\enddemo

\vskip7pt

\proclaim{Lemma 10.15} \, Let  $ V $  be a  $ \Bbbk $--vector  space,
and  $ \, \psi \in \text{\sl Hom}_{\,\Bbbk}\big(V, V \wedge V) \, $.
Let  $ \L(V) $  be the free Lie algebra over  $ V $,  \, and  $ \,
\psi_{d\L} \in \text{\sl Hom}_{\,\Bbbk}\big(\L(V), \L(V) \wedge \L(V)
\big) \, $  the unique extension of  $ \, \psi $  from  $ V $  to
$ \L(V) $  by derivations, i.e.~such that  $ \, \psi_{d\L}{\big|}_V
= \psi \, $  and  $ \; \psi_{d\L}\big([x,y]\big) =
\big[ x \otimes 1 + 1 \otimes x,
 \allowbreak
\, \psi_{d\L}(y) \big] + \big[ \psi_{d\L}(x),
\, y \otimes 1 + 1 \otimes y \, \big]
= x.\psi_{d\L}(y) - y.\psi_{d\L}(x) \; $
in the  $ \L(V) $--module  $ \L(V) \wedge \L(V) \, $,
$ \, \forall \, x, y \in \L(V) \, $.
   \hbox{Let  $ \, K := \text{\sl Ker}\,(\psi) \, $:  then
$ \, \text{\sl Ker}\,\big(\psi_{d\L}\big) = \L(K) \, $,  \,
the free Lie algebra over  $ K \, $.}
\endproclaim

\demo{Proof} For each  $ \, z \in \L(V) \, $  set  $ \, z^\otimes :=
z \otimes 1 + 1 \otimes z \, $.  Let  $ \, I \, $  be a complement
of  $ K $  inside  $ V $,  so that  $ \, V = K \oplus I \, $,  and
$ \psi{\big|}_I $  is injective while  $ \, \psi{\big|}_K = 0_K
\, $.  Let  $ B_K $  and  $ B_I $  be bases of  $ K $  and  $ I $
respectively; then there is a basis of  $ \L(V) $  made of Lie
monomials of the form  $ \, x_{\underline{i}} := \big[ \big[ \cdots
\big[ \cdots \big[ \big[ x_{i_1}, x_{i_2} \big], x_{i_3} \big] \dots,
x_{i_s} \big] \dots,  x_{i_{k-1}} \big], x_{i_k} \big] \, $  (with
$ \, \underline{i} = (i_1,i_2,i_3,\dots,i_s,\dots,i_{k-1},i_k) \, $)
for some  $ \, x_{i_r} \in B_K \cup B_I \; $:  \; for these Lie
monomials definitions yield  $ \; \psi_{d\L}(x_{\underline{i}\,}) =
\sum_{x_{i_s} \in B_I} \big[ \big[ \cdots \big[ \cdots \big[ \big[
x_{i_1}^{\,\otimes}, x_{i_2}^{\,\otimes} \big], x_{i_3}^{\,\otimes}
\big] \dots, \psi(x_{i_s}) \big] \dots, x_{i_{k-1}}^{\,\otimes} \big],
x_{i_k}^{\,\otimes} \big] \, $.  In addition, since  $ \psi{\big|}_I $
is injective the set  $ \, {\big\{\, z^\otimes \,\big\}}_{z \in B_I
\cup B_K} \cup \, {\big\{ \psi(b) \,\}}_{b \in B_I} \, $  is linearly
independent.  Then the set of all Lie monomials  $ \, y_{\underline{i}}
:= \big[ \big[ \cdots \big[ \cdots \big[ \big[ y_{i_1}^{\,\bullet},
y_{i_2}^{\,\bullet} \big], y_{i_3}^{\,\bullet} \big] \dots, y_{i_s}^{\,
\bullet} \big] \dots, y_{i_{k-1}}^{\,\bullet} \big], y_{i_k}^{\,\bullet}
\big] \, $  with the same  $ \underline{i} \, $'s  which give the basis
of  $ \L(V) $  and with  $ \, y_j^{\,\bullet} \in \big\{ x_j^{\,
\otimes}, \psi(x_j) \big\} \, $  is again a linearly independent set
inside  $ \, \L(V) \wedge \L(V) \, $.  Therefore, for a general  $ \,
x = \sum_{\underline{i}} c_{\underline{i}} \, x_{\underline{i}} \in
\L(V) \, $  we have  $ \; \psi_{d\L}(x) \, = \, \sum_{\underline{i}}
c_{\underline{i}} \sum_{x_{i_s} \in B_I} \big[ \big[ \cdots \big[
\cdots \big[ \big[ x_{i_1}^{\,\otimes}, x_{i_2}^{\,\otimes} \big],
x_{i_3}^{\,\otimes} \big] \dots, \psi(x_{i_s}) \big] \dots,
x_{i_{k-1}}^{\,\otimes} \big], x_{i_k}^{\,\otimes} \big] \, $;
\; thus if  $ \, \psi_{d\L}(x) = 0 \, $  we necessarily have  $ \,
c_{\underline{i}} = 0 \, $  for all  $ \, {\underline{i}} $
which sport at least one  $ \, x_{i_s} \in B_I \, $.  The outcome
is that  $ \, x \in \text{\sl Ker}\,(\psi_{d\L}) \, $  implies
$ \, x = \sum_{\underline{i} \; : \, x_{i_s} \in B_K \,
(\,\forall\, s)} c_{\underline{i}} \, x_{\underline{i}}
\in \L(K) \, $;  \, thus  $ \, \text{\sl Ker}\,(\psi_{d\L})
\subseteq \L(K) \, $,  \, and the converse inclusion is
clear because  $ \psi_{d\L} $  is a derivation.   \qed
\enddemo

\vskip5pt

\proclaim{Lemma 10.16}  The Lie cobracket  $ \delta $  of
$ \, U(\L_\nu) $  preserves  $ \tau $.  That is, for each  $ \,
\vartheta \in U(\L_\nu) \, $  in the expansion  $ \, \delta_2
(\vartheta) = \sum_{\underline{b}_1, \underline{b}_{\,2} \in \Bbb{B}}
c_{\underline{b}_1,\underline{b}_{\,2}} \, \boldalpha_{\underline{b}_1}
\otimes \boldalpha_{\underline{b}_{\,2}} \, $  (w.r.t.~the basis  $ \,
\Bbb{B} \otimes \Bbb{B} \, $,  \, where  $ \Bbb{B} $  is a PBW basis
as in \S 10.2 w.r.t.~some total order of  $ B_\nu $)  we have
$ \, \tau \big( \underline{\hat{b}}_{\,1} \big) + \tau \big(
\underline{\hat{b}}_{\,2} \big) = \tau(\vartheta) \, $  for some
$ \underline{\hat{b}}_{\,1} $,  $ \underline{\hat{b}}_{\,2} $  with
$ \, c_{\underline{\hat{b}}_1, \underline{\hat{b}}_{\,2}} \not= 0
\, $,  \, so
   \hbox{$ \, \tau \big( \delta(\vartheta) \big) := \max \big\{
\tau(\underline{b}_{\,1}) + \tau(\underline{b}_{\,2}) \;\big|\;
c_{\underline{b}_1,\underline{b}_{\,2}} \not= 0 \big\} =
\tau(\vartheta) \, $  if  $ \, \delta(\vartheta) \not= 0 \, $.}
 \vskip-5pt   
\endproclaim

\demo{Proof}  It follows from Proposition 10.13 that  $ \, \tau \big(
\delta(\vartheta) \big) \leq \tau(\vartheta) \, $;  \, so  $ \; \delta
\, \colon \, U(\L_\nu) \longrightarrow {U(\L_\nu)}^{\otimes 2} \; $
is a morphism of filtered algebras, hence it naturally induces a
morphism of graded algebras  $ \; \overline{\delta} \, \colon \,
G_{\underline{\varTheta}}\big(U(\L_\nu)\big) \! \llongrightarrow
{G_{\underline{\varTheta}}\big(U(\L_\nu)\big)}^{\otimes 2} \; $
(notation of \S\S 5.3--4).  Therefore proving the claim is equivalent
to showing that  $ \, \text{\sl Ker}\,\big(\overline{\delta}\,\big) =
G_{\underline{\varTheta} \, \cap {Ker}(\delta)} \big( \text{\sl Ker}
\,(\delta) \big) =: \overline{\text{\sl Ker}\,(\delta)} \, $,
\, the latter being thought of as naturally embedded into
$ G_{\underline{\varTheta}}\big(U(\L_\nu)\big) \, $.
                                       \par
   By construction,  $ \, \tau(x \, y - y \, x) = \tau\big([x,y]\big)
< \tau(x) + \tau(y) \, $  for  $ \, x $,  $ y \in U(\L_\nu) \, $,
\, so  $ G_{\underline{\varTheta}}\big(U(\L_\nu)\big) $  is
commutative: indeed, it is clearly isomorphic   --- as an algebra
---   to  $ S(V_\nu) $,  the symmetric algebra over  $ V_\nu \, $.
Moreover,  $ \delta $  acts as a derivation, that is  $ \; \delta(x\,y)
= \delta(x) \, \Delta(y) + \Delta(x) \, \delta(y) \; $  (for all
$ \, x $,  $ y \in U(\L_\nu) \, $),  \, thus the same holds for
$ \overline{\delta} $  too.  Like in Lemma 10.15, since  $ \,
G_{\underline{\varTheta}}\big(U(\L_\nu)\big) \, $  is generated
by  $ \; G_{\underline{\varTheta} \cap \L_\nu}(\L_\nu) =:
\overline{\L_\nu} \; $  it follows that  $ \, \text{\sl Ker}\,
\big(\overline{\delta}\,\big) \, $  is the free (associative
sub)algebra over  $ \, \text{\sl Ker}\, \Big( \overline{\delta}\,
{\big|}_{\overline{\L_\nu}} \Big) \, $,  \, in short  $ \; \text{\sl
Ker}\,\big(\overline{\delta}\,\big) = \Big\langle \text{\sl Ker}\,
\Big( \overline{\delta}\,{\big|}_{\overline{\L_\nu}} \Big)
\Big\rangle \; $.
                                             \par
   Now, by definition  $ \; \delta(x_n) = \sum_{\ell=1}^{n-1} (\ell+1)
\, x_\ell \wedge x_{n-\ell} \; $  (cf.~Theorem 10.6) is a sum of
$ \tau $--{\,}homogeneous  terms of  $ \tau $--{\,}degree  equal to
$ \, (n-1) = \tau(x_n) \, $.  Since in addition  $ \delta $  enjoys
$ \, \delta\big([x,y]\big) = \big[x \otimes 1 + 1 \otimes x, \delta(y)
\big] + \big[\delta(x), y \otimes 1 + 1 \otimes y \,\big] \, $
(for all  $ \, x $,  $ y \in \L_\nu \, $)  we have that  $ \,
\delta{\big|}_{\L_\nu} \, $  is even  $ \tau $--{\,}homogeneous,
which means that  $ \delta\big(\tau(z)\big) $  either is zero or can
be written as a sum whose summands are all  $ \tau $--{\,}homogeneous
terms of  $ \tau $--degree  equal to  $ \, \tau(z) \, $,  \, for any 
$ \tau $--homogeneous  $ \, z \in \L_\nu \, $;  \, this implies that
the induced map  $ \, \overline{\delta}\,{\big|}_{\overline{\L_\nu}}
\, $  enjoys  $ \; \overline{\delta}\,{\big|}_{\overline{\L_\nu}} \big(
\overline{\vartheta} \,\big) = \overline{0} \iff \delta(\vartheta)
= 0 \; $  for any  $ \, \vartheta \in \L_\nu \, $,  \, whence  $ \;
\text{\sl Ker}\, \big(\, \overline{\delta}\,{\big|}_{\overline{\L_\nu}}
\,\big) \, = \, \overline{\text{\sl Ker}\, \big(\delta{\big|}_{\L_\nu}
\big)} \, $.
%
%
%
%
  On the upshot we get  $ \; \text{\sl Ker}\, \big(
\overline{\delta} \,\big) = \Big\langle \text{\sl Ker}\, \Big( \overline{\delta}\,{\big|}_{\overline{\L_\nu}} \Big) \Big\rangle
= \Big\langle \overline{\text{\sl Ker}\, \big(\delta{\big|}_{\L_\nu}
\big)} \, \Big\rangle = \overline{\text{\sl Ker}\,(\delta)} \; $,
\; q.e.d.   \qed
\enddemo
%
%
 \eject   

\proclaim{Proposition 10.17}  $ \, \underline{D} = \underline{\varTheta}
\, $,  \, that is  $ \, D_n = \varTheta_n \, $  for all  $ \, n \in \N
\, $,  or  $ \, \kappa = \tau \, $.  Therefore, given any total order
$ \preceq $  in  $ B_\nu \, $,  \, the set  $ \, \calA_{\leq n} =
\calA \cap \varTheta_n = \calA \cap D_n \, $  of ordered monomials
  $$  \calA_{\leq n} \; = \; \Big\{\, \boldalpha_{\underline{b}} =
\boldalpha_{b_1} \cdots \boldalpha_{b_k} \,\Big|\; k \in \N \, ,
\; b_1, \dots, b_k \in B_\nu \, , \; b_1 \preceq \cdots \preceq
b_k \, ,  \; \tau(\underline{b}\,) \leq n \,\Big\}  $$
is a\/  $ \Bbbk $--basis  of  $ \, D_n \, $,  \, and  $ \, \calA_n
:= \big(\, \calA_{\leq n} \! \mod D_{n-1} \big) \, $  is a\/
$ \Bbbk $--basis  of  $ \, D_n \big/ D_{n-1} \; $  ($ \,
\forall \; n \in \N \, $).
\endproclaim

\demo{Proof} Clearly the claim about the  $ \calA_{\leq n} $'s
and the claim about the  $ \calA_n $'s  are equivalent, and
either of these claims is equivalent to  $ \, \underline{D}
= \underline{\varTheta} \, $.  Note also that  $ \, \calA_n
:= \big(\, \calA_{\leq n} \! \mod D_{n-1} \big) = \big(\,
\calA_{\leq n} \setminus \calA_{\leq n-1} \! \mod D_{n-1}
\big) \, $,  \, where clearly  $ \, \calA_{\leq n} \setminus
\calA_{\leq n-1} = \big\{\, \boldalpha_{\underline{b}} \in
\calA \;\big|\; \tau(\underline{b}\,) = n \,\big\} \, $.
                                             \par
   By  Lemma 10.14{\it (f)\/}  we have  $ \, \calA_{\leq n} =
\calA \bigcap \varTheta_n \subseteq \calA \bigcap D_n \subseteq
D_n \, $;  \, since  $ \calA $  is a basis,  $ \calA_{\leq n} $  is
linearly independent and is a  $ \Bbbk $--basis  of  $ \varTheta_n $
(by definition): so  $ \, \varTheta_n \subseteq D_n \, $  for all
$ \, n \in \N \, $.
 \vskip3pt
   $ \underline{n=0} \, $:  \, By definition  $ \, D_0 := \hbox{\sl
Ker}(\delta_1) = \Bbbk \cdot 1_{\scriptscriptstyle \calH} =:
\varTheta_0 \, $,  \, spanned by  $ \,\calA_{\leq 0} =
\{1_{\scriptscriptstyle \calH}\} \, $,  \, q.e.d.
 \vskip3pt
   $ \underline{n=1} \, $:  \, Let  $ \, \eta' \in D_1 := \hbox{\sl
Ker}(\delta_2) \, $.  Let  $ \Bbb{B} $  be a PBW-like basis of  $ \,
{\calH_\h}^{\!\vee} = U(\L_\nu) \, $  as mentioned in Lemma 10.16;
expanding  $ \eta' $  w.r.t.~the basis  $ \calA $  we have  $ \, \eta'
= \sum_{\boldalpha_{\underline{b}} \in \calA} c_{\underline{b}} \,
\boldalpha_{\underline{b}} = \sum_{\underline{b} \in \Bbb{B}}
c_{\underline{b}} \, \boldalpha_{\underline{b}} \, $.  Then we
have also  $ \; \eta \, := \, \eta' \, - \sum_{\tau(\underline{b}\,)
\leq 1} c_{\underline{b}} \, \boldalpha_{\underline{b}}  \, = \,
\sum_{\tau(\underline{b}\,) > 1} c_{\underline{b}} \,
\boldalpha_{\underline{b}} \, \in \, D_1 \; $  because  $ \,
\boldalpha_{\underline{b}} \in \calA_1 \subseteq \varTheta_1
\subseteq D_1 \, $  whenever  $ \, \tau(\underline{b}\,)
\leq 1 \, $.
                                              \par
   Now,  $ \, \boldalpha_1 := \a_1 \, $  and  $ \, \boldalpha_s :=
\a_s - {\a_1}^{\!s} = \h \, \big( \x_s + \h^{s-1} {\x_1}^{\!s} \big)
\, $  for all  $ \, s \in \N_\nu \setminus \{1\} \, $  yield
  $$  \eta \; = \,  
%
%
 {\textstyle \sum_{\underline{b} \in \Bbb{B} \, , \,
\tau(\underline{b}\,) > 1}} 
c_{\underline{b}} \, \boldalpha_{\underline{b}}  \; = \;
%
%
 {\textstyle \sum_{\underline{b} \in \Bbb{B} \, , \,
\tau(\underline{b}\,) > 1}} \, 
\h^{\,g(\underline{b}\,)}
\, c_{\underline{b}} \, \big(
\x_{\underline{b}} + \h \, \chi_{\underline{b}} \big)
\; \in \;  {\calH_\h}^{\!\vee}  $$   
for some  $ \, \chi_{\underline{b}} \in {\calH_\h}^{\!\vee} \, $:
\, hereafter we set  $ \, g(\underline{b}\,) := k \, $  for each
$ \, \underline{b} = b_1 \cdots b_k \in \Bbb{B} \, $  (i.e.~$ \,
g(\underline{b}\,) \, $  is the degree of  $ \underline{b} $  as
a monomial in the  $ b_i $'s).  If  $ \, \eta \not= 0 \, $,  \,
let  $ \, g_0 := \min \big\{\, g(\underline{b}\,) \,\big|\,
\tau(\underline{b}\,) > 1 \, , \, c_{\underline{b}} \not= 0 \,\big\}
\, $;  \, then  $ \, g_0 > 0 \, $,  $ \, \eta_+ := \h^{-g_0} \,
\eta \in {\calH_\h}^{\!\vee} \setminus \h \, {\calH_\h}^{\!\vee}
\, $  and
  $$  0 \;\; \not= \;\; \overline{\,\eta_+} \;\; = \, {\textstyle
\sum\limits_{g(\underline{b}\,) = g_0}} \, c_{\underline{b}}
\, \overline{\x_{\underline{b}}} \;\; = \, {\textstyle
\sum\limits_{g(\underline{b}\,) = g_0}} \, c_{\underline{b}}
\, x_{\underline{b}} \;\; \in \;\, {\calH_\h}^{\!\vee} \Big/
\h \, {\calH_\h}^{\!\vee} \; = \; U(\L_\nu) \; .  $$   
Now  $ \, \delta_2(\eta) = 0 \, $  yields  $ \, \delta_2
\big( \overline{\,\eta_+} \,\big) = 0 \, $,  \, thus  $ \,
\sum_{g(\underline{b}\,) = g_0} \, c_{\underline{b}} \,
x_{\underline{b}} = \overline{\,\eta_+} \in P\big(U(\L_\nu)\big)
= \L_\nu \, $;  \, therefore all PBW monomials occurring in the
last sum do belong to  $ B_\nu $  (and  $ \, g_0 = 1 \, $).  In
addition,  $ \, \delta_2(\eta) = 0 \, $  also implies  $ \,
\delta_2(\eta_+) = 0 \, $  which yields also  $ \, \delta \big(
\overline{\,\eta_+} \,\big) = 0 \, $  for the Lie cobracket
$ \delta $  of  $ \L_\nu $  arising as semiclassical limit of
$ \Delta_{{\calH_\h}^{\!\vee}} $  (see Theorem 10.6); therefore
$ \; \overline{\,\eta_+} \, = \sum_{b \in B_\nu} c_b \, x_b \; $  is
an element of  $ \L_\nu $  killed by the Lie cobracket  $ \delta $,
i.e.~$ \, \overline{\,\eta_+} \in \text{\sl Ker}\,(\delta) \, $.
                                             \par
   Now we apply Lemma 10.15 to  $ \, V = V_\nu \, $,  $ \, \L(V) =
\L(V_\nu) =: \L_\nu \, $  and  $ \, \psi = \delta{\big|}_{V_\nu}
\, $,  \, so that  $ \, \psi_{d\L} = \delta \, $.  From the
formulas for  $ \delta $  in Theorem 10.6 we see that  $ \, K :=
\text{\sl Ker}\,(\psi) = \text{\sl Ker}\,\Big( \delta{\big|}_{V_\nu}
\Big) = \text{\sl Span}\,\big(\{x_1,x_2\}\big) \, $,  \, hence
$ \, \L(K) = \L\big(\text{\sl Span}\,\big(\{x_1,x_2\}\big)\big)
\, $:  \, by definition the last space is nothing but  $ \, \text{\sl
Span}\,\Big( \big\{\, x_b \,\big|\, b \in B_\nu \, ; \, \tau(b) = 1
\,\big\} \Big) \, $,  \, thus eventually via Theorem 10.6 we get
$ \, \text{\sl Ker}\,(\delta) = \L(K) = \text{\sl Span}\,\big(
\big\{\, x_b \,\big|\, b \in B_\nu \, ; \, \tau(b) = 1 \,\big\}
\big) \, $.
                                             \par
   Since  $ \, \overline{\,\eta_+} \in \text{\sl Ker}\,(\delta) =
\text{\sl Span}\,\big( \big\{\, x_b \,\big|\, b \in B_\nu \, ; \,
\tau(b) = 1 \,\big\} \big) \, $  we have  $ \, \overline{\,\eta_+}
= \sum_{\Sb  b \in B_\nu  \\   \tau(b)=1  \endSb} \hskip-3pt c_b \,
x_b \, $;  \, but  $ \, c_b = 0 \, $  whenever  $ \, \tau(b) \leq 1 \, $,
\, by construction of  $ \eta \, $:  thus  $ \, \overline{\,\eta_+} = 0
\, $,  \, a contradiction.  The outcome is  $ \, \eta = 0 \, $,  \,
whence finally  $ \, \eta' \in \varTheta_1 \, $,  \, q.e.d.
 \vskip3pt
   $ \underline{n>1} \, $:  \, We must show that  $ \, D_n = \varTheta_n
\, $,  \, while assuming by induction that  $ \, D_m = \varTheta_m \, $
for all  $ \, m < n \, $.  Let  $ \, \eta = \sum_{\underline{b} \in
\Bbb{B}} c_{\underline{b}} \, \boldalpha_{\underline{b}} \in D_n \, $;
\, then  $ \, \tau(\eta) = \max \big\{\, \tau(\underline{b}\,) \,\big|\,
c_{\underline{b}} \not= 0 \,\big\} \, $.  If  $ \, \delta_2(\eta) =
0 \, $  then  $ \, \eta \in D_1 = \varTheta_1 \, $  by the previous
analysis, and we're done.  Otherwise,  $ \, \delta_2(\eta) \not= 0
\, $  and  $ \, \tau\big(\delta_2(\eta)\big) = \tau(\eta) \, $  by
Lemma 10.16.  On the other hand, since  $ \underline{D} $  is a Hopf
algebra filtration we have  $ \, \delta_2(\eta) \in \sum_{\Sb  r+s=n
\\   r, s > 0  \endSb} D_r \otimes D_s = \sum_{\Sb  r+s=n  \\   r, s
> 0  \endSb} \varTheta_r \otimes \varTheta_s \, $,  \, thanks to
the induction; but then  $ \, \tau\big(\delta_2(\eta)\big) \leq n \, $,
\, by definition of  $ \tau $.  Thus  $ \, \tau(\eta) = \tau \big(
\delta_2(\eta) \big) \leq n \, $,  \, which means  $ \, \eta \in
\varTheta_n \, $.   \qed
\enddemo

\vskip7pt

\proclaim{Theorem 10.18} \, For any  $ \, b \in B_\nu \, $  set
$ \; \widehat{\boldalpha}_b := \h^{\,\kappa(\boldalpha_b)} \,
\boldalpha_b = \h^{\,\tau(b)} \, \boldalpha_b \; $.
                                 \hfill\break
   \indent   (a) \, The set of ordered monomials
  $$  \widehat{\Cal{A}}_{\leq n} \; := \; \Big\{\,
\widehat{\boldalpha}_{\underline{b}} := \widehat{\boldalpha}_{b_1}
\cdots \widehat{\boldalpha}_{b_k} \;\Big|\;  k \in \N \, , \, b_1,
\dots, b_k \in B \, , \, b_1 \preceq \cdots \preceq b_k \, , \,
\kappa\,(\boldalpha_{\underline{b}}\,) = \tau(\underline{b}\,)
\leq n \,\Big\}  $$
is a  $ \, \Bbbk[\h\,] $--basis  of  $ \, D'_n = D_n \big(
{\calH_\h}^{\!\prime} \big) = \h^n D_n \, $.  So  $ \; \widehat{\Cal{A}}
:= \bigcup_{n \in \N} \widehat{\Cal{A}}_{\leq n} \; $  is a  $ \,
\Bbbk[\h\,] $--basis  of  $ \, {\calH_\h}^{\!\prime} \, $.
                                 \hfill\break
   \indent   (b)  $ \; \displaystyle{ {\calH_\h}^{\!\prime}  \;
= \;  \Bbbk[\h\,] \, \Big\langle {\big\{\, \widehat{\boldalpha}_b
\,\big\}}_{b \in B_\nu} \Big\rangle \! \Bigg/ \!\! \bigg( \! \Big\{\,
\big[ \widehat{\boldalpha}_{b_1}, \widehat{\boldalpha}_{b_2} \big]
- \h \, \widehat{\boldalpha}_{[b_1,b_2]} \;\Big|\; \forall \; b_1,
b_2 \in B_\nu \,\Big\} \bigg) } \, $.
                                 \hfill\break
  \indent   (c) \,  $ {\calH_\h}^{\!\prime} $  is a graded Hopf
$ \; \Bbbk[\h\,] $--subalgebra  of  $ \, {\calH_\h} \, $.
                                        \hfill\break
  \indent   (d) \,  $ {\calH_\h}^{\!\prime}{\Big|}_{\h=0}
:= {\calH_\h}^{\!\prime} \Big/ \h \, {\calH_\h}^{\!\prime}
= \widetilde{\calH} = F \big[ {\varGamma_{\!\L_\nu}
\phantom{|}}^{\hskip-8pt \star} \big] \, $,  \, where
$ \, {\varGamma_{\!\L_\nu}\phantom{|}}^{\hskip-8pt\star} \, $
is a connected Poisson algebraic group with cotangent Lie bialgebra
isomorphic to  $ \L_\nu $  (as a Lie algebra) with the graded Lie
bialgebra structure given by  $ \, \delta(x_n) = (n-2) \, x_{n-1}
\wedge x_1 \, $  (for all  $ \, n \in \N_\nu $).  Indeed,  $ \,
{\calH_\h}^{\!\prime}{\Big|}_{\h=0} $  is the free Poisson
(commutative) algebra over  $ \N_\nu \, $,  generated by all
the\/  $ \bar{\boldalpha}_n := \widehat{\boldalpha}_n{\big|}_{\h=0} $
($ \, n \in \N_\nu \, $)  with Hopf structure given (for all  $ \,
n \in \N_\nu $)  by
  $$  \displaylines{
   \Delta\big(\bar{\boldalpha}_n\big)  \; = \;\,
\bar{\boldalpha}_n \otimes 1 \, + \, 1 \otimes
\bar{\boldalpha}_n \, + \, \sum_{k=2}^{n-1} {n \choose k} \,
\bar{\boldalpha}_k \otimes {\bar{\boldalpha}_1}^{\,n-k} \,
+ \, \sum_{k=1}^{n-1} \, (k+1) \, {\bar{\boldalpha}_1}^{\,k}
\otimes \bar{\boldalpha}_{n-k}  \cr
   S\big(\bar{\boldalpha}_n\big)  \; = \;\,
- \, \bar{\boldalpha}_n \, - \, \sum_{k=2}^{n-1} {n \choose k} \,
S\big(\bar{\boldalpha}_k\big) \, {\bar{\boldalpha}_1}^{\,n-k}
\, - \, \sum_{k=1}^{n-1} \, (k+1) \, S \big( \bar{\boldalpha}_1
\big)^k \, \bar{\boldalpha}_{n-k} \; ,  \qquad
\epsilon\big(\bar{\boldalpha}_n\big) \, = \, 0  \; .  \cr }  $$
Thus  $ \, {\calH_\h}^{\!\prime}{\Big|}_{\h=0} $  is the polynomial
algebra  $ \, \Bbbk \big[ {\{\, \eta_b \,\}}_{b \in B_\nu} \big] \, $
generated by a set of indeterminates  $ \, {\{\, \eta_b \,\}}_{b \in
B_\nu} \, $  in bijection with  $ B_\nu \, $,  so  $ \; {\varGamma_{\!
\L_\nu} \phantom{|}}^{\hskip-8pt \star} \cong \Bbb{A}_\Bbbk^{B_\nu}
\, $
%
%
as algebraic varieties.
                                        \hfill\break
  \indent   Finally,  $ \, {\calH_\h}^{\!\prime}{\Big|}_{\h=0} =
F \big[ {\varGamma_{\!\L_\nu} \phantom{|}}^{\hskip-8pt \star} \big]
= \Bbbk \big[ {\{\, \eta_b \,\}}_{b \in B_\nu} \big] \, $  is a
{\sl graded Poisson Hopf algebra}  w.r.t.~the grading
$ \, \partial(\bar{\boldalpha}_n) = n \, $  (inherited from
$ {\calH_\h}^{\!\prime} $)  and w.r.t.~the grading induced
from  $ \, \kappa = \tau \, $  (on  $ \calH $),  and a  {\sl
graded algebra}  w.r.t.~the  ({\sl polynomial})  grading  $ \,
d(\bar{\boldalpha}_n) = 1 \, $  (for all  $ \, n \in \N_+ $).
                                        \hfill\break
  \indent   (e) \,  The analogues of statements (a)--(d) hold with
$ \calK $  instead of  $ \, \calH \, $,  \, with  $ X^+ $  instead
of  $ X $  for all  $ \, X = \L_\nu, B_\nu, \N_\nu \, $,  and with
$ {\varGamma_{\!\L_\nu^+}\phantom{|}}^{\hskip-8pt \star} $  instead
of  $ {\varGamma_{\!\L_\nu}\phantom{|}}^{\hskip-8pt \star} \, $.
\endproclaim

\demo{Proof} \, {\it (a)} \, This follows from Proposition 10.17 and
the characterization of  $ {\calH_\h}^{\!\prime} $  in \S 10.10.
                                              \par
   {\it (b)} \, This is a direct consequence of claim  {\it (a)\/}
and  Lemma 10.14{\it (g)}.
                                              \par
   {\it (c)} \, Thanks to claims  {\it (a)\/}  and  {\it (b)},  we
can look at  $ {\calH_\h}^{\!\prime} $  as a Poisson algebra, whose
Poisson bracket is given by  $ \, \{x,y\,\} := \h^{-1} [x,y] = \h^{-1}
(x \, y - y \, x) \, $  (for all  $ \, x $,  $ y \in {\calH_\h}^{\!
\prime} \, $);  then  
   $ {\calH_\h}^{\!\prime} $  itself\break
 \eject
\noindent   
 is the free associative Poisson algebra generated by  $ \big\{\,
\widehat{\boldalpha}_n \,\big|\, n \in \N \,\big\} $.  Clearly
$ \Delta $  is a Poisson map, therefore it is enough to prove that
$ \, \Delta\big(\widehat{\boldalpha}_n\big) \in {\calH_\h}^{\!\prime}
\otimes {\calH_\h}^{\!\prime} \, $  for all  $ \, n \in \N_+ \, $.
This is clear for  $ \boldalpha_1 $  and  $ \boldalpha_2 $  which
are primitive; as for  $ \, n > 2 \, $,  \, we have, like in
Proposition 10.13,
  $$  \hbox{ $ \eqalign{
   \hskip-2pt   \Delta \big( \widehat{\boldalpha}_n  &  \big)  \;
= \;  {\textstyle \sum_{k=2}^n} \, \h^{k-1} \boldalpha_k \otimes
\h^{n-k} Q^k_{n-k}(\a_*) \, + \, {\textstyle \sum_{k=0}^{n-1}} \, \h^k
{\boldalpha_1}^{\!k} \otimes \h^{n-k-1} Z^k_{n-k}(\boldalpha_*)
\; =   \hfill  \cr
   {} \hfill   &  \hskip-3pt = \;  {\textstyle \sum_{k=2}^n} \,
\widehat{\boldalpha}_k \otimes \h^{n-k} Q^k_{n-k}(\a_*) \, + \,
{\textstyle \sum_{k=0}^{n-1}} \, {\widehat{\boldalpha}_1}^{\,k} \otimes
\h^{n-k-1} Z^k_{n-k}(\boldalpha_*) \; \in \; {\calH_\h}^{\!\prime}
\otimes {\calH_\h}^{\!\prime}  \cr } $ }   \hskip-1pt (10.8)  $$  
thanks to Lemma 10.12 (with notations used therein).  In addition,
$ \, S\big({\calH_\h}^{\!\prime}\big) \subseteq {\calH_\h}^{\!\prime}
               \, $  also\break
%
%
%
\noindent
follows by induction from (10.8) because Hopf algebra
axioms along with (10.8) give
  $$  S\big(\widehat{\boldalpha}_n\big)  \; = \;
- \widehat{\boldalpha}_n \, - \, {\textstyle \sum_{k=2}^{n-1}} \,
S\big(\widehat{\boldalpha}_k\big) \, \h^{n-k} Q^k_{n-k}(\a_*) \, - \,
{\textstyle \sum_{k=1}^{n-1}} \, S \big( {\widehat{\boldalpha}_1}^{\,k}
\big) \, \h^{n-k-1} Z^k_{n-k}(\boldalpha_*) \; \in \; {\calH_\h}^{\!
\prime}  $$
for all  $ \, n \in \N_\nu \, $  (using induction).  The claim
follows.
                                              \par
   {\it (d)} \, Thanks to  {\it (a)\/}  and  {\it (b)},  $ \,
{\calH_\h}^{\!\prime} {\Big|}_{\h=0} $  is a polynomial
$ \Bbbk $--algebra  as claimed, over the set of indeterminates
$ \big\{\, \bar{\boldalpha}_b := \widehat{\boldalpha}_b
{\big|}_{\h=0} \big( \in {\calH_\h}^{\!\prime}{\big|}_{\h=0} \big)
\,\big|\, b \in B_\nu \,\big\} $.  Furthermore, in the proof of
{\it (c)\/}  we noticed that  $ {\calH_\h}^{\!\prime} $  is also
the free Poisson algebra generated by  $ \big\{\, \widehat{\boldalpha}_n
\,\big|\, n \in \N \,\big\} $;  therefore  $ \, {\calH_\h}^{\!\prime}
{\Big|}_{\h=0} $  is the free commutative Poisson algebra generated
by  $ \big\{\, \bar{\boldalpha}_n := \check{\boldalpha}_{x_n}
{\big|}_{\h=0} \,\big|\, n \in \N \,\big\} $.  Then formula (10.8)
--- for all  $ \, n \in \N_\nu \, $  ---   describes uniquely the
Hopf structure of  $ {\calH_\h}^{\!\prime} $,  hence the formula
it yields at  $ \, \h = 0 \, $  will describe the Hopf structure
of  $ \, {\calH_\h}^{\!\prime}{\big|}_{\h=0} $.
                                              \par
   Expanding  $ \, \h^{n-k} Q^k_{n-k}(\a_*) \, $  in (10.8) w.r.t.~the
basis  $ \widehat{\calA} $  in  {\it (a)\/}  we find a sum of terms
of  $ \tau $--degree  less or equal than  $ (n-k) $,  and the sole
one achieving equality is  $ \, {\widehat{\boldalpha}_1}^{\,n-k} \, $,
\, which occurs with coefficient  $ {n \choose k} $:  \, similarly,
when expanding  $ \, \h^{n-k-1} Z^k_{n-k}(\boldalpha_*) \, $  in (10.8)
w.r.t.~$ \widehat{\calA} $  all summands have  $ \tau $--degree  less
or equal than  $ (n-k-1) $,  and equality holds only for  $ \,
\widehat{\boldalpha}_{n-k} \, $,  \, whose coefficient is  $ \,
(k+1) \, $.  Therefore for some  $ \, \boldsymbol{\eta} \in
{\calH_\h}^{\!\prime} {\big|}_{\h=0} \; $  we have
  $$  \Delta\big(\widehat{\boldalpha}_n\big) \; = \;
{\textstyle \sum_{k=2}^n} \, \widehat{\boldalpha}_k \otimes
{n \choose k} \, {\widehat{\boldalpha}_1}^{\,n-k} \, + \, {\textstyle
\sum_{k=0}^{n-1}} \, (k+1) \, {\widehat{\boldalpha}_1}^{\,k} \otimes
\widehat{\boldalpha}_{n-k} \, + \, \h \; \boldsymbol{\eta} \; ;  $$
this yields the formula for  $ \Delta $,  from which the formula
for  $ S $  follows too as usual.
                                              \par
   Finally, let  $ \, \varGamma := \text{\sl Spec}\,\big(
{\calH_\h}^{\!\prime}{\big|}_{\h=0} \big) \, $  be the
algebraic Poisson group such that  $ \, F\big[\varGamma\big]
= {\calH_\h}^{\!\prime}{\big|}_{\h=0} \, $,  \, and let  $ \,
\boldsymbol{\gamma}_\nu := \text{\sl co$\L$ie}\,(\varGamma) \, $
be its cotangent Lie bialgebra.  Since  $ \, {\calH_\h}^{\!
\prime}{\big|}_{\h=0} $  is Poisson free over  $ \big\{
\bar{\boldalpha}_n \big\}_{n \in \N_\nu} \, $,  \,
as a Lie algebra  $ \boldsymbol{\gamma}_\nu $  is free over
$ \, \big\{\, d_n := \bar{\boldalpha}_n \mod \germ^2
\,\big\}_{n \in \N_\nu} \, $  (where  $ \, \germ :=
J_{{\calH_\h}^{\!\prime}{|}_{\h=0}} \, $),  so  $ \,
\boldsymbol{\gamma}_\nu \cong \L_\nu \, $,  via  $ \,
d_n \mapsto x_n (n \in \N_+) \, $  as a Lie algebra.
The Lie cobracket is   
  $$  \displaylines{
   {} \quad   \delta_{\boldsymbol{\gamma}_\nu}\big(d_n\big)  \,
= \,  (\Delta - \Delta^{\text{op}})\big(\bar{\boldalpha}_n\big)
\mod \germ_\otimes  \; =   \hfill  \cr
   {} \hfill   = \;\,  {\textstyle \sum\limits_{k=2}^{n-1}} {n \choose
k} \, \bar{\boldalpha}_k \wedge {\bar{\boldalpha}_1}^{\,n-k}
\, + \, {\textstyle \sum\limits_{k=1}^{n-1}} \, (k+1) \,
{\bar{\boldalpha}_1}^{\,k} \wedge \bar{\boldalpha}_{n-k}
\mod \germ_\otimes  \; =  \cr
   {} \quad   = \;  {n \choose n-1} \, \bar{\boldalpha}_{n-1}
\wedge \bar{\boldalpha}_1 \, + \, 2 \, \bar{\boldalpha}_1 \wedge
\bar{\boldalpha}_{n-1}  \mod \germ_\otimes  \; = \;   \hfill {}  \cr
   {} \hfill   = \;  (n-2) \, \bar{\boldalpha}_{n-1} \wedge
\bar{\boldalpha}_1  \mod \germ_\otimes  \; = \; (n-2) \,
d_{n-1} \wedge d_1  \; \in \;  \boldsymbol{\gamma} \otimes
\boldsymbol{\gamma}  \cr }  $$
 where  $ \, \germ_\otimes := \Big( \germ^2 \otimes {\calH_\h}^{\!
\prime}{|}_{\h=0} + \germ \otimes \germ + {\calH_\h}^{\!\prime}
{|}_{\h=0} \otimes \germ^2 \Big) \, $,  \, whence  $ \, \varGamma
= {\varGamma_{\!\L_\nu}\phantom{|}}^{\hskip-8pt \star} \, $  as
claimed in  {\it (d)}.
                                              \par
   Finally, the statements about gradings of  $ \, {\calH_\h}^{\!\prime}
{\Big|}_{\h=0} = F \big[ {\varGamma_{\!\L_\nu} \phantom{|}}^{\hskip-8pt
\star} \big] \, $  hold by construction.
%
%
 \eject
   {\it (e)} \, This should be clear from the whole discussion, since
all arguments apply again   --- {\sl mutatis mutandis\/} ---   when
starting with  $ \calK $  instead of  $ \, \calH \, $;  \, we leave
details to the reader.   \qed
\enddemo

\vskip7pt

  {\bf 10.19 Drinfeld's algebra  $ \big( {\calH_\h}^{\!\prime}
\big)^{\!\vee} $.} \, I look now at the other Drinfeld's functor at
$ \h \, $, and consider  $ \; \big( {\calH_\h}^{\!\prime} \big)^{\!
\vee} := \sum_{n \in \N} \h^{-n} {J^{\,\prime}}^n \, $,  \; where
$ \, J^{\,\prime} := J_{{\calH_\h}^{\!\prime}} \, $.  Theorem 2.2
tells us that  {\it  $ \big( {\calH_\h}^{\!\prime} \big)^{\!\vee} $
is a Hopf\/  $ \Bbbk[\h\,] $--subalgebra  of  $ \calH_\h $,
and the specialization of  $ \, \big( {\calH_\h}^{\!\prime}
\big)^{\!\vee} $  at  $ \, \h = 0 \, $,  \, i.e.~$ \; \big(
{\calH_\h}^{\!\prime} \big)^{\!\vee}{\Big|}_{\h=0} := \, \big(
{\calH_\h}^{\!\prime} \big)^{\!\vee} \! \Big/ \h \, \big( {\calH_\h}^{\!
\prime} \big)^{\!\vee} \, $,  \; is the universal enveloping algebra
of the cotangent Lie bialgebra of the connected algebraic Poisson
group which is the spectrum of  $ \, {\calH_\h}^{\!\prime}
{\Big|}_{\h=0} $,  \, that is exactly}  $ {\varGamma_{\!\L_\nu}
\phantom{|}}^{\hskip-8pt \star} \, $.  Thanks to Theorem 10.18,
this means  $ \, \big( {\calH_\h}^{\!\prime} \big)^{\!\vee}
{\Big|}_{\h=0} = U(\L_\nu) \, $  as co-Poisson Hopf
$ \Bbbk $--algebras,  the Lie cobracket of  $ \L_\nu $
being the one given in  Theorem 10.18{\it (d)}.
                                             \par
   Therefore we must show that  $ \big( {\calH_\h}^{\!\prime}
\big)^{\!\vee}{\Big|}_{\h=0} $  is a cocommutative Hopf
$ \Bbbk $--algebra,  it is generated by its primitive elements,
and the latter set inherits a Lie bialgebra structure isomorphic
to that of  $ \, \boldsymbol{\gamma}_\nu := \text{\sl co$\L$ie}\,
\big({\varGamma_{\!\L_\nu}\phantom{|}}^{\hskip-8pt \star}\big) \, $.
We prove all this directly, via an explicit description of  $ \big(
{\calH_\h}^{\!\prime} \big)^{\!\vee} $  and its specialization at
$ \, \h = 0 \, $,  \, provided in the following

\vskip7pt

\proclaim{Theorem 10.20} \, For any  $ \, b \in B_\nu \, $  set
$ \; \check{\boldalpha}_b := \h^{\,\kappa(\boldalpha_b)-1}
\, \boldalpha_b = \h^{\,\tau(b)-1} \, \boldalpha_b = \h^{-1}
\, \widehat{\boldalpha}_b \; $.
                                 \hfill\break
   \indent   (a)  $ \; \displaystyle{ {\big( {\calH_\h}^{\!\prime}
\big)}^{\!\vee}  \; = \;  \Bbbk[\h\,] \, \Big\langle {\big\{\,
\check{\boldalpha}_b \,\big\}}_{b \in B_\nu} \Big\rangle \!
\Bigg/ \!\! \bigg( \! \Big\{\, \big[ \check{\boldalpha}_{b_1},
\check{\boldalpha}_{b_2} \big] - \, \check{\boldalpha}_{[b_1,b_2]}
\;\Big|\; \forall \; b_1, b_2 \in B_\nu \,\Big\} \bigg) } \, $.
                                 \hfill\break
  \indent   (b) \,  $ {\big( {\calH_\h}^{\!\prime} \big)}^{\!\vee} $
is a graded Hopf  $ \; \Bbbk[\h\,] $--subalgebra  of  $ \, \calH_\h
\, $.
                                        \hfill\break
  \indent   (c) \,  $ {\big( {\calH_\h}^{\!\prime} \big)}^{\!\vee}
{\Big|}_{\h=0} := {\big( {\calH_\h}^{\!\prime} \big)}^{\!\vee} \Big/
\h \, {\big( {\calH_\h}^{\!\prime} \big)}^{\!\vee} \cong U \big( \L_\nu
\big) \; $  as co-Poisson Hopf algebra, where  $ \, \L_\nu \, $  bears
the Lie bialgebra structure given by  $ \, \delta(x_n) = (n-2) \,
x_{n-1} \wedge x_1 \, $  (for all  $ \, n \in \N_\nu $).
                                        \hfill\break
   \indent   Finally, the grading  $ \, d $  given by  $ \, d(x_n)
:= 1 \;\, (n \in \N_+) \, $  makes  $ \, {\big( {\calH_\h}^{\!\prime}
\big)}^{\!\vee}{\Big|}_{\h=0} \! = U(\L_\nu) \, $  into a  {\sl graded
co-Poisson Hopf algebra},  and the grading  $ \, \partial $  given
by  $ \, \partial(x_n) := n \;\, (n \in \N_+) \, $  makes  $ \, {\big(
{\calH_\h}^{\!\prime} \big)}^{\!\vee}{\Big|}_{\h=0} \! = U(\L_\nu)
\, $     into a  {\sl graded Hopf algebra}  and  $ \L_\nu $  into
a  {\sl graded Lie bialgebra.}
                                        \hfill\break
  \indent   (d) \,  The analogues of statements (a)--(c) hold
with  $ \calK \, $,  $ \L_\nu^+ $,  $ B_\nu^+ $  and  $ \N_\nu^+ $
respectively instead of  $ \, \calH \, $,  $ \L_\nu^+ $,  $ B_\nu $
and  $ \N_\nu^+ \, $.
\endproclaim

\demo{Proof} \, {\it (a)} \, This follows from  Theorem 10.18{\it
(b)\/}  and the very definition of  $ {\big( {\calH_\h}^{\!\prime}
\big)}^{\!\vee} $  in \S 10.19.
                                              \par
   {\it (b)} \, This is a direct consequence of claim  {\it (a)\/}
and  Theorem 10.18{\it (c)}.
                                              \par
   {\it (c)} \, It follows from claim  {\it (a)\/}  that mapping
$ \, \check{\boldalpha}_b{\big|}_{\h=0} \mapsto b \, $  ($ \forall
\, b \in B_\nu \, $)  yields a well-defined algebra isomorphism
$ \, \Phi \, \colon \, {\big( {\calH_\h}^{\!\prime}
\big)}^{\!\vee}{\Big|}_{\h=0} {\buildrel \cong \over
{\lhook\joinrel\relbar\joinrel\relbar\joinrel\twoheadrightarrow\,}}
U\big(\L_\nu) \, $.  In addition, when expanding  $ \, \h^{n-k}
Q^k_{n-k}(\a_*) \, $  in (10.8) w.r.t.~the basis  $ \calA $  (see
Proposition 10.17) we find a sum of terms of  $ \tau $--degree  less
than or equal to  $ (n-k) $,  and equality is achieved only for 
$ \, \boldalpha_1^{\,n-k} \, $,  \, which occurs with coefficient
$ {n \choose k} $:  \, similarly, the expansion of  $ \, \h^{n-k-1}
Z^k_{n-k}(\boldalpha_*) \, $  in (10.8) yields a sum of terms whose
$ \tau $--degree  is less or equal than  $ (n-k-1) $,  with equality
only for  $ \, \boldalpha_{n-k} \, $,  \, whose coefficient is  $ \,
(k+1) \, $.   Thus using the relation  $ \, \widehat{\boldalpha}_s
= \h \, \check{\boldalpha}_s \, $  ($ \, s \in \N_+ \, $)  we get
  $$  \displaylines{
   \Delta \big( \check{\boldalpha}_n \big)  \, = \,
\check{\boldalpha}_n \otimes 1  \, + \,  1 \otimes \check{\boldalpha}_n
\, + \,  {\textstyle \sum_{k=2}^{n-1}} \, \check{\boldalpha}_k \otimes
\h^{n-k} Q^k_{n-k}(\a_*)  \, + \, {\textstyle \sum_{k=1}^{n-1}} \,
{\check{\boldalpha}_1}^{\,k} \otimes \h^{n-1} Z^k_{n-k}(\boldalpha_*)
\, =   \hfill {}  \cr
   = \,  \check{\boldalpha}_n \otimes 1  \, + \,  1 \otimes
\check{\boldalpha}_n  \, + \, {\textstyle \sum_{k=2}^{n-1}} \,
\h^{n-k} \, \check{\boldalpha}_k \otimes {n \choose k} \,
{\check{\boldalpha}_1}^{\,n-k}  \, + \, {\textstyle \sum_{k=1}^{n-1}}
\, \h^k \, (k+1) \, {\check{\boldalpha}_1}^{\,k} \otimes
\check{\boldalpha}_{n-k}  \, + \,  \h^2 \; \boldsymbol{\eta} \, =  \cr
   {} \hfill   = \, \check{\boldalpha}_n \otimes 1  \, + \,  1 \otimes
\check{\boldalpha}_n  \, + \,  \h \, \big( n \, \check{\boldalpha}_{n-1}
\otimes \check{\boldalpha}_1  \, + \, 2 \, \check{\boldalpha}_1 \otimes
\check{\boldalpha}_{n-1} \big)  \, + \,  \h^2 \; \boldsymbol{\chi}
\cr }  $$
for some  $ \, \boldsymbol{\eta}, \boldsymbol{\chi} \in {\big(
{\calH_\h}^{\!\prime} \big)}^{\!\vee} \otimes {\big( {\calH_\h}^{\!
\prime} \big)}^{\!\vee} \, $.  It follows that  $ \, \Delta
\big( \check{\boldalpha}_n{\big|}_{\h=0} \big)  \, = \,
\check{\boldalpha}_n{\big|}_{\h=0} \otimes 1  \, + \,
1 \otimes \check{\boldalpha}_n{\big|}_{\h=0} \, $  for
all  $ \, n \in\N_\nu \, $.  Similarly we have  $ \, S \big(
\check{\boldalpha}_n{\big|}_{\h=0} \big) = - \check{\boldalpha}_n
{\big|}_{\h=0} \, $  and  $ \, \epsilon\big(\check{\boldalpha}_n
{\big|}_{\h=0}\big) = 0 \, $  for all  $ \, n \in \N_\nu \, $,
\, thus  $ \Phi $  {\sl is an isomorphism of Hopf algebras\/}  too.
In addition, the Poisson cobracket of  $ {\big( {\calH_\h}^{\!\prime}
\big)}^{\!\vee}{\Big|}_{\h=0} $  inherited from  $ {\big( {\calH_\h}^{\!
\prime} \big)}^{\!\vee} $  is given by
  $$  \displaylines{
   \delta\big(\check{\boldalpha}_n{\big|}_{\h=0}\big)
\, = \,  \Big( \h^{-1} (\Delta - \Delta^{\text{op}}) \big(
\check{\boldalpha}_n \big) \Big)  \mod \h \, {\big( {\calH_\h}^{\!
\prime} \big)}^{\!\vee} \otimes {\big( {\calH_\h}^{\!\prime}
\big)}^{\!\vee}  \; =   \hfill  \cr
   {} \hfill   = \;\,  \big( n \, \check{\boldalpha}_{n-1} \wedge
\check{\boldalpha}_1  \, + \, 2 \, \check{\boldalpha}_1 \wedge
\check{\boldalpha}_{n-1} \big)  \mod \h \, {\big( {\calH_\h}^{\!
\prime} \big)}^{\!\vee} \! \otimes {\big( {\calH_\h}^{\!\prime}
\big)}^{\!\vee}  \, = \;  (n-2) \, \check{\boldalpha}_{n-1}
{\big|}_{\h=0} \wedge \check{\boldalpha}_1{\big|}_{\h=0}  \cr }  $$
hence  $ \Phi $  is also an isomorphism of  {\sl co-Poisson\/}
Hopf algebras, as claimed.
                                              \par
   The statements on gradings of  $ \, {\big( {\calH_\h}^{\!\prime}
\big)}^{\!\vee}{\Big|}_{\h=0} = U(\L_\nu) \, $  should be clear by
construction.
                                              \par
   {\it (d)} \, This should be clear from the whole discussion, as
all arguments apply again   --- {\sl mutatis mutandis\/} ---   when
starting with  $ \calK $  instead of  $ \, \calH \, $;  \, details
are left to the reader.   \qed
\enddemo

\vskip7pt

   {\bf 10.21 Specialization limits.} \, So far, Theorem 10.18{\it
(d)\/}  and  Theorem 10.20{\it (c)\/}  prove the following
specialization results for  $ {\calH_\h}^{\!\prime} $  and
$ {\big( {\calH_\h}^{\!\prime} \big)}^{\!\vee} $  respectively:
  $$  {\calH_\h}^{\!\prime} \;{\buildrel \h \rightarrow 0
\over \llongrightarrow}\; \widetilde{\calH} \cong F \big[
{\varGamma_{\!\L_\nu}\phantom{|}}^{\hskip-8pt\star} \big] \; ,
\qquad \qquad  {\big( {\calH_\h}^{\!\prime} \big)}^{\!\vee}
\;{\buildrel \h \rightarrow 0 \over \llongrightarrow}\;
U(\L_\nu)  $$
as graded Poisson or co-Poisson Hopf  $ \Bbbk $--algebras.
In addition,  Theorem 10.18{\it (b)\/}  implies that  $ \;
{\calH_\h}^{\!\prime} \,{\buildrel \h \rightarrow 1 \over
\llongrightarrow}\, \calH' = \calH \; $  as graded Hopf
$ \Bbbk $--algebras.  Indeed, by  Theorem 10.18{\it (b)\/}
$ \, \calH $  (or even  $ \calH_\h $)  embeds as an algebra
into  $ {\calH_\h}^{\!\prime} $,  via  $ \; \boldalpha_n
\mapsto \widehat\boldalpha_n \; $  (for all  $ \, n \in
\N_\nu \, $):  \, then
  $$  [\boldalpha_n, \boldalpha_m] \, \mapsto \, \big[
\widehat\boldalpha_n, \widehat\boldalpha_m \big]  \, =
\,  \h \, \widehat{\boldalpha}_{[x_n,x_m]} \, \equiv \,
\widehat{\boldalpha}_{[x_n,x_m]} \mod (\h \! - \! 1) \,
{\calH_\h}^{\!\prime}   \eqno \big(\, \forall \; n, m
\in \N_\nu \big)  $$
thus, thanks to the presentation of  $ {\calH_\h}^{\!\prime} $  by
generators and relations in  Theorem 10.18{\it (b)\/},  $ \calH $
is isomorphic to  $ \; {\calH_\h}^{\!\prime}{\Big|}_{\h=1}  \; := \;
{\calH_\h}^{\!\prime} \Big/ (\h \! - \! 1) \, {\calH_\h}^{\!\prime}
\; = \;  \Bbbk \big\langle {\widehat\boldalpha}_1{\big|}_{\h=1},
{\widehat\boldalpha}_2{\big|}_{\h=1}, \dots, {\widehat\boldalpha}_n
{\big|}_{\h=1}, \ldots \big\rangle \, $,  \, as a  $ \Bbbk $--algebra,
via  $ \, \boldalpha_n \mapsto {\widehat\boldalpha}_n{\big|}_{\h=1}
\, $.  Moreover, the Hopf structure of  $ {\calH_\h}^{\!\prime}
{\Big|}_{\h=1} $  is given by
  $$  \Delta \big( \widehat{\boldalpha}_n{\big|}_{\h=1} \big)
=  {\textstyle \sum_{k=2}^n} \, \widehat{\boldalpha}_k \otimes
\h^{n-k} Q^k_{n-k}(\a_*)  +  {\textstyle \sum_{k=0}^{n-1}} \,
{\widehat{\boldalpha}_1}^{\,k} \otimes \h^{n-1} Z^k_{n-k}
(\boldalpha_*) \hskip-6pt  \mod (\h-1) {\calH_\h}^{\!\prime}
\otimes {\calH_\h}^{\!\prime} \, .  $$
   \indent   Now,  $ \, Q^k_{n-k}(\a_*) = Q^k_{n-k}(\boldalpha_*
+ {\boldalpha_1}^*) = {\Cal Q}^k_{n-k}(\boldalpha_*) \, $  for
some polynomial  $ {\Cal Q}^k_{n-k}(\boldalpha_*) $  in the
$ \boldalpha_i $'s;   let  $ \, {\Cal Q}^k_{n-k}(\boldalpha_*)
= \sum_s {\Cal T}^{s,k}_{n-k}(\boldalpha_*) \, $  be the splitting
of  $ {\Cal Q}^k_{n-k} $  into  $ \tau $--homogeneous  summands
(i.e., each  $ {\Cal T}^{s,k}_{n-k}(\boldalpha_*) $  is a
homogeneous polynomial of  $ \tau $--degree  $ s \, $):  then
  $$  \h^{n-k} Q^k_{n-k}(\a_*)  \, = \,  \h^{n-k} {\Cal Q}^k_{n-k}
(\boldalpha_*)  \, = \,  \h^{n-k} {\textstyle \sum_s}
{\Cal T}^{s,k}_{n-k}(\boldalpha_*)  \, = \,  {\textstyle \sum_s}
\h^{n-k-s} {\Cal T}^{s,k}_{n-k}(\widehat{\boldalpha}_*)  $$  
with  $ \, n-k-s > 0 \, $  for all  $ s $  (by construction).  Since
clearly  $ \, \h^{n-k-s}{\Cal T}^{s,k}_{n-k}(\widehat{\boldalpha}_*)
\equiv {\Cal T}^{s,k}_{n-k}(\widehat{\boldalpha}_*) \mod (\h-1) \,
{\calH_\h}^{\!\prime} \, $,  \, we find  $ \; \h^{n-k} \, Q^k_{n-k}
(\a_*) \, = \, \h^{n-k} \, {\Cal Q}^k_{n-k} (\boldalpha_*) \,
= \, {\textstyle \sum_s} \h^{n-k-s} \, {\Cal T}^{s,k}_{n-k}
(\widehat{\boldalpha}_*) \, \equiv \, {\textstyle \sum_s}
{\Cal T}^{s,k}_{n-k}(\widehat{\boldalpha}_*)  \mod (\h-1) \,
{\calH_\h}^{\!\prime} \, = \, {\Cal Q}^k_{n-k}(\widehat{\boldalpha}_*)
\, $,  \, for all  $ k $  and  $ n \, $.  Similarly we deduce
that  $ \, \h^{n-1} Z^k_{n-k}(\boldalpha_*) \equiv Z^k_{n-k}
(\widehat{\boldalpha}_*) \mod (\h-1) \, {\calH_\h}^{\!\prime} \, $,
\, for all  $ k $  and  $ n \, $.  The outcome is that
  $$  \displaylines{
   \Delta \big( \widehat{\boldalpha}_n{\big|}_{\h=1} \big)  \! = \!
{\textstyle \sum_{k=2}^n} \, \widehat{\boldalpha}_k \otimes \h^{n-k}
\! {\Cal Q}^k_{n-k}(\boldalpha_*)  +  {\textstyle \sum_{k=0}^{n-1}} \,
{\widehat{\boldalpha}_1}^{\,k} \otimes \h^{n-1} Z^k_{n-k}(\boldalpha_*)
\hskip-4pt  \mod \hskip-2pt (\h-1) \, {\calH_\h}^{\!\prime} \otimes
{\calH_\h}^{\!\prime}  =  \cr
   {} \hfill   = \,  {\textstyle \sum_{k=2}^n} \,
\widehat{\boldalpha}_k \otimes {\Cal Q}^k_{n-k}(\widehat{\boldalpha}_*)
+  {\textstyle \sum_{k=0}^{n-1}} \, {\widehat{\boldalpha}_1}^{\,k}
\otimes Z^k_{n-k}(\widehat{\boldalpha}_*) \mod (\h-1) \,
{\calH_\h}^{\!\prime} \otimes {\calH_\h}^{\!\prime} \, .  \cr }  $$
   \indent   On the other hand, we have  $ \; \Delta(\boldalpha_n)
\, = \,  {\textstyle \sum_{k=2}^n} \, \boldalpha_k \otimes {\Cal
Q}^k_{n-k}(\boldalpha_*) \, + \, {\textstyle \sum_{k=0}^{n-1}}
\, \boldalpha_1^{\,k} \otimes Z^k_{n-k}(\boldalpha_*) \, $  in 
$ \calH $.  Thus the  {\sl graded algebra\/}   isomorphism 
$ \, \Psi \, \colon \, \calH \,{\buildrel \cong \over
{\lhook\joinrel\relbar\joinrel\relbar\joinrel\twoheadrightarrow}}\,
{\calH_\h}^{\!\prime}{\Big|}_{\h=1} \, $  given by  $ \,
\boldalpha_n \mapsto {\widehat\boldalpha}_n{\big|}_{\h=1} \, $
preserves the coproduct too.  Similarly,  $ \Psi $  respects
the antipode and the counit, hence it is a graded Hopf algebra
isomorphism.  In a nutshell, we have (as graded Hopf
$ \Bbbk $--algebras)
  $$  {\calH_\h}^{\!\prime} \;{\buildrel {\h \rightarrow 1}
\over \llongrightarrow}\; \calH'= \calH \; .  $$
   \indent   Similarly, Theorem 10.20 implies that  $ \; {\big(
{\calH_\h}^{\!\prime} \big)}^{\!\vee} \,{\buildrel \h \rightarrow
1 \over \llongrightarrow}\, \calH \, $  as graded Hopf
$ \Bbbk $--algebras.  Indeed,  Theorem 10.20{\it (a)\/}  shows
that  $ \, {\big( {\calH_\h}^{\!\prime} \big)}^{\!\vee} \cong
\Bbbk[\h\,] \otimes_\Bbbk U(\L_\nu) \, $  {\sl as graded associative
algebras},  via  $ \, \check\boldalpha_n \mapsto x_n \, $  ($ \, n
\in \N_\nu \, $),  \, in particular  $ {\big( {\calH_\h}^{\!\prime}
\big)}^{\!\vee} $  is the free associative  $ \Bbbk[\h\,] $--algebra
over  $ \big\{ \check{\boldalpha}_n \big\}_{n \in \N_\nu} $;  \, then
specialization yields a graded algebra isomorphism
  $$  \Omega \, \colon \, {\big( {\calH_\h}^{\!\prime} \big)}^{\!\vee}
{\Big|}_{\h=1} := \, {\big( {\calH_\h}^{\!\prime} \big)}^{\!\vee} \!
\Big/ (\h \! - \! 1) \, {\big( {\calH_\h}^{\!\prime} \big)}^{\!\vee}
\,{\buildrel \cong \over
{\lhook\joinrel\relbar\joinrel\relbar\joinrel\twoheadrightarrow}}\,
\calH \; ,  \qquad  {\check\boldalpha}_n{\big|}_{\h=1} \mapsto
\boldalpha_n \;\; .  $$
As for the Hopf structure, in  $ {\big( {\calH_\h}^{\!\prime}
\big)}^{\!\vee}{\Big|}_{\h=1} $  it is given by
  $$  \Delta \big( \check{\boldalpha}_n{\big|}_{\h=1} \big)  =
{\textstyle \sum_{k=2}^n} \, \check{\boldalpha}_k{\big|}_{\h=1}
\otimes \h^{n-k} {\Cal Q}^k_{n-k}(\boldalpha_*){\big|}_{\h=1}  +
{\textstyle \sum_{k=0}^{n-1}} \, {\check{\boldalpha}_1}^{\,k}{\big|}_{\h=1}
\otimes \h^{n-2} Z^k_{n-k}(\boldalpha_*){\big|}_{\h=1} \; .  $$
   \indent   As before, split  $ {\Cal Q}^k_{n-k}(\a_*) $  as
$ \, {\Cal Q}^k_{n-k}(\a_*) = \sum_s {\Cal T}^{s,k}_{n-k}
(\boldalpha_*) \, $,  \, and split each  $ {\Cal T}^{s,k}_{n-k}
(\widehat{\boldalpha}_*) $  into homogeneous components w.r.t.~the
total degree in the  $ \widehat{\boldalpha}_i $'s,  say  $ \,
{\Cal T}^{s,k}_{n-k}(\widehat{\boldalpha}_*) = {\textstyle \sum_r}
{\Cal Y}^{s,k}_{r,n}(\widehat{\boldalpha}_*) \, $:  \, then  $ \,
\h^{n-k-s} {\Cal T}^{s,k}_{n-k}(\widehat{\boldalpha}_*) = \h^{n-k-s}
{\textstyle \sum_r} {\Cal Y}^{s,k}_{r,n}(\widehat{\boldalpha}_*)
= {\textstyle \sum_r} \h^{n-k-s+r} {\Cal Y}^{s,k}_{r,n}
(\check{\boldalpha}_*) \, $,  \, because  $ \, \widehat{\boldalpha}_*
= \h \, \check{\boldalpha}_* \, $.  As  $ \, \h^{n-k-s+r}
{\Cal Y}^{s,k}_{r,n}(\check{\boldalpha}_*) \equiv {\Cal Y}^{s,k}_{r,n}
(\check{\boldalpha}_*) \! \mod (\h-1) \, {\big( {\calH_\h}^{\!\prime}
\big)}^{\!\vee} \, $,  \, we eventually get
  $$  \h^{n-k} {\Cal Q}^k_{n-k}(\a_*)  =  {\textstyle \sum_{s,r}}
\h^{n-k-s+r} {\Cal Y}^{s,k}_{r,n}(\check{\boldalpha}_*)  \equiv
{\textstyle \sum_{s,r}} {\Cal Y}^{s,k}_{r,n}(\check{\boldalpha}_*)
\hskip-5pt \mod \hskip-2pt (\h-1) \, {\big( {\calH_\h}^{\!\prime}
\big)}^{\!\vee}  \! =  Q^k_{n-k}(\a_*) \, .  $$
for all  $ k $  and  $ n \, $.  Similarly  $ \, \h^{n-1} Z^k_{n-k}
(\boldalpha_*) \equiv Z^k_{n-k}(\boldalpha_*) \hskip-1pt \mod (\h-1)
\, {\big( {\calH_\h}^{\!\prime} \big)}^{\!\vee} $  ($ \forall
\, k \, $,  $ n $).  Thus
  $$  \displaylines{
   \Delta \big( \check{\boldalpha}_n{\big|}_{\h=1} \big)  =
{\textstyle \sum_{k=2}^n} \, \check{\boldalpha}_k{\big|}_{\h=1}
\otimes \h^{n-k} {\Cal Q}^k_{n-k}(\boldalpha_*){\big|}_{\h=1}  +
{\textstyle \sum_{k=0}^{n-1}} \, {\check{\boldalpha}_1}^{\,k}{\big|}_{\h=1}
\otimes \h^{n-2} Z^k_{n-k}(\boldalpha_*){\big|}_{\h=1}  =  \cr
   {} \hfill   = \;  {\textstyle \sum_{k=2}^n} \,
\check{\boldalpha}_k{\big|}_{\h=1} \otimes {\Cal Q}^k_{n-k}
(\boldalpha_*){\big|}_{\h=1}  \, + \,  {\textstyle \sum_{k=0}^{n-1}}
\, {\check{\boldalpha}_1}^{\,k}{\big|}_{\h=1} \otimes Z^k_{n-k}
(\boldalpha_*){\big|}_{\h=1}  \; .  \cr }  $$
   \indent   On the other hand, one has  $ \; \Delta(\boldalpha_n)
\, = \,  {\textstyle \sum_{k=2}^n} \, \boldalpha_k \otimes {\Cal
Q}^k_{n-k}(\boldalpha_*) \, + \, {\textstyle \sum_{k=0}^{n-1}}
\, \boldalpha_1^{\,k} \otimes Z^k_{n-k}(\boldalpha_*) \, $  in 
$ \calH $,  thus the algebra isomorphism  $ \, \Omega \, \colon \, {\big(
{\calH_\h}^{\!\prime}}\big)^{\!\vee}{\Big|}_{\h=1} \, {\buildrel \cong
\over {\lhook\joinrel\relbar\joinrel\relbar\joinrel\twoheadrightarrow}}
\, \calH \, $  given by  $ \, {\widehat\boldalpha}_n{\big|}_{\h=1}
\mapsto \boldalpha_n \, $  also preserves the coproduct; similarly,
it also respects the antipode and the counit, hence it is a graded
Hopf algebra isomorphism.  In a nutshell, we have (as graded Hopf
$ \Bbbk $--algebras)
  $$  {\big({\calH_\h}^{\!\prime}\big)}^{\!\vee} \;{\buildrel
{\h \rightarrow 1} \over \llongrightarrow}\; \calH \; .  $$
 \eject   
   \indent   Therefore we have filled in the top part of the diagram
of deformations (5.5), corresponding to (5.4), for  $ \, H = \calH
\; (:= \calH_\nu) \, $:  \, it reads
  $$  F \big[ {\varGamma_{\!\L_\nu}\phantom{|}}^{\hskip-8pt \star}
\big]  \, = \,  {\calH_\h}^{\!\prime}{\Big|}_{\h=0}
\underset{{\calH_\h}^{\!\prime}} \to {\overset{0 \leftarrow \h
\rightarrow 1} \to{\longleftarrow\joinrel\llongrightarrow}}
{\calH_\h}^{\!\prime}{\Big|}_{\h=1}  = \,  \calH  \, = \,
{\big( {\calH_\h}^{\!\prime} \big)}^{\!\vee}{\Big|}_{\h=1}
\underset{{({\calH_\h}^{\!\prime})}^{\!\vee}} \to {\overset{1 \leftarrow
\h \rightarrow 0} \to {\longleftarrow\joinrel\llongrightarrow}}
{\big( {\calH_\h}^{\!\prime} \big)}^{\!\vee}{\Big|}_{\h=0}
= \,  U(\L_\nu)  $$
or simply  $ \; F \big[ {\varGamma_{\!\L_\nu}\phantom{|}}^{\hskip-8pt
\star} \big] \underset{{\calH_\h}^{\!\prime}}
\to {\overset{0 \leftarrow \h \rightarrow 1} \to
{\longleftarrow\joinrel\relbar\joinrel\relbar\joinrel\llongrightarrow}}
\, \calH \underset{{({\calH_\h}^{\!\prime})}^{\!\vee}}  \to
{\overset{1 \leftarrow \h \rightarrow 0} \to
{\longleftarrow\joinrel\relbar\joinrel\relbar\joinrel\llongrightarrow}}
U(\L_\nu) \; $,  \; where  $ \L_\nu $  is given the Lie bialgebra
structure of Theorem 10.18/20 and  $ {\varGamma_{\!\L_\nu}
\phantom{|}}^{\hskip-8pt \star} $  is the corresponding dual
Poisson group mentioned in Theorem 10.18.  Thus  $ \, \calH \, $
{\sl is intermediate between the (Poisson-type) ``geometrical
symmetries''  $ F \big[ {\varGamma_{\!\L_\nu}\phantom{|}}^{\hskip-8pt
\star} \big] $  and  $ U(\L_\nu) $},  so their geometrical meaning
should shed some light on it; conversely, the physical meaning of
$ \calH $  should have some reflect on the physical meaning of both
$ F \big[ {\varGamma_{\!\L_\nu} \phantom{|}}^{\hskip-8pt \star}
\big] $  and  $ U(\L_\nu) \, $.

\vskip7pt

   {\sl $ \underline{\hbox{\it Remark}} $:}  \; The analysis in \S 10.9
and \S 10.21 yields a complete description of the deformation features
of  $ \calH $  via  $ \calH_\h $  and Drinfeld's functors drawn in
(5.5).  In particular
  $$  G_+ = {\varGamma_{\!\L_\nu} \phantom{|}}^{\hskip-8pt \star} \; ,
\quad  \gerk_+ \equiv \gerg_+^{\,\times} = \big(\L_\nu, \delta_*\big)
\; ,  \qquad  \gerg_- = \big(\L_\nu, \delta_\bullet\big) \; ,  \quad
K_- \equiv G_-^\star = {G_{\!\L_\nu}\phantom{|}}^{\hskip-8pt \star}
\eqno (10.9)  $$
(as  $ \, \Char(\Bbbk) = 0 \, $) where  $ \, \delta_* \, $  and  $ \,
\delta_\bullet \, $  denote the Lie cobracket on  $ \L_\nu $  defined
respectively in Theorems 10.18/20 and in Theorem 10.6.  Next result
shows that the four objects in (10.9) are really different, though
they share some common feature.

\vskip7pt

\proclaim{Theorem 10.22}
                                         \hfill\break
   \indent   (a)  $ \, \big( {\calH_\h}^{\!\vee} \big)' \cong
{\calH_\h}^{\!\prime} \, $  as Poisson  $ \, \Bbbk[\h\,] $--algebras,
but  $ \, \big( {\calH_\h}^{\!\vee} \big)' \not\cong {\calH_\h}^{\!
\prime} \, $  as Hopf  $ \, \Bbbk[\h\,] $--algebras.
                                         \hfill\break
   \indent   (b)  $ \, \big(\L_\nu, \delta_*\big) \cong \big(\L_\nu,
\delta_\bullet\big) \, $  as Lie algebras, but  $ \, \big(\L_\nu,
\delta_*\big) \not\cong \big(\L_\nu, \delta_\bullet\big) \, $  as
Lie bialgebras.
                                         \hfill\break
   \indent   (c)  $ \, {G_{\!\L_\nu}\phantom{|}}^{\hskip-8pt \star}
\cong {\varGamma_{\!\L_\nu}\phantom{|}}^{\hskip-8pt \star} \, $
as (algebraic) Poisson varieties, but  $ \, {G_{\!\L_\nu}
\phantom{|}}^{\hskip-8pt \star} \not\cong {\varGamma_{\!\L_\nu}
\phantom{|}}^{\hskip-8pt \star} \, $  as (algebraic) groups.
                                        \hfill\break
  \indent   (d) \,  The analogues of statements (a)--(c) hold
with  $ \calK $  and  $ \L_\nu^+ $  instead of  $ \, \calH $
and  $ \L_\nu \, $.
\endproclaim

\demo{Proof}  It follows from  Theorem 10.8{\it (a)\/}  that
$ {\big( {\calH_\h}^{\!\vee} \big)}' $  can be seen as a Poisson
Hopf algebra, with Poisson bracket given by  $ \, \{x,y\,\} := \h^{-1}
[x,y] = \h^{-1} (x \, y - y \, x) \, $  (for all  $ \, x $,  $ y
\in {\big( {\calH_\h}^{\!\vee} \big)}' \, $);  then  $ {\big(
{\calH_\h}^{\!\vee} \big)}' $  is the free Poisson algebra generated
by  $ \Big\{\, \widetilde{\b}_{x_n} \! = \widetilde{\x}_n = \a_n
\,\Big|\, n \in \N \,\Big\} \, $;  \, since  $ \, \a_n = \boldalpha_n +
(1 - \delta_{1,n}) \, {\boldalpha_1}^{\!n} \, $  and  $ \, \boldalpha_n
= \a_n - (1 - \delta_{1,n}) \, {\a_1}^{\!n} \, $  ($ \, n \in \N_+
\, $)  it is also (freely) Poisson-generated by  $ \big\{\,
\boldalpha_n \,\big|\, n \in \N \,\big\} $.  We also saw that
$ {\calH_\h}^{\!\prime} $  is the free Poisson algebra over  $ \,
\big\{\, \widehat{\boldalpha}_n \,\big|\, n \in \N \,\big\} \, $;
\, thus mapping  $ \, \boldalpha_n \mapsto \widehat{\boldalpha}_n
\, $  ($ \, \forall \, n \in \N \, $)  does define a unique Poisson
algebra isomorphism  $ \; \Phi \, \colon \, {\big( {\calH_\h}^{\!\vee}
\big)}' \,{\buildrel \cong \over \longrightarrow}\, {\calH_\h}^{\!
\prime} \, $,  \, given by  $ \; \widetilde{\boldalpha}_b := \h^{-d(b)}
\boldalpha_b \mapsto \widehat{\boldalpha}_b \, $,  \; for all  $ \,
b \in B_\nu \, $.  This proves the first half of  {\it (a)},  and
then also (taking semiclassical limits and spectra) of  {\it (c)\/}.
                                               \par
   The group structure of either  $ {G_{\!\L_\nu}
\phantom{|}}^{\hskip-8pt \star} $  or  $ {\varGamma_{\!\L_\nu}
\phantom{|}}^{\hskip-8pt \star} $  yields a Lie cobracket onto
the cotangent space at the unit point of the above, isomorphic
Poisson varieties: this cotangent space identifies with
$ \L_\nu $,  and the two cobrackets are given respectively
by  $ \, \delta_\bullet(x_n) = \sum_{\ell=1}^{n-1} (\ell+1)
\, x_\ell \wedge x_{n-\ell} \, $  for  $ {G_{\!\L_\nu}
\phantom{|}}^{\hskip-8pt \star} $  (by Theorem 10.8) and
by  $ \, \delta_*(x_n) = (n-2) \, x_{n-1} \wedge x_1 \, $  for
$ {\varGamma_{\!\L_\nu}\phantom{|}}^{\hskip-8pt \star} $  (by
Theorem 10.18), for all  $ \, n \in \N_\nu \, $.  It follows that
$ \, \text{\sl Ker}\,(\delta_\bullet) = \{0\} \not= \text{\sl Ker}\,
(\delta_*) \, $,  \, which implies that the
two Lie coalgebra structures on  $ \L_\nu $  are not isomorphic.
This proves  {\it (b)},  and also means that  $ \, {G_{\!\L_\nu}
\phantom{|}}^{\hskip-8pt \star} \not\cong {\varGamma_{\!\L_\nu}
\phantom{|}}^{\hskip-8pt \star} \, $  as (algebraic) groups, hence
$ \, F \big[ {G_{\!\L_\nu} \phantom{|}}^{\hskip-8pt \star} \big]
\not\cong F \big[ {\varGamma_{\!\L_\nu} \phantom{|}}^{\hskip-8pt
\star} \big] \, $  as Hopf  $ \Bbbk $--algebras,  and so  $ \,
\big( {\calH_\h}^{\!\vee} \big)' \not\cong {\calH_\h}^{\!\prime}
\, $  as Hopf  $ \Bbbk[\h\,] $--algebras,  which ends the
proof of  {\it (c)\/}  and  {\it (a)\/}  too.
                                           \par
   Finally, claim  {\it (d)\/}  should be clear: one applies
the like arguments  {\sl mutatis mutandis},  and everything
follows as before.   \qed
\enddemo

\vskip7pt

   {\bf 10.23 Generalizations.} \, Plenty of features of  $ \, \calH
= \calH^\dif \, $  are shared by a whole bunch of  {\sl graded\/}
Hopf algebras, which usually arose in connection with some physical
problem or some (co)homological construction, and all bear a nice
combinatorial content; essentially, most of them can be described
as ``formal series'' over indexing sets   --- replacing  $ \N $  ---
of various (combinatorial) nature: planar trees (with or without
labels), forests, graphs, Feynman diagrams, etc.  Besides the
ice-breaking examples given by Connes and Kreimer (cf.~[CK1--3]),
which are all commutative or cocommutative Hopf algebras, other
non-commutative non-cocommutative examples (like the one of
$ \calH^\dif \, $)  are introduced in [BF1--2], roughly through a
``disabelianization process'' applied to the commutative Hopf
algebras of Connes and Kreimer.  The most general analysis and
wealth of examples in this context is due to Foissy (see [Fo1--3]),
who also makes   --- in other terms ---   an interesting (although
less deep than ours) study of the operators  $ \delta_n $'s  and of
the functor  $ \, H \mapsto H' \, $  ($ H \in \HA_\Bbbk $).  Other
examples, issued out of topological motivations, can be found in
the works of Loday et al.: see e.g.~[LR], and references therein.
                                             \par
   When performing the like analysis, as we did for  $ \calH $,
for a graded Hopf algebra  $ H $  of the afore mentioned type,
{\sl the arguments used for  $ \calH $  apply essentially the
same, up to minor changes, and give much the same results.}
To give an example, the Hopf algebras considered by Foissy
are non-commutative polynomial, say  $ \, H = \Bbbk \big\langle
\{x_i\}_{i \in \Cal{I}} \big\rangle \, $  for some index set
$ \Cal{I} \, $:  \, then one finds  $ \, \widehat{H} = {\calH_\h}^{\!
\vee}{\big|}_{\h=0} = U(\gerg_-) = U(\L_{\Cal{I}}) \, $  where
$ \L_{\Cal{I}} $  is the free Lie algebra over  $ \Cal{I} \, $.
                                             \par
   This opens the way to apply the crystal duality principle
to all these graded Hopf algebras of great interest for their
applications in mathematical physics or in topology (or whatever),
with the simplest case of  $ \calH^\dif $  playing the role of
a toy model which realizes a clear and faithful pattern for
many common features of all Hopf algebras of this kind.

\vskip2,3truecm

\vskip83pt

\vfill

\enddocument